\documentclass[10pt]{article}
\usepackage{times}
\usepackage[ruled,vlined]{algorithm2e}
\usepackage{cl2emono}
\usepackage{amsfonts, amssymb, amsmath}
\usepackage{geneva}
\usepackage{placeins}
\usepackage{imakeidx}
\usepackage{float}
\usepackage[hidelinks]{hyperref}
\usepackage{bibentry}
\usepackage{graphicx}
\usepackage{psfrag}
\usepackage{epsfig}
\usepackage{comment}
\usepackage{color}
\newcommand{\C}{\mathbb C}
\newcommand{\R}{\mathbb R}

\newcommand{\cM}{\mathcal{M}}
\newcommand{\cS}{\mathcal{S}}
\newcommand{\cR}{\mathcal{R}}

\newcommand{\V}{{\mathcal V}}
\def\e{{\rm e}}

\def\s{{\varphi_\star}}
\def\eps{\varepsilon}
\def\oeps{\eps_\star}

\def\bigo{{\mathcal O}}

\def\iu{\mathrm{i}}
\newcommand{\x}{y}
\newcommand{\xx}{w}
\newcommand{\y}{x}

\renewcommand{\diag}{{\rm diag}}
\newcommand\pin{^\dagger}
\newcommand{\F}{\ophi}

\newcommand{\Lameps}{\Lambda_\varepsilon}
\newcommand{\aleps}{\alpha_\varepsilon}
\newcommand{\rhoeps}{\rho_\varepsilon}

\renewcommand{\Re}{{\mbox{\rm Re}}}
\renewcommand{\Im}{{\mbox{\rm Im}}}

\newcommand{\matlab}{{\sc Matlab}}

\newcommand{\Id}{{I}}

\newcommand{\conj}[1]{\overline{#1}}

\newcommand{\goes}{\rightarrow}
\newcommand{\clambda}{\conj{\lambda}}
\newcommand{\la}{\langle}
\newcommand{\ra}{\rangle}

\newcommand{\tp}{\top}
\def\wt{\widetilde}
\def\wh{\widehat}
\newcommand{\bcl}{\color{black}}
\newcommand{\bcltwo}{\color{black}}

\newcommand{\ecl}{\color{black}}
\newcommand{\ecltwo}{\color{black}}

\newcommand{\bng}{\color{black}}
\newcommand{\eng}{\color{black}}
\newenvironment{dedication}
  {\clearpage           
   \thispagestyle{empty}
   \vspace*{\stretch{1}}
   \itshape             
   \raggedleft          
  }
  {\par 
   \vspace{\stretch{3}} 
   \clearpage           
  }

\makeindex

\begin{document}
\title{Matrix nearness problems and
eigenvalue optimization}

%

\author{Nicola Guglielmi and Christian Lubich\\[70mm]
\large Preprint version}

\maketitle
\thispagestyle{empty}

\renewcommand\thepage{\roman{page}}
\setcounter{page}{5}

\begin{dedication}
Dedicated to our children.
\end{dedication}

\markboth{Preface}{Preface}

\noindent
\bf\Large
Preface  
\rm\normalsize

\bigskip\bigskip
\begcitation{A matrix nearness problem consists of finding, for an arbitrary matrix $A$, a nearest member of some given class of matrices, where distance is measured in a matrix norm.}
\signed{(Nicholas J. Higham, \cite{Hig89})}
\endcitation

\bigskip\noindent
This book is about solving matrix nearness problems that are related to eigenvalues or singular values or pseudospectra. These problems arise in great diversity in various fields, be they related to dynamics, as in questions of robust stability and robust control, or 
related to graphs, as in questions of clustering and ranking. Algorithms for such problems work with matrix perturbations that drive eigenvalues or singular values or Rayleigh quotients to desired locations. 

Remarkably, the optimal perturbation matrices are typically of rank one or are projections of rank-1 matrices onto a linear structure, e.g.~a prescribed sparsity pattern. In the approach worked out here, these optimal rank-1 perturbations will be determined in a two-level iteration \bcl that consists of eigenvalue optimization and root-finding. An eigenvalue optimization problem with equality or inequality constraints on the perturbation size is to be solved via gradient-based rank-1 matrix differential equations. This amounts to numerically driving a rank-$1$ matrix, which is represented by two vectors, to a stationary point. 
The root-finding part determines the optimal perturbation size by solving a scalar nonlinear equation. The two algorithmic parts can be nested or can alternate.\ecl 

A wide variety of matrix nearness problems, as outlined in the introductory Chapter~I, will be tackled by such an approach and its nontrivial extensions. In Chapter~II, we study a basic eigenvalue optimization problem and its numerical solution via rank-1 matrix differential equations, which are norm- and rank-1 constrained gradient systems. In Chapter~III, this approach yields algorithms for computing extremal points and boundary curves of pseudospectra. \bcl In Chapter IV, we present algorithms for matrix nearness problems, which in a nested or alternating way combine the rank-1 eigenvalue optimization algorithms of Chapter~II with a Newton--bisection method. \ecl This is illustrated by the problem of computing a nearest unstable complex matrix to a given stable matrix. In Chapter~V, the rank-1 approach is extended to nearness problems for matrices with a prescribed complex-linear or real-linear structure, e.g.~real matrices or matrices with a given sparsity pattern or a Toeplitz or Hankel structure. The approach is applied to various structured matrix nearness problems, which include finding the nearest singular structured matrix to a given invertible structured matrix and, for asymptotically stable linear differential equations, finding the largest norm of structured perturbations of the matrix that still ensures a prescribed transient bound. In Chapter~VI, we propose and analyze algorithms for matrix nearness problems that go beyond those of Chapters IV and~V, among them matrix stabilization and finding a nearest defective real or complex matrix. In~Chapter VII, we discuss algorithms for exemplary nearness problems in the area of systems and control, and in Chapter VIII we rephrase clustering and ranking problems from graph theory as structured matrix nearness problems and extend our algorithmic approach to such graph problems.

In Chapters II to VIII, the inclusion of references to the existing literature is done in Notes at the end of each chapter, and only exceptionally  we add references to the running text. References are not numbered, but are addressed by the names of authors and the year of publication, for example Higham~(\cite{Hig89}) and Lewis \& Overton~(\cite{LeO96}).

\matlab\ codes implementing the algorithms presented in this book are freely
available for non-commercial, academic purposes and can be obtained upon request from the authors. 

\medskip
We are grateful for motivating and critical -- in any case helpful -- comments and discussions with students and colleagues, whom we name here in alphabetical order:
Eleonora Andreotti,
Michele Benzi,
Paolo Butt\`a,
Gianluca Ceruti,
Dominik Edelmann,
Mark Embree,
Antonio Fazzi,
Nicolas Gillis, 
Miryam Gnazzo, 
Stefano Grivet Talocia,
Mert G\"{u}rb\"{u}zbalaban,
Ernst Hairer,
Des Higham,
Michael Karow,
Bal\'azs Kov\'acs,
Daniel Kressner,
Caroline Lasser,
Yoann Le Hénaff,
Maria Lopez Fernandez,
Manuela Manetta,
Mattia Manucci,
Ivan Markovsky, 
Volker Mehr\-mann,
Emre Mengi,
Wim Michiels,
Tim Mitchell,
J\"org Nick, 
Vanni Noferini,
Silvia Noschese,
Michael Overton, 
Federico Poloni,
Vladimir Protasov,
Anton Savostianov,
Carmen Scalone, 
Stefano Sicilia,
Valeria Simoncini,
Pete Stewart,
Nick Trefethen,
Luca Trevisan,
Francesco Tu\-di\-sco,
Bart Vandereycken,
Matthias Voigt.

We express our special thanks to Michael Overton who first introduced us to eigenvalue optimization and matrix nearness problems.

\bigskip\medskip
\noindent
L'Aquila and T\"ubingen, May 2026 


\tableofcontents



\setcounter{page}{1}
\renewcommand\thepage{\arabic{page}}

\chapter{Introduction by examples}
\label{chap:intro}

\rm\normalsize

\noindent
\begcitation{All happy families resemble one another, each unhappy family is unhappy in its own way.}
\signed{(Lev Tolstoy, \cite{T1877})}
\endcitation
\noindent
Analogous to the quote, simple matrix nearness problems tend to look alike, while difficult ones are each difficult in their own distinct way.
In this short introductory chapter we present some of the matrix nearness problems that will be considered in more detail later in this book. We formulate the problems, give some background, and sketch ideas for numerical approaches that will be worked out in later chapters.

\index{matrix nearness problem}

\subsection*{Nearest singular matrix, distance to singularity}
\begcitation{It will be shown that if the least-squares criterion of approximation
be adopted, this problem has a general solution which is relatively
simple in a theoretical sense, though the amount of numerical work
involved in applications may be prohibitive.}
\signed{(Carl Eckart and Gale Young, \cite{EY36})}
\endcitation
\noindent
We begin with the matrix nearness problem that may be most basic and have the longest history. The problem is atypical in that its solution can be directly read off from a matrix factorization, which here is the singular value decomposition.

\medskip\noindent
{\bf Problem.} {\it Given an invertible matrix, find a nearest singular matrix.}

\medskip\noindent
The term ``nearest'' refers to a metric given by a matrix norm. The problem can thus be rephrased as follows when general complex $n\times n$ matrices are considered: Given an invertible matrix $A\in \C^{n,n}$, find a perturbation $\Delta\in\C^{n,n}$ of minimal norm such that $A+\Delta$ is not invertible. The norm of the minimal perturbation $\Delta$ is called the {\it distance to singularity} of the matrix $A$. 

The solution to the problem depends on the choice of matrix norm. Here we consider, for a matrix $M\in \C^{n,n}$ with singular values $\sigma_1\ge \sigma_2 \ge \ldots \ge \sigma_n\ge 0$,
the matrix 2-norm 
$$
\| M \|_2 = \sup_{0\ne x \in\C^n} \frac{\|Mx\|_2}{\|x\|_2} = \sigma_1
$$
(where the 2-norm of vectors is the Euclidean norm)
and the Frobenius norm, which is the Euclidean norm of the vector of matrix entries,
$$
\| M \|_F = \biggl( \sum_{i,j=1}^n |m_{ij}|^2 \biggr)^{1/2} = \biggl( \sum_{i=1}^n \sigma_i^2 \biggr)^{1/2}.
$$

\noindent
For these two norms, the nearest singular matrix to the given invertible matrix $A\in \C^{n,n}$ is obtained directly 
from the singular value decomposition of $A$,
$$
A = \sum_{i=1}^n \sigma_i u_i v_i^*
$$
(with the conjugate transpose $v^*=\overline v^\top$), where $u_i\in\C^n$ and $v_i\in\C^n$ are the left and right singular vectors, respectively, to the $i$th singular value $\sigma_i>0$ of $A$. A nearest singular matrix to $A$ is found by truncating the last term in the sum, i.e.~by adding the {rank-1} perturbation matrix $\Delta=-\sigma_n u_n v_n^*$ to $A$, and the distance to singularity thus equals the smallest singular value $\sigma_n$, both for the matrix 2-norm and the Frobenius norm. This result is the special case of rank $r=n-1$ of the Eckart--Young theorem (\cite{EY36}), going back in its origins to Schmidt (\cite{Schm1907}). The theorem states that a nearest matrix to $A$ of rank~$r$ is obtained by truncating the singular value decomposition after $r$ terms. This holds true for both the matrix 2-norm and the Frobenius norm, and in fact for every unitarily invariant matrix norm, as shown by Mirsky (\cite{Mir60}). For a general matrix norm induced by a vector norm,  the distance to singularity equals $1/\|A^{-1}\|$, i.e.~the relative distance is the reciprocal of the condition number cond$(A)=\|A\|\cdot\|A^{-1}\|$; moreover, a minimal perturbation matrix is still of rank~1. Kahan (\cite{Kah66}) attributes this result to Gastinel; see also Wilkinson (\cite{Wil86}).

\subsubsection*{Structured problem.}
While the distance to singularity with respect to a unitarily invariant matrix norm is obtained directly from the singular value decomposition for a general complex or real  invertible matrix $A$, the situation changes when $A$ is a {\it structured} invertible matrix and we ask for the nearest singular matrix with the same given linear structure, for example having the same sparsity pattern or belonging to the class of Hamiltonian or Toeplitz or Hankel matrices. The nearest singular structured matrix cannot be determined from the singular value decomposition of $A$, and the structured distances to singularity are in general not the same for the matrix 2-norm and the Frobenius norm. The minimal structured perturbation $\Delta$ is no longer a rank-1 matrix in general, but for the Frobenius norm we will show in Section~\ref{sec:sing-S} that a minimal $\Delta$ is the orthogonal projection of a rank-1 matrix onto the structure. 

This rank-1 property leads us to propose a two-level iterative algorithm, where we solve a system of differential equations for two vectors depending on a distance parameter in the inner iteration -- which requires only matrix-vector products with structured matrices and inner products of vectors -- and solve a scalar nonlinear equation for the structured distance to singularity in the outer iteration.

In a numerical example in Section~\ref{sec:sing-S} we use this algorithm to compute the distance of a given pair of polynomials to a pair of polynomials with a common zero. This distance equals the structured distance to singularity of a Sylvester matrix, and the coefficients of the nearest pair of polynomials with a common zero can be read off from the nearest singular Sylvester matrix.

\vfill
\pagebreak

\subsection*{Nearest unstable matrix, distance to instability}
\begcitation{\hfill How~near~is~a~stable~matrix~to~an~unstable~matrix?}
\signed{(Charles Van Loan, \cite{VL85})}
\endcitation

\noindent
Many matrix nearness problems considered in this book are motivated by dynamical systems. In particular, the question posed by Van Loan in the title of his 1985 paper has received much attention in the literature; see the Notes at the end of Chapters \ref{chap:pseudo} and~\ref{chap:two-level}. A linear differential equation $\dot x(t)=Ax(t)$ is asymptotically stable, i.e.~every solution $x(t)$ tends to $0$ as $t\to\infty$, if and only if all eigenvalues of the matrix $A$ have negative real part. Van Loan's question addresses the robustness of asymptotic stability under perturbations of $A$. 

\medskip\noindent
{\bf Problem.} {\it Given a matrix having all eigenvalues in the open complex left half-plane, find a nearest matrix having some eigenvalue on the imaginary axis.}

\medskip\noindent
Various algorithms have been proposed in the literature to solve this problem for general complex matrices (see again the Notes of Chapters \ref{chap:pseudo} and~\ref{chap:two-level}). The problem can be rephrased in terms of the $\eps$-pseudospectrum of $A\in\C^{n,n}$, which for $\eps> 0$ is defined as
$$
\Lambda_\eps(A) = \{ \lambda \in \C \,:\, 
\lambda \in \Lambda(A+\Delta) \text{ for some $\Delta\in \C^{n,n}$ with $\|\Delta\|\le \eps$} \},
$$
where $\Lambda(A)$ is the spectrum of $A$, that is, the set of eigenvalues of $A$. Equivalently, $\Lambda_\eps(A)$ is the set of all complex $\lambda$ for which the distance to singularity of $A-\lambda I$ is at most $\eps$.
In the following, the norm is chosen as the Frobenius norm or the matrix 2-norm. Both yield the same $\eps$-pseudospectrum, since the distance to singularity is the same for both norms.

The problem is then to find the smallest $\eps> 0$
with  $\Lambda_\eps(A)\cap \iu \R\ne \emptyset$. This optimal $\eps$ is called the {\it distance to instability} or {\it stability radius} of $A$.
A minimal perturbation $\Delta\in\C^{n,n}$ with $\Lambda(A+\Delta)\cap \iu \R\ne \emptyset$ is known to have rank~1. This fact has been used  to derive numerical algorithms for the above matrix nearness problem; see Section~\ref{sec:psa} and the Notes of Chapters~\ref{chap:pseudo} and~\ref{chap:two-level}. 

\subsubsection*{Structured problem.} Consider now the same problem posed for real matrices, i.e. finding $\Delta\in\R^{n,n}$ of minimal norm having $\Lambda(A+\Delta)\cap \iu \R\ne \emptyset$, or, more generally, posed for perturbation matrices $\Delta$ in a complex-linear or real-linear subspace $\cS\subset\C^{n,n}$ (the {\it structure space}), e.g. complex or real matrices with a given sparsity pattern or special matrix classes, often named after 19th-century mathematicians. The $\cS$-structured $\eps$-pseudo\-spectrum $\Lambda_\eps^\cS(A)$, where $\C^{n,n}$ in the definition of $\Lambda_\eps(A)$ is replaced by the subspace $\cS$, 
is usually not the same for the Frobenius norm and the matrix 2-norm. With respect to both norms,
the {\it $\cS$-structured distance to instability}, which is the smallest $\eps>0$ such that $\Lambda_\eps^\cS(A)\cap \iu \R\ne \emptyset$, is in general strictly larger than the complex-unstructured distance to instability, including in the real case $\cS=\R^{n,n}$.

A minimal structured perturbation $\Delta\in\cS$ is no longer a rank-1 matrix in general, but for the Frobenius norm we will show in Section~\ref{sec:proto-structured} that a minimal $\Delta\in\cS$ is the orthogonal projection of a rank-1 matrix onto the structure space $\cS$. 
(For the real case $\cS=\R^{n,n}$, this is the real part of a complex rank-1 matrix.) 
Based on this rank-1 property, we propose a two-level algorithm, where we solve a system of differential equations for two vectors depending on a distance parameter in the inner iteration -- which requires computing eigenvalues and eigenvectors of perturbed matrices -- and solve a scalar nonlinear equation for the structured distance to instability in the outer iteration. 

\subsection*{Robustness of transient bounds}
\begcitation{This phenomenon has traditionally been investigated by linearizing the equations
of flow and testing for unstable eigenvalues of the linearized problem, but the results of such
investigations agree poorly in many cases with experiments. Nevertheless, linear effects
play a central role in hydrodynamic instability.}
\signed{(Lloyd N. Trefethen, Anne E. Trefethen, Satish C. Reddy,}
\signed{and Tobin A. Driscoll, \cite{TreTRD93})}
 \endcitation
 \noindent
As has been emphasized by Trefethen et al. (\cite{TreTRD93}, \cite{TreE05}), eigenvalues of non-normal matrices $A$ give no indication of transient bounds of solutions of linear differential equations with $A$, whereas pseudospectra do. This relies on the characterization of the $\eps$-pseudospectrum in terms of resolvent norms, 
\begin{align*}
\Lambda_\eps(A) &=  \{ \lambda \in \C \,:\, 
\lambda \in \Lambda(A+\Delta) \text{ for some $\Delta\in \C^{n,n}$ with $\|\Delta\|\le \eps$} \}
\\
&= \{ \lambda \in \C \,:\, 
\| (A-\lambda I)^{-1} \|_2 \ge \eps^{-1} \},
\end{align*}
which is valid for both the matrix 2-norm and the Frobenius norm as matrix norm $\|\cdot\|$.
This fundamental equality was apparently first noted by Wilkinson (\cite{Wil86}), who ``would like to emphasize the sheer economy of this theorem''; see also the Notes of Chapter~\ref{chap:pseudo}. The result follows directly from the fact that the distance to singularity of the shifted matrix $A-\lambda I$ is equal to $1/\| (A-\lambda I)^{-1} \|_2$ in the cases of both the matrix 2-norm and the Frobenius-norm pseudospectra.

If $\Lambda_\eps(A)$ lies in the closed complex left half-plane, resolvent norms can be used to derive bounds reciprocal to $\eps$ of solutions of homogeneous and inhomogeneous linear differential equations $\dot x(t)=Ax(t)+f(t)$ that are valid for all times $t\ge 0$; see the Notes of Chapter~\ref{chap:pseudo}. To assess the robustness of such transient bounds under perturbations of~$A$, we therefore consider the following pseudospectral version of the previous problem. Here, $\eps$-pseudoeigenvalues are the points in the $\eps$-pseudospectrum $\Lambda_\eps(A)$.

\medskip\noindent
{\bf Problem.} {\it Given a matrix having all $\eps$-pseudoeigenvalues in the open complex left half-plane, find a nearest matrix having some $\eps$-pseudoeigenvalue on the imaginary axis.}

\medskip\noindent
The task here is to find $\Delta\in\C^{n,n}$ of minimal norm having $\Lambda_\eps(A+\Delta)\cap \iu \R\ne \emptyset$.
This problem is readily solved in terms of the distance to instability $\oeps$ of $A$: for $\eps<\oeps$, the distance is simply $\oeps-\eps$ and the minimum perturbation matrix $\Delta$ equals $(1-\eps/\oeps)$ times the minimum perturbation matrix in the nearest unstable matrix problem.

\pagebreak\noindent
{\bf Structured problem.} 

\medskip\noindent
\begcitation{The figure shows that
the unstructured spectral value set can be a misleading indicator of the robustness of
stability.}
\signed{(Diederich~Hinrichsen and Anthony~J.~Pritchard, \cite{HinP05})}
\vskip 1mm
\begcitation{By contrast, the eigenvalues that arise from structured perturbations do not bear as close a relation to the resolvent norm and may not provide much information about matrix behavior.}
\signed{(Lloyd~N.~Trefethen and Mark~Embree, \cite{TreE05})}
\endcitation
\noindent
Consider now the same problem posed for real matrices, i.e. finding $\Delta\in\R^{n,n}$ of minimal norm having $\Lambda_\eps(A+\Delta)\cap \iu \R\ne \emptyset$, or, more generally, posed for perturbation matrices $\Delta$ in a complex-linear or real-linear structure space $\cS\subset\C^{n,n}$.
Here, the situation is more intricate.  The key observation is that pseudoeigenvalues that arise from structured perturbations do bear a close relation to the resolvent norm.
Given $\eps<\oeps$, we introduce the joint unstructured--structured pseudo\-spectrum \bcl as a means to bridge the apparent discrepancy between the two quotes,\ecl
\begin{align*}
\Lambda_{\delta,\eps}^\cS(A) &= \{ \lambda \in \C \,:\, 
\lambda \in \Lambda_\eps(A+\Delta) \text{ for some  $\Delta\in \cS$ with $\|\Delta\|_F \le \delta$} \}
\\
&= \{ \lambda \in \C \,:\, 
\lambda \in \Lambda(A+\Delta+\Theta) \text{ for some  $\Delta\in \cS$ with $\|\Delta\|_F \le \delta$ }
\\
&\hspace{5.2cm} \text{and $\Theta\in \C^{n,n}$ with $\|\Theta\|_F\le\eps$} \}
\\
&= \{ \lambda \in \C \,:\, 
\|(A+\Delta-\lambda I)^{-1}\|_2 \ge \eps^{-1} \text{ for some  $\Delta\in \cS$ with $\|\Delta\|_F \le \delta$} \}.
\end{align*}
We ask for the smallest $\delta>0$ such that $\Lambda_{\delta,\eps}^\cS(A)\cap\iu\R\ne\emptyset$. We call this minimal $\delta$ the {\em $\cS$-structured $\eps$-stability radius} and denote it $\delta_\eps$. Then, transient bounds reciprocal to $\eps$ are  valid uniformly for all structured perturbations $\Delta\in \cS$ with $\|\Delta\|_F \le \delta_\eps$; see Section~\ref{sec:transient-S}.  A minimal perturbation $\Delta\in \cS$ again turns out to be the orthogonal projection onto $\cS$ of a rank-1 matrix.
A two-level algorithm for computing the $\cS$-structured $\eps$-stability radius is based on a {rank-1} matrix differential equation in the inner iteration and is given in Section~\ref{sec:eps-stab}. It takes essentially the same computational steps as the rank-1 algorithm for computing the unstructured stability radius. See also the Notes of Chapter~\ref{chap:struc}.


\subsection*{Matrix stabilization}
\begcitation{The question is to find the smallest perturbation that stabilizes
a given unstable matrix $A$, or, equivalently, to find the closest
stable matrix $X$ to a given unstable matrix $A$. This kind of problem
occurs in system identification when one needs to identify a
stable system from observations.}
\signed{(Fran\c cois-Xavier Orbandexivry, Yurii Nesterov, and Paul Van Dooren, \cite{OrbNVD13})}
\endcitation
\noindent
How near is an unstable matrix to a stable matrix? This formal converse of Van Loan's question has turned out to be challenging from a computational perspective. The problem arises in system identification as well as in model reduction of a stable system, which can yield a smaller but unstable matrix. The question then is how to stabilize it.

\medskip\noindent
{\bf Problem.} {\it Given a matrix having some eigenvalues of positive real part, find a nearest matrix having all eigenvalues of non-positive real part (or below a negative threshold).}

\medskip\noindent
A variety of numerical algorithms for this problem in the general complex case have been proposed in the literature; see the Notes of Chapter~\ref{chap:mnp-mix}. For the Frobenius-norm setting, we describe in Section~\ref{sec:mat-stab} two complementary new approaches,  which can also be used for real and structured versions of the problem. 
The problem can be made more meaningful for transient behaviour -- and potentially less computationally challenging -- by replacing eigenvalues with $\eps$-pseudoeigen\-values in the problem statement.

\bigskip\medskip
\begcitation{\hfill Se non \`e vero, \`e molto ben trovato.\\}
\signed{(Giordano Bruno, \cite{Bru1585})}
\endcitation
\noindent
Contrary to computing the distance to instability considered before, different algorithms tend to compute different stabilized matrices with different distances to the given unstable matrix; see the numerical example in Section~\ref{sec:mat-stab}. Rigorously speaking, the algorithms only yield upper bounds of the distance to stability via locally minimal stabilizing perturbations. Nonetheless, providing a fairly tight bound together with a stabilized matrix can be useful when one considers that the problem posed is a non-smooth, non-convex optimization problem that can have many local minima.

\subsection*{Nearest correlation matrix}
\begcitation{Given a symmetric matrix what is the nearest correlation matrix, that is, the
nearest symmetric positive semidefinite matrix with unit diagonal? This prob-
lem arises in the finance industry, where the correlations are between stocks.}
\signed{(Nicholas J. Higham, \cite{Hig02})}
\endcitation

\noindent
An important special case of the structured matrix stabilization problem is to find, for a given symmetric real matrix, a nearest symmetric positive semidefinite matrix with a prescribed additional structure. A case of particular interest is the following.

\medskip\noindent
{\bf Problem.} {\it Given a symmetric real matrix, find a nearest correlation matrix.}

\medskip\noindent
Here we look for a minimal-norm perturbation of the form $D + \Delta$ with 
$D=I - \text{diag}(A)$ and $\Delta\in \cS_{\rm cor}$ for 
$\cS_{\rm cor}=\{ \Delta \in \R^{n,n}\,:\, \Delta \text{ is symmetric with zero diagonal}\}$
such that all eigenvalues of $A+D+\Delta$ are nonnegative. In the Frobenius-norm setting, we observe $\|D+\Delta\|_F^2 = \|D\|_F^2 + \|\Delta\|_F^2$, and hence $\|\Delta\|_F$ (instead of $\|D+\Delta\|_F$) needs to be minimal. The problem can thus be viewed as the previous matrix stabilization problem for the matrix $-(A+D)$ with structured perturbations $\Delta \in\cS_{\rm cor}$ (unless $A+D$ is already positive semidefinite). A variety of dedicated algorithms, which are not based on this interpretation, have been given in the literature; see
Higham (\cite{Hig02}, \cite{Hig25}) and references therein. On the other hand, the structured matrix stabilization algorithms of Section~\ref{sec:mat-stab} offer an alternative that can be used for the structured problem 
with perturbation matrices $\Delta \in \cS_{\rm cor}$ having a prescribed sparsity pattern.


\subsection*{Nearest defective matrix}
\begcitation{{\it Summary.}\/ Gives a bound for the distance of a matrix having an ill-conditioned eigenvalue problem from a matrix having a multiple eigenvalue which is generally sharper than that which has been published hitherto.}
\signed{(James H.~Wilkinson, \cite{Wil72}, complete abstract)}
\endcitation
\begcitation{For many problems of
numerical analysis, there is a simple relationship between the condition
number of a problem and the shortest distance from that problem to an
ill-posed one: the shortest distance is proportional to the reciprocal of the
condition number (or bounded by the reciprocal of the condition number).
This is true for matrix inversion, computing eigenvalues and eigenvectors,
finding zeros of polynomials, and pole assignment in linear control systems.}
\signed{(James W. Demmel, \cite{Dem87}, from the abstract)}
\endcitation

\pagebreak[3]
\noindent
A matrix is called {\it defective} if 
its Jordan canonical form has a Jordan block of dimension greater than 1. For the corresponding defective eigenvalues (and only for them), the eigenvalue condition number becomes infinite. Here, the reciprocal of the
eigenvalue condition number is the absolute value of the inner product of left and right normalized eigenvectors.

\medskip\noindent
{\bf Problem.} {\it Given a matrix having distinct eigenvalues, find a nearest defective matrix.}

\medskip\noindent
The problem was originally motivated by ill-conditioned eigenvalue problems in numerical linear algebra; see Alam, Bora, Byers \& Overton (\cite{AlaBBO11}) and the Notes of Chapter~\ref{chap:mnp-mix} for references to the history.
It has recently found unexpected practical interest in non-hermitian optics and photonics, where the desired `exceptional points' are nearest defective matrices, possibly with a linear structure and under further constraints; see the references in the Notes of Chapter~\ref{chap:mnp-mix}.

For the general complex case, there exist various algorithms based on a characterization by Alam \& Bora (\cite{AlaB05}), which reduces the problem to finding the smallest $\eps>0$ for which
two components of $\Lambda_\eps(A)$ coalesce in a point; see again the Notes of Chapter~\ref{chap:mnp-mix}.

\subsubsection*{Structured problem.} In the unstructured complex case, every matrix with multiple eigenvalues is arbitrarily close to a matrix with a defective eigenvalue. This need not hold true in structured cases, as is obvious when the structure space consists of normal matrices.
We are not aware of any algorithm in the literature for finding nearest real or structured defective matrices. A two-level algorithm, which is suited also for such cases, is given in 
Section~\ref{sec:defective}. The algorithm minimizes the reciprocal of the eigenvalue condition number of perturbed matrices with a fixed perturbation size via a rank-2 matrix differential equation in the inner iteration. It determines the distance to defectivity in the outer iteration  
by requiring that the reciprocal of the eigenvalue condition number become~$0$.

\subsection*{Nearest singular matrix pencil}
\begcitation{This is a first step toward the goal of finding a way to calculate a smallest norm de-regularizing perturbation of a given square matrix pencil.}
\signed{(Ralph Byers, Chunyang He, and Volker Mehrmann, \cite{ByeHM98})}
\endcitation

\noindent
A matrix pencil $\{ A - \mu B \,:\, \mu \in \C \}$ or more briefly $(A,B)$ (with square matrices $A$ and $B$)  is called {\it singular} if $A-\mu B$ is a singular matrix for all  $\mu \in \C$
(and is called a regular matrix pencil otherwise). This notion is fundamental in the theory of linear differential-algebraic equations $B\dot x(t) = Ax(t)$, which arise in descriptor systems and as linearizations of nonlinear differential-algebraic equations $F(x,\dot x)=0$ near a stationary point. The matrix pencil $(A, B)$ is singular if and only if no initial value $x(0)$ exists such that the corresponding initial value problem has a unique solution. 

\medskip\noindent
{\bf Problem.} \ {\it Given a regular matrix pencil, find a nearest singular matrix pencil.}

\medskip\noindent
Byers, He \& Mehrmann (\cite{ByeHM98}) derived characterizations of minimal perturbations and upper and lower bounds for the distance.
Algorithms for this problem appeared later; see the Notes of Chapter~\ref{chap:mnp-mix}. In Section~\ref{sec:matrix-pencils}, we extend the algorithm for computing a nearest (structured) singular matrix, as given in Section~\ref{sec:sing-S}, to compute the nearest singular matrix pencil to $(A,B)$ under (possibly structured) perturbations of $A$ when $B$ is a fixed singular $n\times n$ matrix. The algorithm extends further to the case where $B$ can also be perturbed. It still uses only matrix--vector products and vector inner products, but now $(n+1)$ times more of them in every iteration. Such an extra factor does not occur for the simpler problem of finding the nearest matrix pencil for which the two matrices have a common null-vector.

There are many more interesting nearness problems for matrix pencils and, more generally, matrix-valued polynomials and even more general nonlinear matrix-valued functions of a complex parameter, not least those related to robust stability of differential-algebraic equations and delay differential equations; see Section~\ref{sec:delay}.

\subsection*{Matrix nearness problems in robust control}
\begcitation{Problems of uncertainty and robustness, which had been forgotten for some time in `modern control'\!, gradually re-emerged and came to the foreground of control theory.}
\signed{(Diederich~Hinrichsen and Anthony~J.~Pritchard, \cite{HinP05})}
\endcitation
\noindent
Consider a linear time-invariant system with input $u(t)\in \C^p$, output $y(t)\in \C^m$ and state vector $x(t)\in\C^n$,
\begin{align*}
\dot x(t) & =  A x(t) + Bu(t)
\\
y(t) & =  C x(t) + Du(t)
\end{align*}
where all eigenvalues of $A$ have negative real part.
In Chapter~\ref{chap:lti} we discuss basic problems of robust control of such systems. We rephrase them as eigenvalue optimization problems and matrix (or operator) nearness problems. The problems considered from this perspective include the following:

\medskip\noindent
{\bf Problems.} {\it
\begin{itemize}
    \item  Compute the $\mathcal{H}_\infty$-norm of the system, which is the operator norm of the input-to-output map $u\in L^2(0,\infty;\C^p) \mapsto y\in L^2(0,\infty;\C^m)$.
    \item  Find an $\mathcal{H}_\infty$-nearest uncontrollable system to a given controllable system under perturbations of the matrix $B$.
    \item  Find an $\mathcal{H}_\infty$-nearest passive system to a given non-passive system 
   under perturbations of the matrix $C$.
    \item  Find a non-contractive system by a minimal structured perturbation to the state matrix $A$ of a given contractive system.
\end{itemize}
}
The distances from systems with undesired properties are important robustness measures of a given control system, whereas algorithms for finding a nearby system with prescribed desired properties are important design tools; see the Notes at the end of Chapter~\ref{chap:lti}. In that chapter, the algorithmic approach to eigenvalue optimization via low-rank matrix differential equations and the two-level approach to matrix nearness problems of previous chapters is extended to problems of robust control and is shown to yield a versatile approach to novel algorithms in that field.

\subsection*{Graph problems as matrix nearness problems}
\begcitation{There is no question
that eigenvalues play a central role in our fundamental understanding of graphs.}
\signed{(Fan Chung, \cite{Chu97})}
\endcitation
\noindent
Eigenvalues or eigenvectors of matrices such as the graph Laplacian or the weighted adjacency matrix characterize basic properties of graphs such as connectivity and centrality. In Chapter~\ref{chap:graphs}, we study problems of clustering and ranking nodes of undirected weighted graphs, which can be viewed as matrix nearness problems that refer to eigenvalues and eigenvectors of sparsity-structured symmetric real matrices under non-negativity constraints for the weights. We consider in detail the following two exemplary problems, which will be reformulated as matrix nearness problems:

\medskip\noindent
{\bf Problem.} \ {\it Given a connected weighted undirected graph, find a minimum cut subject to constraints such as must-link, cannot-link, and cardinality constraints.}

\medskip\noindent
{\bf Problem.} \ {\it Given a connected weighted undirected graph where node $i$ is ranked higher 
than node~$j$, find a minimal perturbation of the weights such that the ranking is reversed.}

\medskip\noindent
For both problems we propose algorithms that do eigenvalue optimization without eigenvalues, but instead use Rayleigh quotients. This requires only  vector inner products and matrix--vector products with matrices having the sparsity pattern of the adjacency matrix of the graph.

\subsection*{Outline}
The matrix nearness problems presented above (and some more) will be discussed in detail in the following chapters. For their numerical solution, a common framework throughout this book is provided by a two-level approach 
\bcl that combines eigenvalue optimization and root-finding in a nested or alternating way:  an eigenvalue optimization problem for a fixed perturbation size is to be solved via a gradient-based low-rank matrix differential equation, and the root-finder adjusts the perturbation size by solving a nonlinear scalar equation, typically by a Newton--bisection method. In the next two chapters we concentrate on the basic eigenvalue optimization problem, and in Chapter~\ref{chap:two-level} on its combination with scalar root-finding. \ecl
Chapter~\ref{chap:struc} extends the programme to structured matrix nearness problems. Chapter~\ref{chap:mnp-mix} discusses diverse matrix nearness problems that go beyond those of previous chapters. Chapters~\ref{chap:lti} and~\ref{chap:graphs} are about matrix nearness problems in the areas of robust control and graph theory, respectively. The appendix, Chapter~\ref{chap:appendix}, presents known basic results on derivatives of eigenvalues and eigenvectors, which are often used in this book.

\newcommand{\lreps}{\mu_\varepsilon}
\newcommand{\pp}{\alpha}
\newcommand{\pq}{\beta}
\newcommand{\cI}{\mathcal{I}}
\newcommand{\cJ}{\mathcal{J}}
\chapter{A basic eigenvalue optimization problem}\label{chap:proto}

\bcl We describe an algorithmic approach to optimization problems where a function is minimized (or maximized) over the $\eps$-pseudospectrum of a given matrix $A$, that is, over all eigenvalues of 
all complex matrices $A+\Delta$ with $\|\Delta\|_F\le\eps$. \ecl This leads to non-convex, non-smooth optimization problems.
The approach taken here uses constrained gradient flows and the remarkable rank-1 property of the optimizers. It leads to an algorithm that numerically follows rank-1 matrix differential equations to their stationary points.

This chapter is basic in the sense that we illustrate essential ideas and techniques on a particular problem class that will be vastly expanded later in this book. The eigenvalue optimization problem considered in this chapter
arises in computing 
pseudospectra or their extremal points. In particular, as will be discussed in Chapter~\ref{chap:pseudo}, right-most points determine the pseudospectral abscissa, and points of largest modulus yield the pseudospectral radius. These eigenvalue optimization problems  (and extensions thereof)  will reappear as a principal building block in the two-level approach to a wide variety of matrix nearness problems to be discussed in later chapters, from Chapter~\ref{chap:two-level} onwards.
In Chapter~\ref{chap:struc} and later we will also consider optimizing eigenvalues over real perturbation matrices and more generally over structured perturbations
that are restricted to a given complex- or real-linear subspace of matrices, for example matrices with a given sparsity pattern, or matrices with prescribed range and co-range, or Hamiltonian or Toeplitz or Hankel matrices. 
In all these cases there is a common underlying rank-1 property of optimizers that will be used to advantage in the algorithms.


\section{Problem formulation}
\label{subsec:proto-problem}




\bcl
Let $A\in\C^{n,n}$ be a given square matrix. The objective is to minimize a function of eigenvalues of matrices $A+\Delta$ over all perturbation matrices $\Delta\in\C^{n,n}$ of a prescribed norm~$\eps$.  

\subsubsection*{Function to be minimized.} We consider a continuously differentiable function
\begin{equation}
\label{ass:f}
f: \C^2 \rightarrow \C \quad\text{with}\quad f\left( \lambda, \clambda \right) = f\left( \clambda, \lambda \right) \in \R \quad\text{for all }\,\lambda\in\C,
\end{equation}
where $\clambda$ is the complex conjugate of $\lambda$. 
\ecl
While our theory applies to general functions $f$ with \eqref{ass:f}, in our examples we often consider specific cases where 
$f$ or $-f$ evaluated at $\left( \lambda, \clambda \right)$ equals
\begin{equation}\label{f-ex}
\Re\,\lambda = \frac{\lambda + \clambda}{2} \quad 
\text{or}\quad| \lambda |^2 = \lambda \clambda.
\end{equation}  
We note that $\Im\,\lambda=\tfrac1{2\iu}(\lambda-\clambda)$ does not satisfy \eqref{ass:f}, but this case can be included in the present setting by first rotating $A$ to $-\iu A$ and then considering the real part.

\bcl
\subsubsection*{Target eigenvalue.}
A {\it target eigenvalue} of a matrix $A+\Delta$ (with respect to the objective function $f$) is an eigenvalue $\lambda$ for which the function value $f(\lambda,\clambda)$ is minimal among all eigenvalues of $A+\Delta$. We then write $\lambda(A+\Delta)$ for a target eigenvalue
of $A+\Delta$. 

For $f$ or $-f$ as in \eqref{f-ex}, a target eigenvalue is an eigenvalue of minimal or maximal real part, or an eigenvalue of minimal or maximal modulus.

A target eigenvalue may not depend continuously on the matrix $\Delta$ when several eigenvalues are simultaneously extremal, but it depends continuously on $\Delta$ if the target eigenvalue is unique. Moreover, the function $\lambda(A+\cdot)$ is arbitrarily differentiable at $\Delta$ if additionally the target eigenvalue $\lambda(A+\Delta)$ is an algebraically simple eigenvalue; see Theorem~\ref{chap:appendix}.\ref{thm:eigderiv} in the Appendix.
\ecl


\bcl
\subsubsection*{Eigenvalue optimization problem.}
\index{eigenvalue optimization!basic problem} 
\index{eigenvalue optimization!complex}
The problem is to find, for a given $\varepsilon>0$, 
\begin{equation} \label{eq:optimiz0-all}
\arg\min\limits_{\lambda\in\Lambda(A+\Delta),\: \Delta \in \C^{n,n} \text{ with } \| \Delta \|_F \le \eps} f \left( \lambda, \clambda   \right), 
\end{equation}
where $\Lambda(A+\Delta)$ is the set of eigenvalues (spectrum) of $A+\Delta$  and
$\| \Delta \|_F$ is the Frobenius norm\index{Frobenius norm} of the matrix $\Delta$. 
This amounts to minimizing $f$ over the $\eps$-pseudo\-spectrum of $A$,
\index{pseudospectrum!complex}
$$
\Lambda_\eps(A)=\{ \lambda \in \C \,:\, \lambda \in \Lambda(A+\Delta) \text{ for some } \Delta \in \C^{n,n} \text{ with } \| \Delta \|_F \le \eps \}.
$$
In many situations of interest, including \eqref{f-ex}, the minimum is attained with equality $\|\Delta\|_F=\eps$, as 
will be discussed in Section~\ref{subsec:ineq}. A sufficient condition for equality is that $\partial f /\partial \lambda$ has no zero in $\Lambda_\eps(A)$. 

We will therefore mostly study the eigenvalue optimization problem with the equality constraint $\|\Delta\|_F=\eps$ in place of the inequality constraint $\|\Delta\|_F \le \eps$ in \eqref{eq:optimiz0-all}.
By definition of a target eigenvalue $\lambda(A+\Delta)$, this 
is equivalent to finding
\begin{equation} \label{eq:optimiz0}
\arg\min\limits_{\Delta \in \C^{n,n},\, \| \Delta \|_F = \eps} f \left( \lambda\left( A + \Delta \right), \clambda \left( A + \Delta \right)  \right).
\end{equation}
\ecl
It is convenient to write
\[
\Delta = \eps E \quad \mbox{ with} \ \| E \|_F = 1
\]
and
\begin{equation} \label{Feps}
\bcl\F_\eps(E)\ecl = f \left( \lambda\left( A + \eps E \right), \clambda\left( A + \eps E \right)  \right)
\end{equation}
\ecl
so that Problem~\eqref{eq:optimiz0} is equivalent to finding
\begin{equation} \label{eq:optimiz}
\arg\min\limits_{E \in \C^{n,n}, \, \| E \|_F = 1} \F_\eps(E).
\end{equation}
Problem \eqref{eq:optimiz0} or \eqref{eq:optimiz} is a nonconvex, nonsmooth optimization problem.
The problem with the inequality constraint $\| E \|_F \le 1$ will equally be considered.
The $\arg\max$ case  is treated in the same way, replacing $f$ by $-f$.

\medskip
There are obvious generalizations with interesting applications, which we will encounter in later chapters \bcl and approach with techniques that extend those of the present chapter.\ecl
\begin{itemize}
\item The perturbation size $\eps$ might vary such that the minimum of $\F_\eps$ takes a prescribed value. 
\bcl (This is the situation encountered with matrix nearness problems in Chapter~\ref{chap:two-level} and later chapters.)\ecl
\item The perturbation matrix $\Delta$ might be restricted to be real or to have a prescribed sparsity pattern or to belong to another complex-linear or real-linear subspace of $\C^{n,n}$. 
\bcl (This is the situation encountered with structured matrix nearness problems in Chapter~\ref{chap:struc} and later chapters.)
\ecl
\item
The objective function $f$
 might depend on several or all eigenvalues of $A+\Delta$ instead of only a single target eigenvalue.
\bcl (This occurs with matrix stabilization in Chapter~\ref{chap:mnp-mix}.)\ecl
\item 
\bcl The minimization in \eqref{eq:optimiz0-all} might be done only over a subset of eigenvalues of $A+\Delta$.
(This occurs with Hamiltonian matrix nearness problems in Chapter~\ref{chap:mnp-mix}.)\ecl
\item The objective function might depend on eigenvectors of $A+\Delta$.
\bcl (This occurs with the problem of finding the nearest defective matrix in Chapter~\ref{chap:mnp-mix} and in the ranking problem of Chapter~\ref{chap:graphs}.)
\item There might be constraints on the perturbed matrices $A+\Delta$ such as nonnegativity of the entries. (This occurs in matrix nearness problems related to graphs in Chapter~\ref{chap:graphs}.)
\ecl

\end{itemize}
\bcl However, in this chapter we will only consider the function $f$ as in \eqref{ass:f} and its minimization as in
\eqref{eq:optimiz0} with a given fixed perturbation size $\eps$.
\ecl

\section{Gradient flow}
\label{sec:proto-complex}

\subsection{Free gradient}
We begin with some notations and normalizations.
Let $x$ and $y$ be  left and right eigenvectors, respectively, associated with a simple eigenvalue $\lambda$ of a matrix $M$: $\ x,y\in\C^n\setminus\{0\}$ with $x^*M=\lambda x^*$ and $My=\lambda y$, where $x^*={\overline x}^\top$. Unless specified differently, we assume that the eigenvectors are normalized such that
\begin{equation} \label{eq:scaling}
\| \y \| = \| \x \| = 1 \quad\text{ and } \quad \y^*\x\, \text{ is real and positive.}
\end{equation}
The norm $\|\cdot\|$ is chosen as the Euclidean norm.
Any pair of left and right eigenvectors $x$ and $y$ can be scaled in this way.
\index{eigenvectors, left and right!scaling}

We denote by 
\begin{equation*}
\langle X,Y \rangle =\sum_{i,j} \conj{x}_{ij}y_{ij} = {\rm trace}(X^* Y)
\end{equation*} 
the inner product in $\C^{n,n}$ that induces the Frobenius norm $\| X \|_F = \langle X,X \rangle^{1/2}$.

The following lemma will allow us to compute the steepest descent direction of the functional $\F_\eps$.

\begin{lemma}[Free gradient] \label{lem:gradient} 
Let $E(t)\in \C^{n,n}$, for real $t$ near $t_0$, be a continuously differentiable path of matrices, with the derivative denoted by $\dot E(t)$.
Assume that $\lambda(t)$ is a simple eigenvalue of  $A+\eps E(t)$ depending continuously on $t$, with associated left and right eigenvectors
$\y(t)$ and $\x(t)$ satisfying \eqref{eq:scaling}, and let the eigenvalue condition number be 
\begin{equation*}
\kappa(t) = \frac1{\y(t)^* \x(t)} > 0.
\end{equation*}
Then, $\F_\eps(E(t))=f \bigl( \lambda(t) , \overline{\lambda(t)}  \bigr)$ 
is continuously differentiable w.r.t. $t$ and we have
\begin{equation} \label{eq:deriv}
\frac1{ \eps \kappa(t) } \,\frac{d}{dt} \F_\eps(E(t)) = \Re \,\bigl\langle  G_\eps(E(t)),  \dot E(t) \bigr\rangle,
\end{equation}
where the (rescaled) gradient of $\F_\eps$ is the rank-1 matrix
\begin{equation} \label{eq:freegrad}
G_\eps(E) = 2 f_{\clambda} \, \y \x^* \in \C^{n,n}
\end{equation}
with\footnote{\bcl We write ${\partial f}/{\partial \clambda}$ or $f_{\clambda}$ for the partial derivative of $f$ with respect to the second argument.} 
$\displaystyle f_{\clambda} = \frac{\partial f}{\partial \clambda}(\lambda, \clambda)$ for the eigenvalue $\lambda=\lambda(A+\eps E)$ and the corresponding left and right eigenvectors $\y$ and $\x$ normalized by \eqref{eq:scaling}.
\end{lemma}
\index{gradient!free}
\begin{proof}
We first observe that \eqref{ass:f} implies $f_{\clambda} = \conj{f_{\lambda}} =\displaystyle\conj{\frac{\partial f}{\partial \lambda}(\lambda, \clambda)}$.
Using Theorem \ref{chap:appendix}.\ref{thm:eigderiv}, we obtain that $\F_\eps(E(t))$ is differentiable with
\begin{eqnarray} \nonumber
\frac{ d }{dt} \F_\eps \left( E(t) \right) & = & f_\lambda \,\dot \lambda + f_{\clambda} \,\overline{\dot \lambda} \\
& = & \frac{\eps}{\y^*\x}  \left( f_\lambda\, \y^* \dot{E} \x + f_{\clambda} \;\conj{\y^* \dot{E} \x} \right) = 
\frac{ \eps}{\y^*\x} \, 2\,\Re \left( f_\lambda\, \y^* \dot{E} \x \right),
\end{eqnarray}
where we omit indication of the ubiquitous dependence on $t$ on the right-hand side.
Noting that
\[
\Re \bigl( f_\lambda \,\y^*\dot{E} \x  \bigr) = \Re\, \bigl\langle \conj{f_{\lambda}}\, \y \x^* , \dot{E}  \bigr\rangle,
\]
we obtain \eqref{eq:deriv}--\eqref{eq:freegrad}.
\qed \end{proof}

\begin{example} \label{ex:G}
For $f(\lambda,\clambda) =-\tfrac12(\lambda+\clambda)=-\Re\,\lambda$ we have $2 f_{\clambda}=-1$ and hence $G_\eps(E)=-xy^*$,
which is nonzero for all $\lambda$.
For $f(\lambda,\clambda) = -\tfrac12 |\lambda|^2 = -\tfrac12 \lambda \clambda$ 
we have 
$2 f_{\clambda} =  -\lambda$.
In this case we obtain
$
G_\eps(E) =  - \lambda \,\y \x^*,
$
which is nonzero whenever $\lambda\ne 0$. 
\end{example}

\subsection{Projected gradient}
\index{gradient!projected}

To satisfy the constraint ${\| E(t) \|_F^2 =1}$, we must have 
\begin{equation} \label{eq:normconstr}
 0 = \frac12\,\frac d{dt}\| E(t) \|_F^2= \Re\, \langle E(t), \dot E(t) \rangle.
\end{equation}
In view of Lemma~\ref{lem:gradient} we are thus led to the following constrained optimization problem for the admissible direction of steepest descent.

\begin{lemma}[Direction of steepest admissible descent]
\label{lem:opt} 
Let $E,G\in\C^{n, n}$ with ${\|E\|_F=1}$. A solution of the optimization problem
\begin{eqnarray}
Z_\star  & = & \arg\min_{\|Z\|_F=1,\,\text{\rm Re}\,\langle E, Z \rangle=0} \ \Re\,\langle  G,  Z \rangle,
\label{eq:opt}
\end{eqnarray}
is given by 
\begin{align}
\mu Z_\star & =  -G + \Re\hspace{1pt}\langle G, E \rangle\, E,
\label{eq:Eopt}
\end{align}
where $\mu$ is the Frobenius norm of the matrix on the right-hand side. 
The solution is unique if $G$ is not a multiple of $E$.
\end{lemma}

\begin{proof}
The result follows on noting that the real part of the complex inner product on $\C^{n,n}$ is a real inner product on $\R^{2n,2n}$, and the real inner product with a given vector (which here is a matrix) is minimized over a subspace by orthogonally projecting the vector onto that subspace. The expression in (\ref{eq:Eopt}) is the orthogonal projection of $-G$ onto the orthogonal complement of the span of $E$, which is the tangent space at $E$ of the manifold of matrices of unit Frobenius norm, i.e. the space of admissible directions.
\qed \end{proof}
\subsection{Norm-constrained gradient flow} \index{gradient flow}\index{gradient flow!norm-constrained}
Lemmas~\ref{lem:gradient} and~\ref{lem:opt} show that the admissible direction of steepest descent of the functional $\F_\eps$ at a matrix $E$ of unit Frobenius norm is given by the positive multiples of the matrix $-G_\eps(E)+ \Re\,\langle G_\eps(E), E \rangle E$. 
This leads us to consider the (rescaled) {\it gradient flow on the manifold of $n\times n$ complex matrices of unit Frobenius norm}:
\begin{equation}\label{ode-E}
\dot E = -G_\eps(E)+ \Re\,\langle G_\eps(E), E \rangle E,
\end{equation}
where we omitted the ubiquitous argument $t$. 

By construction of this ordinary differential equation, we have $\Re \langle E, \dot E \rangle =0$ along its solutions, and so the Frobenius norm $1$ is conserved.
Since we follow the admissible direction of steepest descent of the functional $\F_\eps$ along solutions $E(t)$ of this differential equation, we obtain the following monotonicity property.

\begin{theorem}[Monotonicity] \label{thm:monotone}
Assume that $\lambda(t)$ is a simple eigenvalue of  $A+\eps E(t)$ depending continuously on $t$. 
Let $E(t)$ of unit Frobenius norm satisfy the differential equation {\rm (\ref{ode-E})}.
Then,
\begin{equation}
\frac{d}{dt} \F_\eps (E(t))  \le  0.
\label{eq:pos}
\end{equation}
\end{theorem}

\begin{proof}
Although the result follows directly from Lemmas~\ref{lem:gradient} and \ref{lem:opt}, we compute the derivative explicitly. We write $G=G_\eps(E)$ for short and take the inner product of \eqref{ode-E} with $\dot E$. Since
$\Re\langle E, \dot E \rangle = \tfrac12 d/dt \| E \|_F^2 = 0$, we find
$$
\| \dot E \|_F^2 = - \Re \langle G - \Re\langle G,E \rangle E, \dot E\rangle =
- \Re \langle G, \dot E\rangle
$$
and hence Lemma~\ref{lem:gradient} and \eqref{ode-E} yield
\begin{equation}\label{c-s}
\frac1{\eps\kappa} \,\frac{d}{dt} \F_\eps(E(t)) = \Re \langle G, \dot E\rangle = - \| \dot E \|_F^2 =
- \| G - \Re\,\langle G, E \rangle E \|_F^2 \le 0,
\end{equation}
which gives the precise rate of decay of $\F_\eps$ along a trajectory $E(t)$ of \eqref{ode-E}.
%
\qed \end{proof}
%
%
The stationary points of the differential equation  \eqref{ode-E} are characterized as follows.

\index{stationary point}
\begin{theorem}[Stationary points] \label{thm:stat} 
Let $E_\star\in\C^{n,n}$ with $\| E_\star\|_F=1$ be such that 
\begin{itemize}
\item[(i)] \  The target eigenvalue $\lambda(A+\eps E)$ is simple at $E=E_\star$ and depends continuously on $E$ in a neighborhood of $E_\star$.
\item[(ii)] \ \ \ The gradient $G_\eps(E_\star)$ is nonzero, 
i.e., $f_{\clambda}(\lambda_\star,\clambda_\star)\ne 0\,$
for $\lambda_\star\!=\lambda(A+\eps E_\star)$.
\end{itemize}
Let $E(t)\in \C^{n,n}$ be the solution of \eqref{ode-E} passing through $E_\star$. 
Then the following are equivalent: 
\begin{itemize}
\item[1.] $\displaystyle\frac{ d }{dt} \F_\eps \left( E(t) \right)  = 0$. \\
\item[2.] $\dot E = 0$. \\[-3mm]
\item[3.] $E_\star$ is a real multiple of $G_\eps(E_\star)$.
\end{itemize}
\end{theorem}

\begin{proof} Clearly, 3.~implies 2., which implies 1. Finally, \eqref{c-s} shows that 1. implies 3.
\qed \end{proof}

\begin{remark}[Degeneracies] In degenerate situations where $G_\eps(E_\star)= 0$, we cannot conclude \bng that 2. implies 3., \eng i.e., that the stationary point is a multiple of  $G_\eps(E_\star)$. For the case
$f(\lambda,\clambda) = -\Re\,\lambda$ we have seen in Example~\ref{ex:G} that $G_\eps(E)=-xy^*\ne 0$, where $x,y$ are normalized eigenvectors to the target eigenvalue $\lambda(A+\eps E)$.
For $f(\lambda,\clambda)=-\frac12 |\lambda|^2$ we have
$G_\eps(E)=-\lambda xy^*\ne 0$ for $\lambda\ne 0$. For other functions $f$ we might encounter $G_\eps(E_\star)= 0$, but such a degeneracy can be regarded as an exceptional situation, which will not be considered further.
\end{remark}

\begin{remark}[Stationary points and optimizers] \label{rem:stat-min} \index{optimizer!vs. stationary point}
Every global minimum is a local minimum, and every local minimum is a stationary point. The converse is clearly not true. Stationary points of the gradient system, that are not a local minimum, are unstable. 
It can thus be expected that generically a  trajectory will end up in a local minimum. Running several different trajectories reduces the risk of being caught in a local minimum instead of a global minimum. 
\end{remark}

\begin{remark} [Multiple and discontinuous eigenvalues] \label{rem:mult-eig}
We mention some situations where the assumption of a smoothly evolving simple eigenvalue is violated. 

--- Along a trajectory $E(t)$, the target eigenvalue  $\lambda(t)=\lambda(A+\eps E(t))$ may become discontinuous. For example, in the case of the eigenvalue of largest real part, a different branch of eigenvalues may come to have the largest real part. 
\bng
In such a case of discontinuity, the computed target eigenvalue has a jump and the differential equation is further solved following the new branch, still with descent of the function $f$ until finally a stationary point is approximately reached. 
\eng

--- A multiple eigenvalue $\lambda(t)$ may occur at some finite $t$ because of a coalescence of eigenvalues.
Even if some continuous trajectory runs into a coalescence \bcl at an isolated finite value of $t$, 
after discretization of the differential equation it appears unlikely that the discrete trajectory will run into a coalescence (as opposed to an almost-coalescence)
on a finite time interval. \ecl

--- A multiple eigenvalue may appear in a stationary point, in the limit $t\to\infty$. The computation will stop before, and items 1.-3. in Theorem~\ref{thm:stat} will then be satisfied approximately, in view of~\eqref{c-s}.

Although the cases above do not affect the time-stepping of the gradient system, close-to-multiple eigenvalues do impair the accuracy of the computed left and right eigenvectors that appear in the gradient.

\bcl
The situation can change dramatically when the minimization is not done over {\it all} complex matrices of a given norm (as is done here) but only over matrices of a given norm in a prescribed subspace (as will be considered in Chapter~\ref{chap:struc}). Then, there are cases where multiple eigenvalues become generic; see Burke, Lewis \& Overton (\cite{BuLeOv00}) for a striking example with Jordan blocks of arbitrary dimension.
\ecl
\end{remark}

\bcl
\subsection{Inequality constraint vs. equality constraint} \ecl
\label{subsec:ineq}
\index{inequality constraint}
When we have the inequality constraint $\|\Delta\|_F \le \eps$ in \eqref{eq:optimiz0} or equivalently $\| E \|_F \le 1$ in \eqref{eq:optimiz}, the situation changes only slightly. If $\|E\|_F < 1$, every direction is admissible, and the direction of steepest descent is given by the negative gradient $-G_\eps(E)$. So we choose the free gradient flow
\begin{equation}\label{ode-E-free}
\dot E = -G_\eps(E) \qquad\text{as long as }\ \| E(t) \|_F < 1.
\end{equation}
When  $\|E(t)\|_F=1$, then there are two possible cases. If 
$\Re\,\langle G_\eps(E), E \rangle \ge 0$, then the solution of \eqref{ode-E-free} has (omitting the argument $t$)
$$
\frac d{dt} \|E(t)\|_F^2 = 2 \,\Re\,\langle \dot E, E \rangle= -2\, \Re\,\langle G_\eps(E), E \rangle \le 0,
$$
and hence the solution of \eqref{ode-E-free} remains of Frobenius norm at most 1. 

Otherwise, if $\Re\,\langle G_\eps(E), E \rangle < 0$, 
the admissible direction of steepest descent is given by the right-hand side of \eqref{ode-E}, i.e. $-G_\eps(E) + \Re\,\langle G_\eps(E), E \rangle E$,
and so we choose that differential equation to evolve $E$. The situation can be summarized as taking, if $\|E(t)\|_F=1$,
\begin{equation} \label{ode-E-mu}
\dot E = -G_\eps(E) + \mu E \quad\text{ with }\ \mu=\min\bigl(0,\Re\,\langle G_\eps(E), E \rangle)\bigr).
\end{equation}
Along solutions of \eqref{ode-E-mu}, the functional $\F_\eps$ decays monotonically, and stationary points of \eqref{ode-E-mu} with
$G_\eps(E)\ne 0$ are characterized, by the same argument as in Theorem~\ref{thm:stat}, as
\begin{equation}\label{stat-neg}
\mbox{$E$ is a {\it negative} real multiple of $G_\eps(E)$.}
\end{equation}
\bcl
If  $G_\eps(E_\star)\ne 0$ at an optimizer $E_\star$ (as in Example~\ref{ex:G}), it can thus be concluded that the optimizer of the problem with inequality constraints is a stationary point of the gradient flow \eqref{ode-E} for the problem with equality constraints. We note that the condition ${G_\eps(E_\star)\ne 0}$ translates into the condition $f_{\clambda}(\lambda_\star,\clambda_\star)\ne 0$ for $\lambda_\star=\lambda(A+\eps E_\star)$.

An important result is the following minimum principle, which allows us to restrict our attention to the equality constraint. Here we refer to the subset of the $\eps$-pseudospectrum $\Lambda_\eps(A)$ that collects all target eigenvalues $\lambda(A+\eps E)$ for ${\|E\|_F \le 1}$:
\begin{equation}\label{Lambda-odot}
    \Lambda_\eps^\odot(A) = \{ \lambda(A+\eps E) \,:\, \|E\|_F \le 1 \} \subset \Lambda_\eps(A).
\end{equation}

\begin{theorem}[Equivalence of equality and inequality constraints]
\label{thm:eq-ineq}
 If 
\begin{equation}
\label{df-nonzero} 
 \frac{\partial f}{\partial\clambda}(\lambda,\clambda)\ne 0 \quad\text{ for all } \ \lambda \in \Lambda_\eps^\odot(A),
\end{equation}
then the eigenvalue optimization problems with equality constraint and with inequality constraint have the same minima. In other words, if $E_\star$ minimizes $\F_\eps$ under the inequality constraint $\| E_\star \|_F \le 1$, then $\| E_\star \|_F=1$.
\end{theorem} 

\begin{proof} If $\| E_\star \|_F$ were strictly smaller than 1, then $E_\star$ would be a stationary point of the free gradient flow, i.e. $G_\eps(E_\star) = 0$, which implies $f_{\clambda}(\lambda_\star,\clambda_\star)= 0$ for $\lambda_\star=\lambda(A+\eps E_\star)\in \Lambda_\eps^\odot(A)$ in  contradiction to assumption \eqref{df-nonzero}.
\qed
\end{proof}

We note that condition \eqref{df-nonzero} is satisfied when $f$ is chosen as $\pm\Re \lambda$ and $-|\lambda|^2$,
and for $+|\lambda|^2$ if $0\notin \Lambda_\eps(A)$.

\ecl 

\section{Rank-1 constrained gradient flow}
\label{sec:rank-1-basic}

\subsection{Rank-1 property of optimizers} \index{optimizer!rank-1 property}
We call an optimizer $E_\star$ of \eqref{eq:optimiz} {\it non-degenerate} if conditions (i) and (ii) of Theorem~\ref{thm:stat} are satisfied.
Since optimizers are necessarily stationary points of the norm-constrained gradient flow~\eqref{ode-E-S},
Theorem~\ref{thm:stat} and Lemma~\ref{lem:gradient} immediately imply the following remarkable property.

\begin{corollary}[Rank of optimizers] \label{cor:rank-1}
If $E_\star$ is a non-degenerate optimizer of the eigenvalue optimization problem \eqref{eq:optimiz}, then $E_\star$ is of rank $1$.
\end{corollary}

Let us summarize how this rank-1 property came about: An optimizer is a stationary point of the norm-constrained gradient flow \eqref{ode-E}. In the non-degenerate case, this implies that the optimizer $E$ is a real multiple of the free gradient $G_\eps(E)$, which is of rank 1 as a consequence of the derivative formula for simple eigenvalues.

This corollary motivates us to search for a differential equation on the manifold of rank-$1$ matrices of norm $1$ with the property that the functional $\F_\eps$ decreases along its solutions and has the same stationary points as the differential equation \eqref{ode-E}. Working with rank-1 matrices $E=uv^*$ given by two vectors $u,v\in \C^n$ instead of general complex $n\times n$ matrices is computationally favourable, especially for high dimensions $n$, for two independent reasons:
\begin{itemize}
\item[(i)]  \quad Storage and computations are substantially reduced when the two $n$-vectors $u,v$  are used instead of the full $n\times n$ matrix $E$.
\item[(ii)] \quad The computation of the target eigenvalue $\lambda(t)$ of $A+\eps E(t)$ using inverse iteration is largely simplified thanks to the Sherman--Morrison formula \index{Sherman--Morrison formula}\index{rank-1 matrix!computational advantages}
$$
(A+\eps uv^* - \mu I)^{-1} = (A-\mu I)^{-1} - \frac{(A-\mu I)^{-1}\eps uv^*(A-\mu I)^{-1}}{1+v^*(A-\mu I)^{-1}\eps u}.
$$
\bng In a numerical integration, in order to approximate the target eigenvalue $\lambda(t_n)$, we can shift the matrix by the previously computed target eigenvalue $\mu=\lambda_{n-1} \approx \lambda(t_{n-1})$. 
\bcltwo
The only drawback is that in this way we follow in time a branch of target eigenvalues (as the rightmost), but would not detect jumps to another branch.
\ecltwo
\eng
\end{itemize}
Moreover, after transforming the given matrix $A\in \C^{n,n}$ to Hessenberg form via a unitary similarity transformation with $O(n^3)$ operations,
 linear systems with the shifted matrix $A-\mu I$ for varying shifts $\mu$ can be solved with $O(n^2)$ operations each.  

For sparse matrices $A$, Krylov subspace methods for the perturbed matrix $A+\eps E$ are advantageous when $E$ is of rank 1, since matrix-vector products with $E=uv^*$ only \bcl require computing a vector inner product, which takes \ecl $\bigo(n)$ operations.
\medskip

Since the matrix is complex, you can write "triangular form" here, and observe that Hessenberg form applies in the real case where only real similarity transformations are permitted

\subsection{Rank-1 matrices and their tangent matrices}\index{rank-1 matrix}
\label{subsec:rank-1}
We denote by $\cR_1=\cR_1(\C^{n,n})$ the manifold of complex rank-1 matrices of dimension $n\times n$ and write $E\in\cR_1$ in a non-unique way as
$$
E=\sigma uv^*,
$$
where $\sigma\in \C\setminus\{ 0 \}$ and $u,v\in \C^n$ have unit norm. The tangent space $T_E\cR_1$ at $E\in\cR_1$ consists of the derivatives of paths in $\cR_1$ passing through $E$.
Tangent matrices $\dot E\in T_E\cR_1$ are then of the form 
\begin{equation}\label{E-dot-1}
\dot E = \dot\sigma uv^* + \sigma \dot u v^* + \sigma u \dot v^*,
\end{equation}
where
$\dot\sigma\in\C$ is arbitrary and $\dot u,\dot v\in \C^n$ are such that $\Re(u^*\dot u)=0$ and $\Re(v^*\dot v)=0$ (because of the norm constraint on $u$ and $v$). They are
uniquely determined by $\dot E$ and $\sigma,u,v$ if we impose the orthogonality conditions
$u^*\dot u=0, \  v^*\dot v=0$. Multiplying $\dot E$ with $u^*$ from the left and with $v$ from the right, we then obtain 
\begin{equation}\label{sigma-u-v-Edot}
\dot\sigma = u^* \dot E v, \quad\sigma \dot u = \dot Ev - \dot \sigma u, \quad
\sigma {\dot v}^* = u^* \dot E - \dot\sigma v^*.
\end{equation} 
Extending this construction, we arrive at a useful explicit formula for the projection onto the tangent space that is orthogonal with respect to the Frobenius inner product $\langle\cdot,\cdot\rangle$.

\begin{lemma}[Rank-1 tangent space projection]
\label{lem:P-formula-1}
The orthogonal projection from $\C^{n, n}$ onto the tangent space $T_E\cR_1$ at $E=\sigma uv^* \in\cR_1$
is given by
\begin{equation}\label{P-formula-1}
P_E(Z) = Z - (I-uu^*) Z (I-vv^*)
\quad\text{ for $Z\in\C^{n, n}$}.
\end{equation}
\end{lemma}
\index{rank-1 matrix!tangent space projection}

\begin{proof}  Let $P_E(Z)$ be defined by \eqref{P-formula-1}. To prove that $P_E(Z)\in T_E\cR_1$, we show that $P_E(Z)$ can be written in the form \eqref{E-dot-1}.
Let $\dot\sigma,\dot u,\dot v$ be defined as in \eqref{sigma-u-v-Edot}, but now with $\dot E\in T_E\cR_1$ replaced by arbitrary $Z\in\C^{n, n}$, i.e.,
\begin{equation}\label{sigma-u-v-Z}
\dot\sigma = u^* Z v, \quad \sigma \dot u = Zv - \dot \sigma u, \quad
\sigma {\dot v}^* = u^* Z - \dot\sigma v^*.
\end{equation} 
We obtain the corresponding matrix in the tangent space $T_E\cR_1$ (see \eqref{E-dot-1}) as
\begin{align*}
&\dot\sigma uv^* + \sigma \dot u v^* + \sigma u \dot v^*
\\
&= \dot\sigma uv^* + (Zv-\dot\sigma u)v^* + u(u^*Z - \dot\sigma v^*)
\\
&= Zvv^* - uu^*Zvv^* + uu^*Z = P_E(Z).
\end{align*}
This shows that
\begin{equation}\label{P-formula-dots}
P_E(Z) =\dot\sigma uv^* + \sigma \dot u v^* + \sigma u \dot v^* \in T_E\cR_1.
\end{equation}
Furthermore,
$$
\langle P_E(Z),\dot E \rangle =   \langle Z, \dot E \rangle \qquad\text{for all }\ \dot E \in T_E\cR_1,
$$
because $\langle (I-uu^*) Z (I-vv^*),\dot E \rangle = \langle Z, (I-uu^*) \dot E (I-vv^*) \rangle = 0$ by \eqref{E-dot-1}. Hence,
$P_E(Z)$ is the orthogonal projection of $Z$ onto $T_E\cR_1$.
\qed \end{proof}
We note that $P_E(E)=E$ for $E\in\cR_1$, or equivalently, $E\in T_E \cR_1$, which will be an often used property.

\subsection{Rank-1 constrained gradient flow}
\label{subsec:rank1-gradient-flow}
\index{gradient flow!rank-1 constrained}
\index{rank-1 matrix differential equation}
\index{gradient!projected}

In the differential equation (\ref{ode-E}) we project the right-hand side to the tangent space $T_E\cR_1$:
\begin{equation}\label{ode-E-1}
\dot E = -P_E\Bigl( G_\eps(E) - \Re \langle G_\eps(E),E \rangle E \Bigr).
\end{equation}
This yields a differential equation on the rank-1 manifold $\cR_1$. In view of Lemma~\ref{lem:gradient}, it is
the (rescaled) gradient flow of the functional $\F_\eps$ constrained to the manifold~$\cR_1$.

Assume now that for some $t$, the Frobenius norm of $E=E(t)$ is 1. Since $P_E(E)=E$, 
we have with $Z=-G_\eps(E) + \Re \langle G_\eps(E),E \rangle E$ that
$$
 \Re\, \langle E, \dot E \rangle = \Re\, \langle E, P_E(Z) \rangle = \Re\, \langle P_E(E),Z \rangle =
\Re\, \langle E,Z \rangle = 0.
$$
Hence, solutions $E(t)$ of  \eqref{ode-E-1} stay of Frobenius norm 1 for all $t$.

The proof of Lemma~\ref{lem:P-formula-1} also provides  the following differential equations for the factors of $E(t)=\sigma(t)u(t)v(t)^*$, 
which can be discretized by standard numerical integrators.

\begin{lemma}[Differential equations for the three factors]
\label{lem:suv-1}
For $E=\sigma uv^* \in \cR_1$ with nonzero $\sigma\in \C$ and with
$u\in\C^{n}$ and $v\in\C^{n}$ of unit norm,
the equation $\dot E=P_E(Z)$ is equivalent to
$
\dot E = \dot\sigma uv^* + \sigma \dot u v^* + \sigma u \dot v^* ,
$
where
\begin{eqnarray}
\dot \sigma &=& u^*  Z v
\nonumber\\
\dot u &=& (I-uu^*) Z v \sigma^{-1}
\label{odes-C-1} \\
\dot v &=& (I-vv^*) Z^*  u {\overline\sigma}^{-1}.
\nonumber
\end{eqnarray}
\end{lemma}

\begin{proof} The result follows immediately from \eqref{sigma-u-v-Z} and \eqref{P-formula-dots}.
\qed
\end{proof}

Since we are only interested in solutions of Frobenius norm 1 of \eqref{ode-E-1}, we can simplify the representation of $E$ to
$E=uv^*$ with $u$ and $v$ of unit norm (without the extra factor $\sigma$ of unit modulus).

\begin{lemma} [Differential equations for the two vectors]
\label{lem:uv-1}
For an initial value $E(0)=u(0)v(0)^*$ with $u(0)$ and $v(0)$ of unit norm, the solution of \eqref{ode-E-1} is given as
$E(t)=u(t)v(t)^*$, where $u$ and $v$ solve the system of differential equations (for $G=G_\eps(E)$)
\begin{equation}\label{ode-uv}
\begin{array}{rcl}
 \dot u &=& -\tfrac \iu2 \, \Im(u^*Gv)u - (I-uu^*)Gv
\\[1mm]
 \dot v &=& -\tfrac \iu2 \, \Im(v^*G^*u)v - (I-vv^*)G^*u,
\end{array}
\end{equation}
which preserves $\|u(t)\|=\|v(t)\|=1$ for all $t$.
\end{lemma}

We note that for $G=G_\eps(E)=2f_{\clambda}\,xy^*$ (see Lemma~\ref{lem:gradient}) and with $\alpha=u^*x$, $\beta=v^*y$
and $\gamma=2f_{\clambda}$ we obtain the differential equations
\begin{equation}\label{ode-uv-short}
\begin{array}{rcl}
 \dot u &=& -\tfrac \iu2 \, \Im(\alpha\conj\beta\gamma)u +\alpha\conj\beta\gamma\, u- \conj\beta\gamma \,x
\\[3mm]
 \dot v &=& -\tfrac \iu2 \, \Im(\conj{\alpha}\beta\conj{\gamma})v +\conj{\alpha}\beta\conj{\gamma}\,v -\conj{\alpha\gamma}\,y.
\end{array}
\end{equation}

\begin{proof} We introduce the projection $\widetilde P_E$ onto the tangent space at $E=uv^*$ of the submanifold of rank-1 matrices of unit Frobenius norm,
$$
\widetilde P_E(G) = P_E(G-\Re\langle G,E \rangle E)=P_E(G)-\Re\langle G,E \rangle E.
$$
We find
\begin{align*}
\widetilde P_E(G) &= Gvv^* - uu^*Gvv^* + uu^*G -\Re\langle G,uv^* \rangle uv^*
\\
&= (I-uu^*)Gvv^* + uu^*G(I-vv^*) + uu^*Gvv^* - \Re(u^*Gv)uv^*
\\
&= (I-uu^*)Gvv^* + uu^*G(I-vv^*) + \iu \,\Im(u^*Gv)uv^*
\\
&= \Bigl(  \tfrac \iu2 \, \Im(u^*Gv)u + (I-uu^*)Gv\Bigr) v^* + u\Bigl( \tfrac\iu2\, \Im(u^*Gv)v^* + u^*G(I-vv^*) \Bigr).
\end{align*}
For $\dot E=\dot u v^* + u {\dot v}^*$ we thus have $\dot E=- \widetilde P_E(G)$ if $u$ and $v$ satisfy \eqref{ode-uv}.
Since then $\Re(u^*\dot u)=0$ and $\Re(v^*\dot v)=0$, the unit norm of $u$ and $v$ is preserved.
\qed
\end{proof}

%
%
%
%
%
%
%
%
%
The projected differential equation \eqref{ode-E-1} has the same monotonicity property as the differential equation \eqref{ode-E}.

\begin{theorem}[Monotonicity]
 \label{thm:monotone-C-1}
Let $E(t)\in \cR_1$ of unit Frobenius norm
be a solution to the differential equation \eqref{ode-E-1}.
If $\lambda(t)$ is a simple eigenvalue of $A+\eps E(t)$, then
\begin{equation}
\frac{ d }{dt} \F_\eps \bigl( E(t) \bigr)  \le  0.
\label{eq:pos-C-1}
\end{equation}
\end{theorem}


\begin{proof} As in the proof of Theorem~\ref{thm:monotone}, we abbreviate
$G=G_\eps(E)$ 
and obtain from \eqref{ode-E-1} and $\dot E \in T_E\cR_1$ and $\Re\langle E, \dot E\rangle =0$ that
$$
\| \dot E \|_F^2 = - \Re \bigl\langle P_E\bigl(G-\Re\langle G,E\rangle E\bigr),\dot E \bigr\rangle
= - \Re \bigl\langle G-\Re\langle G,E\rangle E,\dot E \bigr\rangle = - \Re\,\langle G, \dot E \rangle.
$$
Hence Lemma~\ref{lem:gradient}, \eqref{ode-E-1} and 
$\langle P_E(G), E \rangle = \langle G, P_E(E) \rangle = \langle G,E \rangle$
yield
\begin{align}
\frac1{ \eps\kappa } \,\frac{d}{dt} \F_\eps(E(t))  &=  \Re\, \langle G, \dot E \rangle = - \| \dot E \|_F^2
\nonumber
\\[-1mm]
&= - \bigl\| P_E\bigl(G-\Re\langle G,E\rangle E\bigr)\bigr\|_F^2
\label{c-s-1}
\\[1mm]
&= - \| P_E(G) \|_F^2 + \Re\langle G,E\rangle^2,
\nonumber
\end{align}
where the second equality yields the monotone decay and the last equality is noted for later use.
\qed \end{proof}

Comparing the differential equations \eqref{ode-E} and \eqref{ode-E-1} immediately shows that every stationary point of \eqref{ode-E} is also a stationary point of the projected differential equation \eqref{ode-E-1}. Remarkably, the converse is also true for the stationary points $E$ of unit Frobenius norm with ${P_E(G_\eps(E))\ne 0}$. Violation of this non-degeneracy condition is exceptional, as we will explain below.

\index{stationary point}
\begin{theorem}[Stationary points]
\label{thm:stat-1}
Let the rank-1 matrix $E\in \cR_1$ be of unit Frobenius norm and assume that $P_E(G_\eps(E))\ne 0$. If $E$ is a stationary point of the rank-1 projected differential equation \eqref{ode-E-1}, then $E$ is  already a stationary point of the differential equation \eqref{ode-E}.
\end{theorem}

\begin{proof} We show that $E$ is a nonzero real multiple of $G_\eps(E)$. By Theorem~\ref{thm:stat}, $E$ is then a stationary point of the differential equation \eqref{ode-E}.

For a stationary point $E$ of \eqref{ode-E-1}, we must have equality in \eqref{c-s-1}, which shows that $P_E(G)$ (again with $G=G_\eps(E)$) is a nonzero real multiple of $E$. Hence, in view of $P_E(E)=E$, we can write $G$ as
$$
G=\mu E + W, \quad\text{ where $\mu\ne 0$ is real and $P_E(W)=0$.}
$$
Since $E$ is of rank 1 and of unit Frobenius norm, $E$ can be written as $E=uv^*$ with $\| u \|=\|v\|=1$. We then have
 $$
 W=W-P_E(W)= (I-uu^*)W(I-vv^*).
 $$
 On the other hand, $G=2 \overline f_\lambda xy^*$  is also of rank 1. So we have
 $$
2 \overline f_\lambda xy^* = \mu uv^* + (I-uu^*)W(I-vv^*).
 $$
 Multiplying from the right with $v$ yields that $x$ is a complex multiple of $u$, and multiplying from the left with $u^*$ yields that $y$ is a complex multiple of $v$. Hence, $G$ is a complex multiple of $E$. Since we already know that $P_E(G)$ is a nonzero real multiple of $P_E(E)=E$, it follows that $G$ is the same real multiple of $E$. By Theorem~\ref{thm:stat}, $E$ is therefore a stationary point of the differential equation \eqref{ode-E}.
\qed \end{proof}
\begin{remark} [Non-degeneracy condition]\label{rem:exceptional}
\rm
\bng 
Let us discuss the possible failure of occurrence of condition 
$P_E(G_\eps(E))\ne 0$. \eng 
We recall that $G=G_\eps(E)$  is a multiple of $xy^*$, where $x$ and $y$ are left and right eigenvectors, respectively, to the simple eigenvalue $\lambda$ of $A+\eps E$. In which situation can we have $P_E(G)= 0$ \bng when \eng $G\ne 0$\,?

For $E=uv^*$, $P_E(G)= 0$ implies $G=(I-uu^*)G(I-vv^*)$, which yields $Gv=0$ and $u^*G=0$ and therefore $y^*v=0$ and $u^*x=0$. So we have $Ey=0$ and $x^*E=0$. This implies that $\lambda$ is already an eigenvalue of $A$ with the same left and right eigenvectors $x,y$ as for $A+\eps E$, which is an exceptional situation.
\end{remark}


\subsection{Numerical integration by a splitting method}\index{splitting method}
\label{subsec:proto-numer}
%
We need to integrate numerically the differential equations (\ref{ode-uv-short}), viz.
\begin{equation}\label{ode-uv-short-2} \nonumber
\begin{array}{rcl}
 \dot u &=& -\displaystyle \tfrac \iu2 \, \Im(\alpha\conj\beta\gamma)u +\alpha\conj\beta\gamma\, u- \conj\beta\gamma \,x
\\[3mm]
 \dot v &=& -\displaystyle \tfrac \iu2 \, \Im(\conj{\alpha}\beta\conj{\gamma})v +\conj{\alpha}\beta\conj{\gamma}\,v -\conj{\alpha\gamma}\,y,
\end{array}
\end{equation}
where $\alpha=u^* x\in \C$, $\beta=v^* y \in \C$ 
and $\gamma=2f_{\clambda} \in \C$.

The objective here is not to follow a particular trajectory accurately, but to arrive quickly at a stationary point.
The simplest method is the normalized Euler method, where the result after an Euler step (i.e., a steepest descent step) is normalized to unit norm for both the $u$- and $v$-component.
This can be combined with an Armijo-type line-search strategy to determine the step size adaptively. 

We found, however, that a more efficient method is obtained with a {\it splitting method} instead of the Euler method.
The splitting method consists of a first step applied to the differential equations
\begin{equation}\label{ode-uv-horiz}
\begin{array}{rcl}
 \dot u &=&  \alpha\conj\beta\gamma\, u- \conj\beta\gamma \,x  
\\[3mm]
 \dot v &=&  \conj{\alpha}\beta\conj{\gamma}\,v -\conj{\alpha\gamma}\,y  
\end{array}
\end{equation}
followed by a second step for the differential equations
\begin{equation}\label{ode-uv-rot} 
\begin{array}{rcl}
 \dot u &=& -\displaystyle \tfrac \iu2 \, \Im(\alpha\conj\beta\gamma)u 
\\[2mm]
 \dot v &=& - \displaystyle \tfrac \iu2 \, \Im(\conj{\alpha}\beta\conj{\gamma})v.
\end{array}
\end{equation}
As the next lemma shows, the first differential equation moves $\lambda$ in the direction of $-\gamma=-2f_{\clambda}$. In particular, the motion is horizontal if $f_{\clambda}$ is always real. The second differential equation is a mere rotation of $u$ and $v$.


%
%
%
%
%

\begin{lemma}[Eigenvalue motion in the direction of $-f_{\clambda}\,$]
\label{lem:gamma-motion}
\ Along a path of simple eigenvalues $\lambda(t)$ of $A+\eps u(t)v(t)^*$, where $u,v$ of unit norm solve \eqref{ode-uv-horiz}, we have that
\[
\text{$\dot\lambda(t)$ is a nonnegative real multiple of\/ $-\frac{\partial f}{\partial{\clambda}}(\lambda(t),\clambda(t))$.}
\]
\end{lemma}

\begin{proof} By Theorem~\ref{chap:appendix}.\ref{thm:eigderiv},
\begin{equation*}
\dot\lambda = \frac{1}{x^*y} \left( x^*\frac d{dt}(A+\eps  u v^*)\,y \right)= \eps \,\frac{x^* \left( {\dot u} v^* + u {\dot v}^*\right) y}{x^*y}. 
\end{equation*}
With $\alpha=u^* x$ and $\beta=v^*y$ and with $x,y$ normalized by \eqref{eq:scaling}, we obtain from \eqref{ode-uv-horiz}
\begin{equation*}
\frac{\dot \lambda}{\gamma} = -\frac{\eps}{x^*y}\Bigl(|\alpha|^2\cdot \left( 1 - |\beta|^2 \right) + 
|\beta|^2\cdot \left( 1 - |\alpha|^2 \right) \Bigr) \in \R,\ \le 0,
\end{equation*} 
which proves the statement, since $\gamma=2f_{\clambda}$.
\qed
\end{proof}
In general, a splitting method does not preserve stationary points. Here, it does.

\index{stationary point}
\begin{lemma}[Stationary points] \label{lem:stat-split}
If $(u,v)$ is a stationary point of the differential equations \eqref{ode-uv-short}, then it is also a stationary point of the differential equations
\eqref{ode-uv-horiz} and \eqref{ode-uv-rot}.
\end{lemma}

\begin{proof} If $(u,v)$ is a stationary point of \eqref{ode-uv-short}, then $u$ is proportional to $x$ and $v$ is proportional to $y$. Hence,
$x=\alpha u$ and $y=\beta v$. This implies that $(u,v)$ is a stationary point of \eqref{ode-uv-horiz}, and hence also of  \eqref{ode-uv-rot}.
\qed
\end{proof}

\subsubsection*{Fully discrete splitting algorithm.}
Starting from initial values $u_k,v_k$, we denote by $x_k$ and $y_k$ the left and right eigenvectors to the target 
eigenvalue $\lambda_{k}$ of $A+\eps u_k v_k^*$, and set 
\begin{equation}\label{alpha-beta-gamma-n}
\alpha_k=u_k^*x_k, \qquad \beta_k=v_k^*y_k, \qquad \gamma_k = 2f_{\clambda_k}. 
\end{equation}
We apply the Euler method with step size $h$  to \eqref{ode-uv-horiz} to obtain
\begin{equation}
\label{eul-horiz}
\begin{array}{rcl}
{\widehat u}(h) &=& u_k + h\left( \alpha_k\conj\beta_k\gamma_k\, u_k - \conj\beta_k\gamma_k \,x_k \right)
\\[2mm]
{\widehat v}(h) &=& v_k + h\left( \conj\alpha_k\beta_k\conj\gamma_k\,v_k -\conj{\alpha_k\gamma_k}\,y_k \right),
\end{array}
\end{equation}
followed by a normalization to unit norm 
\begin{equation} \label{eq:normal}
\widetilde u(h)=\frac{\widehat u(h)}{\|\widehat u(h)\|},\quad
\widetilde v(h)=\frac{\widehat v(h)}{\|\widehat v(h)\|}.
\end{equation}

Then, as a second step, we integrate the  rotating differential equations \eqref{ode-uv-rot} by setting, 
with $ \vartheta = -\displaystyle \tfrac12\, \Im \left( \alpha_k\conj{\beta_k}\gamma_k \right)$,
\begin{equation} \label{eq:rotate}
u(h)=\e^{\iu \vartheta h} \, \widetilde u(h), \qquad
v(h)=\e^{{}-\iu \vartheta h} \, \widetilde v (h),
\end{equation} 
and compute the target eigenvalue $\lambda(h)$ of $A+\eps u(h)v(h)^*$.
We note that this fully discrete algorithm still preserves stationary points. 
%
%

One motivation for choosing this method is that near a stationary point, the motion is almost 
rotational since $x \approx \alpha u$ and $y \approx \beta v$. 
The dominating term determining the motion is then the rotational term on the right-hand side of 
(\ref{ode-uv-short}), which is integrated by a rotation in the above scheme (the integration would be exact if $\alpha,\beta,\gamma$ were constant).

This algorithm requires in each step one computation of \bcl target \ecl eigenvalues and their left and right eigenvectors
of rank-$1$ perturbations to the matrix $A$; \bcl see the Notes at the end of this chapter. \ecl 

We also tried a variant where $\alpha,\beta,\gamma$ in the rotation step are updated from $( \widetilde u(h),\widetilde v(h))$ and the left and right eigenvectors to the target 
eigenvalue $\widetilde \lambda (h)$ of $A+\eps \widetilde u(h)\widetilde v(h)^*$. 
In our numerical experiments we found, however, that the slight improvement in the speed of convergence to the stationary state does not justify the nearly doubled
computational cost per step.

\medskip
\begin{algorithm}[H] \label{alg_prEul}
\DontPrintSemicolon
\KwData{$A, \eps, \theta > 1, u_k \approx u(t_k), v_k \approx v(t_k)$, $h_{k}$ (proposed step size)}
\KwResult{$u_{k+1}, v_{k+1}$, $h_{k+1}$}
\Begin{
\nl Initialize the step size by the proposed step size, $h=h_{k}$\; 
\nl Compute the value $f_k=f(\lambda_k,\conj{\lambda_k})$\;
\nl Compute left/right eigenvectors 
$x_k, y_k$ of $A + \eps u_k v_k^*$ to $\lambda_k$ such that $\| x_k \| = \| y_k \| = 1, x_k^* y_k > 0$\; 
\nl Compute $\alpha_k,\beta_k,\gamma_k$ by \eqref{alpha-beta-gamma-n} and $g_k$ by \eqref{g-n-formula}\;
\nl Initialize $f(h) = f_k$\;
\While{$f(h) \ge f_k$}{
\nl Compute $u(h), v(h)$ according to \eqref{eul-horiz}-\eqref{eq:rotate}\;
\nl Compute $\lambda(h)$ target eigenvalue of $A + \eps u(h) v(h)^*$\; 
\nl Compute the value $f(h) = f\bigl( \lambda(h), \conj{\lambda(h)} \bigr)$\;
\If{$f(h) \ge f_k$}{Reduce the step size, $h:=h/\theta$}
}
\uIf{$f(h) \ge f_k- (h/\theta) g_k$}{Reduce the step size for the next step, $h_{\rm next}:=h/\theta$}
\uElseIf{$h=h_k$}{Set $h_{\rm next} := \theta h_k$ (augment the stepsize if no rejection has occurred)}
\Else{Set $h_{\rm next} := h_k$}
\nl Set $h_{k+1}=h_{\rm next}$, $\lambda_{k+1}= \lambda(h)$, and the starting values for the next step as 
$u_{k+1}=u(h)$, $v_{k+1}=v(h)$\;
\Return
}
\caption{Integration step for the rank-1 constrained gradient system}
\end{algorithm}

\medskip

\subsubsection*{Step size selection.}\index{step-size selection}\index{Armijo line search}
We use an Armijo-type line search strategy to determine a step size that reduces the functional $f(\lambda,\clambda)$. For the non-discretized differential equation \eqref{ode-E-1}, we know from \eqref{c-s-1} that the decay rate is given by
$$
\frac1{\eps\kappa(t)}\frac{d}{dt} \F_\eps(E(t)) = - \Bigl( \| P_E(G) \|_F^2 - \bigl( \Re\,\langle G, E \rangle \bigr)^2 \Bigr) \le 0.
$$
Here we note that for $E=uv^*$ and again with $\alpha=u^*x$, $\beta=v^*x$, $\gamma=2f_{\clambda}$, so that
$G=\gamma xy^*$, we have
$$
\Re\langle G,E\rangle = \Re(\alpha\conj\beta\gamma)
$$
and 
$$
P_E(G)=\gamma(\alpha uy^* +\conj\beta xv^* - \alpha\conj\beta uv^*).
$$
A calculation shows that the squared Frobenius norm equals
$$
\| P_E(G)\|_F^2 = \langle P_E(G),P_E(G) \rangle = |\gamma|^2 \bigl(|\alpha|^2 + |\beta|^2 - |\alpha|^2  |\beta|^2\bigr) .
$$
We set
$$
g_k = \eps\kappa  
\Bigl( \| P_E(G) \|_F^2 - \bigl( \Re\,\langle G, E \rangle \bigr)^2 \Bigr) \ge 0
$$
for the choice $E=E_k=u_k v_k^*$ and $G=G_\eps(E_k)= 2f_{\clambda}(\lambda_k,\conj{\lambda_k})x_k y_k^*$.
In view of the above formulas, $g_k$ is computed simply as
\begin{equation}\label{g-n-formula}
g_k = \eps\kappa_k 
\Bigl(  
|\gamma_k|^2 \bigl(|\alpha_k|^2 + |\beta_k|^2 - |\alpha_k|^2  |\beta_k |^2\bigr) -  (\Re(\alpha_k\conj\beta_k\gamma_k))^2
\Bigr).
\end{equation}
Let 
$$
f_k = f(\lambda_k,\conj{\lambda_k}), \qquad f(h) = f(\lambda(h),\conj{\lambda(h)}).
$$
We accept the result of the step with step size $h$ if
$$
f(h) < f_k.
$$
If for some fixed $\theta>1$,
$$
f(h) > f_k - (h/\theta) g_k,
$$
then we reduce the step size for the next step to $h/\theta$. If the step size has not been reduced in the previous step, we try for a larger step size.
Algorithm \ref{alg_prEul} describes
the step from $t_{k}$ to $t_{k+1} = t_{k}+h_{k}$. 
\begin{figure}[!ht]\label{fig:split-1}
\vskip -20mm
\centerline{\hskip -5mm
\includegraphics[scale=0.4]{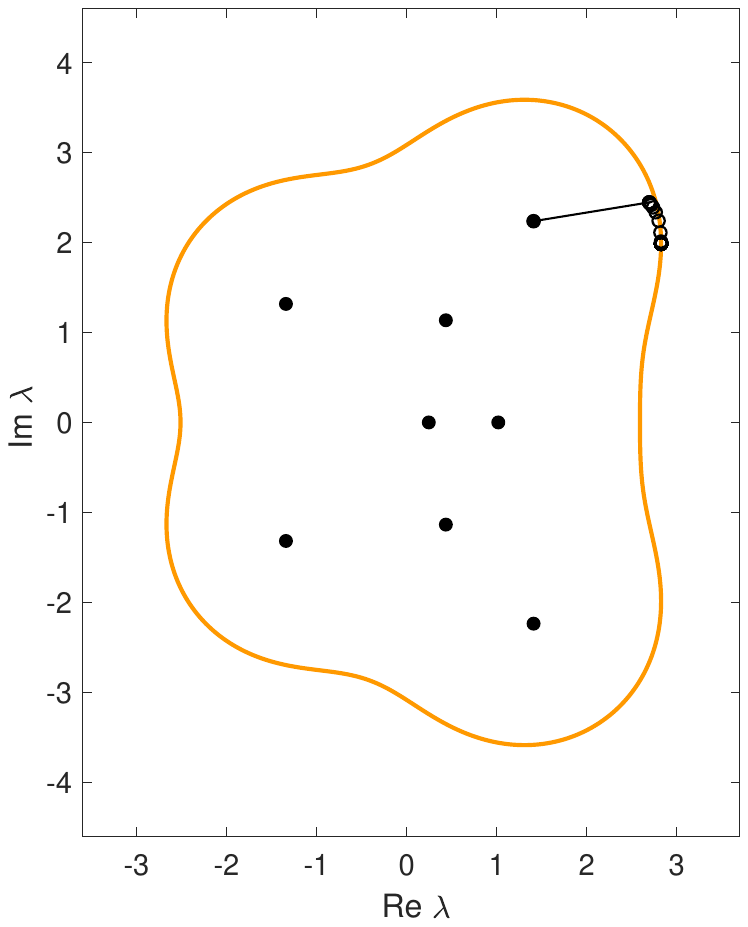} \hskip -20mm
\includegraphics[scale=0.4]{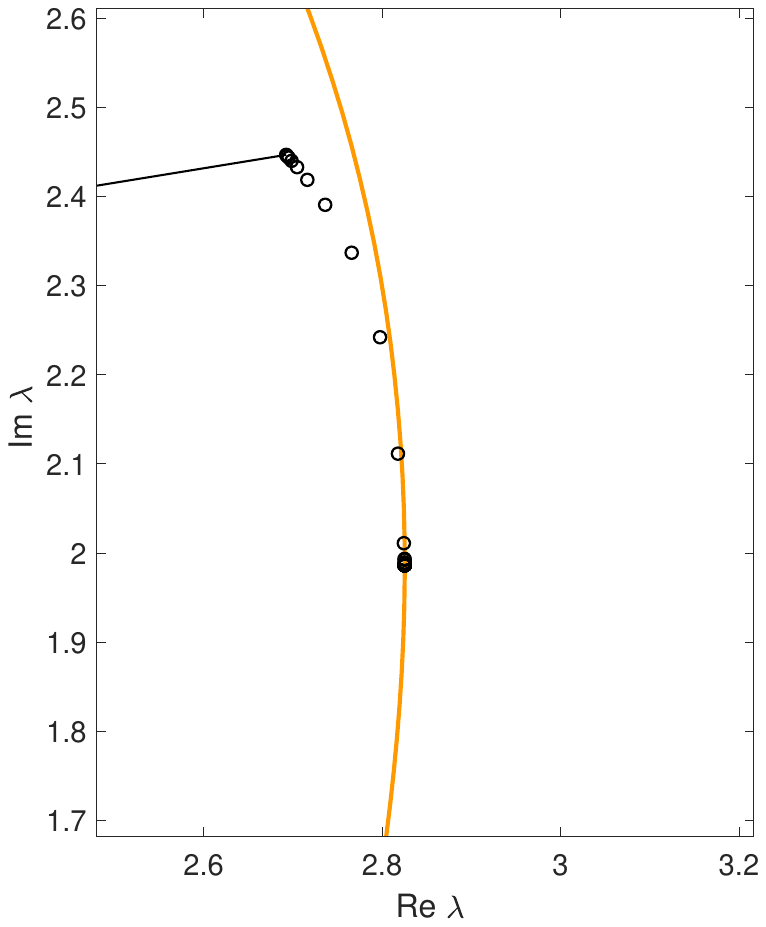}
}
 \vskip -20mm
\caption{Left: in black circles, iterates $\lambda_k$ of the splitting integrator (Algorithm \ref{alg_prEul}) applied to the matrix $A$ of \eqref{eq:example} with $f(\lambda,\clambda) = \Re\, \lambda$ and $\eps=1$. Right: zoom close to the stationary point.}
\end{figure}

\subsubsection*{Numerical example.}

An illustration is given in Fig.~3.1 
for $f(\lambda,\clambda) = - \frac12(\lambda + \clambda)=-\Re\,\lambda$
and for the randomly chosen $8 \times 8$ matrix
\begin{equation}
A = \left( 
\begin{array}{rrrrrrrr}
 0.91 &  1.17 & -0.80 &  0.34 &  0.52  &    0 & -1.39 & -0.28 \\
-0.05 &  0.54 &  1.91 &  1.68 &  1.67 &  1.38 &  1.62 &  2.50 \\
 1.03 & -1.35 & -1.29 &  0.55 & -1.37 & -0.26 &  0.33 & -0.89 \\
-0.27 & -1.05 & -0.87 &  0.99 & -1.23 &  0.04 & -0.11 & -0.62 \\
-0.68 &  0.65 &  1.01 &  0.65 &  0.78 &  0.80 & -0.18 & -0.24 \\
-0.16 & -0.52 &  0.26 & -0.61 & -0.10 & -0.04 &  0.22 &  0.37 \\
-0.67 &  0.17 & -0.69 &  2.23 & -0.23 &  0.94 &  0.19 & -0.22 \\
-1.43 &  0.13 & -0.89 &  0.06 &  1.26 &  0.28 &  0.05 &  0.03
\end{array}
\right).
\label{eq:example}
\end{equation}
The curve in Figure \ref{fig:split-1} is the set, with $\eps=1$ and the target eigenvalue $\lambda(M)$ the rightmost eigenvalue of a matrix $M$,
\begin{equation*}
\{  \lambda\left( A + \eps E \right) \in \C \,: \;  \text{$E \in \C^{n,n}$ with $\| E \|_F = 1$} \}.
\end{equation*}
With our choice of $f(\lambda)=-\Re\,\lambda$ we aim to find a rightmost point of this set. (This set is the boundary of the $\eps$-pseudospectrum of $A$. The problem of computing the real part of a rightmost point will be discussed in detail in Chapter~\ref{chap:pseudo}.)
\bng
The behaviour of Algorithm \ref{alg_prEul} applied to the matrix
\eqref{eq:example} is illustrated in Figure \ref{fig:split-1}.
The initial step size is set to $h=0.1$.
\eng



\bcl
\subsection{Choosing an initial perturbation}
\label{subsec:init}
\index{initial perturbation}

The choice of the initial perturbation $E(0)=u(0)v(0)^*$ may affect to which stationary point the trajectory of the rank-1 and norm-constrained  gradient flow converges. Some choices appear better than others when we aim to arrive at a global minimum or at least at a local minimum with a near-minimal value of the functional $\F_\eps$. 
We describe four possible choices of~$E(0)$ of increasing algorithmic complexity and reliability.

\medskip
(i) In a first approach, the target eigenvalue $\lambda$ of $A$ is computed and, assuming that $\lambda$ is a simple eigenvalue, also its left and right eigenvectors $x$ and $y$ with the normalization \eqref{eq:scaling} and the eigenvalue condition number $\kappa=1/(x^*y)$. Lemma~\ref{lem:gradient} (now with $\eps$ in the role of $t$) yields that for every matrix $E\in \C^{n,n}$ (of Frobenius norm~1),
$$
f\left(\lambda(A+\eps E),{\clambda(A+\eps E)} \right) = 
f\left(\lambda(A),\clambda(A) \right) + \eps\kappa \,\Re\langle G, E \rangle + O(\eps^2)
$$
with the scaled gradient $G= 2 f_{\clambda} xy^*$ as in \eqref{eq:freegrad}. If we momentarily ignore the $O(\eps^2)$ remainder term,  then
the miminum is attained for $E$ in the negative direction of the gradient, and so we choose
\[
E(0) = -\frac{G}{\| G \|_F}, \qquad \mbox{with} \quad
G = 2  \,\frac{\partial f}{\partial\clambda}(\lambda,\clambda) \, x y^*.
\]

(ii) Considering only the target eigenvalue of $A$ as in (i) may not be the best choice,
since a different eigenvalue that has a higher function value but is more sensitive to perturbations might result in a better choice of the initial perturbation. We proceed as in (i), but now for $m$ eigenvalues $\lambda_i$ ($i=1,\dots,m$) of $A$ for which $f(\lambda_i,\clambda_i)$ take values that are not far from that of the target eigenvalue of $A$, with corresponding normalized left and right eigenvectors $x_i$ and $y_i$ and eigenvalue condition number $\kappa_i=1/(x_i^*y_i)>0$. We set
\[
E_i = -\frac{G_i}{\| G_i \|_F}, \qquad \mbox{with} \quad
G_i = 2 \,\frac{\partial f}{\partial\clambda}(\lambda_i,\clambda_i) \,x_i y_i^*.
\]
We then choose 
\begin{equation} \label{eq:choice}
E(0) = E_j  \qquad \mbox{with} \quad
j= \arg\min_{i=1,\dots,m} \left( f\left( \lambda_i, \clambda_i \right)  - \eps\kappa_i \, 2\Bigl| \frac{\partial f}{\partial\clambda}(\lambda_i,\clambda_i) \Bigr|  \right),
\end{equation}
where we note that $ \Re\langle G_i, E_i \rangle = - \|G_i \|_F = - 2 |{\partial f}/{\partial\clambda}(\lambda_i,\clambda_i)|  $.

Approach (ii) is an excellent choice when $\eps$ is small so that the linearization gives a good approximation.

\medskip
\bng
(iii) An alternative approach uses the rank-1 constrained gradient flow without norm constraint, 
\begin{equation}\label{unconstrained-rank-1-flow}
 \dot E = -P_E G_\eps(E),
\end{equation}
starting with a small multiple of $E(0)$ of (ii) and integrating up to the first time $\bar t$ at which $E(\bar t)$ has Frobenius norm~1. 
We then take $E(\bar t)$ as the initial value of the rank-1 and norm-constrained gradient flow.

Note that with this approach, we track the eigenvalue
$\lambda_j$ in \eqref{eq:choice} from time $0$ until time $\bar t$. 
This means, for example, that if the target eigenvalue in the optimization problem is the rightmost one, here -- in the initial interval $[0,\bar t]$ -- we compute instead the path associated with the eigenvalue $\lambda(t)$ which originates from $\lambda_j(A)$, even if this is not the rightmost eigenvalue of $A$; once we reach the unit sphere of the Frobenius norm at $t=\bar t$, we switch to the constrained gradient system \eqref{ode-uv-short} and compute a trajectory of rightmost eigenvalues, originating at the rightmost eigenvalue of $A + \eps E(\bar t)$. 
\eng

\medskip
\bcl
(iv) In a further approach, we put a grid on a complex domain that covers the $\eps$-pseudo\-spectrum $\Lambda_\eps(A)$ (or a part of $\Lambda_\eps(A)$ in which we expect the eigenvalue minimizing~$f$ to occur).  Using the characterization $\Lambda_\eps(A) =\{ \lambda \in \C \,:\, 
 \sigma_{\min}(A-\lambda I) \le \eps \}$ (cf.~Chapter~\ref{chap:intro}), we determine those points on the grid which are in $\Lambda_\eps(A)$, as is done in the software package EigTool (Wright \cite{Wri02}, Wright \& Trefethen \cite{WriT01}). Among those points, we choose one (or a few as in (ii)) for which $f$ assumes the smallest (or nearby) value. For such a point $\lambda_0$, we compute the smallest singular value $\sigma_0$ of $A-\lambda_0 I$ and the associated left and right singular vectors $u_0$ and $v_0$. Then $\Delta_0 = - \sigma_0 u_0v_0^*$ is the matrix of smallest Frobenius norm such that $A-\lambda_0 I + \Delta_0$ is invertible (see again Chapter \ref{chap:intro}), or equivalently, such that $\lambda_0$ is an eigenvalue of $A+\Delta_0$. Since $\lambda_0\in \Lambda_\eps(A)$, we have
$\| \Delta_0 \|_F \le \eps$. We then take $E(0)= \Delta_0/\eps$, which has norm at most 1 as the initial perturbation matrix, and proceed with the concatenated unconstrained and constrained rank-1 gradient flows as in (iii). 

With a sufficiently fine grid, the approach (iv) eliminates the risk of getting stuck in a local instead of global minimum, since the initial value $E(0)$ is then already close to a global minimum, and so the gradient flow converges to this global minimum. On the other hand, for computational efficiency a very fine grid is to be avoided. 

Unlike (i)--(iii), the approach (iv) does not extend to structured eigenvalue optimization problems as considered in Chapter~\ref{chap:struc} and later chapters, since the singular value decomposition is not applicable to determine the nearest structured singular matrix. 
\ecl

\bng
\subsubsection*{An illustrative example.}
Consider $f(\lambda,\clambda)=-\Re\,\lambda$ to be minimized and the $3 \times 3$ complex matrix
\begin{equation} \label{ex:init}
A = \left( \begin{array}{rrr}
-1-\iu & \iu & 0 \\[2mm]
-2 + \iu & \frac12 & 1 + \iu \\[2mm]
0 & -\iu & \frac12 + 2 \iu
\end{array}
\right)
\end{equation}
having eigenvalues, ordered by decreasing real part,
\begin{equation} \nonumber
\begin{array}{l}
\lambda_1 = 0.34279 - 1.2522 \iu, \\[1mm]
\lambda_2 = 0.32691 + 0.94073 \iu, \\[1mm]
\lambda_3 = -0.66970 + 1.3115 \iu.
\end{array}
\end{equation}
Furthermore we have the associated condition numbers  
\[
\kappa_1 = 1.8478, \quad 
\kappa_2 = 2.9556.
\]
We set $\eps=10^{-0.4} = 0.39811$ and 
obtain:
\[
-\Re(\lambda_1) - \frac{\eps}{x_1^* y_1} = -1.0784 > 
-\Re(\lambda_2) - \frac{\eps}{x_2^* y_2} = -1.5035
\]
Then we integrate numerically
(\ref{ode-uv-short}) 
\bcl with the two different initial data as provided by the approaches (i) and (ii).

\begin{figure}[!h]\label{fig:init}
\vskip -20mm
\centerline{
\includegraphics[scale=0.45]{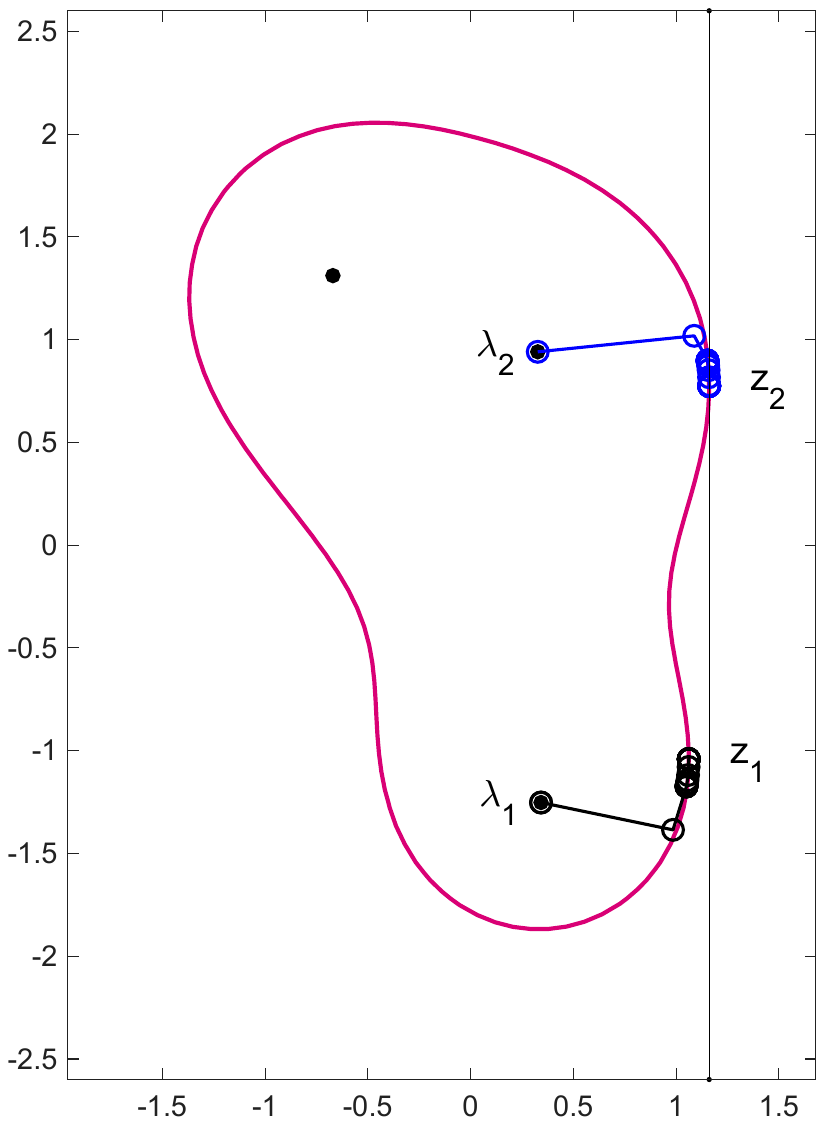} 
}
 \vskip -20mm
\caption{\bcl Eigenvalues computed by the splitting integrator (Algorithm \ref{alg_prEul}) applied to the matrix $A$ of \eqref{ex:init} with two different initial perturbations chosen according to (i) and (ii). Note that $\Re\, \lambda_1> \Re \,\lambda_2$ but $\Re \, z_2 > \Re \, z_1$.}
\end{figure}

\bng
In case (i) we obtain 
\ecl 
\bng
convergence to a stationary point $u_1 v_1^*$
such that the rightmost eigenvalue of $A+\eps u_1 v_1^*$ (which 
minimizes $f$) equals $z_1=1.0616 - 1.0405 \iu$.

\bcl
In case (ii) we obtain 
\ecl 
\bng
convergence to a stationary point $u_2 v_2^*$
such that the rightmost eigenvalue of $A+\eps u_2 v_2^*$ equals 
$z_2=1.1611 + 0.77508 \iu$.
Here the real part of ${}-\lambda$ is indeed minimized. 

\bng
In Figure \ref{fig:init} we show the boundary of the set
\[
\{ \lambda \in \C : \lambda \mbox{ is an eigenvalue of $A + \eps E$
for some \,} E \in \C^{3,3}, \| E \|_F = 1 \}
\]
and the two trajectories obtained by integrating system
\eqref{ode-uv-short} with initial data 
from the approaches (i) and (ii). It is seen that here (i) is not the best option.

Finally let us illustrate a few steps of approach (iii). Here we choose 
$E(0)= h E_2$ (see \eqref{eq:choice}), $h=0.01$ and with $E_2 = x_2 y_2^*$, where $x_2$ and $y_2$ are the normalized left and right eigenvectors associated to the eigenvalue $\lambda_2$ of $A$, which is ranked as second for largest real part. 
In this case we get for the spectrum of $A+\eps E(t)$ at times $0,h,2h$ 
\bcl and $\bar t \approx 1.01$
\bng the following:
\[
\begin{array}{cccc}
   t=0 & t=h & t=2h & t = \bar t \\
   0.3428 - 1.2522 \iu &  0.3394 - 1.2542 \iu &  0.3360 - 1.2562 \iu & 
   0.0461 - 1.3905 \iu \\
\bf   0.3269 + 0.9407 \iu &  \bf 0.3386 + 0.9407 \iu &  \bf 0.3501 + 0.9408 \iu & 
\bf 1.1399 + 0.9515 \iu \\
\rm   -0.6697 + 1.3115 \iu & -0.6766 + 1.3135 \iu & -0.6834 + 1.3154 \iu &
-0.9809 + 1.4390i 
\end{array} 
\]
In bold we denote the eigenvalue that we are tracking, which becomes
the rightmost one at $t=2h$.
\eng

\section{Notes}
The review article by Lewis and Overton (\cite{LeO96}) remains a basic reference on eigenvalue optimization, including a fascinating account of the history of the subject. There is, however, only a slight overlap of problems and techniques considered here and there.

The book by Absil, Mahony \& Sepulchre (\cite{AbMS09}) on optimization on matrix manifolds discusses alternative gradient-based methods to those considered here, though not specifically for eigenvalue optimization nor for low-rank matrix manifolds.

\subsubsection*{Rank-1 property of optimizers.}  \index{optimizer!rank-1 property}
The rank-1 structure of optimizers in an eigenvalue optimization problem was first used by Guglielmi \& Overton (\cite{GO11}) who devised a rank-1 matrix iteration
to compute the complex $\eps$-pseudospectral 
abscissa and radius; see Section~\ref{sec:psa} below.

The approach to eigenvalue optimization via a norm-constrained gradient flow and the associated rank-1 dynamics was first proposed and studied by Guglielmi \& Lubich (\cite{GL11}, \cite{GL12}), where it was used to compute the $\eps$-pseudospectral abscissa and radius as well as sections of the boundary of the $\eps$-pseudospectrum (see Chapter~\ref{chap:pseudo}). 
Our discussion of rank-1 dynamics in Section~\ref{subsec:rank1-gradient-flow} is based on Koch \& Lubich (\cite{KL07}). 

\subsubsection*{Frobenius norm vs.~matrix 2-norm.} \index{Frobenius norm!vs.~matrix 2-norm}
In the approach described in this chapter (and further on in this work), perturbations are measured and constrained in the Frobenius norm.  This choice is made because the Frobenius norm, unlike the matrix 2-norm, is induced by an inner product, which simplifies many arguments. Not least, it allows us to work with gradient systems. However, the approach taken here, with functional-reducing differential equations and their associated rank-1 dynamics, is also relevant for the matrix 2-norm, because the optimizers with respect to the Frobenius norm are of rank 1, and thus their Frobenius norm equals their 2-norm. Since the 2-norm of a matrix does not exceed its Frobenius norm, it follows that the rank-1 Frobenius-norm optimizers constrained by $\|\Delta\|_F\le \eps$ are simultaneously the 2-norm optimizers constrained by $\|\Delta\|_2\le \eps$. 

\bng
\subsubsection*{Computation of rightmost eigenvalues.}

For algorithms addressing the computation of rightmost eigenvalues and corresponding left and right eigenvectors of a large sparse matrix, we refer the reader to a few methods proposed in the literature.
All methods can be extended to the case of arbitrary target eigenvalues by using maps that transform the computational problem into that of computing the largest (in modulus) eigenvalues and associated eigenvectors.
Moreover, all methods can potentially exploit the sparse plus rank-one structure of the unstructured problem we consider in this chapter.

Meerbergen \& Roose (\cite{MeR96}) and Meerbergen, Spence \& Roose (\cite{MeSR94}) consider the Cayley transformation, which transforms the problem into that of computing the eigenvalues of largest modulus of a matrix, and give a review of many iterative methods.

Lehoucq,  Sorensen \& Yang (\cite{LeSY98}) propose 
 an implicitly restarted Arnoldi iteration 
 (which is implemented in ARPACK) to compute a small number of eigenvalues, characterized by specific properties, including the largest real part. They 
 propose a deflation procedure that is designed to improve the convergence of the Arnoldi iteration. In particular, as the iteration progresses, the Ritz value approximations of the eigenvalues converge at different rates.
 The implicitly restarted Arnoldi algorithm is highly effective, but suffers from two problems: the persistence of unwanted Ritz vectors caused by instabilities in the implicit QR step, and the challenge of deflating converged Ritz vectors due to the need to preserve the Arnoldi structure. 
 \index{implicitly restarted Arnoldi iteration}
 
 Stewart (\cite{Ste02}) proposes a generalized Krylov-Schur decomposition that overcomes the structural constraints that limit the transformations that can be applied to the Arnoldi decomposition, and
 efficiently resolves the two above-mentioned key difficulties of implicitly restarted Arnoldi methods.
 Stewart's algorithm is currently used in the MATLAB function {\it eigs}.
\index{Krylov--Schur method}
 
 \eng

%
%

\chapter{Pseudospectra}
\label{chap:pseudo}

\section{Complex $\eps$-pseudospectrum}
\label{sec:ps}
\index{pseudospectrum!complex}

\subsection{Motivation and definitions}
\label{subsec:ps-motivation}
As a motivating example for the $\eps$-pseudospectrum of a matrix $A\in\C^{n,n}$, we consider the linear dynamical system $\dot x(t)=Ax(t)$. The system is asymptotically stable, i.e., solutions $x(t)$ converge to zero as $t\goes\infty$ for all
initial data, if and only if all eigenvalues of $A$ have negative real part. We now ask for the robustness of asymptotic stability under (complex unstructured) perturbations $\Delta\in \C^{n,n}$ of norm bounded by a given $\eps>0$. This clearly depends on the choice of norm, and here we consider the Frobenius norm:
$$
\|\cdot\|=\|\cdot\|_F.
$$
For a {\it normal} matrix, the spectral decomposition yields that the perturbed system remains asymptotically stable for an arbitrary complex perturbation of norm at most $\eps$ if for each eigenvalue $\lambda$ of $A$, the real part is bounded by $\Re\,\lambda + \eps <0$. This condition is, however, not sufficient for {\it non-normal} matrices $A$. 

The question posed is thus: Is the following real number negative?
\begin{align*}
\alpha_\eps(A) = \max \{ \Re\,\lambda : \  \ &\text{There exists $\Delta\in\C^{n,n}$ with $\| \Delta \| \le \eps$ such that}
\\
&\text{$\lambda$ is an eigenvalue of $A+\Delta$} \}.
\end{align*}
This question is answered by solving a problem (\ref{chap:proto}.\ref{eq:optimiz0}) 
with the function to be minimized given by $f(\lambda,\clambda)=-\tfrac12(\lambda+\clambda)=-\Re\,\lambda$ (i.e., we maximize $\Re\,\lambda$).

It is useful to rephrase the question in terms of the $\eps$-pseudospectrum, for which we recall the definition; see also Chapter I and the notes in Section~\ref{sec:ps-notes}.

\begin{definition} \label{def:ps}
The  complex $\eps$-{\it pseudospectrum} of the matrix $A$ is the set
\begin{equation}\label{eq:epsps}
\Lambda_\eps(A) = \{ \lambda \in \C : \ \ \text{$\lambda \in \Lambda(A+\Delta)$ for some $\Delta\in\C^{n,n}$ with $\| \Delta \| \le \eps$}\},
\end{equation}
where $\Lambda(M)\subset \C$ denotes the spectrum (i.e., set of eigenvalues) of a square matrix $M$.
\end{definition}



\index{pseudospectral abscissa}
The above quantity $\alpha_\eps(A)$ can be rewritten more compactly as
\begin{equation}\label{stability-abscissa}
\alpha_\eps(A) = \max \{ \Re\,\lambda : \ \lambda \in \Lambda_\eps(A)\}.
\end{equation}
It is known as the  $\eps$-{\it pseudo\-spectral abscissa} of the matrix $A$. 

\index{pseudospectral radius}
An analogous quantity, of interest for discrete-time linear dynamical systems $x_{k+1}=Ax_k$, is the 
$\eps$-{\it pseudospectral radius} of the matrix $A$,
\begin{equation}\label{stability-radius}
\rho_\eps(A) = \max \{ | \lambda | : \ \lambda \in \Lambda_\eps(A)\}.
\end{equation}

\subsection{Pseudospectrum, singular values, and resolvent bounds}
The complex $\eps$-pseudospectrum can be characterized in terms of singular values. 
The singular
value decomposition of a matrix $M\in\C^{n,n}$ is $M=U\Sigma V^*$, with unitary matrices $U=(u_1,\dots,u_n)$ and $V=(v_1,\dots,v_n)$ formed by the left and right singular vectors $u_k$ and $v_k$, respectively, and with the real diagonal matrix $\Sigma =\diag(\sigma_1,\ldots,\sigma_n)$ of the singular values
$\sigma_1\geq\ldots\geq\sigma_n\ge 0$. 
We use the notation $\sigma_k(M)$ for the $k$th singular value of $M$ when we wish to indicate the dependence on $M$, and we write
$\sigma_{\min}(M)=\sigma_n(M)$ for the smallest singular value.

\begin{theorem}[Singular values and eigenvalues]
\label{thm:ps-sv} 
The complex $\eps$-pseudospectrum of $A \in \C^{n,n}$ is characterized as
\begin{align}\label{ps-sv}
\Lambda_\eps(A) 
&=\{ \lambda \in \C : \ \, \text{$\lambda \in \Lambda(A+\Delta)$\, for some $\Delta\in\C^{n,n}$ with $\| \Delta \| \le \eps$}\} \nonumber
\\
&= \{ \lambda \in \C : \ \sigma_{\min}(A-\lambda I) \le \eps \}. 
\end{align}
Moreover, the perturbation matrix $\Delta$ can be restricted to be of rank 1.
\end{theorem}

\begin{proof} The result relies on the fact that the distance to singularity 
\index{distance to singularity}
of a matrix $M$ equals its smallest singular value:
\begin{equation}\label{dist-sing}
\sigma_{\min}(M) = \min \{ \| \Delta \|: \ \Delta\in\C^{n,n} \text{ is such that }M+\Delta \text{ is singular}\}.
\end{equation}
The perturbation of minimal norm is then the rank-1 matrix
\begin{equation}\label{Delta-star}
\Delta_\star=-\sigma_n u_n v_n^*
\end{equation}
(unique if $\sigma_n<\sigma_{n-1}$), where $\sigma_n=\sigma_{\min}(M)$ and $u_n$, $v_n$ are the left and right $n$th singular vectors. This perturbation is such that $M+\Delta_\star$ has the same singular value decomposition as $M$ except that the smallest singular value is replaced by zero.

Choosing $M=A-\lambda I$ for $\lambda\in\C$ thus shows that $\sigma_{\min}(A-\lambda I) \le \eps$ if and only if there exists a matrix $\Delta\in\C^{n,n}$ of norm at most $\eps$ such that $A-\lambda I+\Delta$ is singular, or equivalently, that $\lambda$ is an eigenvalue of $A+\Delta$.
\qed
\end{proof}

Since $\sigma_{\min}(A-\lambda I)$ depends continuously on $\lambda$, Theorem~\ref{thm:ps-sv} implies that the boundary of the $\eps$-pseudospectrum of $A$ is given as
\begin{equation}\label{ps-sv-bdy}
\partial\Lambda_\eps(A) = \{ \lambda \in \C : \ \sigma_{\min}(A-\lambda I) = \eps \}.
\end{equation}

\begin{remark}[Frobenius norm and matrix 2-norm]\label{rem:ps-norms}\index{Frobenius norm!vs.~matrix 2-norm}
Since for  rank-1 matrices, the Frobenius norm and the matrix 2-norm are the same, Theorem~\ref{thm:ps-sv} and its proof show that the complex $\eps$-pseudospectra defined with respect to these two norms are identical. 
\end{remark}

Since $1/\sigma_{\min}(A-\lambda I)= \sigma_{\max}\bigl((A-\lambda I)^{-1}\bigr) = \| (A-\lambda I)^{-1} \|_2$,
we can reformulate \eqref{ps-sv} in terms of resolvents $(A-\lambda I)^{-1}$ as
\begin{equation}\label{ps-res}
\Lambda_\eps(A) = \{ \lambda \in \C : \ \| (A-\lambda I)^{-1} \|_2 \ge 1/\eps \} .
\end{equation}
This allows us to characterize the $\eps$-pseudospectral abscissa \eqref{stability-abscissa} as 
$$
\alpha_\eps(A)=\max \{ \Re\,\lambda\,:\, \| (A-\lambda I)^{-1} \|_2 \ge 1/\eps \}, 
$$
which implies
\begin{equation} \label{eps-res-bound}
    \frac 1\eps = \max_{\mathrm{Re}\,\lambda \ge \alpha_\eps(A)} \| (A-\lambda I)^{-1} \|_2.
\end{equation}
\index{stability radius}
\index{distance to instability}
If all eigenvalues of $A$ have negative real part, we define the {\it stability radius} (or {\it distance to instability}) as 
$$
\oeps>0\quad\text{such that}\quad  \alpha_{\oeps}(A)=0,
$$
i.e., there exists a perturbation $\Delta \in \C^{n,n}$  of Frobenius norm $\oeps$ such that $A+\Delta$ has an eigenvalue on the imaginary axis, as opposed to all perturbations of smaller norm. The above formula then yields that the inverse stability radius $1/\oeps$ is the smallest upper bound of the resolvent norm $\| (A-\lambda I)^{-1} \|_2$ for $\lambda$ in the complex right half-plane:
\index{resolvent bound}
\begin{equation} \label{oeps-res-bound}
    \frac 1\oeps = \max_{\mathrm{Re}\,\lambda \ge 0} \| (A-\lambda I)^{-1} \|_2.
\end{equation}
As we discuss next, the $\eps$-pseudospectral abscissa $\alpha_\eps(A)$ and the stability radius $\oeps$ are important quantities in bounding solutions of linear differential equations.

\subsection{Transient bounds for linear differential equations}
\label{subsec:ps-exp}
\index{transient bound}
We describe two approaches to bounding solutions to linear differential equations, one for the matrix exponential $\e^{tA}$, which corresponds to the homogeneous initial value problem $\dot x(t)=Ax(t)$ with an arbitrary initial value $x(0)=x_0$, and the other approach for the inhomogeneous problem $\dot x(t)=Ax(t) + f(t)$ with zero initial value.

\subsubsection*{Bounds for the matrix exponential.}
Via Theorem~\ref{thm:ps-sv}, the transient behaviour of $\| \e^{tA} \|_2$ can be bounded in terms of the complex pseudospectrum.
Here we illustrate this with a simple robust bound: Let $\Gamma_\eps\subset\C$ be the boundary curve of a piecewise smooth domain (or several non-overlapping domains) whose closure covers $\Lambda_\eps(A)$, and assume further that the real part of the rightmost point of $\Gamma_\eps$ equals the pseudospectral abscissa $\alpha_\eps(A)$. In particular, we may take $\Gamma_\eps=\partial\Lameps(A)$ when this is a piecewise regular curve. Using the Cauchy integral representation
$$
\e^{tA} = \frac1{2\pi\iu} \int_{\Gamma_\eps} \e^{t \lambda}\, (\lambda I -A)^{-1}\, d\lambda
$$
and noting that by \eqref{ps-res},
$\| (\lambda I -A)^{-1} \|_2  \le  1/\eps$  on $\Gamma_\eps$,
we find by taking norms that
\begin{equation}\label{transient-bound}
\| \e^{tA} \|_2 \le \frac{\gamma_\eps(t)}{\eps} \quad\text{ with}\quad  
\gamma_\eps(t)=  \frac1{2\pi} \int_{\Gamma_\eps} |\e^{t \lambda}|\, | d\lambda| \le \frac{|\Gamma_\eps|}{2\pi}\, \e^{t \alpha_\eps(A)} ,
\end{equation}
where $|\Gamma_\eps|$ is the length of $\Gamma_\eps$. This bound holds for every $\eps>0$.

The same argument can be applied to a perturbed matrix $A+\Delta$ with $\Delta\in\C^{n,n}$ bounded by $\|\Delta\|_2\le\delta< \eps$. 
\index{Weyl inequality}
The Weyl inequality $\sigma_{i+j+1}(B+C)\le \sigma_{i+1}(B)+\sigma_{j+1}(C)$ used with $i+1=n,\;j=0$
and $B=A-\lambda I+\Delta, \ C=-\Delta$ yields the lower bound
$$
\sigma_{\min}(A+\Delta-\lambda I)\ge \sigma_{\min}(A-\lambda I)- \delta \ge \eps-\delta \quad\text{for }\lambda\in\Gamma_\eps.
$$
We then obtain the robust transient bound
\index{transient bound!robust}
\begin{equation}\label{transient-perturbed}
\| \e^{t(A+\Delta)} \|_2 \le \frac{\gamma_\eps(t)}{\eps-\delta} \qquad\text{for  $\,\|\Delta\|_2\le\delta\ $ and for all $\,\eps>\delta$}.
\end{equation}
This bound can be optimized over $\eps>\delta$, provided that a bound for $|\Gamma_\eps|$ and an algorithm for computing the pseudospectral abscissa $\alpha_\eps(A)$ are available. Choosing $\eps$ as the stability radius, $\eps=\oeps$, we have $\alpha_{\oeps}(A)=0$, and so we obtain the time-uniform bound, for every $\delta<\oeps$,
\begin{align}\label{transient-perturbed-oeps}
\| \e^{t(A+\Delta)} \|_2 \le \frac{|\Gamma_{\oeps}|}{2\pi}\,\frac{1}{\oeps-\delta} \quad\ &\text{for all $t>0\ $ and
for all $\Delta\in \C^{n,n}$ with $\,\|\Delta\|_2\le\delta\ $.}
\end{align}


\subsubsection*{Bounds for linear inhomogeneous differential equations.}
 We consider the differential equation $\dot x(t)=Ax(t) + f(t)$ with zero initial value for
inhomogeneities $f\in L^2(0,\infty;\C^n)$, where we assume that all eigenvalues of $A$ have negative real part. We extend  $x(t)$ and $f(t)$ to $t<0$ by zero. Their Fourier transforms $\widehat x$ and $\widehat f$ are then related by
$\iu\omega\,\widehat x(\omega)=A\widehat x(\omega) + \widehat f(\omega)$ for all $\omega\in\R$, i.e.,
$$
\widehat x (\omega)=(\iu\omega I - A)^{-1}\widehat f(\omega), \qquad \omega \in \R,
$$
and hence the Plancherel formula yields
\begin{align*}
    &\int_\R \| x(t) \|^2 \, dt = \int_\R \| \widehat x (\omega)\|^2\, d\omega = 
    \int_\R \| (\iu\omega I - A)^{-1}\widehat f(\omega) \|^2 \, d\omega
    \\
    &\le \max_{\omega\in\R} \| (\iu\omega I - A)^{-1} \|_2^2 \, \int_\R \| \widehat f (\omega)\|^2\, d\omega =
    \max_{\omega\in\R} \| (\iu\omega I - A)^{-1} \|_2^2 \,  \int_\R \| f(t) \|^2 \, dt.
\end{align*}
Using \eqref{oeps-res-bound} and the causality property that $x(t)$, for $0\le t \le T$, only depends on $f(\tau)$ with $0\le \tau \le t \le T$ (which allows us to extend $f(t)$ by $0$  for $t>T$), we thus obtain the $L^2$ transient bound
\begin{equation} \label{oeps-L2-bound}
  \biggl(  \int_0^T \| x(t) \|^2 \, dt \biggr)^{1/2} \le \frac1\oeps
   \biggl(  \int_0^T \| f(t) \|^2 \, dt \biggr)^{1/2} , \qquad 0 \le T \le \infty,
\end{equation}
where $\oeps$ is the stability radius of $A$.

For perturbed differential equations $\dot x_\Delta(t)=(A+\Delta)x_\Delta(t) + f(t)$ with zero initial value with 
perturbation size $\|\Delta\|_2\le\delta<\oeps$ we obtain, by the same argument and using the Weyl inequality as before, the robust $L^2$ transient bound
\index{transient bound!robust}
\begin{equation} \label{oeps-L2-bound-robust}
  \biggl(  \int_0^T \| x_\Delta(t) \|^2 \, dt \biggr)^{1/2} \le \frac1{\oeps-\delta}
   \biggl(  \int_0^T \| f(t) \|^2 \, dt \biggr)^{1/2} , \qquad 0 \le T \le \infty.
\end{equation}


\medskip

\subsection{Extremal perturbations}
\index{extremal perturbation}
The proof of Theorem~\ref{thm:ps-sv} also yields the following result on the perturbations $\Delta$ of minimal norm $\eps$ such that $A+\Delta$ has a prescribed eigenvalue $\lambda$ on the boundary $\partial \Lambda_\eps(A)$ of the complex $\eps$-pseudospectrum of $A$.

\pagebreak[3]

\begin{theorem}
[Extremal complex perturbations] \label{thm:Delta-C}
Let $\lambda\in\partial \Lambda_\eps(A)$, and let  $\Delta\in\C^{n,n}$ of norm $\eps$ be such that $A+\Delta$ has the eigenvalue $\lambda$. Then, $\Delta$ is of rank $1$.
%

Assume now that $\eps$ is a {\em simple} singular value of $A-\lambda I$ and that the corresponding left and right singular vectors are not orthogonal to each other.  Then,
$$
\Delta = \eps\e^{\iu\theta} x y^*,
$$
where $\e^{\iu\theta}$ is the outer normal to $\partial \Lambda_\eps(A)$ at $\lambda$, which is uniquely determined, and $x$ and $y$ are left and right eigenvectors of $A+\Delta$ to the eigenvalue~$\lambda$, of unit norm and with $x^*y>0$.
\end{theorem}

\begin{proof} By \eqref{ps-sv-bdy}, we have $\sigma_{\min}(A-\lambda I)=\eps$. The proof of Theorem~\ref{thm:ps-sv} then shows that
$
\Delta=-\eps uv^*,
$
where $u$ and $v$ are left and right singular vectors of $A-\lambda I$, with
$$
(A-\lambda I+\Delta)v=0 \quad\text{ and } \quad u^*(A-\lambda I+\Delta)=0,
$$
or equivalently,
$(A+\Delta)v=\lambda v$ and $u^*(A+\Delta)=\lambda u^*$.
This shows that $u$ and $v$ are left and right eigenvectors of $A+\Delta$.

Assume now that $\eps$ is a simple singular value of $A-\lambda I$ and that $u^*v\ne 0$. We show that $\partial \Lambda_\eps(A)$ has the outer normal $\e^{\iu\theta}$ at $\lambda$, where the angle $\theta$ is determined by
$$
u^*v = - \rho \,\e^{-\iu\theta}, \quad\rho>0.
$$
Let $\gamma(t)$, for $t$ near $0$, be a path in the complex plane with $\gamma(0)=\lambda\in \partial \Lambda_\eps(A)$. 
With $\nu=\dot\gamma(0)$ we have, by the derivative formula of simple singular values (see Corollary~\ref{chap:appendix}.\ref{lem:singderiv}),
\begin{align*}
\frac d{dt}\bigg|_{t=0} \sigma_{\min} \bigl(A-\gamma(t) I\bigr) &= 
\bcl
\Re \bigl(u^* (-\dot \gamma(0) I ) v \bigr) =
- \Re(\nu \,u^*v) = \rho\, \Re(\nu \e^{-\iu\theta}).
\ecl
\end{align*}
This shows that $\nu=\e^{\iu\theta}$ is the unique direction of steepest ascent, which is orthogonal to the level set $\partial \Lambda_\eps(A)$ and points out of $\Lameps(A)$. Hence,
$\e^{\iu\theta}$ is the outer normal to $\partial \Lambda_\eps(A)$ at $\lambda$.

We set $x=-\e^{-\iu\theta}u$ and $y=v$, which gives us a pair of left and right eigenvectors of $A+\Delta$ with $x^*y=\rho>0$. We then have
$$
\Delta = -\eps uv^* = \eps \e^{\iu\theta} xy^*,
$$
which proves the result.
\qed
\end{proof}

As we will see, Theorems~\ref{thm:ps-sv}  and~\ref{thm:Delta-C} motivate different approaches to computing the boundary of the $\eps$-pseudospectrum (or just extremal points such as a rightmost point): methods that steer the smallest singular value of $A-\lambda I$ to $\eps$, and methods that iterate on rank-1 matrices.

\section{Computing the pseudospectral abscissa} \label{sec:psa}
\index{pseudospectral abscissa}
As we have seen in the previous section, the $\eps$-pseudospectral abscissa $\aleps(A)$ 
is important for ensuring robust stability of a linear
dynamical system.
Various intriguing algorithms based on different ideas have been proposed to compute the pseudospectral abscissa:
\begin{itemize}
\item  the criss-cross algorithm of Burke, Lewis \& Overton (\cite{BuLeOv03}), which is based on Theorem~\ref{thm:ps-sv} and on Byers' Lemma given below;
\item  the rank-1 iteration of Guglielmi \& Overton (\cite{GO11}), which is based on Theorem~\ref{thm:Delta-C};
\item  the rank-1 constrained gradient flow algorithm of Guglielmi \& Lubich (\cite{GL11}); this is the approach presented in~Chapter~\ref{chap:proto} for $f(\lambda,\clambda)=-\tfrac12(\lambda+\clambda)=-\Re\,\lambda$;
\item  the subspace method of Kressner \& Vandereycken (\cite{KV14}).
\end{itemize}

\subsection{Criss-cross algorithm}
\label{subsec:criss-cross}
\index{criss-cross algorithm}

This remarkable algorithm was proposed and analysed by  Burke, Lewis \& Overton (\cite{BuLeOv03}). It uses a sequence of vertical and horizontal
searches in the complex plane to identify the intersection of a
given line with $\partial \Lambda_\eps(A)$. Horizontal searches yield updates to the approximation of $\alpha_\eps(A)$ while vertical searches find
favourable locations for the horizontal searches.

The criss-cross algorithm computes a monotonically growing sequence $(\alpha^k)$ that converges to the complex $\eps$-pseudospectral abscissa $\alpha_\eps(A)$. In its basic form, it can be written as follows:

0. Initialize $\alpha^0=\max\{\Re\,\lambda\,:\, \lambda\in \Lambda(A) \}$.

1. For $k=0,1,2,\dots$ iterate

1.1 (Vertical search) 
\begin{align}
\nonumber
&\text{Find all real numbers $\beta_j$, in increasing order for $j$ from $0$ to $m$,} 
\\
\label{cc-vertical}
&\text{such that $\alpha^k+\iu\beta_j \in \partial \Lambda_\eps(A)$.}
\end{align}

1.2 (Horizontal search) For $j=0,\dots,m-1$,
\begin{align} 
\nonumber
&\text{let the midpoint $\beta_{j+1/2}=\tfrac12(\beta_j+\beta_{j+1})$;}
\\
\nonumber
&\text{if $\alpha^k+\iu \beta_{j+1/2} \in \Lambda_\eps(A)$,
find the largest real number $\widehat\alpha_{j+1/2}$}
\\
\label{cc-horizontal}
&\text{such that $\widehat\alpha_{j+1/2}+\iu \beta_{j+1/2} \in \partial\Lambda_\eps(A)$ .}
\end{align}

1.3 Take $\alpha^{k+1}$ as the maximum of the $\widehat\alpha_{j+1/2}$.


\bigskip\noindent
To turn this into a viable algorithm, \eqref{cc-vertical} and \eqref{cc-horizontal} need to be computed efficiently. This becomes possible thanks to the following basic lemma, applied to $A$ for \eqref{cc-vertical} and to the rotated matrix $\iu A$ for  
\eqref{cc-horizontal}. 

\index{Byers' lemma}
\begin{lemma}[Byers' Lemma]
\label{lem:byers} Let $A\in\C^{n,n}$.
For given real numbers
$\alpha$ and $\beta$, the number $\eps > 0$ is
a singular value of the matrix
\[
A - (\alpha + \iu \beta ) \Id
\]
if and only if\/ $\iu \beta$ is an eigenvalue of the Hamiltonian matrix
\begin{equation}\label{eq:H} 
H(A,\alpha) =
   \left( \begin{array}{cc}   -(A-\alpha I)^* & \eps \Id \\
  -\eps \Id & A- \alpha \Id \end{array} \right).
\end{equation}
\end{lemma}

\begin{proof} 
\bng Without loss of generality we can assume $\alpha=0$ in the following.
\eng
The imaginary number $\iu \beta$ is an eigenvalue of the Hamiltonian matrix \eqref{eq:H} if and only if there exist nonzero vectors $u$ and $v$ such that
\begin{equation} \label{eq:eigH}
\left( \begin{array}{cc} -A^* & \eps \Id \\
  {}-\eps \Id & A  \end{array} \right)
\left( \begin{array}{c} u \\ v \end{array} \right) =
\iu \beta \left( \begin{array}{c} u \\ v \end{array} \right). 	
\end{equation}
This is equivalent to
\begin{equation}
\left( A - \iu \beta \Id \right)^* u = \eps v, \qquad
\left( A - \iu \beta \Id \right) v = \eps u,
\label{eq:uv}
\end{equation}
which expresses that $\eps$ is a singular value of $A - \iu \beta \Id$.
\qed
\end{proof}

Using Lemma~\ref{lem:byers}, the {\it vertical search} \eqref{cc-vertical} is done by computing all purely imaginary eigenvalues $\iu\beta$ of the Hamiltonian matrix and discarding those eigenvalues among them for which $\eps$ is not the smallest singular value of $A-(\alpha_k+\iu\beta)\Id$.
The {\it horizontal search} \eqref{cc-horizontal} is done by computing the purely imaginary eigenvalue $\iu\widehat\alpha_{j+1/2}$ of largest imaginary part of the Hamiltonian matrix that corresponds to the matrix $\iu A$ in the role of $A$ and $- \beta_{j+1/2}$ in the role of $\alpha$, \bng as is described in \eng 
Algorithm \ref{alg_cc}.

The computational cost of an iteration step of the criss-cross algorithm is thus determined by computing imaginary eigenvalues of  complex Hamiltonian $2n\times 2n$ matrices. All imaginary eigenvalues are needed in the vertical search and the ones of largest imaginary part in the horizontal search. In addition, the smallest singular values of complex $n\times n$ matrices need to be computed to decide if a given complex number is in $\Lameps(A)$.

\begin{algorithm}
\DontPrintSemicolon
\KwData{Matrix $A$, $\eps>0$, ${\rm tol}$ a given positive tolerance} 
\KwResult{$\aleps(A)$}
\Begin{
\nl Set $k=1$, $\alpha^0=\alpha(A)$\;
\While{$\alpha^k-\alpha^{k-1} > {\rm tol}$}{
\nl Vertical iteration.\; 
\nl Find all imaginary eigenvalues $\{\widehat \beta_j\}$ of $H(A,\alpha^k)$ (see \eqref{eq:H})
and collect those for which $\eps$ is the {\em smallest} singular value of $A - (\alpha^k +\iu \widehat \beta_j)$
as $\{\beta_j \}_{j=0}^{m}$, 
\begin{equation*}
\beta_0 \le \beta_1 \le \beta_2 \le \ldots \le \beta_{m-1} \le \beta_{m}.
\end{equation*}
\nl \ \hspace{-6mm} Horizontal iteration.\;
\For{$j=0,\ldots,m-1$}{ 
\nl Compute the midpoints $\displaystyle{\beta_{j+1/2} = \frac{\beta_{j}+\beta_{j+1}}{2}}$\; 
Find the highest imaginary eigenvalue (i.e. with largest imaginary part) $\iu \alpha_{j+1/2}$ 
of $H(\iu A,\beta_{j+1/2})$ (see \eqref{eq:H})\;
}
\nl Set $k=k+1$\;
\nl Set $\alpha^{k+1} = \max\limits_{j=0,\ldots,m-1} \alpha_{j+1/2}$\;
}
}
\caption{Criss-cross algorithm}
\label{alg_cc} 
\end{algorithm}

\begin{figure}[ht]\label{fig:cc}
\centerline{
\includegraphics[scale=0.36,keepaspectratio]{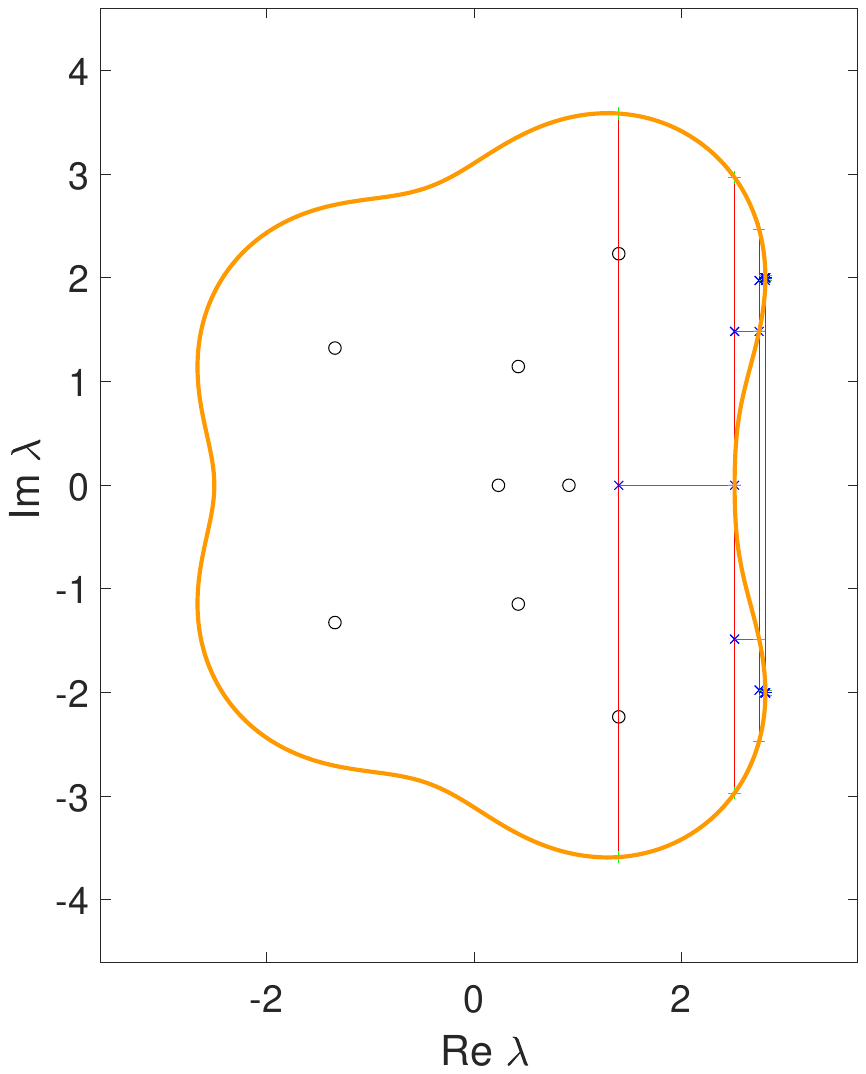}
\hskip -10mm
\includegraphics[scale=0.36,keepaspectratio]{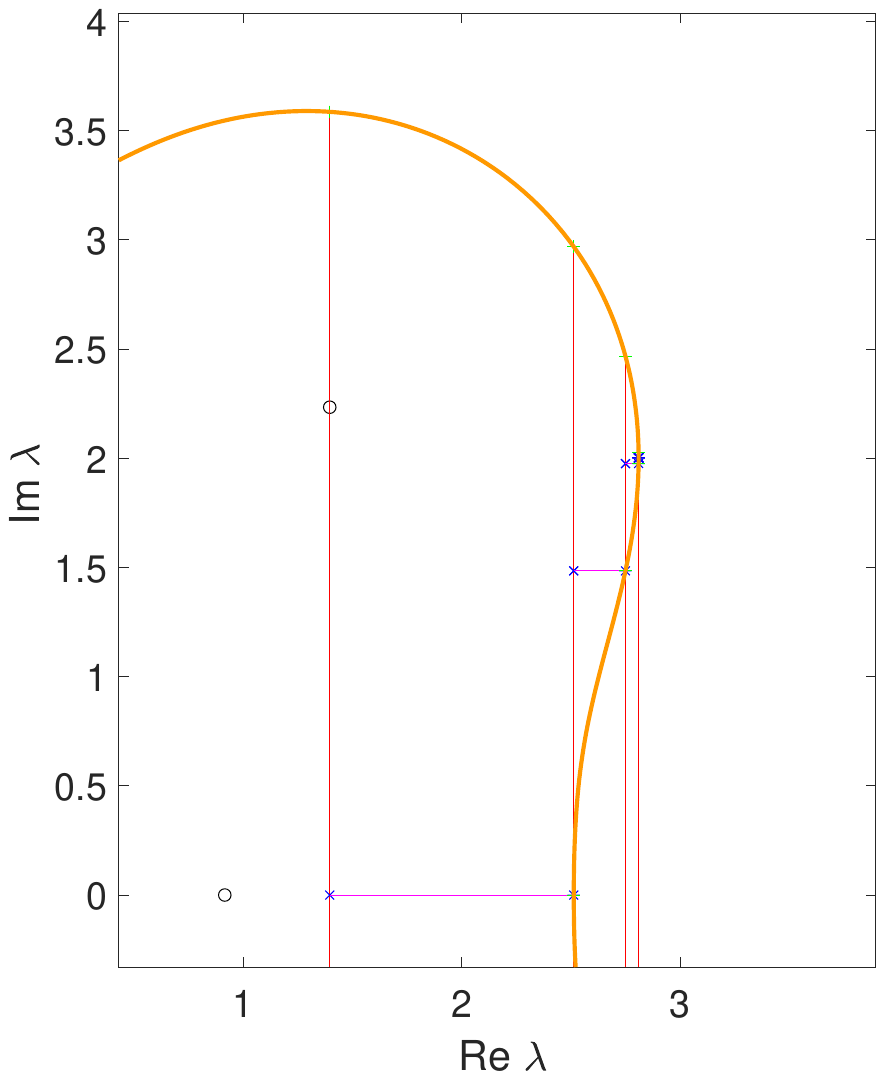}
}
\vskip -10mm
\caption{Iterates of the criss cross algorithm for computing $\alpha_\eps(A)$ for the matrix $A$ of (\ref{chap:proto}.\ref{eq:example})
and $\eps=1$. Right: zoom close to a rightmost point}
\end{figure}

An illustration is given in Fig.~\ref{fig:cc} for the  $8 \times 8$ matrix in (\ref{chap:proto}.\ref{eq:example}).

\subsubsection*{Unconditional convergence.} As the following theorem by Burke, Lewis \& Overton (\cite{BuLeOv03}) shows, the sequence generated by the criss-cross algorithm always converges to the 
$\eps$-pseudospectral abscissa.

\begin{theorem}[Convergence of the criss-cross algorithm]
\label{thm:cc-conv}
For every matrix $A\in\C^{n,n}$, the sequence $(\alpha^k)$ of the criss-cross algorithm converges to the pseudospectral abscissa $\alpha_\eps(A)$.
\end{theorem}
\bcl We give an alternative proof to that of Burke, Lewis \& Overton (\cite{BuLeOv03}).
\ecl

\begin{proof} By construction, the sequence $(\alpha^k)$ is a monotonically increasing sequence of real parts of points on $\partial\Lameps(A)$, and it is bounded as $\Lameps(A)$ is bounded. Therefore, $(\alpha^k)$ converges to a limit $\alpha^\star$, 
\bcl 
which is a fixed-point of the criss-cross iteration and is
\ecl
the real part of some point on $\partial\Lameps(A)$.
 Hence, $\alpha^\star \le \alpha_\eps(A)$. It remains to show that actually 
$\alpha^\star = \alpha_\eps(A)$.

To this end, we use the fact that every path-connected component of $\Lameps(A)$ contains an eigenvalue of $A$. This is readily seen as follows: For any $\lambda_1\in\Lameps(A)$, there exists a matrix $\Delta\in\C^{n,n}$ of norm at most $\eps$ such that $\lambda_1$ is an eigenvalue of~$A+\Delta$. Consider now the path $A+\theta\Delta$,  $0\le\theta\le 1$. By the continuity of eigenvalues, to this path corresponds a path of eigenvalues
$\lambda(\theta)$ of $A+\theta\Delta$ with $\lambda(1)=\lambda_1$, which connects $\lambda_1$ with the eigenvalue $\lambda(0)$ 
of~$A$.

Suppose $\alpha^\star < \alpha_\eps(A)$. We show that this leads to a contradiction. Let $\lambda_1\in\Lameps(A)$ be such that $\Re\,\lambda_1=\alpha_\eps(A)$. Then there is a path $\lambda(\theta)$ to an eigenvalue $\lambda_0$ of $A$, which by construction has a real part that does not exceed $\alpha^0$ and hence is smaller than $\alpha^\star$. So there exists $\theta^\star \in (0,1)$ such that $\lambda(\theta^\star)\in \Lambda_{\theta^\star\!\eps}(A)$ has real part $\alpha^\star$, that is, $\lambda(\theta^\star)=\alpha^\star+\iu\beta^\star$ for some real $\beta^\star$. There exists a smallest interval $[\beta_0,\beta_1]$ that contains $\beta^\star$ and has boundary points such that
$\alpha^\star+\iu \beta_0, \alpha^\star+\iu \beta_1\in\partial\Lameps(A)$. Then, the points $\alpha^\star+\iu \beta$ with $\beta_0<\beta<\beta_1$
are in the interior of $\Lameps(A)$, and in particular this holds true for the midpoint $\beta=\tfrac12(\beta_0+\beta_1)$. Hence there exists a maximal $\wh \alpha>\alpha^\star$ such that $\wh \alpha+\iu\beta\in \partial\Lameps(A)$. 
\bcl 
We thus find an $\wh \alpha>\alpha^\star$ that is equal to or bounded from above by the result of a criss-cross iteration starting from $\alpha^\star$. This contradicts the fact that $\alpha^\star$ is a fixed-point of the criss-cross iteration. 
\bcltwo
So, the assumption $\alpha^\star < \alpha_\eps(A)$, has led to a contradiction and
\ecltwo
\ecl
we must have $\alpha^\star = \alpha_\eps(A)$.
\qed 
\end{proof}

\subsubsection*{Locally quadratic convergence.}
We here show that the criss-cross algorithm converges locally quadratically under the following regularity assumption:
\begin{equation}\label{right-reg}
\begin{aligned}
&\text{At every right-most point of the $\eps$-pseudospectrum of $A$, the} 
\\[-1mm]
&\text{boundary curve $\partial\Lameps(A)$ is smooth with nonzero curvature.}
\end{aligned}
\end{equation}
This condition is stronger than the condition of a simple smallest singular value of $A-\lambda I$ at right-most points $\lambda\in\Lameps(A)$ imposed by Burke, Lewis \& Overton (\cite{BuLeOv03}), but it allows for a short proof of locally quadratic convergence.

\begin{theorem}[Locally quadratic convergence of the criss-cross algorithm]
\label{thm:cc-conv-quad}
Under condition \eqref{right-reg}, the sequence $(\alpha^k)$ of the criss-cross algorithm converges locally quadratically to $\alpha^\star=\alpha_\eps(A)$\,:
$$
0 \le \alpha^\star - \alpha^{k+1} \le C (\alpha^\star - \alpha^{k})^2,
$$
where $C$ is independent of $k$, provided that $\alpha^{k}$ is sufficiently close to $\alpha^\star$.
\end{theorem}

\begin{proof}
Near a right-most boundary point $\alpha^\star+\iu\beta^\star$ of $\Lameps(A)$, boundary points $\alpha+\iu\beta$ are related by
$$
\alpha^\star-\alpha = f(\beta)= c^2(\beta^\star-\beta)^2 + O((\beta^\star-\beta)^3),
$$
where $c\ne 0$ by condition \eqref{right-reg}. For the variables $\eta=\alpha^\star-\alpha$ and $\xi=c(\beta^\star-\beta)$ this relation becomes
$$
\eta= \phi(\xi)=\xi^2 + O(\xi^3).
$$
For a small $\delta>0$, let now $\xi_+=\delta$ and choose $\xi_-=-\delta+O(\delta^2)$ such that $\phi(\xi_-)=\phi(\xi_+)=\delta^2(1+O(\delta))$. Then,
$$
\tfrac12(\xi_+ + \xi_-)=O(\delta^2) \quad\text{ and hence } \quad
\phi\bigl(\tfrac12(\xi_+ + \xi_-)\bigr) = O(\delta^4),
$$
which yields
$$
\phi\bigl(\tfrac12(\xi_+ + \xi_-)\bigr) \le C \phi(\xi_+)^2.
$$
Translated back to the original variables, this yields the stated result.
\qed
\end{proof}

\subsection{Iteration on rank-1 matrices}
\label{subsec:r1-iteration}
Guglielmi and Overton (\cite{GO11}) proposed a strikingly simple iterative algorithm for computing the pseudospectral abscissa that 
uses a sequence of rank-1 perturbations of the matrix. Working with rank-1 perturbations appears natural in view of Theorem~\ref{thm:Delta-C}. Moreover, this theorem (with $\theta=0$) shows that at a point $\lambda\in\partial\Lameps(A)$ such that $\Re \lambda =\alpha_\eps(A)$, where the outer normal is horizontal to the right, the corresponding matrix perturbation $\Delta$ of norm $\eps$ is such that $\Delta=\eps xy^*$, where $x$ and $y$ are left and right eigenvectors, of unit norm and with $x^*y>0$, to the eigenvalue $\lambda$ of $A+\Delta$. This motivates the following fixed-point iteration.

\subsubsection*{Basic rank-1 iteration.} 
The basic iteration starts from two vectors $u_0$ and $v_0$ of unit norm and runs as follows for $k=0,1,2,\dots$: 

Given a rank-$1$ matrix
$E_k=u_k v_k^*$ of unit norm, compute the rightmost eigenvalue 
$\lambda_k$ of $A+\eps E_k$ and left and right eigenvectors $x_k$ and $y_k$, of unit norm and with $x_{k}^*y_{k}>0$, and set 
 $E_{k+1} = u_{k+1}v_{k+1}^* := x_{k} y_{k}^*$. 
 
Algorithm \ref{alg_GO} gives a formal description.
This algorithm requires in each step one computation of rightmost eigenvalues and associated eigenvectors
of rank-$1$ perturbations to the matrix $A$, which can be computed at relatively small computational cost for 
large sparse matrices $A$; \bng see the references in the Notes of Chapter \ref{chap:proto}. \eng 

\begin{algorithm}
\DontPrintSemicolon
\KwData{Matrix $A$, $\eps>0$, ${\rm tol}$ a given positive tolerance} 
\KwResult{$r \le \aleps(A)$, $x,y$}
\Begin{
\nl Compute $x_0$ and $y_0$ left and right eigenvectors to the rightmost eigenvalue of $A$, 
both normalized to unit norm and with $x_0^*y_0>0$.\;
\nl Let $r_{-1} = -\infty$\;
\nl Set $r_0 = \Re(\lambda_0)$\;
\nl Set $k=0$\;
\While{$r_k-r_{k-1} > {\rm tol}$}{
\nl Compute $x_{k+1},y_{k+1}$, left and right eigenvectors to the rightmost eigenvalue 
$\lambda_{k+1}$ of $A+\eps E_{k}$, $E_{k} = x_{k} y_{k}^*$, with $x_k, y_k$ of unit 
norm such that $x_{k}^*y_{k}>0$\;
\nl Set $k=k+1$\;
\nl Set $r_k=\Re(\lambda_k)$\;
}
\nl Set $r=r_k$\;
\nl Set $x=x_k, y=y_k$\;
}
\caption{Rank-1 iteration}
\label{alg_GO} 
\end{algorithm}

\begin{figure}\label{fig:goit}
\vskip -10mm
\centerline{
\includegraphics[scale=0.36]{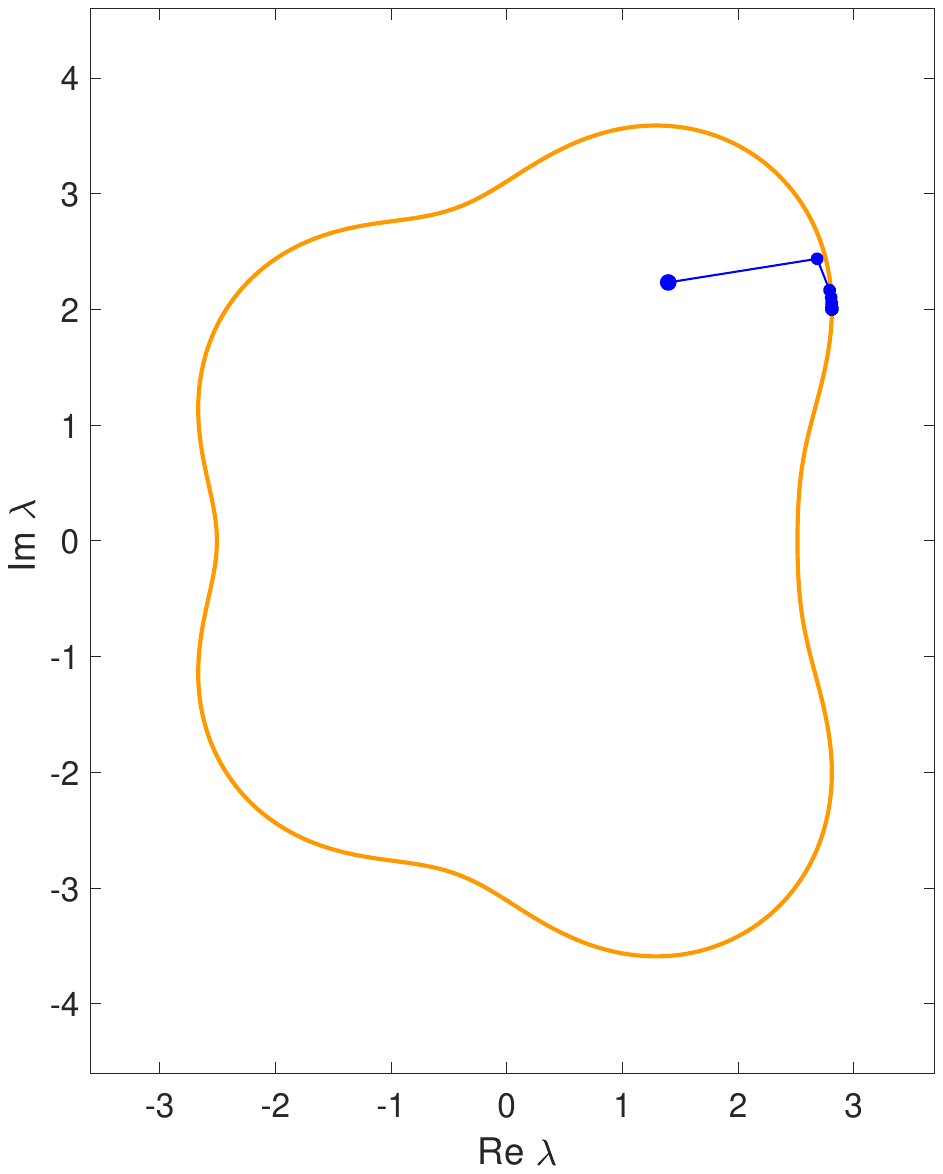}
}
\vskip -1cm
\caption{Iterates $\lambda_k$ of the rank-1 iteration (Algorithm~\ref{alg_GO}) applied to the matrix $A$ of (\ref{chap:proto}.\ref{eq:example})
and $\eps=1$.} 
\end{figure}

An illustration is given in Fig.~\ref{fig:goit} for the $8 \times 8$ matrix
(\ref{chap:proto}.\ref{eq:example}).

The hope is that $\Re\, \lambda_k$ converges to the $\eps$-pseudospecral abscissa $\alpha_\eps(A)$, as is frequently observed in numerical experiments.
Indeed, {\it if\/} the iteration converges, $\lambda_k\to\lambda$ and $x_k\to x$, $y_k \to y$, then $x,y$ are of unit norm with $x^*y>0$ and
\begin{equation}\label{stat-lim}
\text{$x,y$ are left and right eigenvectors to the rightmost eigenvalue $\lambda$ of $A+\eps x y^*$.}
\end{equation}
This implies $(A-\lambda I)y = - \eps x$ and $x^*(A-\lambda I) = -\eps y^*$, which shows that $A-\lambda I$ has the singular value $\eps$
(as is required for having $\lambda\in\partial\Lameps(A)$ by \eqref{ps-sv-bdy}) --- though $\eps$ is here not known to be the smallest singular value. Furthermore, the gradient of the associated singular value is $x^*y>0$, 
that is, the gradient  is horizontal to the right in the complex plane. By Theorem~\ref{thm:Delta-C}, this implies that $\lambda\in\partial\Lameps(A)$ with outer normal 1 if $\eps$ is indeed the {\it smallest} singular value of $A-\lambda I$.

Moreover, in the interpretation of Theorem~\ref{chap:proto}.\ref{thm:stat} and Section~\ref{subsec:ineq}, the property \eqref{stat-lim} implies that $E=xy^*$ is a {\it stationary point} (though not necessarily a maximum) of the
eigenvalue optimization problem
\begin{equation}\label{eigopt-psa}
\arg\max_{\| E \|_F \le 1} \Re\,\lambda(A+\eps E),
\end{equation}
i.e. Problem (\ref{chap:proto}.\ref{eq:optimiz0}) with $f(\lambda,\clambda)=-\Re\,\lambda$.

There exist no results about global convergence of the rank-1 iteration.  Local linear convergence can be shown for a sufficiently small 
ratio of the two smallest
singular values, $\eps/\sigma_{n-1}(A-\lambda I)$,
by studying the derivative of the iteration map at a stationary point. This requires bounds of derivatives of eigenvectors using appropriate representations of the group inverse, as laid out in the Appendix \bng 
(see Theorem \ref{thm:eigvecderiv}) 
\eng


\subsubsection*{Monotone rank-1 iteration.} The simple rank-1 iteration described above is not guaranteed to yield a monotonically increasing sequence $(\Re\,\lambda_k)$. Guglielmi and Overton (\cite{GO11}) also proposed a monotone variant that is described in the following \bng with  slight modifications  to the presentation in the article. \eng 

For given vectors $u,v$ of unit norm, we start from the rank-1 perturbation $A+\eps uv^*$ with rightmost eigenvalue $\lambda_0$, assumed to be simple. Let $x,y$ be left and right eigenvectors associated with $\lambda_0$, of unit norm and with $x^*y>0$. We still have a further degree of freedom in scaling $x$ and $y$, i.e. choosing the argument of the complex numbers $\alpha=u^*x$ and $\beta=v^*y$ of fixed modulus.

--- If $|\alpha|\ge |\beta|$, then we scale $x$ such that $\alpha$ is real and positive. Since we require $x^*y>0$, this also determines $y$ uniquely.

--- If $|\alpha| < |\beta|$, then we scale $y$ such that $\beta$ is real and positive. Since we require $x^*y>0$, this also determines $x$ uniquely.

With this particular scaling, we consider, for $0\le t \le 1$, a family of matrices 
\begin{equation} \label{eq:B}
B(t) = A + \eps p(t) q(t)^*, \qquad 0\le t \le 1,
\end{equation}
that interpolates between
$A+\eps uv^*$ at $t=0$ and $A+\eps xy^*$ at $t=1$:
\begin{equation} \label{eq:xyt}
p(t) = \frac{t x + (1-t) u}{\| t x + (1-t) u \|},
\quad q(t) = \frac{t y + (1-t) v}{\| t y + (1-t) v\|}.
\end{equation}
The following lemma will allow us to formulate a rank-1 iteration with monotonically increasing $\Re\,\lambda_k$.

\begin{lemma} [Monotonicity near $t=0$] \label{lem:loc-mon}
Let $B(t)$, $0\le t \le 1$, be defined as above with the stated scaling of the eigenvectors.
Let $\lambda(t)$, $0\le t \le 1$, be the continuous path of eigenvalues of
$B(t)$ with $\lambda(0)=\lambda_0$. If $\lambda_0$ is a simple eigenvalue of $A+\eps uv^*$, then $\lambda(t)$ is  differentiable at $0$ and
$$
\Re\,\dot\lambda(0)\ge 0.
$$
The inequality is strict except in the following two cases:
\begin{enumerate}
\item $\alpha=\beta=1$;
\item $\alpha$ and $\beta$ are both real, of equal modulus and opposite sign.
\end{enumerate}
\end{lemma}

\begin{proof}
By Theorem~\ref{chap:appendix}.\ref{thm:eigderiv} we have, with $\kappa=1/(x^*y)>0$, 
$$
\dot\lambda(0) = \frac{x^*\dot B(t)y}{x^*y}= \eps \kappa \Bigl(x^* \bigl(\dot p(0)q(0)^*+p(0)\dot q(0)^*\bigr) y \Bigr).
$$
We find $p(0)=u$, $q(0)=v$ and 
$$
\dot p(0)=(x-u)-u\Re(u^*(x-u)) = x - u \Re(u^*x)=x-u\Re\,\alpha, \quad \dot q(0)=y-v\Re\,\beta
$$
This yields
\begin{align*}
\Re\,\dot\lambda(0) &= \eps\kappa \Re\Bigl( x^*(x-u\Re\,\alpha)v^*y + x^* u (y-v \Re\,\beta)^*y \Bigr)
\\
&=\eps\kappa \bigl( \Re\,\beta -\Re\,\alpha\,\Re(\conj\alpha\beta) + \Re\,\alpha - \Re,\beta\,\Re(\conj\alpha\beta) \bigr),
\end{align*}
that is,
\begin{equation}\label{dot-lambda-ab}
\Re\,\dot\lambda(0) = \eps\kappa \left(1-\Re(\conj\alpha\beta)\right) \left( \Re\,\alpha+\Re\,\beta \right).
\end{equation}
With our scaling, the right-hand side is positive except in Case 1. or 2., where it vanishes.
%
\qed
\end{proof}

Lemma~\ref{lem:loc-mon} guarantees that $\Re~\lambda(t)>\Re~\lambda(0)$ for sufficiently 
small $t$. 
Hence the idea is to perform zero or more times the computation of
the rightmost eigenvalue $\lambda(t)$ of \eqref{eq:B}-\eqref{eq:xyt}
until 
$$
\Re~\lambda(t) > \Re~\lambda(0)+ \frac t2   \,\Re~\dot \lambda(0)
$$ 
with $\Re~\dot \lambda(0)$ given by \eqref{dot-lambda-ab}, replacing $t$ by $t/2$
until the inequality is fulfilled. 

In the $k$th iteration step, starting from $(u,v)=(u_k,v_k)$ and $\lambda(0)=\lambda_k$, we determine in this way $t>0$ such that
$\Re~\lambda(t)$ satisfies the above condition and then set $\lambda_{k+1}=\lambda(t)$ and $(u_{k+1},v_{k+1})=(p(t),q(t))$. 
The sequence of the real parts of the eigenvalues $\Re \, \lambda_k$ is then monotonically increasing.
In more detail, the variant is
formulated in Algorithm \ref{alg_GOM}.
\pagebreak[3]
%
\begin{algorithm}
\DontPrintSemicolon
\KwData{Matrix $A$, $\eps>0$, ${\rm tol}$ a given positive tolerance} 
\KwResult{$r \le \aleps(A)$, $x,y$}
\Begin{
\nl Compute $x_0$ and $y_0$ left and right eigenvectors to the rightmost eigenvalue of $A$, 
both normalized to unit norm and with $x_0^*y_0>0$.\;
\nl Let $r_{-1} = -\infty$\;
\nl Set $r_0 = \Re(\lambda_0)$\;
\nl Set $k=0$\;
\While{$r_k-r_{k-1} > {\rm tol}$}{
\nl Compute $x_{k+1},y_{k+1}$, left and right eigenvectors to the rightmost eigenvalue 
$\lambda_{k+1}$ of $A+\eps E_{k}$, $E_{k} = x_{k} y_{k}^*$, with $x_k, y_k$ of unit 
norm such that $x_{k}^*y_{k}>0$\;
\nl Set $\lambda = \lambda_{k+1}$\;
\nl Set $t=1$\;
\Repeat{$\Re \lambda > \Re \lambda_k$}{
\nl Set $t=t/2$\;
\nl Compute $x(t), y(t)$ according to \eqref{eq:xyt}\;
\nl Compute $\lambda(t)$\;
}
\nl Set $x_{k+1} = x(t)$, $y_{k+1}=y(t)$\;
\nl Set $k=k+1$\;
\nl Set $r_k=\Re(\lambda_k)$\;
}
\nl Set $r=r_k$, $x=x_k, y=y_k$\;
}
\caption{Rank-1 iteration: monotone version}
\label{alg_GOM} 
\end{algorithm}
%

%

 

\begin{theorem}[Convergence of the monotone rank-1 iteration] \label{th:convmon}
\bng
If Case 2 in Lemma~\ref{lem:loc-mon} does not occur, \eng then the monotone rank-1 iteration  converges to a stationary point of the eigenvalue optimization problem \eqref{eigopt-psa}, i.e., the limits 
$$
\lambda = \lim_{k\to\infty} \lambda_k,\quad u=\lim_{k\to\infty} u_k, \quad v=\lim_{k\to\infty} v_k
$$
exist,  and the stationarity condition \eqref{stat-lim} is satisfied. In particular, $\eps$ is a singular value of $A-\lambda I$ with left and right singular vectors $u$ and $v$. If $\eps$ is the smallest singular value, then $\lambda\in\partial\Lameps(A)$ with horizontal rightward outer normal.
\end{theorem}

We note that near a stationary point \eqref{stat-lim}, where $\alpha=\beta=1$, Case 2 in Lemma~\ref{lem:loc-mon} cannot occur.

\begin{proof}
Since the sequence $(\Re\,\lambda_k)$ is monotonically increasing and bounded, it converges. This implies that in the limit,
both sides of \eqref{dot-lambda-ab} are zero, and hence one of the two cases 1 or 2 in Lemma~\ref{lem:loc-mon} must occur in the limit. By assumption, we have excluded the exceptional Case 2. In the remaining Case 1, $\alpha=\beta=1$ in the limit, i.e.~$u=x$ and $v=y$ in the limit, and hence the iteration converges and the stationarity condition \eqref{stat-lim} is fulfilled. As noted before, this implies the further statements.
\qed
\end{proof}

%

\begin{figure}[h]
\begin{center}
\vspace{-10mm}
\includegraphics[scale=0.33]{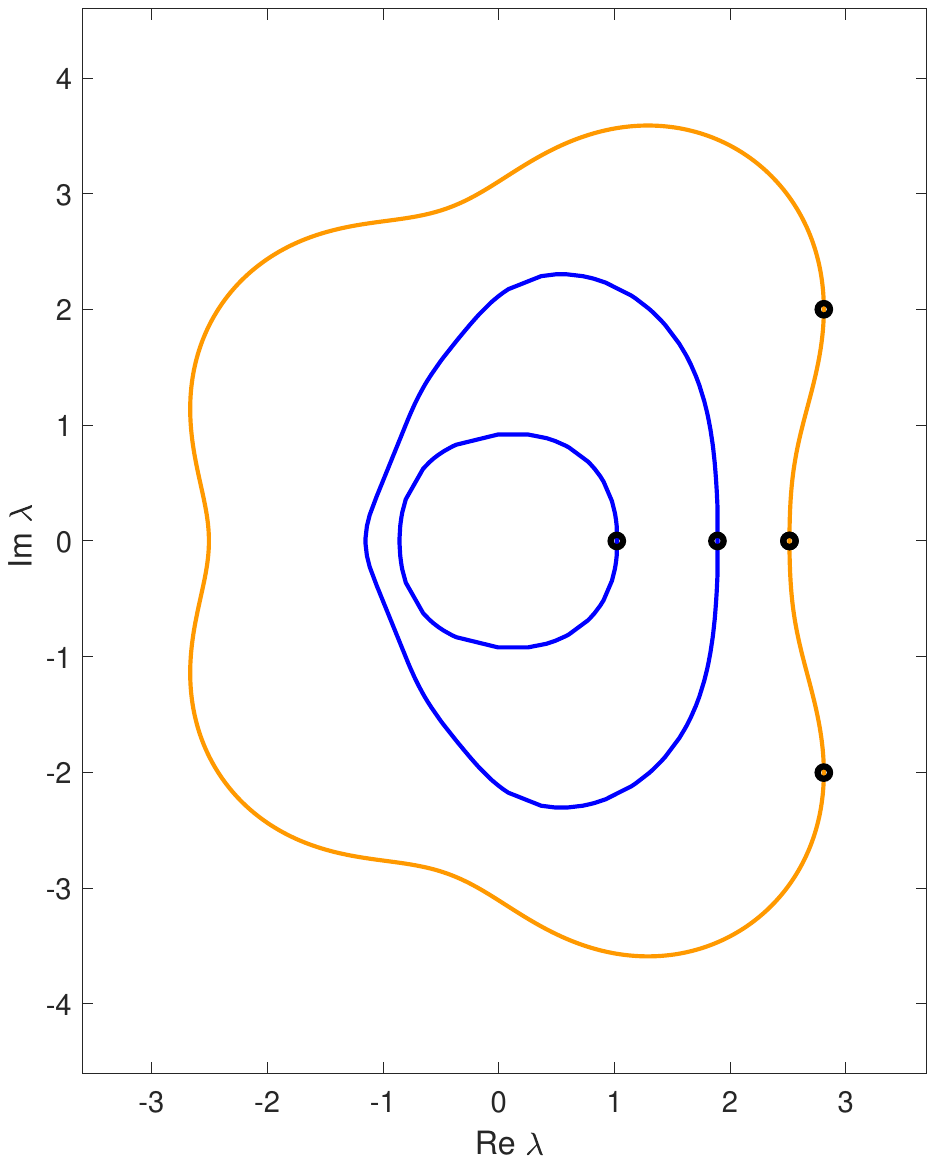}
\vspace{-10mm}
\caption{The outer curve is the boundary of the pseudospectrum
$\Lameps(A)$ and the inner curves are the $\eps$-level sets of
$\sigma_{n-1}(A-zI)$ and $\sigma_{n-2}(A-zI)$ for the matrix $A$ of (\ref{chap:proto}.\ref{eq:example}) and $\eps=1$. The black points are fixed points of the algorithm by Guglielmi \& Overton, as well as stationary points of the differential equations (\ref{chap:proto}.\ref{ode-uv}) and (\ref{chap:proto}.\ref{ode-uv-short}). \label{fig:extsv}}
\end{center}
\end{figure}

Figure \ref{fig:extsv} shows the $\eps$-pseudospectrum of the matrix $A$
in \eqref{eq:example} for $\eps=1$ as well as the $\eps$-level set of 
$\sigma_{n-1}(A-zI)$ and $\sigma_{n-2}(A-zI)$ (the inner curves).

The points to which the algorithms may converge are given by the $5$ points indicated by bullets.
The locally rightmost ones are locally attractive while the other ones turn out to 
be locally unstable.

\subsection{Discretized rank-1 matrix differential equation}
\index{rank-1 matrix differential equation}
\index{gradient flow!rank-1 constrained}

A different iteration on rank-1 matrices results from the rank-1 projected gradient system of Section~\ref{subsec:rank1-gradient-flow}  for the minimization function $f(\lambda,\clambda)=-\tfrac12(\lambda+\clambda)=-\Re\,\lambda$ after discretization as in Section~\ref{subsec:proto-numer} (see Algorithm~\ref{chap:proto}.\ref{alg_prEul}), as was first proposed similarly by Guglielmi \& Lubich~(\cite{GL11}), though with a different time stepping method. The so obtained rank-1 iteration yields a sequence of rank-1 matrices $E_k=u_k v_k^*$  of unit norm and a sequence of eigenvalues $\lambda_k\in\Lameps(A)$ of $A+\eps E_k$ with monotonically growing real part,  which converges to a stationary point~\eqref{stat-lim}; see Theorems~\ref{chap:proto}.\ref{thm:stat-1} and Remark~\ref{chap:proto}.\ref{rem:exceptional}, and also Lemma~\ref{chap:proto}.\ref{lem:stat-split} for the splitting method. The computational cost per step is essentially the same as in the rank-1 iterations of the previous subsection. A numerical example was already presented in Section~\ref{subsec:proto-numer}.

Conceptually, the approach of first deriving a suitable differential equation and then using an adaptive time-stepping to arrive at a stationary point is different from directly devising an iteration. Different tools are available and used in the two approaches. For example, the tangent space of the manifold of rank-1 matrices is a natural concept in the time-continuous setting though not so in the time-discrete setting. This enhanced toolbox results in efficient algorithms that would not be obtained from a purely discrete viewpoint.

\subsection{Acceleration by a subspace method} 
Kressner and Vandereycken (\cite{KV14}) proposed a subspace method to accelerate the basic {rank-1} iteration described in Subsection~\ref{subsec:r1-iteration}. For the subspace expansion, they essentially do a step of Algorithm~\ref{alg_GO}  and add the obtained eigenvector to the subspace. Along the way they compute orthonormal bases $V_k \in \C^{n\times k}$ of nested subspaces. A key element is the computation of the rightmost point of the $\eps$-pseudospectrum of the rectangular matrix pencil $AV_k -\lambda V_k$ in place of $A-\lambda I$, 
$$
\Lambda_\eps(AV_k,V_k)=\{ \lambda \in \C\,:\, \sigma_{\min}(AV_k-\lambda V_k) \le \eps \}.
$$
These pseudospectra are nested: $\Lambda_\eps(AV_k,V_k) \subset \Lambda_\eps(AV_{k+1},V_{k+1}) \subset \Lambda_\eps(A)$.
The rightmost point of $\Lambda_\eps(AV_k,V_k)$ is computed by a variant of the criss-cross algorithm. The basic algorithm is given in Algorithm~\ref{alg:KV-subspace}.

\begin{algorithm}
\DontPrintSemicolon
\KwData{Matrix $A$, $\eps>0$} 
\KwResult{Approximation $\mu_\alpha$ to a locally rightmost point of $\Lambda_\eps(A)$.}
\Begin{
\nl Compute the rightmost eigenvalue $\lambda_0$ and normalized right eigenvector $y_0$ of $A$.\;
\nl Set $\widehat V_1 = y_0$.\;
\nl \For {$k=1,2,\dots$ {\rm until converged}} {
\nl Compute the rightmost point $\mu_k$ of $\Lambda_\eps(AV_k,V_k)$.\;
\nl Compute left/right singular vectors $u_k$ and $v_k$  to $\sigma_{\min}(A-\mu_k I)$.  Set $E_k = -u_k v_k^*$.\;
\nl Compute the rightmost eigenvalue $\lambda_k$ and right eigenvector $y_k$ of $A + \eps E_k$. Compute $V_{k+1} =\, $orth$([V_k, y_k])$.\;
}
\nl Set $\mu_\alpha = \mu_k$.
}
\caption{Subspace method.}
\label{alg:KV-subspace} 
\end{algorithm}

\noindent
Kressner and Vandereycken (\cite{KV14}) show that the sequence $(\mu_k)$ grows monotonically, as a consequence of the growth of the nested subspaces.
A simplified version of the algorithm,  where the right singular vector $v_k$ instead of the right eigenvector $y_k$ is added to the subspace, is shown to converge locally superlinearly to the pseudospectral abscissa.

\section{Tracing the boundary of the pseudospectrum}
\label{sec:ps-tracing}
\index{pseudospectrum boundary}

In this section we describe two algorithms for boundary tracing. While there exist path-following methods to
obtain pseudospectral contours (e.g. those implemented in Eigtool), the computation becomes expensive for large matrices. 
Here we use instead the low-rank structure of the extremal perturbations, which allows us to treat also large sparse matrices efficiently; see Section~\ref{subsec:rank-1}.

We present two algorithms. 
The first algorithm, to which we refer as the {\it tangential--transversal algorithm}, makes use of a combination of the differential equation (\ref{chap:proto}.\ref{ode-uv}) and a similar differential equation that moves eigenvalues horizontally to the boundary. The second algorithm, which we call the {\it ladder algorithm},  aims to compute, for an iteratively constructed sequence of points outside the $\eps$-pseudospectrum, the corresponding nearest points in the  $\eps$-pseudospectrum.
Both algorithms require repeatedly the computation of the eigenvalue of  a rank-1 perturbation
to~$A$ that is nearest to a given complex number. This can be done efficiently also for large matrices, using inverse power iteration combined with the Sherman-Morrison formula.

The ladder algorithm extends readily to real and structured pseudospectra, 
for which the few algorithms proposed in the literature 
\bng (such as implemented in SEigTool, see 
Karow, Kokiopoulou \& Kressner \cite{KarKK10}) \eng are restricted to just a few structures and turn out to be extremely demanding from a computational point of
view.

\subsection{Tangential--transversal algorithm}
\index{tangential--transversal algorithm}

The algorithm alternates between a time step for the system of differential equations (\ref{chap:proto}.\ref{ode-uv}) and the following system of differential equations for vectors $u(t)$ and $v(t)$ of unit norm. This second system is a simplified variant of 
(\ref{chap:proto}.\ref{ode-uv}) for $G=-xy^*$, where $x$ and $y$ are left and right eigenvectors, respectively, 
both of unit norm and with $x^*y>0$,
of an eigenvalue $\lambda$ of the rank-1 perturbed matrix $A+\eps uv^*$:
\begin{equation}\label{ode-hor}
\begin{array}{rcl}
 \dot u &=&  (I-uu^*)xy^*v
\\[1mm]
 \dot v &=&  (I-vv^*)yx^*u.
\end{array}
\end{equation}
The system preserves the unit norm of $u$ and $v$, since $u^*\dot u=0$ and $v^*\dot v=0$. As we show in the next lemma, the system
has the property that for a path of simple eigenvalues $\lambda(t)$ of $A+\eps u(t)v(t)^*$, the derivative $\dot\lambda(t)$ is real and positive, continuing to a stationary point where $A-\lambda I$ has the singular value $\eps$. By Theorem~\ref{thm:ps-sv} it therefore stops at the boundary $\partial\Lambda_\eps(A)$ if $\eps$ is the smallest singular value. While in theory, a trajectory might stop at an interior point where the singular value $\eps$ is not the smallest one, this appears to be an unstable case that is not observed in computations.

\begin{lemma}[Horizontal motion of an eigenvalue]
\label{lem:hor-motion}
Along a path of simple eigenvalues $\lambda(t)$ of $A+\eps u(t)v(t)^*$, where $u,v$ of unit norm solve \eqref{ode-hor}, we have that
$$
\text{$\dot\lambda(t)$ is real and positive for all\, $t$.}
$$
In the limit $\lambda_\star=\lim_{t\to\infty} \lambda(t)$, the matrix $A-\lambda_\star I$ has the singular value $\eps$.
\end{lemma}

\begin{proof}
The perturbation theory of eigenvalues (see Theorem~\ref{chap:appendix}.\ref{thm:eigderiv}) shows that
\begin{equation*}
\dot\lambda = \frac{x^*\frac d{dt}(A+\eps  u v^*)\,y}{x^*y}= \eps \,\frac{x^* \left( {\dot u} v^* + u {\dot v}^*\right) y}{x^*y}. 
\end{equation*}
With $\alpha=u^* x$ and $\beta=v^*y$, we obtain
\begin{equation*}
\dot \lambda = \frac{\eps}{x^*y}\Bigl(|\alpha|^2\cdot \|y-\beta v\|^2 + 
 |\beta|^2\cdot \|x-\alpha u \|^2\Bigr) \in\R, \ \ge 0.
\end{equation*} 
In a stationary point of \eqref{ode-hor}, $u$ and $x$ are collinear, and so are $v$ and $y$. It follows that
$uv^*=\e^{\iu\theta}xy^*$ for some real $\theta$. We thus have
$$
(A+\eps \e^{\iu\theta}xy^*)y=\lambda y, \qquad x^*(A+\eps \e^{\iu\theta}xy^*)=\lambda x^*,
$$
or equivalently
$$
(A-\lambda I)y = \eps \e^{\iu\theta}x, \qquad (A-\lambda I)^* \e^{\iu\theta}x = \eps y,
$$
which states that $\e^{\iu\theta}x$ and $y$ are left and right singular vectors to the singular value $\eps$ of~$A-\lambda I$.
\qed
\end{proof}
While \eqref{ode-hor} moves eigenvalues horizontally to the right, the differential equation (\ref{chap:proto}.\ref{ode-uv}) moves an eigenvalue on the boundary along a path that starts tangentially to the boundary, as is shown by the following lemma.

\begin{lemma}[Tangential motion of an eigenvalue from the boundary]
\label{lem:tang-motion}
Let $\lambda_0\in\partial \Lambda_\eps(A)$ be on a smooth section of the boundary, with outer normal $\e^{\iu\theta}$ at $\lambda_0$ for $0<|\theta|\le\pi/2$.
Let  $u_0$ and $v_0$ of unit norm be such that $A+\eps u_0v_0^*$ has $\lambda_0$ as a simple eigenvalue. Let $u(t)$ and $v(t)$ be solutions of the system of differential equations (\ref{chap:proto}.\ref{ode-uv}) with $G=-xy^*$ with initial values $u_0$ and $v_0$. Then, the path of eigenvalues $\lambda(t)$ of $A+\eps u(t)v(t)^*$ with $\lambda(0)=\lambda_0$ has $\dot\lambda(0)\ne 0$ and
$$
\text{$\dot\lambda(0)$ is tangential to $\partial\Lambda_\eps(A)$ at $\lambda_0$.}
$$
\end{lemma}

\begin{proof} By Theorem~\ref{thm:Delta-C}, we have
$
u_0v_0^* = \e^{\iu\theta} x_0y_0^*.
$
We then find, inserting (\ref{chap:proto}.\ref{ode-uv}) for $\dot u$ and~$\dot v$,
\begin{align*}
\dot \lambda(0) &=  \frac{x_0^*\,\frac d{dt}\big|_{t=0}(A+\eps  u(t) v(t)^*)\,y_0}{x_0^*y_0}
\\
&= \frac\eps{x_0^*y_0} \, x_0^*\bigl( \dot u(0) v(0)^*+u(0)\dot v(0)^*\bigr) y_0
\\
&=  \frac\eps{x_0^*y_0} \, \iu\, \Im(u_0^*x_0y_0^*v_0) x_0^* u_0 v_0^* y_0
\\
&=  \frac\eps{x_0^*y_0} \, \iu \,\Im(\e^{-\iu\theta}) \,\e^{\iu\theta}.
\end{align*}
Since $x_0^*y_0>0$, we find that $\dot \lambda(0)$ points into the tangential direction $-\iu\,\e^{\iu\theta}\,\mathrm{sign}(\theta)$.
\qed
\end{proof}


\subsubsection*{Description of the algorithm.}
Lemmas~\ref{lem:hor-motion} and~\ref{lem:tang-motion} motivate the following algorithm for tracing the boundary of the $\eps$-pseudospectrum $\Lambda_\eps(A)$.
Suppose that $\lambda_0$ is a simple eigenvalue of $A+\eps u_0 v_0^*$ lying on a smooth section of the boundary 
$\partial\Lambda_\eps(A)$, with {\it a priori} unknown outer normal $\e^{\iu\theta_0}$. By Theorem~\ref{thm:Delta-C}, $u_0 v_0^*= \e^{\iu\theta_0}x_0y_0^*$, where $x_0$ and $y_0$ are left and right eigenvectors of $A+\eps u_0v_0^*$, both of unit norm and with $x_0^*y_0>0$. This fact allows us to determine the outer normal as 
\begin{equation} \label{normal-theta}
\e^{\iu\theta_0}=\frac{|u_0^*v_0|} {u_0^*v_0}.
\end{equation}
We first consider the tangential differential equation (\ref{chap:proto}.\ref{ode-uv}) with $G=-xy^*$ 
for the rotated matrix $\iu\e^{-\iu\theta_0}A$,
which leads us to the case $\theta=\pi/2$ in Lemma~\ref{lem:tang-motion}. 
We make a time step of stepsize $h$ with the splitting method presented in Chapter \ref{chap:proto}, 
followed by normalization of $u$ and $v$; 
at $t_1=t_0+h$ this yields a rank-1 matrix $\widetilde E_1=\tilde u_1 \tilde v_1^*$ of 
unit norm. (Note that this step does not require to actually compute the rotated matrix.) 
Since $\dot\lambda(t_0)$ is tangential to the boundary $\partial\Lambda_\eps(A)$ at $\lambda_0$, the 
eigenvalue $\tilde \lambda_1$ of $A+\eps \widetilde E_1$ lies in $\Lambda_\eps(A)$ and is $O(h^2)$
close to the boundary.


With initial values $ \tilde u_1,\tilde v_1$ we then consider the differential equation
\eqref{ode-hor} for the rotated matrix $\e^{-\iu\theta_0}A$
in order to reach the boundary with a horizontal trajectory.
We integrate \eqref{ode-hor} until we stop at a stationary point $E_1=u_1 v_1^*$. 
There, a singular value of $A-\lambda I$ equals $\eps$ (by Lemma~\ref{lem:hor-motion}), which takes
us to the boundary $\partial\Lambda_\eps(A)$ 
provided we start sufficiently close to it.

We then continue the above alternating integration from $\lambda_1\in \partial\Lambda_\eps(A)$ and the 
associated vectors $u_1,v_1$.

The algorithm as described computes a part of $\partial\Lambda_\eps(A)$ to the right of $\lambda_0$. To go to the left,
we change the direction of time in the tangential differential equation~(\ref{chap:proto}.\ref{ode-uv}), i.e., we take a negative stepsize $h$. The algorithm is summarized in Algorithm~\ref{alg_tang-trans}.


\begin{algorithm}
\DontPrintSemicolon
\KwData{Matrix $A$, initial vectors $u,v$ of unit norm such that the eigenvalue $\lambda$ of $A+\eps u v^*$ lies on $\partial\Lambda_\eps(A)$,
stepsize $h$, number $N$ of desired boundary points,
${\rm tol}$ a given positive tolerance} 
\KwResult{vector  $\Gamma$ of $N$ consecutive boundary points}
\Begin{
\For{i=1,\dots,N}{
\nl Set $\zeta= u^*v / |u^*v|$\;
\nl Compute the approximate solution $\widetilde u_1,\widetilde v_1$ of (\ref{chap:proto}.\ref{ode-uv}) for the rotated matrix $\iu\zeta A$
with initial data $\iu\zeta u,v$ doing a single normalized Euler step of size $h$\;
\nl Let $u_0 = -\iu \widetilde u_1$, $v_0 = \widetilde v_1$\; 
\For{$k=1,2,\ldots$ until convergence}{ 
\nl Compute the approximate solution $u_k,v_k$ of (\ref{ode-hor}) for the rotated matrix $\zeta A$ with adaptive stepsize $h_k$
(as in Chapter~\ref{chap:proto})\; 
\nl Compute the rightmost eigenvalue $\lambda_k$ of $\zeta A + \eps u_k\,v_k^*$\; 
}
\nl Set $u=u_k/\zeta$, $v=v_k$\;
\nl Store $\lambda_k/\zeta$ into $\Gamma$\;
}
}
\caption{Tangential--transversal algorithm for tracing the boundary of the $\eps$-pseudospectrum}
\label{alg_tang-trans} 
\end{algorithm}

\begin{figure}[h!]
\centerline{
\includegraphics[scale=0.3]{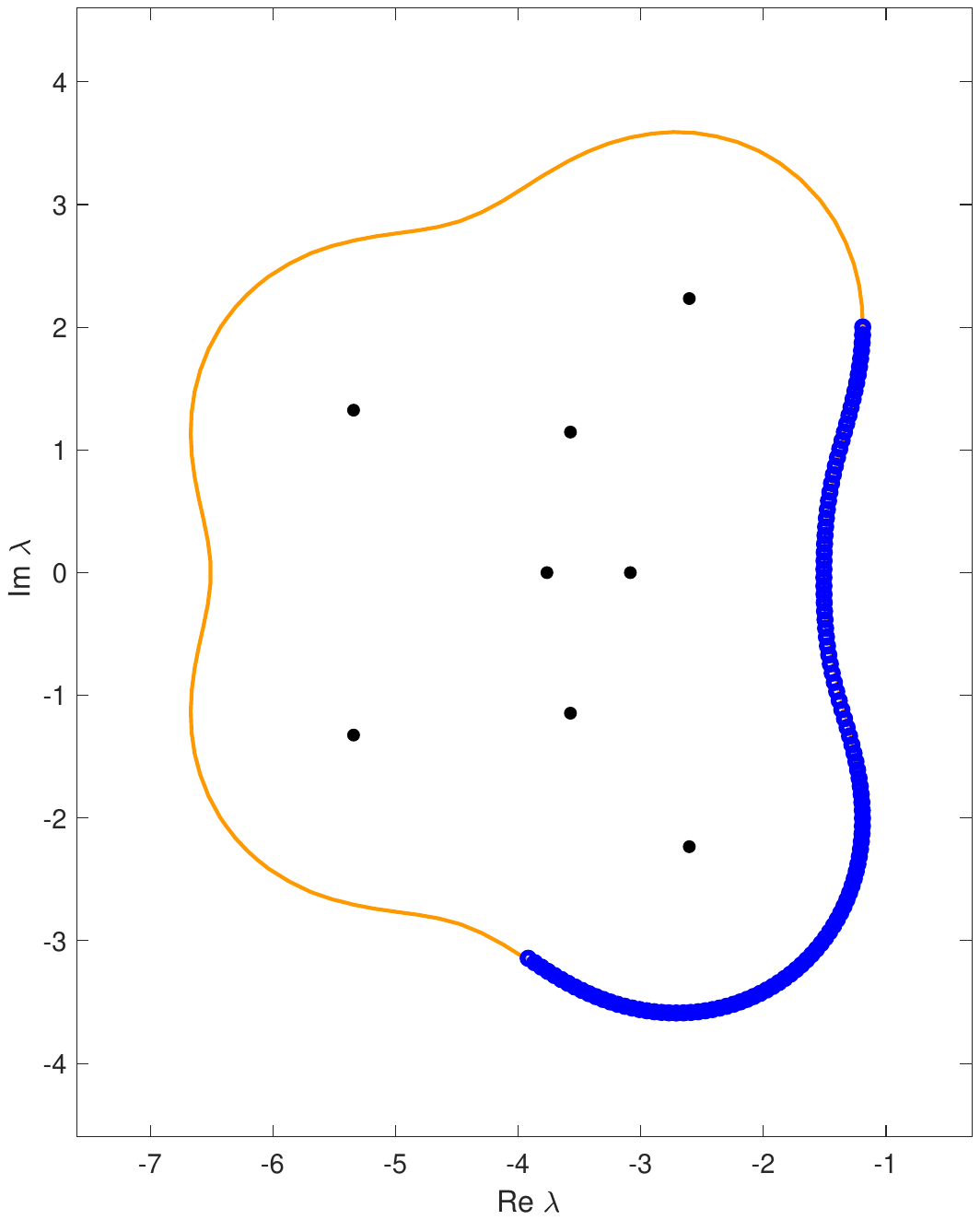}
}
\vspace{-0.3cm}
\caption{The $\eps$-pseudospectrum and a section
of its boundary (small circles) determined by the tangential--transversal Algorithm~\ref{alg_tang-trans} 
for the  matrix (\ref{chap:proto}.\ref{eq:example}) and $\eps=1$.
\bng The algorithm starts at the top rightmost point and continues until it is stopped at $\Re\,\lambda=-4$ to make the boundary curve of the $\eps$-pseudospectrum as computed by Eigtool (thin line) visible. \eng
\label{fig:tt}}
\end{figure}

\subsubsection*{Numerical experiment.}
We consider again the $8 \times 8$  matrix (\ref{chap:proto}.\ref{eq:example})
with $\eps=1$. Applying the tangential--transversal algorithm, we obtain the result in Figure \ref{fig:tt},
which is superimposed 
on the boundary of the $\eps$-pseudospectrum (computed to high accuracy using Eigtool by Wright, 
\cite{Wri02}).

\subsection{Ladder algorithm}
\label{subsec:ladder}
\index{ladder algorithm}


As in the previous algorithm, let $\lambda_0\in\partial\Lambda_\eps(A)$ be a simple eigenvalue of $A+\eps u_0v_0^*$ (with $u_0$ and $v_0$ of unit norm) that lies on a smooth section of $\partial\Lambda_\eps(A)$, with outer normal $\e^{\iu\theta_0}$ at $\lambda_0$
obtained from \eqref{normal-theta}.

With a small distance $\delta>0$, we define the nearby point $\mu_0$ on the straight line normal to $\partial\Lambda_\eps(A)$ at $\lambda_0$,
$$
\mu_0=\lambda_0 + \delta \e^{\iu\theta_0}.
$$
We add a tangential component, either to the left ($+$) or to the right ($-$),
$$
\mu_1= \mu_0 \pm \iu \delta \e^{\iu\theta_0}.
$$
We then apply the eigenvalue optimization algorithm of Chapter~\ref{chap:proto} for the function
$$
f(\lambda,\clambda)=(\lambda -\mu_1)(\clambda -\conj{\mu_1})=|\lambda -\mu_1|^2,
$$
choosing $E_0=u_0v_0^*$ as the starting iterate. That algorithm aims to compute $u_1$ and $v_1$ of unit norm such that $A+\eps u_1 v_1^*$ has the boundary point $\lambda_1\in\partial\Lameps(A)$ nearest to $\mu_1$ as an eigenvalue. At the point $\lambda_1$ we have the outer normal $(\mu_1 - \lambda_1)/| \mu_1 - \lambda_1 | $.

We continue from $\lambda_1$ and $u_1$, $v_1$ in the same way as above, constructing a sequence $\lambda_k$ ($k\ge 1$) of points on the boundary of the pseudospectrum $\Lameps(A)$ with approximate spacing $\delta$. 
This is done in Algorithm~\ref{alg_ladder}.

\begin{algorithm}
\DontPrintSemicolon
\KwData{Matrix $A$, initial vectors $u,v$ of unit norm such that the eigenvalue $\lambda$ of $A+\eps u v^*$ lies on $\partial\Lambda_\eps(A)$,
step size $\delta$, number $N$ of desired boundary points, tolerance
${\rm tol}>0$} 
\KwResult{vector  $\Gamma$ of $N$ consecutive boundary points}
\Begin{
\For{i=1,\dots,N}{
\nl Set $\zeta= u^*v / |u^*v|$\;
\nl Set $\mu=\lambda + (1-\iu)\delta/\zeta$\;
\nl Compute the nearest point to $\mu$ on $\partial\Lambda_\eps(A)$ by the rank-1 eigenvalue optimization algorithm of Chapter~\ref{chap:proto} for minimizing the function
$f(\lambda,\clambda)=|\lambda-\mu|^2$,
with $u$ and $v$ as the starting iterate and with the tolerance parameter ${\rm tol}$. This yields an update of $u$, $v$ and $\lambda$.\;
\nl Store $\lambda$ into $\Gamma$.\;
}
}
\caption{Ladder algorithm for tracing the boundary of the {$\eps$-pseudospectrum}.}
\label{alg_ladder} 
\end{algorithm}

\begin{figure}[h!]
\vspace{-8mm}
\centerline{
\hspace{-5mm}
\includegraphics[scale=0.33]{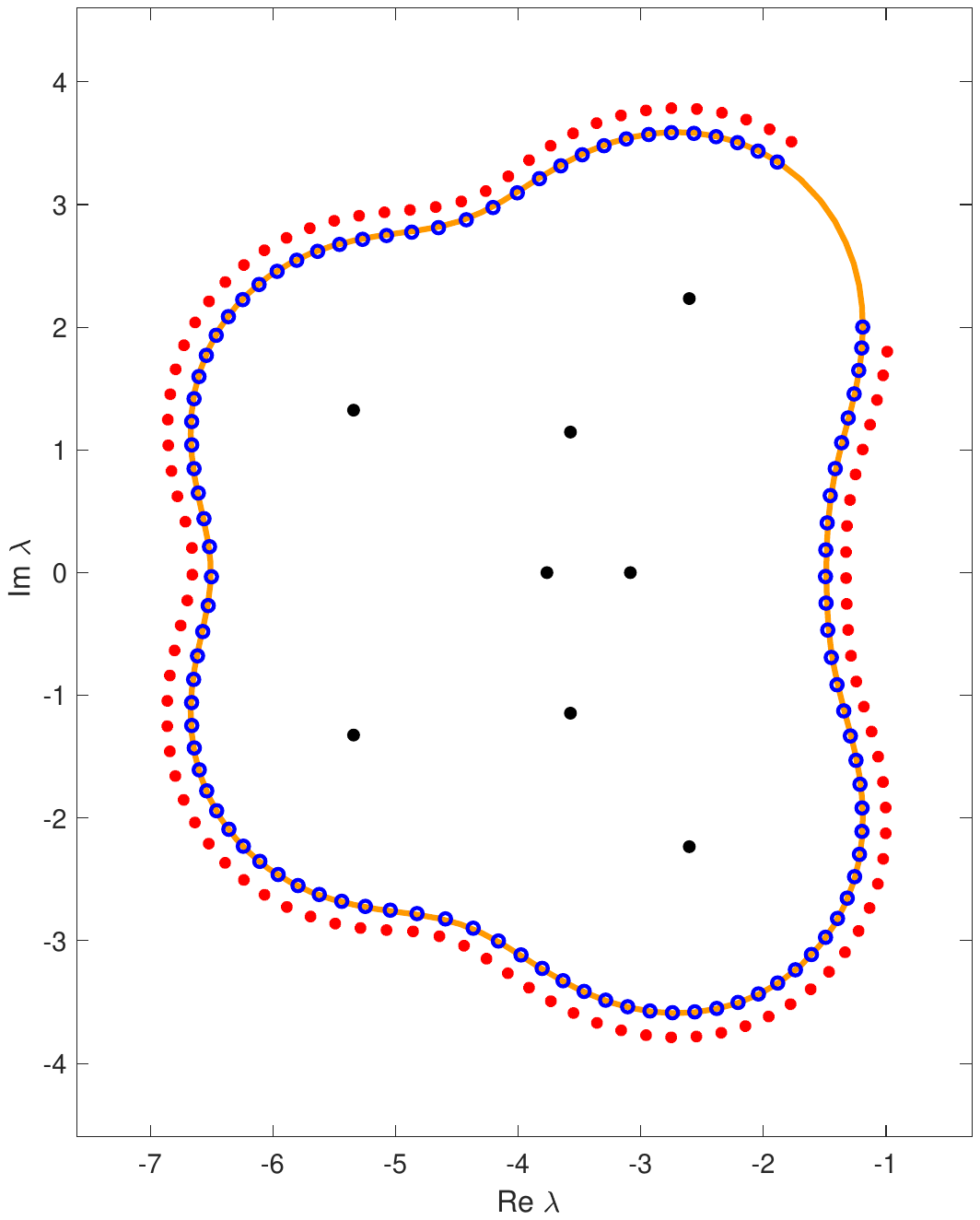}
\hspace{-5mm}
\includegraphics[scale=0.33]{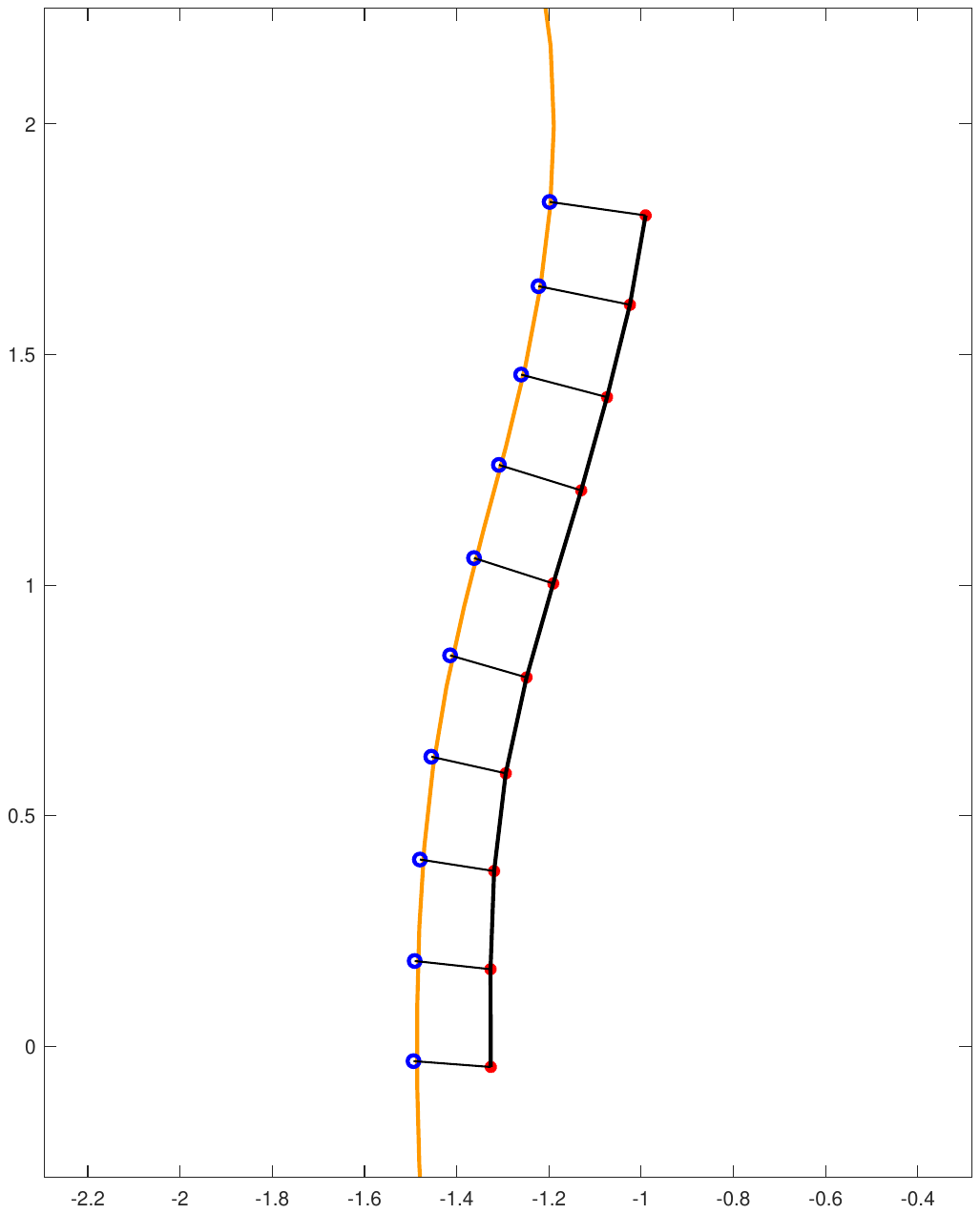}
}
\vspace{-0.3cm}
\caption{Boundary of the $\eps$-pseudospectrum (continuous line)  
with $\eps=1$ for the matrix $A$ in (\ref{chap:proto}.\ref{eq:example})
with boundary points (circles) and control points (bullets) computed by the ladder algorithm (Algorithm~\ref{alg_ladder}).
Right: zoom.
\label{fig:L1}}
\end{figure}

\begin{remark} The points $\lambda_k$ and $\mu_k$ and the straight lines between them form the ``rope ladder'' along the boundary of the pseudospectrum to which the name of the algorithm refers. We climb up or down on this ladder to construct the sequence of boundary points.
\end{remark}

\begin{remark} If the curvature of the boundary is larger than $1/\delta$, then it may happen that $\mu_k$ gets to lie inside $\Lameps(A)$,
and the algorithm will find $\lambda_k=\mu_k$. In such a situation, the step size $\delta$ needs to be reduced.
\end{remark}

\subsubsection*{Numerical experiment.}
We apply the ladder algorithm to the matrix $A$ in (\ref{chap:proto}.\ref{eq:example}).
Algorithm \ref{alg_ladder} produces the curve of circles in Figure \ref{fig:L1} which is superimposed 
on the boundary of the $\eps$-pseudospectrum (computed to high accuracy by Eigtool by Wright, 
\cite{Wri02}).
The average number of steps for each horizontal computation is $5.5$ for an accuracy tolerance
${\rm tol}=10^{-7}$.

\section{Notes}
\label{sec:ps-notes}

\subsubsection*{Pseudospectra.}
The standard reference for (complex unstructured) pseudospectra is the book by Trefethen \& Embree (\cite{TreE05}). 
The name ``pseudospectrum'' was coined by Trefethen (\cite{Tre92}), but as noted there, the concept was used before under different names. The ``approximate eigenvalues'' of Varah (\cite{Var67}) appear to be among the earliest precursors. 
Interpreting approximate eigenvalues or pseudospectral values as exact eigenvalues of a perturbed matrix is in the spirit of backward error analysis as pioneered by Wilkinson (\cite{Wil65}). In one of his last papers, Wilkinson (\cite{Wil86}) considered pseudospectra under the 
name ``fundamental inclusion domains''.
Pseudospectra are of interest not only for matrices but more generally for linear operators on Hilbert or Banach spaces; see Trefethen (\cite{Tre97}) and also the concise account in the book by Davies (\cite{Dav07}), Chapter~9.

The motivating example of robust stability of a matrix (see Subsection~\ref{subsec:ps-motivation}) was already studied by Van Loan (\cite{VL85}) and Hinrichsen \& Pritchard (\cite{HinP86a},\cite{HinP86b},\cite{HinP90}). They aim at determining the minimal norm of complex, real or structured perturbations that turn a stable matrix into an unstable one. This yields the (complex, real or structured) distance to instability, or stability radius, which equals the smallest perturbation size $\eps$ for which the $\eps$-pseudospectral abscissa becomes non-negative; this will be taken up in the next two chapters. 

The fundamental Theorem~\ref{thm:ps-sv} and its proof via the distance to singularity was given by
Wilkinson (\cite{Wil86}).
Theorem~\ref{thm:Delta-C} on extremal complex perturbations is related to Guglielmi \& Overton (\cite{GO11}). 

The use of pseudospectra together with the Cauchy integral formula to better understand the transient behaviour of  dynamical systems was suggested by Trefethen (\cite{Tre92}). Subsection~\ref{subsec:ps-exp} follows this approach, emphasizing the role of the pseudospectral abscissa and the stability radius. There are close connections to the Kreiss matrix theorem (Kreiss \cite{Kr62}, LeVeque \& Trefethen \cite{LeVT84} and Spijker \cite{Spi91}) and to further stability bounds based on resolvent bounds as given, e.g., by Lubich \& Nevanlinna (\cite{LubN91}), Reddy \& Trefethen (\cite{RedT92}) and
van Dorsselaer, Kraaijevanger \& Spijker (\cite{DorKS93}); see also Eisner (\cite{Eis10}) and references therein.

 Pseudospectra of matrix pencils $A-\lambda B$ were studied by van Dorsselaer (\cite{Dor97}) and, e.g.,
 Ahmad, Alam \& Byers (\cite{AhmAB10}), structured pseudospectra for polynomial eigenvalue problems by Tisseur \& Higham (\cite{TisH01}), and pseudospectra for rectangular matrices by Wright \& Trefethen (\cite{WriT02}).
 
An extension of the concept of the $\eps$-pseudospectrum of interest in control theory is the {\em spectral value set} that consists of the eigenvalues of all matrices $A+B\Delta (I-D\Delta)^{-1}C$ with given system matrices $A, B,C,D$  of compatible dimensions and with varying (complex or real or structured) matrices $\Delta$ of 2-norm at most~$\eps$. This was first considered (for $D=0$) by Hinrichsen and Kelb (\cite{HinK93});
 see also Karow (\cite{Kar03}) for a detailed study and the book of Hinrichsen \& Pritchard~(\cite{HinP05}), Chapter 5. 
Apart from eigenvalues of $A$, the unstructured complex spectral value set contains all $\lambda \in \C$ for which the transfer matrix $H(\lambda)= C(\lambda I-A)^{-1}B+D$ has 2-norm at least $\eps^{-1}$; see Guglielmi, G\"urb\"uzbalaban \& Overton (\cite{GugGO13}).
This characterization will be essential in Chapter~\ref{chap:lti}.

\bng
For the Weyl inequality, which we have used in Subsection~\ref{subsec:ps-exp}  to derive bounds for the matrix exponential, we refer the reader to Horn \& Johnson (\cite{HJ90}) and to Marshall, Olkin \& Arnold (\cite{MO11}). 

The Plancherel formula is standard in Fourier theory; see e.g. the classical textbook by Stein \& Shakarchi (\cite{SS03}).
\eng

\subsubsection*{Algorithms for computing the pseudospectral abscissa and radius.}
The criss-cross algorithm of Burke, Lewis \& Overton (\cite{BuLeOv03}) for computing the pseudospectral abscissa relies on Byers' lemma (Lemma~\ref{lem:byers}), due to Byers (\cite{Bye88}). This lemma is fundamental in that it relates singular values of a matrix shifted along a line parallel to the imaginary axis to the eigenvalues of a Hamiltonian matrix. This Hamiltonian connection has been put to important and enduring use in control systems starting with the work by Boyd, Balakrishnan \& Kabamba (\cite{BoyBK89}); 
see e.g.~Grivet-Talocia \& Gustavsen (\cite{GriG15}) and references therein. We will encounter Hamiltonian eigenvalue optimization problems in Chapters~\ref{chap:mnp-mix} and~\ref{chap:lti}. 
Mengi \& Overton (\cite{MeOv05}) extended the criss-cross algorithm to computing the pseudospectral radius (using radial and circular searches), and Lu \& Vandereycken (\cite{LuV17}) developed a criss-cross type algorithm for computing the real pseudospectral abscissa.
Benner \& Mitchell (\cite{BenM19})
\bng made improvements to the original criss-cross algorithms in the root-finding in the horizontal searches. They also generalized the criss-cross techniques to spectral value sets.
\eng

The rank-1 iteration of Guglielmi \& Overton (\cite{GO11}) for approximating the pseudospectral abscissa and radius appears to be the first algorithm that uses the low-rank property of extremal perturbations as described by Theorem~\ref{thm:Delta-C}. It motivated the rank-1 projected gradient flow algorithm of Guglielmi \& Lubich (\cite{GL11}) for approximating the pseudospectral abscissa and radius, which opened up an approach to a wide range of eigenvalue optimization and matrix nearness problems, as discussed throughout this book. 

Subspace acceleration approaches as used by Kressner \& Vandereycken (\cite{KV14}) also have a much wider scope than merely computing the pseudospectral abscissa; see, e.g., Kangal, Meerbergen, Mengi \& Michiels (\cite{KanMMM18}) and
Kressner, Lu \& Vandereycken (\cite{KreLV18}) for their use in other eigenvalue optimization problems.

In retrospect it appears remarkable how the modest aim of computing the pseudospectral abscissa led to the discovery of diverse classes of efficient algorithms that find use in a wide variety of other problems.

\subsubsection*{Computing the boundary of  pseudospectra.} Trefethen (\cite{Tre99}) gave a survey of computing pseudospectra (as of 1999), which was then accompanied by the software package EigTool (Wright \cite{Wri02}, Wright \& Trefethen \cite{WriT01}). The basic algorithm is based on a contour plot of the smallest singular value of $A-zI$ for $z$ on a grid. Algorithms for tracing boundary curves of complex pseudospectra were developed by Br\"uhl (\cite{Bru96}), Mezher \& Philippe (\cite{MezP02}), and
Boulton \& Lancaster (\cite{BL10}).
Bekas \& Gallopoulos (\cite{BekG01}) combined curve-tracing and grids. 


Guglielmi \& Lubich (\cite{GL12}) exploited the rank-1 property of extremal perturbations in curve-tracing algorithms. Those algorithms are related to the tangential--transversal algorithm and the ladder algorithm described in Section~\ref{sec:ps-tracing}, which have not previously appeared in the literature. 
\chapter{Two-level approaches to matrix nearness problems}
\label{chap:two-level}

\index{matrix nearness problem!two-level algorithm}
We formulate a two-level approach to matrix nearness problems whose wide scope will become apparent in the following chapters. In this short chapter we discuss the two-level algorithm in a setting related to the class of eigenvalue optimization problems of Chapter~\ref{chap:proto}. We then apply it 
to the particular problem of computing the stability radius (or distance to instability): 

\medskip\noindent
{\bf Problem. }
{\it Given a matrix with all its eigenvalues having negative real part (a Hurwitz matrix), find the nearest complex matrix with some eigenvalue on the imaginary axis.}

\medskip\noindent
{\bf Problem. }
{\it Given a matrix with all its eigenvalues having modulus smaller than $1$ (a Schur matrix), find the nearest matrix with some eigenvalue on the complex unit circle.}

\medskip\noindent
Here, ``nearest'' will be understood to have minimal distance in the Frobenius norm.
These problems can be conveniently rephrased in terms of 
pseudospectra: 
The stability radius of a Hurwitz matrix is the smallest $\eps>0$ such that the $\eps$-pseudospectrum has some point on the imaginary axis; and the stability radius of a Schur matrix is the smallest $\eps>0$ such that the $\eps$-pseudospectrum has some point on the unit circle. 
More generally, for a given matrix having all eigenvalues in an open set $\Omega$ in the complex plane, the problem is to find the smallest $\eps>0$ such that the 
$\eps$-pseudospectrum has some point on the boundary of $\Omega$.

We consider a two-level approach that uses an {\it inner iteration} to compute the solution of an eigenvalue optimization problem as studied in Chapter~\ref{chap:proto} for a fixed perturbation size $\eps$, and then determines the optimal perturbation size $\oeps$ in an {\it outer iteration}. \bcl The hybrid expansion-contraction (HEC) algorithm developed by Mitchell \& Overton (\cite{MitO16}) and generalized by Mitchell \& Van Dooren (\cite{MitVD23}) uses the same building blocks but alternates between optimizing the perturbation matrix $E$ for a fixed $\eps$ and finding a root $\eps$ for a fixed perturbation matrix $E$. 
\ecl

Neither algorithm is guaranteed to find the global optimum of these nonconvex and nonsmooth optimization problems, but the algorithm computes a matrix with the desired spectral property which is locally nearest and often, as observed in our numerical experiments, has a distance close to or equal to the minimal distance. At convergence, it provides a rigorous upper bound to the minimal distance, and often a tight one.
 Running the algorithm with several different
starting values reduces the risk of getting stuck in a local optimum.

\section{Problem setting}
We consider  matrix nearness problems that are closely related to the eigenvalue optimization problems of Chapter~\ref{chap:proto}. 
For a given matrix $ A \in \C^{n,n}$, let $\lambda( A)\in\C$ be a target eigenvalue of~$ A$. We again consider the  smooth function
$f(\lambda,\clambda)$ satisfying (\ref{chap:proto}.\ref{ass:f}) that is to be minimized.
For a prescribed real number $\s$ in the range of $f$ we assume that
$$
f(\lambda( A),\clambda( A)) > \s,
$$
so that for sufficiently small $\eps>0$ we have $\phi(\eps)>\s$, where
$$
\phi(\eps) :=\min\limits_{\Delta \in \C^{n,n},\, \| \Delta \|_F = \eps} f \left( \lambda\left(  A + \Delta \right), \clambda \left(  A + \Delta \right)  \right).
$$
The objective is to find the smallest $\eps>0$ such that  $\phi(\eps)=\s$: 
\begin{equation}
\oeps = \min\bigl\{\eps > 0 \,:\, 
\phi(\eps) \le \s \bigr\}.
\label{eq:mnpb}
\end{equation}
Note that then $\phi(\oeps)=\s$. Determining $\oeps$ is a one-dimensional root-finding problem for the function $\phi$ that is defined by  eigenvalue optimization problems as studied in Chapter~\ref{chap:proto}.

\begin{example}[Stability radius of a Hurwitz matrix]
With the function 
$
f ( \lambda, \clambda) = -\tfrac12(\lambda + \clambda) = -\Re\, \lambda
$
and $\s=0$ and the target eigenvalue $\lambda(M)$ chosen as an eigenvalue of largest real part of a matrix $M$, we arrive at the classical problem of computing the stability radius (or distance to instability) $\oeps$ of a Hurwitz matrix $ A$, i.e. a matrix with negative spectral abscissa
$\alpha( A)= \max\{\Re\, \lambda\,:\, \lambda \text{ is an eigenvalue of $A$} \} < 0$\,:
$$
\oeps>0\quad\text{such that}\quad  \alpha_{\oeps}(A)=0,
$$
where $\alpha_\eps( A) = \max\limits_{E \in \C^{n,n}, \| E \|_F = 1} \alpha( A + \eps E)$  is the $\eps$-pseudospectral abscissa (see the previous chapter).
\end{example}

\begin{example}[Stability radius of a Schur matrix] With $f( \lambda, \clambda) = -\lambda\clambda = -|\lambda|^2$  and $\s=-1$
and the target eigenvalue $\lambda(M)$ chosen as an eigenvalue of largest modulus of a matrix $M$, we arrive at the problem of computing the stability radius of a Schur matrix $ A$, i.e. a matrix with spectral radius 
$\rho( A)=\max\{ |\lambda| \,:\,  \lambda  \text{ is an eigenvalue of $A$} \} < 1$\,:
$$
\oeps>0\quad\text{such that}\quad  \rho_{\oeps}(A)=1,
$$
where $\rho_\eps( A) = \max\limits_{E \in \C^{n,n}, \| E \|_F = 1} \rho( A + \eps E)$  is the $\eps$-pseudospectral radius.
\end{example}

\section{Nested iteration}
\label{sec:two-level}
\index{two-level iteration}

The approach is summarized by the following two-level method:
\begin{itemize}
\item {\bf Inner iteration:\/} Given $\eps>0$, we aim to compute a  matrix $E(\eps) \in\C^{n,n}$  
of unit Frobenius norm,  such that 
$\F_\eps(E) = f \left( \lambda\left(  A + \eps E \right), \clambda \left(  A + \eps E \right)  \right)$ 
is minimized, i.e. 
\begin{equation} \label{E-eps-2l}
E(\eps) = \arg\min\limits_{E \in \C^{n,n}, \| E \|_F = 1} \F_\eps(E).
\end{equation}


\item {\bf Outer iteration:\/} For the given value $\s$, we compute the smallest positive value $\oeps$ with
\begin{equation} \label{eq:zero}
\phi(\oeps)=\s,
\end{equation}
where $\phi(\eps)= \F_\eps\left(E(\eps) \right) = f \left( \lambda\left(  A + \eps E(\eps) \right), \clambda \left(  A + \eps E(\eps) \right)  \right)$.
\end{itemize}

\subsection{Inner iteration: Rank-1 constrained gradient flow}
\index{gradient flow!rank-1 constrained}
\index{optimizer!rank-1 property}
\index{rank-1 matrix differential equation}

The eigenvalue optimization problem \eqref{E-eps-2l} is of the type studied in Chapter~\ref{chap:proto}, which has rank-1 matrices as optimizers (see Corollary~\ref{chap:proto}.\ref{cor:rank-1}).
To compute $E(\eps)$ for a given $\eps>0$, we make use of a constrained gradient system for the functional $\F_\eps(E)$
under the constraints of unit Frobenius norm of $E \in \C^{n,n}$, which is further constrained to rank-1 matrices. We directly use the rank-1 constrained gradient flow approach developed in Chapter~\ref{chap:proto}, where we follow a suitable discretization of a rank-1 matrix differential equation (i.e., of a system of differential equations for two vectors) into a stationary point.

\bng
Note that, as discussed in Remark 
\ref{chap:proto}.\ref{rem:mult-eig}, even if it is convenient to assume that the target eigenvalue is simple along the trajectory and at the stationary points, this is not strictly necessary for the convergence of the inner iteration to a stationary point.
\eng

\bcl On the other hand note that, also if the inner iteration converges to a stationary point, that point need not be a global minimum. For the non-convex minimization problems that we consider here, the algorithm may get stuck in a local minimum. This can be made less likely by running several trajectories with different starting values, but there is in general no guarantee to arrive at the global minimum.
\ecl

\subsection{Outer iteration: standard Newton--bisection method} \label{sec:Newton-bisection}

In the outer iteration we compute $\oeps$, the smallest positive solution of the one-dimensional root-finding problem \eqref{eq:zero}.
This can be solved by a variety of methods, such as bisection. We aim for a locally quadratically convergent Newton-type method, which can be justified under regularity assumptions that are generic.
If these assumptions are not met, we can always resort to bisection. The algorithm proposed in the next subsection in fact uses a combined 
Newton--bisection approach.

\subsubsection*{Derivative for the Newton step.}
 In the following, an important role is again played by the  gradient $G_\eps(E)= 2f_{\clambda} xy^*$ as defined in
Lemma~\ref{chap:proto}.\ref{lem:gradient}. 

\begin{assumption} \label{ass:E-eps}
For $\eps$ close to $\oeps$ and $\eps<\oeps$, 
we assume the following for the optimizer $E(\eps)$ of \eqref{E-eps-2l}:
\begin{itemize}
\item The eigenvalue $\lambda(\eps)=\lambda( A+\eps E(\eps))$ is a simple eigenvalue.
\item The map $\eps \mapsto E(\eps)$ is continuously differentiable.
\item The gradient $G(\eps)=G_\eps(E(\eps))$ is nonzero.
\end{itemize}
\end{assumption}

Under this assumption, the branch of eigenvalues $\lambda(\eps)$ and its corresponding eigenvectors $\y(\eps), \x(\eps)$ with the scaling (\ref{chap:proto}.\ref{eq:scaling}) are also continuously differentiable functions of $\eps$ in a left neighbourhood of $\oeps$. 
We denote the eigenvalue condition number by
$$
\kappa(\eps) = \frac 1 {\y(\eps)^*\x(\eps)}>0.
$$
The following result gives us an explicit and easily computable expression for the derivative of 
$\phi(\eps)= \F_\eps(E(\eps)) = f( \lambda( \eps), \clambda( \eps)  )$ 
with respect to $\eps$ in terms of the gradient $G(\eps)$.
\begin{theorem}[Derivative for the Newton iteration] 
\label{thm:phi-derivative}
Under Assumption~\ref{ass:E-eps}, the function $\phi$ is continuously differentiable in a left neighbourhood of $\oeps$ and its derivative is 
\begin{equation} \label{eq:dereps}
\phi'(\eps)  = 
- \kappa(\eps) \,\| G(\eps) \|_F < 0.
\end{equation}
\end{theorem}
\begin{proof}
By Lemma~\ref{chap:proto}.\ref{lem:gradient} and (\ref{chap:proto}.\ref{eq:deriv-S}) we obtain, indicating by $'$ differentiation w.r.t.~$\eps$ and noting that $\frac d{d\eps}(\eps E(\eps)) = E(\eps) + \eps E'(\eps)$,
\begin{equation} \label{eq:deriveps}
\frac{1}{\kappa(\eps)} \,\frac{d}{d \eps} \F_\eps(E(\eps)) = 
\Re \bigl\langle  G(\eps),  E(\eps) + \eps E'(\eps) \bigr\rangle.
\end{equation}
By Theorem~(\ref{chap:proto}.\ref{thm:stat}), 
we know that in the stationary point $E(\eps)$, there exists a real $\mu(\eps)$ such that
\begin{equation} \label{E-mu-G}
E(\eps) = \mu(\eps) G(\eps).
\end{equation}
Since $\| E(\eps) \|_F=1$ for all $\eps$, we find  $1= |\mu(\eps)| \, \| G(\eps)\|_F $ (in particular $\mu(\eps)\ne 0$) and
$$
0= \frac12 \frac d{d\eps} \| E(\eps) \|^2 = \Re \langle E(\eps),E'(\eps) \rangle = \mu(\eps)\, \Re \langle G(\eps), E'(\eps) \rangle,
$$
so that 
$$
\Re \langle G(\eps), E'(\eps) \rangle =0.
$$
Inserting this relation into  \eqref{eq:deriveps} and using once again \eqref{E-mu-G}, we obtain
$$
\frac{1}{\kappa(\eps)} \, \phi'(\eps) = \Re \bigl\langle  G(\eps),  E(\eps) \bigr\rangle = 
\frac1{\mu(\eps)} \| E(\eps) \|_F^2 = \frac1{\mu(\eps)}= {\rm sign}(\mu(\eps))\, \| G(\eps) \|_F.
$$
Since for $\eps < \oeps$, we have $\phi(\eps) > \phi(\oeps)=\s$, and since the above formula shows that $\phi'$ cannot change sign, we must have $\phi'(\eps)<0$ and hence $\mu(\eps)<0$. This yields the stated result.
\qed 
\end{proof}\smallskip

\subsubsection*{Newton--bisection iteration.}
\label{subsec:Newton--bisection}
\index{Newton--bisection method}
In view of Theorem~\ref{thm:phi-derivative}, applying Newton's method to the equation $\phi(\eps)=\s$ yields the following iteration:
\begin{equation}\label{CNM1} \eps_{k+1} = \eps_{k}  +
\frac{\y(\eps_k)^* \x(\eps_k)}{\|G(\eps_k)\|_F} \left(\phi(\eps_k)-\s\right),
\end{equation}
where the right-hand side uses the optimizer $E(\eps_k)$ computed by the inner iteration in the $k$-th step.
\begin{figure}[ht]\label{fig:Hur}
\vskip -3cm
\centerline{
\includegraphics[scale=0.4]{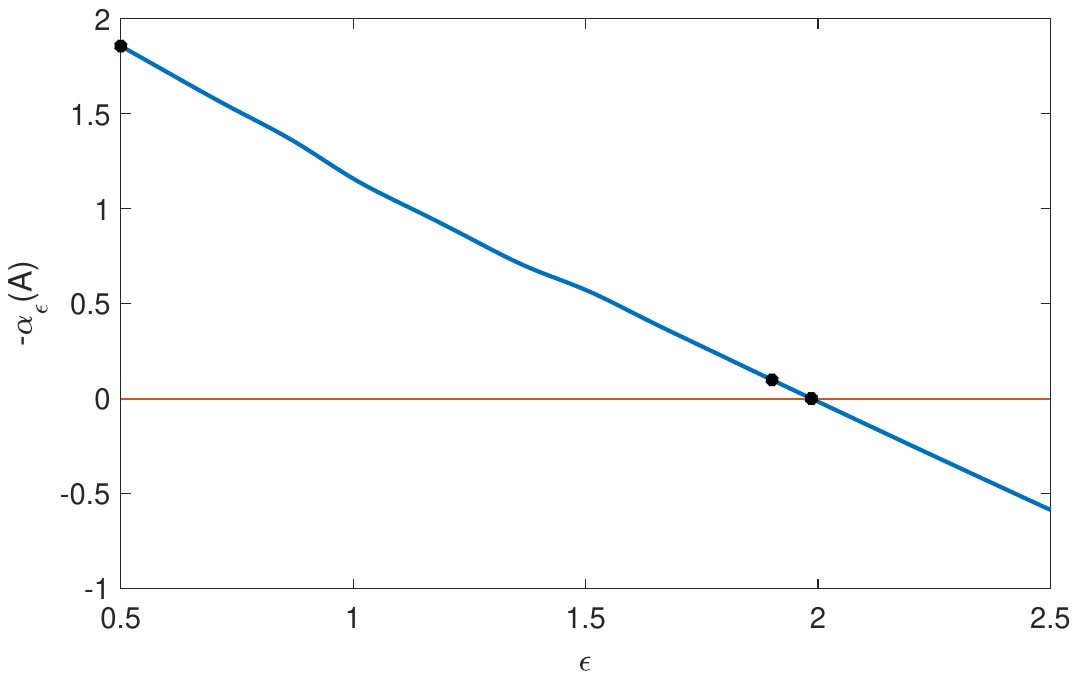}
}
\vskip -3cm
\caption{The function $\phi(\eps) = -\alpha_\eps(A)$ for the matrix
$A$ in (\ref{chap:proto}.\ref{eq:example}). The bullets are the points computed by the outer 
iteration.}
\end{figure}

\begin{algorithm}
\DontPrintSemicolon
\KwData{Matrix $ A$, 
$\s\in\R$, ${\rm tol}_0$ (initial tolerance), $k_{\max}$ (max number of iterations), \\
$\eps_{M}$ (initial upper bound for  $\oeps$)}
\KwResult{$\oeps$ if the inner iteration provides the global minimum and the outer Newton-bisection iteration converges; otherwise an upper bound $\wh \eps \ge \oeps$}
\Begin{
\nl Set $\lambda(0)$ target eigenvalue of $ A$, $\y(0)$ and $\x(0)$ the corresponding left and right
eigenvectors of unit norm with $\y(0)^*\x(0)>0$.\; 
\nl Initialize $E(\eps_0)$ according to the setting. \; 
\nl Initialize $\eps_0$ according to the setting.\; 
    Set $k=0$.\;
\nl Initialize lower and upper bounds: $\eps_{\rm lb}=0$, $\eps_{\rm ub}=\eps_M$.\; 
\nl Set $\phi(\eps_{-1})=\phi(\eps_0)+ 2\,{\rm tol}_0$\;
\While{$|\phi(\eps_k)-\phi(\eps_{k-1})| < {\rm tol}_k$}{
\nl Compute $E(\eps_k)$, $\phi(\eps_k)$ by integrating the constrained gradient system 
with initial datum $E(\eps_{k-1})$. 
(This is the {\bf inner iteration}).\;
\nl Update upper and lower bounds $\eps_{\rm lb}$, $\eps_{\rm ub}$.\; 
\nl \eIf{$\phi(\eps_k) < \s$} 
{Set $ \eps_{\rm ub} = \min(\eps_{\rm ub},\eps_k)$.} 
{Set $ \eps_{\rm lb} = \max(\eps_{\rm lb},\eps_k)$.}
\nl Compute $G(\eps_k)$.\;
\nl Compute $\eps_{k+1} = \eps_{k}+
\displaystyle \frac{\y(\eps_k)^* \x(\eps_k)}{\|G(\eps_k)\|_F} \left( \phi(\eps_k) -\s \right)$.\;
\nl Set $k=k+1$.\; 
\nl \If{$\eps_{k} \not\in [\eps_{\rm lb},\eps_{\rm ub}]$}
{\bng Set $\eps_{k} = \left(\eps_{\rm lb} + \eps_{\rm ub} \right)/2$. \eng}
\eIf{$k=k_{\max}$}
{goto $12$.}
{Set ${\rm tol}_k = \max \{ 10^{-2}\,{\rm tol}_{k-1}, 10^{-8} \}$.}\;
}
\nl \eIf{$k < k_{\max}$}
{Set $\wh\eps = \eps_k$. Return $\wh\eps$.}{Print {\em max number iterations reached}}
}
\caption{Outer iteration: Newton--bisection method}
\label{alg_SR} 
\end{algorithm}

Algorithm \ref{alg_SR} implements a hybrid Newton--bisection method that maintains
an interval known to contain the root, bisecting when the Newton step is outside
the interval $[ \eps_{\rm lb}, \eps_{\rm ub} ]$.

Step {\small \bf 5} (in the {\bf while} loop) gives the computational core of Algorithm \ref{alg_SR}; it implements the inner iteration and is not presented in detail since it depends on the possible structure of the matrix and on whether the low-rank structure is exploited.
The inner iteration performs the algorithm to compute the extremal perturbation $E(\eps_k)$, as described in Chapter~\ref{chap:proto}.
As input to the $k$-th iteration we use
the factors of the final matrix $E$ computed for the previous value of $\eps$ (this explains
the choice of the initial datum at step {\small \bf 5}).

The {\bf while} loop after step {\small \bf 4} implements the outer iteration and makes use of a variable tolerance which decreases as $k$ increases, when the method is expected to approach convergence.
A typical choice of ${\rm tol}_0$ is the norm of the difference of the first two iterates divided by $10$.

The factor $10^{-2}$ between two subsequent tolerances is fruit of an empirical experimentation and is motivated by the fact that we expect convergence - on the average - in about $4$ to $5$ iterates so that we reach the limit tolerance which is set here to $10^{-8}$. These numbers can naturally be considered as parameters of the code. 
Their choice here is only based on the experience with
the numerical experiments we performed.    
When integrating numerically the gradient system with variable stepsize we are assured to fulfill
the termination condition $|\phi(\eps_k)-\phi(\eps_{k-1})| < {\rm tol}_k$ of the {\bf while} loop because we are able to approximate the stationary point with prescribed accuracy. 
Due to the possible convergence of the inner method to a local instead of global minimum, the final value $\wh\eps$
computed by Algorithm \ref{alg_SR} might be larger than the minimal one. 
\bng
\subsection{Starting values} 
\label{sec:starteps}
\index{initial perturbation}

As in subsection~\ref{subsec:init} we consider two options: (i) we make the following choice of the initial perturbation matrix $E_0$ and the initial perturbation size
$\eps = \eps_0$ for the first iteration: let $G$ be the  gradient of $\F_\eps$ at $\eps=0$, i.e.,
$G= 2 f_{\clambda}(\lambda,\clambda) xy^*$ with the target eigenvalue $\lambda$ of the unperturbed matrix $A$ and its normalized left and right eigenvectors $x$ and $y$ with positive inner product. We set
$$
E_0 = -\frac{G}{\|G\|_F},
$$
which is the steepest descent direction at the unperturbed matrix for the functional $\F_\eps(E)$.
At least for $\eps_0$ not too large, this yields that
$\F_{\eps_0}(E_0)$ is smaller than $f(\lambda,\clambda)$, the function value assumed at the target eigenvalue of the unperturbed matrix $A$.

We formally apply 
a Newton step using formula (\ref{CNM1}) for $k=-1$ with $\eps_{-1}=0$ and set
\begin{equation} \label{eq:eps0F}
\eps_{0}=  \frac{\y^*\x}{\|G\|_F} (f(\lambda,\clambda)-\s).
\end{equation}
We check if $\F_{\eps_0}(E_0) < f(\lambda,\clambda)$, and if this is not satisfied we halve $\eps_0$ and check and halve again until this is satisfied.

A better choice -- consistent with (ii) in Section~\ref{subsec:init} -- selects $m$ eigenvalues $\lambda_i$ ($i=1,\dots,m$) of $A$ for which $f(\lambda_i,\clambda_i)$ take values that are not far from that of the target eigenvalue of $A$, with corresponding normalized left and right eigenvectors $x_i$ and $y_i$ with positive inner product, and $\kappa_i=1/(x_i^*y_i)$. We set, for $i=1,\dots,m$,
\begin{align*}
E_i &= -\frac{G_i}{\| G_i \|_F} \quad \mbox{ with} \quad
G_i = 2 \,f_{\clambda}(\lambda_i,\clambda_i) \,x_i y_i^*, 
\\
\eps_{i} &=  \frac{x_i^*y_i}{\|G_i\|_F} (f(\lambda_i,\clambda_i)-\s).
\end{align*}
We then choose 
\begin{equation} \label{eq:choicesr}
E(0) = E_j, \ \ \eps_0=\eps_j  \quad \mbox{ with} \quad
 j= \arg\min_{i=1,\dots,m} \eps_i.
\end{equation}
We note that here the right-hand terms that appear in (\ref{chap:proto}.\ref{eq:choice}) are all equal to $\s$: 
$$
f\left( \lambda_i, \clambda_i \right)  - \eps_i\kappa_i \, 2| f_{\clambda}(\lambda_i,\clambda_i)| =
\s .
$$
We therefore choose $\eps_0$ as the smallest among the $\eps_i$.

Extending the more elaborate approaches (iii) and (iv) in Section~\ref{subsec:init} is yet another option.
\medskip

While the above choice of starting values is reasonable when only a single trajectory is computed, it might in some problems be necessary to run several trajectories to reduce the risk of getting trapped in a local minimum.
\eng

\subsection{Possible failure of convergence} \label{sec:breakdown}
Convergence of the Newton--bisection iteration in combination with the inner iteration is not guaranteed, even if in our numerical experiments it rarely failed. 

\bcl
Failure of convergence can arise if
the inner iteration does not arrive at a (unique) global minimum.
So far in this chapter, we assumed that
the inner iteration computes
$$
E(\eps) = \arg\min\limits_{E \in \C^{n,n}, \| E \|_F = 1} \F_\eps(E).
$$
There are, however,  situations in which the inner iteration
is only able to compute some local minimum $\wh E(\eps)$ of $\F_\eps$ instead of a desired global minimum $E(\eps)$ and the outer iteration switches between different paths of local or global minima, so that $\wh E(\eps)$ does not depend continuously on~$\eps$. In such cases, the Newton--bisection method with the computed $\wh E(\eps)$ instead of $E(\eps)$ need
not converge.

{\it If\/} the Newton--bisection method converges, then its limit $\wh\eps$ provides an upper bound of $\oeps$: For a subsequence for which $\wh E(\eps_k)$ converges to a limit denoted as $\wh E(\wh \eps)$, we have
\begin{equation}\label{wheps-inequality}
 \s = \lim_{k\to\infty} \F_{\eps_k} (\wh E(\eps_k)) = \F_{\wh\eps}(\wh E(\wh \eps)) \ge \F_{\wh\eps}( E(\wh \eps)).
\end{equation}
On the other hand, $\oeps$ is the smallest $\eps>0$ for which $\F_{\eps}( E(\eps))\le \s$, and so we conclude that $\wh \eps$ is an upper bound of $\oeps$:
\begin{equation}\label{oeps-upper-bound}
\oeps \le \wh \eps.
\end{equation}
Equality holds if and only if $\wh E(\wh\eps)$ is a global minimum of $\F_{\wh\eps}$.

\bcltwo
\subsubsection*{Alternative approaches, which have monotone convergence.}
In the following we describe two alternatives to the previously stated Newton--bisection method that are guaranteed to yield a convergent sequence $( \eps_k )_{k \ge 0}$.  
\begin{itemize}
    \item One alternative approach imposes monotone convergence of the sequence $( \eps_k )_{k \ge 0}$ {\it\ from the left} within a Newton--bisection framework. This is described in the next subsection.
  \item A fundamentally different approach is given by the  hybrid expansion--contraction  algorithm introduced by  Mitchell and Overton (\cite{MitO16}) and generalized by Mitchell \& Van Dooren (\cite{MitVD23}), which yields monotone convergence of a sequence $( \eps_k )_{k \ge 0}$ {\it\ from the right}.
This is described in Section~\ref{sec:HEC}.
\end{itemize}
In both cases, the limit $\wh\eps$ thus provides a rigorous upper bound $\wh\eps\ge\oeps$. While the sequence $(\wh E(\eps_k))$ need not necessarily converge, any accumulation point $\wh E(\wh\eps)$ is a stationary point provided that the function 
$(\eps,E)\mapsto f(\lambda(A+\eps E),\conj{\lambda(A+\eps E)})$ is continuously differentiable at $(\wh\eps,\wh E)$, which holds true if the target eigenvalue $\lambda(A+\wh\eps\wh E)$ is a simple eigenvalue.

\begin{remark}
For many matrix nearness problems, either of the two approaches can be equally used. In some problems, however, such as the structured distance to singularity in Section~\ref{sec:sing-S}, $\varphi(\eps)$ assumes the constant value $\s$ for all $\eps\ge \oeps$. In such a case, the HEC method stagnates at the starting value $\eps_0>\oeps$, whereas the monotone Newton--bisection converges from the left to an $\wh\eps \ge \oeps$, with equality if the stationary point $\wh E(\eps)$ from the inner iteration equals the global minimizer $E(\eps)$ for all $\eps<\oeps$ near $\oeps$.
\end{remark}
\ecltwo
\bcl

\subsection{Newton--bisection method with monotone convergence}
\label{subsec:Newton-bisection-monotone}
\index{Newton--bisection method!monotone convergence}
Let $\eps_k>0$ and $E_k\in\C^{n,n}$ of Frobenius norm 1 be given, and set $\phi_k=\F_{\eps_k}(E_k)$.
Assume $\phi_{k} > \s$ and let $\eps_{k+1}$ be the result of a Newton step as given by \eqref{CNM1},
 which evidently yields $\eps_{k+1} > \eps_k$. Then we proceed as follows. 
\medskip

\begin{itemize}
\item[1. ] 
In the inner iteration we first integrate the  rank-1 or full-rank gradient system without the norm constraint, see \eqref{unconstrained-rank-1-flow}, starting with the initial value for $E$ given by
\[
E(0)=\frac{\eps_{k}}{\eps_{k+1}}  E_{k}, 
\]
which is such that $\| E(0) \|_F < 1$ and $\F_{\eps_{k+1}}(E(0))=\F_{\eps_{k}} (E_{k})$.
We integrate up to the first time $\bar t$ at which $E(\bar t)$ is of Frobenius norm~1, before continuing with the integration of the norm-constrained rank-1 gradient system \eqref{ode-E-1} into a stationary point $E_{k+1}$. 

\hskip 5mm 
In this way we obtain that the concatenation of the two trajectories of the target eigenvalue for $\eps=\eps_{k}$ and $\eps=\eps_{k+1}$ is continuous. We then have
\begin{equation}\label{phi-monotone}
\phi_{k+1} \le \phi_{k},
\end{equation}
because the functional $\F_{\eps_{k+1}}$ decays monotonically along the concatenated trajectory and hence $\phi_{k+1} =\F_{\eps_{k+1}} ( E_{k+1}) \le \F_{\eps_{k+1}}(E(0))=\F_{\eps_{k}} (E_{k}) = \phi_{k}$.
\item[2. ] If $\phi_{k+1} \ge \s$, 
then we accept $\eps_{k+1}$  and proceed by computing $\eps_{k+2}$ according to the Newton formula \eqref{CNM1}. We then continue to 1. with $k+1$ in the role of $k$. 
\item[3. ] 
Otherwise, we do not accept the $\eps$-step. We just update $\eps_{k+1}$,  replacing $\eps_k$ by $\frac12(\eps_k+\eps_{k+1}),$
and repeat 1. 
\end{itemize}

\bcltwo
\medskip\noindent
While we described Step 1. in terms of unconstrained and constrained gradient flows, in the following it is only important to use an algorithm in the inner iteration that yields \eqref{phi-monotone} and a bounded sequence $(\eps_k)$.

\ecltwo

\begin{algorithm}
\DontPrintSemicolon
\KwData{$\s\in\R$, matrix $A$ with $f(\lambda(A),\overline{\lambda(A)})>\s$, 
initial perturbation size $\eps_0 > 0$ and initial perturbation matrix $E_0$ of Frobenius norm $1$ such that 
$\phi_0 = \F_{\eps_0}(E_0) > \s$ }
\KwResult{$\wh\eps$ (upper bound for $\oeps$) and $\wh E$ (stationary point of $\F_{\wh\eps}$ 
with $\F_{\wh\eps}(\wh E)=\s$)}
\Begin{
\nl Compute $\eps_1$ by a Newton step \eqref{CNM1}.\qquad \% note $\eps_1>\eps_0$\;
\For{$k=0,1,2,\ldots$} 
{
\nl Initialize ${\bf left}=0$\;
\While{${\bf left}=0$}{
\nl Determine a stationary point $E_{k+1}$ (of Frobenius norm 1) of $\F_{\eps_{k+1}}$ such that 
$\phi_{k+1}= \F_{\eps_{k+1}}(E_{k+1})\le\F_{\eps_{k}}(E_{k})=\phi_k$,
\\ \hspace{2mm}e.g.~using both unconstrained and norm-constrained (rank-$1$ or full-rank)  
\\ \hspace{2mm}gradient systems as described in Step 1. above.\;
\nl \eIf{$\phi_{k+1} < \s$}{Set $\eps_{k+1} = 
\displaystyle
\frac{\eps_{k+1} + \eps_k}{2}$ \qquad\ \ \ \ \ \ \% note $\eps_{k+1} > \eps_{k}$}{Set ${\bf left}=1$} 
}
\nl Compute $\eps_{k+2}$ by a Newton step \eqref{CNM1}. \quad\% note $\eps_{k+2} > \eps_{k+1}$ unless $\phi_{k+1}=\s$\;
\nl {\bf If} $\eps_{k+2}=\eps_{k+1}$ (up to some tolerance) {\bf then return} $\wh\eps=\eps_{k+1}$ and $\wh E =  E_{k+1}$.\;
}
}
\caption{Monotone Newton--bisection method}
\label{alg:SRC} 
\end{algorithm}

\bcl
Algorithm \ref{alg:SRC} summarizes the procedure described above.
We state properties of the algorithm in the following theorem. 

\bcltwo
\begin{theorem} [Newton--bisection with monotone convergence]
\label{thm:newton-bisection-monotone}
Let $\eps_0>0$ and $E_0$ of Frobenius norm 1 be given and $\phi_0= \F_{\eps_0}(E_0) > \s$.
Then Algorithm~\ref{alg:SRC} yields the following.

(a) The algorithm generates a monotonically increasing sequence $(\eps_k)$, a sequence of matrices $(E_k)$ of Frobenius norm~1, and a monotonically decreasing sequence $(\phi_k)$ bounded from below by $\s$. 

(b) If the sequence $(\eps_k)$ is bounded, then the associated sequence $(\phi_k)$ converges to $\s$.
Furthermore, let $\wh\eps=\lim_{k\to\infty} \eps_k$  and let $\wh E$ be an accumulation point of the sequence of matrices $(E_k)$ 
for which 
$\lambda(A+\wh\eps \wh E)$ is a simple eigenvalue. Then $\wh E$ is a stationary point of $\F_{\wh\eps}$ and $\F_{\wh\eps}(\wh E)=\s$.
\end{theorem}

\begin{proof}  (a) The lower bound $\phi_k\ge \s$ is ensured by the assumption on the initial data and the construction of $\eps_k$.  With this bound, the
monotonicity of the sequence $(\eps_k)$ then follows immediately from the Newton iteration formula \eqref{CNM1}.

(b) Under the assumption that the monotone sequence $(\eps_k)$ is bounded, it converges to some limit $\wh\eps$.
Taking the limit in the Newton iteration \eqref{CNM1} shows that $\phi_k \to \s$. By the continuity of the function $(\eps,E)\mapsto \F_\eps(E)$, this shows that for every accumulation point $\wh E$ of $(E_k)$ we have $\F_{\wh\eps}(\wh E)=\s$. 
This implies $\wh \eps \ge \oeps$ by the definition of $\oeps$. Finally, the limit $\wh E$ of a sequence of stationary points is again a stationary point if the function $(\eps,E)\mapsto \F_\eps(E)$ is continuously differentiable in a neighbourhood of the limit, which is the case if $\lambda(A+\wh\eps \wh E)$ is a simple eigenvalue. 
\qed
\end{proof}

Next we show that the monotone sequence $(\eps_k)$ indeed remains bounded when the unconstrained and constrained gradient flows as described in Step 1.~above or their rank-1 counterparts are used in line 3 of Algorithm~\ref{alg:SRC}.

\begin{lemma} [Boundedness of the perturbation sizes $\eps_k$]
    Let $\eps_0>0$ and $E_0$ of Frobenius norm 1 be given and $\phi_0= \F_{\eps_0}(E_0) > \s$.
Assume that there exists $\alpha>0$ such that 
\begin{equation}\label{df-lower-bound}
\biggl| \frac{\partial f}{\partial\clambda}(\lambda,\clambda) \biggr| \ge \alpha \quad\text{ for all }\ 
\lambda \in \C \text{ such that } f(\lambda,\clambda) \le f_0,
\end{equation}
where $f_0=f(\lambda(A+\eps_0 E_0), \clambda(A+\eps_0 E_0))$ is the function value at a target eigenvalue of $A+\eps_0 E_0$.
Then, the monotone Newton--bisection algorithm using full-rank gradient systems in line 3 of Algorithm~\ref{alg:SRC} yields a sequence $(\eps_k)$ bounded by $\bar\eps = \eps_0 + (\phi_0-\s)/\alpha$.
Under the additional assumption
\begin{equation} \label{nb-ass-rank-1}
   f_0 < f(\lambda(A), \clambda(A)),
\end{equation} 
the same result with a possibly larger $\bar\eps$ holds true when instead the rank-1 projected gradient is used.
\end{lemma}

\begin{proof}  The trajectory of the unconstrained gradient flow either runs into a stationary point or is unbounded. 
Since the condition on $f$ excludes a vanishing gradient by Lemma~\ref{chap:proto}.\ref{lem:gradient}, the trajectory is unbounded and hence there is a finite time $\bar t$ at which $E(\bar t)$ has Frobenius norm~1. Hence the sequence $(E_k)$ of matrices of Frobenius norm 1 is well-defined. 

Assume $\eps_k \le \bar\eps$. We will show that $\eps_{k+1}\le \bar\eps$.
According to the algorithm, we first integrate the free gradient system
\begin{eqnarray} 
\dot{E}(t) & = & - \frac{G(t)}{\|G(t)\|_F}, \qquad t \ge 0,
\label{ode-free}
\\
E(0) & = & \frac{\eps_{k}}{\eps_{k+1}}  E_{k}.
\nonumber
\end{eqnarray} 
Here $E_k$ (of Frobenius norm~1) is a stationary point of the inner iteration with $\eps=\eps_k$, and 
$G(t)=2f_{\clambda}(t) \, x(t) y(t)^*$ is the rescaled gradient matrix as in Lemma~\ref{chap:proto}.\ref{lem:gradient},
where $x(t)$ and $y(t)$ are normalized left and right eigenvectors with positive inner product associated with the target eigenvalue $\lambda(t)$ of $A + \eps_{k+1} E(t)$, and $f_{\clambda}(t)= (\partial f/\partial\clambda)(\lambda(t),\clambda(t))$.
For $t \le t_{k+1} := \displaystyle \frac{\eps_{k+1}-\eps_{k}}{\eps_{k+1}}$ we have 
\begin{align*}
\| E(t) \|_F & = \| E(0) - \int\limits_{0}^{t} \frac{G(t)}{\|G(t)\|_F}\, {\rm d}s \|_F
\le  \| E(0) \|_F +  \int\limits_{0}^{t} 1   \,{\rm d}s  
\\
& =  \frac{\eps_{k}}{\eps_{k+1}} +  t \le
\frac{\eps_{k}}{\eps_{k+1}} + \frac{\eps_{k+1}-\eps_{k}}{\eps_{k+1}} = 1.
\nonumber
\end{align*}
Moreover, by Lemma~\ref{chap:proto}.\ref{lem:gradient} and \eqref{ode-free} we have $\displaystyle \frac{d}{dt} f(\lambda(t),\clambda(t))  = - \eps_{k+1} \kappa(t) \| G(t) \|_F$ with $\kappa(t)=1/(x(t)^*y(t))$, and hence
\begin{equation} \nonumber
f(\lambda(t),\clambda(t)) = f(\lambda(0),\clambda(0)) - \int\limits_{0}^{t} \eps_{k+1} \kappa(s) \| G(s) \|_F \, {\rm d} s.
\end{equation}
Since $\kappa(s) \ge 1$ for all $s$ and $\| G(s) \|_F = |f_\lambda(s) | \ge \alpha>0$ by assumption \eqref{df-lower-bound}, we obtain
$$
f(\lambda(t_{k+1}),\clambda(t_{k+1})) \le f(\lambda(0),\clambda(0)) - t_{k+1} \alpha \eps_{k+1} =
f(\lambda(0),\clambda(0))  - \alpha\left( \eps_{k+1}-\eps_{k} \right) ,
$$
and further we have, with $\lambda_k$ a target eigenvalue of $A+\eps_k E_k$,
\begin{align*}
&\phi_{k+1} = f(\lambda_{k+1},\clambda_{k+1}) \le  f(\lambda(t_{k+1}),\clambda(t_{k+1})) 
\ \text{ and } 
\\ 
&f(\lambda(0),\clambda(0))  
= f(\lambda_{k},\clambda_{k}) = \phi_k ,
\end{align*}
so that
\begin{equation} \label{eq:lambda-mono}
\s \le \phi_{k+1} \le \phi_k - \left( \eps_{k+1}-\eps_{k} \right)\alpha ,
\end{equation}
where the first inequality is ensured by the algorithm.
By induction this implies 
$$
\alpha(\eps_{k+1}-\eps_0) \le \phi_0-\phi_{k+1} \le \phi_0- \s,
$$ 
that is,
$$
\eps_{k+1}\le \eps_0 + \frac{\phi_0-\s}\alpha = \bar\eps.
$$

\medskip
{\it Case of the rank-1 algorithm.} Here, the gradient $G=f_{\clambda}xy^*$ is replaced by the projected gradient $P_{E} (G)$ in \eqref{ode-free} and in the subsequent lines in part (b).  $P_E$ is the orthogonal projection onto the tangent space at $E=uv^*$ of the manifold of rank-1 matrices, i.e. $P_E(G) = G - (I-uu^*)G(I-vv^*)$ by (\ref{chap:proto}.\ref{P-formula-1}). If for some $\beta>0$,
\begin{equation} \label{PG-lower-bound}
   \| P_{E(t)} G(t) \|_F \ge \beta > 0 \qquad\text{for all $t\in [0,t_k]$ in every step $k$}, 
\end{equation}
then the above arguments yield $\eps_{k+1} \le \eps_0 + ({\phi_0-\s})/\beta$. It remains to prove \eqref{PG-lower-bound}.

 In the following we emphasize the dependence on $k$ in the notation in writing $E_k(t)=u_k(t)v_k(t)^*$ and $G_k(t)= f_{\clambda,k}(t) x_k(t) y_k(t)^*$.  The monotonic decrease of $\F_{\eps_k}(E_k(t))$ with growing $t$ and $k$ shows that under assumption \eqref{nb-ass-rank-1}, $A+\eps_{k}E_k(t)$ and $A$ cannot have the same target eigenvalue (recall the definition of a target eigenvalue in Section~\ref{subsec:proto-problem}) and the same associated eigenvectors.
Remark~\ref{chap:proto}.\ref{rem:exceptional} then shows that $P_{E_k(t)} (G_k(t))\ne 0$ for all $t$, or equivalently by \eqref{df-lower-bound}, 
$$
P_{E_k(t)}(x_k(t)y_k(t)^*)\ne 0.
$$
By continuity of $u_k(t),v_k(t)$ and $x_k(t),y_k(t)$, we have in every step $k$ that
$$
\beta_k :=\min_{0\le t \le t_{k}} \|P_{E_k(t)}(x_k(t)y_k(t)^*)\|_F >0.
$$
Let $s_k\in\arg\min_{0\le t \le t_{k}} \|P_{E_k(t)}(x_k(t)y_k(t)^*)\|_F$.

Assume by contradiction $\lim_{k\to\infty} \beta_k=0$.
Then there exist convergent subsequences $E^k=E_k(s_k)\to \wh E$ and normalized eigenvectors $x^k=x_k(s_k)\to\wh x$, $y^k=y_k(s_k)\to \wh y$ 
such that $\lim_{k\to\infty} \| P_{E^k}(x^k (y^k)^*)\|_F =0$ and again by continuity, 
$$
P(\wh E)(\wh x \wh y^*)=0.
$$
By Remark~\ref{chap:proto}.\ref{rem:exceptional}, this implies that $A$ and $A+ \eps \wh E$  have the same target eigenvalue $\lambda(A)$ and corresponding eigenvectors for every $\eps>0$, which contradicts the fact that by \eqref{nb-ass-rank-1},
$$
f(\lambda(A+\eps_k E^k), \clambda(A+\eps_k E^k)) < f_0 < f(\lambda(A), \clambda(A))\quad\text{  for all $k\ge 0$.}
$$
This contradiction proves that $\beta = \inf_{k\ge 0} \beta_k>0$ and \eqref{PG-lower-bound} is satisfied with this $\beta$.
\qed
\end{proof}
\ecltwo

\ecl

\bcl
\section{Alternating between optimization and root-finding: HEC algorithm}
\label{sec:HEC}
\index{Hybrid expansion-contraction method}
\index{HEC method}
 The hybrid expansion--contraction (HEC) method was developed by Mitchell \& Overton (\cite{MitO16}) and 
 \bcltwo
 generalized by Mitchell \& Van Dooren (\cite{MitVD23}) to a general class of root-max and root-min problems, which in particular includes the problems of eigenvalue optimization combined with root-finding for the optimal distance $\eps$ as considered here.
 \ecltwo
  \bcl
In the present context, the HEC method converges monotonically in $\eps$ and yields a rigorous upper bound also in cases where the standard Newton--bisection method of the outer iteration does not converge. While the monotone Newton--bisection method converges from the left with a monotonically increasing sequence of iterates, the HEC method converges from the right with a monotonically decreasing sequence. The HEC method is formulated in Algorithm~\ref{alg:hec} in a semi-abstract form without practical stopping criteria, in the notation used in this chapter.

\begin{algorithm}
\DontPrintSemicolon
\KwData{$\s\in\R$, matrix $A$ with $f(\lambda(A),\overline{\lambda(A)})>\s$, 
initial perturbation size $\eps_0 > 0$ and initial perturbation matrix $E_0$ of Frobenius norm 1 such that $\F_{\eps_0}(E_0) \le \s$ }
\KwResult{$\wh\eps$ (upper bound for $\oeps$) and $\wh E$ (stationary point  of $\F_{\wh\eps}$ with $\F_{\wh\eps}(\wh E)=\s$)}
\Begin{
\nl \For {$k=0,1,2,\ldots$} 
{
\nl Determine $\eps_{k+1}\le \eps_k$ such that $\F_{\eps_{k+1}}(E_k)=\s$
\\ \hspace{2mm} (e.g.~by Newton--bisection iteration).\;
\nl {\bf If} $\eps_{k+1}=\eps_k$ (up to some tolerance) {\bf then return} $\wh\eps=\eps_k$ and $\wh E = E_k$.\;
\nl Determine a stationary point $E_{k+1}$ (of Frobenius norm 1) of $\F_{\eps_{k+1}}$ such that 
$\F_{\eps_{k+1}}(E_{k+1})\le\F_{\eps_{k+1}}(E_{k})=\s$ 
\\ \hspace{2mm}(e.g.~by a rank-1 gradient method as in Chapter~\ref{chap:proto}).\;
}
}
\caption{HEC algorithm (formulated in the present context)}
\label{alg:hec} 
\end{algorithm}

Mitchell \& Overton (\cite{MitO16}) call the root-finding step in line 2 the {\it contraction phase} and the optimization step in line 4 the {\it expansion phase}. They consider maximization, in which case the inequality in line 4 is reversed and the functional is expanded (if $\s=0$, as they choose without loss of generality). With minimization as considered here, we nevertheless use the established acronym HEC even if there is no expansion.

The HEC method can still be viewed as a two-level method, but it does not fit into the framework of the 
inner--outer iteration method outlined in Section~\ref{sec:two-level}. Instead, it {\it alternates} between eigenvalue optimization for a fixed perturbation size (e.g.~via the {rank-1} gradient method of Chapter~\ref{chap:proto}) and root-finding of $\F_\eps(E)-\s$ for a fixed perturbation matrix $E$ (e.g.~via a Newton--bisection method). There are the same two building blocks or ``levels'' of eigenvalue optimization and root-finding as previously, but there are no inner and outer levels. Here the same matrix $E_k$ is used in several Newton--bisection iteration steps until a convergence criterion is met, whereas $E$ is updated after every Newton--bisection iteration step in the nested iteration of Section~\ref{sec:two-level}.

It is important to note that the equation $\F_{\eps_{k+1}}(E_k)=\s$ in line 2 always has a solution $\eps_{k+1}\le \eps_k$, because $\F_\eps(E)$ depends continuously on $\eps$ and
$$
\F_{0}(E_k) = f(\lambda(A),\overline{\lambda(A)}) > \s \ge \F_{\eps_k}(E_k),
$$
where the last inequality is assumed to hold for $k=0$ and follows from the construction of $E_k$ for $k\ge 1$ in line 4. The intermediate value theorem then yields the existence of an $\eps_{k+1}\in(0,\eps_k]$ with $\F_{\eps_{k+1}}(E_k)=\s$, and Newton--bisection iteration with bracketing converges. 

This observation leads to the following result given in Theorem 3.7 of Mitchell \& Van Dooren (\cite{MitVD23}), which we rephrase in the context of this chapter.

\begin{theorem} [Convergence of the HEC algorithm]
\label{thm:hec}
The HEC algorithm starting with $\eps_0>\oeps$ and $\F_{\eps_0}(E_0) < \s$ generates a monotonically decreasing sequence $(\eps_k)$ that converges to a limit $\wh\eps\ge\oeps$. 
Furthermore, let $\wh E$ be an accumulation point of the sequence of matrices $(E_k)$ 
for which 
$\lambda(A+\wh\eps \wh E)$ is a simple eigenvalue. Then $\wh E$ is a stationary point of $\F_{\wh\eps}$ and $\F_{\wh\eps}(\wh E)=\s$.
\end{theorem}

\begin{proof} Since the sequence $(\eps_k)$ is monotonically decreasing and bounded by 0 from below, it converges to a limit $\wh\eps$. Since $\F_{\eps_{k+1}}(E_k)=\s$ , we have for every accumulation point $\wh E$ of the bounded sequence $(E_k)$ that by continuity
$ \F_{\wh \eps}(\wh E)=\s$, which implies $\phi(\wh \eps) = \s$ and hence $\wh \eps \ge \oeps$ by the definition of $\oeps$. Finally, the limit $\wh E$ of a sequence of stationary points is again a stationary point if the function $(\eps,E)\mapsto \F_\eps(E)$ is continuously differentiable in a neighbourhood of the limit, which is the case if 
$\lambda(A+\wh\eps \wh E)$ is a simple eigenvalue.
\qed
\end{proof}

Moreover, Mitchell \& Overton (\cite{MitO16}), Theorem 4.4, and Mitchell \& Van Dooren (\cite{MitVD23}), Theorem 3.10, prove quadratic convergence of the HEC algorithm under an additional regularity assumption.

\ecl


\section{Computing the stability radius}
\label{sec:dist-instab}
\index{stability radius}
The standard Newton-bisection iteration, the monotone Newton-bisection iteration, and the HEC method each apply directly to computing the stability radius (or distance to instability) of a Hurwitz-stable matrix $A$. We choose the target eigenvalue $\lambda(M)$ as an eigenvalue of largest real part (among these, the one with largest imaginary part), and we take $f(\lambda,\clambda)=-\,\Re\,\lambda$ as the function to be minimized. The eigenvalue optimization problem
\eqref{E-eps-2l}
then becomes the maximization problem
\begin{equation} \label{E-eps-stab-radius}
E(\eps) = \arg\max\limits_{E \in \C^{n,n}, \| E \|_F = 1} \Re\,\lambda(A+\eps E),
\end{equation}
and the optimal perturbation size $\oeps$ is determined from Equation \eqref{eq:zero}, which here reads
\begin{equation} \label{zero-stab-radius}
\Re\,\lambda(A+\oeps E(\oeps)) =0.
\end{equation}
In the inner iteration, the eigenvalue optimization problem \eqref{E-eps-stab-radius} is of the class studied in Chapter~\ref{chap:proto} 
with $f(\lambda,\clambda)=-\Re\,\lambda$ and is solved with the rank-1 constrained gradient flow approach developed there.

\bcltwo
In the outer iteration we compute the optimal perturbation size $\oeps$ by the monotone Newton--bisection algorithm of Section~\ref{subsec:Newton-bisection-monotone} for $f(\lambda,\clambda)=-\Re\, \lambda$, using the derivative formula of Theorem~\ref{thm:phi-derivative} with  $G=-xy^*$.

We consider a numerical example for the matrix $A$ in (\ref{chap:proto}.\ref{eq:example}) shifted by $\: -4 I$ ($4$ times the negative identity), which makes it Hurwitz,
with $\phi(\eps) = -\alpha_{\eps}(A)=-\Re\,\lambda(A+\eps E(\eps))$ (see Table~\ref{tab:Hur}). 
We note, however, that all values of $\phi(\eps_k)$ are positive
and so the criterion at line 4 of Algorithm \ref{alg:SRC} was never applied in this example.
\begin{table}[hbt]
\begin{center}
\begin{tabular}{|l|l|l|l|}\hline
  $k$ & $\eps_k$ & $\phi(\eps_k)$ & $\#$ iters \\
 \hline
\rule{0pt}{9pt}
\!\!\!\! 
    $0$         & $0.5$ & $1.85659$ & $27$  \\
	$1$         & $1.828513855227206$ & $1.83031\,10^{-1}$ & $51$  \\
	$2$         & $1.985028328768070$ & $9.93050\,10^{-4}$ & $45$  \\
	$3$         & $1.985886608426606$ & $2.71296\,10^{-9}$ & $31$  \\
    $4$         & $1.985886631875649$ & $2.01231\,10^{-15}$ & $1$ \\ 
 \hline
\end{tabular}
\vspace{2mm}
\caption{Computation of the stability radius for the matrix $A -4 I$ (with $A$ the matrix in (\ref{chap:proto}.\ref{eq:example}): 
computed values $\eps_k$, $\phi(\eps_k) =  -\Re\,\lambda (A + \eps_k E_k)$ 
and number of eigenvalue computations of the inner rank-$1$ algorithm.\label{tab:Hur}}
\end{center}
\begin{center}
\begin{tabular}{|l|l|l|l|}
\hline
$k$ & $\eps_k$ & $\phi(\eps_k)$ & $\#$ iters (E+C) \\
\hline
\rule{0pt}{9pt}
$0$ &
$3.200385406490764$ &
$1.36461$ &
-- \\
$1$ &
$2.005376560653111$ &
$2.25352\,10^{-2}$ &
$31+5$ \\
$2$ &
$1.985895420074233$ &
$1.01676\,10^{-5}$ &
$27+3$ \\
$3$ &
$1.985886631877472$ &
$2.10987\,10^{-12}$ &
$10+3$ \\
$4$ &
$1.985886631875651$ &
$3.10862\,10^{-15}$ &
$3+4$ \\
\hline
\end{tabular}
\vspace{2mm}
\caption{Computation of the stability radius for the matrix
$A-4I$ (with $A$ the matrix in (\ref{chap:proto}.\ref{eq:example}))
by the HEC algorithm: computed values $\eps_k$,
$\phi(\eps_k)=-\Re\,\lambda(A+\eps_k E_k)$,
and the numbers of iterations performed in the expansion 
and contraction (Newton--bisection) phases.}
\label{tab:HurHEC}
\end{center}
\end{table}

The HEC iteration converges to
\[
\eps_\ast = 1.985886631875649,
\]
which agrees with the value computed by the nested iteration.
Table \ref{tab:HurHEC} has been obtained by applying Tim Mitchell's code Rostapack which implements the HEC algorithm (http://www.timmitchell.com/software/ROSTAPACK/).
\ecltwo

\section{Notes}

An early survey of matrix nearness problems, with emphasis on the properties of
symmetry, positive definiteness, orthogonality, normality, rank-deficiency and instability, was given by Higham~(\cite{Hig89}). His review is a source of continuing interest in view of the choice of topics and the references to the older literature. 

\subsubsection*{Distance to instability (stability radius) under complex unstructured perturbations.}
\bng
Van Loan (\cite{VL85}) and
Hinrichsen \& Pritchard (\cite{HinP86a}) appear to be first to address the problem how to compute a nearest unstable matrix to a given stable matrix.
\eng
Van Loan considered both complex and real perturbations and came up with heuristic algorithms for approximating the smallest perturbations that shift an eigenvalue to the imaginary axis. His starting point was the characterization of the distance to instability under complex unstructured perturbations of the matrix $A$ as 
\begin{equation}\label{stab-rad}
\beta(A) = \min_{\omega\in\R}\, \sigma_{\min}(A-\iu\omega I)
\end{equation}
and an intricate characterization of the distance to instability under real perturbations.

For the {\it complex} case, Byers (\cite{Bye88}) showed that the Hamiltonian matrix 
$$
H(\sigma)=\begin{pmatrix}  A  & -\sigma I \\ \sigma I & A^* \end{pmatrix} 
$$
has a purely imaginary eigenvalue if and only if $\sigma \ge \beta(A)$; cf.~Lemma~\ref{chap:pseudo}.\ref{lem:byers}. Based on this result, he proposed a bisection method for computing the distance to the nearest complex matrix with an eigenvalue on the imaginary axis (the complex stability radius). Each step of the method requires the solution of an eigenvalue problem of the Hamiltonian matrix $H(\sigma)$ for varying $\sigma> 0$.
Byers  (\cite{Bye88}) also gave an extension of the algorithm to compute the distance to the nearest complex matrix with an eigenvalue on the unit circle. 

Conceptually related Hamiltonian eigenvalue methods by Boyd \& Balakrishnan (\cite{BoyB90}) and Bruinsma \& Steinbuch (\cite{BruS90}) for the more general problem of computing the $\mathcal{H}_\infty$-norm of a transfer function also apply to computing the distance to stability. These methods converge locally quadratically.

He \& Watson (\cite{HeW99}) developed a method for computing the distance to instability that is better suited for large sparse matrices~$A$.
They use a method based on inverse iteration for singular values to compute a stationary point of the function
$f(\omega)= \sigma_{\min}(A-\iu\omega I)$. They then check whether the stationary point reached is a global minimum by solving an eigenvalue problem for $H(\sigma)$.  An alternative method for large sparse matrices was devised by Kressner (\cite{Kre06}) who worked with inverse iterations using sparse LU factorizations of imaginary shifts of Hamiltonian matrices $H(\sigma)$.

For $\sigma=\beta(A)$, the Hamiltonian matrix $H(\sigma)$ has an eigenvalue of even multiplicity on the imaginary axis. Generically, it is expected to be a defective double eigenvalue. Freitag \& Spence (\cite{FreS11}) used a Newton-based method 
to find the parameters $\sigma$ and $\omega$ such that $H(\sigma) - \iu\omega I$ has a zero eigenvalue corresponding to a two-dimensional Jordan block.

\subsubsection*{Two-level iteration.}
A different approach is to combine an algorithm for computing the $\eps$-pseudospectral abscissa $\alpha_\eps(A)$ 
(see Section~\ref{sec:psa}) with a root-finding algorithm such as a Newton--bisection method for determining $\oeps>0$ such that $\alpha_{\oeps}(A)=0$. Then,
$\oeps$ is the distance to instability. \bcl Such a two-level approach was apparently first proposed by Guglielmi, G\"urb\"uzbalaban \& Overton (\cite{GugGO13}) in the closely related, more general setting of computing the $\mathcal{H}_\infty$-norm of a linear time-invariant system; see also Section~\ref{sec:Hinf}. The two-level approach
with inner and outer iterations can be used efficiently to approximate the distance to instability \ecl for large sparse matrices using the rank-1 iteration of Guglielmi \& Overton (\cite{GO11}), or the subspace method of Kressner \& Vandereycken (\cite{KV14}),
or with the discretized rank-1 differential equation of Guglielmi \& Lubich (\cite{GL11}) in the inner iteration. The latter, differential equation based approach is described here in Section~\ref{sec:two-level}. It extends in a direct way to computing the distance to instability under real or structured perturbations, as will be described in the next chapter.

\bng
The initialization \eqref{eq:eps0F} was proposed in the articles by Mitchell \& Overton (\cite{MitO16}) and, in the case of a real structure, by Guglielmi, G\"urb\"uzbalaban, Mitchell \& Overton (\cite{GGMO17}). 
\eng

\bcl
\subsubsection*{HEC method.}
\index{Hybrid expansion-contraction method}
\index{HEC method}
Possible failure of convergence of the standard Newton--bisection method in the outer iteration, in cases where the inner iterations provide local minima that switch paths, was first observed and discussed by Mitchell \& Overton (\cite{MitO16}). As a remedy, they proposed the hybrid expansion--contraction (HEC) method and showed that it has guaranteed monotone convergence to an upper bound of the distance to instability under very weak assumptions and has quadratic convergence under suitable regularity assumptions. The HEC method was later generalized and studied in a more abstract framework by Mitchell \& Van Dooren (\cite{MitVD23}), which extended and clarified the scope of the approach. Instead of optimizing to convergence for a fixed perturbation size in the inner iteration and using the result in one Newton--bisection step in the outer iteration, the HEC method alternates between optimizing to convergence and root-finding to convergence for a fixed perturbation size and fixed perturbation matrix, respectively. The algorithmic building blocks are the same in both methods, but they are arranged differently.

The monotone variant of the Newton--bisection method given in Subsection~\ref{subsec:Newton-bisection-monotone} is another remedy to the possible failure of convergence of the standard Newton--bisection method.
The method and its convergence analysis have not previously appeared in the literature.
\ecl
\chapter{Real and structured perturbations}
\label{chap:struc}

In this chapter we extend the programme and algorithms of the previous chapters to the case where the matrix perturbations are restricted to be real or, more generally, to lie in a structure  space $\cS$, which is a given complex-linear or real-linear subspace of $\C^{n,n}$. For example, $\cS$  may be a space of complex or real matrices with a given sparsity pattern, or of matrices with given range and co-range, or special matrices such as Hamiltonian or Toeplitz or Sylvester matrices. Remarkably, the rank-1 property of optimizers extends to general linear structures: optimizers of the structured eigenvalue optimization problem \bng -- although often of full rank -- \eng
are the orthogonal projection of complex rank-1 matrices onto the structure space. This still enables us to work with rank-1 matrices to compute quantities related to structured pseudospectra such as the structured $\eps$-pseudospectral abscissa,  to trace the boundary of structured pseudospectra, and to adapt 
the two-level approach of the previous chapter to solve structured matrix nearness problems. In particular, we consider determining the following:
\begin{itemize}
    \item the {\em $\cS$-structured stability radius} (distance to instability) of a Hurwitz matrix~$A$, i.e. finding a structured perturbation $\Delta\in\cS$ of minimal norm such that $A+\Delta$ has some eigenvalue on the imaginary axis.
    This quantity assesses the robustness of asymptotic stability of linear differential equations under structured perturbations of the matrix; 
    \index{structured stability radius} 
    \index{stability radius!structured}
    \index{distance to instability!under structured perturbations}
\item the {\em $\cS$-structured $\eps$-stability radius} of a Hurwitz matrix~$A$, i.e. finding a structured perturbation $\Delta\in\cS$ of minimal norm such that the (complex unstructured) pseudospectrum $\Lambda_\eps(A+\Delta)$ has some point on the imaginary axis. This quantity will be shown to assess the robustness of transient bounds of linear differential equations under structured perturbations of the matrix;
\index{structured $\eps$-stability radius}
\item the {\em $\cS$-structured dissipativity radius} of a matrix $A$ that has the numerical range in the open complex left half-plane. The aim is to find a structured perturbation $\Delta\in\cS$ of minimal norm such that the numerical range of $A+\Delta$ has some point on the imaginary axis.
\index{structured dissipativity radius}
The norm of $\Delta$  assesses the robustness of dissipativity $\| \exp(t(A+\Delta)) \|_2 \le 1$ under structured perturbations $\Delta\in\cS$;
\item the {\em $\cS$-structured distance to singularity} of an invertible matrix, i.e. finding a structured perturbation $\Delta\in\cS$ of minimal norm such that $A+\Delta$ is singular. This is illustrated by the problem of finding a nearest pair of polynomials with a nontrivial common divisor to a given pair of polynomials that is coprime, which can be restated as a matrix nearness problem for Sylvester matrices.
\index{structured distance to singularity}
\index{distance to singularity!under structured perturbations}
\end{itemize}

For the first two items, we give two-level algorithms that use rank-1 matrix differential equations in the inner iteration and compute eigenvalues and eigenvectors of perturbed matrices in every step.

For the last two items, we present two-level algorithms with iterations that {\it do not} compute eigenvalues and eigenvectors (or singular values and singular vectors) but which still use rank-1 matrix differential equations in the inner iteration. In the third item we maximize the real part of Rayleigh quotients of perturbed matrices, and in the last item we minimize the squared Euclidean norm of perturbed matrices multiplying vectors. The resulting algorithms only need to compute matrix--vector products and inner products of vectors.

\bcl This chapter shows the versatility of the two-level approach that uses rank-constrained gradient systems for the perturbation matrix 
 and Newton--bisection methods for the perturbation size in a nested or alternating way. For ease of presentation we will here only formulate simple inner--outer iterations, but it is understood that monotonically increasing Newton--bisection iterations and monotonically decreasing HEC iterations  as discussed in Chapter~\ref{chap:two-level} can be used instead (unless noted otherwise) and may be favourable.\ecl

\section{Real version of the eigenvalue optimization problem}
\label{sec:real} 
\index{eigenvalue optimization!real}

\subsection{Problem formulation}

We now consider Problem (\ref{chap:proto}.\ref{eq:optimiz0}) for a {\it real} matrix $A\in\R^{n,n}$ and {\it real} perturbations $\Delta\in \R^{n,n}$: find
\begin{equation} \label{eq:optimiz0rF}
\arg\min\limits_{\Delta \in \R^{n,n}, \| \Delta \|_F = \eps} f \left( \lambda\left( A + \Delta \right), \clambda \left( A + \Delta \right)  \right),
\end{equation}
where again $\lambda(A+\Delta)$ is the target eigenvalue of the perturbed matrix $A + \Delta$, and $f$ satisfies (\ref{chap:proto}.\ref{ass:f}).
As in the complex case, it is convenient to write
\[
\Delta = \eps E \quad \mbox{ with } \ \| E \|_F = 1
\]
and to use the notation
\begin{equation}
\F_\eps(E) = f \left( \lambda\left( A + \eps E \right), \conj\lambda\left( A + \eps E \right) \right)
\end{equation}
so that \eqref{eq:optimiz0rF} can be rewritten as
\begin{equation} \label{eq:optimizrF}
\arg\min\limits_{E \in \R^{n, n}, \| E \|_F = 1} \F_\eps(E).
\end{equation}

\subsection{Norm-constrained gradient flow and rank of optimizers}
\label{subsec:gradient-flow-real} 

The programme of Section~\ref{sec:proto-complex} extends to the real case with minor but important modifications.

\medskip\noindent 
{\bf Real gradient.}
\index{gradient!free}
Consider a smooth path of {\it real} matrices $E(t)\in\R^{n,n}$. Since $\dot E(t)$ is then also a real matrix, we have by Lemma~\ref{lem:gradient}
\begin{equation} \label{eq:deriv-real}
\frac1{ \eps \kappa(t) } \,\frac{d}{dt} \F_\eps(E(t)) = \bigl\langle  G_\eps^\R(E(t)),  \dot E(t) \bigr\rangle
\end{equation}
with the rescaled real gradient
\begin{equation}\label{gradient-real}
G_\eps^\R(E) := \Re\, G_\eps(E) = \Re( 2 f_{\clambda} \,xy^*) \in \R^{n,n},
\end{equation}
which is the real part of a complex rank-1 matrix and hence has rank at most $2$ (as a sum of two rank-1 matrices). As $\eps$ is fixed and only the real case is considered in this section, we often write for short
$$
G(E) := G_\eps^\R(E) .
$$
Lemma~\ref{lem:opt} on the direction of steepest admissible descent extends without ado to the real case: Consider real matrices and take everywhere the real inner product instead of the real part of the complex inner product.

\medskip\noindent
{\bf Norm-constrained gradient flow.}
\index{gradient flow!norm-constrained}
We consider the  gradient flow on the manifold of {\it real}  $n\times n$ matrices of unit Frobenius norm,
\begin{equation}\label{ode-ErF}
\dot E = -G(E) + \langle G(E), E \rangle E.
\end{equation}
{\bf Monotonicity.} 
Assuming simple eigenvalues along the trajectory,
we then still have the monotonicity property of Theorem~\ref{thm:monotone},
\begin{equation}
\frac{d}{dt} \F_\eps (E(t)) =- \| G(E) - \langle G(E), E \rangle E \|_F^2 \le  0,
\label{eq:pos-real}
\end{equation}
with essentially the same proof (the real inner product replaces the real part of the complex inner product). 

\medskip\noindent
{\bf Stationary points.}
\index{stationary point}
Also the characterization of stationary points as given in Theorem~\ref{thm:stat} extends with the same proof: Let
$E\in\R^{n,n}$ with $\| E\|_F=1$ be such that the eigenvalue $\lambda(A+\eps E)$ is simple
and $G_\eps^\R(E)\ne 0$. Then, 
\begin{equation}\label{stat-real}
\begin{aligned}
&\text{$E$ is a stationary point of the differential equation \eqref{ode-ErF}}
\\[-1mm]
&\text{if and only if $E$ is a real multiple of $G_\eps^\R(E)$.}
\end{aligned}
\end{equation}
As a consequence, optimizers of \eqref{eq:optimizrF} have rank at most $2$. We can determine the precise rank as follows.

\begin{theorem} [Rank of optimizers] \label{thm:rank}
For $A\in\R^{n, n}$ and $\eps>0$, let $E\in\R^{n,n}$ with $\| E\|_F=1$  be a stationary point of the differential equation \eqref{ode-ErF}
such that the eigenvalue $\lambda=\lambda(A+\eps E)$ is simple and $G_\eps^\R(E)\ne 0$. Then, $E$ is the real part of a complex rank-1 matrix and the rank of $E$ is as follows.
\begin{itemize}
\item[(a)] $\quad$If\/ $\lambda$ is real, then $E$ has rank $1$.
\item[(b)] $\quad$If\/ $\Im\,\lambda \ne 0$, then $E$ has rank $2$.
\end{itemize}
\end{theorem}
\index{optimizer!rank-2}

\begin{proof}
%
(a) If $\lambda$ is real, then $f_\lambda$ is real and $\y$ and $\x$ can be chosen real,  hence $G(E)=\Re(\overline f_\lambda \y \x^*)=f_\lambda \y \x^\top$ is of 
rank~1, and so is every nonzero real multiple, in particular $E$.

(b) We set $\xx=f_\lambda \x$ and separate the real and imaginary parts in $\y=\y_R+\iu \y_I$ and $\xx=\xx_R+\iu \xx_I$. 
If $\Re(\y \xx^*)=\y_R \xx_R^\top + \y_I \xx_I^\top$ is of rank 1, 
then $\y_R$ and $\y_I$ are linearly dependent or $\xx_R$ and $\xx_I$ are linearly dependent. 
Let us first assume the former. In this case there is a real $\alpha$ such that $\y=\cos(\alpha) v+ \iu \sin(\alpha) v$ for some nonzero real vector $v$. 
Rotating both $\y$ and $\x$ by ${\rm e}^{-\iu \alpha}$ does not change the required property $\y^*\x>0$, so we can assume without loss of generality that $\y$ is a real 
left eigenvector of the real matrix $A+\eps E$, which implies that the corresponding eigenvalue $\lambda$ is real. The argument is analogous when $\xx_R$ and $\xx_I$ 
are linearly dependent.
\qed \end{proof}

\subsection{Low-rank matrices and their tangent matrices}
\label{subsec:low-rank}
\index{rank-$r$ matrix}

Theorem~\ref{thm:rank} motivates us to search for a differential equation on the manifold of real rank-2 matrices
that has the same stationary points with non-real target eigenvalue as \eqref{ode-ErF}. For the stationary points with real target eigenvalues we use a differential equation for rank-1 matrices as in Section~\ref{subsec:rank1-gradient-flow}.

In the following we consider differential equations on the manifold of real $n\times n$ matrices of rank $r$, denoted 
$$
\cM_r = \cM_r(\R^{n,n})= \{ E \in \R^{n,n}\,:\, \text{rank}(E) = r\}.
$$
While only ranks 1 and 2 are of interest for the optimization problem studied in this section, we now consider the case of a general fixed rank $r$, since it is not more complicated than rank 2 and will be useful later. We proceed similarly to the rank-1 case considered in Section~\ref{subsec:rank-1}.

%
%

Every real rank-$r$ matrix $E$ of dimension $n\times n$ can be written in
the form
\begin{equation}\label{USV}
E = USV^\top
\end{equation}
where $U\in\R^{n, r}$ and $V\in\R^{n, r}$ have orthonormal
columns, i.e.,
\begin{equation}\label{UV-orth}
U^\top  U = I_r, \quad\ V^\top  V =I_r,
\end{equation}
(with the identity matrix $I_r$ of dimension $r$),
and $S\in \R^{r,r}$ is nonsingular. The singular value decomposition yields
$S$ diagonal, but here we will not assume a special form of $S$.
The representation (\ref{USV}) is not unique: replacing $U$ by
$\widetilde U=UP$ and $V$ by $\widetilde V=VQ$
with orthogonal matrices $P,Q\in\R^{r,r}$,
and correspondingly $S$ by $\widetilde S=P^\top SQ$, yields the same matrix
$E=USV^\top =\widetilde U \widetilde S \widetilde V^\top $.

Every tangent matrix $\dot  E \in T_E\cM_r$ is of the form
\begin{equation}\label{Edot-r}
 \dot  E = \dot  U S V^\top  + U\dot  S V^\top  + U S \dot  V^\top ,
\end{equation}
where  $\dot  S\in\R^{r,r}$,  and $U^\top \dot  U$ and
$V^\top  \dot  V$ are skew-symmetric (as results from differentiating $U^\top U=I_r$ and $V^\top V=I_r$).
The matrices
$\dot  S, \dot  U, \dot  V$ are uniquely determined
by $\dot  E$ and $U,S,V$ if we impose the orthogonality conditions
\begin{equation}\label{orth}
U^\top \dot  U = 0, \quad\ V^\top  \dot  V=0.
\end{equation}
Multiplying $\dot E$ with $U^\top$ from the left and with $V$ from the right, we then obtain 
$$
\dot S = U^\top \dot E V, \quad \dot U S = \dot E V - U \dot S, \quad
S {\dot V}^\top = U^\top \dot E - \dot S V^\top,
$$ 
which yields $\dot S,\dot U,\dot V$ in terms of $\dot E$.
Extending this construction, we arrive at an explicit formula for the orthogonal projection onto the tangent space. Here, orthogonality refers to the real Frobenius inner product.

\begin{lemma}[Rank-$r$ tangent space projection]
\label{lem:P-formula-r}
\index{rank-$r$ matrix!tangent space projection}
The orthogonal projection from $\R^{n, n}$ onto the tangent space $T_E\cM_r$ at $E=USV^\top \in\cM_r$
is given by
\begin{equation}\label{P-formula-r}
P_E(Z) = Z - (I-UU^\top) Z (I-VV^\top)
\quad\text{ for $Z\in\C^{n, n}$}.
\end{equation}
\end{lemma}

\begin{proof}  The proof is a direct extension of the proof of Lemma~\ref{lem:P-formula-1}.
Let $P_E(Z)$ be defined by \eqref{P-formula-r}. Determining (similarly to above) $\dot S,\dot U,\dot V$ by
\begin{equation*}
\dot S = U^\top Z V, \quad \dot U S = Z V - U \dot S, \quad
S {\dot V}^\top = U^\top Z - \dot S V^\top,
\end{equation*} 
we obtain $P_E(Z)=UU^\top ZVV^\top - ZVV^\top - UU^\top Z$ in the form \eqref{Edot-r} with \eqref{orth} and hence
$$
P_E(Z) \in T_E\cM_r.
$$
Furthermore,
$$
\langle P_E(Z),\dot E \rangle =   \langle Z, \dot E \rangle \qquad\text{for all }\ \dot E \in T_E\cM_r,
$$
because $\langle (I-UU^\top) Z (I-VV^\top),\dot E \rangle = \langle Z, (I-UU^\top) \dot E (I-VV^\top) \rangle = 0$ by \eqref{Edot-r}. Hence,
$P_E(Z)$ is  the orthogonal projection of $Z$ onto $T_E\cM_r$.
\qed \end{proof}
We note that $P_E(E)=E$ for $E\in\cM_r$, or equivalently, $E\in T_E \cM_r$.

\subsection{Rank-constrained gradient flow}
\label{subsec:rank-r-gradient-flow}
\index{gradient flow!rank-$r$ constrained}

In the differential equation (\ref{ode-ErF}) we replace the right-hand side by its orthogonal projection onto $T_E\cM_r$, so that solutions starting with rank $r$ 
will retain rank $r$ for all times:
\begin{equation} \label{ode-ErF-2}
\dot E = P_E \Bigl( -G(E) + \langle G(E), E \rangle E \Bigr).
\end{equation}

Since $E\in T_E\cM_r$, we have $P_E(E)=E$ and $\langle E,Z \rangle = \langle E,P_E(Z) \rangle$, and hence the differential equation can be rewritten as
\begin{equation}\label{ode-ErF-2-v2}
\dot E = -P_E(G(E) ) + \langle E,  P_E(G(E) )\rangle E,
\end{equation}
which differs from (\ref{ode-ErF}) only in that $G(E)$ is replaced by its orthogonal projection to $T_E\cM_r$. 
This shows that $\langle E,\dot E \rangle=0$, so that the unit norm of $E$ is conserved along solutions of (\ref{ode-ErF-2}).

To obtain the differential equation in a form that uses the factors in $E=USV^\top $ rather than the full $n\times n$ matrix $E$, we use the following result, which follows directly from the proof of Lemma~\ref{lem:P-formula-r}.

\begin{lemma}[Differential equations for the factors]
\label{lem:USV}  
For $E=USV^\top \in \cM_r$ with nonsingular $S\in\R^{r, r}$ and with
$U\in\R^{n, r}$ and $V\in\R^{n, r}$ having orthonormal
columns, the equation $\dot E=P_E(Z)$ is equivalent to
$
\dot E = \dot U SV^\top  + U \dot S V^\top  + US\dot V^\top ,
$
where
\begin{eqnarray}
\dot S &=& U^\top  Z V
\nonumber\\
\dot U &=& (I-UU^\top) Z V S^{-1}
\label{odesrFr} \\
\dot V &=& (I-VV^\top) Z^\top  U S^{-\top}.
\nonumber
\end{eqnarray}
\end{lemma}

\noindent
With $Z=-G(E) + \langle G(E), E \rangle E$ and $r=2$, this yields that the differential equation (\ref{ode-ErF-2})  for 
$E=USV^\top$ is equivalent to a system of differential equations for $S,U,V$. On the right-hand side appears the inverse of the matrix $S$, which may be ill-conditioned. In the present context, this appears when the target eigenvalue gets close to the real axis (see
Theorem~\ref{thm:rank}) so that $E$, and hence $S$, becomes almost of rank $1$. In such a situation of a small singular value in $S$, standard numerical integrators become unstable or yield plainly wrong numerical solutions unless used with a tiny stepsize proportional to the smallest nonzero singular value. Later in this section we will describe a numerical integrator that is robust to small singular values.


\subsubsection*{Monotonicity.}
Assuming simple eigenvalues almost everywhere along the trajectory,
we still have the monotonicity property of Theorem~\ref{thm:monotone-C-1} along solutions $E(t)$ of~\eqref{ode-ErF-2},
\begin{equation}
\frac{d}{dt} \F_\eps (E(t)) = - \| P_E(G(E) ) - \langle E,  P_E(G(E) )\rangle E \|_F^2 \le  0,
\label{eq:pos-real-r}
\end{equation}
with essentially the same proof (the real inner product replaces the real part of the complex inner product).
 
\subsubsection*{Stationary points.}
Comparing the differential equations \eqref{ode-ErF} and \eqref{ode-ErF-2} immediately shows that every stationary point of \eqref{ode-E} is also a stationary point of the projected differential equation \eqref{ode-ErF-2}. As in Theorem~\ref{thm:stat-1}, the converse is also true for the stationary points $E$ of unit Frobenius norm with ${P_E(G(E))\ne 0}$. This shows that the low-rank projection does not create spurious stationary points.

\begin{theorem}[Stationary points]
\label{thm:stat-r}
\index{stationary point}
Let $E\in \cM_2$ be of unit Frobenius norm and assume that $P_E(G(E))\ne 0$. If $E$ is a stationary point of the projected differential equation \eqref{ode-ErF-2}, then $E$ is  already a stationary point of the differential equation \eqref{ode-ErF}.
\end{theorem}

\begin{proof} The proof extends the proof of Theorem~\ref{thm:stat-1}.
We show that $E$ is a real multiple of $G_\eps(E)$. By \eqref{stat-real}, $E$ is then a stationary point of the differential equation \eqref{ode-E}.

For a stationary point $E$ of \eqref{ode-ErF-2}, we must have equality in \eqref{c-s-1}, which shows that $P_E(G)$  (with $G=G(E)$ for short) is a nonzero real multiple of $E$. Hence, in view of $P_E(E)=E$, we can write $G$ as
$$
G=\mu E + W, \quad\text{ where $\mu\ne 0$ is real and $P_E(W)=0$.}
$$
With $E=USV^\top$ as above, we then have
 $$
 W=W-P_E(W)= (I-UU^\top)W(I-VV^\top).
 $$
 Since $G$ is of rank at most $2$, it can be written in the form $G=XRY^\top$, where $X,Y\in \R^{n,2}$ have orthonormal columns and $R\in \R^{2,2}$.
So we have
 $$
XRY^\top = \mu USV^\top + (I-UU^\top)W(I-VV^\top).
 $$
 Multiplying from the right with $V$ yields $X(RY^\top V) = \mu US$, which shows 
 that $X$ has the same range  as $U$, and multiplying from the left with $U^\top$ yields that $Y$ has the same range as $V$. Hence, $G$ has the same range and corange as $E$, which implies that $P_E(G)=G$. Since we already know that $P_E(G)$ is a nonzero real multiple of $P_E(E)=E$, it follows that $G$ is the same real multiple of $E$. By \eqref{stat-real}, $E$ is therefore a stationary point of  \eqref{ode-ErF}.
\qed \end{proof}
As in Remark~\ref{rem:exceptional}, it is shown that if $G(E)$ is of rank 2 and $P_E(G(E))=0$, then it follows that $Ey=0$ and $x^*E=0$,
which implies that $\lambda$ is already an eigenvalue of the unperturbed matrix $A$ with the same eigenvectors $x$ and $y$, which is an exceptional situation.

\subsection{Time-stepping for the low-rank differential equation}
\label{subsec:low-rank-integrator}
\subsubsection*{A robust integrator.}
The following method adapts the low-rank integrator of Ceruti \& Lubich (\cite{CeL22}) to the norm-constrained situation considered here. It first updates the basis matrices $U$ and $V$ with orthonormal columns in parallel and then uses a Galerkin approximation to the differential equation \eqref{ode-ErF-2} in the updated basis. This basis update and Galerkin (BUG) integrator has been shown to be robust to the presence of small singular values, which would here appear in the case of a target eigenvalue near the real axis; \bcl cf.~Theorem~\ref{thm:rank}. \ecl
\index{BUG low-rank integrator}

	One time step of integration from time $t_k$ to $t_{k+1}=t_k+h$  starting from a factored rank-$r$ matrix 
	$E_k=U_kS_kV_k^\top$ of unit Frobenius norm computes an updated rank-$r$ factorization $E_{k+1}=U_{k+1}S_{k+1}V_{k+1}^\top$ of unit Frobenius norm as follows.
	
	\begin{enumerate}
		\item 
		Update the basis matrices $ U_k \rightarrow U_{k+1}$ and $ V_k \rightarrow V_{k+1}$:
		\\[1mm]
		Integrate from $t=t_k$ to $t_{k+1}=t_k+h$ the $n \times r$ matrix differential equation
		$$ \dot{K}(t) = -G( K(t) V_k^\top) V_k  , \qquad K(t_k) = U_k S_k.$$
		Perform a QR factorization $K(t_{k+1}) = U_{k+1} {R}_{k+1}$ and compute the $r\times r$ matrix $M= U_{k+1}^\top U_k$.
				\\[2mm]
		Integrate from $t=t_k$ to $t_{k+1}$ the $n \times r$ matrix differential equation
		$$  \dot{L}(t) =-G(U_k L(t)^\top)^\top  U_k, \qquad L(t_k) = V_k {S}_k^\top. $$
		Perform a QR factorization $L(t_{k+1}) = V_{k+1} \widetilde{R}_{k+1}$ and compute the $r\times r$ matrix $N= V_{k+1}^\top V_k$.\\[-2mm]
		
		\item
		Update ${S}_k \rightarrow {S}_{k+1}$\,: \\[1mm]
		Integrate from $t=t_k$ to $t_{k+1}$ the $r \times r$ matrix differential equation
		$$ \dot{S}(t) =  -U_{k+1}^\top G(U_{k+1} S(t) V_{k+1}^\top) V_{k+1}, 
		\qquad S(t_k) = \frac{M {S}_k N^\top}{ \| M {S}_k N^\top \|_F},
		$$
		and set ${S}_{k+1} =S(t_{k+1})/\| S(t_{k+1}) \|_F$.
	\end{enumerate} 
The differential equations in the substeps are solved approximately by a step of some standard numerical integrator, e.g.~the explicit Euler method or better a low-order explicit Runge--Kutta method such as the second-order Runge method. We denote the result of the fully discrete method with stepsize $h$ as $U(h)$, $V(h)$, $S(h)$.

\subsubsection*{Step-size selection.}
\index{step-size selection}
As in Section~\ref{subsec:proto-numer}, the stepsize selection is done by an Armijo-type strategy. Along solutions $E(t)\in \mathcal{M}_2$ of \eqref{ode-ErF-2}, we have by \eqref{eq:deriv-real}--\eqref{gradient-real}
$$
\frac{d}{dt} \F_\eps(E(t)) = - g(E(t)) \quad\text{ with } \quad 
g(E)= \eps \kappa \bigl( \| P_E(G)\|_F^2 - \langle G,E \rangle^2 \bigr) \ge 0
$$
where $\kappa=1/(x^*y)>0$, $G=G_\eps^\R(E)=\Re(2f_{\clambda} xy^*)$ with the normalized left and right eigenvectors $x$ and $y$ to the eigenvalue $\lambda(A+\eps E)$.

We note that on separating real and imaginary parts in $2f_{\clambda} x=w_R +\iu w_I$ and $y=y_R+\iu y_I$ and defining the $n\times 2$ real matrices 
$W=(w_R,w_I)$ and $Y=(y_R,y_I)$, we have the real factorization
$$
G=WY^\top.
$$
With the rank-2 matrix $E=USV^\top$ in factorized form as above, we can then compute $g(E)$ without actually forming the $n\times n$ matrices $E$ and $G$: noting that with the $2\times 2$ matrices $P=U^\top W$ and $Q=V^\top Y$ we have
$$
\langle G,E \rangle = \langle WY^\top, USV^\top\rangle = \langle PQ^\top,S \rangle_{\R^{2\times 2}}
$$
and 
$$
P_E(G)= UPY^\top - UPQ^\top V^\top + WQ^\top V^\top,
$$
which yields after a straightforward computation
$$
\| P_E(G) \|_F^2 = \| PY^\top \|_F^2 + \|WQ^\top \|_F^2 - \| PQ^\top \|_F^2,
$$
so that finally $g=g(E)$ is given by
\begin{equation}\label{g-n-formula-r}
g= \eps \kappa \bigl( \| PY^\top \|_F^2 + \|WQ^\top \|_F^2 - \| PQ^\top \|_F^2 -  \langle PQ^\top,S \rangle^2 \bigr).
\end{equation}
With this quantity $g$, the Armijo-type stepsize selection is then done as in Section~\ref{subsec:proto-numer}. A time step of the method is summarized in Algorithm~\ref{alg:real}.

\medskip
\begin{algorithm}[H] \label{alg:real}
\DontPrintSemicolon
\KwData{$A, \eps, \theta > 1, U_k \approx U(t_k),  V_k \approx V(t_k)\in \R^{n\times 2}$ with orthonormal columns, $S_k \approx S(t_k)\in \R^{2\times 2}$ of unit Frobenius norm, target eigenvalue $\lambda_k=\lambda(A+\eps U_kS_kV_k^\top)$, $h_{k}$ (proposed step size)}
\KwResult{$U_{k+1}, V_{k+1},S_{k+1}$, $\lambda_{k+1}$, $h_{k+1}$}
\Begin{
\nl Initialize the step size by the proposed step size, $h=h_{k}$\; 
\nl Compute $f_k=f(\lambda_k,\conj{\lambda_k})$ \;
\nl Compute the left/right eigenvectors $x_k, y_k$ to $\lambda_k$ such that $\| x_k \| = \| y_k \| = 1, x_k^* y_k > 0$\;
\nl Compute $g_k$ by \eqref{g-n-formula-r}\;
\nl Initialize $f(h) = f_k$\;
\While{$f(h) \ge f_k$}{
\nl Compute $U(h), V(h), S(h)$ by the above rank-2 integrator\;
\nl Compute $\lambda(h)$ target eigenvalue of $A + \eps U(h) S(h)V(h)^\top$\; 
\nl Compute $f(h) = f\bigl( \lambda(h), \conj{\lambda(h)} \bigr)$\;
\If{$f(h) \ge f_k$}{Reduce the step size, $h:=h/\theta$}
}
\nl Initialize $h_{\rm next}=h$\;
\If{$f(h) \ge f_k- (h/\theta) g_k$}{Reduce the step size for the next step, $h_{\rm next}:=h/\theta$}
\If{$h_{\rm next}=h_{k}$}
{\nl Compute $U(\theta h), V(\theta h), S(\theta h)$ by the above rank-2 integrator\;
\nl Compute $\lambda(\theta h)$ target eigenvalue of $A + \eps U(\theta h) S(\theta h)V(\theta h)^\top$\; 
\nl Compute $f(\theta h) = f\bigl(\lambda(\theta h),\conj{\lambda(\theta h)} \bigr)$\;
\If{$f(h) > f(\theta h)$}
{Enlarge the step size for the next step, $h:= \theta h$ and then $h_{\rm next} := h$}
}
\nl Set $h_{k+1}=h_{\rm next}$, $\lambda_{k+1}= \lambda(h)$, and the starting values for the next step as 
$U_{k+1}=U(h)$, $V_{k+1}=V(h)$, $S_{k+1}=S(h)$\;
\Return
}
\caption{Integration step for the rank-2 constrained real gradient system}
\end{algorithm}

\section{Structured eigenvalue optimization problems}
\label{sec:proto-structured}
\index{eigenvalue optimization!structured}


\subsection{Problem formulation. Complex- and real-linear structures}
\index{structure space}
Let the structure space $\cS$ be a complex-linear or real-linear subspace of $\C^{n,n}$, e.g. a space of complex or real matrices with a prescribed sparsity pattern,  or with prescriped range and co-range, or Toeplitz matrices, or Hamiltonian matrices, etc.  

As before, we set with the target eigenvalue $\lambda(M)$
\begin{equation} \label{eq:optimiz0S}
\F_\eps(E) = f \left( \lambda\left( A + \eps E \right), \conj\lambda\left( A + \eps E \right) \right).
\end{equation}
We now restrict the admissible perturbations $\eps E$ to be in $\cS$ and consider the {\it structured} eigenvalue optimization problem to find
\begin{equation} \label{eq:optimizS}
\arg\min\limits_{E \in \cS, \| E \|_F = 1} \F_\eps(E).
\end{equation}

\subsection{Orthogonal projection onto the structure}\label{subsec:proj-structure}
\index{structure space!orthogonal projection}
Let $\Pi^\cS$ be the orthogonal projection (w.r.t. the Frobenius inner product)  onto~$\cS$\/: for every $Z\in \C^{n,n}$, 
\begin{equation}\label{Pi-S}
\Pi^\cS Z \in \cS \quad\text{ and } \quad \Re\langle \Pi^\cS Z, W \rangle = \Re\langle Z,W \rangle \quad \text{for all }\, W\!\in\cS.
\end{equation}
For a complex-linear subspace $\cS$, taking the real part of the complex inner product can be omitted (because with $W\in\cS$, then also $\iu W\in\cS$), but taking the real part is needed for real-linear subspaces. 
In the following examples, the stated action of $\Pi^\cS$ is readily verified. 

\begin{example}[Real matrices]
For $\cS=\R^{n,n}$,  we have $\Pi^\cS Z=\Re\, Z$ for all $Z\in\C^{n,n}$.
\end{example}

\begin{example}[Sparse matrices] If $\cS$ is the space of complex matrices with a prescribed sparsity pattern, then $\Pi^\cS Z$ leaves the entries of $Z$ on the sparsity pattern unchanged and annihilates those outside the sparsity pattern.
If $\cS$ is the space of real matrices with a prescribed sparsity pattern, then $\Pi^\cS Z$ takes the real part of the entries of $Z$ on the sparsity pattern and annihilates those outside the sparsity pattern.
\end{example}
\index{sparse matrix}

\begin{example} [Matrices with prescribed range and co-range]
An example of particular interest in control theory is the perturbation space
$$
\cS = \{ B \Delta C \,:\, \Delta \in \R^{k,l} \},
$$
where $B\in\R^{n,k}$ and $C\in\R^{l,n}$ with $k,l<n$  are given matrices of full rank. Here, $\Pi^\cS Z = B B^\dagger Z C^\dagger C$,
where $B^\dagger$ and $C^\dagger$ are the Moore--Penrose inverses of $B$ and $C$, respectively. 
\end{example}

\begin{example}[Toeplitz matrices] If $\cS$ is the space of complex $n\times n$ Toeplitz matrices, then $\Pi^\cS Z$ is obtained by replacing in  each diagonal all the entries of $Z$ by their arithmetic mean. For real Toeplitz matrices, the same action is done on $\Re\, Z$. 
\end{example}
\index{Toeplitz matrix}

\begin{example} [Hamiltonian matrices] If $\cS$ is the space of $2d\times 2d$ real Hamiltonian matrices, then $\Pi^\cS Z = J^{-1}\mathrm{Sym}(\Re(JZ))$, where $\mathrm{Sym}(\cdot)$ takes the symmetric part of a matrix and
$$
J=\begin{pmatrix}
0 & I_d \\ -I_d & 0
\end{pmatrix},
$$
for which $J^{-1}=J^\top=-J$. We recall that a real matrix $A$ is Hamiltonian if $JA$ is symmetric.
\end{example}
\index{Hamiltonian matrix}

\subsection{Structure- and norm-constrained gradient flow}
\label{subsec:gradient-flow-S}

The programme of Chapter~\ref{chap:proto} extends to  structured cases as follows.

\subsubsection*{Structured gradient.}
\index{gradient!structured}
Consider a smooth path of {\it structured} matrices $E(t)\in\cS$. Since then also $\dot E(t)\in\cS$, we have by Lemma~\ref{chap:proto}.\ref{lem:gradient}
\begin{equation} \label{eq:deriv-S}
\frac1{ \eps \kappa(t) } \,\frac{d}{dt} \F_\eps(E(t)) = \Re\,\bigl\langle  G_\eps^\cS(E(t)),  \dot E(t) \bigr\rangle
\end{equation}
with the rescaled structured gradient
\begin{equation}\label{gradient-S}
G_\eps^\cS(E) := \Pi^\cS  G_\eps(E) = \Pi^\cS(2 f_{\clambda} \,xy^*) \in \cS,
\end{equation}
where $x,y$ are the left and right eigenvectors, normalized to unit norm and with positive inner product, associated with a simple eigenvalue $\lambda$ of $A+\eps E$, and $f_{\clambda} =(\partial f / \partial\clambda)(\lambda,\clambda)$.

We note that $G_\eps^\cS(E)$ is the orthogonal projection onto $\cS$ of a rank-1 matrix. 

Lemma~\ref{chap:proto}.\ref{lem:opt} on the direction of steepest admissible descent extends immediately to the structured case: If $E,G\in\cS$ in Lemma~\ref{chap:proto}.\ref{lem:opt}, then also $Z_\star$ of (\ref{chap:proto}.\ref{eq:Eopt}) is in $\cS$.

\subsubsection*{Norm- and structure-constrained gradient flow.}
\index{gradient flow!norm-constrained}
\index{gradient flow!structure-constrained}
We consider the  gradient flow on the manifold of {\it structured}  $n\times n$ matrices in $\cS$ of unit Frobenius norm,
\begin{equation}\label{ode-E-S}
\dot E = -G_\eps^\cS(E) + \Re \langle G_\eps^\cS(E), E \rangle E.
\end{equation}

\subsubsection*{Monotonicity.} 
Assuming simple eigenvalues along the trajectory,
we then still have the monotonicity property of Theorem~\ref{thm:monotone},
\begin{equation}
\frac{d}{dt} \F_\eps (E(t))  \le  0,
\label{eq:pos-S}
\end{equation}
with essentially the same proof.

\subsubsection*{Stationary points.}
\index{stationary point}
Also the characterization of stationary points as given in Theorem~\ref{thm:stat} extends with the same proof: Let
$E\in\cS$ with $\| E\|_F=1$ be such that the eigenvalue $\lambda(A+\eps E)$ is simple
and $G_\eps^\cS(E)\ne 0$. Then, 
\begin{equation}\label{stat-S}
\begin{aligned}
&\text{$E$ is a stationary point of the differential equation \eqref{ode-E-S}}
\\[-1mm]
&\text{if and only if $E$ is a real multiple of $G_\eps^\cS(E)$.}
\end{aligned}
\end{equation}

Before we proceed in this direction, we take a look at the condition $G_\eps^\cS(E)\ne 0$. We have the following result.

\begin{theorem} [Non-vanishing structured gradient]
\label{thm:nonzero-gradient-S}
Let   $A,E\in \cS$ and $\eps>0$, and let $\lambda$ be a simple target eigenvalue of $A+\eps E$.

(i) Complex case: $\cS$ is a complex-linear subspace of $\C^{n,n}$. Then, 
$$
G_\eps^\cS(E)\ne 0 \quad\text{ if } \quad \clambda  f_{\clambda} \ne 0.
$$

(ii) Real case: $\cS$ is a real-linear subspace of $\R^{n,n}$. Then, 
$$
G_\eps^\cS(E)\ne 0 \quad\text{ if } \quad \Re(\clambda  f_{\clambda}) \ne 0.
$$
\end{theorem}

\noindent
We emphasize that also $A$ needs to be in $\cS$. The result does not hold true when $A\notin\cS$.

\begin{proof}
    We give the proof for the real case. The complex case is analogous but slightly simpler. We take the real inner product of $G_\eps^\cS(E)$ with $A+\eps E\in \cS$ and use the definition \eqref{gradient-S} of $G_\eps^\cS(E)$: 
    \begin{align*}
    &\langle G_\eps^\cS(E), A+\eps E \rangle = \Re \langle \Pi^\cS (2 f_{\clambda} \,xy^*),A+\eps E \rangle =
    \Re \langle 2 f_{\clambda} \,xy^*,A+\eps E \rangle 
    \\
    &= \Re\bigl( 2 f_{\lambda}\, x^*(A+\eps E)y \bigr) = \Re\bigl( 2 f_{\lambda}\lambda\,x^*y \bigr) 
    = 2\, \Re\bigl( f_{\clambda}\clambda \bigr) \, (x^*y),
    \end{align*}
    where $x^*y>0$. This yields the result.
    \qed
\end{proof}

If the identity matrix $I$ is in $\cS$, then the condition for $G_\eps^\cS(E)\ne 0$ can be weakened: 

-- In the complex case, it then suffices to have $f_{\clambda}\ne 0$. This is seen by taking the inner product with $A+\eps E - \mu I \in \cS$ for an arbitrary $\mu\in \C$.

-- In the real case, if $\lambda$ is real, then it suffices to have $\Re\, f_{\clambda}\ne 0$. If $\lambda$ is non-real, then it even suffices to have $f_{\clambda}\ne 0$. In both cases this is seen by taking the inner product with $A+\eps E - \mu I \in \cS$ for an arbitrary $\mu\in \R$.

\subsection{Rank-1 matrix differential equation}
\label{subsec:rank-1-S}
\index{rank-1 matrix differential equation}
As a consequence of \eqref{stat-S}, optimizers of \eqref{eq:optimizS} are projections onto $\cS$ of rank-1 matrices. This motivates us to search for a differential equation on the manifold of rank-$1$ matrices that leads to the same stationary points.

\label{subsec:struc-rank-1-ode}
Solutions of \eqref{ode-E-S} can be written as 
 $E(t)=\Pi^\cS Z(t)$ where $Z(t)$ solves
 \begin{equation}\label{ode-E-S-Z}
\dot Z = -G_\eps(E) + \Re \langle G_\eps(E), E \rangle Z
\quad\text{ with }\ E=\Pi^\cS Z.
\end{equation}
We note that $\Re\langle E, \dot E \rangle=0$ if $\| E\|_F=1$, so that the unit Frobenius norm of $E(t)$ is conserved for all $t$. As every solution tends to a stationary point of rank 1, we project the right-hand side onto the tangent space $T_Y\cM_1$ at $Y$ of the manifold of complex rank-1 matrices $\cM_1=\cM_1(\C^{n,n})$ and consider instead the projected differential equation with solutions of rank 1:
 \begin{equation}\label{ode-E-S-1}
\dot Y = -P_Y G_\eps(E) + \Re \langle P_Y G_\eps(E), E \rangle Y \quad\text{ with }\ E=\Pi^\cS Y.
\end{equation}
Note that then
\begin{equation}\label{ode-E-S-1-Pi}
\dot E = -\Pi^\cS P_Y G_\eps(E) + \Re \langle \Pi^\cS P_Y G_\eps(E), E \rangle E \quad\text{ with }\ E=\Pi^\cS Y,
\end{equation}
which differs from the gradient flow \eqref{ode-E-S} only in that the gradient $G_\eps(E)$ is replaced by the rank-1 projected gradient $P_Y G_\eps(E)$.

For $E=\Pi^\cS Y$ of unit Frobenius norm,
$$
\Re\langle E, \dot E \rangle =  \Re\langle  E, \dot Y \rangle 
=-\Re\langle E, P_Y G_\eps(E)\rangle + \Re \langle P_Y G_\eps(E), E \rangle \,\Re\langle E,Y\rangle =0,
$$
where we used that $\Re\langle E,Y\rangle= \Re\langle E,\Pi^\cS Y\rangle=\Re\langle E,E\rangle =\| E \|_F^2 =1$. So we have
$$
\| E(t) \|_F =1 \qquad\text{for all }t.
$$
We  write a rank-1 matrix $Y\in\cM_1$ in a non-unique way as
\[
Y=\rho uv^*,
\]
where $\rho\in \R,\ \rho>0$ and $u,v\in \C^n$ have unit norm. The following lemma extends Lemma~\ref{chap:proto}.\ref{lem:uv-1} to the structured situation. It shows how  the rank-1 differential equation \eqref{ode-E-S-1}
can be restated in terms of differential equations for the  factors $u, v$ and an explicit formula for~$\rho$.

\begin{lemma}[Differential equations for the factors]
\label{lem:uv-1-S}
Every  solution $Y(t)\in \cM_1$ of the rank-1 differential equation \eqref{ode-E-S-1} with $\| \Pi^\cS Y(t) \|_F=1$
can be written as $Y(t)=\rho(t)u(t)v(t)^*$ where $u(t)$ and $v(t)$ of unit norm satisfy the differential equations 
\begin{align*}
\rho \dot u &=  - \displaystyle \tfrac\iu 2 \Im(u^*Gv)u -(I-uu^*) Gv ,
\qquad
\\
\rho \dot v &=   - \displaystyle \tfrac\iu 2 \Im(v^*Gu)v -(I-vv^*) G^*u ,
\end{align*}
where $G=G_\eps(E)$ for $E=\Pi^\cS Y = \rho\, \Pi^\cS (uv^*)$
and $\rho=1/\| \Pi^\cS(uv^*)\|_F$.
\end{lemma}

We find that with the exception of the additional positive factor $\rho$ on the left-hand side, these differential equations are of the same form as in Lemma~\ref{chap:proto}.\ref{lem:uv-1}. Note that $\rho$ is only related to the speed with which a trajectory is \bcl traversed \ecl, but does not affect the trajectory itself. However, here $G=G_\eps(E)$ for a different matrix $E=\Pi^\cS(\rho uv^*) $ instead of $E=uv^*$ in (\ref{chap:proto}.\eqref{ode-uv}).

\begin{proof}
The equation for $\rho$ is obvious because
$1 = \| E \|_F = \| \Pi^\cS(\rho uv^*)\|_F =\rho \| \Pi^\cS(uv^*)\|_F$.
\\
We write the right-hand side of \eqref{ode-E-S-1} and use \eqref{P-formula-1} to obtain for $Y=\rho uv^*$
\begin{align*}
\dot Y&= -P_Y G + \Re\langle P_Y G,E \rangle Y \\
&= -\ (I-uu^*)Gvv^* - uu^* G(I-vv^*) - uu^*G vv^* + \Re\Big\langle P_Y G,E \Big\rangle Y\\
&= -\  \Bigl( (I-uu^*)Gvv^* + \tfrac\iu 2 \Im(u^*Gv)u \Bigr)v^*  
- u \Bigl( u^*G(I-vv^*) + \tfrac\iu 2 \Im(u^*Gv) v^* \Bigr) 
\\
&\quad\, - \Bigl( \Re (u^*Gv) +  \Re\langle P_Y G,E \rangle \rho \Bigr) uv^*.
\end{align*}
Since this is equal to $\dot Y = (\rho \dot u) v^* + u (\rho \dot v^*) + \dot \rho uv^*$, we can equate $\rho \dot u$, $\rho \dot v^*$ and $\dot\rho$ with the three terms in big brackets. So we obtain the stated differential equations for $u$ and~$v$ (and another one for $\rho$, which will not be needed). Further we have $(d/dt)\| u\|^2 =2\,\Re(u^*\dot u)=0$ and analogously for $v$, which yields that  $u$ and $v$ stay of unit norm.
\qed
\end{proof}

We note that for $G=G_\eps(E)=2f_{\clambda}\,xy^*$ (see Lemma~\ref{chap:proto}.\ref{lem:gradient}) and with $\alpha=u^*x$, $\beta=v^*y$
and $\gamma=2f_{\clambda}$, we obtain differential equations that differ from
(\ref{chap:proto}.\ref{ode-uv-short}) only in the additional factor $\rho$ on the left-hand side:
\begin{equation}\label{ode-uv-short-S}
\begin{array}{rcl}
\rho \dot u &=&  \alpha\conj\beta\gamma\, u- \conj\beta\gamma \,x -\tfrac \iu2 \, \Im(\alpha\conj\beta\gamma)u
\\[3mm]
\rho \dot v &=&  \conj{\alpha}\beta\conj{\gamma}\,v -\conj{\alpha\gamma}\,y -\tfrac \iu2 \, \Im(\conj{\alpha}\beta\conj{\gamma})v.
\end{array}
\end{equation}

\subsubsection*{Stationary points.}
\index{stationary point}
The following theorem states that under some non-degeneracy conditions (see Remark~\ref{chap:proto}.\ref{rem:exceptional}),
the differential equations  \eqref{ode-E-S} and~\eqref{ode-E-S-1} yield the same stationary points.
\begin{theorem}[Relating stationary points]
\label{thm:stat-S}  
(a) Let $E\in\cS$ of unit Frobenius norm be a stationary point of the gradient system \eqref{ode-E-S} that satisfies $\Pi^\cS G_\eps(E)\ne 0$. Then, $E=\Pi^\cS Y$
for an $Y\in\C^{n,n}$ of rank~$1$ that is a stationary point of the differential equation \eqref{ode-E-S-1}.

(b) Conversely, let $Y\in\C^{n,n}$ of rank~$1$  be a stationary point of the differential equation \eqref{ode-E-S-1} such that $E=\Pi^\cS Y$
has unit Frobenius norm and $P_Y G_\eps(E)\ne 0$. Then, $E$ is a stationary point of the gradient system \eqref{ode-E-S}.
\end{theorem}

\begin{proof} Let $G=G_\eps(E)$ in this proof for short.

(a) By \eqref{stat-S}, $E=\mu^{-1} \Pi^\cS G$ for some nonzero real $\mu$. Then, $Y:=\mu^{-1}G$ is of rank 1 and we have $E=\Pi^\cS Y$. We further note that $P_Y G = \mu P_Y Y = \mu Y = G$. We thus have
$$
-P_Y  G + \Re\langle P_Y G, E \rangle Y  = -G + \Re\langle G,E \rangle Y.
$$
Here we find that
$$
 \Re\langle G,E \rangle= \Re\langle \Pi^\cS G,E \rangle =  \Re\langle \mu E,E \rangle = \mu \| E \|_F^2 = \mu.
$$
So we have
$$
-G + \Re\langle G,E \rangle Y = -G + \mu Y =0
$$
by the definition of $Y$. This shows that $Y$ is a stationary point of \eqref{ode-E-S-1}.

(b) By the argument used in the proof of Theorem~\ref{thm:stat-1}, 
stationary points $Y\in\cM_1$ of the differential equation \eqref{ode-E-S-1} are characterized as real multiples of $G$. Hence,
$E=\Pi^\cS Y$ is a real multiple of $\Pi^\cS G$, and by \eqref{stat-S}, $E=\Pi^\cS Y$ is a stationary point of \eqref{ode-E-S}.
\qed \end{proof}

\subsubsection*{Possible loss of global monotonicity and preservation of local monotonicity near stationary points.}
Since the projections $\Pi^\cS$ and $P_Y$ do not commute, we cannot guarantee the monotonicity \eqref{eq:pos-S} along solutions of \eqref{ode-E-S-1}. However, in all our numerical experiments we observed a monotonic decrease of the functional in all steps except possibly (and rarely) in the first step. In the following we will explain this monotonic behaviour locally near a stationary point, but we have no theoretical explanation for the numerically observed monotonic behaviour far from stationary points.

The first observation, already made in the proof of Theorem~\ref{thm:stat-S}, is that at a stationary point $Y$ of \eqref{ode-E-S-1}, we have $P_Y G_\eps(E)=G_\eps(E)$ for $E=\Pi^\cS Y$. Therefore, close to a stationary point, 
$P_Y G_\eps(E)$ will be close to $G_\eps(E)$.
It turns out that it is even {\it quadratically} close. This is made more precise in the following lemma.

\begin{lemma}[Projected gradient near a stationary point] \label{lem:loc-S}
\begin{samepage}
Let $Y_\star\in \cM_1$ with $E_\star=\Pi^\cS Y_\star \in \cS$ of unit Frobenius norm. Let $Y_\star$ 
be a stationary point of the rank-1 projected differential equation \eqref{ode-E-S-1}, 
with an associated target eigenvalue $\lambda$ of $A+\eps E_\star$ that is simple.
Then, there exist $\bar \delta>0$ and a real $C$ such that for all positive $\delta\le \bar\delta$ and all $Y \in \cM_1$ with $\| Y-Y_\star \| \le \delta$ and associated $E=\Pi^\cS Y$ of unit norm,
we have
\begin{equation}
\| P_Y G_\eps\left(E \right) - G_\eps\left(E \right) \| \le C \delta^2. 
\end{equation}
\end{samepage}
\end{lemma}

\begin{proof}
We consider a smooth path $Y(\tau) = u(\tau) v(\tau)^* \in \cM_1$ (with $u(\tau), v(\tau) \in \C^n$) and associated $E(\tau)=\Pi^\cS Y(\tau)$ of unit Frobenius norm with initial value
\begin{align*}
& Y(0) = Y_\star = \mu^{-1} G_\star  \quad \mbox{\rm for some nonzero real} \ \mu \quad  \mbox{and} 
\\
& G_\star = G_\eps \left( E_\star) \right) = 2 \overline f_\lambda xy^*,
\end{align*}
where $E_\star=\Pi^\cS Y_\star$ is of unit Frobenius norm and $(\lambda,x,y)$ is the target eigentriplet
of $A + \eps E_\star$ associated with the target eigenvalue~$\lambda$. 

A direct calculation of the first-order terms in the Taylor expansions of
$P_{Y(\tau)}G_\eps(E(\tau))$ and $G_\eps(E(\tau))$, which uses the formula (\ref{chap:proto}.\ref{P-formula-1}) for the projection $P_{Y(\tau)}$ and the formula (\ref{chap:proto}.\ref{eq:freegrad}) of the rescaled gradient $G_\eps(E(\tau))$, surprisingly yields that the
Taylor expansions of $P_{Y(\tau)}G_\eps(E(\tau))$ and $G_\eps(E(\tau))$ at $\tau=0$ coincide up to $O(\tau^2)$. This gives the stated result.
\qed
\end{proof}

As a direct consequence of this lemma, a comparison of the differential equations \eqref{ode-E-S-1-Pi} and \eqref{ode-E-S} yields that $\delta$-close to a stationary point, the functional decreases monotonically along solutions of \eqref{ode-E-S-1} up to $O(\delta^2)$, and even with the same negative derivative as for the gradient flow \eqref{ode-E-S} up to $O(\delta^2)$. Note that the derivative of the functional is proportional to $-\delta$ in a $\delta$-neighbourhood of a strong local minimum.
Guglielmi, Lubich \& Sicilia (\cite{GuLS23}) use Lemma~\ref{lem:loc-S} 
to prove a result on local convergence as $t\to\infty$ to strong local minima of the functional $\F_\eps$ of \eqref{eq:optimiz} for $E(t)=\Pi^\cS Y(t)$ of unit Frobenius norm associated with solutions $Y(t)$ of the rank-1 differential equation \eqref{ode-E-S-1}.

\subsection{Discrete algorithm} \label{subsec:S-discrete}
Since the differential equations for $u$ and $v$ have essentially the same form as in the unstructured case, the splitting algorithm of Subsection~\ref{subsec:proto-numer} extends in a straightforward way, and also the stepsize selection is readily extended. We refer to Guglielmi, Lubich \& Sicilia (\cite{GuLS23}) for details of the algorithm and for numerical experiments with sparse matrices and matrices perturbed with matrices of given rank and co-rank.

\bcl
\subsection{Initial perturbation} 
\label{subsec:S-init}
\index{initial perturbation}
The choices (i) and (ii) of the initial perturbation matrix $E(0)$ in Section~\ref{subsec:init} are readily adapted to the structured situation, simply replacing the free gradient $G$ by the projected gradient $\Pi^\cS G$.

A structured analogue of the more elaborate approach (iii) of Section~\ref{subsec:init} is obtained by first solving the rank-1 differential equation without norm constraint 
$$
\dot Y = -P_Y G_\eps(E) \quad\text{ with }\ E=\Pi^\cS Y
$$
up to the first time $\bar t$ where $\| E(\bar t) \|_F=1$, starting from a small multiple of the rank-1 matrix determined by (i) or (ii). We then continue with the norm-constrained rank-1 differential equation \eqref{ode-E-S-1} starting from $Y(\bar t)$.
\ecl

\section{Structured pseudospectra}
\label{sec:ps-S}

\subsection{Motivation and definitions}
\label{subsec:ps-motivation-S}
We take up the motivating example for the $\eps$-pseudospectrum of a matrix $A$ and consider the linear dynamical system $\dot x(t)=Ax(t)$. We now ask if asymptotic stability is preserved under real or structured perturbations $\Delta\in\cS$ of norm bounded by a given $\eps>0$.  
This question is answered by the sign of the {\em $\cS$-structured $\eps$-pseudospectral abscissa} defined as
\index{pseudospectral abscissa!structured}
\begin{align*}
\alpha_\eps^\cS(A) = \max \{ \Re\,\lambda : \  \ &\text{There exists $\Delta\in\cS$ with $\| \Delta \| \le \eps$ such that}
\\
&\text{$\lambda$ is an eigenvalue of $A+\Delta$} \}.
\end{align*}
This number is computed by solving a problem  \eqref{eq:optimiz0rF} for $\cS=\R^{n,n}$, or 
\eqref{eq:optimiz0S}--\eqref{eq:optimizS} for an arbitrary complex-linear or real-linear structure space $\cS$, in each case with the function to be minimized given by $f(\lambda,\clambda)=-\tfrac12(\lambda+\clambda)=-\Re\,\lambda$.

It is useful to rephrase the question in terms of the $\eps$-pseudospectrum, which is defined as follows.

\begin{definition} \label{def:ps-S} \index{pseudospectrum!structured}
The $\cS$-structured $\eps$-{\it pseudospectrum} of the matrix $A$ is the set
\begin{equation}\label{eq:epsps-S}
\Lambda_\eps^\cS(A) = \{ \lambda \in \C : \ \ \text{$\lambda \in \Lambda(A+\Delta)$ for some $\Delta\in\cS$ with $\| \Delta \| \le \eps$}\},
\end{equation}
where $\Lambda(M)\subset \C$ denotes the spectrum of a square matrix $M$.
\end{definition}

For $\cS=\C^{n,n}$, $\Lambda_\eps(A)=\Lambda_\eps^\C(A)=\Lambda_\eps^\cS(A) $ is the {\it complex} $\eps$-pseudospectrum, and for $\cS=\R^{n,n}$, $\Lambda_\eps^\R(A)=\Lambda_\eps^\cS(A) $ is known as the {\it real} $\eps$-pseudospectrum.

The  {\em $\cS$-structured $\eps$-pseudospectral abscissa} of the matrix $A$  can then be rewritten more compactly as
\begin{equation}\label{stability-abscissa-S}
\alpha_\eps^\cS(A) = \max \{ \Re\,\lambda : \ \lambda \in \Lambda_\eps^\cS(A)\}.
\end{equation}
An analogous quantity for discrete-time linear dynamical systems $x_{k+1}=Ax_k$ is the {\em $\cS$-structured 
$\eps$-pseudospectral radius} of the matrix $A$,
\index{pseudospectral radius!structured}
\begin{equation}\label{stability-radius-S}
\rho_\eps^\cS(A) = \max \{ | \lambda | : \ \lambda \in \Lambda_\eps^\cS(A)\}.
\end{equation}

\subsection{Extremal structured perturbations}
\index{extremal perturbation}
For structured pseudospectra, there is no characterization in terms of singular value decompositions or resolvent bounds. There is, however, an analogue of Theorem~\ref{chap:pseudo}.\ref{thm:Delta-C} which characterizes extremal perturbations, at least at boundary points $\lambda\in\partial\Lambda_\eps^\cS(A)$ where the boundary is differentiable and thus admits an outer normal $\e^{\iu\theta}$.

\begin{theorem}[Extremal structured perturbations] \label{thm:Delta-S}
Let $\lambda\in\partial \Lambda_\eps^\cS(A)$ be on a smooth section of the boundary, with outer normal $\e^{\iu\theta}$ at $\lambda$.
Let  $\Delta\in\cS$ of Frobenius norm $\eps$ be such that $A+\Delta$ has $\lambda$ as a simple eigenvalue. Let $x$ and $y$ be left and right eigenvectors of $A+\Delta$ to the eigenvalue $\lambda$, of unit norm and with $x^*y>0$. If the matrix $\Pi^\cS(\e^{\iu\theta}xy^*)$ is non-zero, where $\Pi^\cS$ is the orthogonal projection \eqref{Pi-S} onto $\cS$, then
$$
\Delta = \eps \,\eta\,\Pi^\cS(\e^{\iu\theta} x y^*),
$$
where $\eta=1/\|\Pi^\cS(\e^{\iu\theta} x y^*)\|_F>0$.
\end{theorem}

\begin{proof} The proof is based on Section~\ref{subsec:gradient-flow-S}, in particular (\ref{stat-S}).
In the proof we denote the given $\lambda\in\partial \Lambda_\eps^\cS(A)$ as $\lambda_0$ and use $\lambda$ to denote a complex variable. Similarly, we denote the matrix $\Delta$ in the statement of the lemma by $\Delta_0$ and use $\Delta$ for a generic matrix in $\cS$. 

We set $\mu=\lambda_0+\delta \e^{\iu\theta}$ with a small $\delta>0$.
For the function
$$
f(\lambda,\clambda)=|\mu-\lambda|^2=(\mu-\lambda)(\conj{\mu}-\clambda)
$$
we consider the optimization problem (\ref{eq:optimizS})
with $\lambda(A+\Delta)$ denoting the eigenvalue of $A+\Delta$ closest to $\lambda_0$, which has $\Delta_0=\eps E_0$ with $\| E_0\|_F=1$ as a solution, and $\lambda(A+\Delta_0)=\lambda_0$. In particular, $E_0$ is then a stationary point of the gradient system  
(\ref{ode-E-S}) and hence satisfies (\ref{stat-S}) with a negative factor as in (\ref{chap:proto}.\ref{stat-neg}), which yields that $E_0$ is a negative real multiple of the structured gradient $G_\eps^\cS(E_0)= \Pi^\cS(f_{\clambda} xy^*)$ with
$f_{\clambda} = \partial f/\partial \clambda(\lambda_0,\clambda_0)=-(\mu-\lambda_0)=-\delta\e^{\iu\theta}$, which yields $\Delta_0$ as stated. This holds under the condition that  
$G_\eps^\cS(E_0)\ne 0$ is non-zero, i.e.,  $\Pi^\cS(\e^{\iu\theta}xy^*)\ne 0$.
\qed
\end{proof}

\begin{remark}
    In the unstructured real case $\cS=\R^{n,n}$, where $\Pi^\cS(Z)=\Re\,Z$, the condition $\Re(\e^{\iu\theta} xy^*)\ne 0$ excludes sections of $\partial \Lambda_\eps^\R(A)$ that consist of intervals on the real line. 
\end{remark}

\subsection{Tracing the boundary of structured pseudospectra}
\label{subsec:tracing-S}
\index{pseudospectrum boundary}

\subsubsection*{Structured ladder algorithm.}
\index{ladder algorithm}
The ladder algorithm of Section \ref{subsec:ladder} is readily extended to compute smooth sections of the boundary of the structured $\eps$-pseudospectrum of a given matrix for  a general linear structure $\cS$:

Let $\lambda_0\in\partial\Lambda_\eps^\cS(A)$ be a simple eigenvalue of $A+\eps \eta_0\Pi^\cS(u_0v_0^*)$ (with $u_0$ and $v_0$ of unit norm and $\eta_0=1/\|\Pi^\cS(u_0v_0^*)\|$) that lies on a smooth section of $\partial\Lambda_\eps^\cS(A)$, with outer normal $\e^{\iu\theta_0}= |u_0^* v_0|/(u_0^* v_0)$ at $\lambda_0$, as follows from Theorem~\ref{thm:Delta-S}. 

With a small distance $\delta>0$, we define the nearby point $\mu_0$ on the straight line normal to $\partial\Lambda_\eps(A)$ at $\lambda_0$,
$$
\mu_0=\lambda_0 + \delta \e^{\iu\theta_0}.
$$
We add a tangential component, either to the left ($+$) or to the right ($-$),
$$
\mu_1= \mu_0 \pm \iu \delta \e^{\iu\theta_0}.
$$
We then apply the structured eigenvalue optimization algorithm of Section~\ref{sec:proto-structured} for the function
$$
f(\lambda,\clambda)=(\lambda -\mu_1)(\clambda -\conj{\mu_1})=|\lambda -\mu_1|^2,
$$
choosing $E_0=\Pi^\cS(u_0v_0^*)$ as the starting iterate. That algorithm aims to compute $u_1$ and $v_1$ of unit norm such that $A+\eps \Pi^\cS(u_1 v_1^*)$ has the boundary point $\lambda_1\in\partial\Lameps^\cS(A)$ nearest to $\mu_1$ as an eigenvalue. At the point $\lambda_1$ we have the outer normal $(\mu_1 - \lambda_1)/| \mu_1 - \lambda_1 | $.

We continue from $\lambda_1$ and $u_1$, $v_1$ in the same way as above, constructing a sequence $\lambda_k$ ($k\ge 1$) of points on the boundary of the structured pseudospectrum $\Lameps^\cS(A)$ with approximate spacing $\delta$.

Let us consider here some illustrative examples, in which we apply this structured version of the ladder algorithm.

\subsubsection{Real pseudospectra}

\begin{figure}[ht]
\centerline{
\includegraphics[scale=0.36]{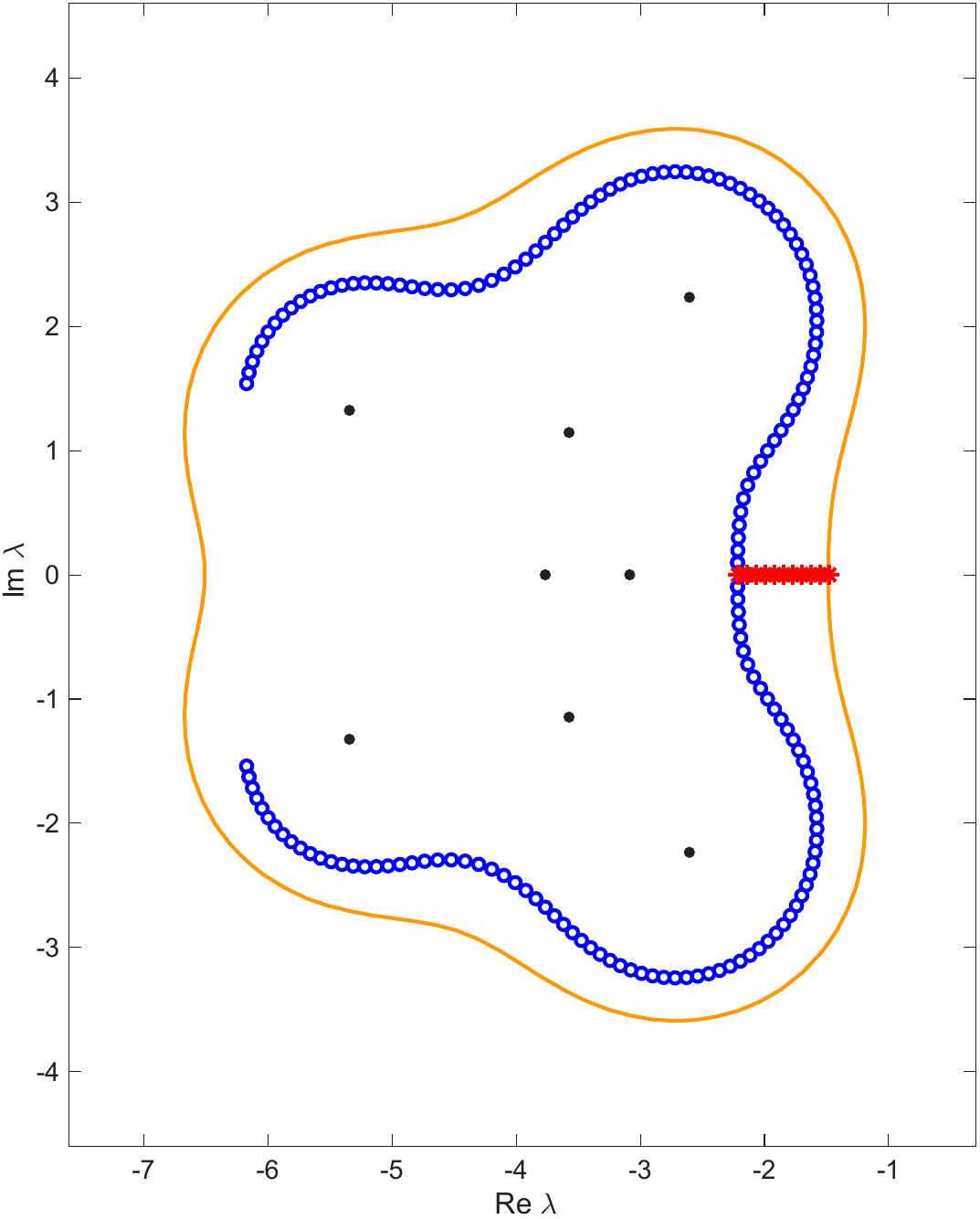}
}
\caption{Section of the boundary of the real pseudospectrum $\Lameps^{\R}(A)$ computed by the ladder algorithm 
compared with the unstructured pseudospectrum $\Lameps(A)$.
\label{fig:R2}}
\end{figure}

We consider the real matrix in (\ref{chap:proto}.\ref{eq:example}) and $\eps=1$ and
are interested in computing the boundary of the real $\eps$-pseudospectrum $\Lameps^{\R}(A)$.
In Figure \ref{fig:R2} we show a section \bcltwo (the right one) \ecltwo of the real $\eps$-pseudospectrum
of $A$ with $\eps=1$, as computed with the ladder algorithm, 
together with the boundary of the complex $\eps$-pseudospectrum $\Lambda_\eps(A)$.
%
%

\subsubsection{Hamiltonian pseudospectra}
\begin{figure}[ht]
\vskip -6cm
\centerline{
\ \hskip -2.5cm \includegraphics[scale=0.5]{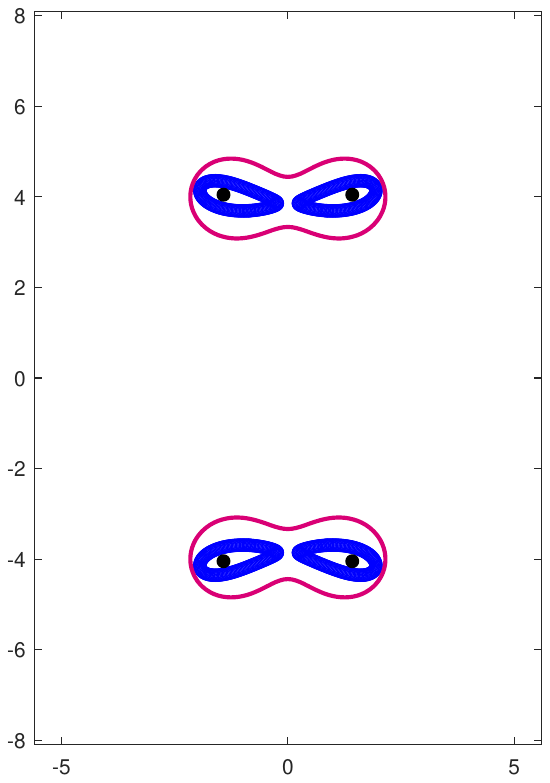}
}
\caption{Boundary points on the real Hamiltonian $\eps$-pseudospectrum $\Lameps^{\mathcal{H}}(H)$ computed by the ladder algorithm versus the complex unstructured pseudospectrum $\Lameps(H)$. 
\label{fig:H3}}
\end{figure}

\index{Hamiltonian matrix}
We show an example of a real Hamiltonian $\eps$-pseudo\-spectrum, $\Lambda_\eps^\mathcal{H}(H)$ for $\mathcal{H}=\text{Ham}(\R^{n,n})$ the space of $n \times n$ real Hamiltonian matrices.
Let
\begin{equation}
H = \left( \begin{array}{rrrr}
    1.0  &  1.6 &   1.2 &   0.4 \\
    2.2  & -0.6 &   0.4 &  -4.4 \\
   -4.0  & -7.4 &  -1.0 &  -2.2 \\
   -7.4  &  6.0 &  -1.6 &   0.6
\end{array}
\right)
\end{equation}
and set $\eps=0.4$.
In Figure \ref{fig:H3} we show the whole set of boundary points of $\Lambda_\eps^\cS(H)$ as computed by the ladder algorithm.

\begin{figure}[h!]
\vskip -1cm
\centerline{
\includegraphics[scale=0.36]{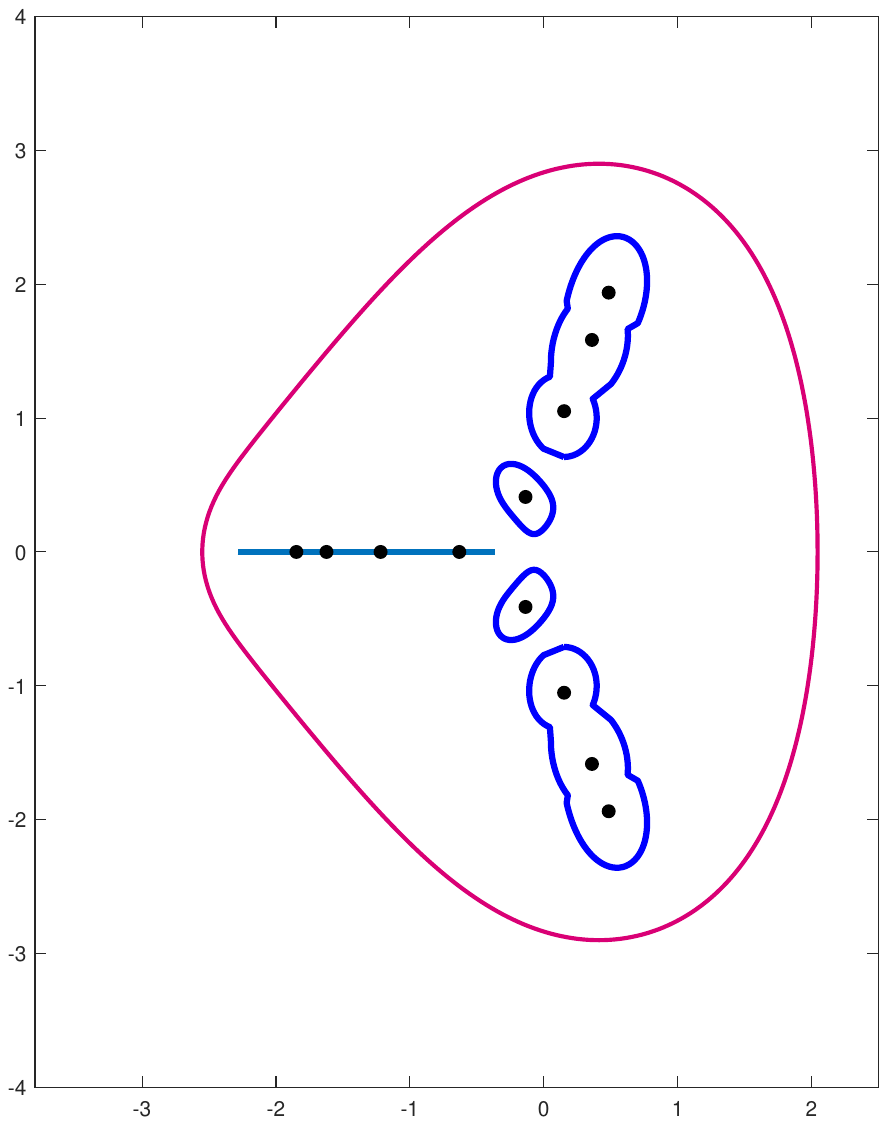}
}
 \vspace{-1.5cm}
\caption{
The penta-diagonal real Toeplitz-structured $\eps$-pseudo\-spec\-trum $\Lameps^{\mathcal{T}_5}(T)$, as computed by the ladder algorithm, inside the unstructured $\eps$-pseudospectrum. 
\label{fig:T2}}
\end{figure}

\subsubsection{Toeplitz pseudospectra}
\label{sec:toepps}
\index{Toeplitz matrix}

We further consider the case of Toeplitz matrices, for which -- as far as we know -- there are no available 
tools to draw structured pseudospectra and compute related quantities.
We denote by $\mathcal T_k$ the space of $k$-diagonal real Toeplitz matrices and $\Lambda_\eps^{{\mathcal T}_k}(A)$
the associated structured $\eps$-pseudospectrum.

We consider the $12 \times 12$ penta-diagonal Toeplitz matrix
\begin{equation}
T =T_{5}(s_1,s_2,d,t_1,t_2) \in \mathcal{T}_5, 
\label{eq:example1}
\end{equation}
that is the matrix with elements
\begin{eqnarray*}
&& a_{i,i} = d, \qquad i=1,\ldots,n
\\
&& a_{i+1,i} = s_1, \qquad i=1,\ldots,n-1
\\
&& a_{i+2,i} = s_2, \qquad i=1,\ldots,n-2
\\
&& a_{i-1,i} = t_1, \qquad i=1,\ldots,n-1
\\
&& a_{i-2,i} = t_2, \qquad i=1,\ldots,n-2.
\end{eqnarray*}
with entries $d=-0.3$, $s_1=-0.1$, $s_2=-0.3$, $t_1=2$, $t_2=0.5$.
We use the ladder algorithm to compute the $\mathcal T_5$-structured $\eps$-pseudospectrum of $T$ for $\eps=0.6$.
In Figure \ref{fig:T2} we show both the structured and the unstructured pseudospectra. The so obtained structured $\eps$ pseudospectrum agrees with what is obtained by intensive sampling of eigenvalues of perturbed matrices $A+\eps E$ with $E\in\cS$ of Frobenius norm at most~1.

\section{Computing the structured stability radius}
\label{sec:stab-rad-S}
\index{stability radius!structured}
\index{structured stability radius}
\index{distance to instability!under structured perturbations}
\subsection{Computing the structured pseudospectral abscissa}
\index{pseudospectral abscissa!structured}
While we are not aware of any algorithm in the literature for computing the $\cS$-structured pseudospectral abscissa for general linear structures~$\cS$,
the projected rank-1 algorithm of Section~\ref{sec:proto-structured} used with the rightmost eigenvalue as target eigenvalue and the functional $f(\lambda,\bar\lambda)= 
-\frac12(\lambda+\bar\lambda)=-\Re\,\lambda$  directly yields such an algorithm. The eigenvalue optimization problem
\eqref{E-eps-2l}
then becomes the maximization problem
\begin{equation} \label{E-eps-stab-radius-S}
E(\eps) = \arg\max\limits_{E \in \cS, \| E \|_F = 1} \Re\,\lambda(A+\eps E),
\end{equation}

\subsection{Two-level iteration}\index{two-level iteration}
The two-level approach of Chapter~\ref{chap:two-level} applies directly to computing the complex, real or structured stability radius of a Hurwitz-stable matrix $A$. 
In the inner iteration, the eigenvalue optimization problem \eqref{E-eps-stab-radius-S} is of the class studied in Section~\ref{sec:proto-structured} 
and is solved with the projected rank-1 approach developed there.

In the outer iteration, the optimal perturbation size $\oeps$ is determined from Equation (\ref{chap:two-level}.\ref{eq:zero}), which here reads
\begin{equation} \label{zero-stab-radius-S}
\Re\,\lambda(A+\oeps E(\oeps)) =0.
\end{equation}
\index{Newton--bisection method}
This scalar equation is solved by the Newton--bisection algorithm of Section~\ref{subsec:Newton--bisection}, where the derivative formula of Theorem~\ref{chap:two-level}.\ref{thm:phi-derivative} is valid also for the projected gradient $G(\eps)=-\Pi^\cS ( x(\eps)y(\eps)^*)$ used here; in particular, $G=-\Re(xy^*)$ in the real unstructured case.

\bcl The HEC algorithm can equally be applied here. The monotone Newton--bisection iteration can be applied in the full-rank case, where a gradient descent is used. However, the flow of the rank-1 matrix differential equation of Section~\ref{subsec:rank-1-S} is not guaranteed to reduce the objective function far from stationary points (although it usually does in our numerical experience). This can possibly impair monotonicity and convergence of the modified Newton--bisection algorithm of Section~\ref{subsec:Newton-bisection-monotone}.
\ecl

\subsubsection{An illustrative example}
We consider the  block tridiagonal sparse matrix  \emph{rdbrusselator} of dimension $n=3200$ from the Matrix Market, and shift it by $\frac12 \Id$ to make it Hurwitz. We choose the structure space $\cS$ as the space of matrices with the same sparsity pattern and use the algorithm described above to compute the structured stability radius $\oeps$.



\begin{table}[hbt]
\begin{center}
\begin{tabular}{|l|l|l|l|}\hline
  $k$ & $\eps_k$ & $\alpha_{\eps_k}^\cS(E(\eps_k))$ & $\#$ eigs \\
 \hline
\rule{0pt}{9pt}
\!\!\!\! 
        $0$         & $0.5$               & $-0.334680920533956$ & $7$   \\
	$1$         & $6.548675639782492$ & $0.024507744851970$  & $1$    \\
	$2$         & $6.400112687385496$ & $3.068774405129098$  & $142$  \\
        $3$         & $1.126792695169999$ & $-0.233912658683879$ & $134$  \\
	$4$         & $1.444080848272038$ & $0.002458991032461$  & $2$    \\
	$5$         & $1.440784407976959$ & $0.000001735157713$  & $4$    \\
	$6$         & $1.440782082056653$ & $0.000000038671781$  & $3$    \\
	$7$         & $1.440782030219339$ & $0.000000002000975$  & $4$    \\   
\hline
\end{tabular}
\vspace{2mm}
\caption{Computation of the structured stability radius: 
computed values $\eps_k$, $\alpha_{\eps_k}(E(\eps_k))$ 
and number of eigenvalue computations of the inner algorithm.\label{tab:Hurdb}}
\end{center}
\end{table}
In Table \ref{tab:Hurdb} we report the values $\eps_k$ computed by the outer Newton--bisection iteration together with the associated structured $\eps$-pseudospectral abscissa and the number of eigenvalue computations in each iteration.
Interestingly, there is a clear jump to a better local extremizer passing from iteration 1 to
iteration 2 in the Newton process in Table~\ref{tab:Hurdb}.

\section{Transient dynamics under structured perturbations}\label{sec:transient-S}
\index{transient bound}
Pseudospectra come with different uses: 
\begin{itemize}
    \item The unstructured and structured pseudospectra provide information on the sensitivity of the spectrum to \bcl unstructured and structured \ecl perturbations;
    \item The complex unstructured pseudospectrum provides information on the transient behaviour of linear differential equations, for which the pseudospectral abscissa $\alpha_\eps(A)$ and the stability radius $\oeps$ with $\alpha_{\oeps}(A)=0$ are key quantities; see the bounds 
    (\ref{chap:pseudo}.\ref{transient-bound})-(\ref{chap:pseudo}.\ref{oeps-L2-bound}) in Section~\ref{subsec:ps-exp}. 
\end{itemize}
The transient bounds in the second item are obtained
from resolvent bounds that are available from the characterization (\ref{chap:pseudo}.\ref{ps-res}), viz.,
\begin{align}
\Lambda_\eps(A)&= \{ \lambda \in \C \,:\, \text{$\lambda$ is an eigenvalue of $A+\Theta$ for some $\Theta\in\C^{n\times n}$ with $\| \Theta \|_F \le \eps$}\}
\nonumber
\\
&= \{ \lambda \in \C \,:\, \| (A- \lambda I)^{-1} \|_2 \ge \eps^{-1} \}.
\label{ps-viz}
\end{align}
In contrast, structured pseudospectra are not directly related to resolvent bounds, and we are not aware of any case where transient bounds have been inferred from knowledge of a structured pseudospectrum.

\subsection{Joint unstructured--structured pseudospectrum}
\index{pseudospectrum!unstructured--structured}
Let again the structure space $\cS$ be a complex-linear or real-linear subspace of $\C^{n\times n}$, and let $\eps>0$ and $\delta>0$. As we will show in the next subsection,
robust transient bounds under structured perturbations in $\cS$ become accessible when the notions of the unstructured $\eps$-pseudospectrum $\Lambda_\eps(A)$, which is characterized by resolvent bounds, and the structured $\delta$-pseudospectrum $\Lambda_\delta^\cS(A)$ are combined in a joint pseudospectrum. 
\index{resolvent bound}
This allows us to
use resolvent bounds as with unstructured pseudospectra and structured
perturbations as with structured pseudospectra. We define
\begin{align}
\nonumber
\Lambda_{\delta,\eps}^{\cS}(A) :\!&= \{ \lambda \in \C \,:\, \text{$\lambda\in \Lambda_\eps(A+\Delta)$ for some $\Delta\in\cS$ with $\| \Delta \|_F \le \delta$} \}
\\
\label{uss-ps}
&= \{ \lambda \in \C \,:\, \text{$\lambda$ is an eigenvalue of $A+\Delta+ \Theta$ for some} 
\\[-1mm]
&\hskip 2cm \text{$\Delta\in\cS$ with $\| \Delta \|_F \le \delta$ and $\Theta \in \C^{n\times n}$ with $\| \Theta \|_F \le \eps$} 
\}
\nonumber
\\
&=  \{ \lambda \in \C \,:\, \text{$\| (A+\Delta-\lambda I)^{-1} \|_2 \ge \eps^{-1}$ for some $\Delta\in\cS$ with $\| \Delta \|_F \le \delta$} \}.
\nonumber
\end{align}
The equalities follow from~\eqref{ps-viz}. 

Let now the matrix $A\in \C^{n\times n}$  be a Hurwitz matrix, and let $\oeps>0$ be its stability radius and $\delta_\star^\cS>0$ its $\cS$-structured stability radius.
For a given $\eps>0$ we define the {\em $\cS$-structured $\eps$-stability radius} of $A$ as
\index{structured $\eps$-stability radius}
\begin{equation} \label{eps-stab-r}
    \delta_\eps^\cS(A) := \max\{ \delta>0 \,:\, \text{$\Lambda_{\delta,\eps}^{\cS}(A)$ has no points with positive real part} \} 
\end{equation}
and note that $\delta_\eps^\cS(A)\nearrow \delta_\star^\cS$ as $\eps \searrow 0$.
Conversely, we define for a given $\delta>0$
\begin{equation} \label{delta-stab-r}
    \eps_\delta^\cS(A) := \max\{ \eps>0 \,:\, \text{$\Lambda_{\delta,\eps}^{\cS}(A)$ has no points with positive real part} \} 
\end{equation}
for which $\eps_\delta^\cS(A) \nearrow \oeps$ as $\delta \searrow 0$.


By the second line of \eqref{uss-ps}, the structured $\eps$-stability radius is characterized as the largest $\delta>0$ such that all eigenvalues of $A+\Delta + \Theta$ have nonpositive real part for every $\Delta\in \cS$ with $\|\Delta\|_F \le \delta$
and every $\Theta\in \C^{n\times n}$ with $\|\Theta\|_F \le \eps$. This characterization is used in the derivation of the algorithm for computing the structured $\eps$-stability radius that will be given in Section~\ref{sec:eps-stab}.
On the other hand, the last line of \eqref{uss-ps} yields that with $\delta=\delta_\eps^\cS(A)$,
\begin{equation}
    \label{robust-res-bound}
\frac 1\eps = \max_{\Delta\in\cS, \| \Delta\|_F \le \delta}\max_{\,\text{Re}\,\lambda \ge 0\,} \| (A+\Delta-\lambda I)^{-1} \|_2,   
\end{equation}
so that the resolvent norm of $A+\Delta$  is bounded in the \bng complex right \eng half-plane by $\eps^{-1}$ uniformly for all
$\Delta\in\cS$ with $\|\Delta\|_F \le \delta$.

In \eqref{robust-res-bound} we considered $\delta=\delta_\eps^\cS(A)$ as a function of $\eps$. Conversely, for a given $\delta$ with $0<\delta<\delta_\star^\cS$, the bound \eqref{robust-res-bound} is valid with $\eps=\eps_\delta^\cS(A)$.
When we determine $\delta_\eps^\cS(A)$ and $\eps_\delta^\cS(A)$, we thus answer the following two questions:
\begin{itemize}
    \item Up to which size of structured perturbations are the resolvent norms of the perturbed matrices within a given bound in the \bng  complex right \eng half-plane? (Given $\eps$, find $\delta$.)
    \item For a given size of structured perturbations, what is the smallest common bound for the resolvent norms of the perturbed matrices in the \bng  complex right \eng half-plane? (Given $\delta$, find $\eps$.)
\end{itemize}

\subsection{Robust transient bounds under structured perturbations}
We give robust bounds of solutions of linear differential equations that follow from the robust resolvent bound \eqref{robust-res-bound}. In the following, $A\in \C^{n\times n}$  is a matrix with all eigenvalues of negative real part, $\oeps>0$ is its stability radius,
and $0<\eps<\oeps$.
The first result is a variant of (\ref{chap:pseudo}.\ref{transient-perturbed-oeps}) with structured perturbations.
\begin{proposition}[Matrix exponential under structured perturbations]
\label{lem:exp-bound-S} \index{transient bound!robust}
For every perturbation $\Delta\in\cS$ with $\|\Delta\|_F\le \delta_\eps^\cS(A)$,
$$
\bigl\| \e^{t(A+\Delta)} \bigr\|_2  \le  
\frac{|\Gamma|}{2\pi\eps} \qquad\text{for all }\ t>0,
$$
where $\Gamma$ is a closed contour in the closed complex left half-plane that is a union of (i) the part in the complex left half-plane of a contour (or union of several contours) that surrounds the pseudospectrum $\Lambda_{\eps+\delta}(A)$ with $\delta=\delta_\eps^\cS(A)$ and (ii) one or several intervals on the imaginary axis that close the contour. Moreover, $|\Gamma|$ is the length of $\Gamma$.
\end{proposition}

\begin{proof} Let $\Delta\in\cS$ with $\|\Delta\|_F\le \delta= \delta_\eps^\cS(A)$.
The bound follows from the Cauchy integral representation
$$
\e^{t(A+\Delta)} = \frac1{2\pi\iu} \int_{\Gamma} \e^{t \lambda}\, (\lambda I -A-\Delta)^{-1}\, d\lambda
$$
on noting that $|\e^{t \lambda}|\le 1$ for all $\lambda\in\Gamma$ and
\begin{equation}\label{res-bound-Gamma}
\| (\lambda I -A-\Delta)^{-1} \|_2  \le  \frac1\eps \quad\ \text{ for all } \lambda\in\Gamma.
\end{equation}
This resolvent bound holds true because 
(i) the inclusion $\Lambda_{\eps}(A+\Delta) \subset \Lambda_{\eps+\delta}(A)$ implies the bound \eqref{res-bound-Gamma}
for $\lambda$ in the closure of $\C \setminus \Lambda_{\eps+\delta}(A)$, and (ii) the resolvent bound \eqref{robust-res-bound} implies the bound \eqref{res-bound-Gamma} for all $\lambda$ on the imaginary axis.
\qed
\end{proof}

The next result is a structured version of the robust bound (\ref{chap:pseudo}.\ref{oeps-L2-bound-robust}). 

\begin{proposition}[Inhomogeneous linear system under structured perturbations]
\label{lem:ode-inhom-bound-S}
For all perturbations $\Delta\in\cS$ with $\|\Delta\|_F\le \delta_\eps^\cS(A)$,
solutions to the inhomogeneous linear differential equations 
$$
\dot x_\Delta(t)=(A+\Delta)x_\Delta(t)+ f(t), \qquad x_\Delta(0)=0, 
$$
have the common bound (with $\|\cdot\|$ the Euclidean norm on $\C^n$)
\begin{equation} \label{robust-L2-bound}
  \biggl(  \int_0^T \| x_\Delta(t) \|^2 \, dt \biggr)^{1/2} \le \frac1\eps \,
   \biggl(  \int_0^T \| f(t) \|^2 \, dt \biggr)^{1/2} , \qquad 0\le T \le \infty.
\end{equation}
\end{proposition}

\begin{proof} We extend  $x_\Delta(t)$ and $f(t)$ to $t<0$ by zero. Their Fourier transforms $\widehat x_\Delta$ and $\widehat f$ are then related by $\iu\omega\,\widehat x_\Delta(\omega)=(A+\Delta)\widehat x_\Delta(\omega) + \widehat f(\omega)$ for all $\omega\in\R$, i.e.,
$$
\widehat x_\Delta (\omega)=(\iu\omega I - A-\Delta)^{-1}\widehat f(\omega), \qquad \omega \in \R.
$$
Applying the Plancherel formula twice yields
\begin{align*}
    &\int_\R \| x_\Delta(t) \| ^2 \, dt = \int_\R \| \widehat x_\Delta (\omega)\| ^2\, d\omega = 
    \int_\R \| (\iu\omega I - A-\Delta)^{-1}\widehat f(\omega) \| ^2 \, d\omega
    \\
    &\le \max_{\omega\in\R} \| (\iu\omega I - A-\Delta)^{-1} \|_2^2 \, \int_\R \| \widehat f (\omega)\| ^2\, d\omega 
    \\
    &=
    \max_{\omega\in\R} \| (\iu\omega I - A-\Delta)^{-1} \|_2^2 \,  \int_\R \| f(t) \| ^2 \, dt.
\end{align*}
Using the resolvent bound \eqref{robust-res-bound} and causality, we obtain the bound
\eqref{robust-L2-bound}.
\qed
\end{proof}
 
Conversely, let $\delta$ be given with $0<\delta<\delta_\star^\cS$, where $\delta_\star^\cS$ is the $\cS$-structured stability radius of $A$. 
The same proofs then  yield that the bounds of Propositions
\ref{lem:exp-bound-S} and \ref{lem:ode-inhom-bound-S} hold true with $\eps=\eps_\delta(A)$.

\section{Computing the structured $\eps$-stability radius}
\label{sec:eps-stab}

In this section we derive an algorithm for computing the structured $\eps$-stability radius $\delta_\eps^\cS(A)$ (as well as the dual quantity $\eps_\delta^\cS(A)$) with essentially the same computational effort as for the stability radius $\oeps(A)$ in Chapter~\ref{chap:two-level}, still using a rank-1 constrained gradient flow within a two-level iteration.

\subsection{Two-level iteration} \index{two-level iteration}
Our numerical approach to computing the structured $\eps$-stability radius $\delta_\eps^\cS(A)$ again uses a two-level iteration. 

The target eigenvalue $\lambda(M)$ of a matrix $M$ is again chosen as an eigenvalue of $M$ of maximal real part (and among those, e.g.~the one with maximal imaginary part).
For $\delta>0$ and $\eps>0$ we introduce the functional
\begin{equation}\label{F-delta}
\F_{\delta,\eps}(E^\cS,E) = - \Re\, \lambda(A+\delta E^\cS + \eps E) 
\end{equation}
for $E^\cS\in \cS$ and $E \in \C^{n,n}$, both of unit Frobenius norm. With this functional we follow the two-level approach of Section~\ref{sec:two-level}:
\begin{itemize}
\item {\bf Inner iteration:\/} We aim to compute  matrices $E^\cS_{\delta,\eps} \in\cS$  
and $E_{\delta,\eps} \in \C^{n,n}$, both of unit Frobenius norm,
that minimize $\F_{\delta,\eps}$:
\begin{equation} \label{E-delta}
(E^\cS_{\delta,\eps},E_{\delta,\eps}) = \arg\min\limits_{E^\cS \in \cS, E\in \C^{n,n} \atop \| E^\cS\|_F= \| E \|_F = 1} \F_{\delta,\eps}(E^\cS,E).
\end{equation}

\item {\bf Outer iteration:\/} Given $\eps>0$, we compute the smallest positive value $\delta_\eps$ with
\begin{equation} \label{zero-delta-S}
\phi(\delta_\eps,\eps)= 0,
\end{equation}
where $\phi(\delta,\eps)=  \F_{\delta,\eps}(E^\cS_{\delta,\eps},E_{\delta,\eps})=-\Re\,\lambda\bigl(A+\delta E^\cS_{\delta,\eps}+ \eps E_{\delta,\eps}\bigr) $.
\end{itemize}

\bcl Then\ecl, 
$\delta_\eps$ is the $\cS$-structured $\eps$-stability radius $\delta_\eps^\cS(A)$.
Conversely, we obtain the inverse common resolvent bound $\eps_\delta^\cS(A)$ on the complex right half-plane under structured perturbations of size $\delta$ as the smallest positive value $\eps_\delta>0$ with
\begin{equation} \label{zero-delta-S-2}
\phi(\delta,\eps_\delta)= 0.
\end{equation}

\subsection{Reduced eigenvalue optimization}
\label{subsec:reduced-opt}

We will significantly reduce the dimension of the optimization problem \eqref{F-delta}--\eqref{E-delta} using the following properties of the optimal perturbation matrices $E^\cS \in \cS$ and $E \in \C^{n,n}$.

\begin{theorem} [Optimizers]
\label{thm:ESE}
    Let $E^\cS \in \cS$ and $E \in \C^{n,n}$, both of Frobenius norm~1, solve the eigenvalue optimization problem \eqref{F-delta}--\eqref{E-delta}. Assume that the
    target eigenvalue $\lambda(A+\delta E^\cS + \eps E)$ is simple and that the left and right eigenvectors $x$ and $y$ with positive inner product satisfy $\Pi^\cS(xy^*)\ne 0$. Then,
    \begin{equation} \label{ESE}
    \begin{aligned}
    &\text{$E$ is a real multiple of $xy^*$ (hence of rank~1) and }
    \\
    &\text{$E^\cS$ is a real multiple of $\Pi^\cS \! E$.}
    \end{aligned}
    \end{equation}
    Moreover, if the equality constraints for the Frobenius norm $1$ are replaced by inequality constraints $\le 1$, then 
    $E^\cS=\eta \Pi^\cS E$ with $\eta=1/\| \Pi^\cS \! E \|_F>0$.
\end{theorem}

\begin{proof}
We combine the procedures of Chapter~\ref{chap:proto} and Section~\ref{sec:proto-structured}. We find that along a differentiable path $(E^\cS(t),E(t))$ in $\cS\times \C^{n,n}$ we have, assuming simple target eigenvalues
$\lambda (A+\delta E^\cS(t)+\eps E(t))$,
\begin{align*}
\frac 1{\kappa(t)} \frac{d}{dt} \F_{\delta,\eps}(E^\cS(t),E(t)) &=
\Re\langle \Pi^\cS G_{\delta,\eps}(E^\cS(t),E(t)), \delta \dot E^\cS(t)\rangle 
\\[-2mm]
&\quad + 
 \Re\langle G_{\delta,\eps}(E^\cS(t),E(t)), \eps \dot E(t)\rangle
\end{align*}
with the rescaled gradient 
$$
G_{\delta,\eps}(E^\cS,E)=-xy^*,
$$
where $x$ and $y$ are left and right eigenvectors, of unit norm and with positive inner product, associated with the simple target eigenvalue $\lambda(A+\delta E^\cS+\eps E)$, and  $\kappa=1/(x^*y)$. This leads us to the norm-constrained gradient flow, with $G=G_{\delta,\eps}(E^\cS,E)=-xy^*$ for short,
\begin{equation} \label{ode-ES-E}
    \begin{aligned}
               \delta \dot E^\cS &= - \Pi^\cS G + \Re\langle \Pi^\cS G, E^\cS \rangle E^\cS
       \\[1mm]
       \eps \dot E &= -  G + \Re\langle G, E \rangle E,
    \end{aligned}
\end{equation}
where the unit norms of $E^\cS(t)\in\cS$ and $E(t)$ are preserved and $\F_{\delta,\eps}(E^\cS(t),E(t))$ decreases monotonically. Provided that $\Pi^\cS G\ne 0$, we have at stationary points that $E^\cS$ is a real multiple of $\Pi^\cS G$ and $E$ is a real multiple of $G$, which yields \eqref{ESE}. In the case of inequality constraints $\le 1$, the argument of Section~\ref{subsec:ineq} shows that then $E^\cS$ is a negative multiple of $\Pi^\cS G$ and $E$ is a negative multiple of $G$, so that $E^\cS$ is a positive multiple of~$\Pi^\cS E$, and both  $E$ and $E^\cS$ are necessarily of norm 1 in a stationary point with $\Pi^\cS G\ne 0$.
\qed
\end{proof}

Because of Theorem~\ref{thm:ESE} we consider instead of \eqref{F-delta} the minimization of the restricted functional
$$
\wt \F_{\delta,\eps}(E) = - \Re \,\lambda\biggl(A+\eps E +  \delta \frac{\Pi^\cS E}{\| \Pi^\cS E\|_F}\biggr),
$$
for which we follow the programme of Chapter~\ref{chap:proto}.
\begin{lemma}[Reduced gradient] 
\label{lem:red-gradient}
\index{gradient!reduced}
Let $E(t)\in \C^{n, n}$, for real $t$ near $t_0$, be a continuously differentiable path of matrices, with the derivative denoted by $\dot E(t)$.
Assume that $\lambda(t)$ is a simple eigenvalue of  $A+\eps E(t)+\delta \eta(t) \Pi^\cS E(t)$ depending continuously on~$t$, with
$\eta(t)=1/\| \Pi^\cS E(t)\|_F$,  with associated left and right eigenvectors
$x(t)$ and $y(t)$ of unit norm and with positive inner product. Let the eigenvalue condition number be 
$\kappa(t) = 1/(x(t)^* y(t)) > 0$.
%
Then, $\wt \F_{\delta,\eps}(E(t))= -\Re \, \lambda(t)$ 
is continuously differentiable w.r.t. $t$ and we have
\begin{equation} \label{deriv-red}
\frac1{ \kappa(t) } \,\frac{d}{dt} \wt \F_{\delta,\eps}(E(t)) = \Re \,\bigl\langle  \wt G_{\delta,\eps}(E(t)),  \dot E(t) \bigr\rangle,
\end{equation}
where the (rescaled) gradient of $\widetilde \F_{\delta,\eps}$ is the matrix, with $G=- xy^*$ and $\eta=1/\| \Pi^\cS E\|_F$,
\begin{equation} \label{red-grad}
\wt G_{\delta,\eps}(E) =  \eps G + \delta\eta\, \Pi^\cS G - \delta \eta \,\Re\bigl\langle G,\eta\,\Pi^\cS E\bigr\rangle \,\eta\,\Pi^\cS E \in \C^{n, n}.
\end{equation} 
\end{lemma}

\begin{proof}
The proof is similar to the proof of Lemma~\ref{chap:proto}.\ref{lem:gradient}, noting in addition that
$$
\dot \eta = \frac d{dt} \bigl\langle \Pi^\cS E , \Pi^\cS E \bigr\rangle^{-1/2} =
-\eta^3 \, \Re\,\bigl\langle \Pi^\cS E , \Pi^\cS \dot E \bigr\rangle.
$$
We obtain 
\begin{align*}
- \frac1\kappa \,\Re\,\dot \lambda & = 
\Re\, \bigl\langle G, \eps \dot E + \delta \eta\, \Pi^\cS \dot E 
- \delta \eta^3 \,\Re\bigl\langle \Pi^\cS E, \Pi^\cS \dot E \bigr\rangle\, \Pi^\cS E \bigr\rangle
\\
&= \Re\, \bigl\langle G, \eps \dot E + \delta \eta\, \Pi^\cS \dot E \bigr\rangle
- \delta \eta^3 \,\Re\bigl\langle \Pi^\cS E, \Pi^\cS \dot E \bigr\rangle \, 
\Re \bigl\langle G, \Pi^\cS E \bigr\rangle
\\
&= \Re\, \bigl\langle \eps G + \delta \eta\, \Pi^\cS G 
- \delta\eta^3 \,\Re \bigl\langle G, \Pi^\cS E \bigr\rangle \,\Pi^\cS E, \dot E \bigr\rangle
\end{align*}
as stated.
\qed
\end{proof}

We consider the norm-constrained gradient flow
\begin{equation}\label{ode-E-red}
\dot E = -\wt G_{\delta,\eps}(E) + \Re \,\langle \wt G_{\delta,\eps}(E), E \rangle\, E.
\end{equation}
We again have 
$$\frac12\,\frac{d}{dt} \| E(t)\|_F^2 = \Re \langle E(t), \dot E(t) \rangle =0,
$$ 
so that the Frobenius norm $1$ is conserved, and with $\wt G=\wt G_{\delta,\eps}(E)$,
\begin{equation}\label{c-s-red}
\frac1{\kappa} \,\frac{d}{dt} \wt \F_\eps(E(t)) = \Re \langle \wt G, \dot E\rangle = - \| \dot E \|_F^2 =
- \| \wt G - \Re\,\langle \wt G, E \rangle E \|_F^2 \le 0,
\end{equation}
so that the reduced functional decays monotonically along solutions of \eqref{ode-E-red}. At stationary points $E$ we have
that $E$ is a real multiple of $\wt G_{\delta,\eps}(E)$, provided this is nonzero. 
 The following result takes us back to rank-1 matrices.

\begin{lemma}[Stationary points] \label{lem:gradprop}
    Under the conditions of Theorem~\ref{thm:ESE} we have at a stationary point $E$ of \eqref{ode-E-red}
    \begin{equation}
        \text{$\wt G_{\delta,\eps}(E )= \eps G_{\delta,\eps}(E)= - \eps xy^*$},       
    \end{equation}
 where $x$ and $y$ are the left and right eigenvectors to a rightmost eigenvalue of the matrix $A+\eps E+ \delta\eta\,\Pi^\cS E$ with $\eta = 1/\| \Pi^\cS E \|_F$. Hence, $E$ is a real multiple of $xy^*$ and thus has rank 1. 
 Moreover, $E$ is a stationary point of \eqref{ode-E-red} if and only if $(E, \eta\,\Pi^\cS E)$  is a stationary point of \eqref{ode-ES-E}.
\end{lemma}

\begin{proof}
    By Theorem~\ref{thm:ESE}, $\wt G = \wt G_{\delta,\eps}(E)$ is a real multiple of $E$, and with $G=- xy^*$,
    \begin{equation} \label{Pi-G}
    \Pi^\cS \wt G =  \eps \Pi^\cS G + \delta \eta \,\Pi^\cS G -
    \delta \eta \, \Re \langle \Pi^\cS G, \eta\Pi^\cS E \rangle \, \eta\Pi^\cS E
    \end{equation}
is a real multiple of $\Pi^\cS E$. This implies that $\Pi^\cS G$ is a real multiple of $\Pi^\cS E$. So
the last two terms in \eqref{Pi-G} cancel, which are the same two terms as in \eqref{red-grad}. This implies 
$\wt G = \eps G$. Moreover, we then know that $E$ is a real multiple of $\wt G$ and thus of $G$, which implies the result. 
\qed
\end{proof}

\subsection{Rank-1 constrained reduced gradient flow} \label{subsec:r1-ode-red}
\index{gradient flow!rank-1 constrained}
\index{rank-1 matrix differential equation}
We proceed as in Section~\ref{subsec:rank1-gradient-flow} with $\wt G_{\delta,\eps}(E)$ in place of $G_\eps(E)$.
In the differential equation (\ref{ode-E-red}) we project the right-hand side to the tangent space $T_E\cR_1$ with the orthogonal projection $P_E$:
\begin{equation}\label{ode-E-1-red}
\dot E = -P_E\Bigl( \wt G_{\delta,\eps}(E) - \Re \langle \wt G_{\delta,\eps}(E),E \rangle E \Bigr).
\end{equation}
As in (\ref{chap:proto}.\ref{ode-E-1}), solutions $E(t)$ of  \eqref{ode-E-1-red} stay of Frobenius norm 1 for all $t$, and the functional $\wt \F_{\delta,\eps}$ decays monotonically along solutions.
The solution of \eqref{ode-E-1-red} is given as $E(t)=u(t)v(t)^*$, where $u$ and $v$ solve the system of differential equations (with $\wt G= \wt G_{\delta,\eps}(E)$)
\begin{equation}\label{ode-uv-red}
\begin{array}{rcl}
 \dot u &=& -\tfrac \iu2 \, \Im(u^*\wt Gv)u - (I-uu^*)\wt Gv
\\[1mm]
 \dot v &=& -\tfrac \iu2 \, \Im(v^*\wt G^*u)v - (I-vv^*)\wt G^*u,
\end{array}
\end{equation}
which preserves $\|u(t)\|=\|v(t)\|=1$ for all $t$. This system is solved numerically into a stationary point using a splitting method with an Armijo-type stepsize selection, similarly to Section~\ref{subsec:proto-numer}. The following result on the preservation of stationary points is proved by the same arguments as Theorem~\ref{chap:proto}.\ref{thm:stat-1}.
\index{stationary point}

\begin{theorem}[Stationary points]
\label{thm:stat-1-red}
Under the assumptions of Theorem~\ref{thm:ESE}, the following holds true:
\begin{enumerate}
\item Every stationary point $E$ of unit Frobenius norm of \eqref{ode-E-red}, which then is a real multiple of $G_{\eps,\delta}(E)$ and satisfies $\wt G_{\eps,\delta}(E)=  G_{\eps,\delta}(E)$,
is also a stationary point of \eqref{ode-E-1-red}.
\item Conversely, let $E$ of rank 1 and unit Frobenius norm be a stationary point of \eqref{ode-E-1-red} that satisfies 
$\wt G_{\eps,\delta}(E)=  G_{\eps,\delta}(E)$ and $P_E(G_{\eps,\delta}(E))\ne 0$. Then, $E$ is also a stationary point of \eqref{ode-E-red}.
\end{enumerate}
\end{theorem}

\subsection{Outer iteration}
\index{Newton--bisection method}
For the solution of the scalar equation $\phi(\delta,\eps)=0$ for given $\eps>0$ or given $\delta>0$
we use a Newton--bisection method as in
Section~\ref{sec:two-level}.  The partial derivatives of $\phi$ for the Newton iteration are obtained with the arguments of the proof of Theorem~\ref{chap:two-level}.\ref{thm:phi-derivative} (under analogous assumptions), yielding 
\begin{align}
\frac{\partial\phi}{\partial\delta}(\delta,\eps )
&= - \kappa\, \| \Pi^\cS (xy^*) \|_F,  
\\[2mm]
\frac{\partial\phi}{\partial\eps}(\delta,\eps )
&= - \kappa \, \| xy^* \|_F  = - \kappa,
\end{align}
where $x$ and $y$ are left and right normalized eigenvectors with positive inner product, associated with the rightmost eigenvalue of
$A+\delta E^\cS_{\delta,\eps}+\eps E_{\delta,\eps}$, and $\kappa=1/(x^*y)>0$.

\medskip
Numerical experiments with the algorithm described in this section are reported by Guglielmi \& Lubich (\cite{GL25}) for Toeplitz and sparse real matrices.

\section{Numerical range under structured perturbations}
\label{sec:num-range}
\index{numerical range} \index{field of values}
Given a matrix that has the numerical range in the open complex left half-plane, we aim to compute a structured perturbation
of minimal Frobenius norm such that the numerical range of the perturbed matrix touches the imaginary axis. In contrast to the algorithms considered so far, the iterative algorithm proposed here will not require eigenvalue and eigenvector computations (except possibly for the starting iterate) but instead work with Rayleigh quotients. 
\index{eigenvalue optimization!without eigenvalues}

\subsection{Numerical range and dissipativity}
The {\it numerical range} (or {\it field of values}) of a matrix $A\in\C^{n,n}$ is the range of the Rayleigh quotients, 
\index{Rayleigh quotient}
\begin{align*}
W(A) &= \biggl\{ \frac{v^* A v}{v^*v}\,:\, v\in \C^n, v\ne 0 \biggr\} 
\\
&= \bigl\{  v^* A v \,:\, v\in \C^n, \| v \|=1 \},
\end{align*}
where $\|\cdot\|$ is the Euclidean norm. 
\bcl The set $W(A)$ is convex and contains the eigenvalues of $A$. For a comparison of numerical range and pseudospectra we refer to Trefethen \& Embree (\cite{TreE05}), Chapter 17.
\ecl

The real part of a rightmost point of $W(A)$ is the {\it numerical abscissa},
\begin{align*}
    \omega(A) &= \max \{ \Re \,z\,:\, z \in W(A) \}
    \\
    &= \text{largest eigenvalue of }\ \tfrac12 (A+A^*).
\end{align*}
The interest in this number relies on the following fact.
\begin{lemma}
    \label{lem:exp-num-absc}
    The numerical abscissa $\omega(A)$ is the smallest real number $\omega$ such that
$$   
\| e^{tA} \|_2 \le e^{t\omega} \quad\text{ for all }\ t\ge 0.
$$
\end{lemma}
\begin{proof}
    Replacing $A$ by $A-\omega(A)I$, we can assume $\omega(A)=0$ without loss of generality. Let $y(t)$ be a solution of the differential equation $\dot y=Ay$, i.e., $y(t)= e^{tA}y(0)$. Then,
    $$
    \frac d{dt}\|y(t)\|^2 = 2\, \Re \langle y,\dot y \rangle = 2 \, \Re \langle y,A y \rangle \le 2\omega(A) \|y\|^2=0,
    $$
    and hence $\|y(t)\|^2 \le \|y(0)\|^2$, that is, $\|e^{tA}y(0)\| \le \|y(0)\|$. Since this is valid for every initial value $y(0)\in \C^n$, we find $\| e^{tA} \|_2 \le 1$, which is the stated bound for $\omega(A)=0$.
    On the other hand, if $y(0)$ is an eigenvector to the eigenvalue $\omega(A)=0$ of $\frac12(A+A^*)$, then the above bound for $\frac d{dt}\|y(t)\|^2$ becomes an equality at $t=0$, and hence there cannot be an $\omega < \omega(A)=0$ such that the inequality of the lemma holds for small $t>0$.
    \qed
\end{proof}

If $\omega(A)\le 0$, the matrix $A$ is called {\em dissipative}.  
\index{dissipative matrix}
Then, $\| e^{tA} \|_2 \le 1$ for all $t\ge 0$.

\subsection{Structured dissipativity radius}
\index{structured dissipativity radius}
\index{dissipativity radius}
We consider the behaviour of the numerical abscissa under structured perturbations (as we did previously in this chapter for the spectral abscissa and pseudospectral abscissa). Given a matrix $A\in\C^{n,n}$, a complex-linear or real-linear structure space $\cS\subset\C^{n,n}$
and a perturbation size $\eps$, we define the {\it $\cS$-structured $\eps$-numerical abscissa}
$$
\omega_\eps^\cS(A) =\max \{ \Re\, z\,:\, z\in W(A+\Delta) \text{ for some $\Delta\in\cS$ with $\|\Delta\|_F \le \eps$} \}.
$$
In the unstructured complex case $\cS=\C^{n,n}$ we simply have $\omega_\eps(A)=\omega(A)+\eps$, and the same holds true in the unstructured real case $\cS=\R^{n,n}$ if also $A$ is real. The extremal perturbation in these cases is $\Delta=\eps\, xx^*$, where
$x$ is an eigenvector to the largest eigenvalue of $\tfrac12(A+A^*)$. For more general structures, however, determining 
$\omega_\eps^\cS(A)$ is not so trivial.

For a matrix $A$ with negative numerical abscissa $\omega(A)<0$ we define the {\it $\cS$-structured dissipativity radius}
\begin{align*}
 \eps_{\text{diss}}^\cS(A) &= \text{smallest $\eps$ such that } \omega_\eps^\cS(A)\ge 0
 \\
  &=  \text{largest $\eps$ for which $\|e^{t(A+\Delta)}\|_2 \le 1$ for all $t\ge 0$ and all}
  \\
  &\quad\ \, \text{structured perturbations $\Delta \in \cS$ of Frobenius norm at most $\eps$.}
\end{align*}
Here the second equality follows from Lemma~\ref{lem:exp-num-absc}.
In the following we will address the computation of this quantity.

\medskip\noindent
{\bf Problem.} {\it Given a matrix $A\in \C^{n,n}$ with negative numerical abscissa and a structure space $\cS$,
compute the $\cS$-structured dissipativity radius $\eps_{\text{\rm diss}}^\cS(A)$.}

\medskip\noindent
We will describe a two-level algorithm for computing the $\cS$-structured dissipativity radius and the
extremal perturbation $\Delta$ of minimal Frobenius norm that has $\omega(A+\Delta)=0$. This algorithm uses a single vector differential equation in the inner iteration and a Newton--bisection method in the outer iteration. While this fits into our framework for dealing with matrix nearness problems considered so far, the decisive difference is that we do not solve eigenvalue problems in each iteration, but instead we work with Rayleigh quotients, which just require matrix--vector products and inner products.

\subsection{Optimizing Rayleigh quotients instead of eigenvalues}
\label{subsec:rayleigh-max}
We aim to minimize the following functional, which contains a Rayleigh quotient instead of eigenvalues:
$$
\F_\eps(E,v) = - \Re\, v^*(A+\eps E)v
$$
is to be minimized over all $E\in\cS$ of Frobenius norm 1 and $v\in\C^n$ of Euclidean norm 1. In the real case where $A\in \R^{n,n}$ and $\cS\subset\R^{n,n}$, we only need to take real vectors $v\in \R^n$.

In the present Hermitian case, our approach with eigenvalues and eigenvectors as considered so far can be interpreted as computing
$$
\min_{(E,v)} \F_\eps(E,v) = \min_E \Bigl( \min_v \F_\eps(E,v) \Bigr),
$$
where the inner minimum over $v$ is given by the negative real part of the rightmost eigenvalue. Now we will first consider the simultaneous minimization of $E$ and $v$, and from the optimality condition we find that in the minimum, $E$ can be expressed as a computationally simple function $E(v)$, and $v$ is indeed an eigenvector of $A+\eps E$. Inserting $E(v)$ into the functional then yields a reduced functional of $v$ only,
$\wt \F_\eps(v)=\F_\eps(E(v),v)$, which we minimize via the (norm-constrained) reduced gradient flow for $v$.

We organize this subsection similarly to Section~\ref{subsec:reduced-opt}.

\subsubsection*{Structured gradient.}
\index{gradient!structured}
Consider a smooth path of structured matrices $E(t)\in\cS$ and vectors $v(t)\in\C^n$. Since then also $\dot E(t)\in\cS$, we have, using that $\dot E(t)= \Pi^\cS \dot E(t) \in \cS$ and denoting the Hermitian part of a matrix $M$ by 
$\text{Herm}(M)=\tfrac12(M+M^*)$
\begin{align*}
   \frac d{dt}\, \F_\eps(E(t),v(t)) &=  - \eps\, \Re\, (v^*\dot E v) -2\, \Re\bigl( v^* \text{Herm}(A+\eps E)\dot v\bigr)
   \\
   &= - \eps \,\Re\, \langle \Pi^\cS(vv^*), \dot E \rangle - 2\,\Re\,\langle \text{Herm}(A+\eps E)v,\dot v \rangle.
\end{align*}

\subsubsection*{Norm- and structure-constrained gradient flow.}
\index{gradient flow!norm-constrained}
\index{gradient flow!structure-constrained}
With a scaling factor $\alpha$, e.g. $\alpha=|\omega(A)|$, we consider the system of differential equations
\begin{align}
    \dot E &= \Pi^\cS(vv^*) - \Re\,\langle \Pi^\cS(vv^*), E \rangle \, E
    \nonumber
    \\
    \alpha\, \dot v &= \text{Herm}(A+\eps E) v - \langle \text{Herm}(A+\eps E) v,v \rangle\, v.
    \label{v-dot-num}
\end{align}
By construction, the Frobenius norm 1 of $E(t)$ and the Euclidean norm 1 of $v(t)$ are preserved, and again we obtain that the functional $\F_{\eps}$ decays monotonically along solutions.

\subsubsection*{Stationary points.}
\index{stationary point}
At stationary points with $\Pi^\cS(vv^*)\ne 0$ we find that
$$
\text{$E$ is a real multiple of $\Pi^\cS(vv^*)$, and $v$ is an eigenvector of $\text{Herm}(A+\eps E)$.}
$$
    Moreover, if the equality constraints for the Frobenius norm $1$ are replaced by inequality constraints $\le 1$, then 
    the same argument as in the proof of Theorem~\ref{chap:proto}.\ref{ESE} shows that
    $E$ is a positive multiple of $\Pi^\cS(vv^*)$, i.e. 
\begin{equation}\label{E-opt-num}
E= \Pi^\cS (vv^*)/\| \Pi^\cS(vv^*) \|_F.
\end{equation}

\subsubsection*{Reduced functional.}
Similar to Section~\ref{sec:eps}, this motivates us to insert  \eqref{E-opt-num} into the functional and to minimize the resulting reduced functional of a single vector $v\in\C^n$ of norm 1 only:
$$
\wt \F_\eps(v) = - \Re \bigl( v^*(A+\eps E)v\bigr) \quad\text{ with }\ E= \frac{\Pi^\cS(vv^*)}{\|\Pi^\cS(vv^*)\|_F}.
$$
\begin{lemma}[Reduced gradient]
\label{lem:red-grad-num}
\index{gradient!reduced}
Along a path $v(t)\in\C^n$, we have (omitting the argument $t$ on the right-hand side)
$$
\frac{d}{dt}\, \wt \F_\eps(v(t)) = - 2\,\Re\,\langle \text{\rm Herm}(A+\eps E)v,\dot v \rangle,
$$
where $E= \Pi^\cS(vv^*)/\|\Pi^\cS(vv^*)\|_F$.
\end{lemma}
\begin{proof}
We find
$$
\frac{d}{dt}\, \Re (v^* A  v) = 2\, \Re\,\langle \text{Herm}(A)v,\dot v \rangle
$$
and with $\eta=1/\|\Pi^\cS(vv^*)\|_F$,
\begin{align*}
\frac{d}{dt}\,\bigl(v^* \eta \Pi^\cS(vv^*) v \bigr) 
= \frac{d}{dt}\,\bigl( \eta \langle vv^*, \Pi^\cS(vv^*)\rangle \bigr) 
= \frac{d}{dt}\,\bigl( \eta \| \Pi^\cS(vv^*)\|_F^2 \bigr)
= \frac{d}{dt}\,\frac 1\eta = - \frac{\dot\eta}{\eta^2}
\end{align*}
and
\begin{align*}
\dot \eta &= -\tfrac12 \eta^3 \, 2\,\Re \,\langle \Pi^\cS(vv^*), \dot v v^* + v \dot v^* \rangle
= -\eta^3 \,\Re \langle 2 \,\text{Herm}(\Pi^\cS(vv^*))v, \dot v \rangle
\end{align*}
so that 
$$
\frac{d}{dt}\,\bigl( v^* \eta\Pi^\cS(vv^*) v \bigr) 
= - \frac{\dot\eta}{\eta^2}= 2\,\Re \langle \text{Herm}(\eta\Pi^\cS(vv^*))v, \dot v \rangle
= 2\,\Re \langle \text{Herm}(E)v,\dot v \rangle,
$$
which yields the result.
\qed
\end{proof}

\subsubsection*{Norm-constrained reduced gradient flow.}
In view of Lemma~\ref{lem:red-grad-num}, the norm-constrained gradient flow of the reduced functional takes the
same form as the differential equation for $v$ in \eqref{v-dot-num}:
\begin{equation}\label{v-dot-red-num}
    \dot v = \text{Herm}(A+\eps E) v - \langle \text{Herm}(A+\eps E) v,v \rangle\, v,
\end{equation}
now for $E= \Pi^\cS(vv^*/\|\Pi^\cS(vv^*)\|_F$. 
The Euclidean norm 1 of $v(t)$ is preserved for all $t$, and the reduced functional $\wt \F_{\eps}$ decays monotonically along solutions of this differential equation.

We  numerically integrate the differential equation \eqref{v-dot-red-num} for the single vector $v(t)$ into a stationary point. Note that this algorithm just requires computing matrix--vector products and inner products of vectors but no computations of eigenvalues and eigenvectors.

In this way we aim to determine
$$
v(\eps)=\arg\min_{v} \;\wt \F_\eps(v) \quad\text{and}\quad E(\eps)= \Pi^\cS(v(\eps)v(\eps)^*) /
\| \Pi^\cS(v(\eps)v(\eps)^*) \|_F,
$$
where the minimum is taken over all $v\in \C^n$ of Euclidean norm 1.

\subsubsection*{Outer iteration.}
The optimal perturbation size $\eps_\text{diss}^\cS$ is determined from the equation
\begin{equation} \label{zero-num}
\Re \Bigl(v(\eps)^*\bigl(A+\eps E(\eps)\bigr)v(\eps)\Bigr) =0.
\end{equation}
\index{Newton--bisection method}
This scalar equation is solved by the Newton--bisection algorithm of Section~\ref{subsec:Newton--bisection}, where the derivative formula of Theorem~\ref{chap:two-level}.\ref{thm:phi-derivative} is valid also for the structure-projected gradient $G(\eps)=-\Pi^\cS ( v(\eps)v(\eps)^*)$ used here.

\subsubsection*{Illustrative example.}

We consider the matrix $A = - I_{10}$, the identity matrix of
dimension $n=10$, and the structure space $\cS$ of matrices with prescribed range and co-range,
\[
\cS = \{ B \Delta C \,:\, \Delta \in \R^{2,3} \},
\]
with
\begin{align*}
& B^\top = \left( \begin{array}{rrrrrrrrrr}
    -2  &  -1  &   1  &  -4  &   0  &  -2  &   2  &  -2  &   0  &  -2 \\
    -1  &  -1  &   0  &  -1  &   2  &   1  &   4  &  -1  &   2  &  -4 
    \end{array} \right)
\\[2mm]    
& C =   \left( \begin{array}{rrrrrrrrrr}
     1  &   1  &  -1  &   1  &   2  &   2  &   1  &   1  &  -1  &  -2 \\
    -1  &   0  &   2  &  -2  &  -1  &   6  &  -1  &   4  &   2  &  -1 \\
     0  &   0  &   1  &   3  &  -3  &   4  &   2  &  -1  &   0  &   1
    \end{array} \right)  
\end{align*}
We set $\eps=1$ and maximize the real part of the rightmost point in the numerical range of $A+\eps E$ both with the method proposed in this section using Rayleigh quotients and in the way considered previously in this book, maximizing the real part of the rightmost eigenvalue of $A+\eps E$ (for which
it is necessary to compute the rightmost eigenpair of the perturbed matrix at each step of the
numerical integrator).

\begin{figure}[t]
\vskip -3cm
\centering
\includegraphics[scale=0.5]{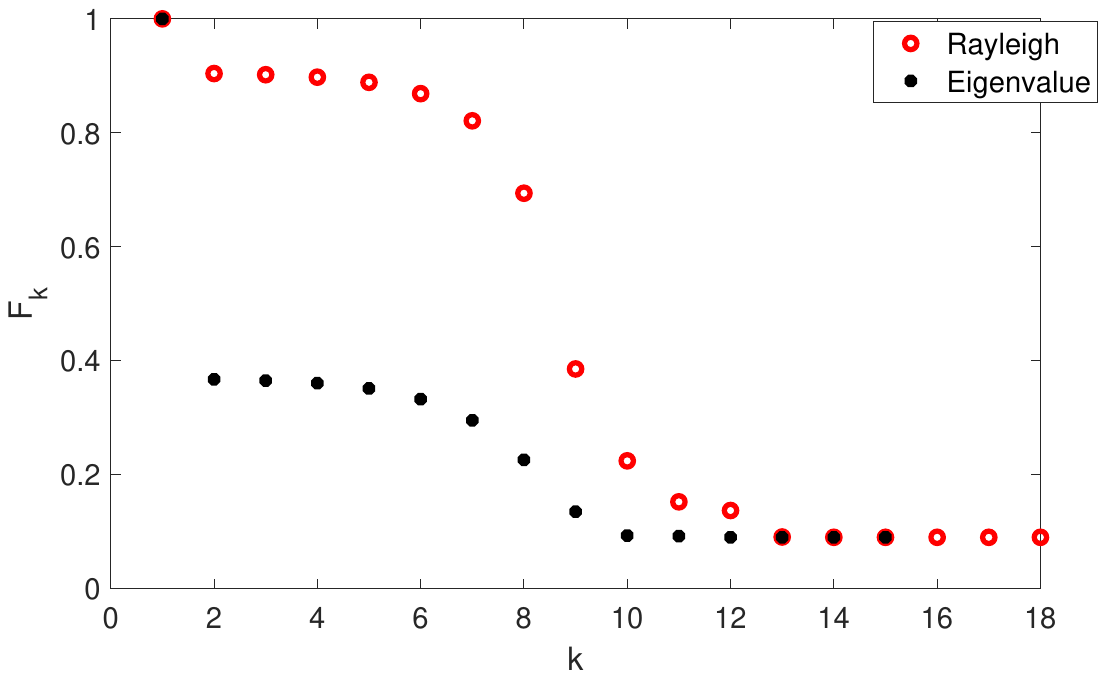}
\vskip -3.7cm
\caption{Behaviour of the algorithm based on the Rayleigh quotient (circles) versus the algorithm based on the eigenvalue computation (bullets).} 
\label{fig:fov}
\end{figure}

Figure \ref{fig:fov} shows the behaviour of the eigenvalue versus that of the Rayleigh quotient.
In the figure, $\F_k$ denotes the value of the functional at the $k$-th step, that is
$\F_k \approx \wt \F_\eps(v(t_k))$.

We note that the final value of $\wt \F_\eps(v) \approx 0.08976$; the difference with respect to the corresponding approximation computed by the method making use of eigenvalues is roughly $5 \times 10^{-7}$
implying an error on the sixth significant digit. The number of steps executed by the numerical integrator is comparable ($18$ for the Rayleigh quotient and $15$ for the eigenvalue method). The structured dissipativity radius of $A=-I$ is approximately computed as  $\oeps \approx 1.098614111289907$, whereas the unstructured dissipativity radius equals 1.

\begin{table}[hbt]
\begin{center}
\begin{tabular}{|l|l|l|l|}\hline
  $k$ & $\eps_k$ & $ -\omega(A+\eps_k E(\eps_k))$ & $\#$ Rayleigh \\
 \hline
\rule{0pt}{9pt}
\!\!\!\! 
        $0$         & $1$ & $1.856586036644719$ & $18$  \\
	$1$         & $1.098615122565759$ & $0.920501421999773 \cdot 10^{-6}$ & $13$  \\
	$2$         & $1.098614111289907$ & $6.104818325601980 \cdot 10^{-12}$ & $10$  \\
 \hline
\end{tabular}
\vspace{2mm}
\caption{Computation of the dissipativity radius: 
computed values $\eps_k$, $-\omega(A+\eps_k E(\eps_k))$
and number of Rayleigh quotient computations of the inner algorithm.\label{tab:Hur-S}}
\end{center}
\end{table}

\section{Structured distance to singularity}
\label{sec:sing-S}
\index{structured distance to singularity}
\index{distance to singularity!under structured perturbations}

Let $A\in \C^{n,n}$ be an invertible matrix, and let $\cS$ be a complex- or real-linear subspace of $\C^{n,n}$ that defines the linear structure that is imposed on perturbations $\Delta\in\cS$ to $A$. We consider the following structured matrix nearness problem.

\medskip\noindent
{\bf Problem.} {\it Given an invertible matrix $A$, find a structured perturbation $\Delta\in \cS$ of minimal Frobenius norm such that $A+\Delta$ is singular.}

\medskip\noindent
The norm of the minimizing $\Delta\in \cS$ is called the {\em $\cS$-structured distance to singularity} of the given matrix $A$. Unlike the complex or real unstructured distance to singularity, the $\cS$-structured distance to singularity (for $\cS\ne \C^{n,n},\R^{n,n}$) cannot be obtained from a singular value decomposition of $A$.

An obvious case of interest is when $\cS$ is a space of complex or real matrices with a prescribed sparsity pattern. At the end of this section we consider a different situation where a structured distance to singularity is of interest: Given two polynomials, compute the nearest pair of polynomials with a common zero. This can be rephrased as the above problem with $\cS$ the space of Sylvester matrices.

\subsection{Optimizing eigenvalues or singular values}
\index{two-level iteration}
For the computation of the structured distance to singularity we can use the two-level iteration of Section~\ref{sec:two-level}. In the inner iteration we could again do eigenvalue optimization using the rank-1 algorithm of Section~\ref{sec:proto-structured} with $f(\lambda,\clambda)=\lambda\clambda =|\lambda|^2$, choosing the eigenvalue of smallest modulus as the target eigenvalue (the smallest singular value could be used in an alternative algorithm). 
We note that here the structured gradient is $G_\eps^\cS(E)=\Pi^\cS(\lambda x y^*)$, for the (simple) eigenvalue $\lambda=\lambda(A+\eps E)$ with left and right eigenvectors $x$ and $y$ of norm 1 and with positive inner product. Alternatively, the smallest singular value could be minimized under structured perturbations of Frobenius norm~$\eps$.
The outer iteration for finding the minimal $\eps$ yielding a zero eigenvalue or singular value uses the Newton--bisection method described in Section~\ref{subsec:Newton--bisection}.

\subsection{Optimizing norms of matrix--vector products}
\label{subsec:mvp-sing-S}
\index{eigenvalue optimization!without eigenvalues}
We now turn to a different approach that requires only matrix--vector products instead of computations of eigenvalues and eigenvectors or singular values and singular vectors.
To avoid some technical complications that would obscure the basic argument, we consider in this subsection the real case:
$$
\text{$A\in \R^{n,n}$ is a real matrix and $\cS\subset \R^{n,n}$ is a real structure space.}
$$
In an approach that does not require the computation of eigenvalues and eigenvectors (or singular values and singular vectors) in every iteration,
we first consider the following functional, which only contains the Euclidean norms of matrix--vector products:
\begin{equation}\label{F-sing}
\F_\eps(E,u,v) = \tfrac12\,\|u^{\top}(A+\eps E)\|^2 + \tfrac12\,\|(A+\eps E)v\|^2
\end{equation}
is to be minimized over all $E\in\cS$ of Frobenius norm 1 and all $u,v\in\R^n$ of Euclidean norm~1. 
We let $(E(\eps),u(\eps),v(\eps))$ be a minimizer of $\F_\eps$ and 
$\oeps$ is the smallest $\eps>0$ such that $\F_\eps(E(\eps),u(\eps),v(\eps))=\min_{(E,u,v)} \F_\eps(E,u,v)=0$.

It will turn out that in a minimum, $u$ and $v$ are left and right singular vectors of $A+\eps E$, and $E=-\Pi^\cS(uv^{\top})/\|\Pi^\cS(uv^{\top})\|$. We will then restrict the functional $\F_\eps$ to such $E$, obtaining a functional of two vectors $u$ and $v$ only. This reduced functional $\wt \F_\eps(u,v)$ will be used in the two-level approach, minimizing $\wt \F_\eps$ for a fixed $\eps$ in the inner iteration and determining the smallest $\eps$ with $\wt \F_\eps(u(\eps),v(\eps))=0$ in the outer iteration.
We closely follow the program of Section~\ref{subsec:rayleigh-max}.

\subsubsection*{Structured gradient.}
\index{gradient!structured}
Consider a smooth path of structured matrices $E(t)\in\cS$ and vectors $u(t),v(t)\in\C^n$. Since then also $\dot E(t)\in\cS$, we have
\begin{align}
   &\frac d{dt}\, \F_\eps(E(t),u(t),v(t)) 
   \nonumber
   \\
   &=  
    \,\langle u^{\top}(A+\eps E), \eps u^{\top}\dot E + \dot u^{\top}(A+\eps E) \rangle +
    \,\langle (A+\eps E)v, \eps \dot E v + (A+\eps E)\dot v \rangle
    \nonumber
   \\
   &=  \eps \, \langle G, \dot E \rangle +
       \,\langle (A+\eps E)(A+\eps E)^{\top}u,\dot u \rangle +
    \,\langle (A+\eps E)^{\top}(A+\eps E)v,\dot v \rangle
        \label{grad-full-sing}
\end{align}
with the structured gradient matrix
$$
G= G_\eps^\cS(E,u,v)=\Pi^\cS\bigl(uu^{\top}(A+\eps E)+(A+\eps E)vv^{\top}\bigr),
$$
which is the projection onto $\cS$ of a matrix of rank at most 2.

\subsubsection*{Norm- and structure-constrained gradient flow.} 
\index{gradient flow!norm-constrained}
\index{gradient flow!structure-constrained}

Taking the norm constraints into account, we arrive at the following system of differential equations:
\begin{align}
\nonumber
    \eps \dot E &= -G +  \,\langle G, E \rangle \, E
    \\
    \label{Euv-ode-S}
    \alpha \, \dot u &= -(A+\eps E)(A+\eps E)^{\top} u + \|(A+\eps E)^{\top} u\|^2\, u
    \\
    \nonumber
    \beta \,\dot v &= -(A+\eps E)^{\top}(A+\eps E) v + \|(A+\eps E) v\|^2\, v,
\end{align}
with scaling factors $\alpha(t)>0$ and $\beta(t)>0$, which we propose to choose as
$\alpha = \|(A+\eps E)^{\top} u\|$ and $\beta=\|(A+\eps E) v\|$ to obtain similar rates of change in $\eps E$, $u$, $v$ near a stationary point.

By construction, the Frobenius norm 1 of $E(t)$ and the Euclidean norm 1 of $u(t)$ and $v(t)$ are preserved, and the functional $\F_{\eps}$ decays monotonically along solutions.

\subsubsection*{Stationary points.}
\index{stationary point}
At stationary points $(E,u,v)$, the vectors $u$ and $v$ are eigenvectors of $(A+\eps E)(A+\eps E)^{\top}$ and $(A+\eps E)^{\top}(A+\eps E)$, respectively. If the sign of $u$ is such that $u^{\top}(A+\eps E)v$ is real and positive, then $u$ and $v$ are left and right singular vectors of $A+\eps E$.
The case of interest is when $u$ and $v$ are left and right singular vectors to the common smallest singular value $\sigma>0$. We then have $\sigma u=(A+\eps E)v$ and
$\sigma v = (A+\eps E)^{\top}u$ and hence a rescaled gradient (rescaled by the factor $(2\sigma)^{-1})$ becomes $G=\Pi^\cS(uv^{\top})$ in the stationary point. This rescaling corresponds to taking the norms instead of the squared norms in the functional $\F_\eps$ of \eqref{F-sing}. If $G\ne 0$, it follows that $E$ is a multiple of  $\Pi^\cS(uv^{\top})$. Moreover, if the equality constraint for the Frobenius norm $1$ is replaced by an inequality constraint $\le 1$, then the same argument as in the proof of Theorem~\ref{chap:proto}.\ref{ESE} shows that $E$ is a {\em negative} multiple of $\Pi^\cS(uv^{\top})$ and still has norm 1, so that 
\begin{equation}\label{E-opt-sing}
E= -\frac{\Pi^\cS(uv^{\top})}{\|\Pi^\cS(uv^{\top})\|_F}.
\end{equation} 

\begin{lemma}[Non-vanishing gradient when $A+\eps E$ is nonsingular] \label{lem:nvg}
    If $A\in\cS$ and $\sigma\ne 0$ is a singular value of $A+\eps E \in\cS$ with associated left and right singular vectors $u$ and $v$, respectively, then $\Pi^\cS(uv^{\top})\ne 0$.
\end{lemma}

\begin{proof}
    The proof uses a similar argument to the proof of Theorem~\ref{thm:nonzero-gradient-S}. We have
    \begin{align*}
    & \langle \Pi^\cS(uv^{\top}),A+\eps E\rangle 
    \\
    &=  \langle uv^{\top},\Pi^\cS(A+\eps E)\rangle = 
     \langle uv^{\top},A+\eps E\rangle 
    \\
    &= \langle u,(A+\eps E)v\rangle =  \langle u,\sigma u\rangle = \sigma \ne 0
    \end{align*}
    and hence $\Pi^\cS(uv^{\top})\ne 0$.
    \qed
\end{proof}

The next lemma characterizes the non-vanishing of the gradient when $\sigma=0$.

\begin{lemma}[Non-vanishing gradient when $A+\eps E$ is singular] \label{lem:nvg-0}
Let $u$ and $v$ be two nonzero vectors. We have
 $\Pi^\cS(uv^{\top})\ne 0$ if and only if there exists a matrix $S\in\cS$ such that $u^\top S v \ne 0$.
\end{lemma}

\begin{proof}
    Let $S\in\cS$ be arbitrary. We have
    \begin{align*}
     \langle \Pi^\cS(uv^{\top}),S\rangle 
    =  \langle uv^{\top},\Pi^\cS(S)\rangle = 
     \langle uv^{\top},S\rangle 
    = \langle u,Sv\rangle = u^\top S v,
    \end{align*}
    and $\Pi^\cS(uv^{\top})\ne 0$ if and only if $\langle \Pi^\cS(uv^{\top}),S\rangle \ne 0$ for some $S\in\cS$.
    \qed
\end{proof}

As a case of particular interest, consider the case where the identity matrix $I\in\cS$, which occurs in many applications. In the situation described before the two lemmas above, in the case when the smallest singular value of $A+\eps E$ is $\sigma=0$, then the associated singular vectors $u$ and $v$ are also eigenvectors to the eigenvalue $0$.  Lemma~\ref{lem:nvg-0} shows the following:
$$
\text{\it If $I\in\cS$, then $\Pi^\cS(uv^{\top})= 0$ only if $u^Tv=0$, that is, if\/ $0$ is a defective eigenvalue.}
$$

The above arguments motivate and yield the following theorem.

\begin{theorem} [Rank-1 property of the minimal deregularizing perturbation]
\label{thm:sing-S}
    Let $A\in \cS$ be invertible and let $\Delta\in\cS$ be a perturbation of minimal Frobenius norm $\oeps$ that makes $A+\Delta$ singular. Assume that the nullspace of $A+\Delta$ is one-dimensional. Let $u$ and $v$ be left and right singular vectors to the simple singular value $0$ of $A+\Delta$. 
    Assume that $\Pi^\cS(uv^{\top})\ne 0$. 
    Then, $\Delta$ is a real multiple of $\Pi^\cS(uv^{\top})$. In particular, $\Delta$ is the orthogonal projection onto the structure space $\cS$ of a rank-1 matrix.
\end{theorem}

\begin{proof} We write $\Delta=\oeps E$ with $E$ of Frobenius norm~1. Then $(E,u,v)$ is a stationary point of
\eqref{Euv-ode-S}, and by the arguments given in the paragraph before Lemma~\ref{lem:nvg}, we have \eqref{E-opt-sing}, 
%
%
which shows that $\Delta=\oeps E$ is a real multiple of $\Pi^S(uv^\top)$.
\qed
\end{proof}

\begin{example}[Noferini \& Poloni's counterexample when $\Pi^\cS(uv^{\top}) = 0$] The assumption $\Pi^\cS(uv^{\top})\ne 0$ cannot be dropped. This was pointed out to us by Vanni Noferini and Federico Poloni, who give the following counterexample:
$$
\cS = \biggl\{ 
\begin{pmatrix}
    p & q \\
    0 & p
\end{pmatrix}\in \R^{2,2}\,:\,p,q \in \R 
\biggr\}, 
\qquad 
A=\begin{pmatrix}
    1 & 1 \\
    0 & 1
\end{pmatrix} \in\cS.
$$
The minimal singularizing perturbation is $\Delta=-I_2 \in\cS$, which is not  the orthogonal projection onto $\cS$ of a rank-1 matrix. The singular vectors of $A+\Delta$ are $u=(0,1)^\top$ and $v=(1, 0)^\top$, for which
$\Pi^\cS(uv^{\top})= 0$ and $v^\top u=0$ in agreement with the defectivity of the eigenvalue $0$ of $A+\Delta=A-I_2$. 

We further remark that for $\eps<\oeps=\sqrt{2}$, the perturbation that minimizes $\F_\eps$ is upper triangular with $\Pi^\cS(u(\eps)v(\eps)^{\top})\ne 0$ in agreement with Lemma~\ref{lem:nvg}. 
\end{example}

\begin{example}[Complementary example when $\Pi^\cS(uv^{\top}) = 0$] 
Consider
\[
\cS = \biggl\{ 
\begin{pmatrix}
    0 & q \\
    q & p
\end{pmatrix}\in \R^{2,2}\,:\,p,q \in \R 
\biggr\}, 
\qquad 
A=\begin{pmatrix}
    0 & 1/2 \\
    1/2 & 1
\end{pmatrix} \in\cS.
\]
The minimal singularizing perturbation is 
\[
\Delta ={}-\begin{pmatrix}
    0 & 1/2 \\
    1/2 & 0
\end{pmatrix} \in\cS,
\]
which is not  the orthogonal projection onto $\cS$ of a rank-1 matrix. 
The singular vectors of $A+\Delta$ are $u=(1,0)^\top$ and $v=(1, 0)^\top$, for which $\Pi^\cS(uv^{\top})= 0$. 
Here, the eigenvalue $0$ of $A+\Delta$ is not defective, but the identity matrix $I_2 \not\in \cS$.
\end{example}



\subsubsection*{Reduced functional.}
We insert  \eqref{E-opt-sing} into the functional \eqref{F-sing} and minimize the resulting reduced functional of two vectors $u,v\in\C^n$ of norm 1 only:
$$
\wt \F_\eps(u,v) = \tfrac12\,\|u^{\top}(A+\eps E)\|^2 + \tfrac12\,\|(A+\eps E)v\|^2 \quad\text{ with }\ 
E= -\frac{\Pi^\cS(uv^{\top})}{\|\Pi^\cS(uv^{\top})\|_F}.
$$
Here we assume that the signs of $u$ and $v$ are such that $u^\top (A+\eps E)v>0$. We will not discuss exceptional situations where this term becomes zero or where the denominator in the expression for $E$ becomes zero.
\begin{lemma}[Reduced gradient]
\label{lem:red-grad-sing}
Along a path $(u(t),v(t))\in\R^n\times \R^n$, we have (omitting the argument $t$) 
$$
\frac{d}{dt}\, \wt \F_\eps(u(t),v(t)) = \langle g,\dot u \rangle + \langle h,\dot v \rangle
$$
with
\begin{align*}
g &= (A+\eps E)(A+\eps E)^\top u - \eps \eta^2 \gamma Ev - \eps \eta Gv
\\ 
h &=  (A+\eps E)^\top(A+\eps E) v - \eps \eta^2 \gamma E^\top u - \eps \eta G^\top u,
\end{align*}
where 
$E=-\eta\Pi^\cS(uv^{\top})$ and $\eta=1/\|\Pi^\cS(uv^{\top})\|_F$, and \\
$G=\Pi^\cS\bigl(uu^{\top}(A+\eps E)+(A+\eps E)vv^{\top}\bigr)$ and $\gamma = u^\top G v$.
\end{lemma}
\index{gradient!reduced}

\begin{proof}
We have
\begin{align*}
    \dot E &= - \dot \eta \, \Pi^\cS(uv^\top) - \eta\, \Pi^\cS(\dot uv^\top+ u\dot v^\top) \quad\text{with }
 \\
  \dot\eta &= - \eta^3 \langle \Pi^\cS(uv^\top), \dot uv^\top+ u\dot v^\top\rangle 
= \eta^2 \langle Ev,\dot u\rangle + \eta^2 \langle E^\top u, \dot v \rangle.
\end{align*}
This yields, using that $G\in\cS$,
\begin{align*}
\langle G,\dot E\rangle &= \eta^3 \langle \Pi^\cS(uv^\top), \dot uv^\top+ u\dot v^\top\rangle \,\gamma
- \eta \langle G,\dot uv^\top+ u\dot v^\top \rangle
\\
&= \langle -\eta^2 \gamma Ev  - \eta Gv, \dot u \rangle + 
\langle -\eta^2 \gamma E^\top u  - \eta G^\top u, \dot v \rangle.
\end{align*}
Inserting this expression into \eqref{grad-full-sing} yields the stated result.
\qed    
\end{proof}

\subsubsection*{Norm-constrained reduced gradient flow.}
\index{gradient flow!norm-constrained}
In view of Lemma~\ref{lem:red-grad-sing}, the norm-constrained gradient flow of the reduced functional takes the
form 
\begin{equation}
\label{uv-ode-sing}
\begin{array}{rcl}
\dot u &=& - g + \langle g,u \rangle u
\\[1mm]
\dot v &=& - h + \langle h,v \rangle v.
\end{array}
\end{equation}
The Euclidean norm 1 of $u(t)$ and $v(t)$ is preserved for all $t$, and the reduced functional $\wt \F_{\eps}$ decays monotonically along solutions of this differential equation.

We  numerically integrate the differential equation \eqref{uv-ode-sing} for the two vectors $u(t)$ and $v(t)$ into a stationary point. In each time step, this algorithm just requires computing matrix--vector products and inner products of vectors but no computations of eigenvalues and eigenvectors or singular values and singular vectors.
The sign of $u$ in the initial value is chosen such that $u^\top (A+\eps E)v>0$, and this sign condition is monitored over the discrete trajectories.

In this way we aim to determine
$$
(u(\eps),v(\eps))=\arg\min_{u,v} \;\wt \F_\eps(u,v) \quad  \text{and} \quad 
E(\eps)= -\frac{\Pi^\cS(u(\eps)v(\eps)^\top)}
{\| \Pi^\cS(u(\eps)v(\eps)^\top) \|_F},
$$
where the minimum is taken over all $u,v\in \R^n$ of Euclidean norm 1.

\begin{remark}[Complex case]
    The approach taken here extends to the complex case, where the phase of $u$ is chosen such that
    $u^*(A+\eps E)v$ is real and positive. This can be ensured by adding the constraint $\Im(u^*(A+\eps E)v)=0$ via an extra term $\iu\omega u$ to the differential equation for $u$, where $\omega$ is determined such that the time-differentiated constraint is satisfied. Note that the functional is invariant under rotations $\e^{\iu\theta}u$.
\end{remark}

\subsubsection*{Outer iteration.}
For a given tolerance $\vartheta\ll 1$,
the optimal perturbation size ${\eps_\vartheta}$ is determined as the smallest $\eps>0$ such that 
$$
\sigma_{\min}(A+\eps E(\eps)) =
\Bigl(\tfrac12\,\|u(\eps)^{\top}(A+\eps E(\eps))\|^2 + \tfrac12\,\|(A+\eps E(\eps))v(\eps)\|^2\Bigr)^{1/2} \le \vartheta \,\|A\|_2.
$$
Here the first equality holds true because $u(\eps)$ and $v(\eps)$ are left and right singular vectors associated with the smallest singular value of $A+\eps E(\eps)$.
This scalar equation (up to $O(\vartheta))$ is again solved by a Newton--bisection algorithm. Here we use the following simple expression for the
derivative.
\index{Newton--bisection method}

\begin{lemma} We have
$$
\frac{d}{d\eps} \sigma_{\min}(A+\eps E(\eps)) = - \|\Pi^\cS(u(\eps)v(\eps)^\top)\|_F.
$$
\end{lemma}

\begin{proof}
By Corollary~\ref{chap:appendix}.\ref{lem:singderiv}, 
$$
\frac{d}{d\eps} \sigma_{\min}(A+\eps E(\eps)) =u(\eps)^\top \bigl(E(\eps)+\eps E'(\eps)\bigr) v(\eps).
$$
We find
\begin{align*}
u(\eps)^\top E'(\eps) v(\eps) &= \langle u(\eps)v(\eps)^\top, E'(\eps) \rangle =
\langle \Pi^\cS(u(\eps)v(\eps)^\top), E'(\eps) \rangle 
\\
&= 
-\mu \langle E(\eps),E'(\eps)\rangle = -\frac \mu 2 \frac{d}{d\eps} \|E(\eps)\|_F^2 =0
\end{align*}
with the irrelevant factor $\mu=\|\Pi^\cS(u(\eps)v(\eps)^\top)\|_F$. So we have
\begin{align*}
&\frac{d}{d\eps} \sigma_{\min}(A+\eps E(\eps)) = u(\eps)^\top E(\eps) v(\eps) 
= - u(\eps)^\top \frac{\Pi^\cS(u(\eps)v(\eps)^\top)}{\|\Pi^\cS(u(\eps)v(\eps)^\top)\|_F}\, v(\eps) 
\\
&= -\frac{\langle u(\eps)v(\eps)^\top,\Pi^\cS(u(\eps)v(\eps)^\top)}{\|\Pi^\cS(u(\eps)v(\eps)^\top)\|_F}
= -\frac{\langle \Pi^\cS (u(\eps)v(\eps)^\top),\Pi^\cS(u(\eps)v(\eps)^\top)\rangle}{\|\Pi^\cS(u(\eps)v(\eps)^\top)\|_F}
\\
&= - \|\Pi^\cS(u(\eps)v(\eps)^\top)\|_F
\end{align*}
as stated in the lemma.
\qed
\end{proof}

\bng
We remark that in this case we cannot use the HEC algorithm, which works with values of $\eps>\oeps$, with $\oeps$ the smallest value of $\eps$ such that
$\sigma_{\min}(A+\eps E(\eps)) = 0$. The problem is that $\sigma_{\min}(A+\eps E(\eps))$ vanishes identically for
$\eps \ge \oeps$. On the other hand, the monotone Newton--bisection algorithm is constructed such that the iterates $\eps_k$ grow monotonically and thus converge from the left.
\eng

\subsection{Example: Nearest pair of polynomials with a common zero} 
\index{coprime polynomials}
We consider polynomials with real coefficients (alternatively, we could allow for complex coefficients).
A pair of polynomials $(p,q)$ is called {\it coprime} if $p$ and $q$ have no nontrivial common divisor, or equivalently, have no common zero. For a pair of polynomials $(p,q)$ that is coprime, the distance to the nearest pair of polynomials with a common zero is of interest. In the following we measure the distance of pairs of polynomials by the Euclidean norm of the difference of the vectors of coefficients in the monomial basis.

\medskip\noindent
{\bf Problem.} {\it Given a pair of polynomials that is coprime, find the nearest pair of polynomials with a nontrivial common divisor.}

\medskip\noindent
Consider polynomials 
\begin{eqnarray}
\begin{array}{rcl}
p(z) & = & a_n z^n + a_{n-1} z^{n-1} + \cdots + a_1 z + a_0
\\[2mm]
q(z) & = & b_m z^m + b_{m-1} z^{m-1} + \cdots + b_1 z + b_0
\end{array}
\label{eq:pq}
\end{eqnarray}
with real coefficients  $a_i$ and $b_i$. We may assume $m\le n$ and $a_n\ne 0$, and  we can take $m=n$ if we allow for $b_n=0$. So we assume $m=n$ in the following. 
With these polynomials we associate the Sylvester matrix of dimension $2n\times 2n$,
\index{Sylvester matrix}
\begin{equation}
S(p,q) := 
\begin{pmatrix}
a_n & & \ldots & & a_0 &  & &   
\\
 & a_n & & \ldots & & a_0 & & 
\\
 & & \ddots &  & & &  \ddots &     
\\
 &  &  & a_{n} & & \ldots & &  a_0 
\\
b_n & & \ldots & & b_0 &  & &   
\\
 & b_n & & \ldots & & b_0 & & 
\\
 & & \ddots &  & & &  \ddots &     
\\
 &  &  & b_n & & \ldots & &  b_0 
\end{pmatrix}.
\label{eq:Sylpq}
\end{equation}
A result from commutative algebra (see Laidacker~(\cite{Lai69})) states:\\
{\it The pair of polynomials $(p,q)$ is coprime if and only if  the associated Sylvester matrix $S(p,q)$  is invertible}.

This allows us to reformulate the problem of nearest non-coprime polynomials as a matrix nearness problem, where the distance is measured by the Frobenius norm.

\medskip\noindent
{\bf Problem.} {\it Given an invertible Sylvester matrix, find the nearest singular Sylvester matrix.}

\medskip\noindent
This amounts to the problem of computing the structured distance to singularity, where the structure is given by the subspace $\cS \subset \R^{2n,2n}$ of Sylvester matrices of the form \eqref{eq:Sylpq}.

Let $S$ be a Sylvester matrix. 
We define the \emph{lower radius} of $S$ as
\[
\mu(S) = \min\{ |\lambda| \, : \, \lambda \text{ is an eigenvalue of S} \},
\]
and note that $S$ is singular if and only if $\mu(S) = 0$. Furthermore, with the structured $\eps$-pseudospectrum $\Lameps^{\cS}(S)$, let
\begin{equation}\label{eq:psa}
\lreps^\cS (S) = \min\{|\lambda|: \lambda \in \Lameps^{\cS}(S) \},
\end{equation} 
which reduces to the lower radius $\mu(S)$ when $\eps=0$. 
We can then express the radius of coprimeness of the pair of polynomials $(p,q)$ as
\[
\rho_{\rm co}(p,q) = \frac{\oeps}{\sqrt{n}}   \quad\text{ with }\quad 
\oeps =  \min \{ \eps>0 \,:\, \lreps^\cS(S) = 0 \}.
\]
(The division by $\sqrt{n}$ is done to account for the fact that each coefficient $a_i$ and $b_i$ appears $n$ times in the Sylvester matrix, which yields a factor $\sqrt{n}$ in the Frobenius norm.)

\subsubsection*{Two-level iteration.}
\index{two-level iteration}
We are thus in the situation of applying the two-level iteration of Section~\ref{sec:two-level} with the functional
$\F_\eps(E)$ (for $E\in \cS$ of unit Frobenius norm) given as
\begin{equation} \label{eq:FepSyl}
\F_\eps(E) = \mu(S+\eps E), 
\end{equation}
which is of the form \eqref{Feps} with $f(\lambda,\clambda)=\sqrt{\lambda\clambda}$ and with the eigenvalue of smallest modulus as target eigenvalue.

To apply the gradient-based algorithm of Section~\ref{sec:proto-structured} in the inner iteration, we need the structured gradient, see \ref{chap:proto}.\ref{eq:deriv-S}), which is the orthogonal projection onto the space $\cS$ of Sylvester matrices of the
(rescaled) gradient
$$
G_\eps(E) = \Re \Bigl(\frac\lambda{|\lambda|}\, xy^*\Bigr),
$$
where $x$ and $y$ are the left and right eigenvectors, normalized to unit norm and with positive inner product, that are associated with the eigenvalue $\lambda$ of smallest modulus of $S+\eps E$, which is assumed to be simple. We note that the structured gradient is nonzero by Theorem~\ref{thm:nonzero-gradient-S}.

The orthogonal projection onto $\cS$ is given in the following lemma.
\index{structure space!orthogonal projection}

\begin{lemma}[Orthogonal projection onto the space of Sylvester matrices] \label{lem:projS}
Let $S \in \cS \subset \R^{2n\times 2n}$ and $Z \in \C^{2n,2n}$; the orthogonal 
 projection $P_\cS$ onto $\cS$, with respect to the Frobenius inner product $\la \cdot, \cdot \ra$, is given by 
\begin{eqnarray}
\Pi^\cS Z & = & S(p,q), \label{eq:PBpq}
\end{eqnarray}
where $p$ and $q$ are the polynomials with coefficients (for $k=0,\ldots,n$)
\begin{equation}
a_{n-k} = \frac{1}{n} \sum\limits_{l=1}^{n} \Re \left( Z_{l,l+k} \right), \qquad
b_{n-k} = \frac{1}{n} \sum\limits_{l=1}^{n} \Re \left( Z_{n+l,l+k} \right). 
\nonumber
\end{equation}
\end{lemma} 
\def\one{\doubleone}
\begin{proof}
We have to find $\arg\min_{S \in \cS} \| Z - S \|_F$.
The result follows directly from the fact that for a complex vector $x \in \C^n$,
\[
\mu_* = \arg\min\limits_{\mu \in \R} \| x - \mu \one \|_F = \frac{1}{n} \sum\limits_{i=1}^{n} \Re(x_i),
\] 
where $\one = \left( 1 \ 1 \ \ldots \ 1 \right)^\tp$.
\qed
\end{proof}

\subsubsection*{Numerical example.}
\label{sec:exill}

Consider the two polynomials of degree $3$,
\begin{equation} \label{ex:pq1}
p(z) = z^3 + 2 z^2 + 2 z + 2, \quad\
q(z) = 2 z^3 + z - 2 ,
\end{equation}
where $p$ is constrained to be monic. Here $a = \left( 1 \ 2 \ 2 \ 2 \right)^\tp$ and $b = \left( 2 \ 0 \ 1 \ -2 \right)^\tp$;
the corresponding Sylvester matrix is given by
\begin{equation}
S(a,b) = 
\left(
\begin{array}{rrrrrr}
1 & 2 & 2 & 2 & 0 & 0 
\\
0 & 1 & 2 & 2 & 2 & 0  
\\
0 & 0 & 1 & 2 & 2 & 2 
\\
2 & 0 & 1 & -2 & 0 & 0 
\\
0 & 2 & 0 & 1 & -2 & 0 
\\
0 & 0 & 2 & 0 & 1 & -2 
\end{array}
\right)
\label{eq:illS}
\end{equation}

\begin{figure}[ht]
\vskip -5.3cm
\centering
\includegraphics[width=\textwidth]{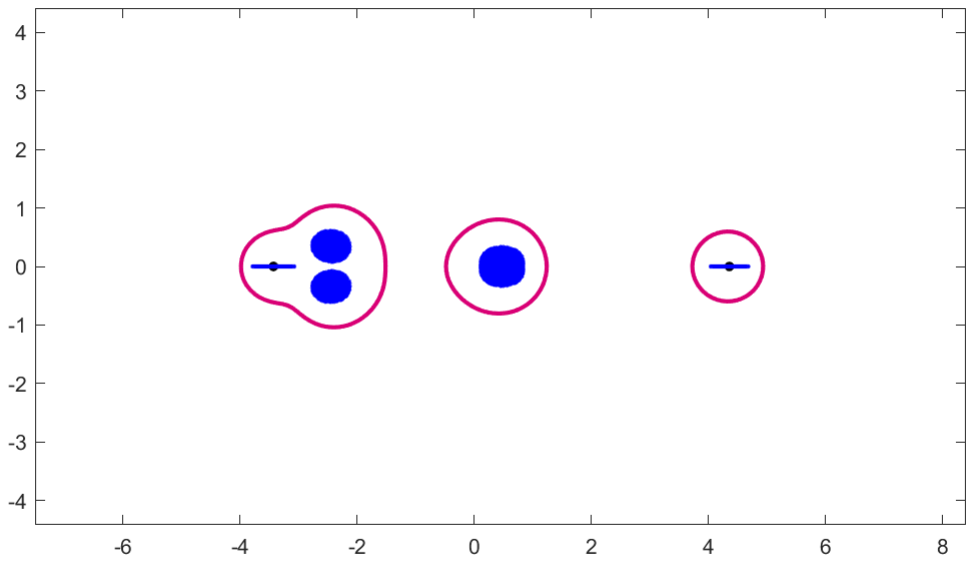}
\vskip -5cm
\caption{The approximated structured Sylvester $\eps$-pseudospectrum 
for  $\eps = \frac12$ for Example \eqref{eq:illS} is filled with blue. 
The red curve represents the boundary of the set of eigenvalues obtained by considering arbitrary 
complex perturbations (that is omitting the constraint of real Sylvester structure) 
of norm bounded by $\frac12$.} 
\label{fig:illS}
\end{figure}

The structured pseudospectrum $\Lameps^{\cS}(S)$ for $\eps=\frac12$ is approximated by dense sampling on the 
set of admissible perturbations and is plotted in blue in Figure \ref{fig:illS}. 

It turns out that for the value 
\begin{equation} \nonumber
\eps = \oeps = 0.618108064
\end{equation}
the functional $\F_\eps(E)$ (see \eqref{eq:FepSyl}) vanishes, while for $\eps < \oeps$ it holds $\F_\eps(E) > 0$.
This gives $\rho_{\rm co}(p,q) = 0.356864857$.

The computed matrix $S + \oeps E(\oeps)$ has rank $2n-2$ due to a double semi-simple zero eigenvalue. The 
coefficients of the perturbed polynomials  $\hat p$, $\hat q$  are shown (with $10$-digit accuracy) in Table \ref{t1}.
The common complex conjugate zeros of $\hat p$, $\hat q$ are 
\[
z_{1,2} = -0.4008686595 \pm 1.03085391659 \iu.
\]
\vskip -7mm
\begin{table*}[htb!]
\centering
\caption{Coefficients of the perturbed polynomials $\hat p$, $\hat q $ in the example 
\ref{ex:pq1}.
} \label{t1}
\(
\begin{array}{llllll}
\hat a_3 =  0.75744188 &
\hat a_2 =  2.10479150 &
\hat a_1 =  2.12724001 &
\hat a_0 =  1.83200184 
\\[1mm]    
\hat b_3 =  1.95430087 &
\hat b_2 = -0.06706025 &
\hat b_1 =  1.08084913 &
\hat b_0 = -1.99883585 
\end{array}
\)
\end{table*}



\section{Notes}

\subsubsection*{Low-rank matrix differential equations for structured eigenvalue optimization.}
The approach to structured eigenvalue optimization via a norm-constrained gradient system and the associated low-rank dynamics was first proposed and studied by Guglielmi \& Lubich (\cite{GL13}), where rank-2 differential equations  were used to compute the real $\eps$-pseudospectral abscissa and radius. Guglielmi, Kressner \& Lubich (\cite{GKL15}) studied Hamiltonian eigenvalue optimization using rank-4 matrix differential equations; see also Section~\ref{sec:Hamilton}.
The structure-projected rank-1 approach in Section~\ref{sec:proto-structured} to eigenvalue optimization with general linear structures $\cS$ was developed by Guglielmi, Lubich \& Sicilia (\cite{GuLS23}).

Our discussion of low-rank dynamics in Sections~\ref{subsec:rank1-gradient-flow} and~\ref{subsec:rank-r-gradient-flow} is based on Koch \& Lubich (\cite{KL07}). Numerical integrators for
low-rank matrix differential equations that are robust to small singular values are given by the projector-splitting integrator of Lubich \& Oseledets (\cite{LuO14}) and the basis-update \& Galerkin (BUG) integrator of Ceruti \& Lubich (\cite{CeL22}), of which a norm-preserving variant is presented in Section~\ref{subsec:low-rank-integrator}.

\subsubsection*{Frobenius norm vs.~matrix 2-norm.}
In the real case and structured cases, where the Frobenius-norm and matrix 2-norm optimizers are in general different, functional-reducing differential equations for the 2-norm constrained eigenvalue optimization problem have been  given that have similar properties to the gradient systems for the Frobenius-norm problem considered here; see Guglielmi \& Lubich (\cite{GL13}) for the computation of the 2-norm real pseudospectral abscissa and Guglielmi, Kressner \& Lubich (\cite{GKL15}) for 2-norm constrained eigenvalue optimization problems for Hamiltonian matrices. 

Graillat (\cite{Gra06}) and Rump (\cite{Rum06}) proved that for complex Toeplitz and Hankel structures, the 2-norm structured pseudospectrum equals the unstructured pseudospectrum. This is not valid in the Frobenius norm, and it is not valid for real Toeplitz and Hankel matrices. 

\subsubsection*{Computing structured pseudospectra.}
The approach of Karow, Kokiopoulou and Kressner (\cite{KarKK10}) to compute structured pseudospectra w.r.t. the matrix 2-norm uses structured singular values
$
\mu^\cS(B)=
(\inf\{ \| \Delta \|_2\,:\, \Delta \in \cS \text{ and $I-\Delta B$ is a singular matrix} \} 
)^{-1},
$
which yield the characterization
$
\Lambda_\eps^\cS(A) = \{ \lambda \in \C\, : \, \mu^\cS\bigl( (A- \lambda I)^{-1}\bigr)\ge 1/\eps \}.
$
Note that $1/\mu^\cS\bigl( (A- \lambda I)^{-1}\bigr)$ equals the structured 2-norm distance to singularity of $A-\lambda I$; 
cf.~Section~\ref{sec:sing-S}, where the structured Frobenius-norm distance to singularity is considered.
For the unstructured complex case we have $\mu^\C(B)=\| B \|_2$, and for the unstructured real case, an expression for
$\mu^\R(B)$ was given by Qiu et al. (\cite{QiuBRDYD95}).
Based on that result, Karow, Kokiopoulou and Kressner (\cite{KarKK10}) gave an algorithm to compute real pseudospectra. 
They also found expressions for structured singular values for complex skew-symmetric, Hermitian and complex Hamiltonian matrices, which allowed them to extend their algorithm to these particular structures. 
For the Frobenius-norm structured pseudospectrum, their basic algorithm can be combined with the algorithm for the structured distance to singularity of Section~\ref{sec:sing-S} to compute the Frobenius-norm structured singular values of the resolvents $(A-\lambda I)^{-1}$ for general linear structures.

Butt\`a, Guglielmi \& Noschese (\cite{ButGN12}) proposed an algorithm to approximate the structured $\eps$-pseudospectral abscissa (w.r.t. the Frobenius norm) and sections of the boundary of structured pseudospectra for complex Toeplitz matrices. This algorithm is related to the rank-1 iteration of Guglielmi \& Overton (\cite{GO11}).

The ladder algorithm and its structured version used in Section~\ref{subsec:tracing-S}, which traces sections of the boundary of structured pseudospectra (w.r.t. the Frobenius norm) for arbitrary linear structures $\cS$, have not previously appeared in the literature.


\subsubsection*{Stability radius under real unstructured perturbations.}
Qiu, Bernhardsson, Rantzer, Davison, Young \& Doyle (\cite{QiuBRDYD95}) characterized the real stability radius $r_\R(A)$ as
$$
\frac1{r_\R(A)} = \sup_{\omega\in\R} \inf_{0<\gamma\le 1} \sigma_2
\begin{pmatrix}
\Re\, M_\omega & - \gamma\, \Im\, M_\omega \\
\gamma^{-1} \Im \, M_\omega & \Re \, M_\omega
\end{pmatrix}
\quad\ \text{with}\quad M_\omega = (A-\iu \omega I)^{-1},
$$
where $\sigma_2(\cdot)$ is the second largest singular value of a matrix.
An algorithm for the computation of $r_\R(A)$ via this formula was proposed by Sreedhar, Van Dooren \& Tits (\cite{SreVDT96}).
Based on a reformulation of this formula and using Byers' connection between singular values and eigenvalues of Hamiltonian matrices,
Freitag \& Spence (\cite{FreS14}) also developed an algorithm to deal efficiently with this two-dimensional optimization.

In a different approach, Lu \& Vandereycken (\cite{LuV17}) developed a criss-cross-type algorithm (cf.~\ref{sec:psa}) for computing the real pseudospectral abscissa. 

Guglielmi \& Manetta (\cite{GM15}) studied an algorithm that is well-suited also for large sparse matrices. It corresponds to the  general two-level approach taken in Section~\ref{sec:two-level}. In the inner iteration, the algorithm computes the real $\eps$-pseudospectral abscissa via rank-2 matrix differential equations of Guglielmi \& Lubich (\cite{GL13}) (which are given there for both the matrix 2-norm and the Frobenius norm). In the outer iteration, it uses a combined Newton--bisection method to optimize the perturbation size $\eps$ to yield $\oeps$ such that the real $\oeps$-pseudospectral abscissa becomes zero. The Newton iteration used the simple derivative formula of Theorem~\ref{thm:phi-derivative} for the particular case of the real gradient $G= \Re (xy^*)$. A related method based on a real version of the iteration method of Guglielmi \& Overton (\cite{GO11}) was proposed by Rostami (\cite{Ros15}) and further developed and analysed by Guglielmi (\cite{Gug16}).

\subsubsection*{Structured stability radii.}
In the control systems literature, Hinrichsen \& Pritchard (\cite{HinP86a},\cite{HinP86b},\cite{HinP90}) considered complex and real stability radii
(i.e. distance to instability under complex and real perturbations) and also structured stability radii
$$
r(A,B,C)= \min \{ \| \Delta \| \, : \, A+B\Delta C \text{ has some eigenvalue of nonnegative real part} \},
$$
where $A$ is a Hurwitz-stable matrix and $B$ and $C$ are given matrices of compatible dimensions. The perturbation matrix $\Delta$ is assumed to be real or complex. Most of the algorithms mentioned above extend to this situation of range- and corange-restricted perturbations. For example, Hinrichsen, Kelb \& Linnemann (\cite{HinKL89}) extended Byers' algorithm to compute the complex structured stability radius $r_\C(A,B,C)$. 

We note that this notion of structured stability radius minimizes the norm of the parameter matrix $\Delta$ and not of the structured perturbation $B\Delta C$ of $A$. The latter would fit directly into the framework of Section~\ref{sec:two-level}, whereas controlling the norm of $\Delta$ requires some (minor) modifications to the algorithm.

We are not aware of algorithms for other structured stability radii (distance to instability under complex or real structured perturbations) in the literature, e.g. for perturbations with a given sparsity pattern and/or symmetry, or Toeplitz perturbations etc. The two-level algorithm of Section~\ref{sec:two-level} addresses such problems with general linear structures.

\subsubsection*{Transient bounds under structured perturbations.}
The joint unstructured--structured pseudospectrum and the structured $\eps$-stability radius were introduced by
Guglielmi \& Lubich (\cite{GL25}) as basic notions to assess the robustness of transient bounds of linear dynamical systems under structured perturbations of the matrix. The rank-1 algorithm of Section~\ref{sec:eps-stab} for computing the structured $\eps$-stability radius was also proposed and studied there.

\subsubsection*{Numerical range under structured perturbations.}
The algorithm in Section~\ref{sec:num-range} for computing the structured dissipativity radius using Rayleigh quotients instead of eigenvalues has not been published before. More generally, replacing eigenvalues by Rayleigh quotients in  real symmetric and complex Hermitian 
\bcl matrix nearness problems can yield a computationally favourable approach; see also Chapter~\ref{chap:graphs}.\ecl

\subsubsection*{Structured distance to singularity.}
While the unstructured 2-norm and Frobenius-norm distance to singularity of a matrix equal its smallest singular value, we are not aware of an algorithm in the literature for computing the distance to singularity under general structured perturbations. 
\bng The important special case of perturbations of the form $B \Delta C$ with given real rectangular matrices $B$ and $C$ is addressed by Qiu et al. (\cite{QiuBRDYD95}). \eng

Approximating the structured distance to singularity by minimizing the norm of matrix-vector products within a two-level approach as is done in Section~\ref{sec:sing-S} appears to be new.



\renewcommand{\Id}{\mathrm{I}}

\chapter{Matrix nearness problems of diverse kinds}
\label{chap:mnp-mix}

In this chapter we discuss extensions of the two-level approach to matrix nearness problems of diverse types that are not covered by the framework of the previous chapters, be it because
\begin{itemize}
    \item the functional in the associated eigenvalue optimization problem depends not only on one target eigenvalue but on a varying number of eigenvalues, as in the problems of computing the nearest stable matrix to a given unstable matrix and computing the nearest correlation matrix to a given symmetric real matrix; or
    \item the problem requires special attention to the structure, as in Hamiltonian matrix nearness problems; or 
    \item the functional in the associated eigenvalue optimization problem depends on eigenvectors, as in the Wilkinson problem of computing the nearest matrix with defective eigenvalues, which amounts to computing the distance to singularity of the eigenvalue condition number; or
    \item the nearness problem deals with matrix pencils, as in the problem of finding the nearest matrix pencil that is singular or additionally has a common null-vector; or
    \item the underlying eigenvalue problems are nonlinear in the eigenvalues, as in computing
    the stability radius of linear delay differential equations; or
    \item the given problem is not itself a matrix nearness problem but is addressed by the same techniques, such as computing
    the Kreiss constant that appears in the Kreiss matrix theorem (which leads to generalized eigenvalue problems similar to the preceding item).
\end{itemize}
These items are exemplary, not exhaustive. They form the sections of this chapter. The sections can be read independently of each other, but we give fewer details in the later sections. 

\bcl Like the previous chapter, also this chapter illustrates the versatility of the two-level approach that uses rank-constrained gradient flows 
 and Newton--bisection methods in a nested or alternating way as discussed in Chapter~\ref{chap:two-level}. For ease of presentation we will here only formulate simple inner--outer iterations, but it is implicitly understood that monotonically increasing Newton--bisection iterations and monotonically decreasing HEC iterations can be used instead and may be favourable.\ecl

This chapter enriches the toolbox for applications in various fields, such as those considered in the final chapters of this book. The proposed algorithms can deal with complex, real and structured perturbations.

\section{Matrix stabilization}
\label{sec:mat-stab}
\index{matrix stabilization}

In this section we extend the two-level approach of Chapter~\ref{chap:two-level} to the problem of moving all eigenvalues of a given matrix into a prescribed closed subset $\overline\Omega$ of the complex plane by a perturbation of minimal Frobenius norm. This spectral recovery problem is complementary to the robustness analysis in the previous two chapters where the original matrix has all its eigenvalues inside $\Omega$ and it was required that one eigenvalue be driven to the boundary of $\Omega$.
While the approach presented here is conceptually applicable to the spectral recovery problem for very general subsets $\Omega$, we will
focus our attention on the guiding problem of {\it Hurwitz stabilization} of an unstable matrix, which corresponds to the case where $\Omega={\C^-}$ is the complex left half-plane:
\index{Hurwitz stabilization}
\index{matrix stabilization!structured}

\medskip
\noindent
{\bf Problem.} {\it Given a square matrix $A$ that has some eigenvalues with positive real part, find a structured perturbation $\Delta\in\cS$ of 
minimal Frobenius norm such that $A+\Delta$ has no eigenvalue with positive real part.} 

\smallskip




The perturbations $\Delta$ are restricted to lie in a structure space $\cS$ that is $\C^{n,n}$ or $\R^{n,n}$ or a linear subspace thereof. 
\bcl
In the latter case, it is not guaranteed that the problem has a solution.
\ecl

We describe two algorithmic approaches to this matrix stabilization problem, which are both of the two-level type considered in Chapter~\ref{chap:two-level}.  
\begin{itemize}
\item The {\it exterior algorithm} moves eigenvalues that lie outside the closed target set $\overline\Omega$ where the eigenvalues of the perturbed matrix $A+\Delta$ should lie ($\overline\Omega=\overline {\C^-}$ for Hurwitz stability and 
$\overline\Omega$ is the closed unit disk for Schur stability).
All eigenvalues outside $\overline\Omega$ are moved towards the boundary of $\Omega$ while increasing the perturbation size. 
\item The {\it interior algorithm} starts from a non-optimal perturbation $\Delta_0$  such that $A+\Delta_0$ has all eigenvalues in $\Omega$ and
moves an eigenvalue inside $\Omega$ to the boundary while reducing the perturbation size.
\end{itemize}

\subsection{Exterior two-level algorithm}
\label{subsec:ext-stab}
\index{matrix stabilization!exterior two-level algorithm}
\index{two-level iteration}

Here we use the following eigenvalue optimization problem: For a given perturbation size $\varepsilon>0$, find 
\begin{equation} \label{f-E-eps-stab}
\arg\min\limits_{\Delta \in \cS,\, \| \Delta \|_F = \eps} \ \sum_{i=1}^n \,
f \left( \lambda_i\left( A + \Delta \right), \clambda_i \left( A + \Delta \right)  \right),
\end{equation}
with
\begin{equation}
f\left( \lambda, \clambda \right) = \sfrac12\,\text{dist}(\lambda,\overline{\C^-})^2 = \sfrac12 \bigl((\Re\,\lambda)_+\bigr)^2 =
\sfrac18\bigl((\lambda + \clambda)_+\bigr)^2,
\label{eq:Hurw} 
\end{equation}
where for $a \in \R$, $a_+ := \max\left\{a,0\right\}$.

The eigenvalues $ \lambda_i(A+\Delta)$ ($i=1,\dots,n$) of the perturbed matrix $A + \Delta$ are
ordered by decreasing size of the real part. Note that only the eigenvalues with positive real part contribute to the sum in \eqref{f-E-eps-stab}, so that the sum only extends from 1 to $m^+(A+\Delta)$, where
$m^+(M)$ is the number of eigenvalues of $M$ with positive real part. 
As in Chapter~\ref{chap:two-level}, this leads us to the following two-level approach to the Hurwitz stabilization problem. Here we additionally introduce a small shift $\delta>0$ that aims for strict Hurwitz stability with all eigenvalues having real part not exceeding $-\delta$.
\begin{itemize}
\item {\bf Inner iteration:\/} Given $\eps>0$, we aim to compute a  matrix $E(\eps) \in\cS$  
of unit Frobenius norm  such that 
\begin{equation}
\F_\eps(E) = \sfrac12 \sum_{i=1}^n \, \bigl((\Re\,\lambda_i( A + \eps E)+\delta)_+\bigr)^2
\label{Feps-stab}
\end{equation}
is minimized, i.e. 
\begin{equation} \label{E-eps-stab}
E(\eps) = \arg\min\limits_{E\in \cS, \| E \|_F = 1} \F_\eps(E).
\end{equation}


\item {\bf Outer iteration:\/} We compute the smallest positive value $\oeps$ with
\begin{equation} \label{eq:zero-stab}
\phi(\oeps)= 0,
\end{equation}
where $\phi(\eps)=  \F_\eps\left(E(\eps) \right) = \tfrac12 \sum_{i=1}^n \, \bigl((\Re\,\lambda_i( A + \eps E(\eps)+\delta)_+\bigr)^2$.
\end{itemize}

\medskip
We remark that the existence of a zero of $\phi$ (i.e., stabilizability) is not guaranteed for arbitrary structure spaces $\cS$, but it is when $\cS$ equals $\C^{n,n}$ or $\R^{n,n}$ or more generally when $\cS$ contains real multiples of the identity matrix $I$ (since then some negative shift of the given matrix moves all its eigenvalues into the complex left half-plane).

As before, the inner iteration uses a norm- and structure-constrained gradient-flow differential equation, possibly further restricted to low-rank dynamics; the outer iteration uses a hybrid Newton--bisection method. We give details in the following subsections.

\subsubsection{Constrained gradient flow for minimizing $\F_\eps(E)$}
\label{subsec:grad-flow-stab}
Here we let $\eps>0$ be fixed.
Let $E(t)\in \cS$, \bng for $t$ in an interval $[0,T]$, \eng be a continuously differentiable path of matrices in the structure space $\cS\subset \C^{n,n}$, and let the eigenvalues $\lambda_i(t)=\lambda_i(A+\eps E(t))$ be simple for $i=1,\dots,n$ and all $t\in [0,T]$. The corresponding left and right eigenvectors $x_i$ and~$y_i$ are assumed to be of unit norm and with $x_i ^* y_i>0$.
Then, applying Theorem \ref{chap:appendix}.\ref{thm:eigderiv} we obtain 
\begin{equation}
\frac{d}{dt} \F_\eps\bigl(E(t)\bigr) = \eps \sum\limits_{i=1}^{n}  \bigl( \Re\,\lambda_i (A + \eps E(t)) + \delta \bigr)_+ \,
\frac{\Re(x_i(t) ^* \dot{E}(t) y_i(t))}{x_i(t) ^* y_i(t)}.
\label{eq:derFeps}
\end{equation}
We introduce the notation
\begin{equation*}
\gamma_i (t) :=\frac{\left( \Re \,\lambda_i (A + \eps E(t)) + \delta \right)_+}{x_i(t)^* y_i(t)} \ge 0 .
\end{equation*}
Here we note that $\gamma_i(t)=0$ for $i>m_\delta(A+\eps E(t))$, which is the number of eigenvalues with real part greater than $-\delta$.
We write \eqref{eq:derFeps} as
$$
\frac{d}{dt} \F_\eps\bigl(E(t)\bigr) = 
\eps\sum\limits_{i=1}^{n} \gamma_i(t) \,\Re\bigl(x_i(t) ^* \dot{E}(t) y_i(t)\bigr) = 
\eps\,\Re\,\bigl\langle G_\eps(E(t)), \dot{E}(t) \bigr\rangle
$$
with the rescaled free gradient
\begin{equation}\label{grad-stab}
G_\eps(E)= \sum\limits_{i=1}^{n} \gamma_i x_i y_i^*.
\end{equation}
This is of rank at most $m_\delta(A+\eps E)$.
In the structured case where $E(t)\in\cS$ for all $t$, and hence also its derivative is in $\cS$ so that $\dot E(t)=\Pi^\cS \dot E(t)$, we further obtain that with the projected gradient
$$
G_\eps^\cS(E)=\Pi^\cS G_\eps(E)\in\cS,
$$ 
we have
\begin{equation}
\frac{d}{dt} \F_\eps\bigl(E(t)\bigr) = 
\eps\,\Re\,\bigl\langle G_\eps^\cS(E(t)), \dot{E}(t) \bigr\rangle.
\label{eq:derFeps2}
\end{equation} 
As in Section~\ref{sec:proto-structured}, see~(\ref{chap:struc}.\ref{ode-E-S}), we consider the constrained gradient flow 
\begin{equation}\label{ode-E-S-stab}
\dot E = -G_\eps^\cS(E) + \Re \langle G_\eps^\cS(E), E \rangle E,
\end{equation}
which again has the properties that 
\begin{itemize}
\item the unit Frobenius norm is conserved along solutions $E(t)$; 
\item $\F_\eps(E(t))$ decays monotonically with growing $t$; 
\item stationary points $E$ are real multiples of $G_\eps^\cS(E)$ provided that $G_\eps^\cS(E)\ne 0$.
\end{itemize}
\index{stationary point}
Therefore, a stationary point $E$ is a projection onto the structure space $\cS$ of a matrix of rank at most $m_\delta=m_\delta(A+\eps E)$. In particular, in the complex unstructured case the rank is at most $m_\delta$, and in the real case at most $2m_\delta$.

\subsubsection{Rank-constrained gradient flow in the complex and real unstructured cases}
\label{subsec:low-rank-stab}
\index{gradient flow!rank-$r$ constrained}

With an expected upper bound $m$ of $m_{\delta,\eps}=m_\delta(A+\eps E(\eps))$ at the minimizer $E(\eps)$, which is of rank at most
$m_{\delta,\eps}$ in the complex unstructured case and at most $2m_{\delta,\eps}$ in the real case,
we can use a rank-constrained gradient flow in the same way as in Section~\ref{subsec:low-rank}, where the chosen rank is now $r=m$ in the complex case and $r=2m$ in the real case. We then consider the rank-$r$ constrained gradient flow, with $\cS=\C^{n,n}$ or $\R^{n,n}$,
\begin{equation}\label{ode-ErF-2-v2-stab}
\dot E = -P_E(G_\eps^\cS(E) ) + \Re \langle E,  P_E(G_\eps^\cS(E) )\rangle E,
\end{equation}
where $P_E$ is the orthogonal projection onto the tangent space at $E$ of the manifold of (complex or real) $n\times n$-matrices of rank~$r$. This differential equation is of the same type as in Section~\ref{subsec:low-rank}  and is treated numerically in the same way as described there. (In the complex case, transposes of matrices are replaced by conjugate transposes.) 

\subsubsection{Rank-constrained matrix differential equation in structured cases}
\label{subsec:low-rank-stab-S}
Similar to Section~\ref{subsec:rank-1-S} we use the projected differential equation with solutions of rank $
m$ with an expected upper bound $m$ of $m_\delta(A+\eps E(\eps))$ at the minimizer $E(\eps)$:
 \begin{equation}\label{ode-E-S-1-stab}
\dot Y = -P_Y G_\eps(E) + \Re \langle P_Y G_\eps(E), E \rangle Y \quad\text{ with }\ E=\Pi^\cS Y,
\end{equation}
where now $P_Y$ is the projection onto the tangent space at $Y$ of the manifold of rank-$m$ matrices. This rank-$m$ matrix differential equation is solved numerically into a stationary point as described in Section~\ref{subsec:low-rank-integrator}. (In the complex case, transposes of matrices are again replaced by conjugate transposes.)


%

%
%

\subsubsection{Iteration for $\eps$}
To solve the one-dimensional root-finding problem \eqref{E-eps-stab}, we use a Newton--bisection method as in Section~\ref{subsec:Newton--bisection}.
We let $E(\eps)$ of unit Frobenius norm be a (local) minimizer of the optimization problem \eqref{Feps-stab} and we denote by 
$\lambda_i(\eps)$, the 
eigenvalues and by~$x_i$ and~$y_i$ corresponding left and right 
eigenvectors of $A+\eps E(\eps)$, of unit norm and with positive inner product.

We denote by $\oeps$ the smallest value of $\eps$ such that $\phi(\eps)=\F_\eps(E(\eps))=0$. 
%
%
%
%
%
%
For a Newton-like algorithm we need an extra assumption that plays the same role as Assumption~\ref{chap:two-level}.\ref{ass:E-eps}.

\begin{assumption}
For $\eps$ close to $\oeps$ and $\eps<\oeps$, 
we assume the following:
\begin{itemize}
\item The eigenvalues  of $A + \eps E(\eps)$  
with positive real part are simple eigenvalues.
\footnote{This assumption is also required in the inner iteration, see the discussion preceding \eqref{eq:derFeps}.}
\item The map $\eps \mapsto E(\eps)$ is continuously differentiable.
\item The structured gradient $G(\eps):=G_\eps^\cS(E(\eps))$ is nonzero.
\end{itemize}
\label{assumpt-stab}
\end{assumption}

%
%
%
%
The following result extends Theorem~\ref{chap:two-level}.\ref{thm:phi-derivative} from one to several eigenvalues and is proved by the same arguments.
\begin{lemma}[Derivative for the Newton iteration]
\label{lem:der-stab}
Under Assumption~{\rm \ref{assumpt-stab}}, the function $\phi(\eps)=\F_\eps(E(\eps))$ 
is differentiable and its derivative equals 
\begin{equation}
\phi'(\eps) =  - \| G(\eps) \|_F.
\label{eq:derFdeps}
\end{equation}
\end{lemma}
With the derivative of $\phi$ at hand, we compute $\oeps$ by a hybrid Newton--bisection algorithm as described in Section~\ref{subsec:Newton--bisection}. In addition, since $\F_\eps(E(\eps))=0$ for $\eps\ge \oeps$, we take a bisection step when all eigenvalues of $A+\eps E(\eps)$ have real part smaller than $-\delta$.
\medskip

\begin{figure}[ht]
\vskip -7cm
\begin{center}
\includegraphics[width=10cm]{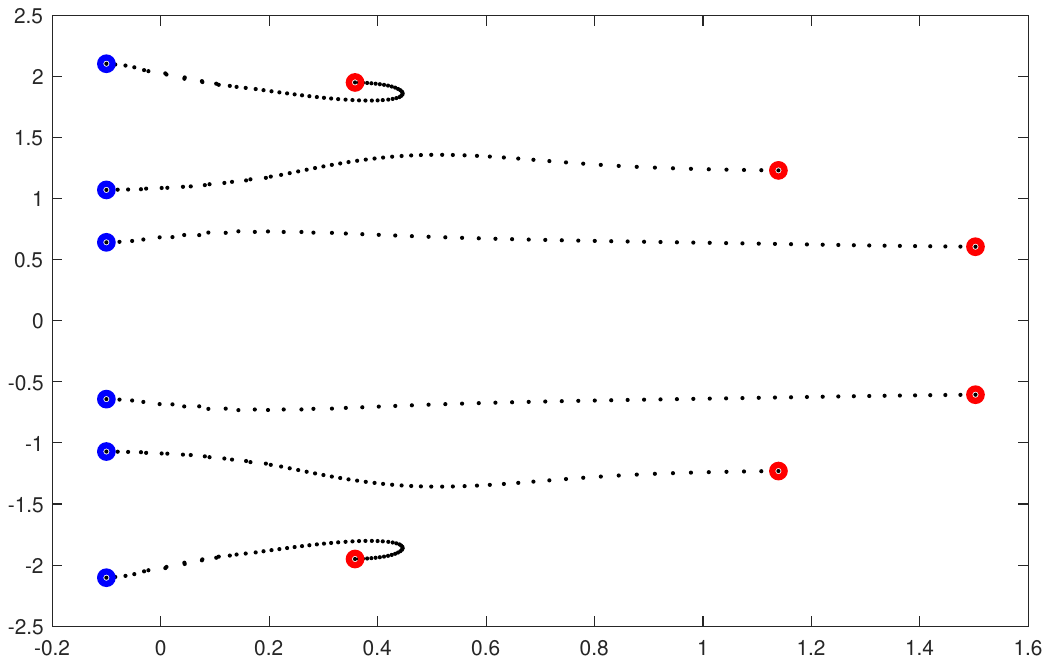} \hskip 0.3cm 
\caption{Paths of the eigenvalues of the matrix $A + \eps E(\eps)$
for $\eps \in [ 0, \oeps ]$.
}
\label{fig:Stab1}
\end{center}
\end{figure}

As an example, consider the Grcar matrix of dimension $n=6$,
\begin{equation}
A = \left(    \begin{array}{rrrrrr}
     1  &   1  &   1    &  1    &  0    & 0 \\
    -1  &   1  &   1    &  1    &  1    & 0 \\
     0  &  -1  &   1    &  1    &  1    & 1 \\
     0  &   0  &  -1    &  1    &  1    & 1 \\
     0  &   0  &   0    & -1    &  1    & 1 \\
     0  &   0  &   0    &  0    & -1    & 1
\end{array}
\right).
\end{equation}
The matrix has all six eigenvalues in the complex right half-plane (red points in Figure \ref{fig:Stab1}).
We set $\delta=0.1$, meaning we intend to push the whole spectrum to the left of the axis $\Re(z) = -\delta$.

In Figure \ref{fig:Stab1} we show the paths of eigenvalues corresponding
to the matrix $A+\eps E(\eps)$ for $\eps \in [0,\oeps]$, where $E(\eps)$
indicates the extremizer of the functional computed at $\eps$.


\subsection{Interior two-level algorithm}
\label{subsec:int-stab}
\index{matrix stabilization!interior two-level algorithm}
In the previous subsection, we worked with perturbed matrices that had some eigenvalues of positive real part, and the algorithm moved them to the left until it terminated with a matrix all of whose eigenvalues had real part at most $-\delta$. In an alternative approach, the given matrix $A$ is first perturbed to a non-optimal matrix $A+\eps_0 E_0$ with $\eps_0>0$ and $E_0\in\cS$ of unit Frobenius norm such that all its eigenvalues have real part smaller than $-\delta$. For example, this can be achieved by a simple shift if $I\in\cS$. More importantly, the stable perturbed initial matrix
$A+\eps_0 E_0$ may result from a non-optimal matrix stabilization algorithm.
We then reduce the perturbation size while moving the rightmost eigenvalue towards the imaginary axis.

Using a two-level iteration starting from $A+\eps_0 E_0$,
we aim to reduce the perturbation size to find $\eps>0$ for which $A+\eps E$ has some eigenvalue of real part at least $-\delta$ for {\it every} matrix $E$ of Frobenius norm 1. This differs from the problem of computing the distance to instability, where the aim was to find the smallest perturbation size $\eps$ for which $A+\eps E$ (with $A$ having only eigenvalues of negative real part) has eigenvalues of nonnegative real part for {\it some} matrix $E$ of Frobenius norm~1. 

As in the algorithm for computing the distance to instability of Section~\ref{sec:dist-instab},  the target eigenvalue $\lambda(M)$ of a matrix $M$ is chosen as an eigenvalue of maximal real part.  The function to be minimized is now $f(\lambda,\clambda)=+\,\Re\, \lambda$, whereas it was
$f(\lambda,\clambda)=-\,\Re\, \lambda$ for computing the distance to instability of a Hurwitz matrix. 
Schematically, the interior two-level algorithm proceeds as follows.
\index{two-level iteration}

\begin{itemize}
\item {\bf Inner iteration:\/} Given $\eps>0$, we aim to compute a  matrix $E(\eps) \in\cS$  
of unit Frobenius norm  such that the rightmost eigenvalue of $A + \eps E$
is minimized, i.e. 
\begin{equation} \label{E-eps-stab-int}
E(\eps) = \arg\min\limits_{E\in \cS, \| E \|_F = 1}  \Re \,\lambda(  A + \eps E ).
\end{equation}
(Note that for computing the distance to instability, we maximized the  same functional; see \eqref{E-eps-stab-radius}.)

\item {\bf Outer iteration:\/} We compute $\oeps>0$ as the solution of the one-dimensional equation
\begin{equation} \label{eq:zero-stab-int}
\phi(\oeps)= -\delta,
\end{equation}
where $\phi(\eps)=  \Re\,\lambda( A + \eps E(\eps))$.
\end{itemize}

This yields the nearest perturbed matrix $A + \oeps E(\oeps)$ with all eigenvalues of real part at most~$-\delta$ and the perturbation matrix in the structure space $\cS$.

In the inner iteration, the structured eigenvalue optimization problem \eqref{E-eps-stab-int} is of the class studied in Section~\ref{sec:proto-structured} and is solved with 
a rank-1 matrix differential equation as derived there.

In the outer iteration we compute the optimal perturbation size $\oeps$ by the Newton--bisection algorithm of Section~\ref{subsec:Newton--bisection} for $f(\lambda,\clambda)=\Re\, \lambda$, using the derivative formula of Theorem~\ref{chap:two-level}.\ref{thm:phi-derivative} with $G(\eps)=\Pi^\cS ( x(\eps)y(\eps)^*)$; in particular, $G=xy^*$ in the complex unstructured case, and $G=\Re(xy^*)$ in the real unstructured case.

\subsubsection{Numerical example}
\label{sec:ill1}

Consider the matrix
\begin{equation}
A = \left( \begin{array}{rrrrrrrrrr}
     0  &  1  &  1  &  1  & -1  &  0  & -1  &  0  &  0  &  0 \\
     1  & -1  &  0  &  1  &  1  &  0  &  1  &  0  &  0  &  0 \\
    -1  &  0  & -1  & -1  & -1  &  1  &  1  &  1  &  0  &  0 \\
     1  &  0  &  0  & -1  &  1  & -1  & -1  &  1  &  0  &  0 \\
     0  &  0  & -1  &  1  &  0  &  1  &  1  & -1  &  0  &  0 \\
     0  & -1  &  1  &  1  & -1  &  0  &  0  &  1  &  1  &  0 \\
    -1  &  1  & -1  &  1  &  1  &  0  & -1  &  0  &  1  &  1 \\
     0  &  0  &  1  & -1  & -1  &  1  &  1  &  1  & -1  &  1 \\
     0  &  0  &  0  &  0  &  0  &  0  &  0  & -1  &  1  & -1 \\
     0  &  0  &  0  &  0  &  0  &  0  &  0  &  0  & -1  &  1
		\end{array}
		\right)
		\label{ex:1}
\end{equation}
The matrix $A$ has $6$ eigenvalues with positive real part. We set $\delta=0.1$.

\subsubsection*{Exterior method.}

The stabilized matrix presents $7$ eigenvalues on the parallel axis to the imaginary axis and abscissa $-\delta=-0.1$; its distance from $A$ is $\oeps \approx 2.56$.
For comparison, the stabilized matrix computed by the algorithm of Orbandexivry, Nesterov \& Van Dooren (\cite{OrbNVD13}) has a much larger 
distance $9.02$ and all eigenvalues are located on the imaginary axis. 
\begin{figure}[ht]
\vskip -7cm
\begin{center}
\includegraphics[width=10cm]{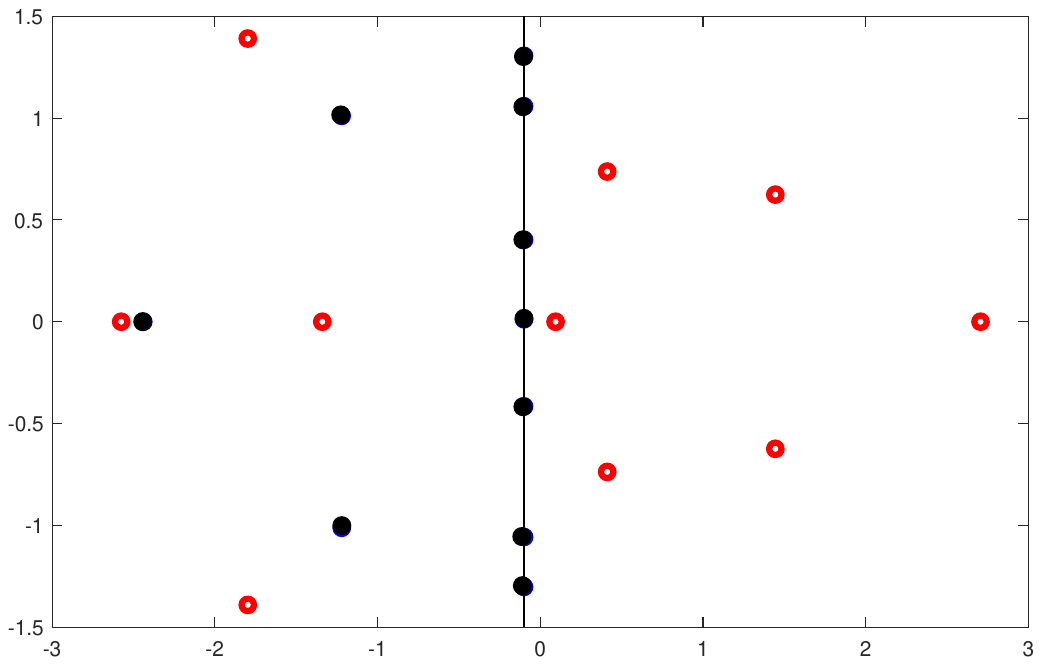} \hskip 0.3cm 
\caption{ 
Spectrum of the matrix
\eqref{ex:1} (red circles) and of the stabilized matrix $A + \oeps E(\oeps)$ (black filled circles) in the exterior method.    
}
\label{fig1-stab}
\end{center}
\end{figure}  
\begin{figure}[h!]
\vskip -7.5cm
\begin{center}
\includegraphics[width=9.9cm]{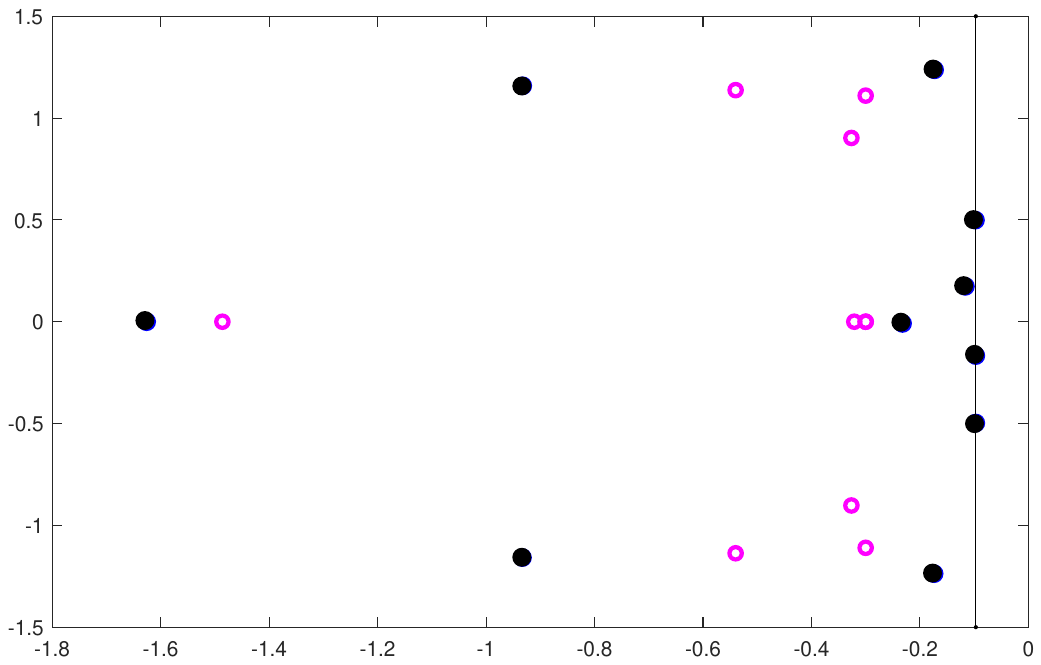} \hskip 0.3cm 
\caption{ 
Spectrum of the initial matrix stabilizing \eqref{ex:1}, computed by
Overton's algorithm (pink circles) and of the 
matrix $A + \oeps E(\oeps)$ (black filled circles) computed by the interior method.    }
\label{fig1-stab-2}
\end{center}
\end{figure}  

\subsubsection*{Interior method.}

We start by computing a stable matrix, with spectrum on the left of the half-line $\Re\, \lambda = -\delta$, using a method developed by Michael Overton, using a BFGS-type method on the penalized functional
\begin{equation*}
    \| X-A \|_F + \rho\,\alpha(X)
\end{equation*}
where $\alpha(X)$ denotes the spectral abscissa of $X$ and $\rho$ is
a penalty parameter.
We obtain in this way an initial matrix for our method, whose distance
from $A$ is approximately $3.14$.
After applying the interior method we obtain a
stabilized matrix which has $3$ eigenvalues on the vertical axis with abscissa $-\delta=-0.1$; its distance from $A$ is $\oeps \approx 2.94$, which is larger than the distance computed by the exterior method, but still a nearby matrix.

The method proposed by Gillis and Sharma (\cite{GilS17}) yields 
a stabilized matrix with approximate distance $1.92$ for $\delta=0$. The matrix has the double eigenvalue $0$ and further $4$ eigenvalues approximately aligned on the imaginary axis. We applied the interior method starting from the 
matrix computed by the algorithm by Gillis and Sharma (\cite{GilS17}). According to our experiments, this matrix turns out to be locally optimal also with respect to the functional of our method.
\bng
This means that if we take it as a starting value for our algorithm, it does not change, since it satisfies the extremality conditions in our approach.
\eng

We refer to Guglielmi \& Lubich (\cite{GL17}) for numerical experiments of matrix stabilization for large sparse matrices, done with a numerical method that is closely related to the exterior method considered here, and for numerical comparisons with other methods.

\subsection{Example: Nearest sparse correlation matrix}
\index{correlation matrix}
As we pointed out in Chapter~\ref{chap:intro}, the problem of computing a nearest correlation matrix (possibly further structured by a prescribed sparsity pattern) to a given symmetric real matrix can be viewed as a special case of the structured matrix stabilization problem considered above. The algorithms of this section can therefore be directly applied to this particular problem.
\medskip

We recall that for a given symmetrix matrix $A$ we look for a minimal-norm perturbation of the form $D + \Delta$ with 
$D=I - \text{diag}(A)$ and $\Delta\in \cS_{\rm cor}$ for 
$$
\cS_{\rm cor}=\{ \Delta \in \R^{n,n}\,:\, \Delta \text{ is symmetric with zero diagonal}\}
$$
such that all eigenvalues of $A+D+\Delta$ are nonnegative.

We consider in addition a sparsity structure of the correction $\Delta$ to $A$. The matrix $A$ might be dense but its correction $\Delta$ is assumed to be sparse. 
This is due to the fact that only a few entries of $A$ are intended to be modified, typically those associated with uncertain or noisy data, while other entries (which are considered sufficiently accurate) are fixed.
This kind of problem occurs in some important cases; for example in the financial market, some   
stock values are updated in time very frequently while some others are updated only a few times in a day. This impacts the accuracy of the computed correlation entries in different ways. In such a case, those entries characterized by poor accuracy can be subject to modifications while those which are more reliable are kept fixed.
Similarly, in statistical applications the data from $m$ observations of $n$ random variables is
collected in an $m \times n$ matrix and it is often the case that some of the observations are
missing (see e.g. Higham and Strabi\'{c} \cite{HigStr16}).
One way to form correlations is via the so-called pairwise deletion method. It calculates the correlation coefficient between a pair of vectors by using only the components available in both vectors simultaneously and the result is a unit diagonal symmetric matrix which is often indefinite.
The problem then is to fill the zero entries so that the matrix becomes positive semidefinite.

Therefore, in this section, we look for
\begin{equation} \nonumber
\Delta\in \cS=\cS_{\rm cor} \cap \cS_{P}
\end{equation}
where $\cS_{P}$ is the space of real symmetric matrices with a given sparsity pattern $P$.

As usual, we write $\Delta = \eps E$, with $E \in \cS_{\rm cor} \cap \cS_{P}$ of Frobenius norm 1 and a given $\eps>0$. We aim to compute a  matrix $E(\eps)$  
of unit Frobenius norm,  such that 
\begin{equation}
\F_\eps(E) = \frac12 \sum_{i=1}^n \, \bigl({}-\lambda_i( A + \eps E)\bigr)_+^2 
\label{Feps-corr}
\end{equation}
is minimized, i.e. 
\begin{equation} \label{E-eps-corr}
E(\eps) = \arg\min\limits_{E\in\cS, \| E \|_F = 1} \F_\eps(E).
\end{equation}
Subsequently we compute the minimal value $\oeps$ such that 
\[
\F_{\oeps}\bigl( E(\oeps) \bigr) = 0.
\]

\paragraph*{Examples.}
We consider first the following simple example of invalid correlation matrix (due to the violation of the unit diagonal constraint) considered by Higham (\cite{Hig02}),
\[
B = \left( \begin{array}{rrrr}
2 & -1 & 0 & 0 \\
-1 & 2 & -1 & 0 \\
0 & -1 & 2 & -1 \\
0 & 0 & -1 & 2
\end{array} \right),
\qquad 
A = \left( \begin{array}{rrrr}
1 & -1 & 0 & 0 \\
-1 & 1 & -1 & 0 \\
0 & -1 & 1 & -1 \\
0 & 0 & -1 & 1
\end{array} \right)
\]
First we change the diagonal of $B$ to the identity matrix, which gives the matrix $A$ we consider; next we run our algorithm.
With full sparsity pattern, meaning 
\[
P = \left( \begin{array}{rrrr}
0 & 1 & 1 & 1 \\
1 & 0 & 1 & 1 \\
1 & 1 & 0 & 1 \\
1 & 1 & 1 & 0
\end{array} \right)
\]
we get the same matrix computed in Higham (\cite{Hig02}) in
$3$ Newton iterations,
\[
A + \oeps E(\oeps) = 
\left( \begin{array}{rrrr}
    1.0000  & -0.8008  &  0.1992  &  0.1231 \\
   -0.8008  &  1.0000  & -0.6777  &  0.1992 \\
    0.1992  & -0.6777  &  1.0000  & -0.8008 \\
    0.1231  &  0.1992  & -0.8008  &  1.0000 \\
    \end{array} \right)
\]
with $\oeps \approx 0.7453$. This provides the same distance $\| B - ( A + \oeps E(\oeps) ) \|_F \approx 2.1343$ as in
Higham \cite{Hig02}.

Next we consider the sparsity pattern in order to maintain the tridiagonal structure,
\[
P = \left( \begin{array}{rrrr}
0 & 1 & 0 & 0 \\
1 & 0 & 1 & 0 \\
0 & 1 & 0 & 1 \\
0 & 0 & 1 & 0
\end{array} \right)
\]
This gives - still in $3$ Newton iterations - 
\[
A + \oeps E(\oeps) = 
\left( \begin{array}{rrrr}
   1.0000 &  -0.6985   &       0   &       0 \\
  -0.6985 &   1.0000   & -0.5121   &       0  \\
        0 &  -0.5121   &  1.0000   & -0.6985 \\
        0 &        0   & -0.6985   &  1.0000
    \end{array} \right)
\]
with $\oeps \approx 0.9163$.
The eigenvalues are:
\[
\lambda_1 = 7.6202 \cdot 10^{-6}, \quad
\lambda_2 = 5.1213 \cdot 10^{-1}, \quad 
\lambda_3 = 1.4879, \quad
\lambda_4 = 2.0000.
\]
Finally let us consider the pattern
\[
P = \left( \begin{array}{rrrr}
0 & 0 & 0 & 1 \\
0 & 0 & 0 & 0 \\
0 & 0 & 0 & 0 \\
1 & 0 & 0 & 0
\end{array} \right)
\]
The eigenvalues of 
\[
A+\eps E = 
\left(
\begin{array}{rrrr}
 1 & -1 & 0 & a \\
 -1 & 1 & -1 & 0 \\
 0 & -1 & 1 & -1 \\
 a & 0 & -1 & 1 \\
\end{array}
\right),
\]
with $a = \displaystyle \frac{\eps \sqrt{2}}{2}$, are given by
\[
\lambda_1 = 1+\sqrt{2}, \quad \lambda_2 = 1-\sqrt{2}, \quad \lambda_3 = 1-\sqrt{a+1}, \quad \lambda_4 = \sqrt{a+1}+1.
\]
This shows that the matrix cannot be corrected to become a correlation matrix by changing only $a$; in this case our algorithm diverges.

Next we consider here a large invalid correlation matrix of stock data considered by
Higham (\cite{Hig02}) and Higham amd Strabi\'{c} (\cite{HigStr16}). The matrix has order 3120 and is dense and of full rank.

We consider a sparse correction with different density, and include preservation of the sparsity pattern as a constraint of the problem.
We give here the results for a sparsity pattern with a density of roughly $70 \%$ non-zero entries.

In our notation we indicate by $m_{-}(\eps)$ the number of negative eigenvalues of the matrix $A + \eps E(\eps)$, and monitor it along the Newton iteration.
\begin{figure}[ht]
\centering
\vspace{-0.2cm}
\includegraphics[width=0.50\textwidth]{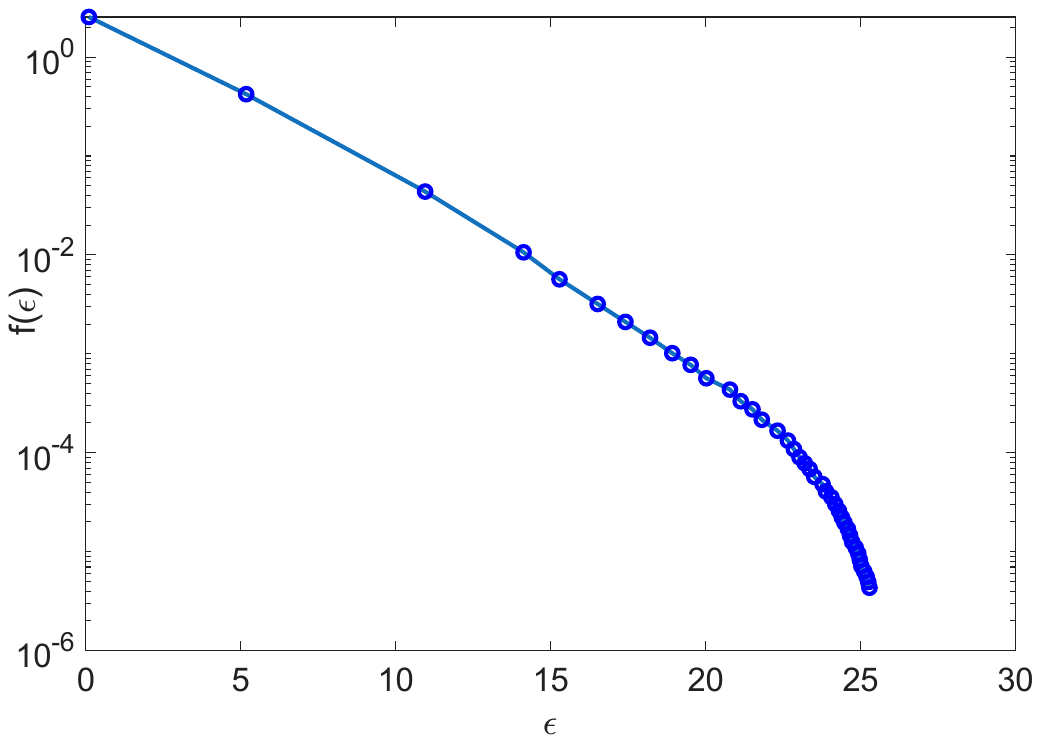}
\hskip -0.1cm 
\includegraphics[width=0.50\textwidth]{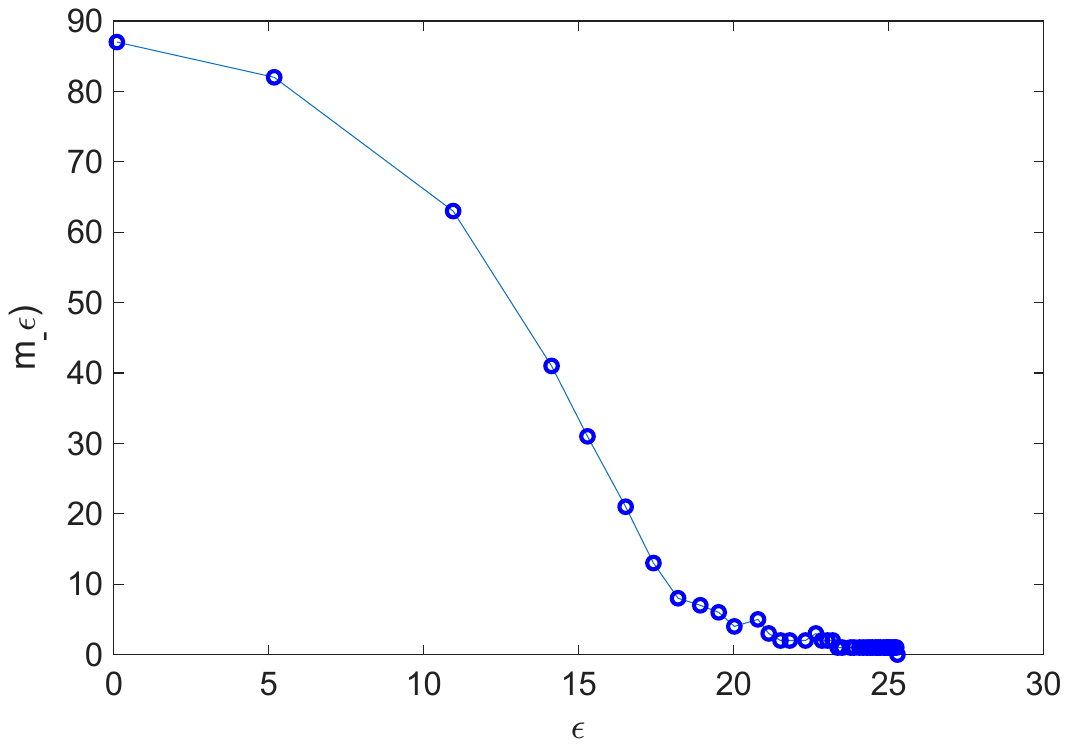}
\vspace{-1cm}
\caption{Left picture: behaviour of the function $\eps_k \mapsto f(\eps_k)=\F_{\eps_k}(E(\eps_k))$ along the Newton iteration. Right picture: behaviour of the function $\eps_k \mapsto m_{-}(\eps_k)$ along the Newton iteration}
\label{fig:Cor}
\end{figure}

In order to accelerate the computation we have adopted the following strategy: with a small target value $\tau>0$, 
we set the functional 
\begin{equation}
\F_\eps(E) = \frac12 \sum_{i=1}^n \, \bigl({}-\lambda_i( A + \eps E)  + 2 \sqrt{\tau})\bigr)_+^2 ,
\label{Feps-corr2}
\end{equation}
which means that we aim to shift the eigenvalues as close as possible to the positive threshold $2 \sqrt{\tau}$. However, we stop the Newton process when all eigenvalues of $A + \eps E(\eps)$ become positive, independently of the effective value of the functional. In this way we are able to avoid oscillating behaviour of the eigenvalues when $\eps$
converges to the minimizer $\oeps$ of $f(\eps)=\F_{\eps}(E(\eps))$.
When all negative eigenvalues approach the value $\sqrt{\tau}$ we have that $\F_\eps(E) \approx  \tfrac{1}{2}c \tau$ where $c$ is the number of eigenvalues on the left of the threshold. 

In Figure \ref{fig:Cor} we show the behaviour of the sequences $f(\eps_k)$ and $m_{-}(\eps_k)$ computed in the Newton process. The iteration is stopped before we meet the tight tolerances set for stopping the Newton iteration, because all eigenvalues leave the negative real axis.



\section{Hamiltonian matrix nearness problems}
\label{sec:Hamilton}
\index{Hamiltonian matrix nearness problem}

In this section we apply the two-level algorithmic approach of Chapter~\ref{chap:two-level} to matrix nearness problems for real Hamiltonian matrices, where `nearest' again refers to the smallest distance in the Frobenius norm:
\begin{itemize}
\item 
{\bf Problem A.} \ {\it Given a Hamiltonian matrix with no eigenvalues on the imaginary axis, find a nearest Hamiltonian matrix having some purely imaginary eigenvalue.}
\item 
{\bf Problem B.} \ {\it Given a Hamiltonian matrix with some eigenvalues on the imaginary axis, find a nearest Hamiltonian matrix such that arbitrarily close to that matrix there exist Hamiltonian matrices with no eigenvalues on the imaginary axis.}
\end{itemize}
Such and related Hamiltonian matrix nearness problems arise 
in the passivation of linear time-invariant control systems 
and in the stabilization of gyroscopic systems; see the references in the notes at the end of this chapter.

We deal here with structured matrix nearness problems with the structure space 
 $\calS$ given as the space $\text{Ham}(\R^{n,n})$ of \emph{real Hamiltonian} matrices (with even dimension $n=2d$), consisting of those matrices $A\in\R^{n,n}$ for which
 \index{Hamiltonian matrix}
$$
JA\, \hbox{ is  real symmetric, where }\,
J=\begin{pmatrix} 0 & I_d \\ -I_d & 0
\end{pmatrix}.
$$
We note that the eigenvalues of a real Hamiltonian matrix lie symmetric to both the real axis and the imaginary axis: with any eigenvalue $\lambda$, also $\clambda,-\lambda,-\clambda$ are eigenvalues. In fact, if $x$ is a left eigenvector of $A$ to the eigenvalue $\lambda$, then $Jx$ is a right eigenvector of $A$ to the eigenvalue $- \clambda$, since
$
AJx = J^{-1}(JA)Jx= -J A^\top J^\top Jx= -J(x^*A)^*=-\clambda Jx.
$

\subsection{Problem A: Moving eigenvalues to the imaginary axis}
\label{subsec:ham-A}

For a real Hamiltonian matrix $M$, we let in the following the target eigenvalue $\lambda(M)$ be the eigenvalue of minimal real part in the first quadrant $\{ \lambda \in \C\,:\,\Re\,\lambda\ge 0, \, \Im\,\lambda\ge 0\}$.
(If this eigenvalue is not unique, we choose the one with minimal imaginary part.) We follow the two-level approach of Chapter~\ref{chap:two-level}:
\index{two-level iteration}
\begin{itemize}
\item {\bf Inner iteration:\/} Given $\eps>0$, we aim to compute a  matrix $E(\eps) \in\cS=\mathrm{Ham}(\R^{n,n})$  
of unit Frobenius norm,  such that 
$ \Re\, \lambda (  A + \eps E)$ 
is minimized:
\begin{equation} \label{E-eps-ham}
E(\eps) = \arg\min\limits_{E \in \mathrm{Ham}(\R^{n,n}), \| E \|_F = 1} \Re\, \lambda (  A + \eps E).
\end{equation}


\item {\bf Outer iteration:\/} We compute the smallest positive value $\oeps$ with
\begin{equation} \label{eq:zero-ham}
\phi(\oeps)= 0,
\end{equation}
where $\phi(\eps)=  \Re \,\lambda(  A + \eps E(\eps) )$.
\end{itemize}

\medskip\noindent
We note that the inner iteration aims to find a leftmost point, in the first quadrant, of the structured $\eps$-pseudospectrum $\Lameps^\calS(A)=\{ \lambda\in \C\,:\, \text{$\lambda$ is an eigenvalue of $A+\eps E$ for some}$ $E\in\cS$ with $\|E\|_F=1 \}$.

We recall from Section~\ref{subsec:proj-structure} that the orthogonal projection $\Pi_\calS$ from $\C^{n,n}$ onto $\calS=\text{Ham}(\R^{n,n})$ is given by
\begin{equation}\label{Pi-Ham-recall}
\Pi^\cS Z = J^{-1}\mathrm{Sym}(\Re\,JZ), \qquad Z\in \C^{n,n}.
\end{equation}
We already know from Section~\ref{subsec:gradient-flow-S} that if
 $\lambda(A+\eps E(\eps))$ is a simple eigenvalue, then the optimizer $E(\eps)$ is a stationary point of the structure-constrained gradient flow given in
(\ref{chap:proto}.\ref{ode-E-S}), viz.,
\begin{align}\nonumber
&\dot E = -G_\eps^\cS(E) + \Re \langle G_\eps^\cS(E), E \rangle E\\
&\text{with the projected gradient }\ G_\eps^\cS(E) = \Pi^\cS (xy^*),   \label{ode-E-S-recall}
\end{align}
where $x$ and $y$ are again left and right eigenvectors of $A+\eps E$ with $\|x\|=\|y\|=1$ and $x^*y>0$. With the given definitions and the orthogonality of $J$,
this becomes
\begin{equation}\label{ode-E-ham}
J\dot E =  -\mathrm{Sym}(\Re\,Jxy^*) + \langle \mathrm{Sym}(\Re\,Jxy^*), JE \rangle JE .
\end{equation}
In a stationary point, $E$ is a real multiple of $G_\eps^\cS(E)=J^{-1}\mathrm{Sym}(\Re\,Jxy^*)$, which is of rank at most 4. The precise rank is as follows.

\begin{theorem} [Rank of optimizers] \label{thm:rank-ham}
For a real Hamiltonian matrix $A$ and $\eps>0$, let $E\in\R^{n,n}$ with $\| E\|_F=1$  be a stationary point of the differential equation \eqref{ode-E-ham}
such that the eigenvalue $\lambda=\lambda(A+\eps E)$ is simple and $\lambda\notin \iu\R$. Then, 
\begin{itemize} 
\item $E$ has  rank $4$ if $\lambda\notin \R$.
\item $E$ has  rank $2$ if $\lambda\in \R$.
\end{itemize}
\end{theorem}
\index{optimizer!rank-4}

\begin{proof} 
If $x$ is a left eigenvector to the eigenvalue $\lambda$, then $Jx$ is a right eigenvector to $-\overline \lambda$. Since $A$ is real, we have in addition that $J\overline x$ and $\overline y$ are eigenvectors to $-\lambda$ and $\overline\lambda$, respectively.

Under the assumption $\lambda\notin \R\cup \iu\R$ we have four different eigenvalues lying symmetric to the real as well as the imaginary axis. The corresponding right eigenvectors $Jx,J\overline x,y,\overline y$ are linearly independent, and so are their real and imaginary parts $Jx_R,y_R,Jx_I,y_I$. Therefore,  the matrix $JG_\eps^\cS(E))$, which is proportional to $JE$, equals
$$
\mathrm{Sym}(\Re\,Jxy^*)
 = \tfrac14 (Jx_R,y_R,Jx_I,y_I) 
\begin{pmatrix} 
0 & 1 & 0 & 0 \\ 
1 & 0 & 0 & 0 \\
0 & 0 & 0 & 1 \\
0 & 0 & 1 & 0
\end{pmatrix} (Jx_R,y_R,Jx_I,y_I)^T,
$$
which has rank 4. Hence, $G_\eps^\cS(E)=J^{-1}\mathrm{Sym}(J\,\Re\,xy^*)$ has rank $4$, and so has its nonzero real multiple $E$.

If $\lambda$ is real and nonzero, then $Jx$ and $y$ are real, and they are linearly independent as eigenvectors to $-\lambda$ and $\lambda$. It follows that $G_\eps^\cS(E)$ and hence $E$ are of rank 2.
\qed
\end{proof}


\subsection{Inner iteration: Rank-4 constrained gradient flow} 
\label{subsec:rank-four-dyn}
\index{gradient flow!rank-4 constrained}
\index{rank-4 matrix differential equation}

In view of Theorem~\ref{thm:rank-ham}, we restrict the gradient flow \eqref{ode-E-S-recall} to Hamiltonian rank-4 matrices, i.e.,
$JE$ will be constrained to lie in the manifold of symmetric rank-4 matrices and hence can be represented as
$$
JE = U S U^\top,
$$
where $U\in \R^{n\times 4}$ has orthonormal columns and $S\in \R^{4\times 4}$ is symmetric and invertible. Like in Section~\ref{subsec:low-rank}, the orthogonal projection onto the tangent space at $JE$ is given as
$$
P_{JE}(Z) = Z - (I-UU^\top)Z(I-UU^\top), \qquad Z\in \R^{n\times n}.
$$
We consider the rank-4 projected gradient flow 
\begin{equation}\label{ode-E-4}
J\dot E = P_{JE} \Bigl( -\mathrm{Sym}(\Re\,Jxy^*) + \langle \mathrm{Sym}(\Re\,Jxy^*), JE \rangle JE \Bigr),
\end{equation}
where $x$ and $y$ are again left and right eigenvectors of $A+\eps E$ with $\|x\|=\|y\|=1$ and $x^*y>0$ to the target eigenvalue $\lambda(A+\eps E)$.
Similarly to Sections~\ref{subsec:rank1-gradient-flow} and~\ref{subsec:rank-r-gradient-flow}, we have the following properties.

\begin{theorem}[Monotonicity]
 \label{thm:monotone-ham}
Let $E(t)$ of unit Frobenius norm
be a solution to the differential equation \eqref{ode-E-4}.
If the eigenvalue $\lambda(t)=\lambda(A+\eps E(t))$ is simple, then
\begin{equation}
\frac{ d }{dt}\, \Re\, \lambda(t)  \le  0.
\label{eq:mon-ham}
\end{equation}
\end{theorem}


\begin{proof} We abbreviate
$G=G_\eps^\cS(E)=J^{-1}\mathrm{Sym}(\Re\,Jxy^*)$ 
and obtain from (\ref{chap:proto}.\ref{eq:deriv-S}) with $f(\lambda,\clambda)=\tfrac12(\lambda+\clambda)=\Re\,\lambda$ and $\kappa=1/(x^*y)>0$, using the orthogonality of $J$,
\begin{align}
\nonumber
\frac1{ \eps\kappa } \,\frac{d}{dt} \, \Re\, \lambda(t) &=  \langle G, \dot E \rangle = \langle JG,J\dot E \rangle
\\
\nonumber
&=  \langle JG,  P_{JE}\bigl(- JG - \langle JG, JE \rangle E\bigr) \rangle
\\
\label{c-s-1-ham}
&=  \Bigl( -\| P_{JE}(JG) \|_F^2 + \bigl( \Re\,\langle P_{JE}(JG), JE \rangle \bigr)^2 \Bigr) \le 0,
\end{align}
where we used  $P_{JE}(JE)=JE$ in the last equality, and $\| JE\|_F=1$ and the Cauchy--Schwarz inequality in
the final inequality.
\qed \end{proof}

\begin{theorem}[Stationary points]
\label{thm:stat-ham}
Let $E$ be a real Hamiltonian rank-4 matrix of unit Frobenius norm and suppose that $P_{JE}(JG_\eps^\cS(E))\ne 0$. If $E$ is a stationary point of the projected differential equation \eqref{ode-E-4}, then $E$ is  also a stationary point of the differential equation \eqref{ode-E-S-recall}.
\end{theorem}
\index{stationary point}

\begin{proof} The proof is similar to the proof of Theorem~\ref{chap:proto}.\ref{thm:stat-r}.
We show that $E$ is a real multiple of $G_\eps^\cS(E)$. By (\ref{chap:proto}.\ref{stat-S}), $E$ is then a stationary point of the differential equation (\ref{chap:proto}.\ref{ode-E-S}), which is the same as~\eqref{ode-E-S-recall}.

For a stationary point $E$ of \eqref{ode-E-4}, we must have equality in the estimate of the previous proof, which shows that $P_{JE}(JG)$  (with $G=G_\eps^\cS(E)$) is a nonzero real multiple of $E$. Hence, in view of $P_{JE}(JE)=JE$, we can write $G$ as
$$
G=\mu E + W, \quad\text{ where $\mu\ne 0$ and $P_{JE}(JW)=0$.}
$$
With $JE=USU^\top$ as above, we then have
 $$
 JW=JW-P_{JE}(JW)= (I-UU^\top)JW(I-UU^\top).
 $$
 Since $G$ is of rank at most $4$ and real Hamiltonian, it can be written in the form $JG=XRX^\top$, where $X\in \R^{n,4}$ has orthonormal columns and $R\in \R^{4,4}$.
So we have
 $$
XRX^\top = \mu USU^\top + (I-UU^\top)JW(I-UU^\top).
 $$
 Multiplying from the right with $U$ yields $X(RX^\top U) = \mu US$, which shows 
 that $X$ has the same range  as $U$. Hence, $JG$ has the same range as $JE$, which implies that $P_{JE}(JG)=JG$. Since we already know that $P_{JE}(JG)$ is a nonzero real multiple of $P_{JE}(JE)=JE$, it follows that $G$ is the same real multiple of $E$. Hence $E$ is a stationary point of  \eqref{ode-E-S-recall}.
\qed \end{proof}

 We further remark that for {\it real} simple eigenvalues $\lambda$ we have an analogous rank-2 dynamics.
 
\subsubsection*{A robust integrator.}
The following time-stepping method is an adaptation to \eqref{ode-E-4} of the low-rank integrator of Ceruti \& Lubich (\cite{CeL22}), similar to the integrator in Section~\ref{subsec:low-rank-integrator}. It first updates the basis matrix $U$ with orthonormal columns and then computes an update of the symmetric $4\times 4$ matrix $S$ by a Galerkin approximation to the differential equation \eqref{ode-E-4} in the updated basis. This integrator is robust to the presence of small singular values, which appear in the case of a target eigenvalue near the real axis, where the rank degenerates from 4 to 2.

	One time step of integration from time $t_k$ to $t_{k+1}=t_k+h$  starting from a factored rank-$4$ matrix 
	$JE_k=U_kS_kU_k^\top$ of unit Frobenius norm computes an updated rank-$r$ factorization $JE_{k+1}=U_{k+1}S_{k+1}U_{k+1}^\top$ of unit Frobenius norm as follows.
	
	\begin{enumerate}
		\item 
		Update the basis matrix $ U_k \rightarrow U_{k+1}$:
		\\[1mm]
		Integrate from $t=t_k$ to $t_{k+1}=t_k+h$ the $n \times r$ matrix differential equation
		$$ \dot{K}(t) = -JG_\eps^\cS( K(t) U_k^\top) U_k  , \qquad K(t_k) = U_k S_k.$$
		Perform a QR factorization $K(t_{k+1}) = U_{k+1} {R}_{k+1}$ and compute the $r\times r$ matrix $M= U_{k+1}^\top U_k$.
		\\[-2mm]		
		\item
		Update the symmetric matrix ${S}_k \rightarrow {S}_{k+1}$\,: \\[1mm]
		Integrate from $t=t_k$ to $t_{k+1}$ the $r \times r$ matrix differential equation
		$$ \dot{S}(t) =  -U_{k+1}^\top JG_\eps^\cS(U_{k+1} S(t) U_{k+1}^\top) U_{k+1}, 
		\qquad S(t_k) = \frac{M {S}_k M^\top}{ \| M {S}_k M^\top \|_F},
		$$
		and set ${S}_{k+1} =S(t_{k+1})/\| S(t_{k+1}) \|_F$.
	\end{enumerate} 
The differential equations in the substeps are solved approximately by a step of some standard numerical integrator, e.g.~the explicit Euler method or a low-order explicit Runge--Kutta method such as the second-order Heun method. The stepsize selection is done as in Section~\ref{subsec:low-rank-integrator}, using an Armijo-type line search.


\subsection{Eigenvalues close to coalescence on the imaginary axis}
\label{sec:eps}
This theoretical section serves as a preparation for the algorithm of the outer iteration that will be presented in the next section.
Let $E(\eps)$ of unit Frobenius norm be a local minimizer of the optimization problem \eqref{E-eps-ham}. We let $\lambda(\eps)$  be the eigenvalue of smallest positive real part (and nonnegative imaginary part) of the Hamiltonian matrix 
$$M(\eps):=A+\eps E(\eps)$$ 
and $x(\eps)$ and $y(\eps)$ are corresponding left and right eigenvectors normalized by $\|x(\eps)\|=\|y(\eps)\|=1$ and $x(\eps)^*y(\eps)>0$.
We let $\oeps$ be the smallest value of $\eps$ such that 
$$
\phi(\eps) = \Re\,\lambda(\eps)
$$
becomes zero:
$$
\phi(\eps)>0 \quad\ \text{for } 0 < \eps < \oeps \quad\ \text{ and } \quad \phi(\eps)=0 \quad\text{for $\eps\ge\oeps$ near $\oeps$.}
$$
Under Assumption~\ref{chap:two-level}.\ref{ass:E-eps}, the function $\phi$ is continuously differentiable in a left neighbourhood of $\oeps$ and its derivative is given by Theorem~\ref{chap:two-level}.\ref{thm:phi-derivative}.
In the following we show that under further assumptions, the function $\phi$ 
has a square-root behaviour $\phi(\eps)\sim \sqrt{\oeps-\eps}$ as $\eps \nearrow \oeps$.

\begin{assumption} \label{assumpt-epsstar}
We assume that the limit $M(\oeps):=\lim_{\eps\nearrow\oeps} M(\eps)$ of the Hamiltonian matrices exists and that the purely imaginary eigenvalue
$\lambda(\oeps)=\lim_{\eps\nearrow\oeps} \lambda(\eps)$ of $M(\oeps)$ has algebraic multiplicity two and is defective (that is, the zero singular value of $M(\oeps)-\lambda(\oeps) I$ is simple).
\end{assumption}

By definition of $\oeps$, the eigenvalue $\lambda(\oeps)$ is on the imaginary axis and has even multiplicity
because of the symmetry of the eigenvalues with respect to the imaginary axis. Here we assume multiplicity two. The defectivity appears to be generic (we have no proof for this but observed defectivity in all our numerical experiments).

Under Assumption \ref{assumpt-epsstar}, the eigenvalue $\lambda(\oeps)$ of $M(\oeps)$ is non-derogatory, that is, only a single Jordan block corresponds to this eigenvalue, and hence its left and right eigenspaces are of dimension 1. Since $\lambda(\oeps)$ is a defective 
eigenvalue, the left and right eigenvectors at $\oeps$ are orthogonal to each other: $x(\oeps)^*y(\oeps) =0$.

We need the following result.

\begin{theorem}  [Eigenvectors at coalescence] \label{thm:yJx}
Let $M(\eps)$, $\eps\in[\eps_0,\oeps]$, be a continuous path of real Hamiltonian matrices, and $\lambda(\eps)$ be a path of eigenvalues of $M(\eps)$ that are simple and not purely imaginary for $\eps<\oeps$ and satisfy Assumption~\ref{assumpt-epsstar} at $\oeps$. Under  a nondegeneracy condition on eigenvectors of $M(\eps)$ stated in \eqref{nondeg-condition} below, 
there exist left and right eigenvectors $x(\eps)$ and $y(\eps)$ to the eigenvalue
$\lambda(\eps)$, normalized to unit norm and with 
 $x(\eps)^*y(\eps)>0$ for $\eps < \oeps$, which depend continuously on $\eps$ in the closed interval $[\eps_0,\oeps]$. In particular, the eigenvectors converge for $\eps\nearrow \oeps$.
In the limit we have
\[
y(\oeps) = \pm Jx(\oeps),
\]
where the sign depends on $x(\eps)$ for $\eps$ near $\oeps$.
\end{theorem}

The important fact here is that $y(\oeps)$ is not just a complex multiple of $Jx(\oeps)$, as would easily be obtained from the symmetry of eigenvalues with respect to the imaginary axis, but that it is a {\it real} multiple.

\begin{proof}
By a result of Paige \& Van Loan~(\cite{PaiVL81}), Theorem 5.1, 
the Hamiltonian matrix $M(\eps)$ with no imaginary eigenvalues (for $\eps<\oeps$) admits a real Schur-Hamiltonian decomposition, that is,  there exists an orthogonal symplectic real matrix 
$S(\eps)$ (i.e., $S(\eps)^\top S(\eps)=I$ and $S(\eps)^\top J S(\eps)=J$) for $\eps > \oeps$) that transforms $M(\eps)$ to a block triangular Hamiltonian matrix
\begin{equation}
M_0(\eps) = S(\eps)^{-1} M(\eps) S(\eps)
= \begin{pmatrix} F(\eps) & H(\eps) \\ 0 & -F(\eps)^\top \end{pmatrix} ,
\end{equation}
where $H(\eps)$ is symmetric and $F(\eps)$ is upper quasi-triangular \bng (that is, it has either scalar or $2 \times 2$ blocks on the diagonal). \eng  

For the eigenvalue $\lambda(\eps)$, the left and right eigenvectors of $M_0(\eps)$ are related to those of $M(\eps)$ by
\begin{equation} \label{x-x0}
x_0(\eps) = S(\eps)^\top x(\eps), \qquad y_0(\eps) = S(\eps)^{-1} y(\eps).
\end{equation}
We assume that $x(\eps)$ and $y(\eps)$ are normalized to norm 1 and such that $x(\eps)^*y(\eps)>0$ for $\eps>\oeps$, and hence we have also 
\begin{equation} \label{x0y0pos}
\text{$x_0(\eps)$ and $y_0(\eps)$ are of norm 1 and\ }\
x_0(\eps)^*y_0(\eps)>0.
\end{equation}
We observe that the lower half of the right eigenvector $y_0(\eps)$ to the block triangular matrix $M_0(\eps)$ consists only of zeros 
and we split the eigenvectors  into the upper and the lower $n/2$-dimensional subvectors as
\begin{equation}\label{yxpsr}
y_0(\eps) = \left( \begin{array}{r}  -p(\eps) \\ 0\ \ \end{array} \right), \qquad 
x_0(\eps) = \left( \begin{array}{r} -s(\eps) \\ r(\eps) \end{array} \right) .
\end{equation}
By the Hamiltonian symmetry,  left and right eigenvectors associated with the eigenvalue $-\conj{\lambda(\eps)}$, with positive inner product, are
$\widetilde x_0(\eps) = J y_0(\eps)$ and $\widetilde y_0(\eps) = J x_0(\eps)$, and so we have
\begin{equation}
\widetilde y_0(\eps) =  \left( \begin{array}{r} r(\eps) \\ s(\eps) \end{array} \right), \qquad 
\widetilde x_0(\eps) = \left( \begin{array}{r} 0\ \ \\ p(\eps) \end{array} \right). 
\end{equation}

By compactness, there exists a sequence $(\eps_m)$ with $\eps_m \nearrow \oeps$ as $m\to\infty$ such that $x_0(\eps_m)$, $y_0(\eps_m)$ and $S(\eps_m)$ converge to vectors $x_{0,\star}$, $y_{0,\star}$ of norm 1 and an  orthogonal symplectic real matrix $S_\star$. By the continuity of $M(\cdot)$ and $\lambda(\cdot)$ at $\oeps$, the limit vectors $x_{0,\star}$, $y_{0,\star}$ are then left and right eigenvectors  corresponding to the purely imaginary eigenvalue $\lambda(\oeps)$ of $M(\oeps)$.

By Assumption~\ref{assumpt-epsstar}, the left and right eigenspaces to $\lambda(\oeps)$ are one-dimensional, and so we have that for some complex $\xi,\eta$ of unit modulus,
\begin{equation}\label{y0-x0-limits}
\lim\limits_{m\to\infty} \widetilde{y}_0(\eps_m) = -\eta \lim\limits_{m\to\infty} y_0(\eps_m), \qquad
\lim\limits_{m\to\infty} \widetilde{x}_0(\eps_m) = \xi \lim\limits_{m\to\infty} x_0(\eps_m) .
\end{equation}   
We thus obtain 
\begin{equation}\label{s-to-zero}
\lim_{m\to\infty}  s(\eps_m) = 0
\end{equation} 
and
\begin{equation} \label{eq:limits}
\lim\limits_{m\to\infty}  r(\eps_m) = \eta \lim\limits_{m\to\infty} p(\eps_m), \qquad
\lim\limits_{m\to\infty}  p(\eps_m) = \xi \lim\limits_{m\to\infty} r(\eps_m) ,
\end{equation} 
so that
\begin{equation}
\xi = \bar\eta.
\end{equation}
By \eqref{x0y0pos},
\begin{equation} \label{sp-pos}
s(\eps)^*p(\eps) \ \text{ is real and positive for } \eps > \oeps,
\end{equation}
and in particular, $s(\eps)\ne 0$ for $\eps > \oeps$ (but recall \eqref{s-to-zero}).

Moreover, from the fact that $x_0(\eps)$ is a left eigenvalue of $M_0(\eps)$, we infer that (omitting the argument $\eps$ in the next few lines)
$s^*F=\lambda s^*$ and $-s^*H-r^*F^\top=\lambda r^*$. Multiplying the second equation with $s$ from the right and using the first equation then yields $-s^*Hs=(\lambda+\clambda)\,r^*s$, which shows that
\begin{equation} \label{sr-real}
s(\eps)^* r(\eps) \ \text{ is real for } \eps > \oeps.
\end{equation}
Under the nondegeneracy condition
\begin{equation}\label{nondeg-condition}
\liminf_{\eps\nearrow\oeps} \,\left| \left(\frac{ s(\eps) }{ \|s(\eps)\| }\right)^* \frac{ r(\eps) }{ \|r(\eps)\| } \right| > 0,
\end{equation}
which states that the normalizations of the vectors $s$ and $r$ are not asymptotically orthogonal,
we conclude that there is a subsequence $(\eps_m')$ of $(\eps_m)$ such that the normalized sequence $\bigl(s(\eps_m')/\|s(\eps_m')\|\bigr)$ is convergent and (on noting that $\|r(\eps_m)\|\to 1$ because of \eqref{yxpsr} and \eqref{s-to-zero})
\begin{equation}\label{sp-conv}
\lim_{m\to\infty} \frac{ s(\eps_m')^*r(\eps_m') }{ \|s(\eps_m')\| } \ne 0.
\end{equation}
As \eqref{eq:limits} implies that this nonzero limit equals
$$
\lim_{m\to\infty} \frac{ s(\eps_m')^*r(\eps_m') }{ \|s(\eps_m')\| } =
\eta \lim_{m\to\infty} \frac{ s(\eps_m')^*p(\eps_m') }{ \|s(\eps_m')\| }  
$$
and the two limits in this formula are real by \eqref{sp-pos} and \eqref{sr-real}, 
it follows that $\eta$ is real and hence $\eta$ equals $1$ or $-1$. In view of $\eqref{sp-pos}$ and \eqref{nondeg-condition}, we actually have
\begin{equation}\label{eta-sr}
\eta=\lim_{\eps\nearrow\oeps}\ \mathrm{sign}( s(\eps)^* r(\eps)) = \pm 1,
\end{equation}
which depends only on the left eigenvector $x_0(\eps)$.
As a consequence, we obtain from \eqref{y0-x0-limits} that
\begin{equation} \label{eq:limy0}
y_{0,\star} = -\eta J x_{0,\star}= \mp J x_{0,\star}.
\end{equation}
By \eqref{x-x0} we have
\begin{equation}
x(\eps) =  \left( S(\eps)^\top \right)^{-1} x_0(\eps) =  - J S(\eps) J x_0(\eps), \qquad y(\eps) = S(\eps) y_0(\eps) = \pm S(\eps) J x_0(\eps)
\end{equation}
and therefore the limits $x_\star = \lim_{m\to\infty} x(\eps_m)$ and $y_\star = \lim_{m\to\infty} y(\eps_m)$ exist and satisfy
\begin{equation} \label{eq:limy}
y_\star  = \pm J x_\star.
\end{equation}
We now use once again that by Assumption~\ref{assumpt-epsstar}, the left and right eigenspaces to $\lambda(\oeps)$ are one-dimensional. 
Hence $x_\star$ is a complex multiple of the unique left eigenvector $x(\oeps)$ of norm 1 for which the first nonzero entry is positive.
If we choose the eigenvectors $x(\eps)$ such that their corresponding entry is also nonnegative, then we find that every convergent subsequence
$(x(\eps_m))$ converges to the same limit $x(\oeps)$ as $m\to\infty$, and hence $x(\eps)$ converges to $x(\oeps)$ as $\eps\nearrow\oeps$. To the left eigenvector $x(\eps)$, there corresponds a unique right eigenvector $y(\eps)$ of norm 1 that satisfies
$x(\eps)^*y(\eps)>0$ for $\eps>\oeps$. By \eqref{eq:limy}, the limit of every convergent subsequence $(y(\eps_m))$ converges to
$\pm J \lim_{m\to\infty} x(\eps_m)=\pm Jx(\oeps)$, and hence the limit $y(\oeps):= \lim_{\eps\nearrow\oeps} y(\eps)$ exists, is a right eigenvector of $M(\oeps)$ to the eigenvalue $\lambda(\oeps)$, and it equals
$$
y(\oeps) = \lim_{\eps\nearrow\oeps} y(\eps) = \pm J \lim_{\eps\nearrow\oeps} x(\eps) = \pm J x(\oeps),
$$
which completes the proof.
\qed
\end{proof}

\begin{remark} If additionally $M(\eps)$, $\eps\in[\eps_0,\oeps]$, is continuously differentiable and 
$\Re\, x(\oeps)^* M'(\oeps) y(\oeps)\ne 0$, then Theorem~\ref{thm:yJx} implies that the eigenvalue $\lambda(\eps)$ approaches the imaginary axis from the normal direction (i.e. horizontally in the complex plane). This is because then we have, with
$\kappa(\eps)=1/(x(\eps)^*y(\eps))>0$,
\begin{align*}
\Im\,\frac{\lambda'(\eps)}{\kappa(\eps)}&=\Im\, x(\eps)^* M'(\eps) y(\eps)
\to \Im\, x(\oeps)^* M'(\oeps) y(\oeps) \\&=
\Im\, (Jx(\oeps))^* JM'(\oeps) y(\oeps) = \pm\, \Im\, y(\oeps)^* JM'(\oeps) y(\oeps) =0
\end{align*}
by the symmetry of $JM'(\oeps)$. By assumption,
$$
\Re\,\frac{\lambda'(\eps)}{\kappa(\eps)}=\Re\, x(\eps)^* M'(\eps) y(\eps) \to \Re\, x(\oeps)^* M'(\oeps) y(\oeps)\ne 0.
$$
Hence, $\Im\,\lambda'(\eps)\,/\,\Re\,\lambda'(\eps)\to 0$ as $\eps\nearrow\oeps$. This can, however, not be concluded when
$M'(\eps)$ has no limit at $\oeps$ and $\| M'(\eps)\|\to \infty$ as $\eps\nearrow\oeps$.
\end{remark}

We are now in a position to characterize the asymptotic behaviour of the function $\phi(\eps)=\Re\,\lambda(\eps)$ as $ \eps \nearrow \oeps$, in the situation described at the beginning of this section.

\begin{theorem}[Square root asymptotics]\label{thm:sqrt} 
Under Assumptions~\ref{chap:two-level}.\ref{ass:E-eps} and~\ref{assumpt-epsstar} and the nondegeneracy condition \eqref{nondeg-condition}, and under the further condition that the eigenvalue $\lambda(\eps)$ of the Hamiltonian matrix $M(\eps)=A+\eps E(\eps)$ does not approach the imaginary axis tangentially as $\eps\nearrow\oeps$,
we have
$$
\Re \,\lambda(\eps)= \gamma \,\sqrt{\oeps-\eps} \;(1+o(1)) \quad\text{ as }\ \eps \nearrow \oeps
$$
for some positive constant $\gamma$.
\end{theorem}

\begin{proof} We split the proof into four parts (a)-(d).
\\
(a) For $\eps<\oeps$, let $x(\eps)$ and $y(\eps)$ be left and right eigenvectors of $M(\eps)=A+\eps E(\eps)$ to the simple eigenvalue $\lambda(\eps)$, of unit norm and with
$x(\eps)^* y(\eps) > 0$ for $\eps<\oeps$ and normalized such that their limits for $\eps \nearrow \oeps$ exist according to Theorem~\ref{thm:yJx}.
We consider the nonnegative function 
$$
\vartheta(\eps) := \frac1{\kappa(\eps)}= x(\eps)^* y(\eps) > 0 \ \text{ for }\ \eps \in (\eps_0, \oeps), \qquad
\vartheta(\oeps) = 0.
$$ 
To compute the derivative of $\vartheta$, we use Theorem~\ref{chap:appendix}.\ref{thm:eigvecderiv} for the left and right eigenvectors of norm 1 and with positive inner product,
\begin{align*}
x'(\eps)^* & =   - x(\eps)^* M'(\eps) Z(\eps) + \Re\bigl( x(\eps)^* M'(\eps) Z(\eps) x(\eps)\bigr) x(\eps)^*
\\
y'(\eps) & =  - Z(\eps) M'(\eps) y(\eps) + \Re\bigl(y(\eps)^* Z(\eps) M'(\eps) y(\eps) \bigr) y(\eps), 
\nonumber
\end{align*} 
where $Z(\eps)$ is the group inverse of $N(\eps):=M(\eps) - \lambda(\eps) I$. 
Since Theorem~\ref{chap:appendix}.\ref{thm:Ginv} shows that $x(\eps)^* Z(\eps)=0$ and $Z(\eps) y(\eps) = 0$,  these 
formulas imply
\begin{equation}
\vartheta'(\eps)  =  \Re\Bigl(x(\eps)^* M'(\eps) Z(\eps) x(\eps) +
y(\eps)^* Z(\eps) M'(\eps) y(\eps)\Bigr) \vartheta(\eps) .
\label{eq:derdelta}
\end{equation} 
By Theorem~\ref{chap:appendix}.\ref{thm:Ginv}, the group inverse is related to
the pseudoinverse $N(\eps)^\dag$ by the formulas 
\begin{eqnarray}
Z(\eps) & = & \frac{1}{\vartheta(\eps)^2} \widehat{Z}(\eps)
\nonumber
\\[1mm]
\widehat{Z}(\eps) & = & \bigl( \vartheta(\eps) I - y(\eps) x(\eps)^* \bigr) 
N(\eps)^\dag
\bigl( \vartheta(\eps) I - y(\eps) x(\eps)^* \bigr).
\label{eq:Ghat}
\end{eqnarray} 
By Assumption \ref{assumpt-epsstar}, the second smallest
singular value $\sigma_{n-1}(\eps)$ of $N(\eps)$
does not converge to zero. Therefore, $
N(\eps)^\dag$ has a finite limit as $\eps \nearrow \oeps$.
We thus have
\begin{eqnarray}
\widehat{Z}(\eps) & = & 
y(\eps) x(\eps)^* N(\eps)^\dag y(\eps) x(\eps)^* + O(\vartheta(\eps)) 
\\
&=& \nu(\eps) y(\eps)x(\eps)^* + O(\vartheta(\eps)) 
\nonumber
\end{eqnarray}
with the factor
$$
 \nu(\eps)  :=  x(\eps)^* N(\eps)^\dag y(\eps),
$$
Furthermore, we set
$$
\mu(\eps) := x(\eps)^* M'(\eps) y(\eps).
$$
We insert the expression for the group inverse $Z(\eps)$ into \eqref{eq:derdelta} and
note the identities 
$N^\dagger(\eps)x(\eps)=0$ and $y(\eps)^*N^\dagger(\eps)=0$, which follow from $x(\eps)^*N(\eps)=0$ and $N(\eps)y(\eps)=0$, respectively.
We then obtain
\begin{align}\nonumber
\vartheta'(\eps)\vartheta(\eps) &= \Re\Bigl(x(\eps)^* M'(\eps) \widehat Z(\eps) x(\eps)  +
y(\eps)^* \widehat Z(\eps) M'(\eps) y(\eps)  \Bigr)
\\
&=  \Re\Bigl(2\nu(\eps) \mu(\eps) + O(\vartheta(\eps)\mu(\eps))\Bigr).
\label{vartheta-prime}
\end{align}
(b) We now study the limit behaviour of $\nu(\eps)$ as $\eps\nearrow\oeps$.
By Theorem~\ref{thm:yJx},  the limits of the left and right eigenvectors for $\eps \nearrow \oeps$ exist and satisfy
$y(\oeps)= \pm Jx(\oeps)$ and further $x(\oeps)^* y(\oeps)=0$.
Since $JN(\oeps)$ is a hermitian matrix, we therefore obtain
\begin{align*}
 \nu(\oeps) &= x(\oeps)^* N(\oeps)^\dag y(\oeps) = x(\oeps)^* (JN(\oeps))^\dag  J y(\oeps) 
 \\
&= \mp \, x(\oeps)^* (JN(\oeps))^\dag x(\oeps) \in \R.
\end{align*}
We next show that $ \nu(\oeps)\ne 0$.
Since $x(\oeps) \,\bot \,y(\oeps)$ and since $x(\oeps)$ spans the nullspace of $N(\oeps)^*$ (by the defectivity condition in Assumption~\ref{assumpt-epsstar}), we obtain
$$
 y(\oeps) \in {\rm Ker} (N(\oeps)^*)^\perp =
{\rm Range} \left( N(\oeps) \right).
$$
Hence there exists $z_1$ such that $y(\oeps)=N(\oeps) z_1$.
Assume, in a proof by contradiction, $N(\oeps)^\dag y(\oeps) \, \bot \, x(\oeps)$, which means 
$$
N(\oeps)^\dag y(\oeps) \in {\rm Range} \left( N(\oeps) \right).
$$
Hence there exists $z_2$ such that $N(\oeps)^\dag y(\oeps)=N(\oeps) z_2$.

%
Multiplying this equation with $N(\oeps)^2$ we obtain (omitting the argument $\oeps$ in the following)
$$
N^3z_2 = N^2 N^\dag y = N \, N N^\dag N z_1 = NNz_1 = Ny = 0.
$$
The null-space of $N^3$ is two-dimensional by Assumption~\ref{assumpt-epsstar} and contains the two nonzero vectors $y$ and $N^\dag y$, since $Ny=0$ and $N^2 N^\dag y =0$. Note that $N^\dag y\ne 0$ because otherwise $y$ would be in the nullspace of $N^\dag$, which is the nullspace of $N^*$, which contradicts the above observation that $y\ne 0$ is in the orthogonal complement of the nullspace of $N^*$.
Moreover, $y$ and $N^\dag y$ are linearly independent, since otherwise the relation $y= c N^\dag y $ would yield, on multiplication with $N$, that
$$
0=Ny = cNN^\dag y = c NN^\dag Nz_1 =c N z_1 = c y,
$$
which contradicts $y\ne 0$. Therefore, the null-space of $N^3$ is {\it spanned} by
$y$ and $N^\dag y$, and since we have shown that it contains $z_2$, we obtain
$$
z_2 = c_1 y + c_2 N^\dag y.
$$ 
Multiplying this equation with $N$ then gives 
$$
N^\dag y = Nz_2 = c_2 NN^\dag Nz_1 = c_2 Nz_1 = c_2 y,
$$
which contradicts the linear independence of $y$ and $N^\dag y$. We have thus  led the assumption $N(\oeps)^\dag y(\oeps) \ \bot \ x(\oeps)$ to a contradiction.
Therefore, $\nu(\oeps) \neq 0$. So we have shown that
\begin{equation}\label{nu-star}
\nu_\star:=
\lim_{\eps \nearrow \oeps} \nu(\eps) \quad\text{exists and is real and nonzero.}
\end{equation}
(c) We next study the limit behaviour of $\mu(\eps)$ as $\eps\nearrow\oeps$.
By Theorem~\ref{chap:appendix}.\ref{thm:eigderiv} we have 
$$
\mu(\eps)=  \lambda'(\eps)\vartheta(\eps). 
$$
By Theorem~\ref{chap:two-level}.\ref{thm:phi-derivative} with the objective function $\phi(\eps)=\Re\,\lambda(\eps)$ and the gradient $G(\eps)=\mathrm{Sym}(\Re\, Jx(\eps)y(\eps)^*)$, we thus have 
$$
\Re\,\mu(\eps)=-\|G(\eps)\|_F.
$$
Since $y(\oeps)= \pm Jx(\oeps)$ by Theorem~\ref{thm:yJx}, we have in the limit $\eps \nearrow \oeps$ that
$$
G(\oeps)= \mathrm{Sym}(\Re\, Jx(\oeps)y(\oeps)^*)=\pm \Re\, y(\oeps)y(\oeps)^* \ne 0.
$$
By assumption, $\lambda(\eps)$ does not approach the imaginary axis tangentially, and hence we have  
\begin{equation}\label{Im-mu}
\frac{|\Im\,\mu(\eps)|}{|\Re\,\mu(\eps)|}=\frac{|\Im\,\lambda'(\eps)|}{|\Re\,\lambda'(\eps)|} \le C
\end{equation}
for some constant $C$ independent of $\eps\in(\eps_0,\oeps)$. This implies that $|\Im \,\mu(\eps)|$ is bounded independently of $\eps$.
In the following we let 
\begin{equation}\label{mu-G}
\rho_\star:= -\lim_{\eps \nearrow \oeps} \Re\,\mu(\eps) = \|G(\oeps)\|_F >0.
\end{equation}
(d) From \eqref{nu-star}--\eqref{mu-G} we conclude that the right-hand side of \eqref{vartheta-prime} has a nonzero finite real limit as $\eps \nearrow \oeps$. 
So we have 
$$
\frac d{d\eps} \vartheta(\eps)^2 = 2 \vartheta'(\eps)\vartheta(\eps) = -4\rho_\star\nu_\star (1+o(1)) \quad\text{ as }\ \eps\nearrow\oeps.
$$
Integrating this relation and using $\vartheta(\oeps)^2=0$ yields
$$
\vartheta(\eps)^2 = (\oeps-\eps) \,4\rho_\star\nu_\star (1+o(1)).
$$
This further allows us to conclude that $\nu_\star$ is not only nonzero and real but actually positive.
We  recall that $\vartheta(\eps)>0$ for $\eps<\oeps$ and take the square root to obtain
\begin{equation}\label{vartheta-asymptote}
\vartheta(\eps) = \sqrt{\oeps-\eps} \ 2\sqrt{ \rho_\star\nu_\star} \,(1+o(1)).
\end{equation}
On the other hand, since $\mu(\eps)=  \lambda'(\eps)\vartheta(\eps)$, we find
$$
\Re\, \lambda'(\eps) = \frac{\Re\, \mu(\eps)}{\vartheta(\eps)} = \frac{-\rho_\star}{\vartheta(\eps)} (1+o(1)) .
$$
Using  \eqref{vartheta-asymptote} and setting $\gamma=\sqrt{\rho_\star/\nu_\star}$, this yields
\begin{equation}\label{phi-prime}
\Re\, \lambda'(\eps) = - \frac{\gamma}{2 \sqrt{\oeps - \eps}} (1+o(1)),
\end{equation}
and integration then implies the stated result for $\Re\, \lambda(\eps)$.
\qed
\end{proof}

\subsection{Outer iteration: Square root model and bisection}
\label{subsec:outer-it-sqrt}
For a small positive parameter $\delta$, we aim to find $\eps_\delta$ as the minimal number such that 
\begin{equation}\label{eps-delta}
\Re\,\lambda(\eps_\delta)=\delta.
\end{equation}
We can use the Newton--bisection method of Section~\ref{subsec:Newton--bisection} to compute $\eps_\delta$, but for small $\delta$, this can lead to many Newton step rejections and bisection steps. In the situation of Theorem~\ref{thm:sqrt},
using a square root model of $\phi(\eps)=\Re\,\lambda(\eps)$ appears more appropriate. 
The algorithm falls back to simple bisection if the local square root model fails.

For $\eps\nearrow \oeps$,  we have in the expected situation of Theorem~\ref{thm:sqrt} the square-root behaviour for $\phi(\eps)=\Re \,\lambda(\eps)$ (and by \eqref{phi-prime} for $\phi'(\eps)$)  
\begin{eqnarray}
\begin{array}{rcl}
\phi(\eps) & = & \gamma \sqrt{\oeps - \eps} \ (1+o(1)) 
\\[2mm]
\phi'(\eps) & = & -\displaystyle{\frac{\gamma}{2 \sqrt{\oeps - \eps}}} \ (1+o(1)) .
\end{array}
\label{eq:eps-sqrt}
\end{eqnarray}
For an iterative process, given $\eps_k$, we use that $\phi'(\eps_k)=-\kappa(\eps)\|G(\eps)\|_F$ by
Theorem~\ref{chap:two-level}.\ref{thm:phi-derivative} and solve (\ref{eq:eps-sqrt}) for $\gamma$ and 
$\oeps$, ignoring the $o(1)$ terms.
We denote the solution as $\gamma_k$ and $\widehat \eps_k$, i.e.,
\begin{eqnarray}
\gamma_k & = & \sqrt{-2 \phi(\eps_k) \phi'(\eps_k)}, \qquad
\widehat \eps_k = \eps_k - \frac{\phi(\eps_k)}{2 \phi'(\eps_k)} .
\label{eq:stepk}
\end{eqnarray}
As a substitute for the equation $\phi(\eps)=\delta$, we solve the equation $\gamma_k \sqrt{\widehat\eps_k-\eps_{k+1}}=\delta$ for $\eps_{k+1}$, which yields
\begin{equation}
\eps_{k+1}  =  \widehat \eps_k + {\delta^2}/{\gamma_k^2} \ge \eps_k.
\label{eq:stepk2}
\end{equation} 
Algorithm \ref{alg:problemA} is based on these formulas.
Here, tol is a tolerance that
controls the desired accuracy of the computed optimal $\eps$ (not to be chosen too small).

\bng
As we have discussed in the context of the Newton--bisection method, there is no guarantee of convergence. The inner iteration could find locally optimal values that lead to invalid bisection updates of lower and upper bounds. A version with monotonically increasing perturbation sizes $\eps_k$ analogous to the monotone Newton--bisection method but with the square-root model is expected to show more favourable behaviour. 
\bcltwo
The HEC method, which is supposed to converge with monotonically decreasing perturbation sizes $\eps_k$, does not appear applicable here, due to the fact that the function $\phi$ is identically zero for $\eps \ge \oeps$.
\ecltwo
\eng
 
\medskip
\begin{algorithm}[H]  \label{alg:problemA}
\DontPrintSemicolon
\KwData{$\delta$, ${\rm tol}$, $\theta$ (default $0.8$), 
and $\eps_0$ (such that $\phi(\eps_0) > {\rm tol}$)} 
\KwResult{$\widehat \eps_{\delta}$, $E(\widehat \eps_{\delta})$}
\Begin{
\nl Set {\rm Reject} = {\rm False} and $k=0$\;
\nl \While{$|\phi(\eps_k) - \delta| \ge {\rm tol}$}{
\nl \eIf{${\rm Reject} = {\rm False}$} {
\nl	Set $\widetilde\eps = \eps_{k}$, \ $\widetilde{\theta}=\theta$\; 
    Compute $\gamma_k$ and $\widehat \eps_k$ by (\ref{eq:stepk})\;
\nl Set $\eps_{k+1} = \widehat \eps_k + {\delta^2}/{\gamma_k^2}$\;
    }{
    Set $\eps_{k+1} = \widetilde\theta\,\eps_{k}+ (1-\widetilde\theta)\,\widetilde\eps$\;
		Set $\widetilde\theta = \theta\widetilde\theta$\;
    }
\nl Set $k=k+1$\;
\nl Compute $\phi(\eps_k)$ by solving the rank-4 differential equation \eqref{ode-E-4} with initial datum $E(\eps_{k-1})$
into a stationary point $E(\eps_k)$ as in Section~\ref{subsec:rank-four-dyn}\;
\nl Compute $\phi'(\eps_k)$ by (\ref{chap:two-level}.\ref{eq:dereps})\;		
\nl \eIf{$\phi(\eps_k) < {\rm tol}$}{ 
    Set ${\rm Reject} = {\rm True}$} {
    Set ${\rm Reject} = {\rm False}$ }
}
\nl Return $\widehat \eps_\delta = \eps_k$\;
}
\caption{Basic algorithm for computing the optimal perturbation for small~$\delta$}
\end{algorithm}
\medskip

\subsection{Problem B: Eigenvalues leaving the imaginary axis}
\label{subsec:ham-B}

We describe two complementary approaches to Problem B. The first approach moves eigenvalues on the imaginary axis, and the second approach moves eigenvalues off the imaginary axis.

\subsubsection{Moving eigenvalues on the imaginary axis to coalescence}
Because of the symmetry of eigenvalues of Hamiltonian matrices with respect to the imaginary axis, paths of eigenvalues can leave the imaginary axis only at multiple eigenvalues. Given a Hamiltonian matrix with some simple eigenvalues on the imaginary axis, it is thus of interest to find its distance to the nearest matrix where two previously adjacent eigenvalues coalesce. This problem is addressed by an extension of the two-level approach considered before.

Let $A$ be a real Hamiltonian matrix with a pair $\lambda_1(A)$ and $\lambda_2(A)$ of adjacent eigenvalues on the imaginary axis, with
$\Im\,\lambda_2(A)>\Im\,\lambda_1(A)$. In the inner iteration we determine, for a fixed perturbation size $\eps>0$, a real Hamiltonian matrix  $E$ of Frobenius norm~$1$ such that the functional
$$
\F_\eps(E)= \Im\,\lambda_2(A+\eps E)-\Im\,\lambda_1(A+\eps E)
$$
is minimized. As in previous sections, this minimization is carried out with a constrained gradient flow. Along a path $E(t)$ of real Hamiltonian matrices with simple eigenvalues $\lambda_k(t)=\lambda_k(A+\eps E(t))$  ($k=1,2$) on the imaginary axis,
corresponding left and right eigenvectors $x_k(t),y_k(t)$ of unit norm with positive inner product, and the eigenvalue condition numbers $\kappa_k(t)=1/(x_k(t)^*y_k(t))$ we find (omitting the ubiquitous argument~$t$)
\begin{align*}
\frac{d}{dt} \F_\eps(E(t)) &= \Im (\dot \lambda_2 - \dot \lambda_1) 
= \Im(\kappa_2 x_2^* \dot E y_2 - \kappa_1 x_1^* \dot E y_1)
\\
&=\langle \Im(\kappa_2 x_2 y_2^* - \kappa_1 x_1 y_1^*),\dot E\rangle 
\\
&= \langle G_\eps(E), \dot E \rangle
\end{align*}
with the Hamiltonian gradient
$$
G_\eps(E) = \Pi_{\cS} \,\Im(\kappa_2 x_2 y_2^* - \kappa_1 x_1 y_1^*) = 
J^{-1} \mathrm{Sym} \bigl( \Im (\kappa_2 Jx_2 y_2^* - \kappa_1 Jx_1 y_1^*)\bigr),
$$
which has rank at most $8$. The corresponding norm-constrained gradient system is then again
$$
\dot E = - G_\eps(E) + \langle G_\eps(E),E \rangle E,
$$
along which $\F_\eps(E(t))$ decreases monotonically. In a stationary point, $E$ is a real multiple of $G_\eps(E)$, which is of rank at most $8$. We can then solve numerically the rank-8 constrained gradient system into a stationary point in the same way as we did with the rank-4 system in Section~\ref{subsec:rank-four-dyn}.\index{gradient flow!rank-8 constrained}
In the outer iteration we aim to determine the smallest zero $\oeps$ of $\phi(\eps)=\F_\eps(E(\eps))$, where $E(\eps)$ is the minimizer corresponding to the perturbation size~$\eps$.
This is again done by a Newton--bisection algorithm (or using a square root model and bisection) as discussed before.

We note, however, that a coalescence on the imaginary axis does not guarantee that the coalescent eigenvalues can be moved off the imaginary axis by an arbitrarily small further perturbation; see Mehrmann \& Xu (\cite{MeX08}), Theorem~3.2. Moreover, mere coalescence on the imaginary axis does not give an answer to the related problem of finding a smallest perturbation that moves the adjacent imaginary eigenvalues to a prescribed positive distance $\delta$ to the imaginary axis.

\subsubsection{Moving non-imaginary eigenvalues of perturbed Hamiltonian matrices back to coalescence on the imaginary axis}
In a complementary approach, we first perturb the given real Hamiltonian matrix $A$, which is assumed to have some eigenvalues on the imaginary axis, to another Hamiltonian matrix $A_0=A+\eps_0 E_0$ (with $\| E_0\|_F=1$) that has no eigenvalues on the imaginary axis, but which is not the one that is closest to $A$. With a sufficiently large perturbation size $\eps_0$, this is always possible; just take $A_0=\text{blockdiag}(B,-B^\top)$, where $B\in\R^{d,d}$ is an arbitrary matrix having no purely imaginary eigenvalues. For example, one might choose $B$ as the left upper block of $A$, if this has no imaginary eigenvalues, or else slightly shifted to have no imaginary eigenvalue. We remark that in our numerical experiments, the choice of $A_0$ was not a critical issue. Starting from $A_0$, we reduce the perturbation size to $\eps<\eps_0$ and in this way drive eigenvalues back to the imaginary axis.

We aim to find the largest perturbation size $\eps$ for which $A+\eps E$ has some eigenvalue on the imaginary axis for {\it every} matrix $E$ of Frobenius norm 1. This differs from Problem~A, where the aim was to find the smallest perturbation size $\eps$ for which $A+\eps E$ (with $A$ having no purely imaginary eigenvalues) has eigenvalues on the imaginary axis for {\it some} matrix $E$ of Frobenius norm~1. 

As in Section~\ref{subsec:ham-A},  the target eigenvalue $\lambda(M)$ of a Hamiltonian matrix $M$
is taken as an eigenvalue of minimal real part in the first quadrant.  
In the inner iteration we use a rank-4-constrained gradient system to compute, for a given perturbation size $\eps>0$, a  real Hamiltonian matrix $E(\eps)$  
of Frobenius norm 1  such that 
$ \Re\, \lambda (  A + \eps E)$ 
is (locally) {\it maximized} (as opposed to {\it minimized} for Problem A):
\begin{equation} \label{E-eps-ham-max}
E(\eps) = \arg\max\limits_{E \in \mathrm{Ham}(\R^{n,n}), \| E \|_F = 1} \Re\, \lambda (  A + \eps E).
\end{equation}
The details of the algorithm are nearly identical to Section~\ref{subsec:ham-A}, except that the sign of the right-hand sides of the differential equations \eqref{ode-E-S-recall}, \eqref{ode-E-ham} and \eqref{ode-E-4} is switched, or in other words, we go backward in time with the same differential equations. 

In the outer iteration we compute, for a given small $\delta>0$, the perturbation size $\eps_\delta$ as the largest $\eps$ with
$ \Re \,\lambda(  A + \eps E(\eps) )=\delta$,  in the same way as in Section~\ref{subsec:outer-it-sqrt}.


\section{Nearest defective matrix}
\label{sec:defective}

The following is a classical matrix nearness problem known as the Wilkinson problem.
\index{Wilkinson problem}
\index{defective eigenvalue}

\medskip
\noindent
{\bf Problem.} {\it Given an $n \times n$ matrix $A$ with $n$ distinct eigenvalues, find a nearest matrix with a defective multiple eigenvalue.}

\medskip\noindent
Here, a multiple eigenvalue is called {\it defective} if 
the Jordan canonical form has a Jordan block of dimension at least 2 with this eigenvalue. Since the left and right eigenvectors of a nonscalar Jordan block are orthogonal to each other (but those of a scalar Jordan block are not), it follows that a matrix is  defective if and only if it has an eigenvalue for which the inner product of the left and right eigenvectors is zero.

We consider the matrix nearness problem with structured perturbations that are in a given structure space $\cS$, which is a complex-linear or real-linear subspace of $\C^{n,n}$ as in Chapter~\ref{chap:struc}. (The unstructured complex and real cases are particular cases with $\cS=\C^{n,n}$ and $\cS=\R^{n,n}$, respectively.) We note that a structured Wilkinson problem need not have a solution. For example, every Hermitian perturbation of a Hermitian matrix is Hermitian and hence its left and right eigenvectors to any eigenvalue are equal.

Our interest is in computing the  nearest matrix w.r.t. the Frobenius norm that differs from the given matrix $A$ by a structured perturbation $\Delta\in \cS$. The distance to defectivity is
\begin{eqnarray}
w^\cS(A) & := & \inf \bigl\{ \| \Delta  \|_F \,:\ \text{$\Delta  \in \cS$ is such that $A + \Delta $ is defective}
\bigr\}.
\label{eq:dist}
\end{eqnarray}
This can be interpreted as the structured distance to singularity of the eigenvalue condition number $\kappa=1/(x^*y)$, where $x$ and $y$ are left and right eigenvectors of unit norm and with positive inner product. This notion of eigenvalue condition number goes back to Wilkinson (\cite{Wil65}). 
\index{distance to singularity!of eigenvalue condition number}


\subsection{Two-level approach}
\index{two-level iteration}
For $0<\eps<w^\cS(A)$ we introduce the functional $\F_\eps(E)$ (for matrices $E\in \C^{n,n}$, which will later be restricted to be structured and of unit Frobenius norm) as follows: For an eigenvalue $\lambda$ of $A+\eps E$, let $x$ and $y$ be corresponding left and right eigenvectors, normalized to unit norm and with positive inner product
$x^*y>0$. \bcl In the present context, a target eigenvalue $\lambda(A+\eps E)$ is an eigenvalue of $A+\eps E$ for which the inner product $x^*y$ of corresponding eigenvectors is minimal among all eigenvalues of $A+\eps E$.
We set
\index{eigenvector optimization}
\begin{equation}\label{F-eps-def}
    \F_\eps(E)= x^*y >0,
\end{equation}
where $x$ and $y$ are eigenvectors to a target eigenvalue $\lambda(A+\eps E)$, normalized to unit norm and with positive inner product.
\ecl
In contrast to previous sections, the functional to be minimized now depends on eigenvectors instead of eigenvalues.

We follow the two-level approach of Chapter~\ref{chap:two-level}:
\begin{itemize}
\item {\bf Inner iteration:\/} Given $\eps>0$, we aim to compute a  matrix $E(\eps) \in\cS $  
of unit Frobenius norm that minimizes $\F_\eps$:
\begin{equation} \label{E-eps-def}
E(\eps) = \arg\min\limits_{E \in \cS , \| E \|_F = 1} \F_\eps(E).
\end{equation}

\item {\bf Outer iteration:\/} For a small threshold $\delta>0$, we compute the smallest positive value $\eps_\delta$ with
\begin{equation} \label{zero-def}
\phi(\eps_\delta)= \delta
\end{equation}
with $\phi(\eps)=  \F_\eps( E(\eps) )=x(\eps)^*y(\eps)$, where $x(\eps)$ and $y(\eps)$ are left and right eigenvectors of $A+\eps E(\eps)$ associated with $\lambda(\eps)$, of unit norm and with positive inner product.
\end{itemize}

\bcltwo
Provided that these computations succeed in computing a global minimum in \eqref{E-eps-def}, we then expect that $\Delta A_\delta= \eps_\delta E(\eps_\delta) \in \cS $ makes $A+\Delta A$ close to a defective matrix and that the limit $\oeps=\lim_{\delta\searrow 0}\eps_\delta$ exists and is equal to the distance $w^\cS(A)$ of the matrix $A$ to the set of defective real matrices. These steps are detailed in the following subsections.
\ecltwo

\subsection{Constrained gradient flow}

As in Chapter~\ref{chap:proto}, we begin by computing the free (complex) gradient of the functional $\F_\eps$.
\index{gradient!free}

\begin{lemma}[Free gradient] \label{lem:gradient-mp-cnv}
Let $E(t)\in \C^{n,n} $, for real $t$ near $t_0$, be a continuously differentiable path of matrices, with the derivative denoted by $\dot E(t)$.
Assume that $\lambda(t)$ is a simple eigenvalue of  $A+\eps E(t)$ depending continuously on $t$,
with corresponding left and right eigenvectors $x(t)$ and $y(t)$, respectively, which are taken to be of unit norm and with positive inner product. Let $\kappa(t)=1/(x(t)^*y(t))$.
Then, $\F_\eps(E(t))=x(t)^* y(t)$ 
is continuously differentiable w.r.t. $t$ and
\begin{equation} \label{eq:deriv-def}
\frac{\kappa(t)}{ \eps  } \,\frac{d}{dt} \F_\eps(E(t)) = \Re \, \bigl\langle  G_\eps(E(t)),  \dot E(t) \bigr\rangle,
\end{equation}
where the rescaled gradient of $\F_\eps$ is a matrix of rank at most 2, given by
\begin{equation} \label{freegrad-def}
G_\eps(E) = xx ^* Z ^* + Z ^*yy ^*.
\end{equation}
Here, $Z$ is the group inverse of $A + \eps E -\lambda I$ (see Section~\ref{sec:eigvec-deriv})
for  the eigenvalue $\lambda$ of $A+\eps E$,  and $x$ and $y$ are the left and right normalized eigenvectors  with positive inner product.
\end{lemma}
\index{group inverse}

\smallskip
\begin{proof}
By (\ref{chap:appendix}.\ref{deigvec}) and using that $Z y=0$, $x ^*Z=0$, we get 
\begin{align*}
\frac{d}{d t} (x ^* y) & = \Re\bigl(\dot x ^* y + x ^*\dot y\bigr) = 
\\ & = 
\eps\,(x ^*y) \,\Re\bigl(x ^*\dot E Z x + y ^* Z \dot E y\bigr) \\ & = 
\eps\,(x ^*y) \, \Re \,\langle xx ^* Z ^* + Z ^*yy ^*, \dot E \rangle = \eps\,(x^*y) \, \Re \,\langle G_\eps(E), \dot E \rangle\;,
\end{align*}
from which \eqref{eq:deriv-def} follows.  
\qed
\end{proof}

\subsubsection*{Structured gradient.} 
\index{gradient!structured}
For a path of {\it structured} matrices $E(t)\in\cS $, also $\dot E(t)$ is in $\cS$, and hence the right-hand side of \eqref{eq:deriv-def} becomes
$\Re \bigl\langle  \Pi^\cS  G_\eps(E(t)),  \dot E(t) \bigr\rangle$,
where $\Pi^\cS$ is the orthogonal projection onto the structure space $\cS$ (see Section~\ref{sec:proto-structured}).
With the structured gradient
\begin{equation}\label{real-grad-def}
G_\eps^\cS(E):=\Pi^\cS G_\eps(E) = \Pi^\cS( xx ^* Z ^* + Z ^*yy ^* ),
\end{equation}
we then have
\begin{equation} \label{eq:deriv-real-def}
\frac{\kappa(t)}{ \eps  } \,\frac{d}{dt} \F_\eps(E(t)) = \Re\,\bigl\langle  G_\eps^\cS(E(t)),  \dot E(t) \bigr\rangle.
\end{equation}
With this structured gradient, we now follow closely the programme of Chapter~\ref{chap:proto} in the structured version of Section~\ref{sec:proto-structured}.

\subsubsection*{Norm-constrained structured gradient flow.}
\index{gradient flow!structure-constrained}
\index{gradient flow!norm-constrained}
We consider the  gradient flow on the manifold of matrices in $\cS $ of Frobenius norm~1,
\begin{equation}\label{ode-E-def}
\dot E = -G_\eps^\cS(E) + \langle G_\eps^\cS(E), E \rangle E.
\end{equation}

\subsubsection*{Monotonicity.} 
Assuming simple eigenvalues along the trajectory,
we again have the monotonicity property of Theorem~\ref{chap:proto}.\ref{thm:monotone} with (\ref{chap:proto}.\ref{eq:pos-real}) and of (\ref{chap:struc}.\ref{eq:pos-S}),
\begin{equation}
\frac{d}{dt} \F_\eps (E(t)) = - \| G_\eps^\cS(E) - \langle G_\eps^\cS(E), E \rangle E \|_F^2 \le  0.
\label{eq:pos-def}
\end{equation}

\subsubsection*{Stationary points.}
\index{stationary point}
Also the characterization of stationary points as given in Theorem~\ref{chap:proto}.\ref{thm:stat} and (\ref{chap:struc}.\ref{stat-S}) extends with the same proof: Let
$E\in\cS $ with $\| E\|_F=1$ be such that the eigenvalue $\lambda(A+\eps E)$ is simple 
and $G_\eps^\cS(E)\ne 0$. Then, 
\begin{equation}\label{stat-def}
\begin{aligned}
&\text{$E$ is a stationary point of the differential equation \eqref{ode-E-S}}
\\[-1mm]
&\text{if and only if $E$ is a real multiple of $G_\eps^\cS(E)$.}
\end{aligned}
\end{equation}
Hence, in this situation an optimizer of \eqref{E-eps-def} is the orthogonal projection onto $\cS$ of a rank-2 matrix.

This raises the question as to whether the structured gradient $G_\eps^\cS(E)=\Pi^\cS( xx ^* Z ^* + Z ^*yy ^* )$ can be the zero matrix. We first note that if the left and right eigenvectors are equal, i.e. $x=y$, then $G_\eps(E)=0$ because we have $Z y=0$ and $x ^*Z=0$. Hence, for a normal matrix $A+\eps E$ the gradient is always zero.

The following two results exclude a vanishing gradient in the unstructured complex and real cases provided that $x\ne y$.
We do not know of a result that guarantees a nonvanishing structured gradient $G_\eps^\cS(E)$ when $\cS$ is a proper subspace of $\C^{n,n}$ or $\R^{n,n}$. The elegant argument of Theorem~\ref{chap:struc}.\ref{thm:nonzero-gradient-S} does not apply here.
\smallskip
\begin{theorem}[Non-vanishing free gradient]
\label{th:ReS-def}
Assume that the complex matrix $B$ has a simple eigenvalue
$\lambda$. Let $x$ and $y$ be left and right eigenvectors of unit norm with positive inner product
associated with $\lambda$ and assume that $x \ne y$. Let $Z$ be the group inverse of $B - \lambda I$ and let
$G=xx ^* Z ^* + Z ^*yy ^*$ as in \eqref{freegrad-def}. Then, $G \ne 0$.
\end{theorem}

\begin{proof}
If $G=0$, then $Zx$ must be a multiple of $y$, and since $y$ is in the null-space of $N=B-\lambda I$, we then have $NZx=0$.
Since the group inverse $Z$ satisfies $NZ=ZN$, this implies $ZNx=0$, and since $NZN=N$, this yields $Nx=0$. Since $\lambda$ is a simple eigenvalue with right eigenvector $y$, $N$ has a 1-dimensional null-space that is spanned by $y$. So $x$ and $y$ are collinear. With the imposed normalization of the left and right eigenvalues this implies $x=y$ contrary to our assumption.
\qed
\end{proof}

For the following result we recall that for the unstructured real case, i.e. $\cS=\R^{n,n}$, we have $\Pi^\cS G = \Re\, G$.
\begin{theorem}[Non-vanishing real gradient]
\label{th:ReS-def-R}
Assume that the real matrix $B$ has a pair of simple complex conjugate eigenvalues
$\lambda$ and $\bar\lambda$. Let $x$ and $y$ be left and right eigenvectors of unit norm with positive inner product
associated with $\lambda$ and assume that $x \ne y$. Let $Z$ be the group inverse of $B - \lambda I$ and let
$G=xx ^* Z ^* + Z ^*yy ^*$ as in \eqref{freegrad-def}. Then, $\Re\,G \ne 0$.
\end{theorem}


\smallskip
\begin{proof}
The proof is done by leading the assumption $\Re\, G= 0$ to a contradiction. So let us assume that $G$ is purely imaginary.
In part (a) of the proof we show that then $\{x,\conj x\}$ and $\{y,\conj y\}$ span the same 2-dimensional invariant subspace. In part (b) we work on this subspace and derive a contradiction.

(a) By definition of the matrix $G$, its range is given by
\[
\mathrm{Ran}(G) = {\rm span}\left\{ x, Z ^* y \right\}.
\]
If $G$ is purely imaginary, then $\mathrm{Ran}(\conj G)= \mathrm{Ran}(G)$ and hence $\conj x\in\mathrm{Ran}(G)$, i.e., we have
$\conj{x}  =  \alpha x + \beta Z ^* y$.
A left premultiplication with $y ^*$  allows us to conclude $\alpha=0$, because
 (i) \ $y ^*\conj{x} = 0$ by the bi-orthogonality of left 
and right eigenvectors corresponding to different eigenvalues (here $\lambda$ and $\clambda$), \ (ii) $Z y = 0$, a property of the group inverse $Z$, and \, (iii)  
$x ^* y \neq 0$ as $\lambda$ is a simple eigenvalue.
This implies $\conj{x} \propto  Z ^*y $. Analogously, we obtain $\bar{y} \propto Z x$. So we have
\begin{equation} \label{eq:d1}
    Z x = \gamma \conj{y}, \qquad Z ^* y = \eta \conj{x},
\end{equation}
with $\gamma \neq 0$ and $\eta \neq 0$.

Since $y$ is in the null-space of $B-\lambda I$ and $B$ is real, it follows that $\conj y $ is a right eigenvector of $B-\lambda I$ to the eigenvalue $\mu = -2\iu \, \Im\,\lambda\ne 0$. For the group inverse $Z$ of $B-\lambda I$, this implies that $\conj y $ is a right eigenvector of $Z$ to the eigenvalue $\nu=1/\mu$. Analogously we find that
$\conj x$ is a left eigenvector of $Z$ to the eigenvalue $\nu$.
So we have
\begin{equation} \label{eq:d2}
  Z \conj{y} = \nu \conj{y} , \qquad Z ^* \conj{x} = \conj\nu \,\conj{x}, 
\end{equation}
with $\nu\ne 0$. Equations \eqref{eq:d1}--\eqref{eq:d2} show that $Z(\gamma^{-1} x - \nu^{-1} \conj y)=0$, and since $y$ spans the null-space of $Z$, we find $y \propto \gamma^{-1} x - \nu^{-1} \conj y$. Analogously we obtain
$x\propto \eta^{-1} y - \conj{\nu}^{-1} \conj x$. We conclude that
\begin{equation}\label{span-x-y}
{\rm span} \left\{ y, \conj{y} \right\} = {\rm span} \left\{ x, \conj{x} \right\}\,.
\end{equation} 

(b) The two-dimensional space $Y:={\rm span} \left\{ y, \conj{y} \right\}$ is thus an invariant subspace of both $B$ and $B^\top$. We choose the real orthonormal basis $(q_1,q_2)$ of $Y$ that is obtained by normalizing the orthogonal vectors $y+\conj y$ and $i(y-\conj y)$. We extend this basis of $Y$  to a real orthonormal basis $Q=(q_1,\dots,q_n)$ of $\C^n$. From \eqref{span-x-y} we infer the block-diagonal structure 
\[
\widetilde{B} = Q^\top B Q =
\left( 
\begin{array}{cc}
B_1 & 0 \\
0 & B_2  
\end{array}
\right) .
\]
The $2\times 2$ matrix 
\begin{equation*}
B_1 = \begin{pmatrix} 
\varrho    & \sigma \\
{}-\tau  & \varrho  
\end{pmatrix}
\end{equation*}
is such that $\varrho = \Re(\lambda)$ and $\sigma > 0, \tau > 0$ with 
$\sigma \tau = \Im(\lambda)^2 > 0$ so that $B_1$ has eigenvalues $\lambda$ 
and $\conj{\lambda}$.

If $\sigma=\tau$ then $B_1$ is normal, which implies that the pair of right 
and left eigenvectors associated with $\lambda$, say $\widetilde{y}, \widetilde{x}$
(scaled to have unit norm and positive inner product)
is such that $\widetilde{x} ^* \widetilde{y} = 1$. 
Since $y = Q \widetilde{y}$ and $x = Q \widetilde{x}$, the orthogonality of $Q$ 
implies $x ^*y=1$, which contradicts the assumption $x^*y<1$.
So we must have $\sigma \ne \tau$.

By the properties of the group inverse we have that
\[
\widetilde{Z} = Q^\top Z Q =
\left( 
\begin{array}{cc}
Z_1 & 0 \\
0 & Z_2  
\end{array}
\right) ,
\]
where $Z_1$ is the group inverse of $B_1-\lambda I$ and 
$Z_2$ is the inverse of the nonsingular matrix $B_2-\lambda I$.
The following formula for the group inverse is verified
by simply checking the three conditions in Definition \ref{chap:appendix}.\ref{def:groupinv}:
\[
Z_1 =
\left(
\begin{array}{rr}
 \frac{\iu}{4 \sqrt{\sigma \tau}} & -\frac{1}{4 \tau} \\
 \frac{1}{4 \sigma} & \frac{\iu}{4 \sqrt{\sigma \tau}}
\end{array}
\right).
\]
It follows that also $Q^\top G Q$ is block diagonal so that we write
\begin{eqnarray*}
\widetilde{G} = Q^\top G Q =
\left( 
\begin{array}{cc}
G_1 & 0 \\
0 & G_2  
\end{array}
\right) 
\quad\ \text{ with }\quad\ 
G_1 = \widetilde{x}_1 \widetilde{x}_1 ^* Z_1 ^* + Z_1 ^*  \widetilde{y}_1 \widetilde{y}_1 ^* ,
\end{eqnarray*}
where $\widetilde{x}_1 \in \C^2$ and $\widetilde{y}_1 \in \C^2$ are the projections onto 
${\rm span}(e_1,e_2)$ (the subspace spanned by the first two vectors of the canonical basis) 
of the eigenvectors of $\widetilde{B}$ associated with $\lambda$, that is 
$\widetilde{y} = Q^\top y$ and $\widetilde{x} = Q^\top x$, 
\begin{eqnarray}
\widetilde{y} & = & \nu_y^{-1} \left( \begin{array}{ccccc}
{} \iu \frac{\sqrt{\sigma}}{\sqrt{\tau}} & 1 & 0 & \ldots & 0 \end{array} \right)^{\top}
\nonumber
\\
\widetilde{x} & = & \nu_x^{-1} \left( \begin{array}{ccccc}
{}-\iu \frac{\sqrt{\tau}}{\sqrt{\sigma}} & 1 & 0 & \ldots & 0 \end{array} \right)^{\top}
\nonumber
\end{eqnarray}
with $\nu_y = \sqrt{\frac{\sigma}{\tau}+1}$ and $\nu_x = \sqrt{\frac{\tau}{\sigma}+1}$  chosen such that 
$\widetilde{x}$ and $ \widetilde{y}$ are of unit norm with positive inner product.
Finally we obtain
\begin{eqnarray}
G_1 & = & \left(
\begin{array}{ll}
 0 & \frac{\tau-\sigma}{2 \sigma \left(\sigma+\tau \right)} \\
 \frac{\tau-\sigma}{2 \tau \left(\sigma+\tau\right)} & 0
\end{array}
\right),
\nonumber
\end{eqnarray}
which is real and cannot vanish due to the fact that $\sigma \neq \tau$.

Recalling that $Q$ is real, if $G$ were purely imaginary then $G_1$ would be purely imaginary as well, which 
would give a contradiction. 
\qed
\end{proof}

\subsection{Inner iteration: Rank-2 dynamics} 
\label{subsec:wilkinson-rank-2-ode}
We proceed in analogy to Section~\ref{subsec:struc-rank-1-ode}.
We project  onto the tangent space $T_Y\cR_2$ at $Y$ of the manifold of complex rank-2 matrices $\cR_2=\cR_2(\C^{n,n})$ and consider the projected differential equation with solutions of rank 2:
\index{rank-2 matrix differential equation}
 \begin{equation}\label{ode-E-S-2}
\dot Y = -P_Y G_\eps(E) + \Re \langle P_Y G_\eps(E), E \rangle Y \quad\text{ with }\ E=\Pi^\cS Y,
\end{equation}
where $P_Y$ is the orthogonal projection onto the tangent space $T_Y\cR_2$. We note that
\begin{equation}\label{ode-E-S-2-Pi}
\dot E = -\Pi^\cS P_Y G_\eps(E) + \Re \langle \Pi^\cS P_Y G_\eps(E), E \rangle E \quad\text{ with }\ E=\Pi^\cS Y,
\end{equation}
which differs from the gradient flow \eqref{ode-E-S} only in that the gradient $G_\eps(E)$ is replaced by the projected gradient $P_Y G_\eps(E)$.
By the same argument as in Section~\ref{subsec:struc-rank-1-ode}, the Frobenius norm 1 of $E(t)=\Pi^\cS Y(t)$ is preserved along solutions $Y(t)\in \cM_2$ of \eqref{ode-E-S-2}. 

As in Theorem~\ref{chap:struc}.\ref{thm:stat-S} we find that the stationary points $Y$ of \eqref{ode-E-S-2} correspond bijectively to the stationary points $E=\Pi^\cS Y$ of the gradient flow \eqref{ode-E-def} provided that $P_E G_\eps^\cS(E) \ne 0$.

The rank-2 differential equation is solved numerically as described in Section~\ref{subsec:low-rank-integrator}.

\subsection{Outer iteration}

We proceed as in Chapter~\ref{chap:two-level} and make an analogous assumption: 
\begin{assumption} \label{ass:E-eps-def}
For $\eps$ close to $\oeps$ and $\eps<\oeps$, 
we assume the following for the optimizer $E(\eps)$ of \eqref{E-eps-def}:
\begin{itemize}
\item The eigenvalue $\lambda(\eps)=\lambda( A+\eps E(\eps))$ is a simple eigenvalue.
\item The map $\eps \mapsto E(\eps)$ is continuously differentiable.
\item The structured gradient $G^\cS(\eps)=G_\eps^\cS(E(\eps))$ is nonzero.
\end{itemize}
\end{assumption}

In the same way as in Theorem~\ref{chap:two-level}.\ref{thm:phi-derivative},  we calculate under this assumption
the derivative of $\phi(\eps)=\F_\eps(E(\eps))=x(\eps)^*y(\eps)$, where $x(\eps)$ and $y(\eps)$ are left and right eigenvectors of $A+\eps E(\eps)$ associated with $\lambda(\eps)$, of unit norm and with positive inner product.
Here we obtain (with $' = d/d\eps$)
\begin{equation}\label{eq:dereps-def}
\phi'(\eps)= - \phi(\eps) \, \| G^\cS(\eps) \|_F < 0.
\end{equation}
Starting from $\eps > 0$ such that $\phi(\eps) > \delta$, we want to compute the smallest root $\eps_\delta>0$ of the equation
$
\phi(\eps) = \delta.
$
It is of interest to study the behaviour of $\phi(\eps)$ as $\eps$ approaches $\oeps = \lim_{\delta\searrow 0}\eps_\delta$, where eigenvalues coalesce to form a Jordan block. We make the following generic assumption.

\begin{assumption} \label{assumpt-epsstar-def}
We assume the following in the limit $\eps\nearrow \oeps$:
\begin{itemize}
    \item The eigenvalue
$\lambda(\eps)$ coalesces with only one other eigenvalue as $\eps\nearrow \oeps$ to form a Jordan block.
 \item The limits  $x_\star= \lim_{\eps\nearrow \oeps} x(\eps)$,
 $y_\star = \lim_{\eps\nearrow \oeps} y(\eps)$, and $E_\star = \lim_{\eps\nearrow \oeps} E(\eps)$ exist.
\end{itemize}
\end{assumption}
We note that if the limit matrix $E_\star$ exists and the matrix $A+\oeps E_\star$ is non-derogatory, i.e., for each distinct eigenvalue there is only one Jordan block, then the existence of the limits $x_\star$ and $y_\star$ of left and right eigenvectors is ensured by a theorem of Conway \& Halmos (\cite{CH80}). On the other hand, if the limits $x_\star$ and $y_\star$ of left and right eigenvectors exist, then also the limit matrix $E_\star$ exists by \eqref{stat-def} and \eqref{real-grad-def}.

\begin{theorem}[Square root asymptotics]\label{thm:sqrt-def} 
Under Assumptions~\ref{ass:E-eps-def} and~\ref{assumpt-epsstar-def} and the non-degeneracy condition that $\gamma\ge 0$ defined in \eqref{gamma-def} below is nonzero, 
we have
$$
\phi(\eps)= \gamma \,\sqrt{\oeps-\eps} \;(1+o(1)) \quad\text{ as }\ \eps \nearrow \oeps.
$$
\end{theorem}

\begin{proof} The result follows if we can show that $\phi(\eps)\phi'(\eps)$ has a finite nonzero limit as $\eps\nearrow \oeps$. By \eqref{eq:dereps-def} we have 
$$
\phi(\eps)\phi'(\eps) = - \phi(\eps)^2 \,\| \Pi^\cS G(\eps) \|_F.
$$
We recall from \eqref{freegrad-def} that 
$$
G(\eps)=x(\eps)x(\eps)^*Z(\eps)^*+ Z(\eps)^* y(\eps)y(\eps)^*,
$$
where $Z(\eps)$ is the group inverse of $N(\eps):= A+ \eps E(\eps)-\lambda(\eps)I$.
By Assumption~\ref{assumpt-epsstar-def}, the rank of the matrix $N(\eps)$ remains equal to $n-1$ also as $\eps\nearrow\oeps$. Therefore, the Moore-Penrose pseudoinverse 
$$
B(\eps):=N(\eps)^\dag
$$ 
has a finite limit $B_\star$ as $\eps\nearrow\oeps$, and by the second part of Assumption~\ref{assumpt-epsstar-def} also
$$
\beta(\eps):= x(\eps)^* B(\eps)y(\eps)
$$
has a finite limit $\beta_\star$. By Theorem~\ref{chap:appendix}.\ref{thm:Ginv}, we have for $\eps\nearrow\oeps$
\begin{align*}
\phi(\eps)^2 Z(\eps) &= (\phi(\eps) I - y(\eps)x(\eps)^*)B(\eps)(\phi(\eps) I - y(\eps)x(\eps)^*) 
\\
&\to y_\star x_\star^*B_\star y_\star x_\star^* = \beta_\star y_\star x_\star^*
\end{align*}
and therefore by \eqref{freegrad-def},
\begin{align*}
\phi(\eps)^2 \,\Pi^\cS G(\eps) &= \Pi^\cS\Bigl(x(\eps)x(\eps)^*\phi(\eps)^2 \,Z(\eps)^*+ \phi(\eps)^2 \,Z(\eps)^* y(\eps)y(\eps)^*\Bigr)
\\
&\to  \Pi^\cS \bigl(2\conj\beta_\star x_\star y_\star^*\bigr),
\end{align*}
so that finally using \eqref{eq:dereps-def},
\begin{equation}\label{gamma-def}
\phi(\eps)\phi'(\eps) = - \phi(\eps)^2 \| \Pi^\cS G(\eps) \|_F \to -\tfrac12 \gamma^2 := - 2 \|\Pi^\cS \bigl(\conj\beta_\star x_\star y_\star^*\bigr)\|_F.
\end{equation}
The stated result then follows in the same way as in part (d) of the proof of Theorem~\ref{thm:sqrt} provided that $\gamma\ne 0$ as is assumed.
\qed
\end{proof}

In view of the expected square root behaviour of Theorem~\ref{thm:sqrt-def}, we use an outer iteration based on a square root model and bisection as described in Section~\ref{subsec:outer-it-sqrt}.
If $\delta$ is not too small, a classical Newton iteration might also be used.

\bcl
\subsection{Initial perturbation}  
\label{subsec:init-defective}
\index{initial perturbation}

The initial perturbation matrix $E(0)$ and perturbation size $\eps_0$ are chosen as in Section \ref{subsec:init} and in (\ref{chap:two-level}.\ref{eq:eps0F}), with the only difference that the gradient $G$ is now given by \eqref{freegrad-def} and \eqref{real-grad-def} for the unstructured and structured cases, respectively.
\ecl

\subsection*{A small illustrative example}
We consider the following $10 \times 10$ matrix 
$A = {\rm Grcar}(10)$ 
from the Eigtool demo,
that is,
\begin{equation} \label{eq:Gr10}
A = 
\left( \begin{array}{rrrrrrrrrr}
   1  &  1  &  1  &  1  &  0  &  0  &  0  &  0   &  0  &  0  \\ 
  -1  &  1  &  1  &  1  &  1  &  0  &  0  &  0   &  0  &  0 \\
   0  & -1  &  1  &  1  &  1  &  1  &  0  &  0   &  0  &  0 \\
   0  &  0  & -1  &  1  &  1  &  1  &  1  &  0   &  0  &  0 \\
   0  &  0  &  0  & \ddots & \ddots & \ddots & \ddots & \ddots & \ddots & \ddots
\end{array}
\right).
\end{equation}
\begin{figure}[ht]
\centering
\vspace{-0.1cm}
\includegraphics[width=0.53\textwidth]{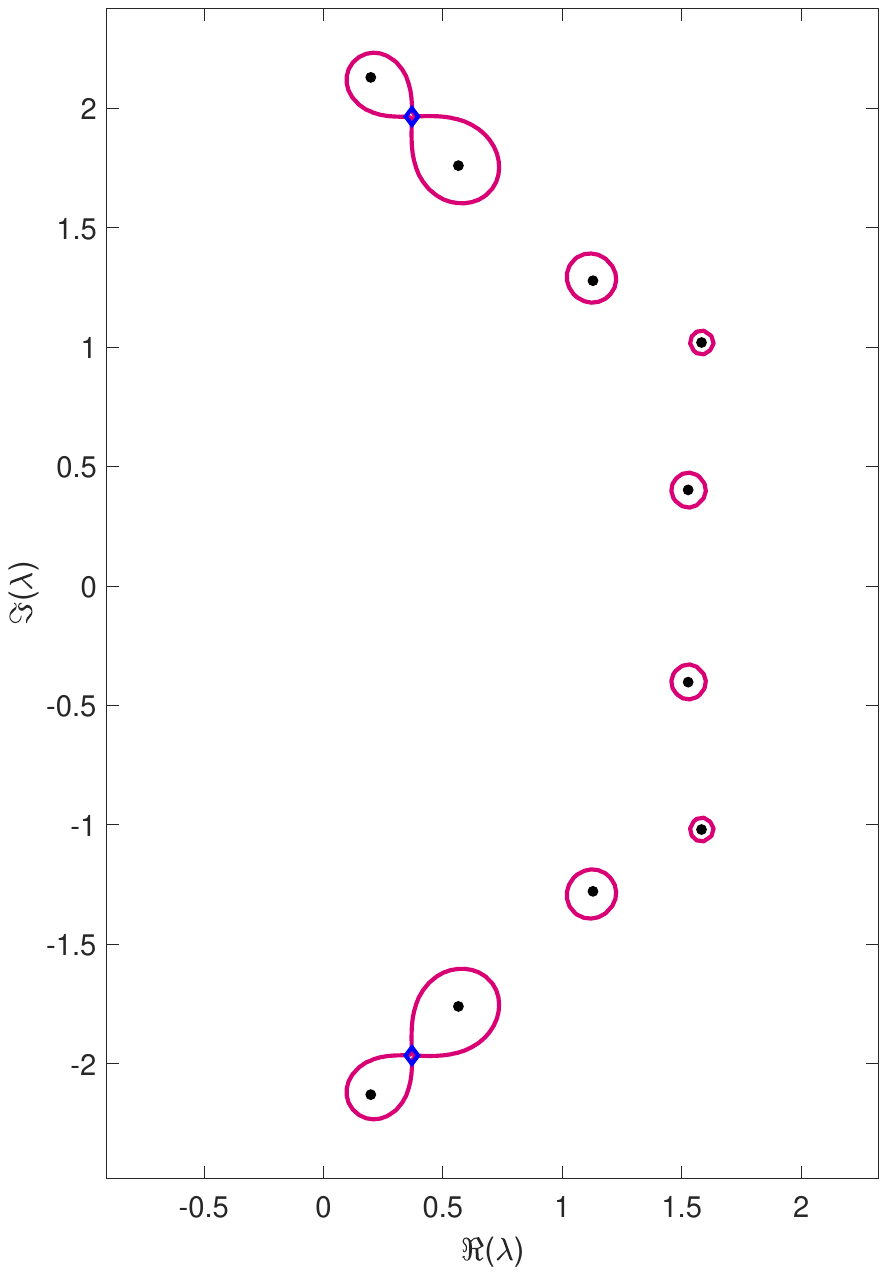} \hskip -10mm
\includegraphics[width=0.53\textwidth]{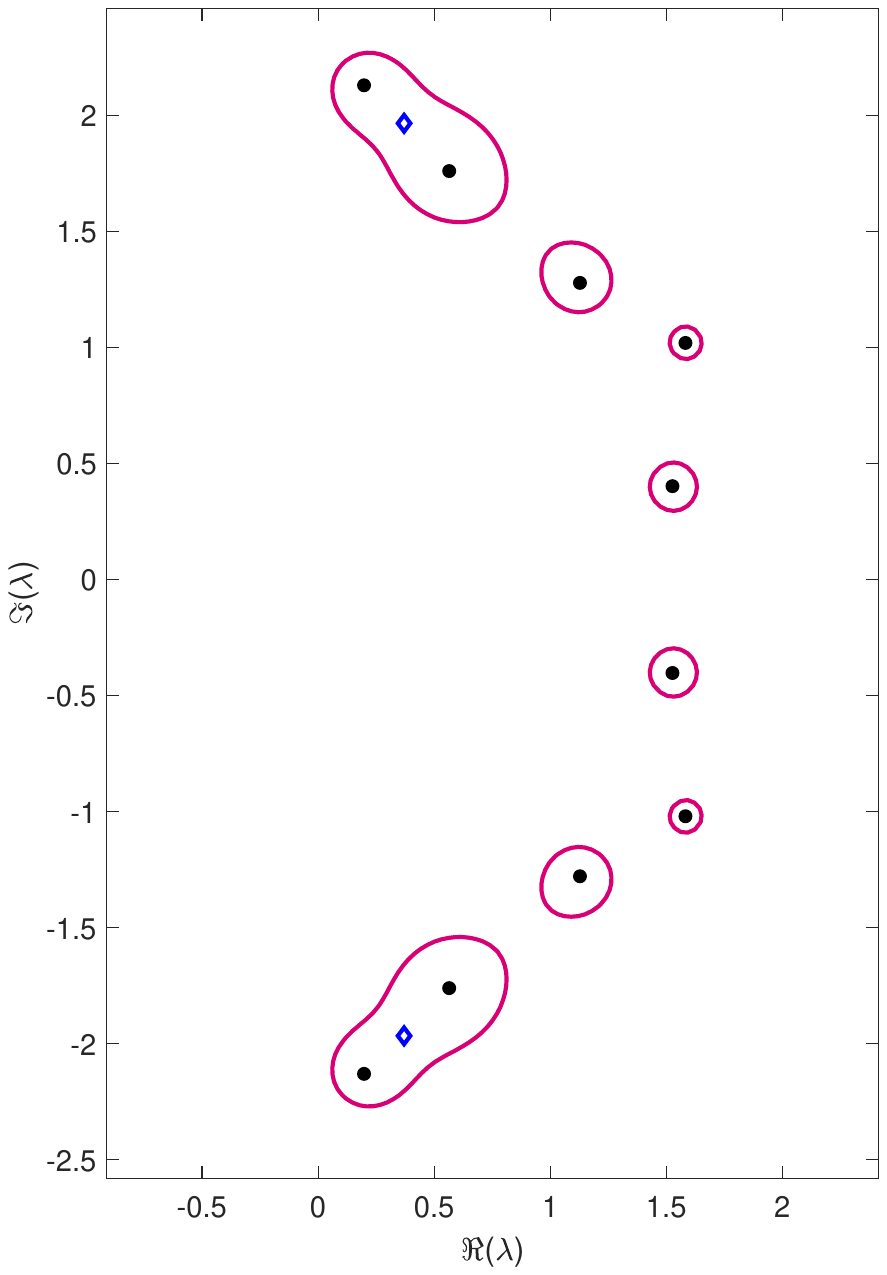}
\vspace{-1cm}
\caption{Left picture: unstructured $\eps$-pseudospectrum of the matrix \eqref{eq:Gr10}
at the smallest value $\eps=\eps_1$ of coalescence 
of a pair of eigenvalues.
Right picture: unstructured $\eps$-pseudospectrum of the matrix \eqref{eq:Gr10}
at the value $\eps=\eps_2$ corresponding to the norm of the smallest real perturbation of $A$ determining the coalescence of a pair of eigenvalues.}
\label{fig:Def1}
\end{figure}
Considering unstructured perturbations to $A$, the smallest value $\eps$ such that
there exists a matrix $\Delta$ of norm $\eps$ such that $A+\Delta$ has a defective pair of eigenvalues is
$\eps_1 = 0.035369524182688$.
The corresponding $\eps$-pseudospectrum with the two pairs of defective eigenvalues (represented by diamond symbols) is plotted in Figure \ref{fig:Def1} (left picture).

\begin{figure}[ht]
\centering
\vspace{-0.1cm}
\includegraphics[width=0.53\textwidth]{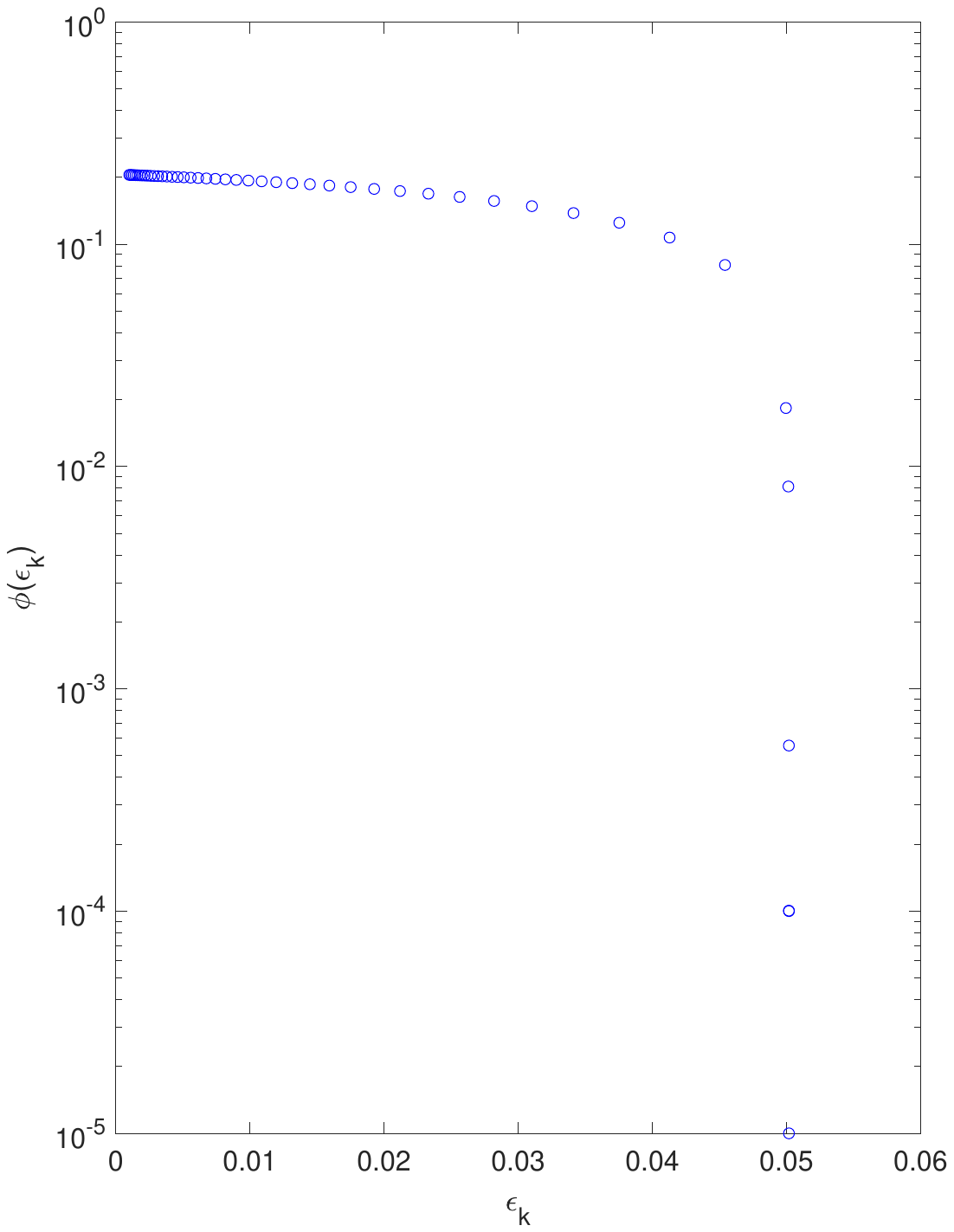} 
\caption{The behaviour of $\phi(\eps)=x(\eps)^*y(\eps)$ 
for the matrix \eqref{eq:Gr10}.}
\label{fig:Def2}
\end{figure}

Considering instead real-structured perturbations to $A$, the smallest value $\eps$ such that
there exists a matrix $\Delta$ of norm $\eps$ such that $A+\Delta$ has a defective pair of eigenvalues is
$\eps_2 = 0.050181307568931$.
The corresponding $\eps$-pseudospectrum with the two pairs of defective eigenvalues (represented by diamond symbols) is plotted in Figure \ref{fig:Def1} (right picture).

It is interesting to observe that the pairs of coalescent eigenvalues determined in the complex case, that is associated to the smallest complex perturbation of $A$, 
are
\[
\lambda_{1,2}^{\C} =
0.370080\ldots \pm 1.96544\ldots \iu.
\]
On the other hand, the two (complex conjugate) pairs of coalescent eigenvalues determined in the real case, associated with the smallest real perturbation of $A$, are
\[
\lambda^{\R}_{1,2} = 
0.370507\ldots \pm 1.965474\ldots \iu.
\]
This indicates that although the values of $\eps_1$ and $\eps_2$ are not very close, the points of coalescence of the eigenvalues are much closer.

\section{Nearest singular matrix pencil}
\label{sec:matrix-pencils}
\index{matrix pencil!singular}
\subsection{Problem setting}
Let $A$ and $B$ be complex $n\times n$ matrices. The matrix pencil $\{ A - \mu B \,:\, \mu \in \C \}$, 
or equivalently the pair $(A,B)$, is called {\it singular} if
\begin{equation}\label{mp-sing}
A-\mu B \ \text{ is singular for all } \mu \in \C.
\end{equation}

This notion is fundamental in the theory of linear differential-algebraic equations $B\dot y(t) = Ay(t) + f(t)$. If and only if the matrix pencil $(A,B)$ is singular, then there exists no initial value $y(0)$ such that the corresponding initial value problem has a unique solution.

A necessary condition for $(A,B)$ to be singular is clearly that both $A$ and $B$ are singular matrices (i.e. non-invertible).
A sufficient condition for $(A,B)$ to be singular is that $A$ and $B$ have a common nonzero vector in their null-spaces. While this is a special case of particular interest, the Kronecker normal form of a matrix pencil (see Gantmacher, \cite{Gan59}) shows that a common null-vector is not a necessary condition for a matrix pencil to be singular.

Given a matrix pencil $(A,B)$ that is not singular, it is of interest to know how far it is from a singular matrix pencil.
In this section we fix $B$, which is assumed to be a singular matrix, and consider structured perturbations  $\Delta A \in \cS$ to $A$, where the structure space $\cS$ is $\C^{n,n}$ or $\R^{n,n}$ or an arbitrary complex- or real-linear subspace, as considered in Section~\ref{sec:proto-structured}. For example, $\cS$ might be a space of real matrices with a given sparsity pattern. We consider the following two structured matrix nearness problems:

\medskip\noindent
{\bf Problem 1.} \ {\it Find $\Delta A\in \cS$ of minimal Frobenius norm such that $(A+\Delta A, B)$ is a singular matrix pencil.}

\medskip\noindent
{\bf Problem 2.} \ {\it Find $\Delta A\in \cS$ of minimal Frobenius norm such that $A+ \Delta A$ and $B$ have a common nonzero vector in their null-spaces.}

\medskip\noindent
If $A$ itself is also in $\cS$, then both problems have a solution, because the trivial choice $\Delta A= -A$ gives us the singular matrix pencil $(0,B)$ and because the set of perturbations $\Delta A$ yielding singular matrix pencils $(A+\Delta A,B)$ is closed. 

We will approach both problems by a two-level method in the spirit of Chapter~\ref{chap:two-level}. We extend the approach of Section~\ref{sec:sing-S}, which does not require computing eigenvalues and eigenvectors, or singular values and singular vectors, in the iterations.
\index{eigenvalue optimization!without eigenvalues}

We note that this approach would equally allow us to treat analogous matrix nearness problems where also $B$ is perturbed, 
but for ease of presentation we have chosen not to do so here. Moreover, in the applications of interest to dynamical systems on networks, $B$ is typically an adjacency matrix or related fixed matrix depending only on the network topology and hence is not subject to perturbations.

\subsection{Distance to structured singular matrix pencils}
Given $n+1$ distinct complex numbers $\mu_0, \dots,\mu_n$,  the fundamental theorem of algebra together with
the fact that a matrix is singular if and only if its determinant vanishes, shows that a matrix pencil $(A,B)$ is singular if and only if 
\begin{equation}\label{mp-sing-d}
\text{the $n+1$ matrices }\ A-\mu_k B \ \text{ are singular for } k=0,\dots,n.
\end{equation}
Our numerical approach for Problem 1 is based on this criterion. The choice of the numbers $\mu_0, \dots,\mu_n$ did not appear to be a critical issue in our numerical experiments. We had good experience with the choice $\mu_k = r e^{2\pi k\iu/(n+1)}$ of modulus $r=\| A \|_F / \| B \|_F$.

We proceed in analogy to Section~\ref{subsec:mvp-sing-S}, working with matrix--vector products instead of eigenvalues or singular values.
For $\eps>0$ we introduce the functional $\F_\eps$ (of matrices $E\in\cS$ of unit Frobenius norm and vectors $u_k,v_k\in\C^n$ of unit Euclidean norm for $k=0,\dots,n$) in the following way: with 
$M_k=A-\mu_k B$,
\begin{equation}\label{F-eps-mp}
 \F_\eps(E,u_0,\dots,u_n,v_0,\dots,v_n) = 
 \sfrac12 \sum_{k=0}^n \Bigl(\| u_k^*(M_k+\eps E)\|^2 + \|(M_k+\eps E)v_k\|^2\Bigr).
 \end{equation}
We follow the two-level approach of Chapter~\ref{chap:two-level}:
\begin{itemize}
\item {\bf Inner iteration:\/} Given $\eps>0$, we aim to compute a minimum of $\F_\eps$. We denote the  minimum as $(E(\eps),u_0(\eps),\dots,u_n(\eps),v_0(\eps),\dots,v_n(\eps))$.

\item {\bf Outer iteration:\/} We compute the smallest positive value $\oeps$ with
\begin{equation} \label{zero-mp}
\phi(\oeps)= 0,
\end{equation}
where $\phi(\eps)$ is the minimal value of $\F_\eps$.
\end{itemize}

Provided that these computations succeed, we then have that $\Delta A_\star= \oeps E(\oeps) \in \cS$ is a solution to Problem 1 above, and $\oeps$ is the distance of the matrix pencil $(A,B)$ to the set of structured singular matrix pencils of the form $(A+\Delta A,B)$ with $\Delta A \in \cS$.

\subsection{Constrained gradient flow for the inner iteration}\label{subsec:gradient-flow-mp}

We closely follow the programme of Section~\ref{subsec:mvp-sing-S}.


\subsubsection*{Structured gradient.}
\index{gradient!structured}
Consider a smooth path of structured matrices $E(t)\in\cS$ and vectors $u_k(t),v_k(t)\in\C^n$ $(k=0,\dots,n+1)$. Since then also $\dot E(t)\in\cS$, we have, similarly to (\ref{chap:struc}.\ref{grad-full-sing}),
\begin{align}
   &\frac d{dt}\,  \F_\eps(E(t),u_0(t),\dots,u_n(t),v_0(t),\dots,v_n(t))  
    \label{grad-full-sing-mp}
   =   \Re \langle G, \eps\dot E \rangle 
      \\ &+\ \sum_{k=0}^n\Bigl(\Re \langle (M_k+\eps E)(M_k+\eps E)^{*}u_k,\dot u_k \rangle +
    \,\Re \langle (M_k+\eps E)^{*}(M_k+\eps E)v_k,\dot v_k \rangle \!\Bigr)
    \nonumber
    \end{align}
with the structured gradient matrix
$$
G= G_\eps^\cS(E,u_0,\dots,u_n,v_0,\dots,v_n)=
\Pi^\cS\sum_{k=0}^n\bigl(u_ku_k^{*}(M_k+\eps E)+(M_k+\eps E)v_kv_k^{*}\bigr),
$$
where $\Pi^\cS$ is again the orthogonal projection onto $\cS$; see Section~\ref{sec:proto-structured}.
\subsubsection*{Norm- and structure-constrained gradient flow.} 
\index{gradient flow!norm-constrained}
\index{gradient flow!structure-constrained}

Taking the norm constraints into account, we arrive at the following system of differential equations; cf.\,(\ref{chap:struc}.\ref{Euv-ode-S}):
\begin{align}
\nonumber
    \eps \dot E &= -G +  \Re\,\langle G, E \rangle \, E
    \\
    \label{Euv-ode-S-mp}
    \alpha_k \, \dot u_k &= -(M_k+\eps E)(M_k+\eps E)^{*} u_k + \|(M_k+\eps E)^{*} u_k\|^2\, u_k
    \\
    \nonumber
    \beta_k \,\dot v_k &= -(M_k+\eps E)^{*}(M_k+\eps E) v_k + \|(M_k+\eps E) v_k\|^2\, v_k,
\end{align}
with scaling factors $\alpha_k(t)>0$ and $\beta_k(t)>0$, which we propose to choose as
$\alpha_k = \|(M_k+\eps E)^{*} u_k\|$ and $\beta_k=\|(M_k+\eps E) v_k\|$ to obtain similar rates of change in $\eps E$, $u_k$, $v_k$ near a stationary point.

By construction, the Frobenius norm 1 of $E(t)$ and the Euclidean norm 1 of $u_k(t)$ and $v_k(t)$ are preserved, and the functional $\F_{\eps}$ decays monotonically along solutions.

\subsubsection*{Stationary points.}
\index{stationary point}
At stationary points $(E,u_0,\ldots,u_n,v_1,\ldots,v_n)$ we find that $u_k$ and $v_k$ are left and right singular vectors to $M_k+\eps E$, respectively, after rescaling $v_k \to \e^{\iu\theta} v_k$ such that $\sigma_k=u_k^{*}(M_k+\eps E)v_k$ is real and positive. 
When we arrive at the stationary point by the gradient flow \eqref{Euv-ode-S-mp}, we can expect that the singular value $\sigma_k$ is the smallest singular value of $M_k$, since other stationary points are not stable.

In the stationary point we have $\sigma_k u_k=(M_k+\eps E)v_k$ and
$\sigma_k v_k = (M_k+\eps E)^{*}u_k$ and hence the functional and the gradient are
\begin{equation}\label{F-G-stat-mp}
\F_\eps= \sum_{k=0}^n \sigma_k^2 \quad\text{ and }\quad
G=\Pi^\cS\sum_{k=0}^n 2\sigma_k u_kv_k^{*}.
\end{equation}
If $G\ne 0$ in the stationary point, it follows that $E$ is a real multiple of $G$. However, in contrast to Section~\ref{subsec:gradient-flow-S}, this now does not imply that non-degenerate optimizers are projections of rank-1 matrices onto the structure space $\cS$, because $G$ is no longer a projected rank-1 matrix. So we cannot work with rank-1 matrices here. In the inner iteration we therefore follow the norm- and structure-constrained gradient flow \eqref{Euv-ode-S-mp} into a stationary point.

 \subsection{Outer iteration} \label{subsec:outer-it-mp}
We proceed as in Section~\ref{subsec:Newton--bisection}. 
Under a differentiability assumption and under the assumption that the gradient  
$G(\eps)=G_\eps^\cS$ is nonzero in the minimum 
of $\F_\eps$,
we find in the same way as in Theorem~\ref{chap:two-level}.\ref{thm:phi-derivative}
a simple
expression for the derivative of $\phi(\eps)=\min \F_\eps$, which here again becomes (with $' = d/d\eps$)
\begin{equation}\label{eq:dereps-mp}
\phi'(\eps)= - \| G(\eps) \|_F < 0.
\end{equation}
This expression can be used in a Newton--bisection method.
\index{Newton--bisection method}
In view of \eqref{F-G-stat-mp},
$\phi(\eps)$ can be expected to behave asymptotically for $\eps\nearrow\oeps$ as
$$
\phi(\eps) \approx c \,(\oeps-\eps)^2.
$$
The unknown quantities $c$ and $\oeps$ can be estimated using $\phi(\eps)$ and $\phi'(\eps)$, which is known from  \eqref{eq:dereps-mp}
for $\eps < \oeps$.
This gives
\begin{eqnarray*}
&& \oeps \approx \eps - 2 \frac{\phi(\eps)}{\phi'(\eps)},
\qquad
c \approx  \frac{\phi'(\eps)^2}{4 \phi(\eps)}.
\end{eqnarray*}
For $\eps = \eps_k < \oeps$, we thus obtain the Newton-type iteration
\begin{equation}
\eps_{k+1} = \eps_k - 2 \frac{\phi(\eps_k)}{\phi'(\eps_k)},
\label{eq:Newton}
\end{equation}
which yields a locally quadratically convergent iteration from the left (if instead $\eps_k > \oeps$ occurs, then we should use bisection, which
would give a linear reduction of the error from the right).

\subsection{Distance to singular pencils with common null-vectors}
We aim to find a perturbation $\Delta A \in \cS$ of minimal Frobenius norm such that there exists a nonzero vector $v$ that is in the null-spaces of both $A+\Delta A$ and $B$. In the unstructured real or complex case ($\cS=\C^{n,n}$ or $\R^{n,n}$),  this problem is solved directly by computing the singular value decompositions $B=U_B\Sigma_B V_B^*$ and $U_B^*A V_B = U \Sigma V^*$ and taking the rank-1 perturbation $\Delta = -\sigma_n U_B u_n v_n^* V_B^*$.

To deal with the structured problem we introduce, for $\eps>0$, a functional $\F_\eps$ (of matrices $E\in\cS$ of unit Frobenius norm and a vector $v$ of unit Euclidean norm) as follows: With $\alpha= \| A \|_2$ and $\beta= \| B \|_2$, we set
\begin{equation}\label{F-eps-mp-cnv}
\F_\eps(E,v) = \sfrac1{2\alpha^2} \| (A + \eps E)v\|^2  + \sfrac1{2\beta^2}  \| Bv \|_2^2.
\end{equation}
With this functional we again follow the two-level approach of Chapter~\ref{chap:two-level}. 
Provided that these computations succeed, we then have that $\Delta A_\star= \oeps E(\oeps) \in \cS$ is a solution to Problem 2 above, and $\oeps$ is the distance of the matrix pencil $(A,B)$ to the set of matrix pencils of the form $(A+\Delta A,B)$ with $\Delta A \in \cS$ for which $A+\Delta A$ and $B$ have a common nonzero null-vector.

The structured gradient is obtained by a straightforward calculation:
For a path $E(t)\in \cS$ and and $v(t)\in\C^n$, we have
$$
 \frac{d}{dt} \F_\eps(E(t),v(t)) = \Re \langle G, \eps \dot E \rangle  + \Re \langle g,\dot v \rangle,
$$
with the projected rank-1 matrix
$$
G=G_\eps^\cS(E,v) = \Pi^\cS \Bigl( \sfrac1{\alpha^2}\, (A+\eps E)^*(A+\eps E)vv^* \Bigr)
$$
and the vector
$$
g=g_\eps(E,v) = \Bigl( \sfrac1{\alpha^2}  (A+\eps E)^*(A+\eps E) + \sfrac1{\beta^2} B^*B \Bigr)v.
$$
With this structured gradient, the full programme of Section~\ref{sec:proto-structured} carries over to the present situation, including the rank-1 matrix differential equation (\ref{chap:struc}.\ref{ode-E-S-1}) and its discretization, which were not available for the (full-rank) gradient of Section~\ref{subsec:gradient-flow-mp}.
\index{rank-1 matrix differential equation}

For the outer iteration, we again use a Newton--bisection method in the same way as in Section~\ref{subsec:outer-it-mp}. Here we have again
$\phi'(\eps)=-\| G(\eps) \|_F$ with $G(\eps)=G_\eps^\cS(E(\eps),v(\eps))$.  

This approach extends directly to the case when also $B$ can be perturbed. The only reason why we have not done so here, is that the notation would be more complicated, requiring two perturbations $E_A\in\cS_A$ and $E_B\in\cS_B$ and
the squared perturbation size $\eps^2=\eps_A^2+\eps_B^2$.



\section{Stability radii for delay differential equations}
\label{sec:delay}
\index{delay differential equation}
\index{stability radius}

\subsection{A simple example}

Consider the scalar delay differential equation
\begin{equation} \label{eq:abdde}
\dot{x}(t) = a x(t) + b x(t-1)
\end{equation}
with $a,b \in \R$.
Looking for solutions $x(t) = c \e^{\lambda t}$ 
gives the characteristic equation \ 
\begin{equation} \label{eq:charab}
\lambda - a - b\,\e^{-\lambda} = 0.
\end{equation}
\begin{figure}
\vspace{-0.5cm}
\centerline{
\includegraphics[width=7.3cm]{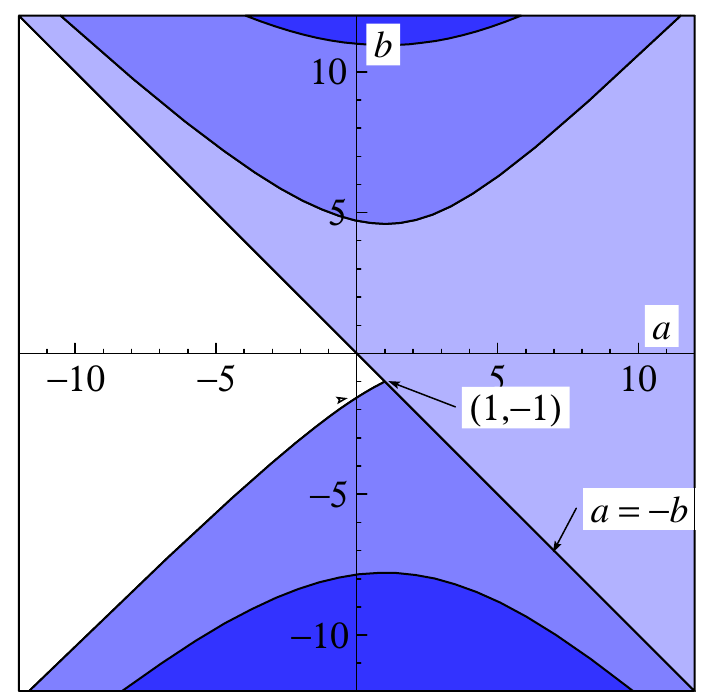} 
}
\label{fig:stabr}
\caption{Region of the $(a,b)$-real plane such that the solution of \eqref{eq:abdde} is asymptotically stable (white set in the right illustration).} 
\end{figure}

A main feature of this problem is that the entire function $\lambda - a - b\,\e^{-\lambda}$
has infinitely many roots and stability can be proved if and only if all the roots lie on the complex left-plane or - in other words - the rightmost root has negative real part. This is what happens when
$a=0.5, b=-1$ as shown in Figure \ref{fig:stabr}
(left illustration).


The stability of a linear system of delay equations 
\begin{eqnarray} \label{delayeig} 
\dot{x}(t) & = & A_1 x(t) + A_2 x\left(t - \tau \right), \qquad t > 0
\\
x(t) & = & g(t), \qquad t \in \left[-\tau,0\right]
\nonumber
\end{eqnarray}
with $A_1, A_2$ given $n \times n$ matrices and $\tau > 0$
a constant delay, can be analyzed in a similar way.

Inserting solutions of the form $x(t) = v\,\e^{\lambda t}$ (with nonzero $v \in \C^n$)
we get 
\begin{equation} \label{eq:charsys}
\left(\lambda \Id - A_1 - A_2 \e^{-\lambda\tau}\right) v=0 ,
\end{equation}
which is a nonlinear eigenvalue problem.
\index{nonlinear eigenvalue problem}

It is known that if the infinitely many eigenvalues strictly lie within the complex left half-plane, then  every solution of \eqref{delayeig} 
is asymptotically stable, independently of the initial data.

The main difficulty here is that in general the matrices $A_1$ and $A_2$ do not commute, so that they cannot be transformed simultaneously to diagonal form, in contrast to the case of linear ordinary differential equations. 
Suitable algorithms for the numerical computation of characteristic roots of linear delay systems are available; see e.g.  Engelborghs and Roose, 2002, Breda, Maset and Vermiglio, 2005, Jarlebring, Meerbergen and Michiels, 2010 (based on Krylov solvers).

%

\subsection{Nonlinear eigenvalue problem and stability radius}
\index{nonlinear eigenvalue problem}

\subsubsection*{Nonlinear eigenvalue problem.}
We consider the problem of determining the eigenvalues $\lambda\in \C$ and corresponding eigenvectors $v\in \C^n\setminus\{0\}$ that satisfy
\begin{equation} \label{nonlin-eig}
\left(\sum_{i=0}^m  f_i(\lambda) A_i\right) v= 0, 
\end{equation}
where $A_0, \ldots,A_m$ are given $n \times n$ matrices and the 
functions $f_0,\ldots, f_m$ are entire complex functions with
\begin{equation*}
f_i(\overline \lambda)=\overline{f_i(\lambda)},\qquad 0\leq i \leq m.
\end{equation*}
We again denote the spectrum by $\Lambda$, i.e.\
\begin{equation}
\Lambda:=\left\{\lambda\in\C:\  \sum_{i=0}^m f_i(\lambda)A_i \  \text{ is a singular matrix}\right\}
\end{equation}
\subsubsection*{Perturbed nonlinear eigenvalue problem.}
We are interested in the effect of bounded perturbations $\Delta A_i$ of $A_i$, which are allowed to be complex.
We consider the  perturbed eigenvalues $\widetilde \lambda$ with
\begin{equation}\label{pert-init}
\left(\sum_{i=0}^m  f_i(\widetilde\lambda)(A_i+\Delta A_i)\right) \widetilde v=0
\end{equation}
for some $\widetilde v\ne 0$. We let
\[
\Delta =\left(  \Delta A_0, \dots , \Delta A_m \right).
\]
%
Introducing weights $w_i > 0$ for $i=0,\ldots, m$, we make use of 
the weighted norm
\begin{equation*}
\|\Delta\|_w := \sqrt{ \sum\limits_{i=0}^{m} w_i^2  \|  \Delta A_i \|_{F}^2}
\end{equation*}
We consider perturbations bounded as
\begin{equation*}
\|\Delta\|_w \le \eps.
\end{equation*}
Taking $w_i=+\infty$ means that the  matrix $A_i$ is not perturbed.
\subsubsection*{Pseudospectra and stability radius.}
\bng
Generalizing the definition given in Chapter \ref{chap:pseudo}, 
the \emph{$\eps$-pseudospectrum} here is the complex set
\eng
\begin{equation*} \label{defps}
\Lambda_{\eps} = \bigcup\limits_{\|\Delta\|_w \le \eps}
\left\{ \lambda\in \C:\ \sum_{i=0}^m  f_i(\lambda)(A_i + \Delta A_i)\  \text{ is a singular matrix} \right\},
\end{equation*}
that is, the set of eigenvalues associated with all perturbed problems having $\|\Delta\|_w \le \eps$.
\begin{figure}[ht]
\vspace{-4cm}
\begin{center}
\includegraphics[width=10.3cm]{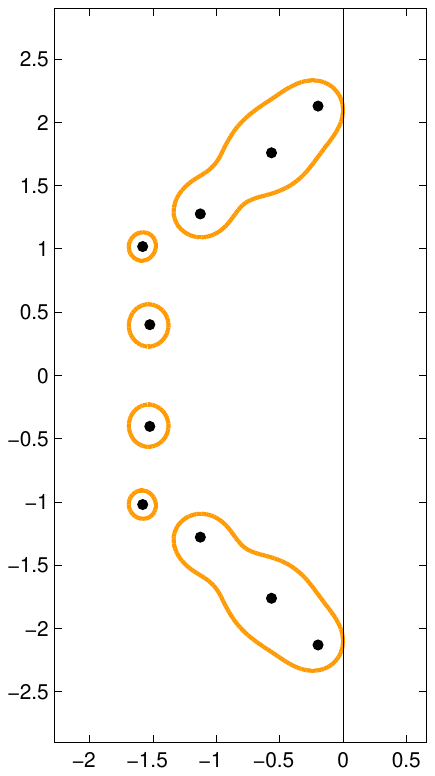}
\vspace{-3cm}
\caption{Rightmost roots of an example of system \eqref{eq:charsys} and $\eps$-pseudospectrum tangential 
to the imaginary axis (showing $\eps=\oeps$).}
\end{center}
\end{figure}

\noindent
The \emph{$\eps$-pseudospectral abscissa} is defined as 
\begin{equation}\label{defpsa}
\alpha_{\eps}:=\sup \left\{\Re\,\lambda\,:\, \lambda \in\Lambda_{\eps}\right\},
\end{equation}
where in this case -- because of the infinitely many eigenvalues -- the $\max$ is replaced by the $\sup$.

Asymptotic stability is associated with the requirement that the spectrum be located in the open left half-plane and bounded away from the imaginary axis. 
The \emph{stability radius} (or \emph{distance to instability}) of an asymptotically  stable system 
is defined as
\index{stability radius}
\begin{equation*}
\oeps:=\inf \left\{\eps>0:\ \alpha_{\eps}\geq 0\right\}.
\end{equation*}

 \medskip\noindent
 {\bf Problem.} {\it Compute the stability radius $\oeps$ for the nonlinear eigenvalue problem~\eqref{pert-init}.}

\subsubsection*{Assumptions.}
\medskip\noindent
For nonlinear eigenvalue problems the pseudospectral abscissa may be equal to infinity (as in differential-algebraic equations), 
\bng
and, even if is finite, a globally rightmost point of the pseudospectrum may not exist.
\eng
To exclude such cases, we assume the following:
\begin{itemize}
\item[\ (i)  ] \ \quad For every $\eps > 0$,\ \ $\alpha_\eps < + \infty$.

\item[(ii) \ ] \ \quad There exists $r>0$ such that the set $\Lambda_{\oeps}\cap\{\lambda\in\C:\ \Re\,\lambda\geq - r\}$ is bounded.
\end{itemize}

As an example which does not fulfill  assumption (ii), consider the following neutral equation
\begin{equation} \label{eq:neut}
\dot x(t) = \dot x(t-\tau) -2x(t)-x(t-\tau).
\end{equation}
\begin{figure}[ht]
\begin{center}
\centerline{
\includegraphics[width=9.3cm]{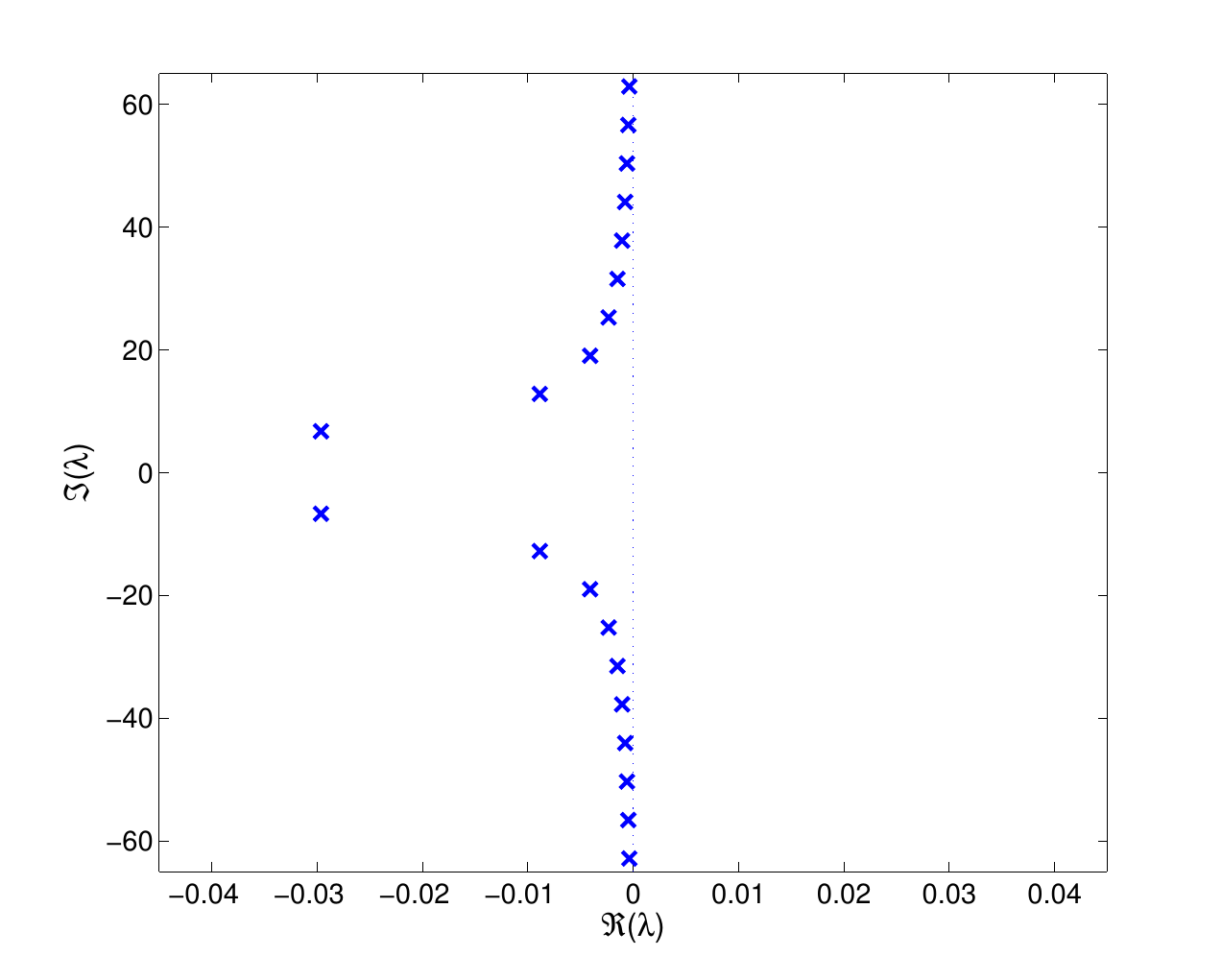}}
\caption{Rightmost roots of the neutral equation.}
\end{center}
\end{figure}
Its characteristic equation is given by
\begin{equation*}\label{quasipol}
\lambda\left(1-\e^{-\lambda}\right)+2+ \e^{-\lambda}=0,
\end{equation*}
which has a sequence of roots approaching the imaginary axis with increasing imaginary part.
Here the spectral abscissa is equal to zero, yet there is no characteristic root with zero real part.

\subsection{Two-level approach}
\index{two-level iteration}
We use a two-level approach as in previous sections to compute the stability radius $\oeps$.
We assume \eqref{nonlin-eig} fulfils assumptions (i) and (ii). Moreover we also assume that
$w_i < \infty$ for $i=0,\ldots,m$ (otherwise the matrices associated to infinite weights
are not perturbed).
We write perturbations as
\begin{equation*}
\Delta A_i = \eps E_i,  \qquad i=0,\dots,m
\end{equation*}
with 
\begin{equation*}
\| (E_0,\dots,E_m) \|_w^2 = \sum_{i=0}^m w_i^2 \| E_i \|_F^2 =1.
\end{equation*}
We denote by $\cS_w^1$ the unit ball for the weighted norm $\| \cdot \|_w$. We look for
\begin{equation*}
\alpha_\eps=\max\limits_{(E_0,  \ldots ,E_m) \in \cS_w^1} 
\biggl\{ \Re \, \lambda: \ \sum_{i=0}^m f_i(\lambda)  (A_i + \eps E_i)   \ \text{ is a singular matrix}\biggr\}.
\end{equation*}
We introduce the functional 
\begin{equation}\label{F-eps-del}
\F_\eps \left(E_0,E_1,\ldots,E_m \right) = 
-\Re\,\lambda,
\end{equation}
where $\lambda$ is the rightmost
eigenvalue of the nonlinear eigenvalue problem \eqref{pert-init} with $\Delta A_i=\eps E_i$. 
We again follow the two-level approach of Chapter~\ref{chap:two-level}:
\begin{itemize}
\item {\bf Inner iteration:\/} Given $\eps>0$, we aim to compute  $(m+1)$-tuples of complex $n\times n$ matrices $\bigl(E_0(\eps),\dots,E_m(\eps)\bigr) \in \cS_w^1$
that minimize $\F_\eps$:
\begin{equation} \label{E-eps-del}
\bigl(E_0(\eps),\dots,E_m(\eps)\bigr) = \arg\min\limits_{(E_0, E_1, \ \ldots ,E_m) \in \cS_w^1} \F_\eps \left(E_0,E_1,\ldots,E_m \right).
\end{equation}

\item {\bf Outer iteration:\/} We compute the smallest positive value $\oeps$ with
\begin{equation} \label{zero-del}
\phi(\oeps)= 0,
\end{equation}
where $\phi(\eps)=  \F_\eps( E_0(\eps),\dots,E_m(\eps) )=-\alpha_\eps$. Then, $\oeps$ is the stability radius.
\end{itemize}

 \subsection{Inner iteration: Rank-1 dynamics}
The programme of Chapter~\ref{chap:proto} extends to  the current situation
 as follows.

\subsubsection*{Free gradient.}
We consider a path $\bigl(E_0(t),\dots,E_m(t)\bigr)\in\cS_w^1$ for $t$ in an open interval around $t_0$ and we assume that the path of eigenvalues $\lambda(t)$ with
$$\det \bigl(f_0(\lambda(t))  (A_0 + \eps E_0(t))  
+ \ldots  + 
f_m(\lambda(t)) (A_m + \eps E_m(t)) \bigr) = 0
$$ 
consists of simple eigenvalues. As in
Lemma~\ref{chap:proto}.\ref{lem:gradient}, we find that
\begin{equation}\label{eq:deriv-S-del}
\frac1{ \eps \kappa(t)} \,\frac{d}{dt} \F_\eps(E_0(t),\dots, E_m(t)) = 
\sum_{i=0}^m w_i^2 \,\Re\bigl\langle  G_i(E_0(t),\dots,E_m(t)),  \dot E_i(t) \bigr\rangle
\end{equation}
with the rescaled partial gradient
\begin{equation}\label{gradient-del}
G_i\left(E_0,E_1,\ldots,E_m \right) = - w_i^{-2} f_i(\clambda) x y^*, \qquad i=0,\dots,m,
\end{equation}
where $x$ and $y$ are left and right eigenvectors corresponding to $\lambda$, chosen of unit norm and with positive inner product, and where $\kappa=1/(x^*y)$.

\subsubsection*{Norm-constrained gradient flow.}
We consider the  gradient flow on the unit sphere $\cS_w^1$,
\begin{equation}\label{ode-E-del}
\dot E_i = -G_i(E_0,\dots,E_m) + \mu E_i, \qquad i=0,\dots,m,
\end{equation}
with the Lagrange multiplier 
$$
\mu =  \sum_{i=0}^m w_i^2 \, \Re \langle G_i(E_0,\dots,E_m) , E_i \rangle.
$$
This ensures that for $(E_0,\dots,E_m)\in \cS_w^1$,
$$
\frac12\,\frac{d}{dt} \| (E_0,\dots, E_m) \|_w^2 = \sum_{i=1} w_i^2 \,\Re\langle \dot E_i, E_i \rangle = 0,
$$
so that the weighted norm 1 is preserved along the flow.

\subsubsection*{Monotonicity.} 
Assuming simple eigenvalues along the trajectory,
we again obtain the monotonicity property as in Theorem~\ref{chap:proto}.\ref{thm:monotone},
\begin{equation}
\frac{d}{dt} \F_\eps (E_0(t),\dots,E_m(t)) = -\eps\kappa\sum_{i=0}^m w_i^2 \| G_i(E_0,\dots,E_m) - \mu E_i \|_F^2
\le  0.
\label{eq:pos-del}
\end{equation}

\subsubsection*{Stationary points.}
\index{stationary point}
Also the characterization of stationary points as given in Theorem~\ref{chap:proto}.\ref{thm:stat} extends with the same proof: Let $\eps>0$ be fixed and let
$(E_0,\dots,E_m)\in\cS_w^1$ be such that the correponding rightmost eigenvalue $\lambda$ is simple 
and $G_i(E_0,\dots,E_m)\ne 0$ for $i=0,\dots,m$. Then, 
\begin{equation}\label{stat-S-mp}
\begin{aligned}
&\text{$(E_0,\dots,E_m)$ is a stationary point of the differential equation \eqref{ode-E-del} if and only if }
\\[-1mm]
&\text{$E_i$ is the same real multiple of $G_i(E_0,\dots,E_m)$ for $i=0,\dots,m$.}
\end{aligned}
\end{equation}
In particular, this shows that at non-degenerate stationary points, $E_i$ is of rank 1 for $i=0,\dots,m$.

\subsubsection*{Rank-$1$ constrained gradient flow.}
\index{gradient flow!rank-1 constrained}
In view of the rank-1 property of optimizers, we proceed further as in Chapter~\ref{chap:proto} and constrain the gradient flow \eqref{ode-E-del} to the rank-$1$ manifold $\cM_1$ by projecting the right-hand side onto the tangent space. This yields the norm- and rank-1 constrained gradient flow
\begin{equation} 
\dot{E}_i  =  - P_{E_i} G_i + \mu E_i, \qquad i=0,\dots,m
\label{ode-E-1-delay}
\end{equation}
with $E_i=\sigma_i u_i v_i^*$ of rank $1$ with $\|u_i\|=\|v_i\|=1$ for all $i$, with
$P_{E_i}$ the orthogonal projection onto the tangent space at $E_i \in \cM_1$, and the Lagrange multiplier $\mu$ as in \eqref{ode-E-del}.

Similarly to the prototype case considered in Chapter \ref{chap:proto}, it can be shown that
the rank-1 projected differential equations preserve norm, monotonicity and stationary points.
The differential equations are integrated numerically into a stationary point  in the way described in Chapter~\ref{chap:proto}, working with the vectors that define the rank-1 matrices $E_i$ and advancing them in time by a splitting method.
This again yields a significant reduction in CPU time and memory.

\subsection{Outer iteration}
\index{Newton--bisection method}
We use again a Newton--bisection method to compute the stability radius, which is the smallest zero $\oeps$ of $\phi(\eps)=\F_\eps(E_0(\eps),\dots,E_m(\eps))$. Under assumptions analogous to Chapter~\ref{chap:two-level}, we find in the same way that $\phi$ is continuously differentiable and the derivative  here becomes, with $' = d/d\eps$
and $G_i(\eps)=G_i^\eps(E_0(\eps),\dots,E_m(\eps))$ and $\kappa(\eps)=1/(x(\eps)^*y(\eps))>0$,
\begin{equation}\label{eq:dereps-delay}
\phi'(\eps)= - \kappa(\eps)\, \| (G_0(\eps),\dots,G_m(\eps)) \|_w < 0,
\end{equation}
which is used in the Newton-based iteration.

\subsection{Numerical examples}

\subsubsection*{A linear system of delay differential equations.}
Consider the linear system of delay equations
\begin{eqnarray} \label{eq:delaysys}
A_0 \dot{x}(t) & = & A_1 x(t) + A_2 x\left(t - 1 \right), \qquad t > 0
\\
x(t) & = & g(t), \qquad t \in \left[-\tau,0\right]
\nonumber
\end{eqnarray}
with $A_0,A_1, A_2$ given $n \times n$ matrices and the delay $\tau =1$.

It is possible to prove that if $A_0$ is nonsingular, \eqref{eq:delaysys} fulfils 
assumptions (i) and (ii).

The associated nonlinear eigenvalue problem is given by
\begin{equation} \label{eq:nlep}
\det \left(\lambda A_0 - A_1 - A_2 \e^{-\lambda}\right) = 0,
\end{equation}
so that - in the considered setting - 
$
f_0(\lambda) = \lambda, \ f_1(\lambda) = -1, \ f_2(\lambda) = -\e^{-\lambda} .
$

We set $A_0 = \Id$ and $w_0 = \infty$ implying this way the identity matrix $A_0$ is not perturbed.
This gives us a constrained gradient system \eqref{ode-E-del} or \eqref{ode-E-1-delay}) with matrices $E_1$ and $E_2$.
Note that $E_0$ is missing because $w_0 = +\infty$.



\subsubsection*{A partial differential equation with a delay.}

Consider the problem (Jarlebring et al., 2010) 
\begin{equation}\label{eq:pdde}
\frac{\partial v(x,t)}{\partial t}=
\frac{\partial^2 v(x,t)}{\partial x^2}+a_0(x)v(x,t)+a_1(x)v(\pi-x,t-1),
\end{equation}
where $a_0(x)=-2\sin(x)$, $a_1(x)=2\sin(x)$, $v_x(0,t)=v_x(\pi,t)=0$.

Space derivatives are approximated by central differences. This gives a delay eigenvalue problem of 
the  form  \eqref{eq:delaysys} with  sparse matrices $A_0$ and $A_1$. 
The dimension is chosen as $n=5000$. Computed pseudospectral abscissae corresponding to different weigths are shown in Table~\ref{tablepde}.

\begin{table}
\begin{center}
{\small
\begin{tabular}{|r|r|l|r|r|r|} \hline
$\alpha$ & $(w_1,w_2)$ & $\eps$ & $\alpha_{\eps}$ & \#steps   \\ \hline
-3.312133337e-1& $(1/2,1/2)$ & 1.e-3 &-3.297978515e-1 &2 \\
&  & 1.e-2 &-3.170982221e-1 &3  \\
&  & 1.e-1 & -1.937166436e-1& 4\\
&  & 1  & 8.647127140e-1 &8\\
  \hline
&  (1/4,$\infty$)& 1.e-3 & -3.300292687e-1  &2 \\
&  & 1.e-2 & -3.193270916e-1 & 3 \\
&  & 1.e-1 & -2.075848221e-1 &4 \\
&  & 1 & +1.641134830e00 & 10 \\ \hline
& ($\infty$,1/4) & 1.e-3 & -3.295667765e-1&2 \\
&  & 1.e-2 & -3.149018637e-1  &3 \\
&  & 1.e-1 & -1.816328080e-1 &4 \\
&  & 1  & +5.599912563e-1  &7 \\ \hline
\end{tabular}
}
\end{center}
\caption{Pseudospectral abscissa $\alpha_{\eps} = -\phi(\eps)$ for the delayed PDE problem.
\label{tablepde}}
\end{table}

\section{Kreiss number}
\label{sec:kreiss}
\index{Kreiss number}

The essential part of the Kreiss matrix theorem (Kreiss \cite{Kr62}) in the sharp version of Spij\-ker (\cite{Spi91}), conjectured by LeVeque \& Trefethen (\cite{LeVT84}), states that if an $n\times n$ complex matrix $A$ has a finite {\it Kreiss number}\footnote{usually called the Kreiss {\it constant}, but in this section we will also consider the dependence on $A$}
\begin{equation}\label{K-kre}
K(A) := \sup_{\mathrm{Re}\, \lambda > 0} \, \Re \,\lambda \ \| (A-\lambda I)^{-1} \|_2 < \infty,
\end{equation}
then its exponential is bounded by
\begin{equation}\label{bound-kre}
K(A) \le \sup_{t\ge 0} \| \e^{tA} \|_2 \le \e\,n\,K(A),
\end{equation}
\bng
where $n$ is the dimension and $\e$ is Euler's number.
\eng
In this section, we consider the following computational problem.

\medskip\noindent
{\bf Problem.} {\it Given a matrix $A\in \C^{n,n}$ having all its eigenvalues in the open complex left half-plane, compute its Kreiss number~$K(A)$.}

\medskip\noindent
Although this is not a matrix nearness problem, we include it in this chapter, because we will address it
with variants of the techniques used for computing the $\eps$-pseudospectral abscissa and the stability radius in Sections~\ref{sec:psa} and~\ref{sec:two-level}. We will describe two such approaches, with emphasis on the second approach.

\subsection{Using pseudospectral abscissae}
\index{pseudospectral abscissa}
For $\eps>0$ let $\alpha_\eps(A)$ 
be the $\eps$-pseudospectral abscissa of $A$. 
We recall formula (\ref{chap:pseudo}.\ref{eps-res-bound}), viz.,
\begin{equation} \label{eps-res-bound-kreiss}
     \max_{\mathrm{Re}\,\lambda \ge \alpha_\eps(A)} \| (A-\lambda I)^{-1} \|_2 = \frac 1\eps.
\end{equation}
Since the maximum is assumed for $\mathrm{Re}\,\lambda = \alpha_\eps(A)$ by the maximum principle, this implies the following formula, which is given by Trefethen \& Embree (\cite{TreE05}), p.\,138.
\begin{lemma}\label{kreiss-number-alpha-eps}
The Kreiss number is determined by the $\eps$-pseudospectral abscissae for ${\eps>0}$ as
\begin{equation}\label{K-psa}
K(A) = \sup_{\eps>0} \alpha_\eps(A)\, \frac 1\eps.
\end{equation}
\end{lemma}
\bcltwo Provided that the pseudospectral abscissa $\alpha_\eps(A)$ is reliably
computed, as is ensured by the criss-cross algorithm described in Section~\ref{sec:psa} (possibly used as a final step after using any of the other, computationally less expensive algorithms described in Section~\ref{sec:psa}),
\ecltwo
it then remains to compute the supremum over all $\eps>0$. While $\alpha_\eps(A)$ grows monotonically with $\eps$, this does not imply that $\alpha_\eps(A)/\eps$ is a unimodal function of $A$, and
little appears to be known about the location and number of local maxima of $\alpha_\eps(A)/\eps$.

\subsection{Using weighted pseudospectra}
\index{pseudospectrum!weighted}
\label{subsec:kreiss-ps}
\bcl In the following we describe an alternative algorithmic approach,
which appears conceptually interesting in that it is based on rephrasing the problem of determining the Kreiss number as a (weighted) matrix nearness problem.
For the nonnegative weight function \ecl 
$$
w(\lambda)=\max(\Re\,\lambda,0),
$$
we consider the {\it weighted $\eps$-pseudospectrum} (for $\eps>0$)
$$
\Lambda_\eps^w(A) = \left\{ \lambda \in \C\setminus \Lambda(A)\,:\, w(\lambda)\, \|(A-\lambda I)^{-1} \|_2 \ge \sfrac1\eps \right\}.
$$
We note that $w(\lambda)\, \|(A-\lambda I)^{-1} \|_2$ is zero in the complex left half-plane, bounded in the right half-plane, and is bounded by 1 asymptotically for ${|\lambda|\to \infty}$. This implies that 
$\Lambda_\eps^w(A)$ lies in the complex right  half-plane or is empty, and $\Lambda_\eps^w(A)$ is a bounded set for every $\eps<1$. There exists $\oeps^w>0$ such that 
$$
\text{$\Lambda_\eps^w(A)=\emptyset\ $ for $\ \eps<\oeps^w\ $ \ but \
$\Lambda_{\oeps^w}^w(A)\ne \emptyset$,}
$$
that is, $\oeps^w$ is the smallest $\eps$ for which $\Lambda_\eps^w(A)$ is not empty.
This means that $1/\oeps^w$ is the smallest upper bound of $w(\lambda)\, \|(A-\lambda I)^{-1} \|_2$ on the complex right half-plane and hence, 
\bcl by definition, \ecl
equals the Kreiss number $K(A)$. We restate this observation as a lemma.

\begin{lemma}\label{lem:kreiss-number}
The Kreiss number is obtained as
\begin{equation} \label{eps-w-kreiss}
     K(A) = \frac 1{\oeps^w}.
\end{equation}
\end{lemma}
We set out to compute $\oeps^w$ using a characterization of $\Lambda_\eps^w(A)$ in terms of eigenvalues of perturbed matrices. We have $\lambda \in \Lambda_\eps^w(A)$ if and only if $\sigma_{\mathrm{min}}(A-\lambda I) \le w(\lambda)\eps$, and the proof of Theorem~\ref{chap:pseudo}.\ref{thm:ps-sv} shows that this is equivalent to the existence of a rank-1 matrix $\Theta$ of 2-norm (or equivalently Frobenius norm) at most $w(\lambda)\eps$ such that $A+\Theta-\lambda I$ is a singular matrix. Writing $\Theta=w(\lambda)\eps E$ with $\|E\|_F=\|E\|_2 \le 1$, we thus obtain the characterization of $\Lambda_\eps^w(A)$ via a nonlinear eigenvalue problem:
\begin{align*}
\Lambda_\eps^w(A) = \bigl\{ \lambda \in \C\setminus \Lambda(A)\,:\,   
& \text{ There exists $E$ with $\|E\|_F \le 1$ such that} \bigr.
\\[-1mm]
\bigl. &\text{ $A+w(\lambda)\eps E - \lambda I$ is singular}\bigr\}.
\end{align*}
This characterization of $\Lambda_\eps^w(A)$ allows us to reformulate the problem of finding $\oeps^w=1/K(A)$.  

\medskip\noindent
{\bf Problem (rephrased as a weighted matrix nearness problem).} {\it Given a matrix $A\in \C^{n,n}$ having all its eigenvalues in the open complex left half-plane, find $\lambda\in \C$ with $\Re\,\lambda\ge 0$ and $\Theta\in \C^{n,n}$ of minimal Frobenius norm such that $A+w(\lambda)\Theta - \lambda I$ is singular.}

\medskip\noindent
Then we have $\|\Theta\|_F = \oeps^w$ and hence $K(A)=1/\|\Theta\|_F$ if $w(\lambda)=\Re\,\lambda$ for $\Re\,\lambda\ge 0$. On the other hand, in the unweighted case $w(\lambda)=1$ the above problem reduces to finding the stability radius of $A$ as $\oeps=\|\Theta\|_F$.

We know that $K(A)\ge 1$ and are mainly interested in the case with strict inequality, which excludes normal matrices.
\begin{assumption} $K(A)>1$.
\end{assumption}
This implies that $\Lambda_\eps^w(A)$ is a nonempty bounded set for $\oeps^w\le\eps<1$.
For such $\eps$ we introduce the right and left weighted $\eps$-pseudospectral abscissae
\begin{align*}
&\alpha_\eps^w(A) = \max \{ \Re\,\lambda\,:\, \lambda \in \Lambda_\eps^w(A) \}<\infty , 
\\
&\beta_\eps^w(A) = \min \ \{ \Re\,\lambda\,:\, \lambda \in \Lambda_\eps^w(A) \} >0.
\end{align*}
We have the weighted variants of \eqref{eps-res-bound-kreiss},
\begin{align*} 
  &   \max_{\mathrm{Re}\,\lambda \ge \alpha_\eps^w(A)} w(\lambda) \,\| (A-\lambda I)^{-1} \|_2 = \frac 1\eps,
     \\
   &  \max_{\mathrm{Re}\,\lambda \le \beta_\eps^w(A)} w(\lambda) \,\| (A-\lambda I)^{-1} \|_2 = \frac 1\eps.
\end{align*}
We now make the further assumption that all global maxima of $w(\lambda) \,\| (A-\lambda I)^{-1} \|_2$ have the same real part. 
\bcl (We expect that generically for a real matrix $A$ there is either a single real global maximum or a pair of conjugate complex maxima, and generically for a complex matrix a single complex global maximum. There are special counterexamples, which we exclude in the following.)
\ecl
\begin{assumption} The set $\{ \Re\,\lambda\,:\, \Re\,\lambda \ \| (A-\lambda I)^{-1} \|_2 = K(A) \}$ consists of a single real number $\rho$.
\end{assumption}
This implies that
$$
 \lim_{\eps\searrow\oeps^w} \alpha_\eps^w(A) = \lim_{\eps\searrow\oeps^w} \beta_\eps^w(A)=\rho .
$$

\subsubsection*{Two-level iteration.}
\index{two-level iteration}
To approximate the Kreiss number $K(A)=1/\oeps^w$, 
\bcl we use one of the two-level iterations discussed in Chapter~\ref{chap:two-level}. We formulate the simple inner-outer iteration although a monotonically convergent variant of the Newton--bisection method or the HEC method might be more reliable and would be guaranteed to give a rigorous lower bound of $K(A)$.
\ecl

\begin{itemize}
\item {\bf Inner iteration:\/} For a given $\eps$ with $\oeps^w\le\eps<1$, we compute the right and left 
weighted $\eps$-pseudospectral abscissae $\alpha_\eps^w(A)$ and $\beta_\eps^w(A)$, analogously to the computation
of the $\eps$-pseudospectral abscissa via a rank-1 gradient flow in Chapter~\ref{chap:proto} or the algorithms in Section~\ref{sec:psa}.

\item {\bf Outer iteration:\/} Given a tolerance parameter $\vartheta>0$, we compute $\eps > \oeps^w$ with
\begin{equation} \label{zero-delta}
\alpha_\eps^w(A)-\beta_\eps^w(A) \le \vartheta,
\end{equation}
and we denote this $\eps$ as $\eps_\vartheta$. This can be done by a Newton--bisection method as in 
Chapter~\ref{chap:two-level}. We then take $K_\vartheta(A)=1/\eps_\vartheta$ as an approximation to $K(A)$.
\end{itemize}

The main difference to the two-level iteration for the stability radius in Chapter~\ref{chap:two-level}
is that we now have to deal with nonlinear eigenvalue problems, where we need to compute eigenvalues $\lambda$ (and corresponding left and right eigenvectors) such that 
\index{nonlinear eigenvalue problem}
$$
A+w(\lambda)\eps E - \lambda I \ \text{ is singular}.
$$
In contrast to the situation considered in the preceding section, here $w$ is not a holomorphic function, and so algorithms for nonlinear eigenvalue problems with holomorphic functions, such as the  method proposed by Beyn (\cite{Bey12}), cannot be applied. See G\"uttel \& Tisseur (\cite{GueT17}) for a review of numerical methods for nonlinear eigenvalue problems. In the present case, where $\eps<1$ (and $\eps\ll 1$ if $K(A)=1/\oeps^w \gg 1$), we can do a (fast) converging iteration where the $(n+1)$st iterate $\lambda_{n+1}$ is determined as a rightmost eigenvalue 
of the perturbed matrix $A+w(\lambda_n)\eps E$. 
\bcl 
This iteration just requires standard eigenvalue algorithms for computing a rightmost eigenvalue of a given matrix as discussed in the Notes of Chapter~\ref{chap:proto}.
\ecl

The whole programme of Chapters~\ref{chap:proto} and~\ref{chap:two-level} then extends in a straightforward way to the present situation.

\subsection{Robustness under structured perturbations}
\label{subsec:kreiss-robust}
So far we considered the Kreiss number of a single given matrix $A$. We now consider the Kreiss numbers of perturbed matrices $A+\Delta$ with (possibly) structured perturbations $\Delta \in \cS$ bounded by
$\|\Delta\|_F\le \delta$ for a given $\delta>0$ and define the supremum of the Kreiss numbers of the collection of such perturbed matrices by
$$
K_\delta^\cS(A) := \sup_{\Delta \in \cS, \|\Delta\|_F\le \delta} K(A+\Delta).
$$
Letting
$$
\eps_\delta^{w,\cS} = \min_{\Delta \in \cS, \|\Delta\|_F\le \delta} \oeps^w(A+\Delta)
$$
with the weight function $w(\lambda)=\Re\,\lambda$ for $\Re\,\lambda\ge 0$, we then have
$$
K_\delta^\cS(A) = \frac 1{\eps_\delta^{w,\cS}}.
$$
The problem of computing $K_\delta^\cS(A)$ can thus be reformulated in the following way.

\medskip\noindent
{\bf Problem (rephrased as a weighted matrix nearness problem).} {\it Given a matrix $A\in \C^{n,n}$ having all its eigenvalues in the open complex left half-plane, and given a structure space $\cS\subset\C^{n,n}$ and $\delta >0$, determine $\eps_\delta^{w,\cS}$ as the smallest $\eps>0$ for which there exist $\lambda\in \C$ with $\Re\,\lambda\ge 0$,  $\Delta \in \cS$ bounded by
$\|\Delta\|_F\le \delta$  and $\Theta\in \C^{n,n}$ bounded by
$\|\Theta\|_F\le \eps$  such that $A+\Delta+w(\lambda)\Theta - \lambda I$ is singular.}

\medskip\noindent
In the unweighted case $w(\lambda)=1$ this problem reduces to finding the inverse common resolvent bound $\eps_\delta^\cS(A)$ studied in Section~\ref{sec:eps-stab}. The weighted case considered here allows for an analogous algorithm. In particular, for the extremal perturbations we again find that $\Delta$ is proportional to $\Pi^\cS \Theta$ and $\Theta$ is of rank~1, which enables us to reduce the problem to that of finding an optimal rank-1 perturbation matrix $\Theta$ as the stationary point of a rank-1 matrix differential equation. The computational cost for computing $K_\delta^\cS(A)$ is then comparable to that for computing just $K(A)$ with the two-level algorithm of the previous subsection.


\section{Notes}


\subsubsection*{Matrix stabilization.}  Finding the smallest (complex, real or structured) stabilizing perturbation to a given matrix is a harder problem than the complementary problem of finding the smallest destabilizing perturbation as discussed in previous chapters.
\bcl
The structured matrix stabilization problem with perturbations of given range and co-range can be NP-hard, as shown by Blondel and Tsitsiklis (\cite{BT97}). 
\ecl

Several conceptually different algorithms for matrix stabilization with complex or real unstructured perturbations have been proposed in the literature.

A black-box approach is to consider the problem of finding the nearest stable matrix as a nonsmooth (but almost everywhere smooth), nonconvex, constrained optimization problem and apply general software for this class of problems, such as given by Curtis, Mitchell \& Overton (\cite{CurMO17}).

Orbandexivry, Nesterov \& Van Dooren (\cite{OrbNVD13}) presented a matrix stabilization algorithm that uses successive convex approximations. They
started from Lyapunov's characterization of stability to reformulate the matrix stabilization problem as finding complex $n\times n$ matrices $X$ and $P$ that give
$$
\inf_{X,P} \tfrac12 \| X-A \|_F^2 \quad\text{ such that $P=P^*$ and $XP+PX^*$ are both positive definite.}
$$
This nonconvex optimization problem is related to the convex problem of finding, for given $X$ and $P$,
$$
\inf_H \tfrac12 \| X+H-A \|_F^2 \quad\text{ such that $H$ is in a suitable ellipsoid defined by $P$ and $X$}.
$$
This update for $X$ is complemented with a procedure that associates an admissible $P$ to~$X$.
With an $O(n^5)$ complexity per iteration, the algorithm is limited to small matrices.

Gillis \& Sharma (\cite{GilS17}) showed that a real square matrix $A$ is stable if and only if it can be written as the matrix of a dissipative Hamiltonian system, i.e. $A = (J-R)Q$, where 
$J$ is skew-symmetric,  $R$ is positive semidefinite and $Q$ is positive definite. This reformulation results in an equivalent nonconvex optimization problem with a convex feasible region onto which points can be projected easily.
The authors proposed a projected gradient method (among other strategies) to solve the problem in the variables $(J,R,Q)$, with $O(n^3)$ complexity per iteration.
Gillis, Karow \& Sharma (\cite{GilKS19}) made an analogous approach to Schur stablilization, based on their characterization of a Schur-stable matrix as being of the form $A = S^{-1}U BS$,
where $S$ is positive definite, $U$ is orthogonal, and $B$ is a positive semidefinite contraction. Choudhary, Gillis \& Sharma (\cite{ChoGS20}) extended the approach to finding the nearest matrix with eigenvalues in more general closed sets $\overline\Omega$ that are a finite intersection of disks, conical sectors and vertical strips.

Noferini \& Poloni (\cite{NofP21}) reformulated matrix stabilization as an optimization problem on the Riemannian manifold of orthogonal or unitary matrices. The problem is then solved using standard methods from Riemannian optimization. The problem of finding the nearest complex Hurwitz-stable matrix is shown to be equivalent to solving
$$
\min_{U\in U(n)} \| L(U^*AU) \|_F^2,
$$
where $L(Z)$ is the lower triangular matrix whose part below the diagonal coincides with that of $Z$, and the diagonal elements are changed to $L(Z)_{ii}=(\Re\, z_{ii})_+$. A related reformulation with orthogonal matrices is given for the real case. The approach is actually formulated for the problem of finding the nearest matrix with eigenvalues in an arbitrary prescribed closed set, thus including Hurwitz- and Schur-stability as special cases.

Guglielmi \& Lubich (\cite{GL17}) studied an ({\it exterior}) two-level approach to matrix stabilization, with a low-rank constrained gradient flow in the inner iteration and a combined Newton--bisection method in the outer iteration. 
This approach moves the eigenvalues in the right half-plane to the imaginary axis. It uses a different functional from the one given in Section~\ref{sec:mat-stab}, which  aims at aligning a fixed number of eigenvalues on the imaginary axis.
The {\it interior} two-level algorithm for matrix stabilization described in Section~\ref{sec:mat-stab} 
is remarkably similar to the two-level algorithm for computing the 
distance to instability. This algorithm has not appeared in the literature before, but it is related to algorithms for Hamiltonian matrix nearness problems and for the passivation of control systems proposed by Guglielmi, Kressner \& Lubich (\cite{GKL15}) and Fazzi, Guglielmi \& Lubich (\cite{FGL21}). In contrast to other methods in the literature, the exterior and interior two-level approaches of Section~\ref{sec:mat-stab} 
can exploit sparsity of the given matrix $A$ in combination with low-rank perturbations and they
apply to matrix stabilization by perturbations with a prescribed linear structure, e.g. for perturbations with a given sparsity pattern.


\subsubsection*{Hamiltonian matrix nearness problems.} Hamiltonian eigenvalue perturbation problems were studied in detail by Mehrmann \& Xu (\cite{MeX08}) and Alam, Bora, Karow, Mehrmann \& Moro (\cite{AlaBKMM11}), motivated by the passivation of linear control systems and the stabilization of gyroscopic mechanical systems, where eigenvalues of Hamiltonian matrices need to be moved to or away from the imaginary axis. In this context, a solution to Problem A considered here yields a lower bound on the distance to non-passivity of a passive system and a solution to Problem B yields a lower bound of the distance to passivity of a non-passive system. The stabilization of a gyroscopic system requires to move all eigenvalues of a Hamiltonian matrix onto the imaginary axis. Those applications have an additional structure of admissible perturbations, which have not been taken into account in this chapter where we consider general real Hamiltonian perturbations. Understanding this general case is, however, basic to addressing the more specific demands of the applications to control systems or mechanical systems. This will become clear in Section~\ref{sec:pass}, where we describe algorithms for finding the nearest passive or non-passive system, which are conceptually close to our treatment of Problems B and A, respectively.

Guglielmi, Kressner \& Lubich (\cite{GKL15}) studied a two-level approach to Problems A and B in the matrix 2-norm instead of the Frobenius norm as considered here. This equally leads to rank-4 differential equations along which the real part of the target eigenvalue decreases (or increases) monotonically. Contrary to the Frobenius norm case, those differential equations cannot be interpreted as constrained gradient systems.

Theorem~\ref{thm:sqrt} on the square root behaviour of the real parts of eigenvalues near a defective coalescence on the imaginary axis was first stated by
Guglielmi, Kressner \& Lubich (\cite{GKL15}). Here we give a corrected proof based on Theorem~\ref{thm:yJx} about the eigenvectors at a defective coalescence, which was first proved by Fazzi, Guglielmi \& Lubich (\cite{FGL21}).

An analogous approach to this chapter can be given for complex Hamiltonian matrices (i.e. matrices $A$ for which $JA$ is hermitian). That case is slightly simpler, as it leads to rank-2 differential equations. Guglielmi, Kressner \& Lubich (\cite{GKL15}) studied the two-level approach to matrix nearness problems in the complex Hamiltonian case for the matrix 2-norm.

Structure-preserving eigenvalue solvers as in the SLICOT library (http://slicot.org/), see  
Benner, Mehrmann, Sima, Van Huffel \& Varga (\cite{BenMSVHV99}),
are distinctly favourable over using a standard general eigenvalue solver, especially for eigenvalues on and close to the imaginary axis as are of interest here; see, e.g., Benner, Losse, Mehrmann \& Voigt
(\cite{BenLMV15}) and references therein.

\subsubsection*{Nearest defective matrix.}
The problem of finding a nearest defective matrix $A'=A+\Delta \in \C^{n,n}$ to a given matrix $A\in \C^{n,n}$ was explicitly addressed by Wilkinson (\cite{Wil84a},\cite{Wil84b}). An instructive review of the history of the problem is given in the introduction to the paper by Alam, Bora, Byers \& Overton (\cite{AlaBBO11}). For the matrix 2-norm, Malyshev (\cite{Mal99}) characterized the distance $d_2(A)$ to the nearest defective matrix as
$$
d_2(A) = \min_{\lambda\in \C} \max_{\gamma\ge 0} \,\sigma_{2n-1} \!
\begin{pmatrix}
    A - \lambda I & \gamma I \\
    0 & A-\lambda I
\end{pmatrix}.
$$
Alam \& Bora (\cite{AlaB05}) showed the following remarkable result for the matrix 2-norm:
\begin{align*}
&\text{\em If two components of the $\eps$-pseudospectrum $\Lambda_\eps(A)$ coalesce at $\lambda_0$, } \\[-1mm]
&\text{\em then $\lambda_0$ is a multiple eigenvalue of a matrix $A'$ such that $\|A - A'\|_2 = \eps$.} 
\end{align*}
Moreover, the matrix $A'$ can be directly computed from a singular value decomposition of $A-\lambda_0 I$. Wilkinson's problem is thus reduced to finding
the smallest $\eps>0$ for which two components of $\Lambda_\eps(A)$ coalesce. It was later shown by Alam, Bora, Byers \& Overton (\cite{AlaBBO11}) that the result of Alam \& Bora (\cite{AlaB05}) extends to the Frobenius norm and that the minimal distance is always attained by a defective matrix.
\bng
Their numerical method is implemented in the code "neardefmat" on Michael Overton's webpage (https://cs.nyu.edu/overton/). The method searches for the lowest saddle point of $\sigma_{\min}(A-zI)$, and it distinguishes between generic and tangential pseudospectral coalescence. 
\bcltwo
Note that no existing method guarantees global convergence to the correct solution.
\ecltwo
\eng

Based on the characterization by Alam \& Bora (\cite{AlaB05}), a computationally efficient Newton method was proposed by
Akinola, Freitag \& Spence (\cite{AkiFS14}).

The {\it real} version of Wilkinson's problem, i.e., to find the nearest defective real matrix to a given real matrix, or more generally {\it structured} versions of Wilkinson's problem, seem to have found little attention in the literature. The two-level algorithm of Section~\ref{sec:defective} with a rank-2 matrix differential equation in the inner iteration is apparently the first algorithm for real and structured versions of Wilkinson's problem.
Finding a small structured perturbation that yields a defective matrix is of substantial current interest in non-hermitian optics and photonics, where the desired ``exceptional points" are synonymous to defective matrices; see e.g.~Chen {\it et al.} (\cite{COZWY17}), Miri \& Alu (\cite{MirA19}), and 
\"Ozdemir, Rotter, Nori \&  Yang (\cite{OezRNY19}). 


\subsubsection*{Nearest singular matrix pencil.}
Byers, He \& Mehrmann (\cite{ByeHM98}) addressed the problem of finding the nearest singular matrix pencil and provided a variety of lower and upper bounds for the distance in the Frobenius norm, both in the complex and real case. Section~\ref{sec:matrix-pencils} modifies the gradient-based approach by Guglielmi, Lubich \& Mehrmann (\cite{GLM17}) in using a different functional that does not depend on eigenvalues. That paper considers the case where both matrices of the pencil $A-\mu B$ are perturbed. This could be done also in the approach of Section~\ref{sec:matrix-pencils} without other than notational complications, but here we allow instead for the restriction to structured perturbations. In a different approach, Riemannian optimization is used by Dopico, Noferini \& Nyman (\cite{DopNN23*}) to compute the nearest (unstructured) singular pencil. Prajaparti \& Sharma (\cite{PraS22}) estimate structured distances to singularity for matrix pencils with symmetry structures.

More general than matrix pencils, it is of interest to compute the nearest
rank-deficient matrix polynomial. This problem has been treated by 
Giesbrecht, Haraldson \& Labahn (\cite{GieHL17}) and 
Das \& Bora (\cite{DasB23}).

The problem of computing the nearest stable matrix pencil is studied by
Gillis, Mehrmann \& Sharma (\cite{GilMS18}) using dissipative Hamiltonian matrix pairs.

\subsubsection*{Stability radii for delay differential equations.}

Stability radii for delay differential equations have been considered in the literature by extending the approaches for ordinary differential equations.

For an extensive discussion we refer the reader to the monograph by Michiels \& Niculescu (\cite{MiNi07}). In particular, in Chapter 1, the robustness of stability and related problems are studied using pseudospectra and stability radii.
Here the technical difficulties are related to the infinite dimensional underlying eigenvalue problem.
We also refer the reader to Michiels, Greeen, Wagenknecht \& Niculescu (\cite{MGWN06}) for pseudospectra and stability radii for delay differential equations.
An algorithm related to the approach proposed in this book has been presented by Michiels and Guglielmi (\cite{MG12}),
based on rank-$1$ iterations.

Finally, a different concept of stability radius, related to the distance to the closest singular matrix-valued analytic function, has been investigated recently by Gnazzo and Guglielmi (\cite{GG23}).
Here the goal is to establish how close is a system of delay differential equations with constant delays to a singular one.

\subsubsection*{Kreiss number of a matrix.}  
Globally convergent algorithms for computing the Kreiss number $K(A)$
have been proposed and studied by Mitchell (\cite{Mit20}, \cite{Mit21}). They are related to algorithms for computing the
distance to uncontrollability of a linear control system and have a computational complexity of $O(n^6)$, which can be reduced to $O(n^4)$ on average by divide-and-conquer variants. 
\bcltwo
The method of Mitchell~(\cite{Mit20}) is related to the distance to uncontrollability and has computational complexity $\bigo(n^6)$. In contrast, the newer methods of Mitchell~(\cite{Mit21}) are fundamentally different: they rely on interpolation-based globality certificates and are significantly faster.
\ecltwo

Apkarian \& Noll (\cite{ApkN20}) take
a different algorithmic approach based on techniques from robust control. They further propose an algorithm for minimizing the Kreiss number in the context of feedback control of transient growth. 
\bng The resulting method, however, is much more expensive. \eng

The reinterpretation as a matrix nearness problem via weighted pseudo\-spectra and the corresponding 
two-level algorithm for computing the Kreiss number as described in Section~\ref{subsec:kreiss-ps} have not appeared in the literature before. Assessing the robustness of the Kreiss number under (possibly structured) perturbations of the matrix, as is done in Section~\ref{subsec:kreiss-robust}, has apparently also not been addressed in the existing literature. 
Analogous algorithms can be given for the time-discrete Kreiss number that is used to bound powers of Schur matrices.

\newcommand{\Hinf}{\mathcal{H}_\infty}
\newcommand{\conjg}[1]{\overline{#1}}
\newcommand{\HinfGc}{\|H\|_{\infty}}
\newcommand{\spec}{\Lambda}
\newcommand{\tfmat}{C(\lambda I-A)^{-1}B+D}
\newcommand{\SVSeps}{\spec_\eps}
\newcommand{\SVSzero}{\spec_0}
\newcommand{\EE}{\Delta}
\newcommand{\xh}{\widetilde x}
\newcommand{\yh}{\widetilde y}
\newcommand{\eigtfmat}{A+B\EE(I-D\EE)^{-1}C}
\newcommand{\beq}{\begin{equation}}
\newcommand{\eeq}{\end{equation}}
\newcommand{\beqs}{\begin{equation*}}
\newcommand{\eeqs}{\end{equation*}}

\chapter{Systems and control}
\label{chap:lti}

In this chapter we reconsider basic problems of robust control of linear time-invariant systems, which we rephrase as eigenvalue optimization problems and matrix (or operator) nearness problems. The problems considered from this perspective include the following:
\begin{itemize}
    \item  computing the $\Hinf$-norm of the matrix transfer function, which is the $L^2$-norm of the input-output map;
    \item  computing the $\Hinf$-distance to uncontrollability of a controllable system;
    \item  passivity enforcement by a perturbation to a system matrix that minimizes
    the $\Hinf$-norm of the perturbation to the matrix transfer function; and
    \item  computing the distance to loss of contractivity under structured perturbations to the state matrix.
\end{itemize}
The distances from systems with undesired properties are important robustness measures of a given control system, whereas algorithms for finding a nearby control system with prescribed desired properties are important design tools. The algorithmic approach to eigenvalue optimization via low-rank matrix differential equations and the two-level approach to matrix nearness problems of previous chapters is extended to a variety of exemplary problems from the area of robust control and is shown to yield versatile and efficient algorithms.

We consider the {\it continuous-time linear time-invariant dynamical system} with inputs $u(t)\in \C^p$, outputs $y(t)\in\C^m$ and states $z(t)\in\C^n$ related by
\index{linear time-invariant system}
\begin{eqnarray}
\dot z(t) & = & A z(t) + Bu(t)
\label{lti}
\\
y(t) & = & C z(t) + Du(t)
\nonumber
\end{eqnarray}
with the real system matrices $A \in \R^{n,n}$, $B \in \R^{n,p}$, $C \in \R^{m,n}$ and $D \in \R^{m,p}$, and with the initial state $z(0)=0$. In this chapter we always assume that all eigenvalues of $A$ have negative real part.

In the final section of this chapter we consider {\it descriptor systems}, 
\index{descriptor system}
where $\dot z(t)$ in \eqref{lti} appears multiplied with a singular matrix $E\in\R^{n,n}$, which yields a differential-algebraic equation instead of the differential equation in \eqref{lti}. Descriptor systems play an essential role in modeling and composing networks of systems. We present an algorithm for computing the $\Hinf$-norm of the associated matrix transfer function, which now needs to be appropriately weighted if the descriptor system has an index (to be defined shortly) higher than 1.
 
\section{$\Hinf$-norm of the matrix transfer function}\label{sec:Hinf}
\index{transfer function}

The matrix-valued {\it transfer function} associated with the system \eqref{lti} is  
\beq
   H(\lambda) = \tfmat \quad  \mathrm{for~} \lambda \in \C\backslash\spec(A) 
   \label{tfmatdef}
\eeq
where $\spec(A)$ denotes the spectrum of $A$.
If all eigenvalues of $A$ have negative real part, as we assumed, then $H$ is a matrix-valued holomorphic function on a domain that includes the closed complex right half-plane $\Re\,\lambda \ge 0$. 

The input--output map $u\mapsto y$ given by the variation-of-constants formula,
\[ 
y(t)= \int_0^t C e^{(t-\tau)A}B u(\tau)\,d\tau + Du(t), \quad t\ge 0,
\]
is the convolution with the inverse Laplace transform of $H$:
\[
y = H(\partial_t)u := (\mathcal{L}^{-1} H)* u.
\]
Taking Laplace transforms, we get the formula that explains the name ``transfer function'',
\beq
\label{Hu}
\mathcal{L}y(\lambda) = H(\lambda)\, \mathcal{L}u(\lambda), \quad \Re\,\lambda \ge 0.
\eeq
The Plancherel formula for the Fourier transform $(\mathcal{F}y)(\omega)=(\mathcal{L}y)(\iu\omega)$ for $\omega\in\R$
(where $y$ is extended by zero to $t<0$) then yields that the operator norm of the input--output map $H(\partial_t): L^2(0,\infty;\C^p) \to L^2(0,\infty;\C^m)$ equals $\sup_{\omega\in\R} \|H(\iu\omega)\|_2$, since
\begin{align*}
\int_0^\infty \|y(t)\|^2 \,dt &=
\int_{\R} \| (\mathcal{L}y)(\iu\omega) \|^2 \,d\omega =
\int_{\R} \| H(\iu\omega) (\mathcal{L}u)(\iu\omega)\|^2 \, d\omega
\\
&\le \sup_{\omega\in\R} \|H(\iu\omega)\|_2^2 \int_{\R}\|(\mathcal{L}u)(\iu\omega)\|^2 \, d\omega = \sup_{\omega\in\R} \|H(\iu\omega)\|_2^2
\int_0^\infty \|u(t)\|^2 \,dt
\end{align*}
and an approximate $\delta$-function argument shows that $\sup_{\omega\in\R} \|H(\iu\omega)\|_2^2$ is the smallest such bound that holds for all square-integrable functions $u$. 
\begin{definition} \label{def:hinfnormcont}
The \emph{$\Hinf$-norm} of the matrix transfer function $H$ is 
\beq
  \HinfGc 
 :=\sup_{\mathrm{Re}\,\lambda \ge 0} \|H(\lambda)\|_2 = \sup_{\omega\in\R} \|H(\iu\omega)\|_2.
                              \label{hinfnormequivcont}
\eeq
\end{definition}
\index{H$_\infty$-norm}
The supremum is a maximum  if $\| D\|_2 = \|H(\infty)\|_2$ is strictly smaller than $\HinfGc$, as will be assumed from now on. The second equation in \eqref{hinfnormequivcont} is a consequence of the maximum principle, which can be applied on noting that $\|H(\lambda)v\|^2$ is a subharmonic function of $\lambda$ for every $v\in\C^n$.

As the $L^2$ operator norm of the input-output map $H(\partial_t)$, the $\Hinf$-norm of $H$ is a fundamental stability measure.
Moreover, using the causality property that $y(t)$ depends only on $u(\tau)$ for $\tau\le t$, we can rewrite the above bound as
\begin{equation}
    \label{yu-bound}
    \biggl( \int_0^T \| y(t) \|^2 \, dt\biggr)^{1/2} \le \| H \|_\infty \, \biggl(\int_0^T \| u(t) \|^2 \, dt\biggr)^{1/2}, \qquad
    0\le T \le \infty,
\end{equation}
and $\| H \|_\infty $ is the smallest such bound.

\pagebreak[3]
\bigskip\noindent
{\bf Problem.} {\it Compute the $\Hinf$-norm of the matrix transfer function $H$.}

\medskip\noindent
Making use of the theory of spectral value sets, which are suitable extensions
of pseudospectra presented in the monograph by Hinrichsen and Pritchard  (\cite{HinP05}), 
we will show that  the $\Hinf$-norm equals the reciprocal of the stability radius, which here is the largest value of $\eps$ such that the associated $\eps$-spectral value set is contained in the complex left half-plane. 

This characterization allows us to extend the algorithmic approach of previous sections, using rank-1 constrained gradient systems, from pseudospectra to spectral value sets, and then use a scalar Newton--bisection method to approximate the $\Hinf$-norm.

\subsection{Matrix transfer function and perturbed state matrices}\label{subsec:tf-psm}

We start with discussing the relationship between the singular vectors of the transfer matrix $H(\lambda)$ and the eigenvectors of a corresponding set of matrices.

Given $A,B,C,D$ defining the linear dynamical system \eqref{lti}, consider the \emph{perturbed state matrix}, for perturbations $\EE\in \C^{p,m}$ such that $I-D\EE$ is invertible,
\beq
    M(\EE) = \eigtfmat \in \C^{n,n} \label{AEdef}
\eeq
and the associated transfer matrix \eqref{lti}.
The next theorem relates the $2$-norm of the transfer matrix, which is its largest singular value, to eigenvalues of perturbed state matrices. This result extends the basic characterization of complex pseudospectra given in Theorem~\ref{chap:pseudo}.\ref{thm:ps-sv} together with (\ref{chap:pseudo}.\ref{ps-res}), to which it reduces for $B=C=I$ and $D=0$.

\begin{theorem}[Singular values and eigenvalues]\label{thm:basicequiv-sv-eig}
Let $\eps > 0$ and $\eps\|D\|_2 < 1$.
Then, for $\lambda\not\in\spec(A)$ the following two statements are equivalent:
\begin{itemize}
    \item[(i)]
    $\ \ \|H(\lambda)\|_2 \geq \eps^{-1}$ 
    \item[(ii)]\ $\ \lambda$ is an eigenvalue of  $M(\EE)$ for some $\EE\in\C^{p,m}$ with $\|\EE\|_2 \leq \eps$. \label{twoequiv}
\end{itemize}
Moreover, $\Delta$ can be chosen to have rank $1$, and the two inequalities can be replaced by equalities in the equivalence.
\end{theorem}

\begin{proof} We first observe that under the condition $\eps\|D\|_2<1$ we have
that $I-D\EE$ is invertible when $\|\EE\|_2\leq\eps$, and hence $M(\EE)$ is then well-defined.

Suppose (i) holds true, with $\rho=\|H(\lambda)\|_2^{-1}\leq \eps$.
Let $u$ and $v$ be right and left singular vectors of $H(\lambda)$, respectively, corresponding to the largest
singular value
$\rho^{-1}$, so that 
\[
\rho H(\lambda)u=v, \quad \rho v^*H(\lambda)=u^*, \quad \mbox{and} \quad \|u\|=\|v\|=1.
\]
Define $\EE=\rho uv^*$ so that $\|\EE\|_2 =\rho\leq\eps$.
We have $H(\lambda)\EE=vv^*$, so
\beq
        (\tfmat)\EE v=v.  \label{tfmatEv}
\eeq
Next define $Y=(I-D\EE)^{-1}C$ and $Z=(\lambda I - A)^{-1}B\EE$, so we have $YZ v=v$.
It follows that $Zv\ne 0$ and $ZY y=y$, with $y:=Zv  =\rho(\lambda I - A)^{-1}Bu$ an eigenvector of $ZY$. 
Multiplying through by $\lambda I-A$, we have
\beq
       B\EE(I-D\EE)^{-1}Cy = (\lambda I - A)y,  \label{xeigvec-sc-1}
\eeq
which is equivalent to $M(\Delta)y=\lambda y$. 
\bng
This proves (ii). 
\eng

Conversely, suppose that (ii) holds true. Then there exists 
$y\not= 0$ such that
\eqref{xeigvec-sc-1} holds. We have $ZYy=y$, so $y$ is an eigenvector of $ZY$ corresponding to
the eigenvalue 1. Consequently, $YZw=w$ where $w=Yy\not = 0$ is an eigenvector of $YZ$.
Multiplying by $I-D\EE$ and rearranging we have
\[
    (\tfmat)\EE w = w, \quad \text{i.e.},\quad
    H(\lambda)\Delta w = w.
\]
This implies
\[
    \eps \|H(\lambda)\|_2 \geq \|H(\lambda)\EE\|_2 \geq 1,
\]
which proves the first statement in \eqref{twoequiv}.

The equivalence \eqref{twoequiv} also holds if we restrict $\EE$ in the second statement
to have rank one. The proof remains unchanged.
\qed
\end{proof}

We reformulate the remarkable relationship between eigenvectors of $M(\EE)$ and singular vectors of $H(\lambda)$ revealed by the previous proof in a separate corollary.

\begin{corollary}[Singular vectors and eigenvectors] \label{thm:evecssvecs}
Let $\eps > 0$ and $\eps\|D\|_2 < 1$, and
let $u\in\C^p$ and $v\in\C^m$ with $\|u\|=\|v\|=1$ be right and left singular vectors
of $H(\lambda)$, respectively, corresponding to a singular value $\eps^{-1}$. 
Then,  the nonzero vectors
\beq
       \yh= (\lambda I - A)^{-1}Bu \quad \mathrm{and} \quad 
       \xh= (\lambda I - A)^{-*}C^*v, \label{xyformulas}
\eeq
 are  (non-normalized)  right and left eigenvectors 
 associated with the eigenvalue~$\lambda$ 
 of $M(\EE)$ for $\EE=\eps uv^*$.
\end{corollary}

\begin{proof}
In the proof of Theorem~\ref{thm:basicequiv-sv-eig} we showed that $\yh$ is a right eigenvector of $M(\EE)$ for $\EE=\eps uv^*$ to the eigenvalue $\lambda$. The proof for the left eigenvector $\xh$ is analogous.
\qed
\end{proof}

\subsection{Spectral value sets}\label{subsec:specvalsets}

We define spectral value sets, which generalize the notion of pseudospectrum of a matrix $A$ to linear control systems with the matrices $(A,B,C,D)$, and we
show their relationship with the 2-norm of the matrix transfer function.

\begin{definition}
Let $\eps \geq 0$ and $\eps\|D\|_2 < 1$, and define the \emph{spectral value set}
\beqs
            \SVSeps(A,B,C,D) = \bigcup \left\{\spec(M(\EE)) : \EE\in\C^{p, m}, \|\EE\|_2\leq \eps\right\}. 
\eeqs 
\end{definition}
\index{spectral value set}
Note that $\SVSeps(A,B,C,D) \supset \SVSzero(A,B,C,D) = \spec(A)$, and note further
that $\SVSeps(A,I,I,0)$ equals the $\eps$-pseudospectrum $\Lambda_\eps(A)$.
The following corollary of Theorem \ref{thm:basicequiv-sv-eig} is immediate.
\begin{corollary}[Characterization of the spectral value set]\label{cor:SVSequiv}
Let $\eps > 0$ and $\eps\|D\|_2 < 1$.  Then,
\begin{eqnarray*}
  \SVSeps(A,B,C,D)\backslash\spec(A) 
                  & = &\bigcup\left\{ \spec(M(\EE)) : \EE\in\C^{p, m}, \|\EE\|_2\leq \eps, \mathrm{rank}(\EE)=1\right\}
                  \\
        &=& \bigcup \left\{\lambda\in \C\backslash\spec(A): \|H(\lambda)\|_2 \geq \eps^{-1} \right\}.
\end{eqnarray*}
\end{corollary}

\subsection{$\Hinf$-norm and stability radius}\label{sec:Hinfcont}

For $\eps\geq 0$ with $\eps\|D\|_2<1$, the \emph{spectral value set abscissa} is
\beq  
         \aleps(A,B,C,D) =  \max\{\mathrm{Re}~\lambda : \lambda \in \SVSeps(A,B,C,D)\} \label{alepsdef}
\eeq
with $\alpha_0(A,B,C,D)=\alpha(A)$, the spectral abscissa of $A$. This definition extends the notion of the pseudospectral abscissa $\alpha_\eps(A)$.



\index{stability radius}
The $\Hinf$-norm can be characterized as the reciprocal of the \emph{stability radius}, which is the largest $\eps$ such
that $\SVSeps(A,B,C,D)$ is contained in the complex left half-plane.
The following theorem states this remarkable equality on which 
our algorithmic approach to computing the $\Hinf$-norm will be based. It extends (\ref{chap:pseudo}.\ref{oeps-res-bound}), to which it reduces for $B=C=I$ and $D=0$.

\begin{theorem}[$\Hinf$-norm via the stability radius]\label{thm:Hinf-oeps}
Assume that all eigenvalues of $A$ have negative real part. Let the stability radius of the system $(A,B,C,D)$ be
$$
\eps_\star:=\inf\{\eps>0\ \text {with }\  \eps \|D\|_2 < 1 \,:\,  \aleps(A,B,C,D) = 0\},
$$
where $\aleps(A,B,C,D)$ is the spectral value set abscissa defined in \eqref{alepsdef}.
Then,
\begin{equation}
  \HinfGc = \frac1{\eps_\star}.
\label{hinfnormdefcont}
\end{equation}
\end{theorem}
\begin{proof}
We first consider the case where $\SVSeps(A,B,C,D)$ does not intersect the
imaginary axis for any $\eps>0$ with $\eps\|D\|_2 < 1$. Then we take the infimum
in \eqref{hinfnormdefcont} to be $1/\|D\|_2$.
By Corollary~\ref{cor:SVSequiv}, 
$\|H(\iu \omega)\|<\eps^{-1}$ for
all $\omega\in\R$ and all $\eps$ with $\eps\|D\|_2<1$, and hence the supremum in \eqref{hinfnormequivcont}
is at least $\|D\|_2$, and therefore equal to $\|D\|_2$ as is seen by letting $\omega\rightarrow\pm\infty$. So we have equality in \eqref{hinfnormdefcont} in this degenerate case.

Otherwise, there exists a smallest $\eps_\star$ with $\eps_\star \|D\|_2<1$ such that $\alpha_{\eps_\star}(A,B,C,D)=0$. So there exists $\omega_\star\in\R$ such that 
$\iu\omega_\star\in \Lambda_{\eps_\star}(A,B,C,D)$. 
By Corollary~\ref{cor:SVSequiv}, this implies 
$\| H(\iu\omega_\star)\|_2\ge 1/\eps_\star$. Here we have actually equality, because  $\| H(\iu\omega_\star)\|_2= 1/\eps$ with $\eps<\eps_\star$ would imply, again by Corollary~\ref{cor:SVSequiv}, that $\iu\omega_\star\in \Lambda_{\eps}(A,B,C,D)$ and hence $\alpha_{\eps}(A,B,C,D)=0$, which contradicts the minimality of $\eps_\star$.
\qed
\end{proof}


It follows
from Corollary \ref{cor:SVSequiv} that, for $\eps>0$ with $\eps\|D\|_2<1$,
the  spectral value set abscissa in \eqref{alepsdef} equals
\beq
     \aleps(A,B,C,D)=\max \left\{\Re~\lambda : \lambda\in\spec(A) \mathrm{~or~}\|H(\lambda)\|_2 \geq \eps^{-1} \right\}.  \label{alepsdef2}
\eeq
The set of admissible $\lambda$ must include $\spec(A)$ because of the possibility that
the spectral value set $\SVSeps(A,B,C,D)$ has isolated points. Excluding such points,
we obtain local optimality conditions for \eqref{alepsdef2}.

In order to proceed we make the following generic assumptions.
\begin{assumption} \label{assumptcont}
Let $\eps > 0$ with $\eps\|D\|_2 < 1$, and
let $\lambda\not\in\spec(A)$ be a locally rightmost point of $\SVSeps(A,B,C,D)$.  We assume:
\begin{enumerate}
\item The largest singular value $\eps^{-1}$ of $H(\lambda)$ is simple. 
\item Letting $u$ and $v$ be corresponding right and left singular vectors and
setting $\EE=\eps uv^*$, the eigenvalue $\lambda$ of $M(\EE)$ is simple. 
\end{enumerate}
\end{assumption}
Here we note that $\eps^{-1}$ equals the largest singular value of $H(\lambda)$, i.e.~$\|H(\lambda)\|_2$, by Corollary~\ref{cor:SVSequiv} and the minimality argument at the end of the proof of Theorem~\ref{thm:Hinf-oeps}, and $\lambda$ is an eigenvalue of $M(\EE)$ by 
Theorem~\ref{thm:basicequiv-sv-eig}.

\begin{lemma}[Eigenvectors with positive inner product]
\label{lem:firstordercont}
Let $\eps > 0$ with $\eps\|D\|_2 < 1$, and
let $\lambda\not\in\spec(A)$ be a locally rightmost point of $\SVSeps(A,B,C,D)$.
Under Assumption \ref{assumptcont}, we then have that
\beq
     \widetilde x^* \widetilde y
     \ \text{ is real and positive}, \label{firstordercont}
\eeq
where $\widetilde x$ and $\widetilde y$ are the (non-normalized) left and right eigenvectors to the eigenvalue $\lambda$ of $M(\Delta)$ with $\Delta=\eps uv^*$ that, via
\eqref{xyformulas}, correspond to the
left and right singular vectors $u$ and $v$ associated with the largest singular value $\eps^{-1}$ of $H(\lambda)$.
\end{lemma}
\begin{proof}
The standard first-order necessary condition
for $\widehat\zeta\in\R^2$ to be a local maximizer of an optimization problem
$
\max\{ f(\zeta): g(\zeta)\leq 0,~\zeta\in \R^2 \},
$
when $f$, $g$ are continuously differentiable and $g(\widehat\zeta)=0$,
$\nabla g(\widehat\zeta)\not=0$, is the existence of a Lagrange multiplier
$\mu\geq 0$ such that $\nabla f(\widehat\zeta)=\mu\nabla g(\widehat\zeta)$.  In our
case, identifying $\lambda\in\C$ with $\zeta\in\R^2$, the gradient of the maximization
objective is $(1,0)^T$, while the constraint function
\[
        \frac{1}{\eps} - \|C\left( \lambda I - A \right)^{-1}B + D\|_2
\]
is differentiable with respect to $\lambda$ because of the first part of Assumption \ref{assumptcont},
and it has the gradient
\[
    \left ( \begin{array}{c}\Re(v^* C( \lambda I - A )^{-2}B u) \\[2mm]
                            \Im(v^* C( \lambda I - A )^{-2}B u)\end{array} \right )
\]
using standard perturbation theory for singular values.
Defining $\Delta=\eps uv^*$ and applying Theorem \ref{thm:evecssvecs} we know that 
$\xh$ and $\yh$
as defined in \eqref{xyformulas} 
are left and right eigenvectors of $M(\EE)$,
with inner product
\beq
        \xh^*\yh = v^* C (\lambda I - A)^{-2} B u.            \label{evecip}
\eeq
By the second part of Assumption \ref{assumptcont}, $\lambda$ is a simple eigenvalue of $M(\EE)$ and so $\xh^*\yh \not = 0$. Therefore, the constraint gradient is nonzero implying that
the Lagrange multiplier $\mu > 0$ exists with
$v^* C( \lambda I - A )^{-2}B u = 1/\mu>0$, and by \eqref{evecip} we thus find $\xh^*\yh>0$.
\qed
\end{proof}

\subsection{Two-level iteration}
\index{two-level iteration}

Like for the matrix nearness problems in Chapter~\ref{chap:two-level}, we approach the computation of the $\Hinf$-norm by a 
two-level method:
\begin{itemize}

\item {\bf Inner iteration:\/} Given $\eps>0$, we aim to compute a  matrix $E(\eps)$  of rank 1 and 
of unit Frobenius norm,  such that the functional
\[
\F_\eps(E) = -\frac{\lambda+\clambda}{2} = {}-\Re(\lambda), \quad\text{ for } \lambda=\lambda(M(\eps E)),
\] 
where $\lambda(M)$ is the rightmost eigenvalue of a matrix $M$, is minimized in the manifold of rank-$1$ matrices of unit norm, i.e. 
\begin{equation} \label{E-epsl}
E(\eps) = \arg\min\limits_{E \in \cM_1, \| E \|_F = 1} \F_\eps(E).
\end{equation}

The obtained optimizer is denoted by $E(\eps)$ to emphasize its dependence on $\eps$,
and the rightmost eigenvalue of $M\left( \eps E(\eps) \right)$ is denoted by $\lambda(\eps)$. We then have
the $\eps$-spectral value set abscissa $\Re\,\lambda(\eps)=\alpha_\eps(A,B,C,D)$.
(Note that the Frobenius norm and the matrix 2-norm are the same for a rank-1 matrix.)

\item {\bf Outer iteration:\/} We compute the smallest positive value $\oeps$ with
\begin{equation} \label{eq:zero-sc}
\phi(\oeps)= 0,
\end{equation}
where $\phi(\eps)= \F_\eps\left(E(\eps) \right) = - \Re \, \lambda(\eps) = - \aleps(A,B,C,D)$ is minus the spectral value set abscissa. 
\end{itemize}
If the numerical result computed by such a two-level iteration were exact, it would yield the $\Hinf$-norm in view of Theorem~\ref{thm:Hinf-oeps},
\[
\HinfGc = \frac{1}{\oeps}.
\]


\subsection{Norm-constrained gradient flow}\label{sec:specvalsetabsc}
\index{gradient flow!norm-constrained}

In this and the next subsection we show how to deal with the inner iteration, following a programme that directly extends the programme of Section~\ref{sec:proto-complex}.

As in \eqref{AEdef}, consider the perturbed matrix, for $\Delta\in\C^{p,m}$,
\[
      M\left( \Delta \right) = A + B \Delta \left( I- D \Delta \right)^{-1}C.
\] 
We consider the eigenvalue optimization problem \eqref{E-epsl}, but for the moment with respect to all complex perturbations $E \in \C^{p,m}$ of unit Frobenius norm, although later
we will restrict to rank-$1$ matrices $E$.
So  we look for
\begin{equation} \label{eq:optim-sys}
\arg\min\limits_{E \in \C^{p,m}, \, \| E \|_F = 1} \F_\eps(E).
\end{equation}
To treat this eigenvalue optimization problem,
we will closely follow the course of Section~\ref{sec:proto-complex} and adapt it to the present situation.

\bigskip\noindent
{\bf Free gradient.}
To derive the gradient of the functional $\F_\eps$, we first state a simple auxiliary result.
\begin{lemma}[Derivative of the perturbed matrix]
Given a smooth matrix valued function $\Delta(t)$ with $\|\Delta(t)\|_2 \| D \|_2  < 1$,
we have
\begin{equation}
\frac{d}{d t} M(\EE(t))=
B\left( I - \Delta(t) D \right)^{-1} \dot \Delta(t)  \left( I - D \Delta(t) \right)^{-1}C.
\end{equation}
\label{lem:derE}
\end{lemma}

\begin{proof}
For conciseness, we omit the dependence on $t$, differentiate and regroup terms as
\begin{align}\label{eq:matrix-Ft-deriv}
\frac{d}{d t} \biggl( \Delta \left( I - D \Delta \right)^{-1} \biggr) &= \dot{\Delta}  \left( I - D \Delta \right)^{-1} +
\Delta \frac{d}{d t} \left( I - D \Delta \right)^{-1} \nonumber \\ 
&= \dot{\Delta} \left( I - D \Delta \right)^{-1} +
\Delta \left( I - D \Delta \right)^{-1} D\dot{\Delta} \left( I - D \Delta \right)^{-1} \nonumber \\
&= \left(I+ \Delta \left( I - D \Delta \right)^{-1} D \right) \dot{\Delta} \left( I - D \Delta \right)^{-1}.
\end{align}
We then observe that
\begin{equation}\label{eq:matrix-inf-series}
I+ \Delta \left( I - D \Delta \right)^{-1} D = I + \Delta \biggl( \sum\limits_{k=0}^{\infty} (D \Delta)^k \biggr) D
= I + \sum\limits_{k=1}^{\infty} (\Delta D)^k = \left( I - \Delta D \right)^{-1}.
\end{equation}
Combining \eqref{eq:matrix-Ft-deriv} and \eqref{eq:matrix-inf-series} yields 
\begin{equation}
\frac{d}{d t} \biggl( \Delta(t) \left( I - D \Delta(t) \right)^{-1} \biggr)=
\left( I - \Delta(t) D \right)^{-1} \dot \Delta(t)  \left( I - D \Delta(t) \right)^{-1},
\end{equation}
which implies the result.
\qed
\end{proof}

We will from now on use a normalization of the eigenvectors of $M(\Delta)$ of \eqref{AEdef}, which we previously considered in this book:
\begin{equation} \label{eq:scaling-sys}
\| x \| = \| y \| = 1 \quad\text{ and } \quad \ x^* y \, \text{ is real and positive.}
\end{equation}
In the following we work with the normalized left and right eigenvectors $x$ and $y$ instead of the non-normalized $\widetilde x$ and $\widetilde y$ of Corollary~\ref{thm:evecssvecs}.


\begin{lemma}  [Derivative of a simple eigenvalue] 
\label{lem:lambdaderiv-sys}  
Let  $\Delta(t)$ be a smooth matrix valued function with $\|\Delta(t)\|_2 \| D \|_2 \le 1$.
Suppose that $\lambda(t)$ is a simple eigenvalue of $M(\Delta(t))$ depending continuously on t, with associated eigenvectors
$x(t)$ and $y(t)$ normalized according to \eqref{eq:scaling-sys}, and let 
$\kappa(t) = 1/(x(t)^*y(t)) > 0$. 
Then, $\lambda(t)$ is differentiable with
\[
            \dot{\lambda}(t) = \kappa(t)\, r(t)^* \dot{\EE}(t) s(t)
\label{eq:optprob2}
\]
with 
\beq \label{eq:rsdef}
      r(t) = \left( I - \Delta(t) D \right)^{-*} b(t), \quad 
      s(t) = \left( I -  D \Delta(t) \right)^{-1} c(t),
\eeq
with $b(t)=B^*x(t) \quad \mbox{and} \quad c(t)=C y(t).$
\end{lemma}

\begin{proof}
Applying Theorem \ref{chap:appendix}.\ref{thm:eigderiv} we get 
\beq \nonumber
\dot{\lambda} = \frac{x^* \dot{M} y}{x^* y} 
\eeq
with
\beq
\dot{M} =
         B \left(I - \Delta D \right)^{-1} \dot{\EE} \left(I -  D \Delta \right)^{-1} C,
\label{eq:optprob}
\eeq
where we omitted the dependence on $t$ for brevity. The result is then immediate.
\qed
\end{proof}
A direct consequence of Lemma~\ref{lem:lambdaderiv-sys}, for $\Delta(t)=\eps E(t)$, is that
\beq
\label{F-der-sys}
\frac1{\eps\kappa(t)}\,\frac{d}{dt} \F_\eps(E(t))=  \Re\langle G_\eps(E(t)), \dot{E}(t)  \rangle
\eeq
with the (rescaled) gradient given by the rank-1 matrix
\beq\label{grad-sys}
G_\eps(E) =  -r s^*.
\eeq
Let 
\begin{equation}
\psi_\eps=\frac{\eps}{1-\eps v^*Du} .   
\label{eq:psieps}    
\end{equation}
For $E = u v^*$, using the formulas 
\begin{equation} \nonumber 
\left( I - \eps u v^* D \right)^{-1} = I + \psi_\eps u v^* D, \qquad  
\left( I - \eps D u v^* \right)^{-1} = I + \psi_\eps D u v^*
\end{equation}
we get (with $\beta = u^* b, \gamma=v^*c$)
\begin{eqnarray} \nonumber
r & = & \left( I + \conj{\psi_\eps} D^* v u^* \right) b = b + \conj{\psi_\eps} \beta D^* v
\\
\nonumber
s & = & \left(  I + \psi_\eps D u v^* \right) c = c +  \psi_\eps \gamma D u.
\end{eqnarray}

\medskip\noindent
{\bf Norm-constrained gradient flow.} 
In the same way as in Section~\ref{sec:proto-complex},
this suggests to consider the following constrained gradient flow on the manifold of  $p\times m$ complex matrices of unit Frobenius norm:
\beq\label{ode-E-sys}
\dot E = -G_\eps(E) + \Re \left\langle G_\eps(E),  E \right\rangle E,
\eeq
where
$(\lambda(t)$, $x(t)$, $y(t))$ is a rightmost eigentriple for 
the matrix $M( \eps E(t))$.
Assume the initial condition $E(0)=E_0$,
a given matrix with unit Frobenius norm, chosen so that $M(\eps E_0)$ has a unique rightmost eigenvalue $\lambda(0)$, which is simple.

We can now closely follow the programme of Section~\ref{chap:proto}.\ref{sec:proto-complex} with straightforward minor adaptations. 

\medskip\noindent
{\bf Monotonicity.} 
Assuming simple eigenvalues along the trajectory of \eqref{ode-E-sys},
we again have the monotonicity property of Theorem~\ref{chap:proto}.\ref{thm:monotone},
\begin{equation}
\frac{d}{dt} \F_\eps (E(t))  \le  0.
\label{monotone-sys}
\end{equation}

\medskip\noindent
{\bf Stationary points.}
\index{stationary point}
Also the characterization of stationary points as given in Theorem~\ref{chap:proto}.\ref{thm:stat} extends with the same proof: Let
$E\in\C^{p,m}$ with $\| E\|_F=1$ be such that the rightmost eigenvalue $\lambda$ of $M(\eps E)$ is simple
and $r,s\ne 0$. Then, 
\begin{equation}\label{stat-sys}
\begin{aligned}
&\text{$E$ is a stationary point of the differential equation \eqref{ode-E-sys}}
\\[-1mm]
&\text{if and only if $E$ is a real multiple of $rs^*$.}
\end{aligned}
\end{equation}



Since local minima of $\F_\eps$ are necessarily 
stationary points of the constrained gradient flow \eqref{ode-E-sys}, this immediately implies the following.

\begin{corollary}[Rank of optimizers] \label{cor:rank-1-sys}
If $E$ is an optimizer of problem \eqref{eq:optim-sys} and we have $r,s \ne 0$, 
then $E$ is of rank $1$.
\end{corollary}
\index{optimizer!rank-1 property}

As in Section~\ref{sec:proto-complex}, 
Corollary~\ref{cor:rank-1-sys} motivates us to project the differential equation \eqref{ode-E-sys} onto the manifold $\cM_1$ of rank-$1$ matrices, which is computationally favourable.

\subsection{Rank-$1$ constrained gradient flow}
\label{subsec:rank-1-sys}
\index{gradient flow!rank-1 constrained}

Since stationary points of \eqref{ode-E-sys} have rank $1$,  we consider - 
as we have done in Chapter \ref{chap:proto}  - 
the differential equation \eqref{ode-E-sys} projected  onto the tangent space $T_E\cM_1$ at $E$ of the rank-$1$ manifold.

We recall that the orthogonal projection with respect to Frobenius inner product $\langle \cdot , \cdot \rangle$ from $\C^{p, m}$ onto the tangent space $T_E\cM_1$ at $E=uv^* \in\cM_1$ (with $u$ and $v$ of Euclidean norm 1) is given by
Lemma~\ref{chap:proto}.\ref{lem:P-formula-1} as
\begin{equation}\label{P-formula-sys}
P_E(Z) = Z - (I-uu^*) Z (I-vv^*)
\quad\text{ for $Z\in\C^{p, m}$}.
\end{equation}
As in Section~\ref{subsec:rank1-gradient-flow} we consider the projected gradient system
\begin{equation} \label{ode-E-1-sys}
\dot{E} =  - P_E \Bigl(G_\eps(E) + \Re\hspace{1pt}\langle P_E(G_\eps(E)), E \rangle\, E   \Bigr)
\end{equation}
and find  the following properties, again by the arguments of Section~\ref{subsec:rank1-gradient-flow}.

\medskip\noindent
{\bf Conservation of unit norm.} Solutions $E(t)$ of \eqref{ode-E-1-sys} have Frobenius norm 1 for all $t$, provided that the initial value $E(0)$ has Frobenius norm 1.

\medskip\noindent
{\bf Differential equations for the two vectors.}
For an initial value $E(0)=u(0)v(0)^*$ with $u(0)$ and $v(0)$ of unit norm, the solution of \eqref{ode-E-1-sys} is given as
$E(t)=u(t)v(t)^*$, where $u$ and $v$ solve the system of differential equations (for $G=G_\eps(E)=-rs^*$)
\begin{equation}\label{ode-uv-sys}
\begin{array}{rcl}
 \dot u &=& -\tfrac \iu2 \, \Im(u^*Gv)u - (I-uu^*)Gv
\\[1mm]
 \dot v &=& -\tfrac \iu2 \, \Im(v^*G^*u)v - (I-vv^*)G^*u,
\end{array}
\end{equation}
which preserves $\|u(t)\|=\|v(t)\|=1$ for all $t$. For its numerical integration we can again use the splitting method of
Section~\ref{subsec:proto-numer}.

\medskip\noindent
{\bf Monotonicity.} 
Assuming simple eigenvalues almost everywhere along the trajectory of \eqref{ode-E-1-sys},
we again have the monotonicity property of Theorem~\ref{chap:proto}.\ref{thm:monotone-C-1},
\begin{equation}
\frac{d}{dt} \F_\eps (E(t))  \le  0.
\label{mon-sys}
\end{equation}

\medskip\noindent
{\bf Stationary points.}
\index{stationary point}
Let $E\in \cM_1$ be of unit Frobenius norm and assume that $P_E(r s^*)\ne 0$. If $E$ is a stationary point of the projected differential equation \eqref{ode-E-1-sys}, then $E$ is  already a stationary point of the differential equation \eqref{ode-E-sys}.

\subsection{Approximating the $\Hinf$-norm (outer iteration)}\label{subsec:hinf-cont}

We wish to compute $\HinfGc$, using the characterization \eqref{hinfnormdefcont}. 
We start by observing that since the spectral value set abscissa $\aleps(A,B,C,D)$
is a monotonically increasing function of $\eps$, we simply need to solve the equation
\beq \label{eq:fdef-sys}
     \aleps(A,B,C,D)=0
\eeq
for $\eps > 0$. The first step is to characterize how $\aleps$ depends on $\eps$.
\begin{theorem}[Derivative of the spectral value set abscissa] \label{thm:lambdaprime}
Let $\lambda(\eps)$ denote the rightmost point of $\SVSeps(A,B,C,D)$
for $\eps > 0$, $\eps\|D\|<1$, and assume that Assumption~\ref{assumptcont} holds for all such $\eps$.
Define $u(\eps)$ and $v(\eps)$ as right and left singular vectors with unit norm
corresponding to $\eps^{-1}$, the largest singular value of $H(\lambda(\eps))$, and applying
Theorem~\ref{thm:evecssvecs} with $\Delta(\eps)=\eps E(\eps) = \eps u(\eps)v(\eps)^*$, define $\xh(\eps)$ and $\yh(\eps)$ by \eqref{xyformulas}.
Furthermore, assume that at $\eps$, the rightmost point
$\lambda(\eps)$ is simple and unique. Then $\lambda$ is continuously differentiable
at $\eps$ and its derivative is real, with
\beq
 \frac{d}{d \eps} \aleps(A,B,C,D) = \frac{d}{d\eps}\lambda({\eps}) =
 \frac{1}{\xh(\eps)^* \yh(\eps)} > 0.
 \label{eq:deralpha}
\eeq
\end{theorem}
\begin{proof}
For the purposes of differentiation, we identify $\lambda\in\C$ with $\zeta\in\R^2$ as in the proof of Lemma~\ref{lem:firstordercont}.
The first part of Assumption \ref{assumptcont} ensures that the largest singular value of $H(\lambda)$ is differentiable with respect to $\lambda$ and that the singular vectors $v(\eps)$ and $u(\eps)$ are well defined up to multiplication of both by a unimodular scalar, and that $E(\eps)$ is not only well defined but differentiable with respect to $\eps$.  The second part ensures that $\xh(\eps)^*\yh(\eps)$ is nonzero,  using standard eigenvalue perturbation theory. As in the proof of Lemma \ref{lem:firstordercont}, observe that
\[
 \frac{1}{\eps} - \|C\left( \lambda(\eps) I - A \right)^{-1}B + D\|_2  = 0
\]
so differentiating this with respect to $\eps$ 
and using the chain rule yields
\[
     \frac{d\lambda(\eps)}{d\eps} = 
     \frac{1}{\eps^2 v^*C(\lambda(\eps)I - A)^{-2}Bu}.
\]
Furthermore, \eqref{evecip} follows (for $\lambda=\lambda(\eps)$) from \eqref{xyformulas}. Combining these
with the first-order optimality conditions for \eqref{alepsdef2} in \eqref{firstordercont} gives the result.
\qed
\end{proof}

\begin{corollary}
Make the same assumptions as in Theorem {\rm \ref{thm:lambdaprime}}, except normalize
$x(\eps)$ and $y(\eps)$ so that they fulfil \eqref{eq:scaling-sys}.  This can be seen to be 
equivalent to scaling $\xh(\eps)$ and $\yh(\eps)$ by $1/\beta(\eps)$ and $1/\gamma(\eps)$
respectively where
\beq
\label{beta-gamma-eps}
\beta(\eps)  =\frac{1-\eps u(\eps) ^*D^*v (\eps) }{u(\eps)^*b(\eps)}, \qquad 
\gamma(\eps)=\frac{1-\eps v(\eps)^*Du(\eps)}{v(\eps)^*c(\eps)}.
\eeq
Hence
\beq
 \frac{d}{d \eps} \aleps(A,B,C,D) = \frac{d}{d\eps}\lambda({\eps}) =
 \frac{1}{\beta(\eps)\gamma(\eps) \big( x(\eps)^* y(\eps)\big)} \in \R^{+}.
 \label{eq:deralphasc}
\eeq
\end{corollary}


If $A,B,C,D$ are all real, then $\SVSeps(A,B,C,D)$ is symmetric with respect to the
real axis and hence its rightmost points must either be real or part of a conjugate pair.
In the latter case, the assumption that $\lambda(\eps)$
is unique does not hold but the result still holds if there is no third rightmost point.

The derivative formula \eqref{eq:deralphasc} naturally leads to a formulation of Newton's
method for computing $\HinfGc$, similar to 
the one previously considered to compute stability radii in 
Section~\ref{subsec:Newton--bisection}.

\subsubsection*{Algorithm.}
We present the basic algorithm in concise form in Algorithm~\ref{alg_Hinf}.

\begin{algorithm}
\DontPrintSemicolon
\KwData{Matrices $A,B,C,D$, initial vectors $u_0,v_0$, $\eps_0 > 0$,
${\rm tol}$ a given positive tolerance} 
\KwResult{$\eps^J \approx \HinfGc$}
\Begin{
\For{$j=0,\dots,j_{\max}$}{
\nl Approximate numerically $\alpha_{\eps^j}(A,B,C,D)$ by integrating numerically
\eqref{ode-uv-sys} into a stationary point\;
\nl return rightmost point $\lambda^j$, and the associated left and right eigenvectors $x^j$, $y^j$ and corresponding scalars $\beta^j$, $\gamma^j$ defined as in \eqref{beta-gamma-eps},
where $b= B^*x^j$ and $c = C y^j$\;
\nl \eIf{$|\Re~\lambda^{j}| < {\rm tol}$}{
\nl Set $J=j$\;
\nl \ Return $\HinfGc \approx 1/\eps^J$\;
}
{\nl Set
\[
       \eps^{j+1}=\eps^j - \big(\Re~\lambda^j\big)\beta^j\gamma^j\big((x^j)^*y^j\big) .
\]
}
}
}
\caption{Basic algorithm to compute the $\Hinf$-norm $\HinfGc$}
\label{alg_Hinf} 
\end{algorithm}

\medskip
\index{Newton--bisection method}
Since Newton's method may not converge, it is standard practice to combine it with a bisection method that maintains an interval known to contain the root, bisecting when the Newton step is either outside the interval or does not yield a sufficient decrease in the absolute function value (in this case $|\alpha_{\eps^j}(A,B,C,D)|=|\Re~\lambda^j|$). 

\bng
We emphasize, however, that this is still an idealized algorithm because (a) there is no guarantee
that the computed stationary point of the differential equation will return the correct value of $\alpha_{\eps^j}(A,B,C,D)$, and (b) the possible breakdown of the Newton--bisection iteration, in analogy to the algorithm to approximate the stability radius, discussed in Section \ref{sec:Newton-bisection}. \bng
\bcl An appropriate choice of initial perturbations, analogous to the choices discussed in Section~\ref{subsec:init}, substantially mitigates (a). The monotone Newton--bisection method and the HEC method (see Chapter~\ref{chap:two-level} for both) are remedies to (b).
\ecl

When $n \gg \max(m,p)$ and the matrix $A$ is sparse,
this algorithm is much faster than the standard
Boyd-Bala\-krish\-nan-Bru\-insma-Stein\-buch algorithm to compute the $\Hinf$-norm.
 \bng However, the latter algorithm is globally convergent in contrast to the locally convergent behavior of the algorithm we consider here.
\eng

\subsubsection*{Numerical example.}

We consider the following example (5.2.1) from Hinrichsen and Pritchard (\cite{HinP86a}):
\begin{align} \nonumber
A  &=  { \left( \begin{array}{rrrrrr} 
0 & 1 & 0 & 0 & 0 & 0 \\
0 & 0 & 1 & 0 & 0 & 0 \\
0 & 0 & 0 & 1 & 0 & 0 \\
0 & 0 & 0 & 0 & 1 & 0 \\
0 & 0 & 0 & 0 & 0 & 1 \\
-1595.48 & -2113.96 & -1361.70 & -518.13 & -122.38 & -15.92 \end{array} \right)}, 
\\
\label{eq:ex1}
B  &=  {\left( \begin{array}{rr}
0 & 0  \\
0 & 0  \\
0 & 0  \\
0 & \frac12  \\
0 & 0  \\
1 & 0  \\
\end{array} \right)}, 
\quad C = I.
\end{align}



\begin{figure}[ht!]
\centering
\vspace{-5mm}
\includegraphics[width=0.45\textwidth]{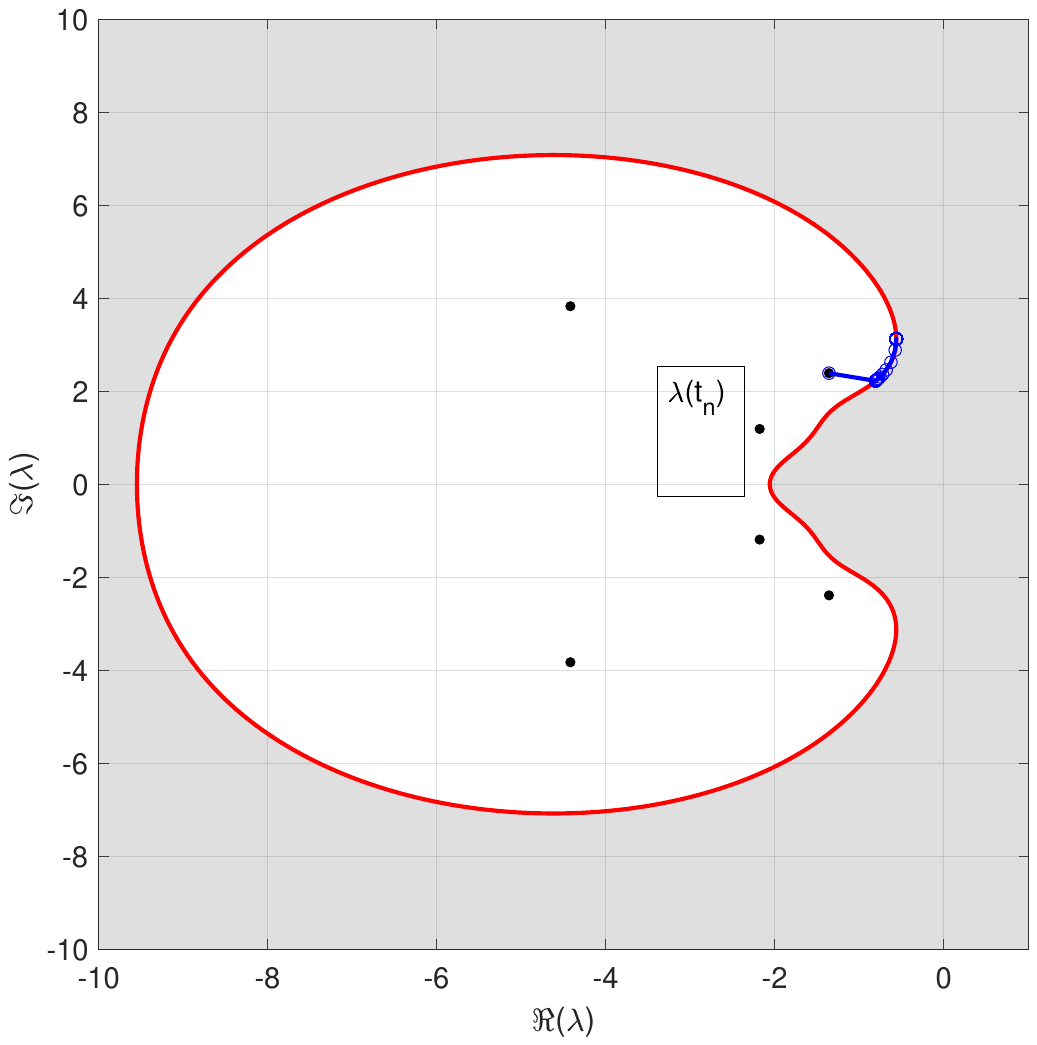} 
\includegraphics[width=0.45\textwidth]{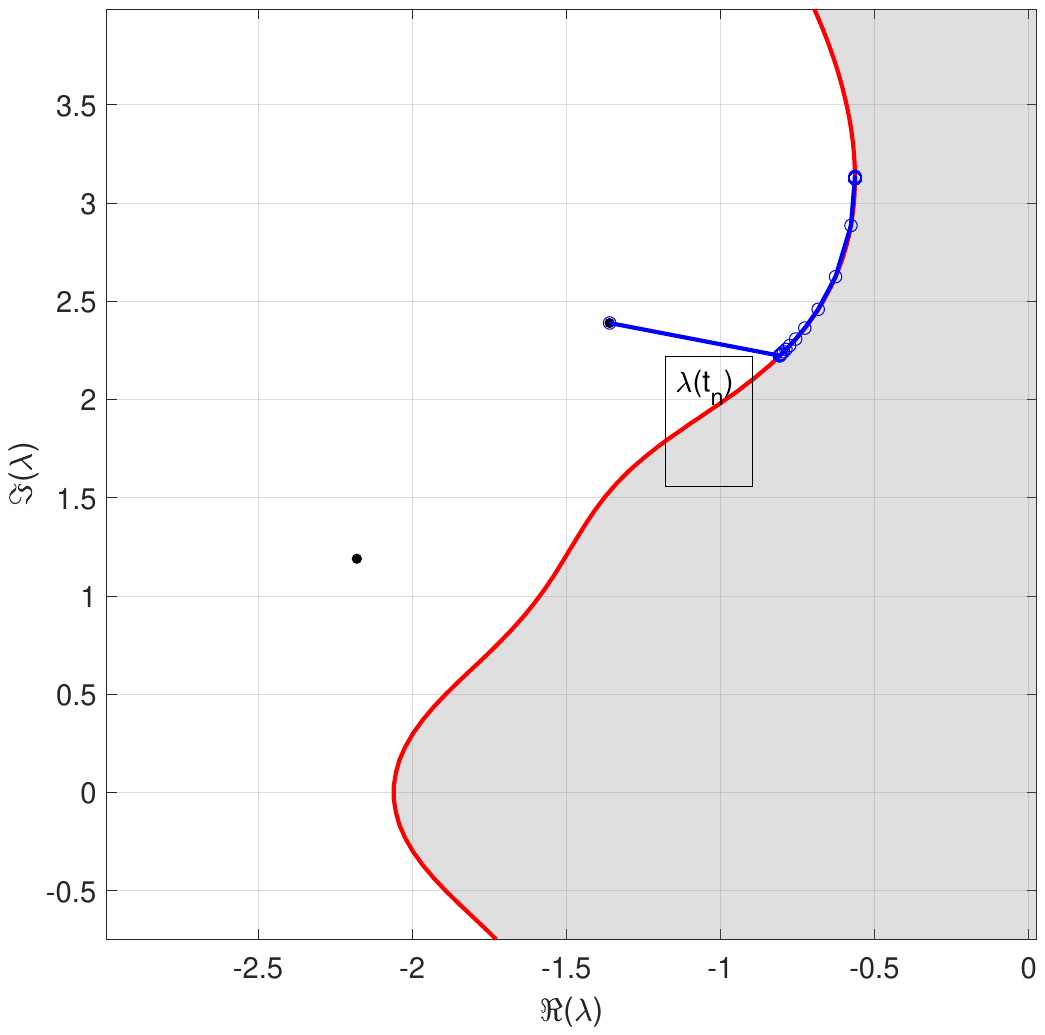}
\\[-15mm]
\includegraphics[width=0.45\textwidth]{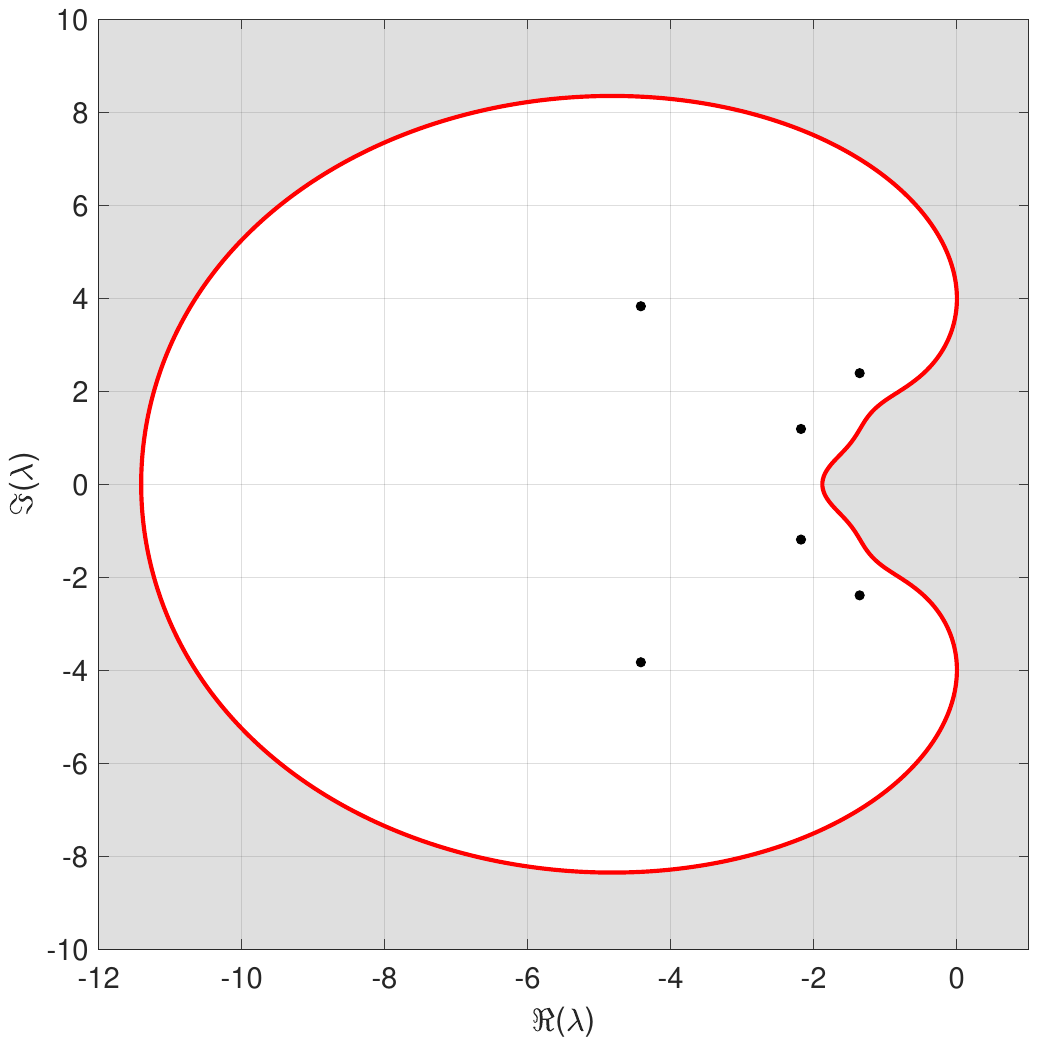} 
\vspace{-1cm}
\caption{Upper left picture: spectral value set of Example \eqref{eq:ex1} and trajectory of the splitting method ($\eps=0.1$).
Upper right picture: zoom of the spectral value set and of the trajectory of the splitting method close to the rightmost point.
Lower picture: Spectral value set  for the value $\eps=\eps_4$ in Table \ref{tab:Hinf}. The picture illustrates that the set is tangential to the imaginary axis.}
\label{fig:Hinf1}
\end{figure}

\begin{table}[hbt]
\begin{center}
\begin{tabular}{|l|l|l|l|}\hline
  $k$ & $\eps_k$ & $\alpha_{\eps^k}$ & $\alpha'_{\eps^k}$  \\
 \hline
\rule{0pt}{9pt}
\!\!\!\! 
        $0$         & $0.1$                     & $-0.563733756826769$ & $7.762584440188086$  \\
	$1$         & $0.172621916214944$       & $ 0.056254763899198$ & $8.944158197488047$  \\
	$2$         & $0.{\bf 1663}32362446206$ & $ 0.000080652047660$ & $8.917214167674935$  \\
	$3$         & $0.{\bf 166323317}912102$ & $ 0.000000000209716$ & $8.917168420440781$  \\
	$4$         & $0.{\bf 1663233178885}83$ & $-0.000000000000036$ & $$  \\
 \hline
\end{tabular}
\vspace{2mm}
\caption{Values generated by applying Algorithm \ref{alg_Hinf}.\label{tab:Hinf}}
\end{center}
\end{table}
We start by computing the spectral value abscissa by applying the splitting method to \eqref{ode-uv-sys}. The spectral value set and the computed trajectory $\lambda(t_n)$ are shown in Figure~\ref{fig:Hinf1}.
Table \ref{tab:Hinf} illustrates the behavior of Algorithm~\ref{alg_Hinf}
and in particular shows quadratic convergence.

\bcltwo
\begin{table}[hbt]
\begin{center}
\begin{tabular}{|l|l|l|l|l|}\hline
  $k$ & $\eps_k$ & $\alpha_{\eps^k}$ & eigs & accept/reject \\
 \hline
\rule{0pt}{9pt}
\!\!\!\! 
    $0$  &   $0.1$                        & $5.65382\,10^{-1}$  &        $20$    & A \\
    $1$  &   $\mathit{0.16229595981186}$  & $\mathit{-8.08079\,10^{-3}}$  & $26$ & R \\
    $2$  &   $0.13361479799059289$        & $2.87473\,10^{-1}$     &  $46$  & A \\    
    $3$  &   $0.{\bf 16}5128673019468$    & $1.06493\,10^{-2}$    &  $26$  & A \\
    $4$  &   $0.{\bf 166}274664829137$    & $4.33842\,10^{-4}$  &  $24$  & A  \\
    $5$  &   $0.{\bf 16632}1327917072$    & $1.77449\,10^{-5}$  &   $24$    & A \\ 
    $6$  &   $0.{\bf 166323}236465875$    & $7.26060\,10^{-7}$  &   $20$    & A \\
    $7$  &   $0.{\bf 16632331}4556989$    & $2.97082\,10^{-9}$  &   $20$    & A \\
    $8$  &   $0.{\bf 166323317888}13$     & $3.44159 \,10^{-13}$ & $4$ & A  \\
\hline
\end{tabular}
\vspace{2mm}
\caption{Values generated by applying the adapted monotone Newton-bisection Algorithm \ref{alg:SRC} of Section \ref{subsec:Newton-bisection-monotone}.\label{tab:Hinf2}}
\end{center}
\end{table}
\ecltwo

\bcltwo
\begin{table}[hbt]
\begin{center}
\begin{tabular}{|l|l|l|l|}
\hline
$k$ & $\eps_k$ & $\phi(\eps_k)$ & $\#$ iters (E+C) \\
\hline
\rule{0pt}{9pt}
$0$ &
$1.279372505046394$ &
$9.66901$ &
-- \\

$1$ &
$0.173574506144666$ &
$6.47765\,10^{-2}$ &
$14+6$ \\

$2$ &
$0.166417083318858$ &
$8.36144\,10^{-4}$ &
$12+3$ \\

$3$ &
$0.166323335133081$ &
$1.53772\,10^{-7}$ &
$9+3$ \\

$4$ &
$0.166323317888565$ &
$8.65018\,10^{-15}$ &
$3+4$ \\
\hline
\end{tabular}
\vspace{2mm}
\caption{Computation of the inverse of the $\Hinf$-norm for the
Hinrichsen and Pritchard example by the HEC algorithm:
computed values $\eps_k$, $\phi(\eps_k)= -\Re\,\lambda(A+\eps_k E_k)$
and the numbers of iterations performed in the expansion (SVSAR)
and contraction (Newton--bisection) phases.}
\label{tab:HPHEC}
\end{center}
\end{table}

Table \ref{tab:HPHEC} has been obtained by applying Tim Mitchell's code Rostapack which implements the HEC algorithm (http://www.timmitchell.com/software/ROSTAPACK/). \vskip 7mm

The HEC iteration converges to
\[
\eps_\ast = 0.166323317888564 .
\]

Finally, we successfully experimented with both the convergent Newton-bisection method and the HEC algorithm on problems CM3 and CM4 in COMPleib, for which the failure of the Newton-bisection method, originally proposed in Guglielmi, G\"urb\"uzbalaban and Overton (\cite{GugGO13}), provided wrong results. 
\ecltwo

\section{$\Hinf$-robustness of controllability}

\label{sec:uncon}
In this section we ask for the robustness of controllability of a given linear time-invariant system \eqref{lti}. We consider the operator nearness problem of finding the distance of a given controllable system to the nearest uncontrollable system, where the distance is taken as the $\Hinf$-norm of the perturbation to the transfer function. This is a natural metric that measures the change in the input-output behaviour due to the perturbation of the system matrices.
We propose and study a two-level algorithm that extends the algorithm of the preceding section.

\subsection{$\Hinf$-nearest uncontrollable system}
\index{linear time-invariant system!controllable}
The linear time-invariant system \eqref{lti} is {\it controllable} if and only if  the $n\times (n+p)$ matrix
\beq
\label{controllable}
\text{$(A-\lambda I, B)\:$ has full row rank for all $\lambda\in \C$.}
\eeq
For this,  it obviously suffices that the condition holds for the eigenvalues of $A$. Hence, a system is uncontrollable if and only if there exists an eigenvalue $\lambda_A$ of $A$ with corresponding left eigenvector $x_A$ such that
$$
x_A^* B=0.
$$

 We consider the problem of finding the distance of a given controllable system \eqref{lti} to the set of uncontrollable systems with the same matrices $A$, $C$ and $D$, but with perturbed input--state matrix
$
\widetilde B= B+ \Delta B.
$

Here an important question arises: which distance? We might minimize the Frobenius norm of $\Delta B$ under the condition that the perturbed system becomes uncontrollable (which would yield $\Delta B = -x_A x_A^*B$ for a left eigenvector $x_A$ to one of the eigenvalues of $A$), but this choice of a distance is not invariant under similarity transformations of $A$,
\beq\label{A-sim}
A \to V^{-1}AV, \quad B\to V^{-1}B, \quad C \to CV,
\eeq
which leave the matrix transfer function $H(\lambda)=C(\lambda I-A)^{-1}B + D$ invariant and hence also the input-output map
$y=H(\partial_t)u$. 
We therefore measure the distance by the $L^2$ operator norm of the difference of the perturbed and unperturbed input-output operators, $\Delta H(\partial_t)= \widetilde H(\partial_t)-H(\partial_t)$. This $L^2$ operator norm is the $\Hinf$-norm of the perturbation $\Delta H(\lambda)=C(\lambda I-A)^{-1}\Delta B$ to the matrix transfer function. Considering {\it real} system matrices $(A,B,C,D)$ and {\it real} perturbations $\Delta B$, we arrive at the following operator nearness problem.

\medskip\noindent
{\bf Problem.} {\it Given a controllable system \eqref{lti}, find a perturbation $\Delta B\in \R^{n,p}$ such that the perturbed input--state matrix $B+\Delta B$ yields an uncontrollable system and such that the perturbation to the transfer function, $\Delta H(\lambda)=C(\lambda I-A)^{-1}\Delta B\/$ for $\Re\,\lambda\ge 0$, is of minimal $\Hinf$-norm.}

\medskip\noindent 



Let $\lambda_A$ be an eigenvalue of $A$ with normalized left eigenvector $x_A$. We consider perturbations $R$ to $B$ that satisfy the uncontrollability condition
$x_A^*(B+R)=0$, i.e.,
\beq
\label{R-cond}
x_A^* R = - x_A^*B.
\eeq
In the following we choose some left eigenvector $x_A$ of $A$ and denote the corresponding admissible affine space by
$$
\mathcal{R} = \{ R\in\R^{n,p}\,:\, R \text{ satisfies \eqref{R-cond} } \}. 
$$
Our proposed algorithm is based on the following characterization.

\index{distance to uncontrollability}
\begin{theorem} [$\Hinf$-distance to uncontrollability] \label{thm:uncontrollability}
With the spectral value set abscissa as in \eqref{alepsdef},
$$
\alpha_\eps(A,R,C,0)=\max_{\|E\|_F=1} \Re\,\lambda(A+\eps REC),
$$
where $\lambda(M)$ is a rightmost eigenvalue of $M$, 
let
\begin{align} \nonumber
&R(\eps) = \arg \min_{R\in \mathcal{R}} \alpha_{\eps}(A,R,C,0) \quad\text{ and }
\\
&\oeps = \inf\{ \eps>0\,:\, \alpha_\eps(A,R(\eps),C,0)=0 \}.
\label{oeps-control}
\end{align}
Then, 
$
\Delta B=R(\oeps)
$ 
is a perturbation to the input--state matrix $B$ such that the perturbed system $(A,B+\Delta B,C,D)$ is not controllable and 
has minimal $\Hinf$-distance between the transfer functions of the unperturbed system and uncontrollable perturbed systems with $\Delta B \in \mathcal{R}\,$. This $\Hinf$-distance equals $1/\oeps$.
\end{theorem}

\begin{proof} For $R\in \mathcal{R}$, we denote the difference of the transfer functions corresponding to the systems
$(A,B+R,C,D)$ and $(A,B,C,D)$ by
$$
\Delta H_R(\lambda)=C(\lambda I-A)^{-1}R.
$$
By Theorem~\ref{thm:Hinf-oeps}, we have 
$$
\| \Delta H_R \|_\infty = \frac1{\eps_R} \quad \text{with} \quad
\eps_R= \min\{\eps>0\,:\, \alpha_\eps(A,R,C,0)=0\} 
$$
and therefore
$$
\min_{R\in \mathcal{R}} \| \Delta H_R \|_\infty = \frac1{\eps_{\max}}
\quad \text{with} \quad
   \eps_{\max} = \max_{R\in\mathcal{R}} \eps_R.
$$
We now prove that $\eps_{\max}$  equals $\oeps$ of \eqref{oeps-control}. Since the spectral value set abscissa depends monotonically on $\eps$ and we have
$\alpha_{\oeps}(A,R,C,0)\ge 0$ for all $R\in \mathcal{R}$ and $\alpha_{\oeps}(A,R,C,0)= 0$ for some $R\in \mathcal{R}$, we obtain
$\oeps \ge \eps_R$ for all $R\in \mathcal{R}$ and $\oeps = \eps_R$ for some $R\in \mathcal{R}$, and hence
$\oeps=\eps_{\max}$.
This yields the stated result.
\qed    
\end{proof}





\subsection{Multi-level iteration}
\index{three-level iteration}
Theorem~\ref{thm:uncontrollability} motivates us to consider the following schematic multi-level iteration, for which we will provide details in the following subsections.



\begin{itemize}
\item {\bf Innermost iteration:\/}
For a given $\eps>0$ and $R\in\mathcal{R}$, compute the spectral value set abscissa $\alpha_\eps(A,R,C,0)$,
e.g., by the rank-1 gradient flow algorithm of the preceding section. This algorithm also yields a rank-1 matrix $E_\eps(R)\in\C^{n,n}$ of Frobenius norm~1 such that
\begin{equation} \label{alpha-uncon}
\alpha_\eps(A,R,C,0) = \max_{E\in\C^{n,n},\,\|E\|_F=1} \Re\,\lambda(A+\eps REC)= \Re\,\lambda(A+\eps RE_\eps(R)C),
\end{equation}
where $\lambda(M)$ is a rightmost eigenvalue of a matrix $M$.
\item {\bf Inner iteration:\/}
For a given $\eps>0$, compute
\begin{equation}
\label{eig-opt-uncon}
R(\eps)=  \arg \min_{R \in\mathcal{R}} \alpha_\eps(A,R,C,0),
\end{equation}
e.g. by a gradient flow as described below.
\item {\bf Outer iteration:\/}
Compute the zero $\oeps$ of the monotonically growing function
\begin{equation}\label{phi-uncon}
\phi(\eps)= \alpha_\eps(A,R(\eps),C,0).
\end{equation}
We use a Newton--bisection method for this scalar equation, similar to the preceding section.
\end{itemize}
The zero $\oeps$ depends on the chosen left eigenvector $x_A$ of $A$, which appears in \eqref{R-cond} and hence in the definition of $\mathcal{R}\,$. While in theory we would have to compute $\oeps=\oeps(x_A)$ for each eigenvector $x_A$ and then choose the maximal value among the $\oeps(x_A)$, we can expect that in practice it suffices to choose just a few eigenvectors $x_A$ for which $x_A^*B$ has smallest norms. (Note that $x_A^*B=(x_A^*V)(V^{-1}B)$ is invariant under transformations \eqref{A-sim}.)

For the maximal $\oeps$ we have that $\Delta B=R(\oeps)$ is a perturbation to the input--state matrix $B$ such that the perturbed system $(A,B+\Delta B,C,D)$ is not controllable and 
has the minimal $\Hinf$-distance $1/\oeps$ between the transfer functions of the unperturbed and perturbed systems.



\subsection{Inner iteration: constrained gradient flow}

We start with the free gradient of the functional in the inner iteration.

\begin{lemma}[Derivative of the spectral value set abscissa] \label{lem:gradient-uncon}
Let $R(t)\in \R^{n,p}$, for real $t$ near $t_0$, be a continuously differentiable path of matrices, with the derivative denoted by $\dot R(t)$. Let $E(t)=E_\eps(R(t))\in\C^{n,n}$ of Frobenius norm 1 be such that
$$
\alpha_\eps(A,R(t),C,0) = \Re\,\lambda(A+\eps R(t) E(t) C)
$$
and the rightmost eigenvalue on the right-hand side is simple, with left and right eigenvectors $x(t)$ and $y(t)$,
of norm 1 and with positive inner product. Let $\kappa(t)=1/(x(t)^*y(t))$. Assume that $x(t)^*R(t)\ne 0$ and $Cy(t)\ne 0$.
Then we have, omitting the argument $t$ on the right-hand side,
$$
\frac1{\eps\kappa(t)}\,\frac{d}{dt}\,\alpha_\eps(A,R(t),C,0) =
\langle G_R, \dot R \rangle \quad \text{ with }\quad G_R=\Re\bigl(x(ECy)^*\bigr).
$$
\end{lemma}

\begin{proof}
    We have
    \begin{align*}
    \frac1{\eps\kappa(t)}\,\frac{d}{dt}\, \Re\,\lambda(A+\eps R(t) E(t) C) &=
    \Re \, x^*(\dot REC + R\dot E C)y 
    \\
    &=\langle \Re\,x(ECy)^*, \dot R \rangle + \Re\,\langle (R^Tx)(Cy)^*, \dot E \rangle
    \\
    &= \langle G_R, \dot R \rangle + \Re\,\langle G_E, \dot E \rangle
    \end{align*}
    with $G_R=\Re\bigl(x(ECy)^*\bigr)$ and $G_E=(R^Tx)(Cy)^*$.
    We show that the last term vanishes.
    Since $E(t)$ is a stationary point of the functional $\F_\eps(E)=\Re\,\lambda(A+\eps R E C)$ (with $R=R(t)$ for a fixed $t$) under the constraint of Frobenius norm 1, and since $G_E$ is nonzero by assumption, we find once again, as in Chapter~\ref{chap:proto}, that $E=\mu G_E$ for some nonzero real factor $\mu$. We then note, omitting the argument $t$ on the right-hand side,
    $$
    0 = \frac12\, \frac{d}{dt}\, \|E(t)\|^2 = \Re \langle E, \dot E \rangle = \mu\, \Re \langle G_E,\dot E \rangle.
    $$
    This yields the stated result.
\qed
\end{proof}

\subsubsection*{Constrained gradient flow.} Choosing $\dot R$ such that it points in the direction of steepest admissible descent under the linear constraint \eqref{R-cond}, we arrive at the projected gradient system  in which we take the constraint $x_A^*\dot R=0$ into account:
\beq
\label{grad-sys-uncon}
\dot R = -(I-P_{x_A})\, G_R \quad\text{ with }\quad G_R = \Re\bigl(x(ECy)^*\bigr),
\eeq
where $P_{x_A}$ is the orthogonal projection onto the range of $(\Re\,x_A,\Im\,x_A)$ and $E=E_\eps(R)$ is the perturbation that yields the spectral value set abscissa (computed in the innermost iteration as in Section~\ref{sec:Hinf}).  






\subsection{Outer iteration}
\index{Newton--bisection method}
For the outer iteration we again use a Newton method, which is justified under additional conditions, 
\bng
and we couple it with bisection, as we have done previously. A breakdown of the standard Newton--bisection iteration might still occur. As in Chapter \ref{chap:two-level}, remedies are given by the monotone Newton--bisection method and the HEC method.
\eng
Under the following assumption, the real function $\phi$ of \eqref{phi-uncon} is differentiable and we will give a simple formula for its derivative.

\begin{assumption} \label{ass:E-eps-uncon}
For $\eps$ in a neighbourhood of $\oeps$, 
we assume the following for the optimizer $R(\eps)$ of the inner iteration and $E(\eps)=E_\eps(R(\eps))$ of the innermost iteration:
\begin{itemize}
\item The eigenvalue $\lambda(\eps)=\lambda( A+\eps R(\eps)E(\eps)C)$ is a simple eigenvalue, with normalized left and right eigenvectors denoted by $x(\eps)$ and $y(\eps)$.
\item The map $\eps \mapsto (E(\eps),R(\eps))$ is continuously differentiable.
\item The partial gradient $G_E(\eps)=(R(\eps)^Tx(\eps))(Cy(\eps))^*$ is nonzero.
\end{itemize}
\end{assumption}

Theorem~\ref{chap:two-level}.\ref{thm:phi-derivative} extends to the present situation.
We again denote the eigenvalue condition number by
$
\kappa(\eps) =  1 /(x(\eps)^*y(\eps))>0.
$
\\
\begin{lemma}[Derivative for the Newton iteration] 
\label{lem:phi-derivative-uncon}
Under Assumption~\ref{ass:E-eps-uncon}, the function $\phi(\eps)=\alpha_\eps(A,R(\eps),C,0)$ is continuously differentiable in a neighbourhood of $\oeps$ and its derivative is given as
\begin{equation} \label{eq:dereps-uncon}
\phi'(\eps)  = 
 \kappa(\eps) \,\| G_E(\eps) \|_F > 0.
\end{equation}
\end{lemma}
\begin{proof}
By Lemma~\ref{lem:gradient-uncon} we obtain, with $G_R(\eps) = x(\eps)(E(\eps)Cy(\eps))^*$ and indicating by $'$ differentiation w.r.t.~$\eps$, 
\begin{align*} 
&\frac{1}{\kappa(\eps)} \,\frac{d}{d \eps} \alpha_\eps(A,R(\eps),C,0) 
\\
&= 
\frac{1}{\kappa(\eps)} \,\frac{d}{d \eps} \alpha_\eps(A,R(\eta),C,0)\Big|_{\eta=\eps} +
\frac{1}{\kappa(\eps)} \,\frac{d}{d \eps} \alpha_\eta(A,R(\eps),C,0)\Big|_{\eta=\eps} 
\\
&= \Re\,\langle G_E(\eps), E(\eps) \rangle + \Re\,\langle G_R(\eps), R'(\eps) \rangle.
\end{align*}
As in the proof of Theorem~\ref{chap:two-level}.\ref{thm:phi-derivative} we find 
$$
\Re \bigl\langle  G_E(\eps),  E(\eps) \bigr\rangle =
\| G_E(\eps) \|_F.
$$
In the stationary point $R(\eps)$ of \eqref{grad-sys-uncon} we have
$(I-P_{x_A})G_R(\eps)=0$, that is,  
$G_R(\eps)=P_{x_A}G_R(\eps)$. On the other hand, because of \eqref{R-cond} we have $P_{x_A} R'(\eps)=0$. Hence,
$$
\bigl\langle G_R(\eps),R'(\eps))\bigr\rangle 
=
\bigl\langle P_{x_A}G_R(\eps),R'(\eps))\bigr\rangle
= \bigl\langle G_R(\eps),P_{x_A} R'(\eps))\bigr\rangle =0,
$$
which yields the result.
\qed 
\end{proof}



\section{Nearest passive system}
\label{sec:pass}
\index{linear time-invariant system!passive}
Consider the linear time-invariant system \eqref{lti}, which we here assume quadratic ($p=m$). The system is called {\it passive} if every input $u\in L^2(\R_+,\R^m)$ and its corresponding output $y$ satisfy the relation
\begin{equation}\label{passive}
\int_0^T y(t)^\top u(t)\, dt \ge 0 \quad\ \text{ for all }\: T>0.
\end{equation} 
Passivity is a fundamental property in control theory. When a given system is not passive, it is often required to enforce passivity by modifying it such that it becomes passive, yet remains `near' the given system. 
One approach, 
as discussed by Grivet-Talocia \& Gustavsen (\cite{GriG15}) and to be adopted here, is to perturb only the state--output matrix $C\in\R^{m,n}$ to $C+\Delta C$ and to minimize a suitable norm of the perturbation such that the perturbed system is passive. 

A favoured choice in the literature is to minimize the Frobenius norm of $\Delta C\, L$, where $L\in \R^{n,n}$ is a Cholesky factor of the controllability Gramian $G_c=LL^T$, which is the unique solution of the Lyapunov equation $AG_c +G_cA^\top= BB^\top$.
This Gramian is symmetric positive definite for a controllable system, as will be assumed in the following. The squared norm $\langle \Delta C \,G_c, \Delta C \rangle = \| \Delta C \, L \|_F^2$ is invariant under transformations \eqref{A-sim} of the system matrices, as is readily checked. Based on the characterization of passivity via Hamiltonian matrices, we present a two-level algorithm for this minimization problem in the Frobenius norm.

In a different direction that appears not to have been addressed in the literature before, we present in Subsection~\ref{subsec:pass-Hinfty} an algorithm for the problem of enforcing passivity by perturbing $C$ in such a way that the $\Hinf$-norm of the difference between the  transfer functions of the passive perturbed system and of the original system is minimized, similar to the approach in Section~\ref{sec:uncon}.

\subsection{Hamiltonian matrix related to passivity}
\index{Hamiltonian matrix}
By the Plancherel formula and the relation \eqref{Hu} of the matrix transfer function $H(\lambda)$, we have (with $u$ extended by 0 outside [0,T])
\begin{align*}
\int_0^T y(t)^\top u(t)\, dt &= \Re \int_{\R} \mathcal{L}y(\iu\omega)^* \mathcal{L}u(\iu\omega)\,d\omega =
\Re \int_{\R} \mathcal{L} \bigl(H(\iu\omega)u(\iu\omega)\bigr))^* \mathcal{L}u(\iu\omega)\,d\omega 
\\
&=
\frac12 \int_{\R} u(\iu\omega)^* 
\bigl(H(\iu\omega) + H(\iu\omega)^*\bigr)u(\iu\omega)\,d\omega,
\end{align*}
and so we find that passivity \eqref{passive} is equivalent to the property of {\it positive realness} of the transfer function:
\index{transfer function!positive real}
\begin{equation}\label{pr}
H(\iu\omega)+H(\iu\omega)^* \text{ is positive semi-definite for all $\omega\in\R$}.
\end{equation}
Since $H(\iu\omega)\to D$ as $\omega\to\infty$, a necessary condition for \eqref{pr} is that $D+D^\top$ is positive semi-definite.

{\it Strict positive realness} is defined in the same way, with `positive definite' instead of `positive semi-definite' for all $\omega\in \R\cup\{\infty\}$. A necessary condition is now that
$D+D^\top$ is positive definite.

The following remarkable result can be found in Chapter 9 of the book by Grivet-Talocia \& Gustavsen (\cite{GriG15}). It goes back to Boyd, Balakrishnan \& Kabamba (\cite{BoyBK89}) and in its conceptual origins further back to Byers (\cite{Bye88}). It characterizes strict positive realness, which is a condition on the transfer function on the whole imaginary axis, by the location of the
eigenvalues of a single Hamiltonian matrix built from the system matrices.

\begin{theorem}
[Passivity and eigenvalues of a Hamiltonian matrix]
\label{thm:passive-ham}
The matrix transfer function $H(\cdot)$ is strictly positive real if and only if $D+D^\top$ is positive definite and the Hamiltonian matrix
\beq
\label{K-pass}
K=  \begin{pmatrix}
A & 0\\ 0 & -A^T
\end{pmatrix}
-
 \begin{pmatrix}
B \\-C^\top
\end{pmatrix}
(D+D^\top)^{-1}
 \begin{pmatrix}
C^\top \\ B
\end{pmatrix}^\top\ 
 \in \R^{2n\times 2n}
\eeq
has no eigenvalues on the imaginary axis.
\end{theorem}

In the following we will consider the Hamiltonian matrix also for systems with perturbed state--output matrices $C+\Delta C$. To indicate the dependence on $C$, we write $K(C)$ for $K$ of \eqref{K-pass}.

\subsection{Two-level iteration}
\index{two-level iteration}
\index{Hamiltonian matrix nearness problem}
Theorem~\ref{thm:passive-ham} leads us to the following matrix nearness problem, which is reminiscent of Problem B in Section~\ref{sec:Hamilton}. Here, the matrix $L\in\R^{n,n}$ is usually chosen as a Cholesky factor of the controllability Gramian, as mentioned above.

\medskip\noindent
{\bf Problem (F).} 
{\it Given a system  \eqref{lti} with positive definite matrix $D+D^\top$
for which the Hamiltonian matrix $K(C)$ has some purely imaginary eigenvalues, 
and given $\delta>0$, compute a perturbed state--output matrix   $C + \Delta C$ with minimal $\|\Delta C \,L\|_F$ such that all eigenvalues of $K(C + \Delta C)$ have a real part of absolute value at least~$\delta$.}

\medskip\noindent
The proposed algorithm is a two-level iterative method similar to the second method of Section~\ref{subsec:ham-B}.

\begin{itemize}
\item {\bf Inner iteration:\/}
For a given $\eps>0$, we use a (low-rank) gradient system to solve the eigenvalue optimization problem, over $E\in \R^{p\times n}$ with $\| E \|_F =1$, 
\begin{equation}\label{eig-opt-pass}
E(\eps) =  \arg \max_{\|E\|_F=1}  \ \Re\, \lambda\bigl(K(C+ \eps E L^{-1})\bigr),
\end{equation}
where
$\lambda(K)$ is an eigenvalue of minimal nonnegative real part (chosen with the largest nonnegative imaginary part) of a real Hamiltonian matrix $K$.

\item {\bf Outer iteration:\/} We compute the  smallest $\eps$ such that
\beq
\label{phi-pass}
\phi(\eps):=\Re\, \lambda\bigl(K(C+ \eps E(\eps) L^{-1})\bigr)= \delta
\eeq 
for the given small threshold $\delta>0$. This uses a   Newton--bisection method and, for very small $\delta$, the asymptotic square-root behavior $ \phi(\eps)\sim \sqrt{\eps-\oeps}$ as the eigenvalue tends to the imaginary axis
at perturbation size $\oeps$; see Theorem~\ref{chap:mnp-mix}.\ref{thm:sqrt}.
\end{itemize}

\subsection{Norm- and rank-constrained gradient flows}\label{subsec:gradientflows-pass}

In this subsection we show how to deal with the inner iteration, once again following and adapting the by now well-trodden path of Chapter~\ref{chap:proto}, here in the real version of Section~\ref{sec:real}.

In the resulting algorithm we do not move eigenvalues of Hamiltonian matrices on the imaginary axis. Instead, like in the second algorithm of Section~\ref{subsec:ham-B}, we work with Hamiltonian matrices all of whose eigenvalues are off the imaginary axis, corresponding to perturbed
matrices $C$ that yield a passive system.
We move
eigenvalues with smallest positive real part toward the imaginary axis, starting from a
non-optimal passive perturbation of the original, non-passive system. This starting perturbation can come from a computationally inexpensive but non-optimal passivity enforcement algorithm. 

\medskip\noindent
\index{gradient!free}
{\bf Free gradient.}
The following lemma will allow us to compute the steepest descent direction of the functional $\F_\eps(E)=-\Re\, \lambda\bigl(K(C+ \eps E L^{-1})\bigr)$.

\begin{lemma}[Real gradient] \label{lem:gradient-pass}
Let $E(t)\in \R^{m,n}$, for real $t$ near $t_0$, be a continuously differentiable path of matrices, with the derivative denoted by $\dot E(t)$.
Assume that $\lambda(t)$ is a simple eigenvalue of  $K(C+ \eps E(t) L^{-1})\bigr)$ depending continuously on $t$, with associated eigenvectors
$\y(t)$ and $\x(t)$ normalized by \eqref{eq:scaling-sys}, and let 
$
\kappa(t) = 1/(\y(t)^* \x(t))>0.
$
\\
Then, the derivative of $\F_\eps(E(t))=\Re\, \lambda\bigl(K(C+ \eps E(t) L^{-1})\bigr)$ 
is given by
\begin{equation} \label{eq:deriv-pass-2}
\frac1{ \eps \kappa(t) } \,\frac{d}{dt} \F_\eps(E(t)) = \,\bigl\langle  G_\eps(E(t)),  \dot E(t) \bigr\rangle
\end{equation}
with the (rescaled) real gradient
\begin{equation}
G_\eps(E)= K'(C+\eps E L^{-1})^\ast[\Re(x y^*)] L^{-\top}\in \R^{m,n},
\label{eq:grad-pass}
\end{equation}
where $K'(C)^\ast:\R^{2n,2n}\to\R^{m,n}$ is the adjoint of the derivative
$K'(C):\R^{m,n}\to \R^{2n,2n}$ defined by $\langle K'(C)^\ast[W], Z \rangle = \langle W, K'(C)[Z]\rangle$ for all  $W\in \R^{2n,2n}$ and $Z\in \R^{m,n}$.
\end{lemma}

\begin{proof} By the derivative formula for simple eigenvalues and the chain rule, we have (omitting the ubiquitous argument $t$ after the first equality sign)
\begin{align*}
&\frac d{dt} \, \Re\, \lambda\bigl(K(C+ \eps E(t) L^{-1})\bigr)
= \kappa \,\Re \Bigl( x^*  K'(C+\eps E L^{-1})[\eps \dot E L^{-1}]y \Bigr)
\\
&= \kappa \bigl\langle \Re(xy^*),K'(C+\eps E L^{-1})[\eps \dot E L^{-1}] \bigr\rangle = 
\kappa \eps \bigl\langle K'(C+\eps E L^{-1})^\ast[\Re(xy^*)]L^{-\top},\dot E \bigr\rangle,
\end{align*}
which yields the stated result.
\qed
\end{proof}

An explicit formula for the adjoint of the derivative of $K$ is given next. 

\begin{lemma}[Adjoint of the derivative]
\label{lem:K-der-adj}
For $W=\begin{pmatrix}
W_{11} & W_{12} \\ W_{21} & W_{22} 
\end{pmatrix}\in \R^{2n\times 2n}$ partitioned according to the $n\times n$ blocks of~$K(C)$, we have  (with $T=D+D^\top$ for short)
\begin{equation}
\label{K-der-adj}
K'(C)^\ast[W] =   -T^{-1} B^\top (W_{11} - W_{22}^\top) +  T^{-1} C\,(W_{21}+W_{21}^\top).
\end{equation}
\end{lemma}

\begin{remark}\label{rem:rank-pass}
A noteworthy consequence of 
\eqref{K-der-adj} is that 
$G_\eps(E)$ of \eqref{eq:grad-pass} has rank at most~8.
(The rank is at most 4 for real eigenvalues.)
\end{remark}

\begin{proof} For any path $C(t)$ we have $\dot K(t) = \frac d{dt} K(C(t)) = K'(C(t))[\dot C(t)]$ and hence
$$
\langle  K'(C)^\ast [W], \dot C \rangle = \langle W, K'(C) [\dot{C}] \rangle = \langle W, \dot{K} \rangle.
$$
Differentiation of \eqref{K-pass} yields
\begin{align*}
\dot K  =&\  
-
 \begin{pmatrix}
0 \\-\dot C^\top
\end{pmatrix}
T^{-1}
 \begin{pmatrix}
C^\top \\ B
\end{pmatrix}^\top
-
 \begin{pmatrix}
B \\-C^\top
\end{pmatrix}
T^{-1}
 \begin{pmatrix}
\dot C^\top \\ 0
\end{pmatrix}^\top.
\end{align*}
We note that
\begin{align*}
&\left\langle
W,
\begin{pmatrix}
0 \\-\dot C^\top
\end{pmatrix}
T^{-1}
 \begin{pmatrix}
C^\top \\ B
\end{pmatrix}^\top  \right\rangle=  
\left\langle
W  \begin{pmatrix}
C^\top \\ B
\end{pmatrix} T^{-1},
 \begin{pmatrix}
0 \\-\dot C^\top 
\end{pmatrix}  \right\rangle 
\\
&\quad = - \bigl\langle   (W_{21} C^\top + W_{22} B) T^{-1}, \dot C^\top \bigr\rangle = -
\bigl\langle   ( (W_{21} C^\top + W_{22} B) T^{-1})^\top , \dot C \bigr\rangle
\\[2mm]
&  \left\langle W,
\begin{pmatrix}
B \\-C^\top
\end{pmatrix}
T^{-1}
 \begin{pmatrix}
\dot C^\top \\ 0
\end{pmatrix}^\top \right\rangle = \langle  T^{-1} (B^\top,-C)W, (\dot C, 0) \rangle
\\
&\quad = \bigl\langle  T^{-1} (B^\top W_{11} - C W_{21}), \dot C \bigr\rangle 
\end{align*}
so that finally  
$$ 
\langle  K'(C)^\ast [W], \dot C \rangle =
\langle W,\dot{K} \rangle  =  \langle  -T^{-1} B^\top (W_{11} - W_{22}^\top) +  T^{-1} C\,2\,\mathrm{Sym}(W_{21}), \dot C \rangle
$$
for all $\dot C\in \R^{m,n}$. This yields $K'(C)^\ast [W]$ as stated.
\qed
\end{proof}

\medskip\noindent
{\bf Norm-constrained gradient flow.} 
\index{gradient flow!norm-constrained}
As in Section~\ref{sec:real},
we consider the  projected gradient flow on the manifold of  $m\times n$ real matrices of unit Frobenius norm:
\beq\label{ode-E-pass}
\dot E = -G_\eps(E) + \left\langle G_\eps(E),  E \right\rangle E,
\eeq
where $G_\eps(E)$ is defined by \eqref{eq:grad-pass} via an eigentriple
$(\lambda$, $x$, $y)$ of 
the Hamiltonian matrix $K(C+ \eps E L^{-1})$ with $\lambda$ the target eigenvalue of minimal nonnegative real part (and among these, the one with largest nonnegative imaginary part). The Frobenius norm 1 is conserved along trajectories. 

We now follow closely the programme of Section~\ref{sec:real}: 
We again have the monotonic decay of $\F_\eps(E(t))$ as in
(\ref{chap:proto}.\ref{eq:pos-real}), and
 the characterization of stationary points as given in Theorem~\ref{chap:proto}.\ref{thm:stat} also extends: Let
$E\in\C^{m,n}$ with $\| E\|_F=1$ be such that the target eigenvalue $\lambda$ of $K(C+\eps E L^{-1})$ is simple
and $G_\eps(E)\ne 0$. Then, $E$ is a stationary point of the differential equation \eqref{ode-E-pass}
if and only if $E$ is a real multiple of $G_\eps(E)$.
Together with Remark~\ref{rem:rank-pass}, this implies the following.



\begin{corollary}[Rank of optimizers] \label{cor:rank-8-pass}
If $E$ is an optimizer of the eigenvalue optimization problem \eqref{eig-opt-pass} and if $G_\eps(E) \ne 0$, 
then $E$ is of rank at most $8$.
\end{corollary}
\index{optimizer!rank-8}

As in Section~\ref{sec:real}, 
this motivates us to constrain the differential equation \eqref{ode-E-pass} to a manifold of real low-rank matrices, which turns out to be computationally favourable for large systems.

\medskip\noindent
{\bf Rank-8 constrained gradient flow.}
\index{gradient flow!rank-8 constrained}
In the same way as in Section~\ref{subsec:rank-r-gradient-flow}, this time with rank $r=8$, we orthogonally 
 project the right-hand side of \eqref{ode-E-pass} onto the tangent space at $E$ of the manifold $\cM_r\subset \R^{m,n}$ of real rank-$r$ matrices, so that solutions starting with rank $r$ retain the rank $r$: 
\begin{equation} \label{ode-ErF-pass}
\dot E = P_E \Bigl( -G_\eps(E) + \langle G_\eps(E), E \rangle E \Bigr).
\end{equation}
Then also the Frobenius norm 1 is conserved (see (\ref{chap:proto}.\ref{ode-ErF-2-v2})), and $\F_\eps(E(t))$ decays monotonically (see (\ref{chap:proto}.\ref{eq:pos-real-r})).

Using the SVD-like factorization $E=USV^\top$, where  $U\in\R^{m,r}$ and $V\in\R^{n,r}$ have orthonormal columns and
$S\in \R^{r,r}$, the seemingly abstract differential equation \eqref{ode-ErF-pass} is solved numerically for the factors $U,V,S$ as described in Section~\ref{subsec:low-rank-integrator}.

\medskip\noindent
{\bf Using the low-rank structure in the eigenvalue computation.}
For the computation of the gradient matrix $G_\eps(E)$, one needs to compute the  eigenvalue of smallest positive real part and the associated left and right eigenvectors of  the Hamiltonian matrix $K(C+\eps E L^{-1})$. Except in the very first step of the algorithm, one can make use of the eigenvalue of smallest real part of the previous step in an inverse iteration (and possibly of the eigenvalues of second and third smallest real part etc.~to account for a possible exchange of the leading eigenvalue).

Moreover, we get from a perturbation $\eps E L^{-1}$ with $E=U\Sigma V^\top$ of rank 8 that $C$ is perturbed by $\Delta C=\eps E L^{-1}=\eps (U\Sigma) (L^{-\top}V)^\top$ of the same rank 8, which yields the perturbed Hamiltonian matrix
$$
K(C+\Delta C)= K(C) + \Delta K, 
$$
where the perturbation 
$\Delta K$ is still of moderate rank in view of \eqref{K-pass}.
This fact can be used in the computation of the required eigenvalues in the case of a high-dimensional system, using the Sherman--Morrison--Woodbury formula in an inverse iteration.
\index{Sherman--Morrison formula}

If $p\ll n$, then $K(C)$ can be viewed as a low-rank perturbation to the matrix
$$
\begin{pmatrix} A & 0 \\ 0 & -A^\top \end{pmatrix}.
$$
With the Sherman--Morrison--Woodbury formula, this can yield an efficient inverse iteration when $A$ is a large and sparse matrix for which shifted linear systems can be solved efficiently.

\subsection{Outer iteration}
For the solution of the scalar nonlinear equation \eqref{phi-pass} we use a Newton--bisection method in the variants discussed in Chapter~\ref{chap:two-level} or, for small $\delta$, the square root model and bisection as in Section~\ref{subsec:outer-it-sqrt}.

\medskip\noindent
Numerical results obtained with the above method for passivity enforcement are given by Fazzi, Guglielmi \& Lubich (\cite{FGL21}).

\subsection{$\Hinf$-nearest passive system} \label{subsec:pass-Hinfty}
We extend the approach of Section~\ref{sec:uncon} to the problem of $\Hinf$-optimal passivity enforcement.

\medskip\noindent
{\bf Problem ($\Hinf$).} 
{\it Given a system  \eqref{lti} with positive definite matrix $D+D^\top$
for which the Hamiltonian matrix $K(C)$ has some purely imaginary eigenvalues, 
and given $\delta>0$, compute a perturbed state--output matrix   $C + \Delta C$ with minimal $\Hinf$-norm of the perturbation to the transfer function, $\Delta H(\lambda)= \Delta C(\lambda I -A)^{-1}B$, such that all eigenvalues of $K(C + \Delta C)$ have a real part of absolute value at least~$\delta$.}



\medskip
We propose a multi-level approach similar to Section~\ref{sec:uncon}. We use the variable $S$ for~$\Delta C$.
\index{three-level iteration}

\begin{itemize}
\item {\bf Innermost iteration:\/}
For given $\eps>0$ and $S\in \R^{m,n}$, compute the spectral value set abscissa $\alpha_\eps(A,B,S,0)$,
e.g., by the rank-1 gradient flow algorithm of the Section~\ref{sec:Hinf}. This algorithm also yields a rank-1 matrix $E_\eps(S)\in\C^{n,n}$ of Frobenius norm~1 such that
\begin{equation} \label{alpha-pass}
\alpha_\eps(A,B,S,0) = \max_{E\in\C^{n,n},\,\|E\|_F=1} \Re\,\lambda(A+\eps BES)= \Re\,\lambda(A+\eps BE_\eps(S)S),
\end{equation}
where $\lambda(M)$ is a rightmost eigenvalue of a matrix $M$.
\item {\bf Inner iteration:\/}
For given $\delta>0$ and $\eps>0$, we use a constrained gradient system to solve the constrained eigenvalue optimization problem, over 
$S \in \R^{m,n}$,
\begin{equation}\label{eig-opt-pass-2}
S(\eps) =   \arg \min_{S} \alpha_\eps(A,B,S,0),
\quad\ \text{ subject to}\quad
\Re\,\lambda_{H}(K(C+S)) \ge \delta,
\end{equation}
where $\lambda_{H}(K)$ is an eigenvalue of minimal nonnegative real part of a Hamiltonian matrix $K\in \R^{2n,2n}$.

\item {\bf Outer iteration:\/} 
Compute the zero $\oeps$ of the monotonically growing function
\begin{equation}\label{phi-pass-2}
\phi(\eps)= \alpha_\eps(A,B,S(\eps),0).
\end{equation}
We use a   Newton--bisection method for this scalar equation.
\end{itemize}


We then have the following analogue of Theorem~\ref{thm:uncontrollability}, which is proved by the same argument based on Theorem~\ref{thm:Hinf-oeps}.

\begin{theorem}[$\Hinf$-distance]
\label{thm:opt-pass} 
Let $\oeps>0$ be the exact solution of the problem \eqref{alpha-pass}--\eqref{phi-pass-2}.
Then, the perturbed state--output matrix 
$C+ S(\oeps)$ yields a passive system having $\Re\,\lambda(K(C+S))\ge \delta$ with minimal $\Hinf$-distance between the transfer functions of the perturbed and unperturbed systems. This $\Hinf$-distance equals $1/\oeps$.
\end{theorem}

For the inner iteration we use a constrained gradient flow, based on the following direct analogue of Lemma~\ref{lem:gradient-uncon} and on Lemma~\ref{lem:gradient-pass}.

\begin{lemma}[Derivative of the spectral value set abscissa] \label{lem:gradient-2-pass}
Let $S(t)\in \R^{m,n}$, for real $t$ near $t_0$, be a continuously differentiable path of matrices, with the derivative denoted by $\dot S(t)$. Let $E(t)=E_\eps(S(t))\in\C^{n,n}$ of Frobenius norm 1 be such that
$$
\alpha_\eps(A,C,S(t),0) = \Re\,\lambda(A+\eps B E(t) S(t))
$$
and the rightmost eigenvalue on the right-hand side is simple, with left and right eigenvectors $x(t)$ and $y(t)$,
of norm 1 and with positive inner product. Let $\kappa(t)=1/(x(t)^*y(t))$. Assume that $S(t)y(t)\ne 0$ and $x(t)^*B\ne 0$.
Then we have, omitting the argument $t$ on the right-hand side,
$$
\frac1{\eps\kappa(t)}\,\frac{d}{dt}\,\alpha_\eps(A,B,S(t),0) =
\langle G_S, \dot S \rangle \quad \text{ with }\quad G_S=\Re\bigl((x^*BE)^*y^*\bigr).
$$
\end{lemma}

\section{Structured contractivity radius}
\index{linear time-invariant system!contractive}
\index{transfer function!contractive}
A linear time-invariant system \eqref{lti} is called {\it contractive} if its transfer function $H$ is bounded in the $\Hinf$-norm by
$$
\|H\|_\infty \le 1,
$$
and it is called {\it strictly contractive} if the above inequality is strict. Contractive systems play an important role as subsystems in large networks, because their composition remains contractive and thus yields well-controlled input-output relations. A strictly contractive system may be susceptible to perturbations (or uncertainties) in the entries of its matrices $(A,B,C,D)$, and it is then of interest to know which size of perturbations still guarantees contractivity. Here, we consider perturbations only in the state matrix $A$ and allow for  structured perturbations $\Delta\in \cS$, where the structure space $\cS\subset \C^{n,n}$ is a given complex- or real-linear subspace,
e.g. real matrices with a prescribed sparsity pattern. We study the following problem in this section.

Let $(A,B,C,D)$ be the system matrices of a strictly contractive linear time-invariant system. In particular, this implies $\|D\|_2<1$.
We consider structured perturbations $A \to A+\Delta$ with $\Delta\in \cS$, which yield perturbed transfer functions 
$$
H_\Delta(\lambda)=   C(\lambda I -A-\Delta)^{-1} B+D.
$$

\noindent
{\bf Problem.} {\it Find the largest possible perturbation size $\delta_\star>0$ such that} 
$$
\| H_\Delta \|_\infty \le 1 \quad\ \text{\it for all $\,\Delta\in \cS\,$ with $\,\|\Delta\|_F \le \delta_\star$\,.}
$$
The number $\delta_\star>0$ measures the robustness of contractivity of a system and is called the $\cS$-structured {\it contractivity radius}. We present and discuss a two-level algorithm that is closely related to that of
Section~\ref{sec:eps-stab} with $\eps=1$, to which it reduces for the special case $B=C=I$ and $D=0$.
\index{contractivity radius}

\subsection{Two-level iteration}
\index{two-level iteration}
We consider the matrix $M(\Theta)$ of \eqref{AEdef} that corresponds to $A+\Delta$ instead of $A$. So we let
$$
M(\Delta,\Theta)= A+ \Delta + B\Theta(I-D\Theta)^{-1}C.
$$
By Theorem~\ref{thm:basicequiv-sv-eig}, $\|H_\Delta\|_\infty = 1$ if and only if there exists $\Theta\in\C^{p,m}$ with
$\|\Theta\|_2=1$ such that $M(\Delta,\Theta)$ has an eigenvalue with nonnegative real part. Moreover, the optimizing matrix $\Theta$ is of rank 1, and hence its Frobenius and 2-norms are the same. We define the functional $\F_\delta(E^\cS,E)$
(for  $E^\cS\in \cS$ and $E\in\C^{p,m}$, both of unit Frobenius norm) by
\begin{equation}\label{F-eps-con}
    \F_\delta(E^\cS,E) = - \Re \,\lambda\bigl(M(\delta E^\cS,E)\bigr),
\end{equation}
where $\lambda(M)$ is the eigenvalue of $M$ of largest real part (and among those, the one with largest imaginary part). With this functional we follow the two-level approach of Section~\ref{sec:two-level}:
\begin{itemize}
\item {\bf Inner iteration:\/} For a given $\delta>0$, we aim to compute  matrices $E^\cS(\delta) \in \cS$ and $E(\delta) \in\C^{p,m}$, both of unit Frobenius norm,
that minimize $\F_\delta$:
\begin{equation} \label{E-theta}
(E^\cS(\delta),E(\delta)) = \arg\min\limits_{ E^\cS \in \cS, E \in \C^{p,m} \atop \| E^\cS \|_F= \| E \|_F = 1} \F_\delta(E^\cS,E).
\end{equation}

\item {\bf Outer iteration:\/} We compute the smallest positive value $\delta_\star$ with
\begin{equation} \label{zero-delta-sc}
\phi(\delta_\star)= 0,
\end{equation}
where $\phi(\delta)=  \F_\delta(E^\cS(\delta),E(\delta))=\alpha_\eps(A+\delta E^\cS(\delta),B,C,D)$ for $\eps=1$.
\end{itemize}
Provided that these computations succeed, we have from Theorem~\ref{thm:Hinf-oeps} that $\Delta_\star= \delta_\star E^\cS(\delta_\star) \in \cS$ is a perturbation matrix with $\| H_{\Delta_\star} \|_\infty =1$, and $\delta_\star$ is the $\cS$-structured contractivity radius.

\subsection{Rank-1 matrix differential equations for the inner iteration}
\index{rank-1 matrix differential equation}

The following lemma is obtained as in the proof of Lemma~\ref{lem:lambdaderiv-sys}.

\begin{lemma}[Structured gradient] \label{lem:gradient-cr}
Let $E^\cS(t)\in \cS$ and $E(t)\in\C^{p,m}$, for real $t$ near $t_0$, be continuously differentiable paths of matrices. 
Assume that $\lambda(t)$ is a simple eigenvalue of  $M (\delta E^\cS(t),E(t))$ depending continuously on $t$, with associated left and right eigenvectors
$\y(t)$ and $\x(t)$ normalized by \eqref{eq:scaling-sys}, and let 
$\kappa(t) = 1/(\y(t)^* \x(t))>0$.
Then, the derivative of $\F_\delta(E^\cS(t),E(t))$ 
is given by
\beq \label{eq:deriv-pass}
\begin{aligned}
\frac1{ \kappa(t) } \,\frac{d}{dt} \F_\delta(E^\cS(t),E(t)) &= \ \Re\,\bigl\langle  G_\Delta^\cS(\delta E^\cS(t),E(t)),  \delta \dot E^\cS(t) \bigr\rangle
\\[-2mm]
&\quad + \Re\,\bigl\langle  G_\Theta(\delta E^\cS(t),E(t)),  \dot E(t) \bigr\rangle
\end{aligned}
\eeq
with the (rescaled) gradient
\begin{equation}
\begin{aligned}
G_\Delta^\cS(\Delta,\Theta) &= \Pi^\cS(xy^*)\in \cS, \\
G_\Theta(\Delta,\Theta) &=  rs^*     \in \C^{p,m},
\end{aligned}
\label{eq:grad-uncon}
\end{equation}
where $r,s$ are obtained from $x,y$ via \eqref{eq:rsdef}.
\end{lemma}
\index{gradient!structured}

\medskip\noindent
\subsubsection*{Norm- and structure-constrained gradient flow.} 
\index{gradient flow!norm-constrained}
\index{gradient flow!structure-constrained}
Similar to Section~\ref{sec:eps-stab},
we consider the  projected gradient flow, with $G_\Delta^\cS=G_\Delta^\cS(\delta E^\cS,E)$ and $G_\Theta=G_\Theta(\delta E^\cS,E)$ for short,
\beq\label{ode-ES-E-contr}
\begin{aligned}
\delta \dot E^\cS &= -G_\Delta^\cS + \Re\left\langle G_\Delta^\cS,  E^\cS \right\rangle E^\cS,
\\
\dot E &= - G_\Theta +  \Re\left\langle G_\Theta,  E \right\rangle E.
\end{aligned}
\eeq
The unit Frobenius norm of $E^\cS$ and $E$ is conserved along trajectories and the functional $\F_\delta(E^\cS(t),E(t))$ decreases monotonically. At a non-degenerate stationary point $(E^\cS,E)$, where $G_\Delta^\cS$ and $G_\Theta$ do not vanish, we find that
$E^\cS$ and $E$ are real multiples of $G_\Delta^\cS$ and $G_\Theta$, respectively. Hence, $E^\cS$ is the projection onto $\cS$ of a rank-1 matrix, and $E$ is a rank-1 matrix.

\subsubsection*{Rank-1 matrix differential equations.} To make use of the rank-1 structure of optimizers, we proceed as in Section~\ref{sec:eps-stab}
and combine the rank-1 approaches of Sections~\ref{sec:proto-complex} and \ref{sec:proto-structured}. We 
consider differential equations for rank-1 matrices $Y(t)$ and $E(t)$, where the former yields $E^\cS(t)=\Pi^\cS Y(t)$. These differential equations are 
obtained from \eqref{ode-ES-E} by replacing $G_\Delta^\cS$ and $G_\Theta$ by their projections $P_Y$ and $P_E$ onto the tangent spaces of the manifold of rank-1 matrices at $Y$ and $E$, respectively:
\begin{equation} \label{ode-ES-E-1-contr}
    \begin{aligned}
       \delta\, \dot Y &= - P_Y G_\Delta^\cS + \Re\,\langle P_Y G_\Delta^\cS, E^\cS \rangle Y \quad\text{ with }\  E^\cS=\Pi^\cS Y,
       \\[1mm]
       \dot E &= -  P_E G_\Theta + \Re\,\langle G_\Theta, E \rangle E.
    \end{aligned}
\end{equation}
These differential equations yield rank-1 matrices $Y(t)$ and $E(t)$ and preserve the unit Frobenius norm of $E^\cS(t)$ and $E(t)$.  As in Sections~\ref{sec:proto-complex} and \ref{sec:proto-structured} it is shown that under a nondegeneracy condition, the stationary points $(Y,E)$ of \eqref{ode-ES-E-1-contr} correspond bijectively to the stationary points $(E^\cS,E)$ of \eqref{ode-ES-E-contr} via $E^\cS=\Pi^\cS Y$ and with the same $E$. The differential equations are integrated numerically into a stationary point $(E^\cS,E)$ as is done in Sections~\ref{sec:proto-complex} and \ref{sec:proto-structured}, working with the vectors that define the rank-1 matrices $Y$ and~$E$.

\subsection{Outer iteration} 
\index{Newton--bisection method}
For the solution of the scalar equation $\phi(\delta)=0$ we use a   Newton--bisection method as in
Section~\ref{sec:two-level}. The derivative of $\phi$ for the Newton iteration is obtained with the arguments of the proof of Theorem~\ref{chap:two-level}.\ref{thm:phi-derivative} (under analogous assumptions), which yields
$$
\phi'(\delta)= - \kappa(\delta) \,\| G_\Delta^\cS(\delta E^\cS(\delta),E(\delta))\|_F =
- \kappa(\delta)\, \| \Pi^\cS (x(\delta)y(\delta)^*) \|_F,  
$$
where $x(\delta)$ and $y(\delta)$ are left and right normalized eigenvectors associated with the rightmost eigenvalue of
$M(\delta E^\cS(\delta), E(\delta))$, and $\kappa(\delta)=1/(x(\delta)^*y(\delta))>0$.

\section{Descriptor systems} \label{sec:descriptor}
\index{descriptor system}
We consider a {\it descriptor system}, which formally differs from the system \eqref{lti} only in that the derivative of the state vector is multiplied with a singular matrix\footnote{In this section only we adhere to the convention in the control literature to denote by $E$ the matrix multiplying the time derivative of the state vector and to work with the matrix pencil $(A,E)$. In the rest of this book $E$ appears as a matrix of unit Frobenius norm when writing a perturbation matrix as $\Delta=\eps E$. In this section we will write instead $\Delta=\eps Z$ with $Z$ of Frobenius norm 1, choosing $Z$ as the letter of last resort.} $E\in\R^{n,n}$: for $t\ge 0$,
\begin{align}
    E \dot z(t) &= A z(t) + B u(t)    \label{descriptor}
    \\
    y(t) &= Cz(t).   \nonumber
\end{align}
We choose zero initial values and assume that the input function $u$ can be extended by zero to a sufficiently differentiable function on the whole real axis. Since $E$ is singular, the equation for the state vector $z$ is now a differential-algebraic equation instead of a pure differential equation. Descriptor systems arise naturally in systems with state constraints and in modelling and composing networks of such systems.

We assume that all finite eigenvalues of the matrix pencil $A-\lambda E$ have negative real part, i.e.,
\begin{equation}
\label{pencil-ass}
   \text{$A-\lambda E\,$ is invertible for all $\lambda\in\C$ with $\Re\,\lambda\ge 0$.}
 \end{equation}
In particular, the matrix $A$ is invertible. 

The matrix transfer function of the descriptor system is
\begin{equation}
\label{H-desc}
H(\lambda)= C(\lambda E-A)^{-1} B , \qquad \Re\,\lambda\ge 0.
\end{equation}

\begin{remark}
In the equation for the output $y$ in \eqref{descriptor} we have set the feedthrough matrix $D=0$ for convenience. The term $Du(t)$ could be added for nonzero $D$. The required changes in the theory and algorithm of this section can be done by combining the constructions and arguments of Section~\ref{sec:Hinf} with those given here. As we wish to concentrate on the effects of the singular matrix $E$, we chose to forego the technical complications resulting from a nonzero feedthrough matrix $D$, which were already dealt with in Section~\ref{sec:Hinf}.
\end{remark}

\subsection{Index and asymptotics of the transfer function at infinity}

In contrast to \eqref{lti}, the transfer function $H(\lambda)$ of \eqref{H-desc} need not be uniformly bounded for $\Re\,\lambda \ge 0$. We show that, aside from exceptional choices of $B$ and $C$, the norm of $H(\lambda)$ grows proportionally to $|\lambda|^{k-1}$ as $\lambda\to\infty$, where $k\ge 1$ is the {\it index} of the differential-algebraic equation $E\dot z = Az+ f$. The index can be determined from the Schur normal form of $A^{-1}E$ as follows.

We premultiply \eqref{descriptor} with $A^{-1}$ and $\lambda^{-1}$ so that the Laplace-transformed state equation $(\lambda E-A)\mathcal{L} z(\lambda) = B\,\mathcal{L} u(\lambda)$
becomes 
$$
\Bigl(A^{-1}E - \frac1\lambda\, I \Bigr) \mathcal{L} z(\lambda) = \frac1\lambda\, A^{-1}B \, \mathcal{L} u(\lambda).
$$
We want to understand the behaviour of the inverse of the matrix in brackets on the left-hand side as $\lambda\to \infty$, which is not obvious as $E$ is singular. To this end we transform to a block Schur normal form 
\begin{equation}
\label{block-schur-dae}
A^{-1}E = Q \begin{pmatrix} G & K \\ 0 & N \end{pmatrix} Q^\top
\end{equation}
with an orthogonal matrix $Q$, an invertible matrix $G$ and a nilpotent matrix $N$. The smallest integer $k\ge 1$ such that 
$$
N^k=0
$$
\index{matrix pencil!index}
\index{descriptor system!index}
is called the {\it index} of the matrix pencil $(A,E)$ (or of the differential-algebraic equation $E\dot z = Az+ f$, or of the descriptor system \eqref{descriptor}).

We have, for $\Re\,\lambda \ge 0$ and $\zeta=1/\lambda$,
$$
(A^{-1}E -\zeta I)^{-1} = Q \begin{pmatrix} (G-\zeta I)^{-1} & -(G-\zeta I)^{-1}K(N-\zeta I)^{-1} \\ 0 & (N-\zeta I)^{-1} \end{pmatrix} Q^\top.
$$
Here we note that 
$$
-\lambda^{-1} (N-\lambda^{-1} I)^{-1} =  (I-\lambda N)^{-1} =  I + \lambda N + \ldots +  \lambda^{k-1} N^{k-1}.
$$
We conclude that the norm of $H(\lambda)=  C (A^{-1}E - \lambda^{-1}I)^{-1} \lambda^{-1} A^{-1}B$ is bounded by a constant times $|\lambda|^{k-1}$ as $\lambda\to \infty$, and for generic $B$ and $C$ the asymptotic growth is actually proportional to $|\lambda|^{k-1}$.

\subsection{Weighted matrix transfer function and its $\Hinf$-norm}
\index{transfer function!weighted}
For a system of index $k>1$, we therefore want to bound the {\it weighted} matrix transfer function
\begin{equation}
\label{H-k-desc}
H^{[k]}(\lambda)=(1+\lambda)^{-(k-1)} H(\lambda), \qquad \Re\,\lambda \ge 0.
\end{equation}
(In the scalar factor, $\lambda$ should be replaced by $\tau\lambda$ with a characteristic time scale $\tau>0$, which we assume to be 1 for ease of presentation.) Note that the Laplace-transformed input-output relation is
\begin{equation}
\label{yu-L-Hk}
\mathcal{L} y(\lambda) = H(\lambda) \, \mathcal{L} u(\lambda) = H^{[k]}(\lambda) \, (1+\lambda)^{k-1}\mathcal{L} u(\lambda), \qquad \Re\,\lambda \ge 0,
\end{equation}
and that 
$$
(1+\lambda)^{k-1}\mathcal{L} u(\lambda)= \bigl(\mathcal{L} (1+d/dt)^{k-1}u\bigr)(\lambda)
$$
under our running assumption that $u$ together with its extension by $0$ to the negative real half-axis is a sufficiently differentiable function.
As in \eqref{yu-bound}, since $\| H^{[k]} \|_\infty = \sup_{\mathrm{Re}\, \lambda \ge 0} \|  H^{[k]}(\lambda) \|_2$ is finite,
these relations imply that the output $y$ is bounded in terms of the input $u$ as
\begin{equation}\label{yu-bound-k}
   \biggl( \int_0^T \| y(t) \|^2 \, dt \biggl)^{1/2} 
   \le \| H^{[k]} \|_\infty \, \biggl(\int_0^T \| (1+ d/dt)^{k-1} u(t) \|^2 \, dt\biggl)^{1/2}, \quad\
    0\le T \le \infty.
\end{equation}
Note that for index $k\ge 2$, the bound depends on derivatives of $u$ up to order $k-1$.
\\
This raises the following problem.

\medskip\noindent
{\bf Problem.}
{\it Compute the $\Hinf$-norm of the weighted matrix transfer function $H^{[k]}$ of the descriptor system.}

\medskip\noindent
In the following we restrict our attention to the case of principal interest where
\begin{equation} \label{ass:H-k-inf}
\| H^{[k]}(\infty) \|_2 <  \sup_{\mathrm{Re}\,\lambda \ge 0} \| H^{[k]}(\lambda) \|_2 \ (= \| H^{[k]} \|_\infty ).
\end{equation}
Then the supremum is a maximum that is attained at a finite $\lambda=\iu\omega$ on the imaginary axis.
We will modify the algorithm of Section~\ref{sec:Hinf} to compute $\| H^{[k]} \|_\infty$.

\subsection{$\Hinf$-norm via a stability radius}

In this subsection we give analogues of Theorems \ref{thm:basicequiv-sv-eig} and \ref{thm:Hinf-oeps} for the descriptor system \eqref{descriptor} with the transfer function $H(\lambda)$ of \eqref{H-desc}. We use the notation
$$
\Lambda(A,E)=\{ \lambda\in \C\,:\, A-\lambda E \,\text{ is singular}\}.
$$
The following result will later be used with $\chi(\lambda)=(1+\lambda)^{-(k-1)}$, where $k$ is the index of the matrix pencil $(A,E)$.

\begin{theorem}[Singular values and eigenvalues]\label{thm:basicequiv-pen}
Let $\eps > 0$, $\lambda\in \C\setminus \Lambda(A,E)$, and nonzero $\chi(\lambda)\in\C$.
The following two statements are equivalent:
\\[2mm]
    \ (i)$\ \ \|\chi(\lambda)H(\lambda)\|_2 \geq \eps^{-1}.$ 
\\[2mm]    
(ii)\ \  There exists $\Delta\in \C^{p,m}$ with $\|\Delta\|_F \le \eps$  such that $\lambda$ is an eigenvalue of the following nonlinear eigenvalue problem: There is an eigenvector $y\in \C^n \setminus \{0\}$ such that
    \begin{equation}
        \label{twoequiv-pen}
    (A + \chi(\lambda) B \Delta C - \lambda E)y=0.
    \end{equation}
Moreover, $\Delta$ can be chosen to have rank $1$, and the two inequalities can be replaced by equalities in the equivalence.
\end{theorem}
\index{nonlinear eigenvalue problem}

\begin{proof} For $\chi(\lambda)=1$, we can repeat the proof of Theorem~\ref{thm:basicequiv-sv-eig}, noting that there the replacement of $A-\lambda I$ by $A-\lambda E$ only leads to obvious changes. For general nonzero $\chi(\lambda)$,
we use the result with $\widetilde \eps= |\chi(\lambda)|\eps$ and $\widetilde \Delta=\chi(\lambda)\Delta$.
\qed
\end{proof}

 
Given a descriptor system such that the matrix pencil $(A,E)$ satisfies \eqref{pencil-ass} and is of index~$k$, 
and choosing $\chi(\lambda)=(1+\lambda)^{-(k-1)}$ so that $H^{[k]}(\lambda)=\chi(\lambda)H(\lambda)$, we proceed as in Section~\ref{sec:Hinf} and define the corresponding spectral value set
\index{spectral value set}
\begin{align*}
\Lambda_\eps^{[k]}&= \{ \lambda \in \C \setminus \Lambda(A,E)\,:\, \| H^{[k]}(\lambda) \|_2 \ge \eps^{-1} \}
\\
&= \{ \lambda \in \C \setminus \Lambda(A,E)\,:\, \lambda\ \text{satisfies (ii) of Theorem~\ref{thm:basicequiv-pen}} \}. \end{align*}
The spectral value abscissa
$$
\alpha_\eps^{[k]} = \sup \{ \Re\,\lambda \,:\, \lambda \in \Lambda_\eps^{[k]} \} 
$$
then yields
$$
  \sup_{\mathrm{Re}\,\lambda \ge \alpha_\eps^{[k]}}  \| H^{[k]}(\lambda) \|_2 = \frac1\eps,
$$
and by \eqref{ass:H-k-inf}, the supremum is a maximum if $\eps$ is so small that $\alpha_\eps^{[k]}\le 0$. With the stability radius
\index{stability radius}
$$
\oeps^{[k]} = \min\{ \eps>0 \,:\, \alpha_\eps^{[k]} = 0 \},
$$
we therefore again characterize the $\Hinf$-norm, which takes the maximum over $\Re\,\lambda\ge 0$, as the inverse stability radius. Since this result is essential for our numerical approach, we formulate it as a theorem.

\begin{theorem}[$\Hinf$-norm via the stability radius]\label{thm:Hinf-oeps-k}
Let the descriptor system be such that the matrix pencil $(A,E)$ has all finite eigenvalues with negative real part
and is of index~$k$. Then, the $\Hinf$-norm of the weighted matrix transfer function $H^{[k]}$ of \eqref{H-k-desc} and the stability radius $\oeps^{[k]}$ are related by
$$
    \| H^{[k]} \|_\infty = \frac1{\oeps^{[k]}}.
$$
\end{theorem}

\subsection{Two-level iteration}
\index{two-level iteration}
We use a two-level iteration similar to that of Section~\ref{sec:Hinf} to compute $\| H^{[k]} \|_\infty$.

\begin{itemize}
\item {\bf Inner iteration:\/} For a given $\eps>0$, compute the spectral value abscissa $\alpha_\eps^{[k]}$ using the nonlinear eigenvalue problem \eqref{twoequiv-pen} with perturbation matrices $\Delta$ of norm $\eps$ and of rank~1, which are determined via a rank-1 projected gradient flow that aims to maximize the real part of the rightmost eigenvalue.

\item {\bf Outer iteration:\/} Compute $\oeps^{[k]}$ as the smallest $\eps>0$ such that $\alpha_\eps^{[k]} = 0$, using some Newton--bisection method as discussed in Chapter~\ref{chap:two-level}.
\end{itemize}
Provided that these computations succeed, we obtain  $ \| H^{[k]} \|_\infty = 1/{\oeps^{[k]}}$ by Theorem~\ref{thm:Hinf-oeps-k}.

\subsection{Inner iteration: constrained gradient flow }
We aim to find $\Delta\in\C^{p,m}$ of Frobenius norm $\eps$ and of rank 1 such that the rightmost eigenvalue $\lambda$  yielding a singular matrix 
$$
M(\Delta,\lambda):=A+\chi(\lambda) B\Delta C - \lambda E
$$
has maximal real part, which equals the $\eps$-spectral value abscissa $\alpha_\eps^{[k]}$.
To this end we extend the norm- and rank-constrained gradient flow approach of Section~\ref{sec:proto-complex} for doing the inner iteration by a discretized rank-1 gradient flow.
Instead of Lemma~\ref{chap:proto}.\ref{lem:gradient} we now have the following gradient.

\begin{lemma}[Free gradient] \label{lem:gradient-desc}
Let $\Delta(t)\in \C^{p,m}$, for real $t$ near $t_0$, be a continuously differentiable path of matrices.
Let $\lambda(t)$ be a unique continuously differentiable path of eigenvalues  that yield singular matrices $M(\Delta(t),\lambda(t)) $ of co-rank~1, 
with associated left and right eigenvectors
$\y(t)$ and $\x(t)$. Assume that
$$
\eta(t):=\y(t)^* \bigl(E-\chi'(\lambda(t))B\Delta(t)C\bigr)\x(t)\ne 0 \quad\text{ and set }\quad \gamma(t) := \frac 1{\eta(t)}.
$$
Then, 
\begin{equation} \label{eq:deriv-desc}
-\Re \,\dot \lambda(t) = \Re \,\bigl\langle  G(\Delta(t)),  \dot \Delta(t) \bigr\rangle
\end{equation}
with the rank-1 matrix (omitting the argument $t$)
\begin{equation} \label{eq:freegrad-desc}
G(\Delta) =  -\bigl(\gamma\chi(\lambda)C y x^* B\bigr)^* \in \C^{p,m}.
\end{equation}
\end{lemma}
\index{gradient!free}

\begin{proof} Differentiating the matrix $A+\chi(\lambda(t)) B\Delta(t) C - \lambda(t) E$ and multiplying with $x(t)^*$ from the left and $y(t)$ from the right yields (omitting the argument $t$)
$$
-\eta{\dot\lambda} + \chi(\lambda)x^*B \dot\Delta Cy =0,
$$
which implies 
$$
\dot\lambda = \gamma \chi(\lambda)x^*B \dot\Delta Cy = \langle x, \gamma\chi(\lambda) B \dot\Delta Cy \rangle= 
\langle \conj{\gamma\chi(\lambda)} B^* xy^* C^*, \dot \Delta \rangle
$$
and hence yields the stated result.
\qed    
\end{proof}

With the gradient \eqref{eq:freegrad-desc}, the whole programme of Section~\ref{sec:proto-complex} carries through, as we briefly sketch in the following.
We write $\Delta\in \C^{p,m}$ of Frobenius norm $\eps$ as 
$$
\Delta = \eps Z \quad\text{with}\quad \|Z\|_F=1
$$
(in previous sections we wrote $\Delta = \eps E$ with $\|E\|_F=1$, but now $E$ is the singular matrix in the descriptor system)  and
$$
G_\eps(Z)= G(\eps Z).
$$
As in \ref{chap:proto}.\ref{ode-E}, we consider the gradient flow on the Frobenius-norm unit sphere of $\C^{p,m}$,
\begin{equation}
    \label{ode-Z-desc}
    \dot Z = - G_\eps(Z) + \Re \langle G_\eps(Z),Z \rangle Z.
\end{equation}
As in Theorem~\ref{chap:proto}.\ref{thm:monotone} we have that $-\Re\,\lambda(t)$ decreases along solutions of \eqref{ode-Z-desc}, where $\lambda(t)$ is a rightmost eigenvalue yielding a singular matrix $M(\lambda(t))$, provided the assumptions of Lemma~\ref{lem:gradient-desc} are satisfied.

\index{stationary point}
Stationary points $Z_\star$ of \eqref{ode-Z-desc} are again real multiples of $G_\eps(Z_\star)$ and are therefore of rank 1, as in Corollary~\ref{chap:proto}.\ref{cor:rank-1}.
As in Section~\ref{subsec:rank1-gradient-flow} we therefore consider the rank-1 constrained gradient flow
\index{gradient flow!rank-1 constrained}
\begin{equation}
    \label{ode-Z-1-desc}
    \dot Z = - P_Z \Bigl( G_\eps(Z) - \Re \langle G_\eps(Z),Z \rangle Z \Bigr),
\end{equation}
where $P_Z$ is the orthogonal projection onto the tangent space at $Z$ of the manifold of rank-1 matrices in $\C^{p,m}$.
This differential equation has the same properties as Equation (\ref{chap:proto}.\ref{ode-E-1}) and is discretized in the same way, as described in Section~\ref{subsec:proto-numer}.

\subsubsection*{Computing eigenvalues of the nonlinear eigenvalue problem.}
\index{nonlinear eigenvalue problem}
What differs from Section~\ref{sec:proto-complex} is the computation of eigenvalues $\lambda$ that yield a singular matrix $M(\lambda)$. If $\eps$ is sufficiently small, this can be done efficiently by a fixed-point iteration. Given an iterate $\lambda_n$, we set $A_n=A+\chi(\lambda_n) B\Delta C$ and compute $\lambda_{n+1}$ as the rightmost eigenvalue of the matrix pencil $A_n - \lambda E$, for which $\zeta_{n+1}=(\lambda_{n+1}+1)/(\lambda_{n+1}-1)$ is the eigenvalue of largest modulus of the matrix $(A_n - E)^{-1}(A_n + E)$.
When $\eps$ is not small, one can use algorithms for general nonlinear eigenvalue problems, such as the method of Beyn (\cite{Bey12}) or other methods as reviewed by G\"uttel \& Tisseur (\cite{GueT17}).

\subsection{Outer iteration}
We need to compute the zero of 
$$
\phi(\eps)=-\alpha_\eps^{[k]}= -\Re \, \lambda_\eps,
$$
where $\lambda_\eps$ is a rightmost eigenvalue of the nonlinear eigenvalue problem for the matrix-valued function $A+\chi(\lambda) B\eps Z_\eps C - \lambda E$ and $Z_\eps$ maximizes the real part of the rightmost eigenvalue among all matrices $Z\in\C^{p,m}$ of Frobenius norm 1. Note that $Z_\eps$ is to be computed in the inner iteration.

As in Chapter~\ref{chap:two-level}, we use some Newton--bisection method, for which we need the derivative $\phi'(\eps)$. As in Theorem~\ref{chap:two-level}.\ref{thm:phi-derivative} we find, under appropriate regularity assumptions, that
$$
\phi'(\eps) = -\| G_\eps(Z_\eps) \|_F.
$$
\index{Newton--bisection method}

\section{Notes}

Robust control is certainly one of the most important applications of the theory presented in this book and a field of research that has received much attention. For a survey we refer the reader to the seminal monograph by Hinrichsen \& Pritchard (\cite{HinP05}) and the references therein.

\subsubsection*{$\Hinf$-norm of the matrix transfer function.}

The $\Hinf$ norm 
associated with a linear time invariant (LTI) system is a basic quantity for measuring
robust stability; see e.g. Hinrichsen \& Pritchard (\cite{HinP86a},\cite{HinP90}) and Zhou, Doyle \& Glover (\cite{ZDG96}).

Numerical methods for computing the $\Hinf$-norm have been known for a long time.
Most of them are based on relations between the $\Hinf$-norm and the spectrum of certain
Hamiltonian matrices or pencils. 
The standard method to compute the $\Hinf$-norm is the 
algorithm due to Boyd \& Balakrishnan (\cite{BoyB90}) and Bruinsma \& Steinbuch (\cite{BruS90}), henceforth called the BBBS algorithm,
which generalizes and improves an algorithm of Byers (\cite{Bye88}, \cite{BN93}) for
computing the distance to instability for $A$.  The method relies on
Definition \ref{def:hinfnormcont}: 
for stable $A$, it needs to maximize $\|H(\iu \omega)\|$
for $\omega\in\R$.  The key idea is that, given any $\delta>0$, it is possible
to determine whether or not $\omega\in\R$ exists such that $\|H(\iu \omega)\|=\delta$ by
computing all eigenvalues of an associated $2n\times 2n$ Hamiltonian matrix and determining whether
any are imaginary.  The algorithm is globally quadratically convergent, but the computation of
the eigenvalues and the evaluation of the norm of the transfer matrix both
require on the order of $n^3$ operations, which is not practical when $n$ is sufficiently large.

The algorithm presented in Section~\ref{sec:Hinf} computes the $\Hinf$-norm as the reciprocal of the stability radius $\oeps$, as given in Theorem~\ref{thm:Hinf-oeps}, by extending the two-level algorithm of Chapter~\ref{chap:two-level} with eigenvalue optimization by the rank-1 matrix gradient flow of Chapter~\ref{chap:proto} in the inner iteration.  

\subsubsection*{Discrete-time systems.}

For discrete-time systems 
\begin{eqnarray}
x_{k+1} & = & Ax_k + Bu_k
\label{lindynsysdiscr}
\\
y_k & = & Cx_k + Du_k
\nonumber
\end{eqnarray}
there are extensions of the BBBS algorithm and also
an analogous approach to the algorithm presented in Section~\ref{sec:Hinf} can be adopted, with $\F_\eps(E) = -|\lambda|^2$ and 
$\phi(\eps) = \F_\eps\left(E(\eps) \right) = \rhoeps(A,B,C,D)$,
where
\beq
         \rhoeps(A,B,C,D) =  \max\{|\lambda| : \lambda \in \SVSeps(A,B,C,D)\} \label{rhoepsdef}
\eeq 
is the spectral value set radius. 
Here one looks for the smallest solution of $\phi(\oeps) = 1$.

\subsubsection*{Distance to uncontrollability.}

Given $A \in \C^{n \times n}, B \in \C^{n \times m}$, the linear control system
$\dot x = Ax + Bu$
is controllable if for every pair of states 
$x_0, x_1 \in \C^n$ there exists a continuous control
function $u(t)$ able to determine a trajectory $x(t)$ such that $x(t_0)= x_0$ and $x(t_1)=x_1$ within finite time $t_1-t_0$.
Equivalently, according to a result by Kalman (\cite{Kal63}), the system is
controllable if the matrix 
$(A - \lambda I,\ B)$ has full row rank for all $\lambda \in \C$.

\bng
Several algorithms have been proposed to compute a nearest uncontrollable system to a given controllable system. The distance to uncontrollability is a measure of robustness of the controllability of the system with respect to perturbations.
The first globally convergent algorithm has been proposed by Gu (\cite{Gu00}), which has complexity $\bigo(n^6)$; afterwards Gu, Mengi, Overton, Xia \& Zhu (\cite{GMO06}) improved it to average-case complexity $\bigo(n^4)$.
Burke, Lewis and Overton (\cite{BLO04}) proposed a trisection variant of Gu's algorithm, which is arbitrarily accurate, in contrast to Gu's algorithm, which computes the distance only within a factor of two.
\eng



An algorithm for measuring the robustness of controllability with respect to the $\Hinf$-norm, as investigated in Section~\ref{sec:uncon}, appears to be a novel contribution to the literature, where other norms are usually considered. The invariance of the $\Hinf$ distance with respect to changes of bases representing the LTI system makes this metric particularly appealing.

\subsubsection*{Nearest passive system.}
Classical references for  passivity of LTI systems are the monographs by Anderson \& Vongpanitlerd (\cite{AndV73}) and  Grivet Talocia \& Gustavsen (\cite{GriG15}).

Passivity enforcement has been addressed in the literature by various authors. 
Grivet-Talocia (\cite{Grivet2004})  used the characterization of passivity in terms of Hamiltonian matrices and applied first-order perturbation theory to move and coalesce imaginary eigenvalues of the associated Hamiltonian matrix. 
Schr\"oder \& Stykel (\cite{SchS07}) extended the approach of Grivet-Talocia, using structure-preserving algorithms to compute the required eigenvalues and eigenvectors. Two-level algorithms based on Hamiltonian eigenvalue optimization for finding a nearest passive or non-passive system with respect to the
Frobenius-norm of the perturbation to the state-to-output matrix $C$, as described in Section~\ref{sec:pass}, have been studied by Fazzi, Guglielmi \& Lubich (\cite{FGL21}).

Br\"ull \& Schr\"oder (\cite{BruS13})  substantially extended the approach to  descriptor systems and to more general notions of dissipativity that they characterized in terms of structured matrix pencils, which they used together with first-order perturbation theory to move imaginary eigenvalues. Gillis \& Sharma (\cite{GilS18}) studied the problem of finding the nearest port-Hamiltonian descriptor system. Such systems are always passive and conversely, with the appropriate terminology, every extended strictly passive system is port-Hamiltonian. 
For related work on port-Hamiltonian systems see also Gillis \& Sharma (\cite{GilS17}) and
Mehl, Mehrmann, \& Sharma (\cite{M16}) and Mehl, Mehrmann, \& Sharma (\cite{MMS17}).

The problem of determining the distance to non-passivity, or passivity radius, was studied by Overton \& Van Dooren (\cite{OveVD05}) and recently
by Mitchell and  \& Van Dooren (\cite{MitVD23})
for complex perturbations of the system matrices, but the case of real perturbations as considered in Section~\ref{sec:pass} has been left open.

Algorithms for finding the nearest passive or non-passive system with respect to the $\Hinf$-distance of the matrix transfer functions, as proposed in Section~\ref{sec:pass}, appear to be new.

\subsubsection*{Structured contractivity radius.}

This is a notion that -- to our knowledge -- is first investigated in this book.

\subsubsection*{Descriptor systems.}

Benner \& Voigt (\cite{BV14}) proposed an approach to the computation of the $\Hinf$-norm for descriptor systems of index 1, still based on the theory of spectral value sets and rank-$1$ iterations. Benner, Lowe \& Voigt (\cite{BLV18}) addressed the case of large descriptor systems, for which an algorithm using structured iterative eigensolvers is proposed.

The computation of the $\Hinf$-norm of weighted transfer functions of descriptor systems of index $k>1$ has apparently not been considered before.

\subsubsection*{Structured singular values.}

Structured singular values (SSV), as defined by Packard \& Doyle (\cite{PD93}), are an important and versatile tool in control, as they allow addressing a central problem in the analysis and synthesis of control systems: to quantify the stability of a closed-loop linear time-invariant systems subject to structured perturbations.
The class of structures addressed by the SSV is very general and allows to cover all types of parametric uncertainties that
can be incorporated into the control system via real or complex linear fractional transformations.
We refer e.g.  to Chen, Fan \& Nett ~(\cite{CFN96a}, \cite{CFN96b}), 
Hinrichsen \& Pritchard (\cite{HinP86a}), Karow (\cite{Ka11}), Karow, Hinrichsen \& Pritchard (\cite{KHP06}), Qiu et al. (\cite{QiuBRDYD95}) and the monograph by Zhou, Doyle \& Glover (\cite{ZDG96}) and the references therein for examples and applications of the SSV.

The versatility of the SSV comes at the expense of being notoriously hard, in fact NP hard to compute; see Braatz, Young, Doyle \& M.~Morari~(\cite{BYDM94}).
Algorithms used in practice thus aim at providing upper and lower bounds, often resulting in a coarse estimate of 
the exact value. An upper bound of the SSV provides sufficient conditions to 
guarantee robust stability, while a lower bound provides sufficient conditions 
for instability and often also allows to determine structured perturbations that 
destabilize the closed loop linear system.

\bng
The function {\tt mussv} in the \matlab\  Control Toolbox 
computes an upper bound of the SSV using diagonal balancing / LMI (Linear Matrix Inequalities) techniques proposed by Young,  Newlin \& Doyle~(\cite{YND92}) and
Fan, Tits \& Doyle (\cite{FTD91}).
\eng

In Guglielmi, Rehman \& Kressner (\cite{GRK17}), a gradient system approach in the spirit of the algorithms presented in this book is proposed,  which often provides tighter bounds than those computed by {\tt mussv}.

\chapter{Graphs}
\label{chap:graphs}

\newcommand{\cE}{\mathcal{E}}
\newcommand\Lap{L}

\newcommand\sym{\mathrm{Sym}}
\newcommand{\diagvec}{\mathrm{vec\,diag}}

\newcommand\Pact {P^+}
\newcommand\Gact{G^+}

\index{graph}
Fundamental properties of graphs such as connectivity and centrality are characterized by eigenvalues or eigenvectors of matrices such as the graph Laplacian or the weighted adjacency matrix. In this chapter we study exemplary problems of clustering and ranking for undirected weighted graphs, which can be viewed as matrix nearness problems that refer to eigenvalues and eigenvectors of structured symmetric real matrices. We propose algorithms that do eigenvalue optimization without eigenvalues but instead work directly with the quadratic form (Rayleigh quotients). This requires only matrix-vector products and vector inner products.

\section{Cutting edges}
\label{sec:cut}

We propose an iterative algorithm for partitioning a given weighted undirected graph under constraints, e.g. under must-link and cannot-link constraints or under cardinality constraints. The algorithm changes the weights such that the second eigenvalue of the graph Laplacian is driven to zero. The connected components of the partitioned graph are then read off from the corresponding eigenvector, known as the Fiedler vector. The proposed algorithm only requires the computation of matrix-vector products with the graph Laplacian and vector inner products but no computations of eigenvalues and eigenvectors. It is a two-level algorithm that follows a graph diffusion equation into a stationary point in the inner iteration and determines the cut in the outer iteration.

\subsection{Problem formulation and algebraic connectivity}
\label{sec:problem}

\subsubsection*{The minimum cut problem as a matrix nearness problem.}
\index{graph!constrained minimum cut}
A graph $\mathcal{G}=(\mathcal{V},\cE)$ is given by a set of vertices (or nodes) $\mathcal{V}=\{1,\dots,n \}$ and an edge set $\cE\subset \mathcal{V}\times\mathcal{V}$. We consider an {\it undirected} graph: if $(i,j)\in\cE$, then $(j,i)\in\cE$. With the undirected graph we associate {\it weights} $a_{ij}$ for $(i,j)\in\cE$ with 
$$
a_{ij}=a_{ji} > 0 \quad\hbox{ for all }\ (i,j)\in\cE.
$$
The graph is {\it connected} if for all $i,j\in\mathcal{V}$, there is a path $(i_0,i_1),(i_1,i_2),\dots,(i_{\ell-1},i_\ell)\in\cE$ of arbitrary length $\ell$, such that $i=i_0$ and $j=i_\ell$ and $a_{i_{k-1},i_k}>0$ for all $k=1,\dots,\ell$.
\index{graph!connected}

Given a connected weighted undirected graph with weights $a_{ij}$, we aim to find a {\it disconnected} weighted undirected graph with the same edge set $\cE$ and modified weights $\wt a_{ij}$ such that 
\begin{equation}\label{dist-2}
\sum_{(i,j)\in\cE} (\wt a_{ij} - a_{ij})^2 \quad\hbox{ is minimized.}
\end{equation}
In the minimizer the weights $\wt a_{ij}$ must either be zero or remain unchanged, i.e. $\wt a_{ij} = a_{ij}$.
Solving the matrix nearness problem \eqref{dist-2} is therefore equivalent to finding a cut  $\mathcal{C}$, i.e. a set of edges that yield a disconnected graph when they are removed from $\cE$, where
$$
\text{the cut $\mathcal{C}$ is such that }\ \sum_{(i,j)\in\mathcal{C}} a_{ij}^2 \quad\hbox{ is minimized.}
$$
Replacing the weights by their square roots yields $a_{ij}$ instead of $a_{ij}^2$ in the above sum, which is the usual
formulation of
the minimum cut problem.

\subsubsection*{Constrained minimum cut problems.}
\label{subsection:constrained-mincut}
\index{graph partitioning}
\index{graph partitioning!must-link constraint}
\index{graph partitioning!cannot-link constraint}
\index{graph partitioning!cardinality constraint}
\index{graph partitioning!membership constraint}
The above problem will be considered with additional constraints such as the following:
\begin{itemize}
\item {\it Must-link constraints:}\/ Pairs of vertices in a given set $\mathcal{P}\subset \mathcal{V} \times \mathcal{V}$ are in the same connected component. This specifies subgraphs that must not be cut.
\item {\it Cannot-link constraints:}\/ Pairs of vertices in a given set $\mathcal{P}\subset \mathcal{V} \times \mathcal{V}$ are in different connected components.
\item {\it Membership constraints:}\/ A given set of vertices $\mathcal{V}^+\subset\mathcal{V}$ is in one team (connected component) and another given set of vertices $\mathcal{V}^-\subset\mathcal{V}$ is in the other team.
\item {\it Cardinality constraint:}\/ Each of the connected components has a prescribed minimum number of vertices.
\end{itemize}

\subsubsection*{Graph Laplacian and algebraic connectivity.}
\index{graph!Laplacian}
\index{graph!algebraic connectivity}
\index{graph!weighted adjacency matrix}
Setting $a_{ij}=0$ for $(i,j)\notin \cE$, we have the {\em weighted adjacency matrix}
$$
A=(a_{ij}) ,
$$
which is a symmetric non-negative real matrix  with the sparsity pattern given by the set of edges $\cE$, i.e. $a_{ij}\ne 0$ only if $(i,j)\in\cE$.
The {\em weighted degrees}  $d_i = \sum_{j=1}^n a_{ij} $ are collected in the diagonal matrix 
$$
D(A) = \diag(d_i)= \diag(A  \one), \qquad\hbox{where $\one:=(1,\ldots,1)^\top \in \R^n$.}
$$
The {\it Laplacian matrix} $\Lap(A)$ is defined by
$$
\Lap(A) = D(A)-A=\diag(A  \one) - A .
$$
By the Gershgorin circle theorem, all eigenvalues of $\Lap(A)$ are non-negative, and 
$\Lap(A)\one=0$, so that $\lambda_1=0$ is the smallest  eigenvalue of $\Lap(A)$ with the eigenvector $\one$. The connectivity of the graph is characterized by the second-smallest eigenvalue of $\Lap(A)$, as the following basic theorem shows.

\begin{theorem}[Fiedler \cite{Fie73}] \label{thm:fiedler}
 Let $A \in \R^{n \times n}$ be the weighted adjacency matrix of an undirected graph and $\Lap(A)$ the graph Laplacian. Let $0 = \lambda_1 \le \lambda_2 \le \ldots \le \lambda_n$ be the eigenvalues of $\Lap(A)$. Then, the graph is disconnected if and only if $\lambda_2 = 0$. Moreover, if  $0=\lambda_2<\lambda_3$, then the entries of the corresponding eigenvector orthogonal to $\one$ assume only two different values, of different sign, which indicate the membership to the two connected components.
\end{theorem}
\index{Fiedler vector}

Because of this result, the second smallest eigenvalue $\lambda_2$ of $\Lap(A)$ is called {\it algebraic connectivity} of $A$. The corresponding eigenvector orthogonal to $\one$ is known as the {\it Fiedler vector}, provided that $\lambda_2<\lambda_3$. There are extensions of this theorem to $k>2$ components ($0=\lambda_2=\ldots=\lambda_k<\lambda_{k+1}$)
\begin{remark}[Normalized Laplacian] 
Instead of the above Laplacian $L=L(A)$, often a normalized Laplacian is considered in the literature (e.g. von Luxburg (\cite{Lux07})), which is given as
$$
\widehat L = \widehat{\Lap }(A)= I - D^{-1/2} A D^{-1/2} \qquad\text{or}\qquad \widetilde L = I - D^{-1} A.
$$
The eigenvalues $\lambda$ of these normalized Laplacians are generalized eigenvalues $Lv=\lambda Dv$ of the Laplacian $L$, and Fiedler's theorem is equally valid for them.
Our approach below extends in a straightforward way to these normalized cases.  
\end{remark}

\begin{remark}[Partitioning into more than two components]
    There are extensions of Theorem~\ref{thm:fiedler} to $k>2$ components ($0=\lambda_2=\ldots=\lambda_k<\lambda_{k+1}$). While this will not be considered in this section, we note that the algorithmic approach taken here for cutting into two components extends to partitioning into $k>2$ components.
\end{remark}

\subsection{Minimum cut problem as a matrix nearness problem}
\label{sec:two-level method}

In this subsection we propose a two-level approach for the unconstrained minimum cut problem \eqref{dist-2}. The extension to minimum cut problems with constraints is presented in the next subsection. Throughout this section, we denote the structure space of the matrices in the problem by
$$
\cS = \sym(\cE),
$$
the space of symmetric real $n\times n$ matrices with the sparsity pattern given by the edges $\cE$ of the graph.

In view of Theorem~\ref{thm:fiedler}, the minimum cut problem can be reformulated as the following matrix nearness problem.

\bigskip\noindent
{\bf Problem.} {\it Given $A\in \cS$ with $A\ge 0$, find a matrix $\Delta\in\cS$ of minimal Frobenius norm such that $A+\Delta\ge 0$ and the Laplacian $\Lap(A+\Delta)$ has a double eigenvalue $0$.}

\bigskip\noindent
\index{two-level iteration}
{\bf Two-level iteration.} We propose a two-level approach in the spirit of Chapter~\ref{chap:two-level}.
\begin{itemize}
\item {\bf Inner iteration:\/}
Given $\eps>0$, we look for a matrix $E=(e_{ij})\in\cS$ of unit Frobenius norm,  with $A+\eps E\ge 0$ (with componentwise inequality) such that the second smallest eigenvalue of $\Lap (A+\eps E)$ is minimal. The minimizer is denoted by~$E(\eps)$.
\item {\bf Outer iteration:\/}
 We determine the smallest value of $\eps>0$ such that the second smallest eigenvalue of 
$\Lap (A+\eps E(\eps))$ equals $0$.
\end{itemize}

\smallskip\noindent
Computing $E(\eps)$ for a given $\eps>0$ is based on a constrained gradient flow for minimizing the functional given by the Rayleigh quotient 
\index{Rayleigh quotient}
\index{eigenvalue optimization!without eigenvalues}
\begin{equation}\label{F-eps-cut}
\F_\eps(E,v) = \tfrac12\,v^\top \Lap (A+\eps E) v,
\end{equation}
for $E\in \cS$ of Frobenius norm~1 with $A+\eps E\ge 0$ and vectors $v\in\R^n$ of Euclidean norm~1 that are orthogonal to the eigenvector $\one$ to the eigenvalue $0$ of $\Lap(A+\eps E)$:
$$
\one^\top v=0.
$$
This will be further reduced to a functional of only the vector $v$.

In the outer iteration we compute the optimal $\eps$ by a combined Newton-bisection method. This optimal $\eps$, denoted $\oeps$, is the smallest $\eps$ such that $\lambda_2(\Lap(A+\eps E))=0$ for some admissible $E$ of unit norm.

The algorithm computes a partition of the graph as provided by the Fiedler vector corresponding to the weight matrix $A+\oeps E(\oeps)$, 
which is the eigenvector orthogonal to $\one$ for the double eigenvalue $0$. 

The iteration can be stopped if, with given tolerance parameters $\theta$ and $\vartheta$,
$$
\text{$\F_\eps(E,v) \le \theta \,\lambda_2(A)$ \ and \ $\| v^\pm - \langle v^\pm \rangle \one \| \le \vartheta$.}
$$ 
Here, $v^+=(\max(v_i,0))$ and $v^-=\min(v_i,0)$ collect the positive and negative entries of $v$, and 
$\langle v^\pm \rangle$ is the arithmetic mean of the nonzero entries in $v^\pm$.
If this condition is satisfied, then the signs of the entries of $v$ mark the membership to the two connected components.


\subsection{Inner iteration: diffusion on the graph}
In this subsection we adapt the programme of Section~\ref{subsec:rayleigh-max} to the current problem.

\subsubsection*{Adjoint of the Laplace map $\Lap(\cdot)$.}
\label{subsec:gradient-cut}
Let $\Pi^\cS$  be the orthogonal projection from $\R^{n,n}$ onto $\cS=\sym(\mathcal{E)}$: for $Z=(z_{ij})$,
\begin{equation}\label{Pi-S-cut}
\Pi^\cS Z\big|_{ij} = \begin{cases}
	\tfrac12(z_{ij}+z_{ji})\,, &\text{if } (i,j)\in\cE \,,\\
	0\,, &\text{otherwise.}\end{cases}
\end{equation}
Considering the dependence on the adjacency matrix $A\in\cS$ in the graph Laplacian $L(A)$ leads us to the Laplace operator $\Lap$ as a linear map
$$
\Lap: \cS\to \R^{n,n}.
$$
Its adjoint with respect to the Frobenius inner product,
$$
\Lap^*: \R^{n,n}\to \cS,
$$
is defined by $\langle \Lap^*(V),W \rangle = \langle V, \Lap(W) \rangle$ for all  $V\in \R^{n,n}$ and $W\in \cS$ and
is given explicitly in the following lemma.
\begin{lemma}[Adjoint of the Laplace operator]
For $V\in \R^{n,n}$,  
$$
\Lap^*(V) = \Pi^\cS (\diagvec(V)\one^\top - V),
$$
where $\diagvec(V)\in \R^n$ is the vector of the diagonal entries of $V$. 
\end{lemma}

\begin{proof}
For all $V\in \R^{n, n}$ and $W\in \cS$,
\begin{align*}
\langle V, \Lap(W) \rangle &= \langle V, \diag(W\one)-W \rangle = \langle \diagvec(V), W\one \rangle - \langle V,W \rangle 
\\
&= \langle \diagvec(V)\one^\top, W \rangle - \langle V,W \rangle = \langle \Pi^\cS(\diagvec(V)\one^\top -V), W \rangle,
\end{align*}
which yields the result.
\qed
\end{proof}

\subsubsection*{Structured gradient.}
\index{gradient!structured}
Consider a smooth path of matrices $E(t)\in \cS$ and vectors $v(t)\in\R^n$. Then, omitting the argument $t$ on the right-hand side,
\begin{align}
\nonumber
    \frac{d}{dt}\F_\eps(E(t),v(t))&=
    \tfrac12\,v^\top \Lap(\eps \dot E) v +
    v^\top \Lap(A+\eps E) \dot v 
    \\
     \nonumber
    &= \tfrac12\,\langle vv^\top,\Lap(\eps\dot E) \rangle +
     \langle \Lap(A+\eps E)v,\dot v\rangle 
 \\
 \label{F-dot-cut}
 &=
    \eps \,\langle G,\dot E \rangle +
     \langle g,\dot v\rangle 
\end{align}
with the (rescaled) gradient
\begin{equation}
\label{grad-cut}
   G=G_\eps(v)=\tfrac12\Lap^*(vv^\top) \in \cS,\quad\ 
g=g_\eps(E,v)= \Lap(A+\eps E)v \in \R^n. 
\end{equation}

\subsubsection*{Constrained gradient flow.}
With the constraints of norm 1 for $E$ and $v$ and the orthogonality constraint $\one^\top v=0$ we have the constrained gradient flow
\begin{align}
    \label{ode-E-cut}
    \dot E &= - G + \langle G,E \rangle E
    \\
    \label{ode-v-cut}
    \dot v &= -g + \langle g,v \rangle v + \frac1n \langle g,\one \rangle \one .
\end{align}

\subsubsection*{Non-negativity constraint.}
To satisfy the non-negativity condition of the perturbed weight matrix $A+\eps E(t)$,  we need that $\dot{e}_{ij}(t)  \ge 0$ for all $(i,j)$ with
$a_{ij} + \varepsilon e_{ij}(t) = 0$. So we replace
\eqref{ode-E-cut} by
$$
\dot E = \wh{P^+} \bigl(- G + \gamma E \bigr)
$$
where $\gamma$ is to be chosen such that $\langle E, \dot E \rangle =0$ and
$$
\wh{P^+}(t)(Z)\big|_{i,j} =
\begin{cases}
	0 &\text{if }\ a_{ij} + \varepsilon e_{ij}(t) = 0 \ \text{ and }\  z_{ij}\le 
 0,\\
  z_{ij}    &\text{otherwise.}
\end{cases}
$$
On every interval where $\cE^+(t):=\{ (i,j)\in\cE\,:\, a_{ij} + \varepsilon e_{ij}(t) > 0 \}$ does not change with~$t$,
the differential equation simplifies to
\begin{equation}
\label{ode-E-plus-cut}
\dot E =-\Pact G + \gamma \Pact E 
\quad\hbox{ with }\quad
\gamma = \frac{\langle \Pact G, \Pact  E \rangle }{\| \Pact  E\|^2},
\end{equation}
where 
$$
P^+Z\big|_{ij} = \begin{cases}
	z_{ij} &\text{if } (i,j)\in\cE^+ \,,\\
	0 &\text{otherwise.}\end{cases}
$$
Along every solution $(E(t),v(t))$ of \eqref{ode-E-plus-cut} together with  \eqref{ode-v-cut}, all constraints $\|E\|_F=1$, $\|v\|=1$, $\one^\top v=0$ and $A+\eps E\ge 0$ remain satisfied if they are at the initial value. The functional \eqref{F-eps-cut} decays monotonically along 
the solution, as is readily verified using \eqref{F-dot-cut}.

\subsubsection*{Stationary points.} 
\index{stationary point}
In a stationary point $(E,v)$ we have cut edges $(i,j)$ with $a_{ij}+\eps e_{ij}=0$ and for the remaining edges we have that $P^+ E$ is a multiple of $P^+G$ provided that $P^+G\ne 0$.

\subsubsection*{Reduced functional.}
Since $G=\Lap^*(vv^\top)$ of \eqref{grad-cut} depends only on $v$ (and not on $E$), we reduce the functional \eqref{F-eps-cut} to a functional of $v$ only:
\begin{equation}
\label{F-red-cut}
\wt \F_\eps(v) = \tfrac12\,v^\top \Lap (A+\eps E) v,
\quad\text {where}\quad 
P^+E = -\eta {P^+\Lap^*(vv^\top)}
\end{equation}
and the other edges are cut: $a_{ij}+\eps e_{ij}=0$ \ for \,$(i,j)\in \cE^0 = \cE \setminus \cE^+$, i.e.
with the complementary projection $P^0=I-P^+$,
\begin{equation}
\label{P-zero-cut}
P^0 E = - \eps^{-1}\,P^0 A.
\end{equation}
The normalizing factor $\eta$ is chosen such that $\|E\|_F^2=\|P^+ E\|_F^2+\|P^0 E\|_F^2=1$, i.e.
\begin{equation}
\label{eta-cut}
\eta =  \frac{(1-\|\eps^{-1} P^0 A\|_F^2)^{1/2}}{\|P^+{\Lap}^*(vv^\top)\|_F} .
\end{equation}
\index{gradient!reduced}
\begin{lemma}[Reduced gradient]
    \label{lem:red-grad-cut}
Along a path $v(t)\in\R^n$, we have (omitting the argument $t$ on the right-hand side)
$$
\frac{d}{dt}\, \wt \F_\eps(v(t)) =  \langle \Lap(A+\eps E)v,\dot v \rangle,
$$
where $E=P^+E+P^0E$ is given by \eqref{F-red-cut}--\eqref{eta-cut}.
\end{lemma}

\begin{proof}
We have
$$
\frac{d}{dt}\, \wt \F_\eps(v(t)) =  \langle \Lap(A+\eps E)v,\dot v \rangle + \tfrac12 \eps \, v^\top \Lap(\dot E)v,
$$
and a direct calculation shows that the last term vanishes.
\qed
\end{proof}

\subsubsection*{Constrained reduced gradient flow: diffusion on the graph.}
\index{graph!diffusion}
With the constraints $\|v\|^2=1$ and $\langle \one, v\rangle =0$ we obtain the constrained reduced gradient flow, which
can be viewed as a diffusion equation on the graph,
with the negative graph Laplacian in place of the Laplacian on a domain:
\begin{align}
\label{red-ode-v-cut}
&\dot v = - \Lap(A+\eps E)v + \langle \Lap(A+\eps E)v,v\rangle  v 
\\
\nonumber
&\text{with $E=P^+E + P^0E$ of \eqref{F-red-cut}-\eqref{eta-cut}}.
\end{align}
Note that here $E$ is a function of $v$, which further depends on the decomposition $\cE=\cE^+\dot\cup\,\cE^0$ of the
edge set into active and cut edges. If $\|v\|^2=1$ and $\one^\top v=0$, then $\langle v,\dot v\rangle =0$ and also
$\langle 1,\dot v \rangle =0$, since 
$\langle 1,\Lap(A+\eps E)v \rangle= \langle \Lap(A+\eps E)1,v \rangle =0$. 
Hence, the two constraints are conserved.
In a stationary point, $v$ is an eigenvector of $\Lap(A+\eps E)$.

\bng
Note that, as is common to the previously considered gradient system-based methods, the inner iteration is not guaranteed to find a global minimum, but only a local one. An appropriate choice of initial values, as discussed in Section~\ref{subsec:init}, and the fact that the gradient system is solved with different values of $\eps$ as needed in the outer iteration, help to mitigate this problem.
\eng



\subsection{Constrained minimum cut problems}
\index{graph!constrained minimum cut}
We add non-negative terms $C(v)$ 
to the functional $\wt \F_\eps(v)$ of \eqref{F-red-cut} such that the augmented functional takes the minimum value zero if and only if the graph with the modified weighted adjacency matrix $A+\eps E$ is disconnected and the imposed constraints are satisfied. In the following we consider must-link, cannot-link, membership and cardinality constraints. They can be trivially combined in our approach by just including them in the constraint functional $C(v)$ and the augmented functional
\begin{equation}\label{F-red-aug-cut}
\wh \F_\eps(v) = \wt \F_\eps(v) + C(v),
\end{equation}
which is minimized under the previously considered constraints $\|v\|=1$, $\langle \one,v \rangle =0$ and $A+\eps E\ge 0$, where $E=E(v)$ is still defined by \eqref{F-red-cut}-\eqref{eta-cut}. The minimization of this augmented functional is done in the same way as in the previous subsection, using a two-level approach that minimizes $\wh \F_\eps(v)$ for a given $\eps$ in the inner iteration and determines the optimal $\oeps$ in the outer iteration as the smallest $\eps>0$ such that $\wh \F_\eps(v(\eps))=0$.
In the inner iteration we have, instead of \eqref{red-ode-v-cut}, the gradient flow
\begin{align}
\nonumber
&\dot v = - \Lap(A+\eps E)v - \nabla C(v) + \langle \Lap(A+\eps E)v+\nabla C(v),v\rangle  v 
+ \frac1n \langle \nabla C(v),\one \rangle \doubleone
\\
&\text{with $E=P^+E + P^0E$ of \eqref{F-red-cut}-\eqref{eta-cut}}.
\label{red-ode-v-con-cut}
\end{align}
We  numerically integrate this differential equation for the single vector $v(t)$ into a stationary point. This just requires computing matrix--vector products and inner products of vectors but no computations of eigenvalues and eigenvectors of varying matrices.
In this way we aim to determine
\begin{equation} \label{stat-eps-cut}
v(\eps)=\arg\min_{v} \;\wh \F_\eps(v) \quad\text{and the associated perturbation matrix}\  E(\eps),
\end{equation}
where the minimum is taken over all $v\in \R^n$ of Euclidean norm 1 that are orthogonal to the vector $\one$.

\subsubsection*{Must-link constraints.}
\index{graph partitioning!must-link constraint}
It is required that  pairs of vertices in a given set $\mathcal{P}\subset \mathcal{V} \times \mathcal{V}$ are in the same connected component. 
Motivated by the special form of the eigenvectors as given in Fiedler's theorem (Theorem~\ref{thm:fiedler}), we consider the non-negative functional
\begin{equation} \label{F-must-link}
C(v) = \frac \alpha 2 \!\sum_{(i,j)\in\mathcal{P}} (v_i - v_j )^2 ,
\end{equation}
where $\alpha>0$ is a scaling factor chosen e.g.~as $\alpha=\|A\|_2$. 
The augmented functional \eqref{F-red-aug-cut} is still non-negative, and it is zero if and only if 
$$
\text{(i) $\:\wt \F_\eps(v)=0\quad$ and \quad(ii) $\: C(v)=0$.}
$$
Since $v$ is constrained to be orthogonal to the eigenvector $\one$ that corresponds to the smallest eigenvalue $0$ of $A+\eps E$, (i) implies that $v$ is an eigenvector to the second-smallest eigenvalue $\lambda_2(A+\eps E)$, and this eigenvalue must be $0$. By Fiedler's theorem, the graph with the weighted adjacency matrix $A+\eps E$ is disconnected, and (assuming that $\lambda_3(A+\eps E)> 0$) the entries of $v$ take only two different values of different sign, which indicate the component of the disconnected graph to which the node $i$ corresponding to the entry $v_i$ belongs. Since (ii) implies that for each pair $(i,j)\in\mathcal{P}$ we have $v_i=v_j$, the nodes $i,j$ are in the same component, as is required by the must-link constraint.

\subsubsection*{Cannot-link constraints.}
\index{graph partitioning!cannot-link constraint}
Here, pairs of vertices in a given set $\mathcal{P}\subset \mathcal{V} \times \mathcal{V}$ must be in different connected components. We consider the non-negative functional
\begin{equation} \label{F-cannot-link}
C(v) = \frac \alpha 2 \!\sum_{(i,j)\in\mathcal{P}} \min(0,v_iv_j)^2 ,
\end{equation}
where again $\alpha>0$ is a scaling factor. By the same arguments as before, we conclude that a graph with a perturbed adjacency matrix $A+\eps E$ is disconnected and satisfies the cannot-link constraints if and only if the augmented functional \eqref{F-red-aug-cut} with \eqref{F-cannot-link} is zero.

\subsubsection*{Membership constraints.}
\index{graph partitioning!membership constraint}

A given set of vertices $\mathcal{V}^+\subset\mathcal{V}$ is required to be in one connected component and another given set of vertices $\mathcal{V}^-\subset\mathcal{V}$ must be in the other connected component.
We consider the functional
\begin{equation} \label{F-member}
C(v)= \frac\alpha 2 \sum_{i\in \V^-} (v_i - \langle v^- \rangle )^2 +
\frac\alpha 2 \sum_{i\in \V^+} (v_i - \langle v^+ \rangle )^2,
\end{equation}
with a scaling factor $\alpha>0$. Here again, $v^+=(\max(v_i,0))$ and $v^-=\min(v_i,0)$ collect the positive and negative entries of $v$, and 
$\langle v^\pm \rangle$ is the arithmetic mean of the nonzero entries in $v^\pm$. By the same arguments as before, we conclude that a graph with a perturbed adjacency matrix $A+\eps E$ is disconnected and satisfies the membership constraints if and only if the augmented functional \eqref{F-red-aug-cut} with \eqref{F-member} is zero.

\subsubsection*{Cardinality constraints.}
\index{graph partitioning!cardinality constraint}
Here, each of the connected components has a prescribed minimum number $\overline n$ of vertices.
We use the same functional $C$ as in the membership constraint, except that the sets $\V^-$ and $\V^+$ are not given {\it a priori}, but are chosen depending on $v$ in the following way: $\V^-$ and $\V^+$ collect the indices of the smallest and largest $\overline n$ components of $v$, respectively, augmented by those indices for which the components of $v$ do not differ by more than a threshold $\vartheta$ from the average of the smallest and largest $\overline n$ components, respectively.

\subsection{Outer iteration}
With $v(\eps)$ of \eqref{stat-eps-cut} we consider the scalar function
$$
\varphi(\eps)=\wh \F_\eps(v(\eps))
$$
and compute the smallest zero $\oeps>0$ of $\varphi$. By construction, $A+\oeps E(\oeps)$ is the nearest weighted adjacency matrix  on the same graph as $A$ for which the Laplacian $L(A+\oeps E(\oeps))$
has a double eigenvalue $0$, and $v(\oeps)$ is the corresponding eigenvector orthogonal to $\one$. By Fiedler's theorem,
the undirected weighted graph with adjacency matrix $A+\oeps E(\oeps)$ is disconnected and, provided that $0$ is not a triple eigenvalue, the components of the graph can be read off from the signs of the entries of the Fiedler vector $v(\oeps)$. To compute $\oeps$, we use a   Newton--bisection method as in previous chapters. Here we use the following expression for the
derivative of $\varphi$.
\index{Newton--bisection method}

\begin{lemma}[Derivative for the Newton iteration] 
On an interval $(\underline\eps,\overline\eps)$ with $\overline\eps\le\oeps$ on which the set $\cE^0$ of cut edges does not change, we have

$$
 \varphi'(\eps) = - \tfrac12 \eta(\eps)^{-1}
$$
with $\eta$ of \eqref{eta-cut}.
\end{lemma}

\begin{proof}
We have (omitting the argument $\eps$ wherever this does not impair clarity)
\begin{align*}
 \varphi'(\eps) &= \frac d{d\eps} \Bigl( \tfrac12 v(\eps)^\top \Lap(A+\eps E(\eps)) v(\eps) + C(v(\eps)) \Bigr)
\\
&= \tfrac12 v^\top \Lap\bigl( \tfrac d{d\eps}(A+\eps E(\eps)) \bigr) v 
+ \langle L(A+\eps E)v,v'\rangle + \langle \nabla C(v),v'\rangle
\\
&=
\tfrac12 \langle vv^\top,\Lap\bigl( \tfrac d{d\eps}(A+\eps E(\eps)) \bigr) \rangle
+ \langle L(A+\eps E)v+\nabla C(v),v'\rangle.
\end{align*}
In the stationary point $v(\eps)$ of \eqref{red-ode-v-con-cut}, $L(A+\eps E)v+\nabla C(v)$ is a linear combination of $\one$ and $v$, and so the last term is a linear combination of $\langle \one,v'\rangle=0$, since $v$ is orthogonal to $\one$, and $2\langle v,v'\rangle =\frac d{d\eps}\|v\|^2=0$, since we normalized $\|v\|=1$. Hence the last term vanishes.
So we have, using that $\eps P^0 E= -P^0 A$ and $P^+E = -\eta P^+\Lap^* (vv^\top)$ with $\eta$ of \eqref{eta-cut},
\begin{align*}
2\, \varphi'(\eps) &= \langle \Lap^* (vv^\top), (A+\eps E(\eps))' \rangle =
\langle \Lap^* (vv^\top), (P^+A+\eps P^+E(\eps))' \rangle 
\\
&= \langle \Lap^* (vv^\top), P^+E + \eps P^+ E' \rangle 
= \langle P^+\Lap^* (vv^\top), P^+E+\eps P^+ E' \rangle
\\
&= -\eta^{-1}\Bigl(\langle P^+E, P^+E \rangle + \eps \langle P^+E, P^+ E' \rangle\Bigr).
\end{align*}
We have $\|P^+E\|_F^2 = 1- \eps^{-2}\|P^0A\|_F^2$ by \eqref{P-zero-cut} and hence
\begin{align*}
\eps \langle P^+E, P^+ E' \rangle &= \eps \,\tfrac12 \, \tfrac d{d\eps} \|P^+E\|_F^2
=\eps \,\tfrac12 \, \tfrac d{d\eps}\bigl( 1- \eps^{-2}\|P^0A\|_F^2 \bigr) = \eps^{-2}\|P^0A\|_F^2,
\end{align*}
so that 
$$
\|P^+E\|_F^2 + \eps \langle P^+E, P^+ E' \rangle = 1.
$$
This yields $ \varphi'(\eps)=-\tfrac12\eta(\eps)^{-1}$ as stated.
\qed
\end{proof}

\begin{figure}[ht]
\centering
\vspace{-4.0cm}
\includegraphics[width=0.8\textwidth]{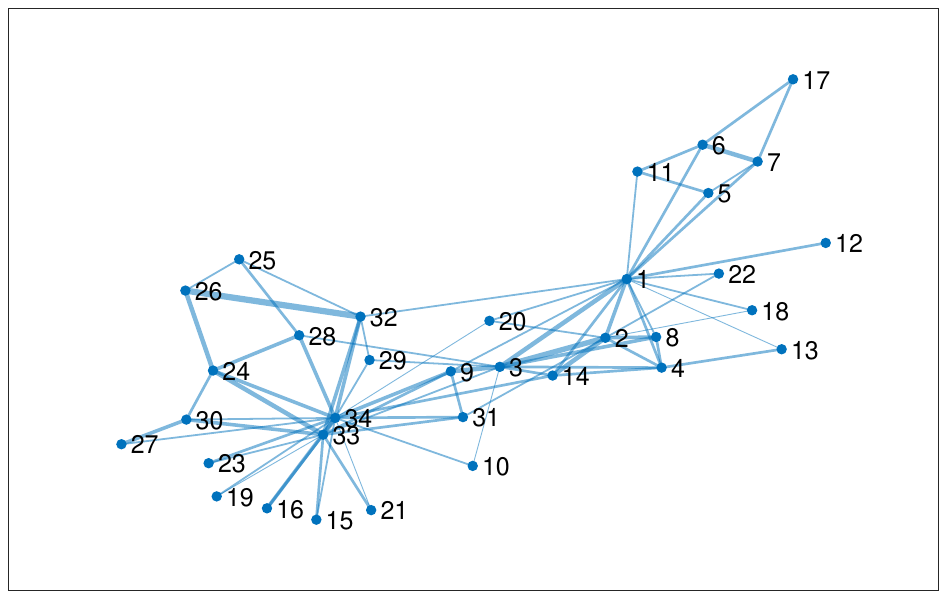} 
\vspace{-3cm}
\caption{The graph associated with Zachary's karate club.}
\label{fig:exKara1}
\end{figure}

\begin{figure}[h!]
\centering
\vspace{-3.3cm}
\includegraphics[width=0.7\textwidth]{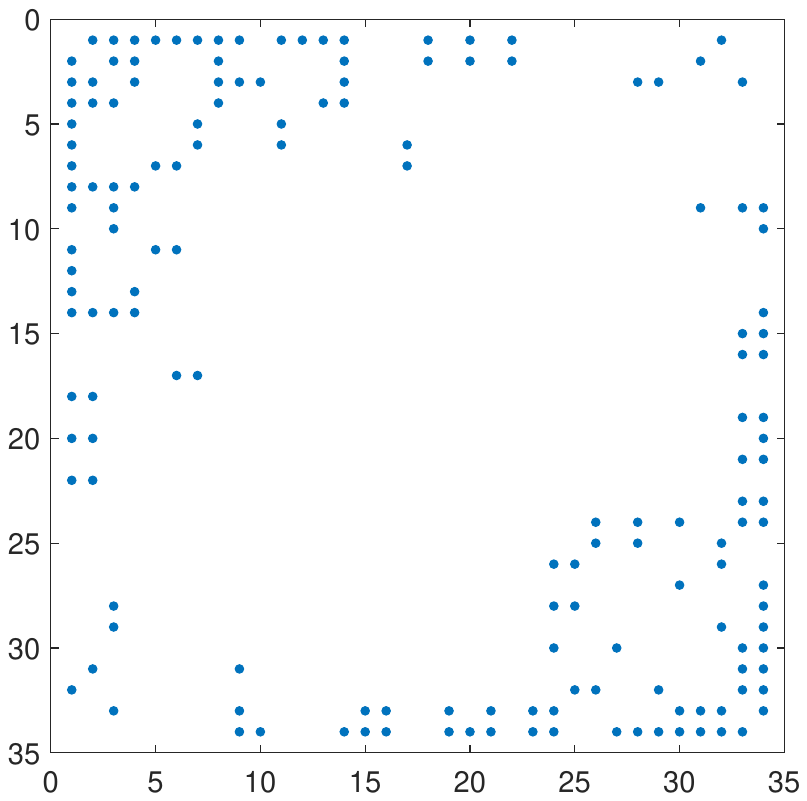} 
\vspace{-2.3cm}
\caption{Sparsity pattern of weighted adjacency matrix associated with Zachary's karate club.}
\label{fig:PN1-Kara}
\end{figure}

\begin{figure}[ht!]
\centering
\vspace{-2cm}
\hspace{-1.9cm}
\includegraphics[width=0.6\textwidth]{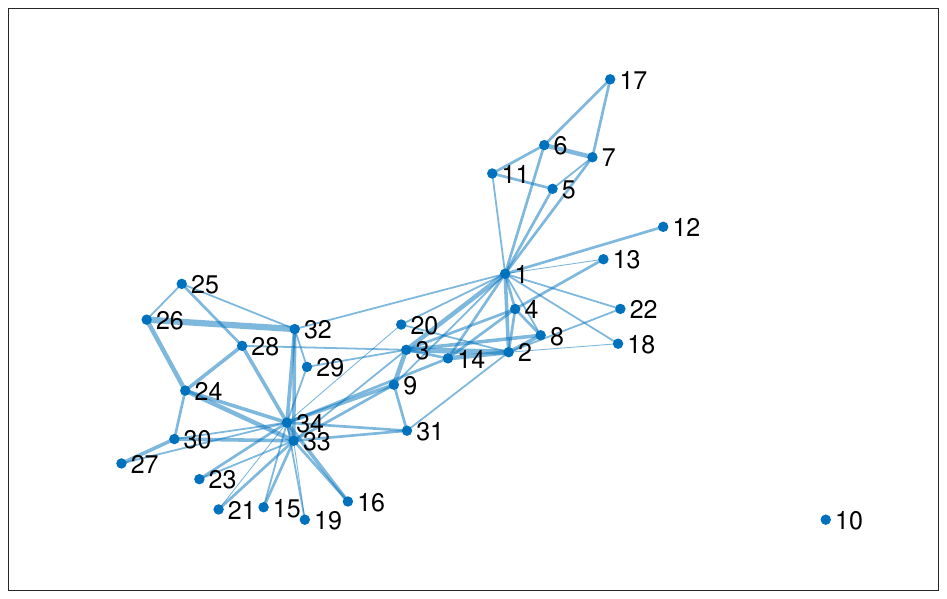} 
\hspace{-1.2cm}
\includegraphics[width=0.6\textwidth]{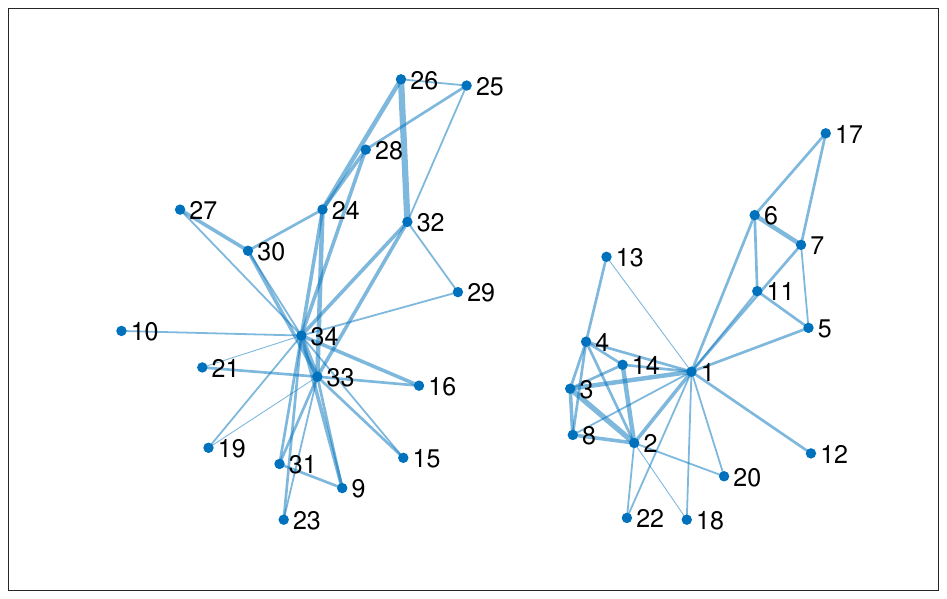} 
\vspace{-2.6cm}
\caption{Left: unconstrained minimum cut. Right: cardinality-constrained minimum cut, both for the graph associated with Zachary's karate club.}
\label{fig:exKara2}
\end{figure}

\begin{remark}
    In order to decide where to cut the graph, it is not necessary to compute $\oeps$ with high accuracy. The iteration can be stopped if, for a given tolerance $\vartheta$, the entries of the vector $v(\eps_k)$ take only values within two disjoint intervals of width $\vartheta$. This marks the membership to the two connected components.
\end{remark}


\subsection{Illustrative examples}

We consider a classical test problem of small dimension (Zachary's karate club). This weighted graph consisting of $34$ vertices; weights describe the relationship between the $34$ members of the karate club; see Zachary (\cite{Zac77}).

The weights vary in the integer range $[1,7]$. The graph is illustrated in Figure 
\ref{fig:exKara1} where the thickness of the
edges increases for larger weights.

The unconstrained minimum cut algorithm determines the graph represented in Figure \ref{fig:exKara2}, where the graph is cut into two components, one made by the singleton node $\{10\}$ and the other made by the rest of the graph; the value $\oeps = \sqrt{10}$, which corresponds to the $2$-norm of the $10$-th row of the adjacency matrix
multiplied by $2$.

Evidently the result is disappointing,
although very predictable, because it does not express a significant clustering of the two communities which attend the Karate club.
Note that alternatively also nodes $\{18\}$ and $\{19\}$
might have been cut out since they would have required a cut with the same value.

In order to obtain a meaningful clustering we make use of a cardinality constraint, imposing that each of the two components should have a cardinality of at least $40\%$, that is, the dimension of each cluster should be at least 
$\lceil 0.4\,n \rceil = 14$, where $n=34$.

The result is shown in Figure \ref{fig:exKara2} on the right-hand side,
which shows two groups of $n_1 = 18$ and $n_2 = 16$
members, respectively.
The result has been obtained with a penalty parameter $\alpha=1$, to which corresponds the
value $\oeps = \sqrt{120} = 10.95445115..$.

\subsubsection{Practical strategy}

 Numerical experiments seem to indicate that when an active edge is \emph{deleted} by our algorithm, which means it is set to $0$, then it remains $0$ without reactivating. This suggests an alternative simpler strategy which makes the set of active vertices monotonically decreasing as $\eps$ increases.

 A second practical trick consists in terminating the clustering in a combinatorial way. This means that when the functional reaches a value below a certain threshold $\theta$, i.e.
 \[
 \wh \F_\eps(v)  \le \theta
 \]
 we proceed as follows.
 We set $B \approx A + \eps E$, where the approximation is obtained, for example, by rounding to the closest integer matrix if the weights are integer or similarly for matrices whose weights are expressed with a $2$-digit fractional part.
 If $\Lap(B) = 0$ then $\oeps = \| B - A \|_F$.

\begin{figure}[ht]
\centering
\vspace{-3.0cm}
\includegraphics[width=0.80\textwidth]{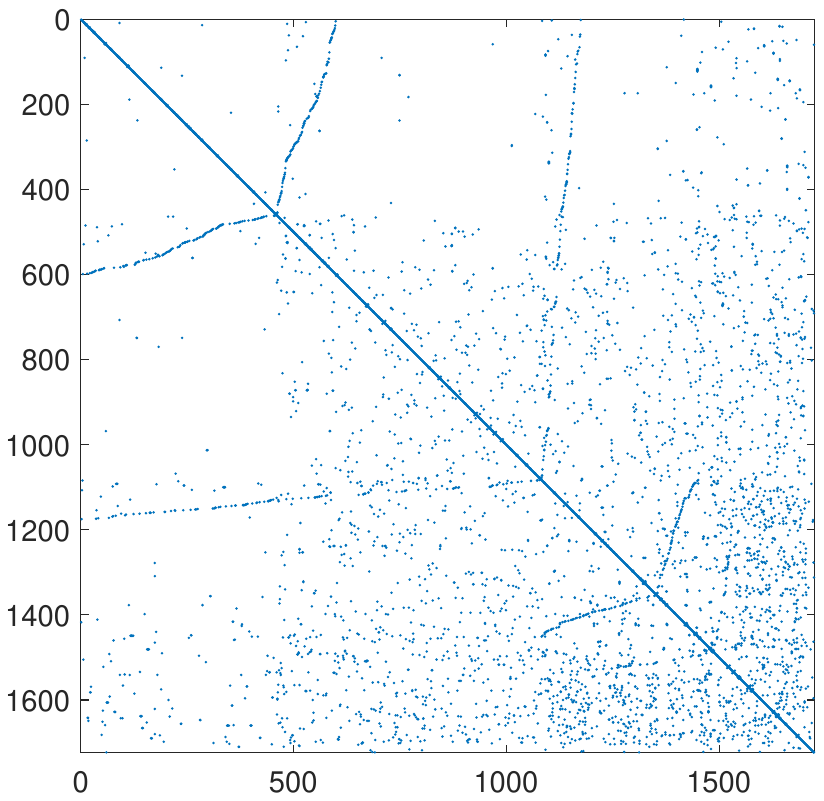} 
\vspace{-3cm}
\caption{Sparsity pattern of matrix \emph{BCSPWR09}.}
\label{fig:PN1}
\end{figure}

 \subsubsection{A larger test example}

We consider the graph represented by the weighted adjacency matrix \emph{BCSPWR09},
collected by B. Dembart and J. Lewis, Boeing Computer Services, Seattle, WA, USA.

It consists of a $1723 \times 1723$
sparse symmetric matrix, whose sparsity pattern is illustrated in Figure \ref{fig:PN1}. 
All weights are equal to $1$.

\begin{figure}[ht]
\centering
\vspace{-4.0cm}
\includegraphics[width=\textwidth]{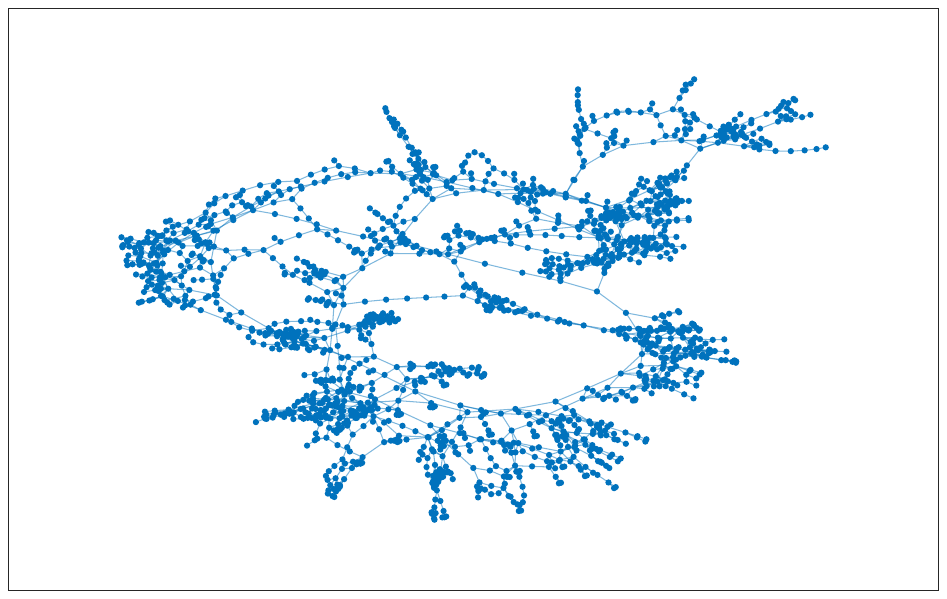} 
\vspace{-4.6cm}
\caption{The graph represented by matrix \emph{BCSPWR09}.}
\label{fig:PN2}
\end{figure}

The graph is illustrated in Figure \ref{fig:PN2}.
 \begin{figure}[h!]
\centering
\vspace{-4.0cm}
\includegraphics[width=\textwidth]{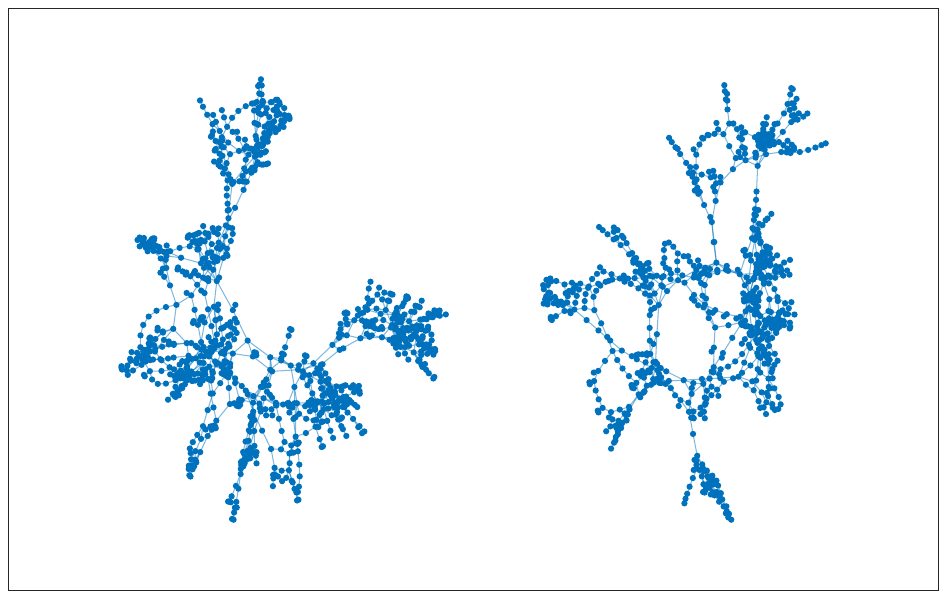} 
\vspace{-4.6cm}
\caption{The obtained partition of the graph \emph{BCSPWR09}, with 10 edges cut.}
\label{fig:PN3}
\end{figure}

Setting a cardinality constraint of $40\%$, that is asking for cutting the graph such that each connected subgraph has at least 40\% of the nodes, we obtain the two clusters illustrated in Figure \ref{fig:PN3}.
The graph is illustrated in Figure \ref{fig:PN2}.

The two components have cardinality $n_1 = 758$ and $n_2 = 965$ and the computed distance is
$\oeps = \sqrt{20} \approx 4.47213595$.
This means that $10$ edges have been cut to provide
the splitting.

\section{Toppling the ranking}
\label{sec:ranking}
\index{graph!eigenvector centrality}
\index{graph!ranking of nodes}
\index{Perron--Frobenius theorem}
Which nodes in a graph are important? How to rank their importance? These questions are addressed by the notion of centrality, of which there are several variants. Here we consider {\it eigenvector centrality}, which is based on the Perron--Frobenius theorem. The weighted adjacency matrix $A\ge 0$ (with entrywise non-negativity) has a largest, positive real eigenvalue with a non-negative eigenvector $x=(x_i)\ge 0$. The nodes are ranked according to the entries of this eigenvector:
$$
\text{Node $i$ is ranked higher than node $j$ if $x_i>x_j$.}
$$
In this section we address a basic question: How robust is the ranking under perturbations of the weights $a_{ij}$ in the graph?

\begin{example} \label{ex:ranking} Assume that node 1 is ranked higher than node 2, i.e. $x_1>x_2$. What is the Frobenius-norm distance to the nearest weighted adjacency matrix $\wt A\ge 0$ on the same graph such that the ranking of nodes 1 and 2 is reversed, i.e. the Perron--Frobenius eigenvector $\wt x\ge 0$ of $\wt A$ has $\wt x_1\le \wt x_2$?
\end{example}

We address such questions for {\em undirected} graphs (i.e. symmetric $A$) using again computationally inexpensive Rayleigh quotients instead of eigenvalues and eigenvectors and their derivatives. The latter could, however, be used in a similar way also for directed graphs, which cannot be done with Rayleigh quotients.

\subsection{Robustness of ranking as a matrix nearness problem} \label{subsec:robust-rk}
We consider an undirected graph with vertex set $\mathcal{V}=\{1,\dots,n\}$ and edge set $\wh\cE$,
and an undirected subgraph with the same vertex set but a possibly reduced edge set $\cE\subseteq\wh\cE$. We let 
$\wh\cS=\mathrm{Sym}(\wh\cE)$ and
$\cS=\mathrm{Sym}(\cE)$ be the spaces of symmetric real matrices with sparsity pattern given by the edge sets $\wh\cE$ and $\cE$, respectively. 
Let $A\in\wh\cS$ with $A\ge 0$ be a weighted adjacency matrix of the graph $(\mathcal{V},\wh\cE)$. 

Let $u\in\R^n$ be a given nonzero vector. We require that $u$ has entries with positive and negative sign.
For Example~\ref{ex:ranking} we would take $u=e_1-e_2$, where $e_i$ is the $i$th standard unit vector. We consider the following matrix nearness problem.
\index{eigenvector optimization}

\medskip\noindent
{\bf Problem.} 
{\em Find a matrix $\Delta\in\cS=\mathrm{Sym}(\cE)$ of minimal Frobenius norm such that ${A+\Delta\ge 0}$ and the Perron--Frobenius eigenvector $x\ge 0$ of $A+\Delta$ is orthogonal to the given vector $u$.}

\medskip\noindent
Note that for Example~\ref{ex:ranking} the orthogonality constraint is 
$$
0=u^\top x = (e_1-e_2)^\top x = e_1^\top x - e_2^\top x = x_1-x_2,
$$
i.e. the orthogonality imposes that nodes 1 and 2 are ranked equally, as opposed to the original ranking $x_1>x_2$. Choosing $u=1.01 \,e_1 - e_2$ would reverse the ranking.

We will propose a two-level approach that does not use eigenvectors and not even eigenvalues in the functional to be minimized. This approach relies on the following lemma.
\index{Rayleigh quotient}

\begin{lemma} [Rayleigh quotients under an orthogonality constraint]\label{lem:u-orth-v-rk}
    Let $A\in\R^{n,n}$ be a symmetric matrix, and let the nonzero vector $u\in\R^n$ be normalized to have Euclidean norm~1. Let
    $v_\star\in\R^n$ maximize the Rayleigh quotient $v^\top A v\,/\,v^\top v$ over all nonzero vectors $v\in\R^n$ that are orthogonal to $u$:
    $$
    v_\star^\top A v_\star = \max_{\| v\|=1, u^\top v = 0} v^\top A v, \qquad\quad \| v_\star\|=1,\ \ u^\top v_\star =0.
    $$
    Then, the following holds true:

    (a) $v_\star$ is an eigenvector to the largest eigenvalue $\lambda_\star = v_\star^\top A v_\star$ of $P_u^\perp A P_u^\perp$, where $P_u^\perp=I-P_u$ and $P_u=uu^\top$ is the orthogonal projection onto the span of $u$.

    (b) $v_\star$ is an eigenvector of $A$ if and only if 
    $$
    u^\top A v_\star=0.
    $$
    In this case, $v_\star$ or $u$ is an eigenvector to the largest eigenvalue of $A$.
\end{lemma}

\begin{proof}
(a) follows directly from the fact that $P_u^\perp v=v$ for every $v$ orthogonal to $u$. For part (b), after changing to an orthonormal basis of $\R^n$ that has $u$ as the first basis vector, we may assume $u=e_1$. We partition the matrix $A$ and the vector $v$ orthogonal to $e_1$ as
$$
A = \begin{pmatrix}
    \alpha & a^\top \\
    a & \check A
\end{pmatrix}, \quad\ v= \begin{pmatrix}
    0 \\
     \check v
\end{pmatrix}, \qquad \text{so that} \quad Av=\begin{pmatrix}
    a^\top  \check v\\
     \check v
\end{pmatrix}.
$$
Hence, $Av=\lambda v$ if and only if 
$$
0=a^\top \check v = e_1^\top Av = u^\top A v \quad\text{ and }\quad \check A \check v =\lambda \check v.
$$
The latter equation is equivalent to $P_u^\perp A P_u^\perp v = \lambda v$, since for $u=e_1$,
$$
P_u^\perp A P_u^\perp = \begin{pmatrix}
    0 & 0 \\
    0 & \check A
\end{pmatrix}.
$$
By the Courant--Fischer minimax principle, 
\index{Courant--Fischer minimax principle}
the second-largest eigenvalue $\lambda_2$ of $A$ equals
$$
\lambda_2 = \min_{\mathcal{U}_1} \max_{0\ne v \perp \mathcal{U}_1} \frac{v^\top Av}{v^\top v},
$$
where the minimum is over all one-dimensional subspaces $\mathcal{U}_1$ of $\R^n$ and is attained only for the span of an eigenvector to the largest eigenvalue of $A$. With the particular choice $\mathcal{U}_1=\text{span}(u)$
we obtain, unless $u$ is an eigenvector to the largest eigenvalue of $A$,
$$
\lambda_2 < \max_{0 \ne v \perp u} \frac{v^\top Av}{v^\top v}  = v_\star^\top A v_\star = \lambda_\star,
$$
so that $\lambda_\star$ is the largest eigenvalue of $A$ if it is some eigenvalue of $A$.
\qed    
\end{proof}

\subsection{Two-level iteration} 
\index{two-level iteration}
Let $A\in\wh\cS=\sym(\wh\cE)$ with $A\ge 0$ and $u\in\R^n$ of norm~1 be given.
The choice of the functional in the following two-level approach is motivated by Lemma~\ref{lem:u-orth-v-rk}. 
Here, $\lambda_{\max}(M)$ is the largest eigenvalue of a symmetric matrix $M$.

\begin{itemize}
\item {\bf Inner iteration:\/}
Given $\eps>0$, we aim to compute a  matrix $E \in\cS=\sym(\cE)$ 
of Frobenius norm~1 with $A+\eps E\ge 0$ and a vector $v\in\R^n$ of norm~1 and orthogonal to $u$  such that 
\begin{equation}\label{F-eps-rk}
\F_\eps(E,v) = \bigl( \lambda_{\max}(A+\eps E) - v^\top (A+\eps E) v \bigr) +  \bigl( u^\top (A+\eps E)v \bigr)^2
\end{equation}
is minimized under the constraints $\|E\|_F=1$, $\|v\|=1$, $u^\top v =0$ and $A+\eps E\ge 0$. We write
\begin{equation} \nonumber
\bigl(E(\eps),v(\eps)\bigr) = \arg\min\limits_{(E,v)} \F_\eps(E,v).
\end{equation}
\item {\bf Outer iteration:\/}
 We compute $\oeps$ as the smallest $\eps>0$ such that $\F_\eps(E(\eps),v(\eps))=0$.
\end{itemize}

Lemma~\ref{lem:u-orth-v-rk} shows that if this iteration succeeds, then $\Delta=\oeps E(\oeps)\in \cS$ \bng is a matrix of locally smallest Frobenius norm  such that $A+\Delta\ge 0$ has a Perron-Frobenius eigenvector that is orthogonal to $u$, and $v(\oeps)$ (taken with the proper sign) is this eigenvector.
\eng

Since the functional $\F_\eps$ does not involve eigenvectors, the corresponding gradient flow system does not require computing derivatives of eigenvectors. It contains, however, a largest eigenvalue, and therefore the gradient system requires the computation of largest eigenvalues and their eigenvectors of varying matrices $A+\eps E$.

\subsubsection*{Constrained gradient flow for the inner iteration.}
Consider a smooth path of matrices $E(t)\in \cS$ and vectors $v(t)\in\R^n$, and let $x(t)$ be the eigenvector of $A+\eps E(t)$ to the largest eigenvalue, which is assumed to be a simple eigenvalue. 
We let $\sigma(t)=2 u^\top (A+\eps E)v(t)$.
Then we find, omitting the argument $t$ on the right-hand side,
\begin{align}
\nonumber
    \frac{d}{dt}\F_\eps(E(t),v(t))&=
    \langle \Pi^\cS(xx^\top - vv^\top + \sigma uv^\top), \eps \dot E \rangle
    \\
    &\quad + \langle -2(A+\eps E)v + \sigma (A+\eps E)u, \dot v \rangle.
\label{F-dot-rk}
\end{align}
The gradient flow subject to the constraints $E\in\cS$, $\|E\|_F=1$, $\|v\|=1$ and $u^\top v=0$ therefore becomes
\begin{equation} \label{Ev-ode-rk}
\begin{aligned} 
    \eps \dot E &= - G + \langle G,E\rangle E \qquad\qquad\text{with}\quad
     G= \Pi^\cS(xx^\top - vv^\top + \sigma\, uv^\top) 
    \\
    \dot v &= -g + \langle g,v\rangle v + \langle g,u\rangle u
    \quad\text{with}\quad
    g = -2(A+\eps E)v + \sigma (A+\eps E)u  .
\end{aligned}
\end{equation}
For simplicity, we do not include the nonnegativity constraint $A+\eps E\ge 0$, which can be enforced in the same way as in the previous section. 

\subsubsection*{An approach without computing eigenvalues and eigenvectors.}
\index{eigenvector optimization!without eigenvectors}
We modify the above approach, using a Rayleigh quotient instead of the largest eigenvalue.  
We replace the eigenvalue $\lambda_{\max}(A+\eps E)$ by a Rayleigh quotient $x^\top (A+\eps E) x$ (with $\|x\|=1$), where $x(t)$ satisfies the gradient system for maximizing $x^\top (A+\eps E) x$ over all $x\in \R^n$ of norm~1 without orthogonality constraint:
with a  time-scaling factor $\rho>0$, we let
    \begin{equation} \label{x-ode-rk}
    \rho \dot x = (A+\eps E) x - \langle (A+\eps E) x, x \rangle x,
    \end{equation}
ideally starting from the Perron--Frobenius eigenvector to $A+\eps E_0$. We combine \eqref{x-ode-rk} with the differential equations \eqref{Ev-ode-rk}. We then have the same stationary point $(E(\eps),v(\eps))$ as before, and $x(\eps)$ is the Perron--Frobenius eigenvector of $A+\eps E(\eps)$.

\subsubsection*{Newton--bisection method for the outer iteration.} 
\index{Newton--bisection method}
As in previous chapters, we use a Newton--bisection method or the HEC method to compute $\oeps$. By the same arguments as in the proof of Theorem~\ref{chap:two-level}.\ref{thm:phi-derivative} we find that under analogous assumptions the univariate function $ \varphi(\eps)=\F_\eps(E(\eps),v(\eps))$ has the derivative
$$
 \varphi'(\eps)= - \| G(\eps) \|_F,
$$
where $G(\eps)$ is the gradient matrix $G$ of \eqref{Ev-ode-rk} that corresponds to the minimizer $(E(\eps),v(\eps))$.
This derivative is used in Newton's method for solving $ \varphi(\eps)=0$.

\subsubsection*{Optimal perturbation.} We show that under plausible assumptions, the matrix nearness problem stated in Section~\ref{subsec:robust-rk} has an optimal perturbation $\Delta\in\cS$ such that 
$$
\text{$\Delta$ is a real multiple of a projected rank-1 matrix $\Pi^\cS(uv^\top)$ for some $v\in\R^n$.} 
$$
Moreover, it turns out that $\pm v$ is then the Perron--Frobenius eigenvector of $A+\Delta$.

\medskip
To arrive at this result, we consider the minimization of a modification of the functional $\F_\eps$ in \eqref{F-eps-rk} 
in which we take the square root of the second term:
$$
\F_\eps(E,v) = \bigl( \lambda_{\max}(A+\eps E) - v^\top (A+\eps E) v \bigr) +  \bigl| u^\top (A+\eps E)v \bigr|.
$$
We let $\eps\nearrow\oeps$, for which the minimum of $\F_{\oeps}$ becomes zero. In the stationary point $(E(\eps),v(\eps))$ of
\eqref{Ev-ode-rk}, $E(\eps)$ is a multiple of 
$$
G(\eps)=\Pi^\cS\bigl(x(\eps)x(\eps)^\top - v(\eps)v(\eps)^\top + \sigma(\eps)uv(\eps)^\top\bigr),
$$
where $x(\eps)$ is the Perron--Frobenius eigenvector of $A+\eps E(\eps)$ and the factor
$\sigma(\eps)$ is now changed to  $\sigma(\eps)=\mathrm{sign}\bigl(u^\top (A+\eps E(\eps)) v(\eps)^\top\bigr)$. 

We assume that $\sigma(\eps)=1$ (or alternatively $-1$) for a sequence of $\eps<\oeps$ that converges to $\oeps$. We further assume that $\F_\eps(E(\eps),v(\eps))\to 0$ as $\eps\nearrow\oeps$.

By compactness, $(E(\eps),v(\eps))$ has an accumulation point, which then equals a minimizer $(E(\oeps),v(\oeps))$ of $\F_{\oeps}$, i.e., $\F_{\oeps}(E(\oeps),v(\oeps))=0$. 
At $\oeps$ we have
$v(\oeps)=x(\oeps)$, and (up to extraction of a subsequence) we have $v(\eps)\to v(\oeps)$ as $\eps\nearrow\oeps$. This implies $G(\eps)\to G_*= \pm\Pi^\cS\bigl(uv(\oeps)^\top\bigr)$, and since $E(\eps)\to E(\oeps)$ (up to extraction of a subsequence), we have that $E(\oeps)$ is a multiple of $G_*$. Since $\Delta=\oeps E(\oeps)$ is an optimal perturbation, this yields the stated result.

\subsubsection*{Example: Matrix $A_{27}$.}
\begin{figure}[ht]
\centering
\vspace{-2.5cm}
\includegraphics[width=.68\textwidth]{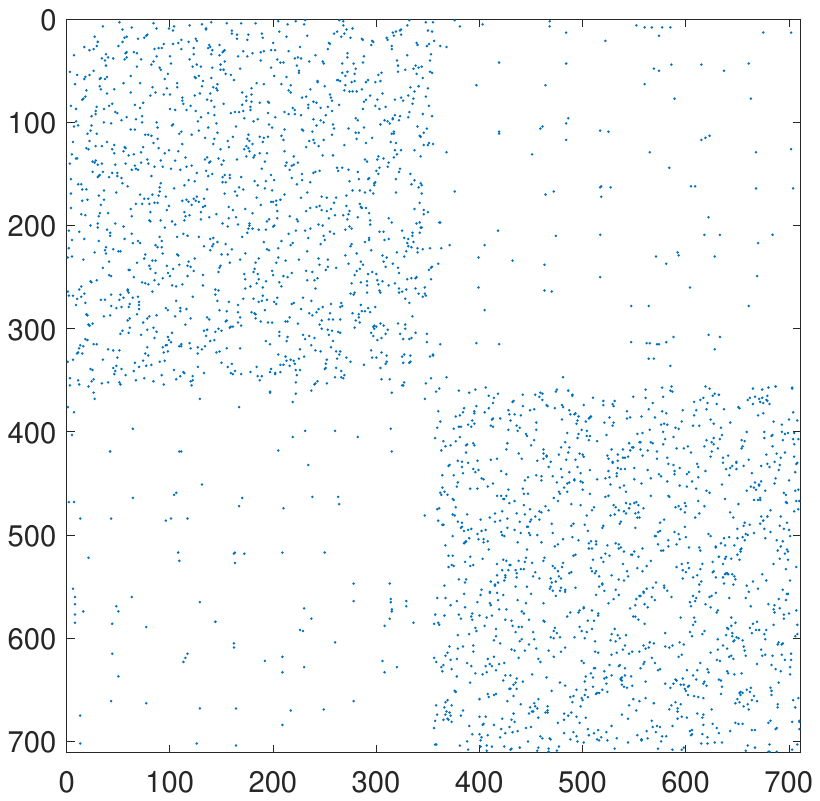} 
\vspace{-2.5cm}
\caption{Pattern of the adjacency matrix for the graph example $A_{27}$.}
\label{fig:ex1_1}
\end{figure}

The considered graph is one of a sequence of graphs arising from molecular dynamics simulations aimed at 
analyzing the structure of water networks under different temperature and pressure conditions; see Faccio, Benzi, Zanetti-Polzi \& Daidone (\cite{FacBZD22}).
It consists of $710$ nodes (corresponding to oxygen atoms), most of which ($613$) have degree $4$. Its pattern is shown in Figure \ref{fig:ex1_1}.

\begin{figure}[ht]
\centering
\vspace{-2cm}
\includegraphics[width=.48\textwidth]{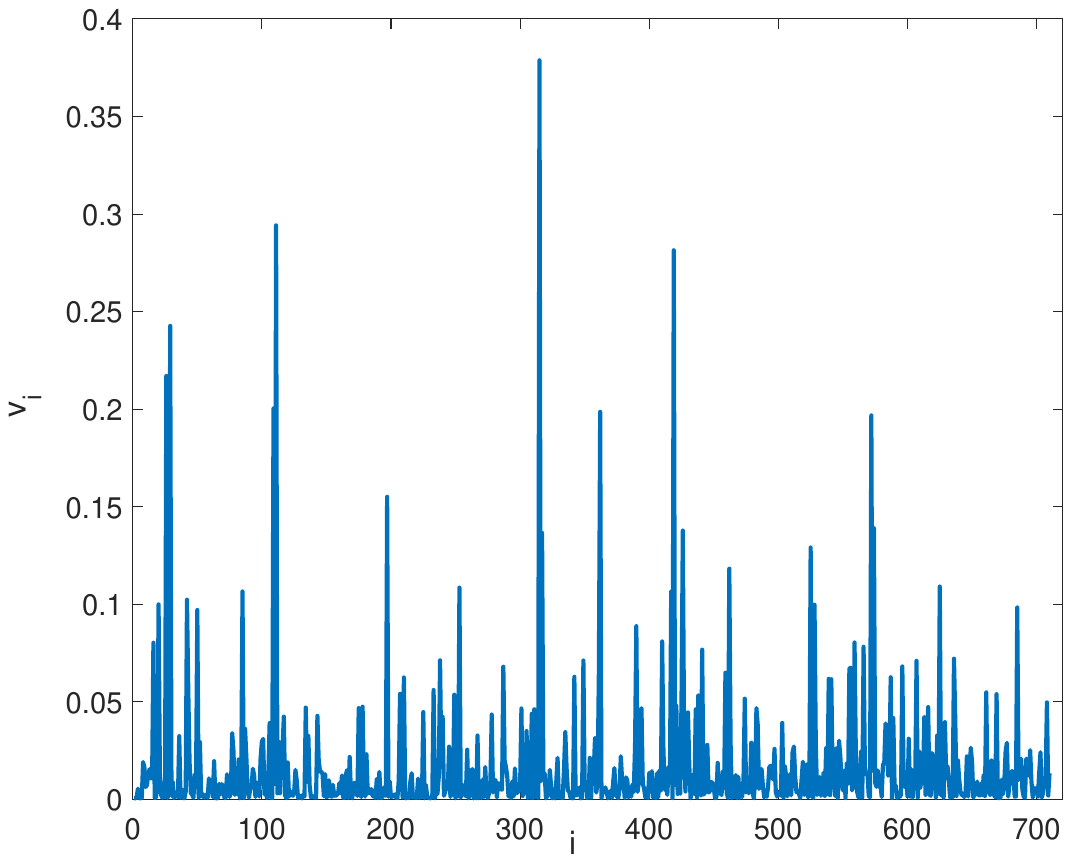} \includegraphics[width=.48\textwidth]{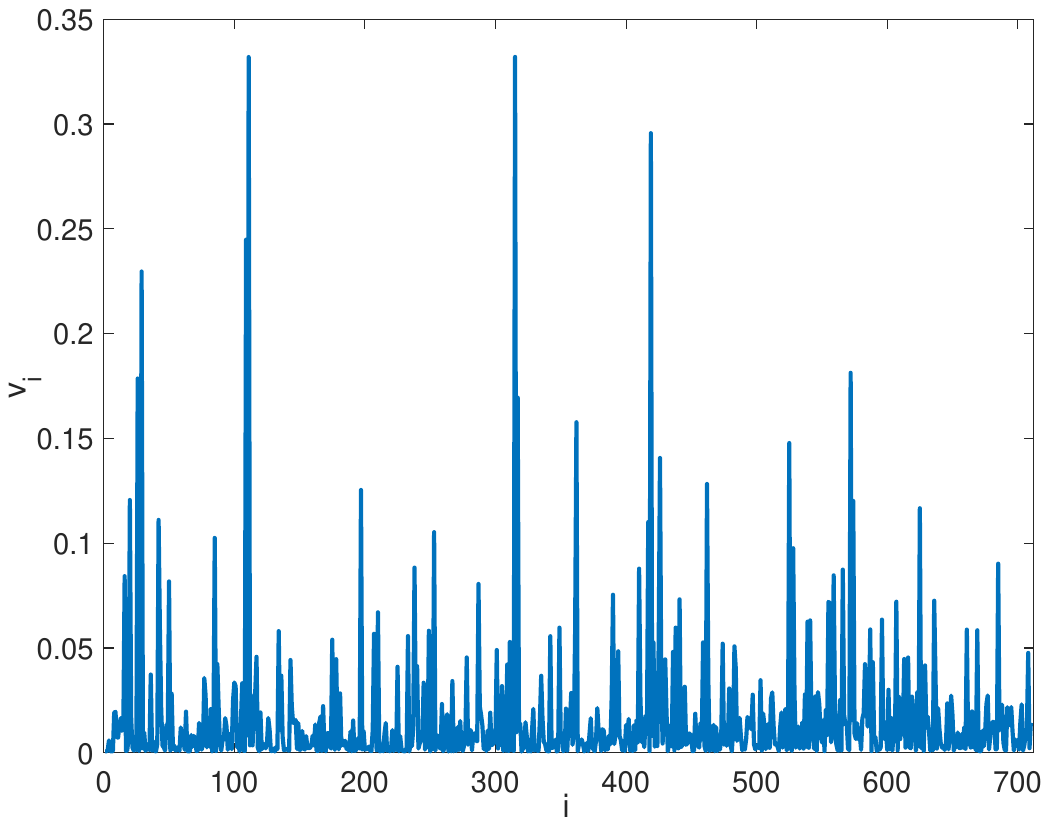}
\vspace{-1.5cm}
\caption{Distribution of the entries of Perron eigenvector for the example $A_{27}$.}
\label{fig:ex1_2}
\end{figure}

In Figure \ref{fig:ex1_2} we show the entries of the Perron eigenvector of the original (left picture) and perturbed (right picture).

The five leading entries of the Perron eigenvector of $A$ are the following:
\begin{align*}
&v_{26}  = 0.21705332, \quad
v_{29}  = 0.24269500, \quad 
v_{419} = 0.28145420, \\
&v_{111} = 0.29422909, \quad 
v_{315} = 0.37879622.
\end{align*}
After scaling the matrix $A$ to have Frobenius norm $1$, we apply the algorithm.
The norm of the perturbation which determines the coalescence of the two main entries of the Perron eigenvector is
$\eps_\star = 0.0092654$, which is below $1\%$. The target vertices are $v_{111}$ and $v_{315}$;  after the perturbation they 
are numerically equal:
\[
v_{111}=0.33101615, \qquad 
v_{315}=0.33136357.
\]
The variations of the two entries  are
$\delta_{111} = +0.1281$ and $\delta_{315} = -0.1235$.

\begin{table}[hbt]
\begin{center}
\begin{tabular}{|l|l|l|}\hline
  $k$ & $\eps_k$ & $\ophi(\eps_k)$  \\
 \hline
\rule{0pt}{9pt}
\!\!\!\! 
  $0$         & $0.005$ &  $0.000375176389506$  \\
	$1$         & $0.005136626814311$ & $0.000351719368954$ \\
	$2$         & $0.008469968038737$ & $0.000013522265277$  \\
	$3$         & $0.008489282220273$ & $0.000012884575588$  \\
	$6$         & $0.009184673253577$ & $0.000000185817345$  \\
	$10$        & $0.009212325217551$ & $0.000000093285281$   \\ 
 \hline
\end{tabular}
\vspace{2.5mm}
\caption{Computation of  $\eps_k$ and $\ophi(\eps_k)$.
 \label{tab:ex1}}
\end{center}
\end{table}

\subsection{Commuting projections: a minimalistic algorithm}
 We assume in the following that the orthogonal projection $P_u=uu^\top$ onto the span of $u$ (for $u\in\R^n$ of norm 1) and the projection $\Pi^\cE$ onto the sparsity pattern given by the edge set $\cE$ commute. The latter projection is simply, for $Z=(z_{ij})\in\R^{n,n}$,
$$
\Pi^\cE(Z) = \begin{cases}
	z_{ij}, &\text{if } (i,j)\in\cE \,,\\
	0\,, &\text{otherwise.}\end{cases}
$$
Note that the orthogonal projection $\Pi^\cS$ onto $\cS=\sym(\cE)$ as defined in \eqref{Pi-S-cut} is given by $\Pi^\cS(Z)=\Pi^\cE(\sym(Z))$,
where $\sym(Z)=\tfrac12(Z+Z^\top)$ is the symmetric part of $Z\in\R^{n,n}$.
We will assume that the projections $P_u$ and $\Pi^\cE$ commute in the following sense:
\begin{equation}
    \label{proj-commute-rk}
P_u \Pi^\cE (Z)  = \Pi^\cE  (P_u Z)    \quad\text{ and } \quad \Pi^\cE (Z) P_u = \Pi^\cE  (Z P_u).
\end{equation}
In Example~\ref{ex:ranking}, this condition is satisfied if with $(1,k)\in\cE$ for some $k$, also $(2,k)\in\cE$, and conversely. In general, this condition is satisfied if with each $(i,k)\in\cE$ with $u_i\ne 0$, also $(j,k)\in\cE$ for each $j$ with $u_j\ne 0$.
If the graph $(\mathcal{V},\mathcal{E})$ does not satisfy this condition, then we can extend the edge set in this way and allow for perturbations of weights on the extended graph with the closure of the edge set $\cE$ with respect to the fixed vector $u$:
\begin{equation}
    \label{E-closure}
    \overline{\cE} = \cE \cup \{ (j,k) \text{ and } (k,j)\,:\text{$u_j \ne 0$ and $(i,k)\in\cE$ for some $i$ with $u_i\ne 0$}\}.
\end{equation}
It will be seen that the commutativity \eqref{proj-commute-rk} leads to a substantial simplification of the problem thanks to orthogonality relations as given in the following lemma.


\begin{lemma}[Orthogonality implied by commutativity] \label{lem:orth-rk}
    Under the commutativity condition \eqref{proj-commute-rk}, we have for all $v\in\R^n$ with ${u^\top v=0}$
    $$
    \langle \Pi^\cS (uv^\top), \Pi^\cS(vv^\top)\rangle =0.
    $$
\end{lemma}

\begin{proof}
Using $P_u^\perp v =v$ and $P_u^\perp u=0$ together with \eqref{proj-commute-rk}, we find
\begin{align*}
    &\langle \Pi^\cS (uv^\top), \Pi^\cS(vv^\top)\rangle
    = \langle \Pi^\cS (uv^\top), vv^\top\rangle = v^\top \Pi^\cS (uv^\top) v
    \\
    &= v^\top P_u^\perp \Pi^\cS (uv^\top) P_u^\perp v = v^\top P_u^\perp \Pi^\cE \bigl(\tfrac12(uv^\top+vu^\top)\bigr) P_u^\perp v
    \\
    &= v^\top  \Pi^\cE \bigl(\tfrac12 P_u^\perp(uv^\top+vu^\top) P_u^\perp \bigr) v = v^\top  \Pi^\cE (0) v =0
\end{align*}
as stated.
\qed
\end{proof}

With this orthogonality we obtain the following partial solution to the problem, ignoring the non-negativity constraint on the perturbed matrix for the moment. A remarkable feature is that the optimal perturbation $\Delta$ is fully characterized by the vector $u$ and the vector $v$ that maximizes the Rayleigh quotient of the original matrix $A$ under the orthogonality constraint $u^\top v=0$. Moreover, the perturbation matrix $\Delta$ changes only entries $(i,j)\in\cE$ that have $u_i\ne0$ or $u_j\ne 0$. In Example~\ref{ex:ranking}, only entries in the first two rows and columns of $A$ are changed.

\begin{theorem}[Matrix nearness problem solved without non-negativity constraint]\label{thm:ranking} 
Let $A\in\cS=\sym(\cE)$ and $u\in\R^n$ of norm $1$ be given. Let $v\in\R^n$ be an eigenvector to the largest eigenvalue $\lambda_u(A)$ of $P_u^\perp A P_u^\perp$, i.e. $v$ maximizes the Rayleigh quotients of $A$ subject to the orthogonality constraint $u^\top v=0$.
Under the commutativity condition \eqref{proj-commute-rk}, the matrix $\Delta$ of smallest Frobenius norm 
such that an eigenvector to the largest eigenvalue of $A+\Delta$ is orthogonal to the vector $u$, is given by
$$
\Delta = -\vartheta \Pi^\cS(uv^\top) \quad\text{with}\quad \vartheta =\frac{u^\top Av}{\|\Pi^\cS(uv^\top)\|_F^2}.
$$
Moreover, $v$ is such an eigenvector to the largest eigenvalue of $A+\Delta$, which is still $\lambda_u(A)$.
\end{theorem}

\begin{proof}
For every $v\in\R^n$ that is orthogonal to $u$ and every $\vartheta\in\R$ we have
$v^\top P_u^\perp A P_u^\perp v = v^\top A v $ and
$$
v^\top (A-\vartheta \Pi^\cS(uv^\top))v = v^\top A v - \vartheta \langle \Pi^\cS (uv^\top), \Pi^\cS(vv^\top)\rangle = v^\top A v
$$
by Lemma~\ref{lem:orth-rk}. An eigenvector $v$ to the largest eigenvalue of $P_u^\perp A P_u^\perp$ is then also an eigenvector to the largest eigenvalue of $P_u^\perp (A -\vartheta \Pi^\cS(uv^\top)) P_u^\perp$ for every $\vartheta\in\R$.

On the other hand, we know from above that the matrix $\Delta$ of a given norm that minimizes $|u^\top (A+\Delta)v|$, is a multiple of $\Pi^\cS(uv^\top)$ as long as the minimum value is positive. By continuity, this remains true also for the smallest norm for which the minimum value becomes zero. We have
$$
u^\top (A-\vartheta \Pi^\cS(uv^\top))v = u^\top A v - \vartheta \langle \Pi^\cS (uv^\top), uv^\top\rangle 
= u^\top A v - \vartheta \|\Pi^\cS(uv^\top)\|_F^2,
$$
which is zero for $\vartheta$ as stated in the theorem. For this $\vartheta$, Lemma~\ref{lem:u-orth-v-rk} yields that an eigenvector $v\perp u$ to the largest eigenvalue of $P_u^\perp (A -\vartheta \Pi^\cS(uv^\top)) P_u^\perp$ is an eigenvector to the largest eigenvalue of $A-\vartheta \Pi^\cS(uv^\top)$ (unless $u$ is).
\qed
\end{proof}

\subsubsection*{Non-negativity constraint.}
Theorem~\ref{thm:ranking} does not take the non-negativity of the matrix into account. It solves the problem only if $A+\Delta\ge 0$, which may or may not happen. 

Algorithm~\ref{alg_rk-1} aims at reducing $|u^\top(A+\Delta)v|$ to zero while respecting the non-negativity constraint $A+\Delta\ge 0$ and retaining the commutativity of projections on subgraphs that exclude those matrix entries that cannot be further changed because of the non-negativity constraint. Note that the reduced edge sets $\cE_{k+1}\subset \cE_{k}$ are constructed in such a way that $\Pi^{\cE_{k+1}}$ still commutes with $P_u$.

When the algorithm finishes after finitely many steps, it returns, in the case of `success', a perturbation matrix $\Delta$ 
that satisfies $A+\Delta\ge 0$ and has the vector $v$ with $u^\top v=0$ of Theorem~\ref{thm:ranking} as an eigenvector to the largest eigenvalue (assuming it is not $u$). If the largest eigenvalue is simple, then $v$ must be a multiple of the non-negative Perron--Frobenius eigenvector of the non-negative matrix $A+\Delta$.
If the algorithm returns with `failure', it has reduced $|u^\top(A+\Delta)v|$ as far as is feasible with commuting projections.

\begin{algorithm}
\DontPrintSemicolon
\KwData{weighted adjacency matrix $A$, edge set $\cE$, vector $u$ of norm 1, \\
$\qquad v$ eigenvector to the largest eigenvalue of $P_u^\perp A P_u^\perp$}
\KwResult{perturbation matrix $\Delta$ in the cases of `success' and `failure'}
\Begin{
\nl Initialize $A_0=A$ and $\cE_0=\{(i,j)\in\overline\cE\,:\, u_i\ne 0 \; \text{ or } \, u_j \ne 0\}$ (see \eqref{E-closure})\;
\nl Initialize $\Theta_{0}=\Pi^{\cE_0}(\sym(uv^\top))$\;
\nl Initialize $k=0$ \ and \ flag = `ongoing'\;
\While{\text{\rm flag = `ongoing'}}{
\nl Set $\vartheta_k = \displaystyle\frac{u^\top A_k v}{\| \Theta_k \|_F^2}$\; 
\If{$A_k - \vartheta_k \Theta_k \ge 0$}{Set $\Delta=A_k - \vartheta_k \Theta_k - A$, flag = `success', \Return}
\nl Set $\rho_k = \max\{\rho\,:\, A_k - \rho\, \mathrm{sign}(\vartheta_k)\, \Theta_k \ge 0\}$\;
\nl Set $A_{k+1} = A_k - \rho_k\, \mathrm{sign}(\vartheta_k)\, \Theta_k$\;
\nl Set $\cE_{k+1}^0 = \{ (i,j)\in \cE_k\,:\, A_{k+1}\big|_{i,j}=0 \}$\;
\nl Set $\cE_{k+1} = \cE_k \setminus \overline{\cE_{k+1}^0}$ (see \eqref{E-closure})\;
\If{$\cE_{k+1}=\emptyset$}{Set $\Delta=A_{k+1}-A$, flag = `failure', \Return}
Set $\Theta_{k+1}=\Pi^{\cE_{k+1}}(\Theta_{k})$\;
Set $k=k+1$
}
}
\caption{Iterative reduction 
preserving non-negativity of weights and commutativity of projections.}
\label{alg_rk-1} 
\end{algorithm}

\section{Notes}

The usefulness of eigenvalues in studying graphs has been explored and amply demonstrated in spectral graph theory; see the book by Chung (\cite{Chu97}) and the introductory articles by Spielman (\cite{Spi07}) and von Luxburg (\cite{Lux07}). In this chapter we put forward
an approach that does not use spectral information of a matrix associated with the {\it given} graph (the graph Laplacian or the weighted adjacency matrix), but of such matrices for graphs with suitably perturbed weights. We reformulate classical problems in graph theory as matrix nearness problems related to eigenvalue optimization. We then adapt algorithms that we developed in previous chapters, using a two-level iteration with a constrained gradient system in the inner iteration. Such an approach to optimization problems for graphs was first taken by Andreotti, Edelmann, Guglielmi \& Lubich (\cite{AndEGL19}), for the problem of constrained graph partitioning. Robustness of spectral clustering to perturbations in the weight matrix was subsequently studied in a related setting by Andreotti, Edelmann, Guglielmi \& Lubich (\cite{AndEGL21}) and
Guglielmi \& Sicilia (\cite{GugS24}).

\subsubsection*{Constrained clustering.}
Constrained partitioning of graphs is a varied research area; see e.g. 
Kamvar, Sepandar, Klein, Dan, Manning, and Christopher (\cite{KamSKDMC03}),
Li, Liu \& Tang (\cite{LiLT09}),
Xu, Li \& Schuurmans (\cite{XuLS09}),
Wang \& Davidson (\cite{WanD10}), 
and B\"uhler, Rangapuram, Setzer \& Hein (\cite{BueRSH13}).
For the 
classical unconstrained minimum cut problem, there exist algorithms with complexity $O(n^2 \log n + n m)$ for graphs with $n$ nodes and $m$ edges; see Stoer \& Wagner (\cite{StoW97}) and references therein. However, cardinality constraints make the minimum cut problem NP-hard; see Bruglieri, Maffioli \& Ehrgott (\cite{BruME04}, \cite{BruEHM06}), who also propose a heuristic combinatorial algorithm for the cardinality-constrained minimum cut problem. Must-link and cannot-link constraints were imposed by Wagstaff, Cardie, Rogers \& Schr\"odl (\cite{WagCRS01}) and in many works thereafter; see e.g. Basu, Davidson, and Wagstaff (\cite{BasDW08}).

Andreotti, Edelmann, Guglielmi \& Lubich (\cite{AndEGL19}) formulated constrained graph partitioning as a matrix nearness problem and gave an iterative algorithm using eigenvalues and eigenvectors of the Laplacian of graphs with suitably perturbed weights. The new algorithm presented in Section~\ref{sec:cut} is similar in spirit but uses instead Rayleigh quotients, which are computationally significantly cheaper.

\subsubsection*{Robustness of spectral clustering.}
The spectral clustering algorithm, see e.g. the review by von Luxburg (\cite{Lux07}), is widely used for partitioning an undirected graph. Robustness of spectral clustering has been addressed from a statistical viewpoint by von Luxburg (\cite{Lux10}). In a different direction, asking for the smallest perturbation of the weight matrix on the given graph that changes the partition, yields a matrix nearness problem that has first been studied by Andreotti, Edelmann, Guglielmi \& Lubich (\cite{AndEGL21}). They devise a two-level iteration that uses a constrained gradient system in the inner iteration. This has been taken further by Guglielmi \& Sicilia (\cite{GugS24}), where the inner iteration works with projected {rank-1} perturbation matrices similarly to the rank-1 approach in Chapter~\ref{chap:struc}. We can directly refer to the algorithms of those two papers and therefore this topic is not included as a separate section in this chapter.

\subsubsection*{Robustness of ranking.}
Ranking the nodes of a graph via the leading eigenvector can be traced back to Landau~(\cite{Lan1895}). After being rediscovered many times, it became widely known through the Google PageRank algorithm; see e.g. Bryan \& Leise (\cite{BryL06}).
Assessing the robustness of a ranking can be viewed as a matrix nearness problem related to eigenvalue and eigenvector optimization. Based on this interpretation, an algorithm for quantifying the robustness of the ranking of nodes --- and simultaneously, for changing the ranking by a minimal perturbation --- in a given directed or undirected weighted graph was given by Benzi \& Guglielmi (\cite{BenG25}), minimizing a functional of the leading eigenvector. The algorithms for undirected graphs in Section~\ref{sec:ranking} are new and do not require to compute eigenvalues and eigenvectors and their derivatives but work instead with Rayleigh quotients.
Related problems have been considered in the literature by other authors.  Most relevant to our work are the papers by Cipolla, Durastante and Meini (\cite{CDM24}) and by Gillis and Van Dooren (\cite{GVD24}).

\subsubsection*{Topological stability of simplicial complexes.}
Finally, let us mention an application 
to topological stability of simplicial complexes with respect to perturbations of the underlying graph structure, where a gradient system approach (similar to that considered in this chapter) is applied to the Hodge Laplacian of a weighted simplicial complex to determine the minimal perturbation of edge weights sufficient to create an additional k-dimensional hole; see Guglielmi, Savostianov and Tudisco (\cite{GST23}).

\chapter{Appendix: Derivatives of eigenvalues and eigenvectors}
\label{chap:appendix}

In this appendix we put together basic results of first order perturbation
theory of eigenvalues and eigenvectors. We refer to Greenbaum, Li \& Overton (\cite{GreLO20}) for a recent review,
which also traces the rich history of this subject area.

\section{Derivative of paths of simple eigenvalues}

We often use the following standard perturbation result for
eigenvalues; see  e.g.  Horn \& Johnson (\cite{HJ90}), Lemma 6.3.10 and Theorem 6.3.12, and Greenbaum, Li \& Overton (\cite{GreLO20}), Theorem~1.

\begin{theorem}[Derivative of simple eigenvalues]\label{thm:eigderiv} 
Consider a continuously differentiable path of square complex matrices $A(t)$ for $t$ in an open interval $I$. Let $\lambda(t)$, $t\in I$, be a continuous path of simple eigenvalues of $A(t)$. Let $x(t)$ and $y(t)$ be left and right eigenvectors, respectively, of $A(t)$ to the eigenvalue $\lambda(t)$. Then, $x(t)^*y(t) \neq 0$ for $t\in I$ and $\lambda$ is continuously differentiable on $I$ with the derivative (denoted by a dot)
\begin{equation}
\dot{\lambda} = \frac{x^* \dot{A} y}{x^* y}\,.
\end{equation}
Moreover, ``continuously differentiable'' can be replaced with ``analytic'' in the assumption and the conclusion.
\end{theorem}
\index{derivative of!simple eigenvalue}
Since we have  $x(t)^*y(t) \neq 0$, we can apply the  normalization
\begin{equation}\label{eq:scalxy}
\| x(t) \|=1, \ \ \|y(t) \| =1, \quad x(t)^* y(t) \text{ is real and positive.}
\end{equation}
%
%
%
%
Clearly, a pair of left and right eigenvectors $\y$ and
$\x$ fulfilling \eqref{eq:scalxy} may be replaced by $\mu \x$ and $\mu \y$ for 
any complex $\mu$ of modulus $1$ without changing the property~\eqref{eq:scalxy}.

Next we turn to singular values. The following result is obtained from Theorem~\ref{thm:eigderiv} by using the
equivalence between singular values of $M$ and eigenvalues of $(0~M;
M^*~0)$; see Horn \& Johnson (\cite{HJ90}), Theorem 7.3.7.

\begin{corollary}[Derivative of singular values]
 \label{lem:singderiv} Consider a continuously differentiable path of matrices $M(t)\in\C^{m,n}$ for $t$ in an open interval $I$. Let $\sigma(t)$, $t\in I$, be a path of simple nonzero singular values of $M(t)$.  Let $u(t)$ and $v(t)$ be left and right singular vectors of $M(t)$ to the singular value $\sigma(t)$, that is, $M(t) v(t) = \sigma(t) u(t)$ and $u(t)^*M(t)= \sigma(t) v(t)^*$ with $\|u(t)\|=\|v(t)\|=1$.
 Then, $\sigma$ is continuously differentiable on $I$ with the derivative
$$
  \dot\sigma = \Re(u^* \dot M v).
$$
\end{corollary}
\index{derivative of!singular value}


\section{Derivative of paths of eigenvectors}
\label{sec:eigvec-deriv}

For the derivative of eigenvectors we need the notion of  group inverse (or reduced resolvent); see Meyer \& Stewart (\cite{MS88}) as well as Kato (\cite{Kat95}), Section I.5.3.

\begin{definition}[Group inverse]
\label{def:groupinv}
Let $N\in {\mathbb C}^{n,n}$ be a singular matrix with a simple zero eigenvalue.
The \emph{group inverse} (or reduced resolvent) of $N$ is the unique matrix
$Z$ with
\begin{equation}\label{group-inv-cond}
NZ =ZN, \qquad ZNZ=Z, \quad \mbox{and} \quad NZN=N.
\end{equation}
\end{definition}
\index{group inverse}
For the left and right eigenvectors $x$ and $y$ associated with the eigenvalue $0$ of $N$,
these equations yield the useful relations 
\begin{equation}\label{Zy-xZ}
Zy = ZNZy = Z^2Ny =0 \quad\text{ and analogously } \quad x^*Z =0.
\end{equation}
It is known from Meyer \& Stewart (\cite{MS88}) that if $N$ is normal, then
its group inverse $Z$ is equal to the better known Moore--Penrose pseudoinverse $N\pin$. In general, the two pseudoinverses are not the same. 
\index{Moore--Penrose pseudoinverse}
They are, however, related by the following result, which is a special case of a more general result in Appendix A of Guglielmi, Overton \& Stewart (\cite{GOS15}) but is also simply verified directly.

\begin{theorem} [Group inverse via Moore--Penrose pseudoinverse]\label{thm:Ginv}
Suppose that the matrix $N$ has the simple eigenvalue $0$ with corresponding left and right eigenvectors 
$\y$ and $\x$ of unit norm and such that $\y^* \x > 0$.
Let $Z$ be the group inverse of $N$, and with $\kappa=1/(\y^*\x)$ define the projection $\Pi=I-\kappa\x\y^*$.
 Then, the group inverse $Z$ of $N$ is related to the Moore--Penrose pseudoinverse $N\pin$ by
\begin{equation}
Z  =  \Pi N\pin \Pi .
\label{groupinvformula}
\end{equation}
\end{theorem}

The Moore-Penrose pseudo-inverse $N\pin$ is obtained from the singular value decomposition $N=U\Sigma V^*$ with unitary matrices $U$ and $V$ and the diagonal matrix $\Sigma=(\Sigma_+ \ 0; \,0\ 0)$, where $\Sigma_+$ is the diagonal matrix of the positive singular values. Then, $N\pin=V\Sigma\pin U^*$, where $\Sigma\pin=((\Sigma_+)^{-1}\ 0;\, 0\ 0)$. It is then a simple exercise to verify that
$\Pi N\pin \Pi $ satisfies the conditions  \eqref{group-inv-cond} that define  the group inverse $Z$. In particular, we find $NZ=ZN=\Pi$, which implies the other two conditions. 


\medskip
The group inverse appears in the following result on the derivative of eigenvectors, which is a variant of Theorem 2 of Meyer \& Stewart (\cite{MS88}).
\begin{theorem}[Derivative of eigenvectors] \label{thm:eigvecderiv}
Consider an analytic path of square complex matrices $A(t)$ for $t$ in an open interval $I$. Let $\lambda(t)$, $t\in I$, be a path of simple eigenvalues of $A(t)$. Then, there exists a continuously differentiable path of associated left and right eigenvectors $\y(t)$ and $\x(t)$, $t\in I$, which are of unit norm with $x^*(t)y(t)>0$ and satisfy the differential equations
\begin{equation}\label{deigvec}
\begin{array}{rcl}
\dot \y^*\!\!\! & = & \!\! - \y^* \dot A Z + \Re(\y^* \dot A Z \y)\y^*,
\\
\dot \x  & = &  \!\! - Z \dot A \x \ \ +  \Re(\x^* Z\dot A \x)\x,
\end{array}
\end{equation}
where $Z(t)$ is the group inverse of $A(t)-\lambda(t) I$.
\end{theorem}
\index{derivative of!eigenvector}

We note that the last terms on the right-hand sides of \eqref{deigvec} are in the direction of $\y^*$ and $\x$. They serve to ensure that the unit norm of $\y(t)$ and $\x(t)$ is conserved and that $\y(t)^*\x(t)$ remains real (and hence positive, since it cannot change sign by Theorem~\ref{thm:eigderiv}).
This is shown by verifying that
$(d/dt)(\y^*\x)$ is real, using the relations \eqref{Zy-xZ}, and that $(d/dt)\|\y\|^2=2\,\Re\, \dot \y^*\y=0\ $ and $\ (d/dt)\|\x\|^2=2\,\Re\, \x^*\dot\x=0\ $ when $x$ and $y$ are of norm~1.

In Theorem 2 of Meyer \& Stewart (\cite{MS88}), the last terms in \eqref{deigvec} appear without taking the real part. While this preserves the unit norm of the left and right eigenvectors, the positivity of their inner product is then not conserved.
Dropping the last terms (which are not analytic) altogether yields an analytic path of non-normalized left and right eigenvectors $\y^*(t)$ and $\x(t)$ with constant inner product $\y^*(t)\x(t)$; see Greenbaum, Li \& Overton (\cite{GreLO20}), Theorem~2 and Section~3.4.

The following result for singular vectors is obtained directly from Theorem~\ref{thm:eigvecderiv}, using that a left singular vector $u$ of a matrix $M$ is an eigenvector of the matrix $M M^*$ and 
that the pseudoinverse $\left( M M^* - \sigma^2 I \right)^\dagger$ has $u$ in its kernel.
The same argument holds true for the right singular vector $v$, which is an eigenvector of the matrix $M^* M$.

\begin{corollary}[Derivative of singular vectors]
\label{lem:singvecderiv} Consider a continuously differentiable path of matrices $M(t)\in\C^{m,n}$ for $t$ in an open interval $I$. 
Let $u(t)$ and $v(t)$ be left and right singular vectors of $M(t)$ to the simple nonzero singular value $\sigma(t)$, that is, $M(t) v(t) = \sigma(t) u(t)$ and $u(t)^*M(t)= \sigma(t) v(t)^*$ with $\|u(t)\|=\|v(t)\|=1$.
Then, $u$ and $v$ are continuously differentiable on $I$ with the derivative
\begin{equation}\label{dsingvec}
\begin{array}{rcl}
  \dot u & = & - \left( M M^* - \sigma^2 I \right)^\dagger \left( \frac{d}{d t} (M M^*) \right) u
\\
  \dot v & = & - \left( M^* M - \sigma^2 I \right)^\dagger \left( \frac{d}{d t} (M^* M) \right) v.
\end{array}
\end{equation}
\end{corollary}
\index{derivative of!singular vector}

\printindex

\end{document}